%% file: sm.tex
\documentclass[a4paper,12pt,twoside]{article}	

\usepackage[colorlinks=true,urlcolor=blue,menucolor=blue,linkcolor=blue]{hyperref} 
\usepackage{titletoc} 
\usepackage[inner=3cm,outer=2.5cm,top=1.5cm,bottom=2.5cm,includeheadfoot]{geometry}  
\usepackage{fancyhdr}
\usepackage[ngerman]{babel}
\usepackage{setspace}

\usepackage[latin1]{inputenc}
\usepackage{amssymb,amsmath,amsfonts}
\usepackage{theoremref} 
\usepackage{stmaryrd} 
\usepackage{scalerel} 

\DeclareFontFamily{U}{mathx}{\hyphenchar\font45}
\DeclareFontShape{U}{mathx}{m}{n}{
      <5> <6> <7> <8> <9> <10>
      <10.95> <12> <14.4> <17.28> <20.74> <24.88>
      mathx10
      }{}
\DeclareSymbolFont{mathx}{U}{mathx}{m}{n}
\DeclareFontSubstitution{U}{mathx}{m}{n}
\DeclareMathAccent{\widecheck}{0}{mathx}{"71}

\let\oldsqrt\sqrt
\def\sqrt{\mathpalette\DHLhksqrt}
\def\DHLhksqrt#1#2{%
\setbox0=\hbox{$#1\oldsqrt{#2\,}$}\dimen0=\ht0
\advance\dimen0-0.2\ht0
\setbox2=\hbox{\vrule height\ht0 depth -\dimen0}%
{\box0\lower0.4pt\box2}}

\usepackage[show]{ed} 
\usepackage{changebar} 
\usepackage{verbatim}	

\usepackage{color}  
\usepackage{graphicx}
\date{} 

\allowdisplaybreaks 

\usepackage{enumerate}
\usepackage[T1]{fontenc}

\begin{document}

\abovedisplayskip=2mm plus 0.5mm minus 0.5mm 
\belowdisplayskip=2mm plus 0.5mm minus 0.5mm 

\definecolor{refkey}{rgb}{1,0,0}
\definecolor{labelkey}{rgb}{1,0,0}
\numberwithin{equation}{section}

\input{macros}

\begin{titlepage}
\begin{center}
\vspace*{1cm}
{\Huge \bf Einige Beiträge zu vollständig nichtdegenerierten matriziellen Momentenproblemen vom $\alpha$-Stieltjes-Typ} \\
\vspace*{1.5cm}
{\Large Der Fakultät für Mathematik und Informatik} \\
{\Large der Universität Leipzig} \\
{\Large eingereichte} \\
\vspace*{1cm}
{\Large D I S S E R T A T I O N} \\
\vspace*{1cm}
{\Large zur Erlangung des akademischen Grades} \\
\vspace*{1cm}
{\Large DOCTOR RERUM NATURALIUM} \\
{\Large (Dr.\,rer.\,nat.)} \\
\vspace*{1cm}
{\Large im Fachgebiet} \\
\vspace*{1cm}
{\Large Mathematik} \\
\vspace*{1cm}
{\Large vorgelegt} \\
\vspace*{1cm}
{\Large von Diplom-Mathematiker Benjamin Jeschke} \\
{\Large geboren am 15. Februar 1988 in Leipzig} \\
\vspace*{1cm}
{\Large Leipzig, den 9. März 2017}
\end{center}
\end{titlepage}
\thispagestyle{empty}
\cleardoublepage

\newpage
\begin{center}
\vspace*{5cm}
Ich bedanke mich herzlichst bei meinem Betreuer Prof.\,Dr. Bernd Kirstein für die zahlreichen Hinweise und Anregungen sowie den ausführlichen historischen Hintergrund. Auch möchte ich mich bei Dr. Conrad Mädler für die Hilfestellung bei vielen mathematischen Fragen und Prof.\,Dr. Bernd Fritzsche für einige hilfreiche Anregungen bedanken. Ein besonderer Dank gilt meinen Eltern, die mich stets unterstützt haben.
\end{center}
\vfill
{\large Version 1.0.0}
\thispagestyle{empty}
\newpage
\thispagestyle{empty}
\hfill

\parskip=0.7mm	
\newpage
\tableofcontents

\parskip=2.5mm 
\newpage

\input{einleitung}
\input{sm0}
\input{sm1}
\input{sm2}
\input{sm3}
\input{sm4}
\input{sm5}

\appendix	
\input{anhang}

\newpage
\fancyhf{}
\fancyfoot[RO,LE]{\thepage}

\input{begriffe}
\input{symbole}
\input{literatur}

\end{document}

%% file: macros.tex
\newtheorem{satz}		{Satz}[section]
\newtheorem{lemma}		[satz]{Lemma}
\newtheorem{bem}		[satz]{Bemerkung}
\newtheorem{defi}		[satz]{Definition}
\newtheorem{theo}		[satz]{Theorem}
\newtheorem{folg}		[satz]{Folgerung}
\newtheorem{beispiel}	[satz]{Beispiel}
\newtheorem{bez}		[satz]{Bezeichnung}

\newcommand{\anf}[1]	{\glqq #1\grqq{}}

\newcommand{\eklam}[1]	{\left[ #1 \right]}
\newcommand{\beklam}[1]	{\big[ #1 \big]}
\newcommand{\Beklam}[1]	{\Big[ #1 \Big]}
\newcommand{\bbeklam}[1]	{\bigg[ #1 \bigg]}
\newcommand{\Bbeklam}[1]	{\Bigg[ #1 \Bigg]}
\newcommand{\gklam}[1]	{\left\{ #1 \right\}}
\newcommand{\bgklam}[1]	{\big\{ #1 \big\}}
\newcommand{\rklam}[1]	{\left( #1 \right)}
\newcommand{\brklam}[1]	{\big( #1 \big)}
\newcommand{\Brklam}[1]	{\Big( #1 \Big)}
\newcommand{\bbrklam}[1]	{\bigg( #1 \bigg)}
\newcommand{\Bbrklam}[1]	{\Bigg( #1 \Bigg)}
\newcommand{\sklam}[1]	{\stretchleftright{\langle}{#1}{\rangle}}
\newcommand{\fklam}[1]	{\raisebox{1pt}{\scaleobj{0.5}{\llfloor}} #1 \raisebox{1pt}{\scaleobj{0.5}{\rrfloor}}}
\newcommand{\fklamo}[1]	{\raisebox{1.5pt}{\scaleobj{0.75}{\llfloor}} #1 \raisebox{1.5pt}{\scaleobj{0.75}{\rrfloor}}}

\newcommand{\abs}[1]	{\left| #1 \right|}
\newcommand{\babs}[1]	{\big| #1 \big|}
\newcommand{\snorm}[1]	{\left|\left| \, #1 \, \right|\right|_{S}}
\newcommand{\bsnorm}[1]	{\big|\big| \, #1 \, \big|\big|_{S}}
\newcommand{\enorm}[1]	{\left|\left| \, #1 \, \right|\right|_{E}}
\newcommand{\benorm}[1]	{\big|\big| \, #1 \, \big|\big|_{E}}

\newcommand{\bbinom}[2]{\Brklam{\mkern-2mu\genfrac{}{}{0pt}{}{#1}{#2}\mkern-2mu}}

\newcommand{\fref}[1]	{(\ref{#1})}
\newcommand{\sref}[1]	{Seite \pageref{#1}}

\newcommand{\prodr}		{\operatorname*{\overrightarrow{\prod}}}
\newcommand{\prodl}		{\operatorname*{\overleftarrow{\prod}}}
\newcommand{\im}		{\operatorname{Im}}
\newcommand{\re}		{\operatorname{Re}}
\newcommand{\rank}		{\operatorname{rank}}
\newcommand{\tr}		{\operatorname{tr}}
\newcommand{\Rstr}		{\operatorname{Rstr.}}
\newcommand{\diag}		{\operatorname{diag}}

\newcommand{\qq}		{q\times q}
\newcommand{\pp}		{p\times p}
\newcommand{\pq}		{p\times q}
\newcommand{\qp}		{q\times p}
\newcommand{\qr}		{q\times r}
\newcommand{\pr}		{p\times r}
\newcommand{\rr}		{r\times r}

\newcommand{\N}			{\mathbb N}
\newcommand{\C}			{\mathbb C}
\newcommand{\R}			{\mathbb R}
\newcommand{\A}			{\mathfrak{A}}
\newcommand{\B}			{\mathfrak{B}}
\newcommand{\M}			{{\cal M}}
\newcommand{\D}			{{\cal D}}
\newcommand{\G}			{{\cal G}}
\newcommand{\E}			{{\cal E}}
\newcommand{\eps}		{\varepsilon}

\newcommand{\Z}[2]		{{\mathbb Z}_{#1,#2}}
\newcommand{\Zok}		{\Z{0}{\kappa}}
\newcommand{\Zokm}		{\Z{0}{\kappa-1}}
\newcommand{\Zokp}		{\Z{0}{\kappa+1}}
\newcommand{\Zon}		{\Z{0}{n}}
\newcommand{\Zonm}		{\Z{0}{n-1}}
\newcommand{\Zonp}		{\Z{0}{n+1}}
\newcommand{\Zozn}		{\Z{0}{2n}}
\newcommand{\Zom}		{\Z{0}{m}}
\newcommand{\Zomm}		{\Z{0}{m-1}}
\newcommand{\Zoq}		{\Z{0}{q}}
\newcommand{\Zekm}		{\Z{1}{\kappa-1}}
\newcommand{\Zen}		{\Z{1}{n}}
\newcommand{\Zenm}		{\Z{1}{n-1}}
\newcommand{\Zem}		{\Z{1}{m}}
\newcommand{\Zemm}		{\Z{1}{m-1}}
\newcommand{\Zek}		{\Z{1}{\kappa}}
\newcommand{\Zekp}		{\Z{1}{\kappa+1}}
\newcommand{\Zep}		{\Z{1}{p}}
\newcommand{\Zeq}		{\Z{1}{q}}
\newcommand{\Zofk}		{\Z{0}{\fklam{\kappa}}}
\newcommand{\Zofkm}		{\Z{0}{\fklam{\kappa-1}}}
\newcommand{\Zofkp}		{\Z{0}{\fklam{\kappa+1}}}
\newcommand{\Zefk}		{\Z{1}{\fklam{\kappa}}}
\newcommand{\Zefkm}		{\Z{1}{\fklam{\kappa-1}}}
\newcommand{\Zefkp}		{\Z{1}{\fklam{\kappa+1}}}
\newcommand{\Zzk}		{\Z{2}{\kappa}}
\newcommand{\Zzfk}		{\Z{2}{\fklam{\kappa}}}
\newcommand{\Zzfkm}		{\Z{2}{\fklam{\kappa-1}}}
\newcommand{\Zzfkp}		{\Z{2}{\fklam{\kappa+1}}}
\newcommand{\Zofm}		{\Z{0}{\fklam{m}}}

\newcommand{\Cp}		{\C^p}
\newcommand{\Cq}		{\C^q}
\newcommand{\Cqq}		{\C^{\qq}}
\newcommand{\Cpp}		{\C^{\pp}}
\newcommand{\Cpq}		{\C^{\pq}}
\newcommand{\Cqp}		{\C^{\qp}}
\newcommand{\Cqr}		{\C^{\qr}}
\newcommand{\Cpr}		{\C^{\pr}}
\newcommand{\Crr}		{\C^{\rr}}
\newcommand{\Crq}		{\C^{r\times q}}
\newcommand{\Cpqq}		{\C^{\qq}_{>}}
\newcommand{\CHqq}		{\C^{\qq}_H}

\newcommand{\Pp}		{\Pi_{+}}
\newcommand{\Pm}		{\Pi_{-}}
\newcommand{\Lap}		{\C_{\alpha,+}}
\newcommand{\Lam}		{\C_{\alpha,-}}

\newcommand{\za}		{\overline{z}}
\newcommand{\omegaa}	{\overline{\omega}}

\newcommand{\No}		{\N_{0}}
\newcommand{\Na}		{\overline{\N}}
\newcommand{\Noa}		{\Na_{0}}

\newcommand{\sj}		{(s_{j})_{j=0}^{\infty}}
\newcommand{\sjk}		{(s_{j})_{j=0}^{\kappa}}
\newcommand{\sjzk}		{(s_{j})_{j=0}^{2\kappa}}
\newcommand{\sjfk}		{(s_{j})^{2\fklam{\kappa}}_{j=0}}
\newcommand{\sjfkm}		{(s_{j})^{2\fklam{\kappa-1}}_{j=0}}
\newcommand{\sjm}		{(s_{j})_{j=0}^{m}}
\newcommand{\sjme}		{(s_{j})_{j=0}^{m+1}}
\newcommand{\sjl}		{(s_{j})_{j=0}^{l}}
\newcommand{\sjn}		{(s_{j})_{j=0}^{2n}}
\newcommand{\sjne}		{(s_{j})_{j=0}^{2n+1}}
\newcommand{\sjnm}		{(s_{j})_{j=0}^{2n-1}}
\newcommand{\sjnz}		{(s_{j})_{j=0}^{2(n+1)}}
\newcommand{\sjo}		{(s_{j})_{j=0}^{0}}

\newcommand{\tjk}		{(t_{j})_{j=0}^{\kappa}}
\newcommand{\tjm}		{(t_{j})_{j=0}^{m}}

\newcommand{\Mqo}		{\M^q_{\geq}(\Omega)}
\newcommand{\Mqco}		{\M^q_{\geq}(\widecheck{\Omega})}
\newcommand{\Mqoa}		{\M^q_{\geq}(\Omega,\A)}
\newcommand{\Mqoae}		{\M^q_{\geq}(\Omega_1,\A_1)}
\newcommand{\Mqoaz}		{\M^q_{\geq}(\Omega_2,\A_2)}
\newcommand{\Mqa}		{\M^q_{\geq}([\alpha,\infty))}
\newcommand{\Mqma}		{\M^q_{\geq}((-\infty,\alpha])}
\newcommand{\Mqmb}		{\M^q_{\geq}((-\infty,\beta])}
\newcommand{\Mqmma}		{\M^q_{\geq}((-\infty,-\alpha])}
\newcommand{\Mqab}		{\M^q_{\geq}([\alpha,\beta])}
\newcommand{\Mqko}		{\M^q_{\geq,\kappa}(\Omega)}
\newcommand{\Mqkco}		{\M^q_{\geq,\kappa}(\widecheck{\Omega})}
\newcommand{\MqR}		{\M^q_{\geq,\infty}(\R)}
\newcommand{\MqkR}		{\M^q_{\geq,\kappa}(\R)}
\newcommand{\Mqmo}		{\M^q_{\geq,m}(\Omega)}
\newcommand{\Mskg}		{\mathrm{M}\eklam{\Omega,\sjk,=}}
\newcommand{\Mqsg}		{\M^q_{\geq}\eklam{\Omega,\sj,=}}
\newcommand{\MqRsg}		{\M^q_{\geq}\eklam{\R,\sj,=}}
\newcommand{\Mqasg}		{\M^q_{\geq}\eklam{[\alpha,\infty),\sj,=}}
\newcommand{\Mqmasg}	{\M^q_{\geq}\eklam{(-\infty,\alpha],\sj,=}}
\newcommand{\Mqsmg}		{\M^q_{\geq}\eklam{\Omega,\sjm,=}}
\newcommand{\MqRsmg}	{\M^q_{\geq}\eklam{\R,\sjm,=}}
\newcommand{\Mqasmg}	{\M^q_{\geq}\eklam{[\alpha,\infty),\sjm,=}}
\newcommand{\Mqmasmg}	{\M^q_{\geq}\eklam{(-\infty,\alpha],\sjm,=}}
\newcommand{\Mqsnmg}	{\M^q_{\geq}\eklam{\Omega,\sjnm,=}}
\newcommand{\Mqmasnmg}	{\M^q_{\geq}\eklam{(-\infty,\alpha],\sjnm,=}}
\newcommand{\Mqskg}		{\M^q_{\geq}\eklam{\Omega,\sjk,=}}
\newcommand{\MqRskg}	{\M^q_{\geq}\eklam{\R,\sjk,=}}
\newcommand{\Mqaskg}	{\M^q_{\geq}\eklam{[\alpha,\infty),\sjk,=}}
\newcommand{\Mqmaskg}	{\M^q_{\geq}\eklam{(-\infty,\alpha],\sjk,=}}
\newcommand{\Mqmbskg}	{\M^q_{\geq}\eklam{(-\infty,\beta],\sjk,=}}
\newcommand{\Msmu}		{\mathrm{M}\eklam{\Omega,\sjm,\leq}}
\newcommand{\Mqsnu}		{\M^q_{\geq}\eklam{\Omega,\sjn,\leq}}
\newcommand{\MqRsnu}	{\M^q_{\geq}\eklam{\R,\sjn,\leq}}
\newcommand{\Mqasnu}	{\M^q_{\geq}\eklam{[\alpha,\infty),\sjn,\leq}}
\newcommand{\Mqmasnu}	{\M^q_{\geq}\eklam{(-\infty,\alpha],\sjn,\leq}}
\newcommand{\Mqmbsnu}	{\M^q_{\geq}\eklam{(-\infty,\beta],\sjn,\leq}}
\newcommand{\Mqmasnmu}	{\M^q_{\geq}\eklam{(-\infty,\alpha],\sjnm,\leq}}
\newcommand{\Mqsmu}		{\M^q_{\geq}\eklam{\Omega,\sjm,\leq}}
\newcommand{\Masmu}		{\mathrm{M}\eklam{[\alpha,\infty),\sjm,\leq}}
\newcommand{\Mqasmu}	{\M^q_{\geq}\eklam{[\alpha,\infty),\sjm,\leq}}
\newcommand{\Mmasmu}	{\mathrm{M}\eklam{(-\infty,\alpha],\sjm,\leq}}
\newcommand{\Mqmasmu}	{\M^q_{\geq}\eklam{(-\infty,\alpha],\sjm,\leq}}
\newcommand{\Mqsku}		{\M^q_{\geq}\eklam{\Omega,\sjk,\leq}}
\newcommand{\MqRsku}	{\M^q_{\geq}\eklam{\R,\sjk,\leq}}
\newcommand{\Msmuu}		{\mathrm{M}\eklam{\Omega,\sjm,\geq}}
\newcommand{\Mqsmuu}	{\M^q_{\geq}\eklam{\Omega,\sjm,\geq}}
\newcommand{\Mqasneu}	{\M^q_{\geq}\eklam{[\alpha,\infty),\sjne,\leq}}
\newcommand{\Mqmasneuu}	{\M^q_{\geq}\eklam{(-\infty,\alpha],\sjne,\geq}}
\newcommand{\Mqmbsneuu}	{\M^q_{\geq}\eklam{(-\infty,\beta],\sjne,\geq}}
\newcommand{\Masmuu}	{\mathrm{M}\eklam{[\alpha,\infty),\sjm,\geq}}
\newcommand{\Mmasmuu}	{\mathrm{M}\eklam{(-\infty,\alpha],\sjm,\geq}}
\newcommand{\Mqmasmuu}	{\M^q_{\geq}\eklam{(-\infty,\alpha],\sjm,\geq}}
\newcommand{\Mabskg}	{\mathrm{M}\eklam{[\alpha,\beta],\sjk,=}}
\newcommand{\Mqabskg}	{\M^q_{\geq}\eklam{[\alpha,\beta],\sjk,=}}
\newcommand{\Mqabsng}	{\M^q_{\geq}\eklam{[\alpha,\beta],\sjn,=}}
\newcommand{\Mqabsneg}	{\M^q_{\geq}\eklam{[\alpha,\beta],\sjne,=}}
\newcommand{\Mqabsog}	{\M^q_{\geq}\eklam{[\alpha,\beta],\sjo,=}}

\newcommand{\Sqa}		{{\cal S}_{q,[\alpha,\infty)}}
\newcommand{\Soqa}		{{\cal S}_{0,q,[\alpha,\infty)}}
\newcommand{\Saskg}		{S\eklam{[\alpha,\infty),\sjk,=}}
\newcommand{\Sqaskg}	{\Soqa[\sjk,=]}
\newcommand{\Sasmu}		{S\eklam{[\alpha,\infty),\sjm,\leq}}
\newcommand{\Sqasmu}	{\Soqa[\sjm,\leq]}
\newcommand{\Sqma}		{{\cal S}_{q,(-\infty,\alpha]}}
\newcommand{\Sqmma}		{{\cal S}_{q,(-\infty,-\alpha]}}
\newcommand{\Sqmb}		{{\cal S}_{q,(-\infty,\beta]}}
\newcommand{\Soqma}		{{\cal S}_{0,q,(-\infty,\alpha]}}
\newcommand{\Soqmma}	{{\cal S}_{0,q,(-\infty,-\alpha]}}
\newcommand{\Soqmb}		{{\cal S}_{0,q,(-\infty,\beta]}}
\newcommand{\Smaskg}	{S\eklam{(-\infty,\alpha],\sjk,=}}
\newcommand{\Sqmaskg}	{\Soqma[\sjk,=]}
\newcommand{\Smasmu}	{S\eklam{(-\infty,\alpha],\sjm,\leq}}
\newcommand{\Sqmasmu}	{\Soqma[\sjm,\leq]}
\newcommand{\Sabskg}	{S\eklam{[\alpha,\beta],\sjk,=}}

\newcommand{\Hq}		{{\cal H}^{\geq}_{q,\infty}}
\newcommand{\Hqo}		{{\cal H}^{\geq}_{q,0}}
\newcommand{\Hqn}		{{\cal H}^{\geq}_{q,2n}}
\newcommand{\Hqne}		{{\cal H}^{\geq}_{q,2(n+1)}}
\newcommand{\Hqnm}		{{\cal H}^{\geq}_{q,2(n-1)}}
\newcommand{\Hqk}		{{\cal H}^{\geq}_{q,2\kappa}}
\newcommand{\Hpq}		{{\cal H}^{>}_{q,\infty}}
\newcommand{\Hpqo}		{{\cal H}^{>}_{q,0}}
\newcommand{\Hpqk}		{{\cal H}^{>}_{q,2\kappa}}
\newcommand{\Hpqfk}		{{\cal H}^{>}_{q,2\fklam{\kappa}}}
\newcommand{\Hpqfkm}	{{\cal H}^{>}_{q,2\fklam{\kappa-1}}}
\newcommand{\Hpqfkzm}	{{\cal H}^{>}_{q,2\fklam{\kappa-2}}}
\newcommand{\Hpqn}		{{\cal H}^{>}_{q,2n}}
\newcommand{\Hpqne}		{{\cal H}^{>}_{q,2n+1}}
\newcommand{\Hpqnm}		{{\cal H}^{>}_{q,2(n-1)}}
\newcommand{\Heqm}		{{\cal H}^{\geq,e}_{q,m}}
\newcommand{\Heqn}		{{\cal H}^{\geq,e}_{q,2n}}
\newcommand{\Heqne}		{{\cal H}^{\geq,e}_{q,2n+1}}
\newcommand{\Heqnm}		{{\cal H}^{\geq,e}_{q,2n-1}}

\newcommand{\Kqa}		{{\cal K}^{\geq}_{q,\infty,\alpha}}
\newcommand{\Kqoa}		{{\cal K}^{\geq}_{q,0,\alpha}}
\newcommand{\Kqna}		{{\cal K}^{\geq}_{q,2n,\alpha}}
\newcommand{\Kqnea}		{{\cal K}^{\geq}_{q,2n+1,\alpha}}
\newcommand{\Kqma}		{{\cal K}^{\geq}_{q,m,\alpha}}
\newcommand{\Kqmea}		{{\cal K}^{\geq}_{q,m+1,\alpha}}
\newcommand{\Kqla}		{{\cal K}^{\geq}_{q,l,\alpha}}
\newcommand{\Kqka}		{{\cal K}^{\geq}_{q,\kappa,\alpha}}
\newcommand{\Kpqa}		{{\cal K}^{>}_{q,\infty,\alpha}}
\newcommand{\Kpqoa}		{{\cal K}^{>}_{q,0,\alpha}}
\newcommand{\Kpqna}		{{\cal K}^{>}_{q,2n,\alpha}}
\newcommand{\Kpqnea}	{{\cal K}^{>}_{q,2n+1,\alpha}}
\newcommand{\Kpqma}		{{\cal K}^{>}_{q,m,\alpha}}
\newcommand{\Kpqmma}	{{\cal K}^{>}_{q,m,-\alpha}}
\newcommand{\Kpqla}		{{\cal K}^{>}_{q,l,\alpha}}
\newcommand{\Kpqka}		{{\cal K}^{>}_{q,\kappa,\alpha}}
\newcommand{\Kpqkma}	{{\cal K}^{>}_{q,\kappa,-\alpha}}
\newcommand{\Keqma}		{{\cal K}^{\geq,e}_{q,m,\alpha}}
\newcommand{\Keqla}		{{\cal K}^{\geq,e}_{q,l,\alpha}}
\newcommand{\Keqka}		{{\cal K}^{\geq,e}_{q,\kappa,\alpha}}

\newcommand{\Lqa}		{{\cal L}^{\geq}_{q,\infty,\alpha}}
\newcommand{\Lqoa}		{{\cal L}^{\geq}_{q,0,\alpha}}
\newcommand{\Lqna}		{{\cal L}^{\geq}_{q,2n,\alpha}}
\newcommand{\Lqnea}		{{\cal L}^{\geq}_{q,2n+1,\alpha}}
\newcommand{\Lqma}		{{\cal L}^{\geq}_{q,m,\alpha}}
\newcommand{\Lqmea}		{{\cal L}^{\geq}_{q,m+1,\alpha}}
\newcommand{\Lqla}		{{\cal L}^{\geq}_{q,l,\alpha}}
\newcommand{\Lqka}		{{\cal L}^{\geq}_{q,\kappa,\alpha}}
\newcommand{\Lpqa}		{{\cal L}^{>}_{q,\infty,\alpha}}
\newcommand{\Lpqoa}		{{\cal L}^{>}_{q,0,\alpha}}
\newcommand{\Lpqna}		{{\cal L}^{>}_{q,2n,\alpha}}
\newcommand{\Lpqnea}	{{\cal L}^{>}_{q,2n+1,\alpha}}
\newcommand{\Lpqma}		{{\cal L}^{>}_{q,m,\alpha}}
\newcommand{\Lpqmma}		{{\cal L}^{>}_{q,m,-\alpha}}
\newcommand{\Lpqla}		{{\cal L}^{>}_{q,l,\alpha}}
\newcommand{\Lpqka}		{{\cal L}^{>}_{q,\kappa,\alpha}}
\newcommand{\Lpqkma}	{{\cal L}^{>}_{q,\kappa,-\alpha}}
\newcommand{\Leqma}		{{\cal L}^{\geq,e}_{q,m,\alpha}}
\newcommand{\Leqnma}	{{\cal L}^{\geq,e}_{q,2n-1,\alpha}}
\newcommand{\Leqla}		{{\cal L}^{\geq,e}_{q,l,\alpha}}

\newcommand{\Iq}		{I_{q}}
\newcommand{\Ip}		{I_{p}}
\newcommand{\Ir}		{I_{r}}
\newcommand{\Inq}		{I_{nq}}
\newcommand{\Inpq}		{I_{(n+1)q}}
\newcommand{\Ijq}		{I_{jq}}
\newcommand{\Ikq}		{I_{kq}}
\newcommand{\Izq}		{I_{2q}}
\newcommand{\Oq}		{0_{\qq}}
\newcommand{\Op}		{0_{\pp}}
\newcommand{\Opq}		{0_{\pq}}
\newcommand{\Onq}		{0_{nq \times q}}
\newcommand{\Oqn}		{0_{q \times nq}}
\newcommand{\Ojq}		{0_{jq \times q}}
\newcommand{\Oqj}		{0_{q \times jq}}
\newcommand{\Okq}		{0_{kq \times q}}
\newcommand{\Oqk}		{0_{q \times kq}}
\newcommand{\Onqn}		{0_{nq \times nq}}

\newcommand{\Hn}		{H_{n}}
\newcommand{\Hk}		{H_{k}}
\newcommand{\Hm}		{H_{m}}
\newcommand{\Hj}		{H_{j}}
\newcommand{\Hnm}		{H_{n-1}}
\newcommand{\Hjm}		{H_{j-1}}
\newcommand{\Hsn}		{\Hn^{\sklam{s}}}
\newcommand{\Hsnm}		{\Hnm^{\sklam{s}}}
\newcommand{\Htn}		{\Hn^{\sklam{t}}}
\newcommand{\dH}		{\widehat{H}}
\newcommand{\dHn}		{\dH_{n}}
\newcommand{\dHnm}		{\dH_{n-1}}
\newcommand{\dHsn}		{\dH^{\sklam{s}}_{n}}
\newcommand{\dHsnm}		{\dH^{\sklam{s}}_{n-1}}
\newcommand{\tH}		{\widetilde{H}}
\newcommand{\tHn}		{\tH_{n}}
\newcommand{\Kn}		{K_{n}}
\newcommand{\Knm}		{K_{n-1}}
\newcommand{\Ksn}		{\Kn^{\sklam{s}}}
\newcommand{\Ktn}		{\Kn^{\sklam{t}}}
\newcommand{\Ksnm}		{\Knm^{\sklam{s}}}
\newcommand{\dK}		{\widehat{K}}
\newcommand{\dKn}		{\dK_{n}}
\newcommand{\dKsn}		{\dK^{\sklam{s}}_{n}}
\newcommand{\wKn}		{\widetilde{K}_{n}}
\newcommand{\wKnm}		{\widetilde{K}_{n-1}}
\newcommand{\wKsn}		{\wKn^{\sklam{s}}}
\newcommand{\yjk}		{y_{j,k}}
\newcommand{\ysjk}		{\yjk^{\sklam{s}}}
\newcommand{\zjk}		{z_{j,k}}
\newcommand{\zsjk}		{\zjk^{\sklam{s}}}

\newcommand{\Vo}		{V_0}
\newcommand{\Ve}		{V_1}
\newcommand{\Vea}		{\Ve^{\ast}}
\newcommand{\Vn}		{V_n}
\newcommand{\Vna}		{\Vn^{\ast}}
\newcommand{\Vnm}		{V_{n-1}}
\newcommand{\Vnma}		{\Vnm^{\ast}}

\newcommand{\ar}		{\alpha \triangleright}
\newcommand{\aro}		{\ar 0}
\newcommand{\arj}		{\ar j}
\newcommand{\arn}		{\ar n}
\newcommand{\arm}		{\ar m}
\newcommand{\ark}		{\ar k}
\newcommand{\arl}		{\ar l}
\newcommand{\sarjk}		{(s_{\arj})_{j=0}^{\kappa-1}}
\newcommand{\Harn}		{H_{\arn}}
\newcommand{\Haro}		{H_{\aro}}
\newcommand{\Harj}		{H_{\arj}}
\newcommand{\Hark}		{H_{\ark}}
\newcommand{\Harl}		{H_{\arl}}
\newcommand{\Hsarn}		{\Harn^{\sklam{s}}}
\newcommand{\Harnm}		{H_{\arn-1}}
\newcommand{\Hsarnm}	{\Harnm^{\sklam{s}}}
\newcommand{\Harnp}		{H_{\arn+1}}
\newcommand{\dHaro}		{\dH_{\aro}}
\newcommand{\dHarn}		{\dH_{\arn}}
\newcommand{\dHsarn}	{\dHarn^{\sklam{s}}}
\newcommand{\tHarn}		{\tH_{\arn}}
\newcommand{\Karn}		{K_{\arn}}
\newcommand{\Ksarn}		{\Karn^{\sklam{s}}}
\newcommand{\dKarn}		{\dK_{\arn}}
\newcommand{\dKsarn}	{\dKarn^{\sklam{s}}}
\newcommand{\yarjk}		{y_{\arj,k}}
\newcommand{\ysarjk}	{\yarjk^{\sklam{s}}}
\newcommand{\zarjk}		{z_{\arj,k}}
\newcommand{\zsarjk}	{\zarjk^{\sklam{s}}}

\newcommand{\al}		{\alpha \triangleleft}
\newcommand{\alo}		{\al 0}
\newcommand{\alj}		{\al j}
\newcommand{\alk}		{\al k}
\newcommand{\aln}		{\al n}
\newcommand{\alm}		{\al m}
\newcommand{\bl}		{\beta \triangleleft}
\newcommand{\bln}		{\bl n}
\newcommand{\saljk}		{(s_{\alj})_{j=0}^{\kappa-1}}
\newcommand{\Haln}		{H_{\aln}}
\newcommand{\Hbln}		{H_{\bln}}
\newcommand{\Halo}		{H_{\alo}}
\newcommand{\Hmaln}		{H_{-\aln}}
\newcommand{\Halj}		{H_{\alj}}
\newcommand{\Hsaln}		{\Haln^{\sklam{s}}}
\newcommand{\Halnm}		{H_{\aln-1}}
\newcommand{\Hsalnm}	{\Halnm^{\sklam{s}}}
\newcommand{\dHaln}		{\dH_{\aln}}
\newcommand{\dHmaln}	{\dH_{-\aln}}
\newcommand{\dHsaln}	{\dHaln^{\sklam{s}}}
\newcommand{\Kaln}		{K_{\aln}}
\newcommand{\Ksaln}		{\Kaln^{\sklam{s}}}
\newcommand{\dKaln}		{\dK_{\aln}}
\newcommand{\dKsaln}	{\dKaln^{\sklam{s}}}
\newcommand{\yaljk}		{y_{\alj,k}}
\newcommand{\ysaljk}	{\yaljk^{\sklam{s}}}
\newcommand{\zaljk}		{z_{\alj,k}}
\newcommand{\zsaljk}	{\zaljk^{\sklam{s}}}

\newcommand{\Qarj}		{(Q_{\arj})_{j=0}^{\infty}}
\newcommand{\Qarjk}		{(Q_{\arj})_{j=0}^{\kappa}}
\newcommand{\Qsarjk}	{(Q^{\sklam{s}}_{\arj})_{j=0}^{\kappa}}
\newcommand{\Qarjm}		{(Q_{\arj})_{j=0}^{m}}
\newcommand{\Qalj}		{(Q_{\alj})_{j=0}^{\infty}}
\newcommand{\Qaljk}		{(Q_{\alj})_{j=0}^{\kappa}}
\newcommand{\Qsaljk}	{(Q^{\sklam{s}}_{\alj})_{j=0}^{\kappa}}
\newcommand{\Qtmaljk}	{(Q^{\sklam{t}}_{-\alj})_{j=0}^{\kappa}}
\newcommand{\Qaljm}		{(Q_{\alj})_{j=0}^{m}}

\newcommand{\En}		{E_{n}}
\newcommand{\Ena}		{\En^{\ast}}
\newcommand{\Enm}		{E_{n-1}}
\newcommand{\Enma}		{\Enm^{\ast}}
\newcommand{\Ek}		{E_{k}}
\newcommand{\Ej}		{E_{j}}
\newcommand{\Ejm}		{E_{j-1}}

\newcommand{\Maro}		{\mathbf{M}_{\ar 0}}
\newcommand{\Mare}		{\mathbf{M}_{\ar 1}}
\newcommand{\Marz}		{\mathbf{M}_{\ar 2}}
\newcommand{\Marn}		{\mathbf{M}_{\ar n}}
\newcommand{\Marnm}		{\mathbf{M}_{\ar n-1}}
\newcommand{\Marnp}		{\mathbf{M}_{\ar n+1}}
\newcommand{\Marj}		{\mathbf{M}_{\ar j}}
\newcommand{\Marjp}		{\mathbf{M}_{\ar j+1}}
\newcommand{\Mmaro}		{\mathbf{M}_{-\ar 0}}
\newcommand{\Mmare}		{\mathbf{M}_{-\ar 1}}
\newcommand{\Mmarn}		{\mathbf{M}_{-\ar n}}
\newcommand{\Mmarnm}	{\mathbf{M}_{-\ar n-1}}
\newcommand{\Mmarnp}	{\mathbf{M}_{-\ar n+1}}
\newcommand{\Mmarj}		{\mathbf{M}_{-\ar j}}
\newcommand{\Laro}		{\mathbf{L}_{\ar 0}}
\newcommand{\Lare}		{\mathbf{L}_{\ar 1}}
\newcommand{\Larn}		{\mathbf{L}_{\ar n}}
\newcommand{\Larnm}		{\mathbf{L}_{\ar n-1}}
\newcommand{\Larj}		{\mathbf{L}_{\ar j}}
\newcommand{\Lmaro}		{\mathbf{L}_{-\ar 0}}
\newcommand{\Lmare}		{\mathbf{L}_{-\ar 1}}
\newcommand{\Lmarn}		{\mathbf{L}_{-\ar n}}
\newcommand{\Lmarnm}		{\mathbf{L}_{-\ar n-1}}
\newcommand{\Lmarj}		{\mathbf{L}_{-\ar j}}
\newcommand{\LMarn}		{[(\Larn)^{\infty}_{n=0}$,""$(\Marn)^{\infty}_{n=0}]}
\newcommand{\LMkarn}	{[(\Larn)^{\fklam{\kappa-1}}_{n=0}$,""$(\Marn)^{\fklam{\kappa}}_{n=0}]}
\newcommand{\LMskarn}	{[(\Larn^{\sklam{s}})^{\fklam{\kappa-1}}_{n=0}$,""$(\Marn^{\sklam{s}})^{\fklam{\kappa}}_{n=0}]}

\newcommand{\Malo}		{\mathbf{M}_{\al 0}}
\newcommand{\Male}		{\mathbf{M}_{\al 1}}
\newcommand{\Maln}		{\mathbf{M}_{\al n}}
\newcommand{\Malnm}		{\mathbf{M}_{\al n-1}}
\newcommand{\Malnp}		{\mathbf{M}_{\al n+1}}
\newcommand{\Malj}		{\mathbf{M}_{\al j}}
\newcommand{\Maljp}		{\mathbf{M}_{\al j+1}}
\newcommand{\Mmalo}		{\mathbf{M}_{-\al 0}}
\newcommand{\Mmaln}		{\mathbf{M}_{-\al n}}
\newcommand{\Lalo}		{\mathbf{L}_{\al 0}}
\newcommand{\Lale}		{\mathbf{L}_{\al 1}}
\newcommand{\Laln}		{\mathbf{L}_{\al n}}
\newcommand{\Lalnm}		{\mathbf{L}_{\al n-1}}
\newcommand{\Lalj}		{\mathbf{L}_{\al j}}
\newcommand{\Lmalo}		{\mathbf{L}_{-\al 0}}
\newcommand{\Lmaln}		{\mathbf{L}_{-\al n}}
\newcommand{\LMaln}		{[(\Laln)^{\infty}_{n=0}$,""$(\Maln)^{\infty}_{n=0}]}
\newcommand{\LMkaln}	{[(\Laln)^{\fklam{\kappa-1}}_{n=0}$,""$(\Maln)^{\fklam{\kappa}}_{n=0}]}
\newcommand{\LMskaln}	{[(\Laln^{\sklam{s}})^{\fklam{\kappa-1}}_{n=0}$,""$(\Maln^{\sklam{s}})^{\fklam{\kappa}}_{n=0}]}
\newcommand{\LMtkmaln}	{[(\Lmaln^{\sklam{t}})^{\fklam{\kappa-1}}_{n=0}$,""$(\Mmaln^{\sklam{t}})^{\fklam{\kappa}}_{n=0}]}

\newcommand{\Cnk}		{(C_{n})^{\kappa}_{n=1}}
\newcommand{\Cnfk}		{(C_{n})^{\fklam{\kappa+1}}_{n=1}}
\newcommand{\Csnfk}		{(C^{\sklam{s}}_{n})^{\fklam{\kappa+1}}_{n=1}}
\newcommand{\Dnk}		{(D_{n})^{\kappa}_{n=0}}
\newcommand{\Dnfk}		{(D_{n})^{\fklam{\kappa}}_{n=0}}
\newcommand{\Dsnfk}		{(D^{\sklam{s}}_{n})^{\fklam{\kappa}}_{n=0}}
\newcommand{\CDk}		{[\Cnk$,""$\Dnk]}
\newcommand{\CDfk}		{[\Cnfk$,""$\Dnfk]}
\newcommand{\CDsfk}		{[\Csnfk$,""$\Dsnfk]}

\newcommand{\Anfk}		{(A_{n})^{\fklam{\kappa-1}}_{n=0}}
\newcommand{\Bnfk}		{(B_{n})^{\fklam{\kappa}}_{n=0}}
\newcommand{\Asnfk}		{(A^{\sklam{s}}_{n})^{\fklam{\kappa-1}}_{n=0}}
\newcommand{\Bsnfk}		{(B^{\sklam{s}}_{n})^{\fklam{\kappa}}_{n=0}}
\newcommand{\ABfk}		{[\Anfk$,""$\Bnfk]}
\newcommand{\ABsfk}		{[\Asnfk$,""$\Bsnfk]}
\newcommand{\Aarnfk}	{(A_{\arn})^{\fklam{\kappa-2}}_{n=0}}
\newcommand{\Barnfk}	{(B_{\arn})^{\fklam{\kappa-1}}_{n=0}}
\newcommand{\Asarnfk}	{(A^{\sklam{s}}_{\arn})^{\fklam{\kappa-2}}_{n=0}}
\newcommand{\Bsarnfk}	{(B^{\sklam{s}}_{\arn})^{\fklam{\kappa-1}}_{n=0}}
\newcommand{\ABarfk}	{[\Aarnfk$,""$\Barnfk]}
\newcommand{\ABsarfk}	{[\Asarnfk$,""$\Bsarnfk]}
\newcommand{\Aalnfk}	{(A_{\aln})^{\fklam{\kappa-2}}_{n=0}}
\newcommand{\Balnfk}	{(B_{\aln})^{\fklam{\kappa-1}}_{n=0}}
\newcommand{\ABalfk}	{[\Aalnfk$,""$\Balnfk]}

\newcommand{\Pn}		{(P_{n})^\infty_{n=0}}
\newcommand{\Pnk}		{(P_{n})^\kappa_{n=0}}
\newcommand{\Pnfk}		{(P_{n})^{\fklam{\kappa+1}}_{n=0}}
\newcommand{\Parnfk}	{(P_{\arn})^{\fklam{\kappa}}_{n=0}}
\newcommand{\Palnfk}	{(P_{\aln})^{\fklam{\kappa}}_{n=0}}
\newcommand{\Ps}		{P^{\sklam{s}}}
\newcommand{\Pars}		{P^{\sklam{\ar s}}}
\newcommand{\Pals}		{P^{\sklam{\al s}}}
\newcommand{\Psn}		{(\Ps_{n})^\infty_{n=0}}
\newcommand{\Psnk}		{(\Ps_{n})^\kappa_{n=0}}
\newcommand{\Psnfk}		{(\Ps_{n})^{\fklam{\kappa+1}}_{n=0}}
\newcommand{\Psarnfk}	{(\Pars_{\arn})^{\fklam{\kappa}}_{n=0}}
\newcommand{\Psalnfk}	{(\Pals_{\aln})^{\fklam{\kappa}}_{n=0}}
\newcommand{\dP}		{\widehat{P}}
\newcommand{\dPs}		{\dP^{\sklam{s}}}
\newcommand{\dPsarnfk}	{(\dPs_{\arn})^{\fklam{\kappa}}_{n=0}}
\newcommand{\dPsalnfk}	{(\dPs_{\aln})^{\fklam{\kappa}}_{n=0}}
\newcommand{\SQark}		{[\Pnfk$,""$\Psnfk$,""$\Parnfk$,""$\dPsarnfk]}
\newcommand{\SQsark}	{[(P_{n,s})^{\fklam{\kappa+1}}_{n=0}$,""$(\Ps_{n,s})^{\fklam{\kappa+1}}_{n=0}$,""$(P_{\arn,s})^{\fklam{\kappa}}_{n=0}$,""$(\dPs_{\arn,s})^{\fklam{\kappa}}_{n=0}]}
\newcommand{\SQalk}		{[\Pnfk$,""$\Psnfk$,""$\Palnfk$,""$\dPsalnfk]}
\newcommand{\SQsalk}	{[(P_{n,s})^{\fklam{\kappa+1}}_{n=0}$,""$(\Ps_{n,s})^{\fklam{\kappa+1}}_{n=0}$,""$(P_{\aln,s})^{\fklam{\kappa}}_{n=0}$,""$(\dPs_{\aln,s})^{\fklam{\kappa}}_{n=0}]}

\newcommand{\Aaro}		{\mathbf{A}_{\aro}}
\newcommand{\Asaro}		{\Aaro^{\sklam{s}}}
\newcommand{\Aare}		{\mathbf{A}_{\ar1}}
\newcommand{\Aarn}		{\mathbf{A}_{\arn}}
\newcommand{\Asarn}		{\Aarn^{\sklam{s}}}
\newcommand{\Aarnm}		{\mathbf{A}_{\arn-1}}
\newcommand{\Aarfm}		{\mathbf{A}_{\ar\fklam{m}}}
\newcommand{\Asarfm}	{\Aarfm^{\sklam{s}}}
\newcommand{\Baro}		{\mathbf{B}_{\aro}}
\newcommand{\Bsaro}		{\Baro^{\sklam{s}}}
\newcommand{\Bare}		{\mathbf{B}_{\ar1}}
\newcommand{\Barn}		{\mathbf{B}_{\arn}}
\newcommand{\Bsarn}		{\Barn^{\sklam{s}}}
\newcommand{\Barnm}		{\mathbf{B}_{\arn-1}}
\newcommand{\Barnp}		{\mathbf{B}_{\arn+1}}
\newcommand{\Barfm}		{\mathbf{B}_{\ar\fklam{m+1}}}
\newcommand{\Bsarfm}	{\Barfm^{\sklam{s}}}
\newcommand{\Caro}		{\mathbf{C}_{\aro}}
\newcommand{\Csaro}		{\Caro^{\sklam{s}}}
\newcommand{\Care}		{\mathbf{C}_{\ar1}}
\newcommand{\Carn}		{\mathbf{C}_{\arn}}
\newcommand{\Csarn}		{\Carn^{\sklam{s}}}
\newcommand{\Carnm}		{\mathbf{C}_{\arn-1}}
\newcommand{\Carfm}		{\mathbf{C}_{\ar\fklam{m}}}
\newcommand{\Csarfm}	{\Carfm^{\sklam{s}}}
\newcommand{\Daro}		{\mathbf{D}_{\aro}}
\newcommand{\Dsaro}		{\Daro^{\sklam{s}}}
\newcommand{\Dare}		{\mathbf{D}_{\ar1}}
\newcommand{\Darn}		{\mathbf{D}_{\arn}}
\newcommand{\Dsarn}		{\Darn^{\sklam{s}}}
\newcommand{\Darnm}		{\mathbf{D}_{\arn-1}}
\newcommand{\Darnp}		{\mathbf{D}_{\arn+1}}
\newcommand{\Darfm}		{\mathbf{D}_{\ar\fklam{m+1}}}
\newcommand{\Dsarfm}	{\Darfm^{\sklam{s}}}
\newcommand{\ABCDark}	{[(\Aarn)^{\fklam{\kappa}}_{n=0}$,""$(\Barn)^{\fklam{\kappa+1}}_{n=0}$,""$(\Carn)^{\fklam{\kappa}}_{n=0}$,""$(\Darn)^{\fklam{\kappa+1}}_{n=0}]}
\newcommand{\ABCDsark}	{[(\Asarn)^{\fklam{\kappa}}_{n=0}$,""$(\Bsarn)^{\fklam{\kappa+1}}_{n=0}$,""$(\Csarn)^{\fklam{\kappa}}_{n=0}$,""$(\Dsarn)^{\fklam{\kappa+1}}_{n=0}]}
\newcommand{\ABCDarm}	{[(\Aarn)^{\fklam{m}}_{n=0}$,""$(\Barn)^{\fklam{m+1}}_{n=0}$,""$(\Carn)^{\fklam{m}}_{n=0}$,""$(\Darn)^{\fklam{m+1}}_{n=0}]}
\newcommand{\ABCDsarm}	{[(\Asarn)^{\fklam{m}}_{n=0}$,""$(\Bsarn)^{\fklam{m+1}}_{n=0}$,""$(\Csarn)^{\fklam{m}}_{n=0}$,""$(\Dsarn)^{\fklam{m+1}}_{n=0}]}

\newcommand{\Aalo}		{\mathbf{A}_{\alo}}
\newcommand{\Asalo}		{\Aalo^{\sklam{s}}}
\newcommand{\Aale}		{\mathbf{A}_{\al1}}
\newcommand{\Aaln}		{\mathbf{A}_{\aln}}
\newcommand{\Asaln}		{\Aaln^{\sklam{s}}}
\newcommand{\Aalnm}		{\mathbf{A}_{\aln-1}}
\newcommand{\Aalfm}		{\mathbf{A}_{\al\fklam{m}}}
\newcommand{\Asalfm}	{\Aalfm^{\sklam{s}}}
\newcommand{\Balo}		{\mathbf{B}_{\alo}}
\newcommand{\Bsalo}		{\Balo^{\sklam{s}}}
\newcommand{\Bale}		{\mathbf{B}_{\al1}}
\newcommand{\Baln}		{\mathbf{B}_{\aln}}
\newcommand{\Bsaln}		{\Baln^{\sklam{s}}}
\newcommand{\Balnm}		{\mathbf{B}_{\aln-1}}
\newcommand{\Balnp}		{\mathbf{B}_{\aln+1}}
\newcommand{\Balfm}		{\mathbf{B}_{\al\fklam{m+1}}}
\newcommand{\Bsalfm}	{\Balfm^{\sklam{s}}}
\newcommand{\Calo}		{\mathbf{C}_{\alo}}
\newcommand{\Csalo}		{\Calo^{\sklam{s}}}
\newcommand{\Cale}		{\mathbf{C}_{\al1}}
\newcommand{\Caln}		{\mathbf{C}_{\aln}}
\newcommand{\Csaln}		{\Caln^{\sklam{s}}}
\newcommand{\Calnm}		{\mathbf{C}_{\aln-1}}
\newcommand{\Calfm}		{\mathbf{C}_{\al\fklam{m}}}
\newcommand{\Csalfm}	{\Calfm^{\sklam{s}}}
\newcommand{\Dalo}		{\mathbf{D}_{\alo}}
\newcommand{\Dsalo}		{\Dalo^{\sklam{s}}}
\newcommand{\Dale}		{\mathbf{D}_{\al1}}
\newcommand{\Daln}		{\mathbf{D}_{\aln}}
\newcommand{\Dsaln}		{\Daln^{\sklam{s}}}
\newcommand{\Dalnm}		{\mathbf{D}_{\aln-1}}
\newcommand{\Dalnp}		{\mathbf{D}_{\aln+1}}
\newcommand{\Dalfm}		{\mathbf{D}_{\al\fklam{m+1}}}
\newcommand{\Dsalfm}	{\Dalfm^{\sklam{s}}}
\newcommand{\ABCDalk}	{[(\Aaln)^{\fklam{\kappa}}_{n=0}$,""$(\Baln)^{\fklam{\kappa+1}}_{n=0}$,""$(\Caln)^{\fklam{\kappa}}_{n=0}$,""$(\Daln)^{\fklam{\kappa+1}}_{n=0}]}
\newcommand{\ABCDsalk}	{[(\Asaln)^{\fklam{\kappa}}_{n=0}$,""$(\Bsaln)^{\fklam{\kappa+1}}_{n=0}$,""$(\Csaln)^{\fklam{\kappa}}_{n=0}$,""$(\Dsaln)^{\fklam{\kappa+1}}_{n=0}]}
\newcommand{\ABCDalm}	{[(\Aaln)^{\fklam{m}}_{n=0}$,""$(\Baln)^{\fklam{m+1}}_{n=0}$,""$(\Caln)^{\fklam{m}}_{n=0}$,""$(\Daln)^{\fklam{m+1}}_{n=0}]}
\newcommand{\ABCDsalm}	{[(\Asaln)^{\fklam{m}}_{n=0}$,""$(\Bsaln)^{\fklam{m+1}}_{n=0}$,""$(\Csaln)^{\fklam{m}}_{n=0}$,""$(\Dsaln)^{\fklam{m+1}}_{n=0}]}

\newcommand{\Uaro}		{U_{\aro}}
\newcommand{\Uarm}		{U_{\arm}}
\newcommand{\Usarm}		{\Uarm^{\sklam{s}}}
\newcommand{\Uarnm}		{U_{\arn-1}}
\newcommand{\Uarnp}		{U_{\arn+1}}
\newcommand{\Uarzn}		{U_{\ar2n}}
\newcommand{\Uarznm}	{U_{\ar2n-1}}
\newcommand{\Uarznp}	{U_{\ar2n+1}}
\newcommand{\Uark}		{U_{\ar\kappa}}
\newcommand{\Usark}		{\Uark^{\sklam{s}}}
\newcommand{\Uarmk}		{(\Uarm)^{\kappa}_{m=0}}
\newcommand{\Usarmk}	{(\Usarm)^{\kappa}_{m=0}}
\newcommand{\tUaro}		{\widetilde{U}_{\aro}}
\newcommand{\tUarm}		{\widetilde{U}_{\arm}}

\newcommand{\Ualo}		{U_{\alo}}
\newcommand{\Ualm}		{U_{\alm}}
\newcommand{\Usalm}		{\Ualm^{\sklam{s}}}
\newcommand{\Ualnm}		{U_{\aln-1}}
\newcommand{\Ualnp}		{U_{\aln+1}}
\newcommand{\Ualzn}		{U_{\al2n}}
\newcommand{\Ualznm}	{U_{\al2n-1}}
\newcommand{\Ualznp}	{U_{\al2n+1}}
\newcommand{\Ualk}		{U_{\al\kappa}}
\newcommand{\Usalk}		{\Ualk^{\sklam{s}}}
\newcommand{\Ualmk}		{(\Ualm)^{\kappa}_{m=0}}
\newcommand{\Usalmk}	{(\Usalm)^{\kappa}_{m=0}}
\newcommand{\tUalo}		{\widetilde{U}_{\alo}}
\newcommand{\tUalm}		{\widetilde{U}_{\alm}}

\newcommand{\Sarn}		{S_{\arn}}
\newcommand{\Sarnm}		{S_{\arn-1}}
\newcommand{\dSaro}		{\widehat{S}_{\aro}}
\newcommand{\dSarn}		{\widehat{S}_{\arn}}
\newcommand{\dSarnm}	{\widehat{S}_{\arn-1}}

\newcommand{\Saln}		{S_{\aln}}
\newcommand{\Salnm}		{S_{\aln-1}}
\newcommand{\dSalo}		{\widehat{S}_{\alo}}
\newcommand{\dSaln}		{\widehat{S}_{\aln}}
\newcommand{\dSalnm}	{\widehat{S}_{\aln-1}}

\newcommand{\un}		{u_n}
\newcommand{\usn}		{\un^{\sklam{s}}}
\newcommand{\uo}		{u_0}
\newcommand{\uso}		{\uo^{\sklam{s}}}
\newcommand{\una}		{\un^{\ast}}
\newcommand{\usna}		{\big(\usn\big)^{\ast}}
\newcommand{\um}		{u_m}
\newcommand{\uma}		{\um^{\ast}}
\newcommand{\unm}		{u_{n-1}}
\newcommand{\usnm}		{\unm^{\sklam{s}}}
\newcommand{\unma}		{\unm^{\ast}}
\newcommand{\usnma}		{\big(\usnm\big)^{\ast}}

\newcommand{\uaro}		{u_{\aro}}
\newcommand{\uarn}		{u_{\arn}}
\newcommand{\usaro}		{\uaro^{\sklam{s}}}
\newcommand{\usarn}		{\uarn^{\sklam{s}}}
\newcommand{\uarna}		{\uarn^{\ast}}
\newcommand{\usarna}	{\big(\usarn\big)^{\ast}}
\newcommand{\uarnm}		{u_{\arn-1}}
\newcommand{\usarnm}	{\uarnm^{\sklam{s}}}
\newcommand{\uarnma}	{\uarnm^{\ast}}
\newcommand{\usarnma}	{\big(\usarnm\big)^{\ast}}

\newcommand{\ualo}		{u_{\alo}}
\newcommand{\ualn}		{u_{\aln}}
\newcommand{\usalo}		{\ualo^{\sklam{s}}}
\newcommand{\usaln}		{\ualn^{\sklam{s}}}
\newcommand{\umalo}		{u_{-\alo}}
\newcommand{\umaln}		{u_{-\aln}}
\newcommand{\ualna}		{\ualn^{\ast}}
\newcommand{\usalna}	{\big(\usaln\big)^{\ast}}
\newcommand{\ualnm}		{u_{\aln-1}}
\newcommand{\usalnm}	{\ualnm^{\sklam{s}}}
\newcommand{\ualnma}	{\ualnm^{\ast}}
\newcommand{\usalnma}	{\big(\usalnm\big)^{\ast}}

\newcommand{\dL}		{\widehat{L}}
\newcommand{\Ln}		{L_n}
\newcommand{\dLn}		{\dL_n}
\newcommand{\Lna}		{\Ln^{\ast}}
\newcommand{\dLna}		{\dLn^{\ast}}

\newcommand{\Rn}		{R_n}
\newcommand{\Rna}		{\Rn^{\ast}}
\newcommand{\Rm}		{R_m}
\newcommand{\Rma}		{\Rm^{\ast}}
\newcommand{\Rnm}		{R_{n-1}}
\newcommand{\Rnma}		{\Rnm^{\ast}}

\newcommand{\Tn}		{T_n}
\newcommand{\Tna}		{\Tn^{\ast}}
\newcommand{\Tnm}		{T_{n-1}}
\newcommand{\Tnma}		{\Tnm^{\ast}}

\newcommand{\vn}		{v_n}
\newcommand{\vna}		{\vn^{\ast}}
\newcommand{\vm}		{v_m}
\newcommand{\vma}		{\vm^{\ast}}
\newcommand{\vnm}		{v_{n-1}}
\newcommand{\vnma}		{\vnm^{\ast}}

\newcommand{\yn}		{y_{0,n}}
\newcommand{\yna}		{\yn^{\ast}}
\newcommand{\ysn}		{y^{\sklam{s}}_{0,n}}
\newcommand{\ysna}		{\big(\ysn\big)^{\ast}}
\newcommand{\ynm}		{y_{0,n-1}}
\newcommand{\ynma}		{\ynm^{\ast}}
\newcommand{\ysnm}		{y^{\sklam{s}}_{0,n-1}}
\newcommand{\ysnma}		{\big(\ysnm\big)^{\ast}}

\newcommand{\yarn}		{y_{\aro,n}}
\newcommand{\yarna}		{\yarn^{\ast}}
\newcommand{\yarnm}		{y_{\aro,n-1}}
\newcommand{\yarnma}	{\yarnm^{\ast}}

\newcommand{\yaln}		{y_{\alo,n}}
\newcommand{\yalna}		{\yaln^{\ast}}
\newcommand{\yalnm}		{y_{\alo,n-1}}
\newcommand{\yalnma}	{\yalnm^{\ast}}

\newcommand{\tJq}		{\widetilde{J}_q}
\newcommand{\bJPp}		{\mathfrak{P}_J(\Pp)}
\newcommand{\btJqPp}	{\mathfrak{P}_{\tJq}(\Pp)}
\newcommand{\tbJPp}		{\widetilde{\mathfrak{P}}_J(\Pp)}
\newcommand{\tbtJqPp}	{\widetilde{\mathfrak{P}}_{\tJq}(\Pp)}
\newcommand{\Jq}		{J_q}
\newcommand{\bJLam}		{\mathfrak{P}_J(\Lam)}
\newcommand{\bJqLam}	{\mathfrak{P}_{\Jq}(\Lam)}
\newcommand{\tbJLam}	{\widetilde{\mathfrak{P}}_J(\Lam)}
\newcommand{\tbJqLam}	{\widetilde{\mathfrak{P}}_{\Jq}(\Lam)}
\newcommand{\jqq}		{j_{qq}}
\newcommand{\bJJPp}		{\mathfrak{P}_{J^{(2)},J^{(1)}}(\Pp)}
\newcommand{\btJJPp}	{\mathfrak{P}_{\jqq,\tJq}(\Pp)}
\newcommand{\tbJJPp}	{\widetilde{\mathfrak{P}}_{J^{(2)},J^{(1)}}(\Pp)}
\newcommand{\tbtJJPp}	{\widetilde{\mathfrak{P}}_{\jqq,\tJq}(\Pp)}

\newcommand{\PtJqCa}	{{\cal P}^{(q,q)}_{-\tJq,\geq}(\C$,""$[\alpha,\infty))}
\newcommand{\dPtJqCa}	{\widehat{\cal P}^{(q,q)}_{-\tJq,\geq}(\C$,""$[\alpha,\infty))}
\newcommand{\PtJqCma}	{{\cal P}^{(q,q)}_{-\tJq,\geq}(\C$,""$(-\infty,\alpha])}
\newcommand{\dPtJqCma}	{\widehat{\cal P}^{(q,q)}_{-\tJq,\geq}(\C$,""$(-\infty,\alpha])}

\newcommand{\tphi}		{\widetilde{\phi}}
\newcommand{\tpsi}		{\widetilde{\psi}}
\newcommand{\phipsi}	{\binom{\phi}{\psi}}
\newcommand{\tphipsi}	{\bbinom{\tphi}{\tpsi}}
\newcommand{\phipsiz}	{\binom{\phi(z)}{\psi(z)}}
\newcommand{\tphipsiz}	{\binom{\tphi(z)}{\tpsi(z)}}
\newcommand{\aphipsi}	{\sklam{\binom{\phi}{\psi}}}
\newcommand{\atphipsi}	{\sklam{\binom{\tphi}{\tpsi}}}

\newcommand{\SpqG}		{{\cal S}_{\pq}(\G)}
\newcommand{\SpqcG}		{{\cal S}_{\pq}(\widecheck{\G})}
\newcommand{\SqPp}		{{\cal S}_{\qq}(\Pp)}
\newcommand{\SqPm}		{{\cal S}_{\qq}(\Pm)}
\newcommand{\Dqa}		{\D_{q,\alpha}}
\newcommand{\Dqma}		{\D_{q,-\alpha}}
\newcommand{\Eqa}		{\E_{q,\alpha}}
\newcommand{\Eqma}		{\E_{q,-\alpha}}

\newcommand{\tWtJqa}	{\widetilde{\cal W}_{-\tJq, \alpha}}

\newcommand{\tTheta}	{\widehat{\Theta}}

\newcommand{\Fn}		{\mathbf{F}_n}
\newcommand{\Fns}		{\mathbf{F}_{n,s}}
\newcommand{\Farn}		{\mathbf{F}_{\arn}}
\newcommand{\Farns}		{\mathbf{F}_{\arn,s}}
\newcommand{\Faln}		{\mathbf{F}_{\aln}}
\newcommand{\Falns}		{\mathbf{F}_{\aln,s}}
\newcommand{\dFn}		{\widehat{\mathbf{F}}_n}
\newcommand{\dFns}		{\widehat{\mathbf{F}}_{n,s}}
\newcommand{\dFarn}		{\widehat{\mathbf{F}}_{\arn}}
\newcommand{\dFarns}	{\widehat{\mathbf{F}}_{\arn,s}}
\newcommand{\dFarnm}	{\widehat{\mathbf{F}}_{\arn-1}}
\newcommand{\dFaln}		{\widehat{\mathbf{F}}_{\aln}}
\newcommand{\dFalns}	{\widehat{\mathbf{F}}_{\aln,s}}
\newcommand{\dFalnm}	{\widehat{\mathbf{F}}_{\aln-1}}
\newcommand{\dFalnms}	{\widehat{\mathbf{F}}_{\aln-1,s}}
\newcommand{\Ffm}		{\mathbf{F}_{\fklam{m}}}
\newcommand{\Farfm}		{\mathbf{F}_{\ar\fklam{m-1}}}
\newcommand{\Falfm}		{\mathbf{F}_{\al\fklam{m-1}}}
\newcommand{\Ffn}		{\Fn^{[f]}}
\newcommand{\Ffns}		{\Fns^{[f]}}
\newcommand{\Ffarn}		{\Farn^{[f]}}
\newcommand{\Ffarns}	{\Farns^{[f]}}
\newcommand{\Ffaln}		{\Faln^{[f]}}
\newcommand{\Ffalns}	{\Falns^{[f]}}
\newcommand{\dFfn}		{\dFn^{[f]}}
\newcommand{\dFfns}		{\dFns^{[f]}}
\newcommand{\dFfarn}	{\dFarn^{[f]}}
\newcommand{\dFfarns}	{\dFarns^{[f]}}
\newcommand{\dFfarnm}	{\dFarnm^{[f]}}
\newcommand{\dFfaln}	{\dFaln^{[f]}}
\newcommand{\dFfalns}	{\dFalns^{[f]}}
\newcommand{\dFfalnm}	{\dFalnm^{[f]}}
\newcommand{\dFfalnms}	{\dFalnms^{[f]}}
\newcommand{\FSn}		{\Fn^{[S]}}
\newcommand{\FSarn}		{\Farn^{[S]}}
\newcommand{\FSaln}		{\Faln^{[S]}}
\newcommand{\dFSn}		{\dFn^{[S]}}
\newcommand{\dFSarn}	{\dFarn^{[S]}}
\newcommand{\dFSaln}	{\dFaln^{[S]}}
\newcommand{\FSfm}		{\Ffm^{[S]}}
\newcommand{\FSarfm}	{\Farfm^{[S]}}
\newcommand{\FSalfm}	{\Falfm^{[S]}}
\newcommand{\FtSfm}		{\Ffm^{[\widetilde{S}]}}
\newcommand{\FtSarfm}	{\Farfm^{[\widetilde{S}]}}
\newcommand{\FtSalfm}	{\Falfm^{[\widetilde{S}]}}

\newcommand{\Varm}		{\mathbf{V}_{\arm}}
\newcommand{\Varn}		{\mathbf{V}_{\ar2n}}
\newcommand{\Varnp}		{\mathbf{V}_{\ar2n+1}}
\newcommand{\Varnm}		{\mathbf{V}_{\ar2n-1}}
\newcommand{\Valm}		{\mathbf{V}_{\alm}}
\newcommand{\Valn}		{\mathbf{V}_{\al2n}}
\newcommand{\Valnp}		{\mathbf{V}_{\al2n+1}}
\newcommand{\Valnm}		{\mathbf{V}_{\al2n-1}}

\newcommand{\Tarn}		{\Theta_{\arn}}
\newcommand{\tTarn}		{\widetilde{\Theta}_{\arn}}
\newcommand{\Yarn}		{\Phi_{\arn}}
\newcommand{\tYarn}		{\widetilde{\Phi}_{\arn}}
\newcommand{\Taln}		{\Theta_{\aln}}
\newcommand{\tTaln}		{\widetilde{\Theta}_{\aln}}
\newcommand{\Yaln}		{\Phi_{\aln}}
\newcommand{\tYaln}		{\widetilde{\Phi}_{\aln}}

\newcommand{\dphi}		{\widehat{\phi}}
\newcommand{\dpsi}		{\widehat{\psi}}
\newcommand{\dphipsi}	{\bbinom{\dphi}{\dpsi}}
\newcommand{\dphipsiz}	{\begin{pmatrix}\dphi(z)\\\dpsi(z)\end{pmatrix}}

\newcommand{\Smin}		{S_{min}}
\newcommand{\Smax}		{S_{max}}
\newcommand{\Sarmmin}	{\Smin^{(\arm)}}
\newcommand{\Sarmmins}	{\Smin^{(\arm,s)}}
\newcommand{\Sarmmax}	{\Smax^{(\arm)}}
\newcommand{\Sarmmaxs}	{\Smax^{(\arm,s)}}
\newcommand{\Salmmin}	{\Smin^{(\alm)}}
\newcommand{\Salmmins}	{\Smin^{(\alm,s)}}
\newcommand{\Salmmax}	{\Smax^{(\alm)}}
\newcommand{\Salmmaxs}	{\Smax^{(\alm,s)}}

\newcommand{\MF}		{{\cal M}_F}
\newcommand{\tF}		{\widetilde{F}}

\newcommand{\Sarm}		{\Sigma_{\arm}}
\newcommand{\Salm}		{\Sigma_{\alm}}

\newcommand{\Rqab}		{{\cal R}_q[\alpha,\beta]}
\newcommand{\Rqabsjk}	{{\cal R}_q\eklam{[\alpha,\beta],\sjk}}

\newcommand{\dgfa}		{Dann gelten folgende Aussagen:}
\newcommand{\egfa}		{Es gelten folgende Aussagen:}
\newcommand{\dsfaa}		{Dann sind folgende Aussagen äquivalent:}
\newcommand{\esfaa}		{Es sind folgende Aussagen äquivalent:}

\newcommand{\bwanf}		{\textit{Beweis:} }
\newcommand{\bwend}		{\hfill $\Box$}

\hyphenation{Stiel-tjes Qua-dru-pel Dyu-ka-rev}	

%% file: einleitung.tex
\newpage
\setcounter{section}{-1}
\section{Einleitung}

\subsection{Historischer Hintergrund}

Diese Arbeit beschäftigt sich mit einer speziellen Art von matriziellen Momentenproblemen, die ihren Ursprung in der fundamentalen klassischen Arbeit \cite{Sur} von Thomas-Joannes Stieltjes (1856-1894) hat. Dort behandelte er die Konvergenz von Kettenbrüchen. Seine Untersuchungen führten ihn zu einem Momentenproblem auf der rechten reellen Halbachse $[0,\infty)$.

Wesentliche Impulse für das Studium von Momentenproblemen gingen von den Arbeiten von H.\,L. Hamburger \cite{Ham} und R. Nevanlinna \cite{Nev} aus, in denen im Gegensatz zu T.-J. Stieltjes nun das Potenzmomentenproblem auf der reellen Achse studiert wurde. Während in den Untersuchungen von H.\,L. Hamburger die Methode der Kettenbrüche noch eine tragende Rolle spielte, wurden in der Arbeit von R. Nevanlinna erstmals Methoden der komplexen Funktionentheorie im Kontext von Momentenproblemen angewendet. Eine fundamentale Rolle in der Geschichte von Momentenproblemen inklusive ihrer matriziellen und operatoriellen Verallgemeinerung spielte Mark Grigorjewitsch Krein (1907-1989). Insbesondere arbeitete er wesentliche Verbindungen zwischen dem Ideengut der Sankt Petersburger Schule von P.\,L. \v{C}eby\v{s}ev und A.\,A. Markov in der zweiten Hälfte des 19. Jahrhunderts und den außerhalb Russlands vollzogenen Entwicklungslinien heraus, welche neben den bereits erwähnten T.-J. Stieltjes, H.\,L. Hamburger und R. Nevanlinna auch mit solchen Namen wie C. Carath\`{e}odory, G. Pick, F. Hausdorff, M. Riesz, E. Hellinger u.\,a. verbunden sind. 

Die Untersuchungen von M.\,G. Krein brachten wichtige neue Erkenntnisse und eine größere Transparenz für das Studium von Momentenproblemen (siehe insbesondere \cite{AK}). M.\,G. Krein war auch derjenige, der erstmals dem Studium von sogenannten finiten Potenzmomentenproblemen Aufmerksamkeit schenkte, bei denen nur eine endliche Anzahl von Momenten vorgeschrieben sind. Hierzu ist insbesondere seine fundamentale Arbeit \cite{Kr5} zu erwähnen, aus der zwei Monographien mit großer Signalwirkung und vielfältiger Anwendung in Analysis und Stochastik entsprangen, nämlich \cite{KS} und \cite{KN}. Insbesondere die Monographie \cite{KN} von M.\,G. Krein und A.\,A. Nudelman lieferte wesentliche Impulse für die Herangehensweise in der vorliegenden Dissertation. Dies betrifft speziell die intensive Heranziehung von Methoden der Theorie orthogonaler Polynome. Hierbei sei erwähnt, dass nun entsprechende matrizielle Verallgemeinerungen vorgenommen werden mussten, was aufgrund des nunmehr behandelten Matrixfalles und der hiermit verbundenen Nichtkommutativität der Matrixmultiplikation den Aufbau verschiedener neuer Konzepte erforderte.

Hinsichtlich weiterer bedeutender Beiträge von M.\,G. Krein zu Potenzmomentenproblemen sei auf die Arbeiten \cite{Kr1}, \cite{Kr2}, \cite{Kr3} und \cite{Kr4} verwiesen (eine Zusammenstellung der Studien von M.\,G. Krein über Potenzmomentenprobleme findet man in den Übersichtsartikeln \cite{Nm1} und \cite{Nm2} von A.\,A. Nudelman). Dabei nutzte er Methoden der Operatortheorie, die später von V.\,M. Adamyan, I.\,M. Tkachenko und M. Urrea in \cite{ATU1}, \cite{ATU2}, \cite{ATU3} und \cite{ATU4} weiterentwickelt wurden. Für eine ausführliche Beschreibung der Geschichte von Potenzmomentenproblemen bis 2001 sei auf die Einleitung der Dissertation von A.\,E. Choque Rivero \cite{CR} verwiesen. 

Vor dem Hintergrund der vorliegenden Dissertation gehen wir nun detaillierter auf die Matrixversion des Potenzmomentenproblems auf $[0,\infty)$ ein. Hierzu sei zunächst bemerkt, dass die Behandlung dieses Problems in der ukrainischen Stadt Charkiw begann. Dies ist untrennbar mit dem Namen V.\,P. Potapov (1941-1980) verbunden. In der ersten Hälfte der 1970er Jahre hatte er bereits in Odessa mit der nunmehr nach ihm benannten Methode der fundamentalen Matrixungleichungen einen sehr leistungsfähigen Apparat zum Studium von Matrixversionen klassischer Interpolations- und Momentenprobleme geschaffen. Im Jahr 1976 wechselte V.\,P. Potapov von Odessa nach Charkiw, wo er nun eine ganze Reihe von Schülern mit der Behandlung konkreter Aufgabenstellungen betraute. So beschäftigte sich I.\,V. Kovalishina mit den Matrixversionen der Interpolationsprobleme von Carath\'{e}odory und Nevanlinna-Pick sowie auch mit der Matrixversion des Hamburgerschen Momentenproblems (siehe \cite{Ko1} und \cite{Ko2}), während V.\,K. Dubovoj seine Aufmerksamkeit der Matrixversion des Interpolationsproblems von I. Schur zuwandte und hierbei stets Verbindungen zur Theorie kontraktiver Operatoren im Hilbertraum anstrebte (siehe \cite{Du1}). 

Für uns besonders interessant ist jedoch das Dissertationsthema, welches V.\,P. Potapov seinem jungen Aspiranten Yu.\,M. Dyukarev stellte. Dies beinhaltete nämlich gerade das Studium der Matrixversion des Stieltjesschen Momentenproblems. Nach dem Tod V.\,P. Potapovs übernahm V.\,E. Katsnelson die weitere Betreuung der Dissertation von Yu.\,M. Dyukarev. Die ersten in diesem Zusammenhang entstandenen Publikationen \cite{DK} sind hierbei nicht direkt mit dem Momentenproblem, sondern mit dessen algebraischen und funktionentheoretischen Umfeld verbunden. Darüber hinaus wurden dort Fragen der Nevanlinna-Pick-Interpolation für Funktionen der matriziellen Stieltjesklasse diskutiert. Das finite matrizielle Stieltjessche Momentenproblem wurde in der Dissertation \cite{dyu1} von Yu.\,M. Dyukarev lediglich für den Fall einer gegebenen Folge von einer ungeraden Anzahl von komplexen quadratischen Matrizen in der \anf{$\leq$}-Version und zudem im vollständig nichtdegenerierten Fall studiert. Der Fall des analogen Momentenproblems für eine gegebene Folge von einer geraden Anzahl von komplexen quadratischen Matrizen wurde von Yu.\,M. Dyukarev in der deponierten Arbeit \cite{dyu2} abgehandelt. Eine simultane Behandlung beider Fälle erwies sich als schwierig. Unter Verwendung eines allgemeineren Schemas, welches auf Operatoridentitäten und Integraldarstellungen basiert, gelang es Yu.\,M. Dyukarev in \cite{dyu3} simultan gewisse Resultate für den geraden oder ungeraden Fall einer vorgegebenen Anzahl von Momenten zu erhalten. Eine systematische Ausarbeitung dieser allgemeinen Konzeption erfolgte dann in der Habilitationsschrift \cite{dyu4} von Yu.\,M. Dyukarev.

Mitte der 1980er Jahre stellte V.\,E. Katsnelson seinem Aspiranten V.\,A. Bolotnikov die Aufgabe, das finite matrizielle Stieltjessche Momentenproblem im allgemeinsten Fall zu untersuchen. Dies erfolgte vor dem Hintergrund der inzwischen geführten Untersuchungen \cite{Du1} von V.\,K. Dubovoj, in denen eine Methode kreiert wurde, welche die Lösung des matriziellen Schurproblems im allgemeinsten Fall ermöglichte. Die Arbeiten an der Dissertation von V.\,A. Bolotnikov wurden dadurch wesentlich erschwert, dass sowohl V.\,E. Katsnelson als auch V.\,A. Bolotnikov Anfang der 1990er Jahre unabhängig voneinander nach Israel emigrierten und an verschiedenen Institutionen landeten. V.\,E. Katsnelson wurde Professor am Weizmann-Institut in Rechovot, während V.\,A. Bolotnikov eine Doktorandenstelle an der Ben-Gurion-Universität des Negev in Be'er Scheva erhielt und dort von D. Alpay betreut wurde. In Be'er Scheva arbeitete er an einer Synthese der Potapovschen Methode der fundamentalen Matrixungleichungen mit dem von D. Alpay und H. Dym geschaffenen Apparat des RKHS-Zugangs zu Matrixinversionen von klassischen Interpolationsaufgaben. Hierbei steht \anf{RKHS} für Hilberträume mit reproduzierendem Kern. Nach Fertigstellung seiner Dissertation griff V.\,A. Bolotnikov nochmals seine früheren Untersuchungen zum degenerierten finiten matriziellen Stieltjesschen Momentenproblem auf und publizierte hierzu die Arbeit \cite{Bo2}, welche interessante neue Ansätze bereitstellt, aber auch einige Unkorrektheiten enthält.

Ende der 1990er Jahre begannen B. Fritzsche und B. Kirstein ihr (bis in die unmittelbare Gegenwart hineinreichendes) Studium von Matrixversionen klassischer Potenzmomentenprobleme auf der reellen Achse. Am Ausgangspunkt dieser Forschungen standen zahlreiche längere Arbeitsaufenthalte von Yu.\,M. Dyukarev in Leipzig, welche neben der Fertigstellung der Dissertation \cite{CR} von A.\,E. Choque Rivero zu einer Reihe von umfassenden Ergebnissen zum finiten matriziellen Hausdorffschen Momentenproblem für beide möglichen Konstellationen führten (siehe \cite{C06} und \cite{C07}). In die gemeinsamen Untersuchungen von Yu.\,M. Dyukarev, B. Fritzsche und B. Kirstein wurden in einem frühen Stadium bereits die Doktoranden H.\,C. Thiele und C. Mädler eingebunden (siehe \cite{13} und \cite{12}). Es sei bemerkt, dass auch einige Untersuchungen zu matriziellen Potenzmomentenproblemen von G.-N. Chen und Y.-J. Hu vorgenommen wurden (vergleiche \cite{CH3} für das Hamburgersche Momentenproblem sowie \cite{CH1} und \cite{CH2} für das Stieltjessche Momentenproblem).

Ein wesentlicher Bestandteil der Untersuchungen von B. Fritzsche, B. Kirstein und C. Mädler über Matrixversionen von finiten Potenzmomentenproblemen auf der reellen Achse besteht in einer gründlichen Analyse der inneren Struktur jener Folgen von komplexen quadratischen Matrizen, für welche die jeweiligen matriziellen Momentenprobleme lösbar sind. Hierbei liegt die Orientierung darauf, die jeweiligen Matrizenfolgen in bijektiver Weise durch innere Parameter zu charakterisieren. Im Fall Hankel-nichtnegativ definiter Folgen wird hierfür die kanonische Hankel-Parametrisierung herangezogen (siehe \cite{13} und \cite{15}), während im Fall von $\alpha$-Stieltjes-nichtnegativ definiten Folgen die $\alpha$-Stieltjes-Parametrisierung (siehe  \cite{12}, \cite{ot226} und \cite{Trans}) Verwendung findet. Für Stieltjes-positiv definite Folgen von komplexen quadratischen Matrizen ist eine weitere Parametrisierung, die im skalaren Fall auf T.-J. Stieltjes \cite{Sur} zurückgeht und im Matrixfall erstmals von Yu.\,M. Dyukarev \cite{dyu} verwendet wurde, von besonderem Interesse. Hierbei handelt es sich um die Dyukarev-Stieltjes-Parametrisierung, welche systematisch erstmals in \cite[Chapter 8]{Trans} studiert und hierzu anschließend wesentlich in \cite{CR1} von A.\,E. Choque Rivero verwendet wurde.

Seit ca. 2010 arbeiten B. Fritzsche, B. Kirstein und C. Mädler an der Fertigstellung einer speziellen Konzeption zur Behandlung des allgemeinsten Falles von finiten matriziellen Potenzmomentenproblemen auf der rellen Achse. Ein Eckpfeiler dieser Konzeption basiert auf der umfassenden Verwendung von Methoden der Schuranalysis. Hierbei lässt sich eine gewisse Dualität betrachten. Es werden nämlich unabhängig voneinander zwei verschiedene Formen von Schur-Typ-Algorithmen ausgearbeitet. Der erste dieser Typen ist von algebraischer Natur und betrifft Folgen von quadratischen Matrizen, welche dadurch gekennzeichnet sind, dass spezielle aus ihnen gebildete Block-Hankel-Matrizen nichtnegativ hermitesch sind. Der zweite Typ von Schur-Typ-Algorithmen ist in der Welt spezieller Klassen von in der offenen oberen Halbebene der komplexen Zahlen oder geeigneter anderer offenen Halbebenen holomorphen Matrixfunktionen mit verschiedenen Zusatzeigenschaften angesiedelt. Die wesentliche Komponente besteht nun in der Synthese dieser beiden Schur-Typ-Algorithmen zu einem einheitlichen Ganzen. Auf der algebraischen Seite des Algorithmus gehen nämlich genau jene endlichen Folgen von komplexen quadratischen Matrizen ein, für die das Momentenproblem lösbar ist, während auf der funktionentheoretischen Seite des Algorithmus genau jene holomorphen Matrixfunktionen agieren, welche als Stieltjes-Transformierte der Lösungen des Momentenproblems in Erscheinung treten. Hinsichtlich der Realisierung dieser Vorgehensweise sei für den Fall des finiten Hamburgerschen Momentenproblems auf \cite{SimH} verwiesen, während der Fall des finiten matriziellen Momentenproblems auf $[\alpha,\infty)$ in \cite{Sim} behandelt wurde.

In der vorliegenden Arbeit konzentrieren wir uns vorwiegend auf den sogenannten vollständig nichtdegenerierten Fall der finiten matriziellen Potenzmomentenprobleme auf Intervallen des Typs $[\alpha,\infty)$ oder $(-\infty,\alpha]$ mit beliebigem reellen $\alpha$. Der allgemeine Fall dieser Momentenprobleme auf Intervallen des Typs $[\alpha,\infty)$ für eine gerade Anzahl von vorgegebenen Momenten wurde in der Dissertation \cite{Maka} von T. Makarevich mit einer speziellen Herangehensweise abgehandelt, welche auf die besondere Spezifik möglicher Degenerierungen eingeht und eine Modifikation einer Methode darstellt, welche auf V.\,K. Dubovoj \cite{Du1} zurückgeht und von ihm für die Behandlung des degenerierten matriziellen Schur-Problems entwickelt wurde. Eine Adaption der Methode von V.\,K. Dubovoj auf das finite matrizielle Stieltjessche Momentenproblem für das Intervall $[0,\infty)$ war bereits von V.\,A. Bolotnikov (siehe \cite{Bo1} und \cite{Bo2}) vorgenommen worden. Die Ziele der Dissertation \cite{Maka} von T. Makarevich bestanden darin, die Resultate von V.\,A. Bolotnikov auf den Fall eines Intervalls $[\alpha,\infty)$ auszudehnen und darüber hinaus einige Lücken und Unkorrektheiten in den Beweisen von V.\,A. Bolotnikov zu beheben. Die Grundstrategie zur Realisierung dieser Zielstellung war hierbei ähnlich wie bei V.\,A. Bolotnikov die Methode der fundamentalen Matrixungleichungen von V.\,P. Potapov unter Einbeziehung der von V.\,K. Dubovoj entwickelten Modifikationen zur Betrachtung degenerierter Situationen.

Das Thema der vorliegenden Dissertation sprengt auf den ersten Blick den Rahmen der in den letzten Jahren verfolgten Konzeption der Forschungen von B. Fritzsche und B. Kirstein zu Matrixversionen von Potenzmomentenproblemen auf der reellen Achse, indem nämlich erneut der vollständig nichtdegenerierte Fall des von T. Makarevich \cite{Maka} in der allgemeinen Situation behandelten $\alpha$-Stieltjes Momentenproblems aufgegriffen wird. Im Unterschied zu \cite{Maka}, wo nur der Fall einer geraden Anzahl von vorgegebenen Momenten behandelt wird, wird es nun allerdings möglich, simultan auch den Fall einer ungeraden Anzahl von vorgegebenen Momenten zu behandeln. Die hier nochmals vorgenommene Behandlung des vollständig nichtdegenerierten Falles erfolgt mit einer alternativen Vorgehensweise, die einen im Vergleich zum allgemeinen Fall einfacheren Zugang ermöglicht und zudem eine ganze Reihe neuer Erkenntnisse von eigenständigem Interesse hervorbringt. Dies betrifft insbesondere die Theorie orthogonaler Matrixpolynome und die damit verbundene Theorie inverser Probleme.

Entscheidende Impulse für die Arbeit an der vorliegenden Dissertation gingen von der 2015 erschienenen Arbeit \cite{CR1} von A.\,E. Choque Rivero aus, welche  nochmals an die früheren Untersuchungen von Yu.\,M. Dyukarev (vergleiche \cite{dyu2}, \cite{dyu1} und \cite{dyu}) zum vollständig nichtdegenerierten Fall des Stieltjesschen Momentenproblems auf dem Intervall $[0,\infty)$ anknüpft. Die wesentlichen neuen Erkenntnisse von A.\,E. Choque Rivero basieren auf einer tiefsinnigen Analyse des von Yu.\,M. Dyukarev verwendeten Apparats von Matrixpolynomen einerseits sowie anderseits auf der Heranziehung einer von B. Fritzsche, B. Kirstein und C. Mädler in \cite[Chapter 8]{Trans} vorgestellten und gründlich studierten Parametrisierung Stieltjes-positiv definiter Folgen von komplexen quadratischen Matrizen. Diese in \cite{Trans} als Dyukarev-Stieltjes-Parametrisierung bezeichnete Parametrisierung kann man im skalaren Fall bereits in der Arbeit \cite{Sur} von T.-J. Stieltjes antreffen, während wesentliche Aspekte der matriziellen Situation von Yu.\,M. Dyukarev in \cite{dyu} erstmals berührt wurden. Insbesondere erkannte Yu.\,M. Dyukarev die Bedeutung der Dyukarev-Stieltjes-Parametrisierung für die Faktorisierung der von ihm zur Beschreibung der Lösungsmenge verwendeten Matrixpolynome. Hieran anknüpfend entwickelte A.\,E. Choque Rivero in \cite{CR1} einen simultanen Zugang zum Studium der Matrixversionen des finiten Stieltjesschen Momentenproblems. Weitere Schlüsselresultate der Arbeit \cite{CR1} betreffen die Herausarbeitung von Zusammenhängen zwischen den von Yu.\,M. Dyukarev eingeführten Matrixpolynomen und der Theorie von orthogonalen Matrixpolynomen. 

Das Hauptziel der vorliegenden Arbeit ist eine Verallgemeinerung der Resultate von Yu.\,M. Dyukarev und A.\,E. Choque Rivero auf den Fall eines Intervalls $[\alpha,\infty)$ mit beliebigem reellen $\alpha$. Hierbei wurde die gleiche Grundstrategie wie in \cite{CR1} verfolgt. Als erster Programmpunkt stand hierbei eine Verallgemeinerung der Dyukarev-Stieltjes-Parametrisierung für $\alpha$-Stieltjes-positiv definite Folgen (vergleiche Kapitel \ref{chapadp}). Auf dieser Grundlage aufbauend werden dann in den weiteren Kapiteln die wesentlichen Resultate der Arbeit \cite{CR1} von A.\,E. Choque Rivero auf den Fall eines beliebigen reellen $\alpha$ verallgemeinert. Die so erzielten Ergebnisse werden dann durch duale Beziehungen auf das Problem für das Intervall $(-\infty,\alpha]$ übertragen.

\subsection[Erste Bezeichnungen und die Formulierung matrizieller Momentenprobleme]{Erste Bezeichnungen und die Formulierung matrizieller \\ Momentenprobleme}

\label{Mengen}
Bevor wir nun die zugrundeliegenden Momentenprobleme formulieren, führen wir zuerst einige Bezeichnungen ein. Wir bezeichnen die Menge der positiven ganzen Zahlen mit $\N$, die Menge aller ganzen Zahlen mit ${\mathbb Z}$, die Menge der reellen Zahlen mit $\R$ und die Menge der komplexen Zahlen mit $\C$. Insbesondere seien dann
\begin{align*}
	\No := \N \cup \gklam{0}, \quad \Na := \N \cup \gklam{\infty}, \quad \Noa := \N \cup \gklam{0,\infty}
\end{align*}
und für $a, b\in{\mathbb Z}$ weiterhin
\begin{align*}
	\Z{a}{b} = \gklam{k\in{\mathbb Z} \;|\; a \leq k \leq b}.
\end{align*}
Für $\kappa \in \Noa$ sei
\label{fklam}
\begin{align*}
 	\fklamo{\kappa} := \max\{n \in \Noa \;|\; 2n \leq \kappa\}.
\end{align*}
Mit
\label{P}
\begin{align*}
	\Pp := \gklam{z\in\C \;|\; \im z >0} \quad \text{bzw.} \quad \Pm := \gklam{z\in\C \;|\; \im z <0}
\end{align*}
bezeichnen wir die obere bzw. untere offene Halbebene von $\C$. \label{La}Für $\alpha\in\R$ seien weiterhin 
\begin{align*}
	\Lap := \gklam{z\in\C \;|\; \re z >\alpha} \quad \text{und} \quad \Lam := \gklam{z\in\C \;|\; \re z <\alpha}.
\end{align*}
\label{Rstr}Seien $X, Y$ und $Z$ nichtleere Mengen mit $Z\subseteq X$ sowie $f:X\rightarrow Y$ eine Abbildung. Dann bezeichne $\Rstr_Z f:Z\rightarrow Y$ definiert gemäß $\eklam{\Rstr_Z f}(z) := f(z)$ die Einschränkung von $f$ auf $Z$.

\label{Cpq}
Im Verlauf der Arbeit sind stets $p,q,r\in\N$. Mit $\Cpq$ bezeichnen wir die Menge der komplexen \textit{p}$\times$\textit{q}-Matrizen. Insbesondere sei $\Cq := \C^{q\times1}$. 
Mit $\Opq$ bzw. $\Iq$ bezeichnen wir die Nullmatrix aus $\Cpq$ bzw. die Einheitsmatrix aus $\Cqq$. 
\label{A}
Sei nun $A\in\Cpq$. Dann bezeichnen $A^{\ast}$ die zu $A$ adjungierte Matrix, $A^{+}$ die Moore-Penrose-Inverse von $A$, $\rank A$ den Rang von A und $\tr A$ die Spur der Matrix A. Mit
\label{reim}
\begin{align*}
	\re A := \frac{1}{2}\rklam{A+A^{\ast}} \quad \text{bzw.} \quad \im A := \frac{1}{2i}\rklam{A-A^{\ast}}
\end{align*}
bezeichnen wir den Real- bzw. Imaginärteil von $A$.
Sei nun $A=(a_{ij})_{i\in\Z{0}{p},j\in\Z{0}{q}}$. Dann bezeichne
\label{norm}
\begin{align*}
	\enorm{A} := \sqrt{\sum^p_{i=0}\sum^q_{j=0}\abs{a_{ij}}^2} \quad \text{bzw.} \quad
	\snorm{A} := \max_{x\in\Cq,\enorm{x}=1}\enorm{Ax}
\end{align*}
die euklidische Norm bzw. Spektralnorm von $A$.
Für $j,k\in\No$ bezeichne
\label{djk}
\begin{align*}
	\delta_{j,k} := \begin{cases} 0 & \text{falls } j\neq k \\
	1 & \text{falls } j=k \end{cases}
\end{align*}
das Kronecker-Delta. Seien nun $n\in\No$ sowie $(p_j)^n_{j=0}$ und $(q_j)^n_{j=0}$ Folgen aus $\N$. Weiterhin sei $A_j\in\C^{p_j\times q_j}$ für alle $j\in\Zon$. Dann sei
\label{diag}
\begin{align*}
	\diag\rklam{A_0, \ldots, A_n} := (\delta_{j,k}A_j)^n_{j=0}.
\end{align*}
\label{A2}Sei nun $A\in\Cqq$. Dann bezeichnen $\det A$ die Determinante von $A$ und im Fall, dass $A$ regulär ist, also $\det A\neq0$ erfüllt, $A^{-1}$ die Inverse von $A$. Weiterhin sei im Fall, dass $A$ regulär ist, $A^{-\ast}:=\eklam{A^{-1}}^{\ast}$.

Seien nun $\G$ eine nichtleere Teilmenge von $\C$ und $B\in\Cpq$. Dann schreiben wir oft zur Vereinfachung statt der konstanten Matrixfunktion ${\cal B}:\G\rightarrow\Cpq$ mit dem Wert $B$ nur $B$. 
\label{f}
Sei weiterhin $f:\G\rightarrow\Cpq$. Dann sei $f^{\ast}:\G\rightarrow\Cqp$ definiert gemäß $f^{\ast}(z):= \eklam{f(z)}^{\ast}$.
Sei nun $f:\G\rightarrow\Cqq$, so dass $\det f:\G\rightarrow\C$ definiert gemäß $\eklam{\det f}(z):=\det\eklam{f(z)}$ nicht die Nullfunktion ist. Dann sei $f^{-1}:\G\rightarrow\Cqq$ bzw. $f^{-\ast}:\G\rightarrow\Cqq$ definiert gemäß $f^{-1}(z):= \eklam{f(z)}^{-1}$ bzw. $f^{-\ast}(z):= \eklam{f(z)}^{-\ast}$.
Wir werden in weiten Teilen der Arbeit holomorphe und speziell meromorphe Matrixfunktionen verwenden. Eine Matrixfunktion $f=(f_{ij})_{i\in\Zep,j\in\Zeq}:\G\rightarrow\Cpq$ heißt in $\G$ holomorph bzw. meromorph, falls $f_{ij}$ für alle $i\in\Zep$ und $j\in\Zeq$ in $\G$ holomorph bzw. meromorph ist. Hierbei sei erwähnt, dass eine skalare Matrixfunktion $g:\G\rightarrow\C$ genau dann meromorph in $\G$ ist, falls eine diskrete Teilmenge $\D$ von $\G$ (d.\,h. $\D$ besitzt keinen Häufungspunkt) existiert, so dass $g$ holomorph in $\G\setminus\D$ ist und $g$ in jedem Punkt von $\D$ einen Pol hat. Eine gute Übersicht über meromorphe skalare Funktionen findet man z.\,B. in \cite[Abschnitt 10.3]{Funk}.

\label{CH}
Mit $\Cqq_H$ bezeichnen wir die Menge aller hermiteschen Matrizen aus $\Cqq$, d.\,h. es ist $A\in\Cqq_H$ genau dann, wenn $A=A^{\ast}$ erfüllt ist.
Mit $\Cqq_>$ bzw. $\Cqq_{\geq}$ bezeichnen wir die Menge der positiv bzw. nichtnegativ hermiteschen Matrizen aus $\Cqq$, d.\,h. es ist $A\in\Cqq_>$ bzw. $A\in\Cqq_{\geq}$ genau dann, wenn $x^{\ast}Ax\in(0,\infty)$ für alle $x\in\Cq\setminus\gklam{0_{q\times1}}$ bzw. $x^{\ast}Ax\in[0,\infty)$ für alle $x\in\Cq$ erfüllt ist. Insbesondere gelten die Inklusionen $\Cqq_> \subseteq \Cqq_{\geq} \subseteq \Cqq_H$. 
Weiterhin verwenden wir die Löwner-Halbordnung, d.\,h. für $A,B\in\Cqq_H$ schreiben wir $A\geq B$ bzw. $A>B$, falls $A-B\in\Cqq_{\geq}$ bzw. $A-B\in\Cqq_{>}$ erfüllt ist. Für $A,B\in\Cqq_H$ mit $A\leq B$ ist
\label{MI}
\begin{align*}
	[A,B] := \gklam{ C\in\Cqq_H \;\big|\; A \leq C \leq B }
\end{align*}
ein Matrixintervall und es heißt nichtdegneriert, falls sogar $A<B$ erfüllt ist. Für $A\in\Cqq_{\geq}$ bezeichne $\sqrt{A}$ die nichtnegative Quadratwurzel von $A$.
In dieser Arbeit werden grundlegende Kenntnisse über hermitesche Matrizen vorausgesetzt, hierbei sei dem Leser z.\,B. auf \cite[Chapters 4, 7, 8]{HJ} verwiesen.

\label{Mq}
Sei $(\Omega,\A)$ ein messbarer Raum. Dann ist $\mu:\A\rightarrow\Cqq_{\geq}$ ein nichtnegativ hermitesches \textit{q}$\times$\textit{q}-Maß auf $(\Omega,\A)$, falls $\mu$ $\sigma$-additiv ist, d.\,h. für jede Folge $(A_j)^{\infty}_{j=0}$ von paarweise disjunkten Teilmengen aus $\A$ stets
\begin{align*}
	\mu\rklam{\bigcup^{\infty}_{j=0}A_j} = \sum^{\infty}_{j=0}\mu(A_j)
\end{align*}
erfüllt ist. Die Menge aller nichtnegativ hermiteschen \textit{q}$\times$\textit{q}-Maße auf $(\Omega,\A)$ bezeichnen wir mit $\Mqoa$.
Sei $\Omega$ eine nichtleere Borel-Teilmenge von $\R$. Dann bezeichnen $\B_{\Omega}$ die zu $\Omega$ gehörige Borel-$\sigma$-Algebra und $\Mqo:=\M^q_{\geq}(\Omega,\B_{\Omega})$.

Seien $(\Omega,\A)$ ein messbarer Raum, $\nu$ ein Maß auf $(\Omega,\A)$ und $\B_{\C}$ die Borel-$\sigma$-Algebra von $\C$. Sei weiterhin $f:\Omega\rightarrow\C$ eine $\A$-$\B_{\C}$-messbare Funktion. Dann heißt $f$ $\nu$-integrierbar, falls
\begin{align*}
	\int_{\Omega}\abs{f} \;d\nu < \infty
\end{align*}
erfüllt ist. Wir schreiben dann $f\in{\cal L}^1(\Omega,\A,\nu;\C)$. Die Funktion $\Phi=(\Phi_{jk})_{j\in\Zep,k\in\Zeq}:\Omega\rightarrow\Cpq$ heißt $\A$-$\B_{\pq}$-messbar, falls für alle $j\in\Zep$ und $k\in\Zeq$ $\Phi_{jk}$ eine $\A$-$\B_{\C}$-messbare Funktion ist. Die Menge $\beklam{{\cal L}^1(\Omega,\A,\nu;\C)}^{\pq}$ umfasst dann alle $\A$-$\B_{\pq}$-messbaren Funktionen $\Phi=(\Phi_{jk})_{j\in\Zep,k\in\Zeq}:\Omega\rightarrow\Cpq$, für die $\Phi_{jk}\in{\cal L}^1(\Omega,\A,\nu;\C)$ für alle $j\in\Zep$ und $k\in\Zeq$ erfüllt ist, und wir schreiben
\begin{align*}
	\int_{\Omega}\Phi \;d\nu := \rklam{\int_{\Omega}\Phi_{jk} \;d\nu}_{j\in\Zep,k\in\Zeq}.
\end{align*}
Sei nun $\mu=(\mu_{jk})^q_{j,k=1}\in\Mqoa$. Für alle $j\in\Zeq$ ist dann $\mu_{jj}$ ein endliches Maß aus $(\Omega,\A)$ und für alle $j,k\in\Zeq$ mit $j\neq k$ ist $\mu_{jk}$ ein komplexes Maß auf $(\Omega,\A)$. Weiterhin ist das Spurmaß von $\mu$
\begin{align*}
	\tau := \sum^q_{j=1} \mu_{jj}
\end{align*}
ein endliches Maß auf $(\Omega,\A)$. Es ist für alle $j,k\in\Zeq$ zudem $\mu_{jk}$ absolut stetig bezüglich $\tau$ und somit existiert die Radon-Nikodym-Ableitung $\frac{d\mu_{jk}}{d\tau}$ von $\mu_{jk}$ bezüglich $\tau$. Dann heißt
\begin{align*}
	\mu^{'}_{\tau} := \rklam{\frac{d\mu_{jk}}{d\tau}}^q_{j,k=1}
\end{align*}
Spurableitung von $\mu$ und es gilt $\mu(A) = \int_A \mu^{'}_{\tau} \;d\tau$ für alle $A\in\A$.
Seien nun $\Phi:\Omega\rightarrow\Cpq$ eine $\A$-$\B_{\pq}$-messbare Funktion und $\Psi:\Omega\rightarrow\Crq$ eine $\A$-$\B_{r\times q}$-messbare Funktion. Dann heißt das Paar $[\Phi,\Psi]$ integrierbar bezüglich $\mu$, falls
\begin{align*}
	\Phi\mu^{'}_{\tau}\Psi^{\ast} \in \beklam{{\cal L}^1(\Omega,\A,\tau;\C)}^{\pr}
\end{align*}
erfüllt ist und wir setzen
\begin{align*}
	\int_{\Omega} \Phi d\mu \Psi^{\ast} := \int_{\Omega} \Phi\mu^{'}_{\tau}\Psi^{\ast} \;d\tau
\end{align*}
sowie für $A\in\A$
\begin{align*}
	\int_A \Phi d\mu \Psi^{\ast} := \int_{\Omega} (1_A\Phi) d\mu (1_A\Psi)^{\ast},
\end{align*}
wobei $1_A:\A\rightarrow\gklam{0,1}$ die Identikatorfunktion von $A$ definiert gemäß
\begin{align*}
	1_A(\omega) := \begin{cases} 0 & \text{falls } \omega\notin A \\
	1 & \text{falls } \omega\in A \end{cases}
\end{align*}
bezeichnet.
\label{L1}
Die Menge ${\cal L}^1(\Omega,\A,\mu;\C)$ umfasst alle $\A$-$\B_{\C}$-messbaren Funktionen $f:\Omega\rightarrow\C$, für die $[f\Iq,1_{\Omega}\Iq]$ integrierbar bezüglich $\mu$ ist, und wir setzen für $A\in\A$
\begin{align*}
	\int_A f \;d\mu := \int_A (f\Iq) d\mu (1_{\Omega}\Iq)^{\ast},
\end{align*}
wofür wir auch gelegentlich $\int_{A}f(\omega) \;\mu(d\omega)$ schreiben. Eine ausführlichere Beschreibung dieses Integrationsbegriffes findet man in \cite{Ka}, \cite{Ro}, \cite{GL} oder \cite[Anhang M]{MP}. Einige weiterführende Resultate behandeln wir in Anhang \ref{chapAA}.

\label{Mqko}
Seien $\kappa \in \Noa$, $m \in \No$, $\Omega\in\B_{\R}$ und für $j\in\No$ weiterhin $f_j:\Omega\rightarrow\C$ definiert gemäß $f_j(t):=t^j$. Dann bezeichne $\Mqko$ die Menge aller $\mu \in \Mqo$, für die $f_j\in{\cal L}^1(\Omega,\B_{\Omega},\mu;\C)$ für alle $j\in\Zok$ erfüllt ist.
Sei nun $\mu\in\Mqko$. Dann bezeichne für $j \in \Zok$
\label{sjs}
\begin{align*}
  s^{(\mu)}_{j}:=\int_{\Omega} t^j \;\mu(dt)
\end{align*}
das $j$-te Moment von $\mu$. Wir betrachten folgende drei Arten von matriziellen Potenzmomentenproblemen:
\label{Ms}
\begin{itemize}
  \item $\Mskg$: Sei $\sjk$ eine Folge aus $\Cqq$. Beschreibe die Menge \linebreak $\Mqskg$ aller $\mu \in \Mqko$, für die $s^{(\mu)}_{j} = s_{j}$ für alle $j \in \Zok$ erfüllt ist.
  \item $\Msmu$: Sei $\sjm$ eine Folge aus $\Cqq$. Beschreibe die Menge \linebreak $\Mqsmu$ aller $\mu \in \Mqmo$, für die $s^{(\mu)}_{m} \leq s_{m}$ und im Fall $m>0$ weiterhin $s^{(\mu)}_{j} = s_{j}$ für alle $j \in \Zomm$ erfüllt sind.
  \item $\Msmuu$: Sei $\sjm$ eine Folge aus $\Cqq$. Beschreibe die Menge \linebreak $\Mqsmuu$ aller $\mu \in \Mqmo$, für die $s^{(\mu)}_{m} \geq s_{m}$ und im Fall $m>0$ weiterhin $s^{(\mu)}_{j} = s_{j}$ für alle $j \in \Zomm$ erfüllt sind.
\end{itemize}
Im Fall $\Omega = \R$ sprechen wir vom Hamburgerschen Momentenproblem, im Fall \linebreak $\Omega = [0,\infty)$ vom Stieltjesschen Momentenproblem und im Fall $\Omega = [\alpha,\beta]$ für $\alpha, \beta\in\R$ mit $\alpha<\beta$ vom Hausdorffschen Momentenproblem. In der vorliegenden Arbeit befassen wir uns hauptsächlich mit den Fällen $\Omega = [\alpha,\infty)$ (rechtsseitiges $\alpha$-Stieltjes Momentenproblem) und $\Omega = (-\infty,\alpha]$ (linksseitiges $\alpha$-Stieltjes  Momentenproblem), wobei $\alpha \in \R$.

Des Studium von finiten Potenzmomentenproblemen, bei denen die Bedingung an das höchste vorgeschriebene Moment die Gestalt einer Ungleichung besitzt, erfolgte erstmals in M.\,G. Kreins richtungsweisender Arbeit \cite[Chapter 1]{Kr5} im Rahmen der Übertragung der Erweiterung der Methode der \v{C}eby\v{s}ev-Systeme auf halbunendliche Intervalle. Eine ausführliche Darlegung dieser Thematik mit eingehender Diskussion der entsprechenden Motivation findet man auch in der Monographie \cite[Chapter 5]{KN} von M.\,G. Krein und A.\,A. Nudelman.

In Kapitel \ref{chapasm} behandeln wir Lösbarkeitsbedingungen für das matrizielle rechts- bzw. linksseitige $\alpha$-Stieltjes Momentenproblem. Es stellt sich heraus, dass ein solches Potenzmomentenproblem genau dann lösbar ist, falls die gegebene Momentenfolge linksseitig bzw. rechtsseitig $\alpha$-Stieltjes"=nichtnegativ definit ist. Hierfür führen wir zunächst den Begriff von Hankel-nichtnegativ definiten Matrizenfolgen ein, die in Verbindung mit dem matriziellen Hamburgerschen Momentenproblem stehen. Weiterhin gehen wir mithilfe der Stieltjes-Transformation vom ursprünglichen Momentenproblem zu einem äquivalenten Interpolationsproblem für eine spezielle Klasse von holomorphen Matrixfunktionen über. Abschließend führen wir die Potapovschen Fundamentalmatrizen für das matrizielle rechts- bzw. linksseitige $\alpha$-Stieltjes Momentenproblem ein, die später in Kapitel \ref{chapdr} eine tragende Rolle spielen.

In Kapitel \ref{chappmp} behandeln wir folgende erste Parametrisierungen einer Matrizenfolge unter Verwendung von gewissen Block-Hankel-Matrizen und deren linken Schur"=Komplemente: Die $\alpha$-Stieltjes-Parametrisierung, die kanonische Hankel-Parametrisierung und das Favard"=Paar. Weiterhin befassen wir uns mit dem monisch links-orthogonalen System von Matrixpolynomen und dem linken System von Matrixpolynomen zweiter Art bezüglich einer Matrizenfolge. Hierbei zeigen wir besonderes Interesse für den Fall, dass die gegebene Folge bis zu einem gewissen Index Hankel-positiv definit ist.

In Kapitel \ref{chapadp} behandeln wir eine für unsere Betrachtungen besonders wichtige Parametrisierung, die $\alpha$-Dyukarev-Stieltjes-Parametrisierung einer rechts- bzw. linksseitig $\alpha$-Stieltjes-positiv definiten Folge, welche im Skalaren für den Fall $\alpha=0$ bereits auf T.-J. Stieltjes \cite{Sur} zurückgeht.

Kapitel \ref{chapdr} nimmt eine besonders zentrale Rolle in der vorliegenden Arbeit ein. Hier erfolgt eine vollständige Beschreibung der Menge der Stieltjes-Transformierten aller Lösungen des finiten vollständig nichtdegenerierten matriziellen $\alpha$-Stieltjes Momentenproblems, bei dem die Bedingung an das höchste vorgeschriebene Moment die Gestalt einer Ungleichung besitzt, in Form einer gebrochen linearen Transformation, deren erzeugende Matrixfunktion ein aus den Ausgangsdaten gebildetes $2$\textit{q}$\times2$\textit{q}-Matrixpolynom ist und als deren Parametermenge eine spezielle Klasse von geordneten Paaren von meromorphen \textit{q}$\times$\textit{q}-Matrixfunktionen fungiert. Die zur Herleitung dieser Beschreibung herangezogene Methode basiert auf der Bestimmung der Lösungsmenge des Systems der beiden zugehörigen fundamentalen Matrixungleichungen vom Potapov-Typ. Die zugrundeliegende erzeugende Matrixfunktion wird auch Resolventenmatrix genannt.

Die Betrachtung der natürlichen \textit{q}$\times$\textit{q}-Blockstruktur der Resolventenmatrix führt uns auf ein Quadrupel von \textit{q}$\times$\textit{q}-Matrixpolynomen, welche mit einer $\alpha$-Stieltjes-positiv definiten Folge verknüpft ist, nämlich das sogenannte $\alpha$-Dyukarev-Quadrupel. Das Studium dieses Quadrupels bildet eine zentrale Rolle von Kapitel \ref{chapdr}. Insbesondere erhalten wir wichtige Informationen über die Verteilung der Nullstellen der Determinanten aller \textit{q}$\times$\textit{q}-Matrixpolynome des Quadrupels.

Weiterhin erfolgt in Kapitel \ref{chapdr} eine ausführliche Diskussion zweier in gewissem Sinne extremaler Lösungen des mithilfe der Stieltjes-Transformation umformulierten Momentenproblems. Diese extremalen Lösungen sind rationale \textit{q}$\times$\textit{q}-Matrixfunktionen, welche interessante Darstellungen besitzen. Die zugehörigen Stieltjes-Maße sind hierbei insbesondere molekular, also auf endlich vielen Punkten der Halbachse $[\alpha,\infty)$ bzw. $(-\infty,\alpha]$ konzentriert.

Im Mittelpunkt von Kapitel \ref{chapfa} steht die multiplikative Zerlegung der Resolventenmatrix für das vollständig nichtdegnerierte matrizielle $\alpha$-Stieltjes Momentenproblem in lineare Matrixpolynome. Hierdurch wird eine Grundlage für einen möglichen Algorithmus vom Schur-Typ zur Lösung des $\alpha$-Stieltjes Momentenproblems geschaffen.

In Kapitel \ref{chapar} wenden wir uns einer alternativen Beschreibung der Menge der Stieltjes-Transformierten aller Lösungen des vollständig nichtdegenerierten matriziellen $\alpha$"=Stieltjes Momentenproblems, bei dem die Bedingung an das höchste vorgeschriebene Moment die Gestalt einer Ungleichung besitzt, zu. Diese Beschreibung erfolgt zwar erneut auf der Grundlage einer gebrochen linearen Transformation von Matrizen. Allerdings verändern wir nun das erzeugende Matrixpolynom der gebrochen linearen Transformation in einer solchen Weise, dass die Menge der Parameterfunktionen nun aus einer gewissen Teilmenge von \textit{q}$\times$\textit{q}-Schur-Funktionen anstelle von speziellen geordneten Paaren von meromorphen \textit{q}$\times$\textit{q}-Matrixfunktionen besteht.

In Kapitel \ref{chapsq} nehmen wir eine weitere Untersuchung der Struktur der in Kapitel \ref{chapdr} eingeführten $2$\textit{q}$\times2$\textit{q}-Dyukarev-Matrixpolynome vor, welche mit einer $\alpha$-Stieltjes-positiv definiten Folge assoziiert sind. Insbesondere stellen wir eine Verbindung zu den in Abschnitt \ref{chapmp} behandelten orthogonalen Polynome her, welche mit einer bis zu einem gewissen Index Hankel-positiv definiten Folge assoziiert werden. Hierbei sei bemerkt, dass mit einer $\alpha$-Stieltjes-positiv definiten Folge zwei jeweils Hankel-positiv definite Folgen verknüpft sind, welche nun für unsere Betrachtungen herangezogen werden.

Wir vertiefen in Kapitel \ref{chapwz} unsere Betrachtungen zur Struktur $\alpha$-Stieltjes-positiv definiter Folgen. Wir studieren nun, in welchem Zusammenhang deren Favard-Paare (vergleiche Abschnitt \ref{chapfp}) zu einigen bereits behandelten Parametrisierungen stehen.

Im abschließenden Kapitel \ref{chapma} gehen wir vor dem Hintergrund der Charakterisierung der Lösbarkeit kurz auf einige Zusammenhänge zwischen finiten matriziellen Momentenproblemen vom Stieltjes-Typ und dem matriziellen Hausdorffschen Momentenproblem ein.

In verschiedenen Anhängen stellen wir oftmals verwendete Hilfsmittel aus Matrixtheorie, Maß- und Integrationstheorie sowie komplexer Funktionentheorie (insbesondere der Theorie verschiedener Klassen meromorpher Matrixfunktionen) zusammen.

%% file: sm0.tex
\newpage
\section[Erste Beobachtungen zu matriziellen \texorpdfstring{$\alpha$}{a}-Stieltjes Momentenproblemen]{Erste Beobachtungen zu matriziellen \texorpdfstring{$\alpha$}{a}-Stieltjes \\ Momentenproblemen} \label{chapasm}

In diesem Kapitel widmen wir uns hauptsächlich Lösbarkeitsbedingungen für das matrizielle rechts- bzw. linksseitige $\alpha$-Stieltjes Momentenproblem. Wir verwenden dafür die Herangehensweise von \cite{12} für den rechtsseitigen Fall. Da wir für diesen Fall nicht ausführlich auf die Beweise eingehen, kann der Leser eine detaillierte Darstellung der Resultate dort finden. Der linksseitige Fall wurde in \cite{ot226} behandelt. Wir werden die dortigen Beweise ausführlicher selbst wiedergeben und die dort erzielten Resultate erweitern.

Als weiterer Bestandteil dieses Kapitels betrachten wir vor dem Hintergrund der Stieltjes"=Transformation (vergleiche Kapitel \ref{chapAS}) eine Umformulierung der matriziellen $\alpha$-Stieltjes Momentenprobleme derart, dass wir statt Maße nun holomorphe Funktionen verwenden. Dies verwenden wir später zur Parametrisierung der Lösungsmenge jener Momentenprobleme (vergleiche Kapitel \ref{chapdr}). Hierfür legen wir mit der Einführung der Potapovschen Fundamentalmatrizen den Grundstein.

\subsection{Lösbarkeitsbedingungen für das matrizielle Hamburgersche Momentenproblem}

Wir wollen zuerst Lösbarkeitsbedingungen für das matrizielle Hamburgersche Momentenproblem aufrufen. Für eine ausführliche Betrachtung jener Lösbarkeitsbedingungen sei dem Leser im Fall einer endlichen gegebenen Momentenfolge auf \cite{13} verwiesen. Der Fall einer unendlichen gegebenen Momentenfolge wurde z.\,B. in \cite{15} behandelt. Die Resultate dieses Abschnitts werden uns im Anschluss helfen, Lösbarkeitsbedingungen für das matrizielle $\alpha$-Stieltjes Momentenproblem zu formulieren. Zunächst führen wir den Begriff der Hankel-nichtnegativ bzw. -positiv Definitheit ein.

\begin{bez}	\thlabel{asmbz1}
	Seien $\kappa \in \Noa$ und $\sjk$ eine Folge aus $\Cpq$.
	Für $j,k \in \No$ mit $j \leq k \leq \kappa$ seien
	\begin{align*}
  		\ysjk := \begin{pmatrix} s_{j} \\ \vdots \\ s_{k} \end{pmatrix} \quad \text{und} \quad \zsjk := \begin{pmatrix} s_{j} & \ldots & s_{k} \end{pmatrix}.
	\end{align*}
	Für $n\in\Zofk$ bezeichnen wir die zu $\sjk$ gehörige \textbf{n-te Block-Hankel-Matrix} mit
	\begin{align*}
  		\Hsn := (s_{j+k})^{n}_{j,k=0}
	\end{align*}
	und das zu $\Hsn$ gehörige \textbf{linke Schur-Komplement} mit
	\begin{align*}
  		\dHsn := \begin{cases} s_{0} & \text{falls } n=0 \\ s_{2n}-z^{\sklam{s}}_{n,2n-1}\big(\Hsnm\big)^{+}y^{\sklam{s}}_{n,2n-1} & \text{falls } n > 0.\end{cases}
	\end{align*}	
	Im Fall $\kappa\geq1$ sei für $n\in\Zofkm$ weiterhin
	\begin{align*}
  		\Ksn := (s_{j+k+1})^{n}_{j,k=0}.
	\end{align*}
	Falls klar ist, von welchem $\sjk$ die Rede ist, lassen wir das \anf{$\sklam{s}$} als oberen Index weg.
\end{bez}

\begin{defi}	\thlabel{asmdef1} \hfill
\begin{itemize}
  \item [\rm{(a)}] Seien $n \in \No$ und $\sjn$ eine Folge aus $\Cqq$. Dann heißt $\sjn$ \textbf{Hankel-nichtnegativ bzw. 
  -positiv definit}, falls $\Hn$ nichtnegativ bzw. positiv hermitesch ist. Mit $\Hqn$ bzw. $\Hpqn$ bezeichnen wir die Menge aller
  Hankel-nichtnegativ bzw. -positiv definiten Folgen $\sjn$ aus $\Cqq$.

  \item [\rm{(b)}] Seien $n \in \No$ und $\sjne$ eine Folge aus $\Cqq$. Dann heißt $\sjne$ \textbf{Hankel-nichtnegativ 
  definit fortsetzbar}, falls ein $s_{2n+2} \in \Cqq$ existiert, so dass $\sjnz \in \Hqne$ erfüllt ist.
  Mit $\Heqne$ bezeichnen wir die Menge aller Hankel-nichtnegativ definit fortsetzbaren Folgen $\sjne$ aus $\Cqq$.

  \item [\rm{(c)}] Seien $n \in \No$ und $\sjn$ eine Folge aus $\Cqq$. Dann heißt $\sjn$ \textbf{Hankel-nichtnegativ 
  definit fortsetzbar}, falls $s_{2n+1},s_{2n+2} \in \Cqq$ existieren, so dass $\sjnz \in \Hqne$ erfüllt ist.
  Mit $\Heqn$ bezeichnen wir die Menge aller Hankel-nichtnegativ definit fortsetzbaren Folgen $\sjn$ aus $\Cqq$. 
\end{itemize}
\end{defi}

Die Definition im Fall einer unendlichen Folge erfolgt vor dem Hintergrund des folgenden Resultats (vergleiche \cite[S. 449]{15}).

\begin{bem}	\thlabel{asmbm1}
  Seien $n \in \N$ und $\sjn \in \Hqn$ bzw. $\sjn \in \Hpqn$. Dann gilt $(s_{j})_{j=0}^{2k} \in {\cal H}^{\geq}_{q,2k}$ bzw.
  $(s_{j})_{j=0}^{2k} \in {\cal H}^{>}_{q,2k}$ für alle $k \in \Zonm$.  
\end{bem}


\bwanf Wegen Teil (a) von \thref{asmdef1} gilt, dass $H_{n}$ nichtnegativ bzw. positiv hermitesch ist.
Sei nun $k \in \Zonm$. Dann liefert ein wohlbekanntes Resultat über nichtnegativ bzw. positiv hermitesche Matrizen,
dass $H_{k}$ als Hauptuntermatrix von $H_{n}$ ebenfalls nichtnegativ bzw. positiv hermitesch ist. Hieraus folgt dann mittels
Teil (a) von \thref{asmdef1} die Behauptung. \bwend

\begin{defi}	\thlabel{asmdef1b}
	Sei $\sj$ eine Folge aus $\Cqq$. Dann heißt $\sj$ \textbf{Hankel"=nichtnegativ bzw. -positiv definit}, 
  falls $\sjn \in \Hqn$ bzw. $\sjn \in \Hpqn$ für alle $n \in \No$ erfüllt ist. Mit $\Hq$ bzw. $\Hpq$ bezeichnen wir die 
  Menge aller Hankel-nichtnegativ bzw. -positiv definiten Folgen $\sj$ aus $\Cqq$.
\end{defi}

Wir geben nun gewisse Inklusionen für die in \thref{asmdef1} und \thref{asmdef1b} eingeführten Mengen an.

\begin{bem}	\thlabel{asmbm6}
	Es gelten $\Hpqn\subseteq\Heqn\subseteq\Hqn$ für alle $n\in\No$ und $\Hpq\subseteq\Hq$.
\end{bem}

\bwanf Sei $n\in\No$. Dann folgt die Inklusion $\Hpqn\subseteq\Heqn$ aus \cite[Remark 2.8]{13} oder \cite[Folgerung 2.10]{Sch}. Die Inklusion $\Heqn\subseteq\Hqn$ folgt direkt aus den Teilen (a) und (c) von \thref{asmdef1} in Verbindung mit \thref{asmbm1}. Aus diesen beiden Inklusionen für beliebige $n\in\No$ folgt dann wegen \thref{asmdef1b} die Inklusion $\Hpq\subseteq\Hq$. \bwend

Folgende Äquivalenz ist eine elementare Schlussfolgerung aufgrund der Struktur der Block-Hankel-Matrix (vergleiche auch \cite[Lemma 2.6(a)]{Maka}).

\begin{bem}	\thlabel{asmbm2}
  Seien $\kappa \in \Noa$ und $\sjk \in \Cqq$. Weiterhin sei $n\in\Zofk$. \dsfaa
  \begin{itemize}
   \item [\rm{(i)}] Es gilt $\Hn \in \C^{(n+1)q \times (n+1)q}_{H}$.
   \item [\rm{(ii)}] Es gilt $s_j \in \Cqq_{H}$ für alle $j \in \Zozn$.
  \end{itemize}
\end{bem}

Es sei zu bemerken, dass im weiteren Verlauf unter Beachtung von \thref{asmbm2} für $n\in\No$ oft $\sjn$ eine Folge aus $\Cqq_H$ ist.
In diesem Fall gilt dann $(\ysjk)^{\ast} = \zsjk$ für alle $j,k \in \No$ mit $j \leq k \leq 2n$.

Wir kommen nun zum Hauptresultat dieses Abschnitts.

\begin{theo}	\thlabel{asmth1} \egfa
\begin{itemize}
  \item [\rm{(a)}] Seien $m \in \No$ und $\sjm$ eine Folge aus $\Cqq$. Dann ist
  $\MqRsmg$ genau dann nichtleer, wenn $\sjm \in \Heqm$ erfüllt ist.
  
  \item [\rm{(b)}] Sei $\sj$ eine Folge aus $\Cqq$. Dann ist
  $\MqRsg$ genau dann nichtleer, wenn $\sj \in \Hq$ erfüllt ist.
  
  \item [\rm{(c)}] Seien $n \in \No$ und $\sjn$ eine Folge aus $\Cqq$. Dann ist
  $\MqRsnu$ genau dann nichtleer, wenn $\sjn \in \Hqn$ erfüllt ist.
\end{itemize}
\end{theo}

\bwanf Zu (a): Siehe \cite[Theorem 4.17]{13} oder \cite[Theorem 9.4]{Sch}.

Zu (b): Siehe \cite[Theorem 6.5]{15}.

Zu (c): Siehe \cite[Theorem 4.16]{13} oder \cite[Theorem 6.6]{Sch}. \bwend

Wir können nun mithilfe von \thref{asmth1} folgende Schlussfolgerung für matrizielle $\alpha$-Stieltjes Momentenprobleme machen:
Seien $\kappa \in \Noa$, $\Omega \in \B_{\R}\setminus\gklam{\emptyset}$, $\tau \in \Mqko$ und $\mu:\B_{\R} \rightarrow \Cqq$ definiert gemäß $\mu(B):=\tau(B \cap \Omega)$. Dann ist augenscheinlich $\mu \in \MqkR$ mit $s^{(\mu)}_{j} = s^{(\tau)}_{j}$ für alle $j \in \Zok$.
Falls für eine Folge $\sjk$ aus $\Cqq$ die Beziehung $\tau \in \Mqskg$ bzw. für ein $n \in \No$ und eine Folge $\sjn$ aus $\Cqq$ die Beziehung $\tau \in \Mqsnu$ erfüllt ist, so gilt auch $\mu \in \MqRskg$ bzw. $\mu \in \MqRsnu$. 
Falls für ein $m \in \No$ und eine Folge $\sjm$ aus $\Cqq$ die Menge $\Mqsmg$ bzw. für eine Folge $\sj$ aus $\Cqq$ die Menge $\Mqsg$ bzw. für ein $n \in \No$ und eine Folge $\sjn$ aus $\Cqq$ ide Menge $\Mqsnu$ nichtleer ist, folgt mithilfe von \thref{asmth1} dann $\sjm \in \Heqm$ bzw. $\sj \in \Hq$ bzw. $\sjn \in \Hqn$. 
Betrachten wir also Lösbarkeitsbedingungen matrizieller $\alpha$-Stieltjes Momentenprobleme führt uns dies auf Teilmengen von $\Heqm$ bzw. $\Hq$ bzw. $\Hqn$.

\subsection[Lösbarkeitsbedingungen für das matrizielle rechtsseitige \texorpdfstring{$\alpha$}{a}"=Stieltjes Momentenproblem]{Lösbarkeitsbedingungen für das matrizielle rechtsseitige \\ \texorpdfstring{$\alpha$}{a}"=Stieltjes Momentenproblem} \label{chaprsm}

Nun befassen wir uns mit Lösbarkeitsbedingungen für das matrizielle rechtsseitige $\alpha$-Stieltjes Momentenproblem. Für eine ausführliche Betrachtung jener Lösbarkeitsbedingungen sei dem Leser im Fall einer endlichen gegebenen Momentenfolge auf \cite{12} verwiesen. Der Fall einer unendlichen gegebenen Momentenfolge wurde z.\,B. in \cite{ot226} behandelt.
Zunächst führen wir den Begriff der rechtsseitig $\alpha$-Stieltjes-nichtnegativ bzw. -positiv Definitheit ein. Der Fall einer unendlichen Folge steht erneut vor dem Hintergrund von \thref{asmbm1}.

\begin{defi}	\thlabel{asmdef2} Sei $\alpha\in \R$.
\begin{itemize}
  \item [\rm{(a)}] Wir setzen $\Kqoa := \Hqo$ sowie für $n \in \N$
  \begin{align*}
    \Kqna &:= \bgklam{\sjn \in \Hqn \;\big|\; (-\alpha s_{j} + s_{j+1})^{2(n-1)}_{j=0} \in \Hqnm} \quad \text{und} \\
    \Kqnea &:= \bgklam{\sjne \text{ Folge aus } \Cqq \;\big|\; \sjn,(-\alpha s_{j} + s_{j+1})^{2n}_{j=0} \in \Hqn}.
  \end{align*}
  Weiterhin bezeichnet $\Kqa$ die Menge aller Folgen $\sj$ aus $\Cqq$, für die $\sjm \in \Kqma$ für alle $m \in \No$ erfüllt ist. Sei $\kappa \in \Noa$. Dann heißt \linebreak $\sjk \in \Kqka$ \textbf{rechtsseitig $\alpha$-Stieltjes-nichtnegativ definit}.
  
  \item [\rm{(b)}] Wir setzen $\Kpqoa := \Hpqo$ sowie für $n \in \N$
  \begin{align*}
    \Kpqna &:= \bgklam{\sjn \in \Hpqn \;\big|\; (-\alpha s_{j} + s_{j+1})^{2(n-1)}_{j=0} \in \Hpqnm} \quad \text{und} \\
    \Kpqnea &:= \bgklam{\sjne \text{ Folge aus } \Cqq \;\big|\; \sjn,(-\alpha s_{j} + s_{j+1})^{2n}_{j=0} \in \Hpqn}.
  \end{align*}
  Weiterhin bezeichnet $\Kpqa$ die Menge aller Folgen $\sj$ aus $\Cqq$, für die $\sjm \in \Kpqma$ für alle $m \in \No$ erfüllt ist. Sei $\kappa \in \Noa$. Dann heißt \linebreak $\sjk \in \Kpqka$ \textbf{rechtsseitig $\alpha$-Stieltjes-positiv definit}.
  
  \item [\rm{(c)}] Seien $m \in \No$ und $\sjm$ eine Folge aus $\Cqq$. Dann heißt $\sjm$ \textbf{rechtsseitig $\alpha$-Stieltjes-nichtnegativ definit fortsetzbar}, falls ein $s_{m+1} \in \Cqq$ existiert mit $\sjme \in \Kqmea$. Mit $\Keqma$ bezeichnen wir die Menge aller rechtsseitig $\alpha$-Stieltjes-nichtnegativ definit fortsetzbaren Folgen $\sjm$ aus $\Cqq$.
\end{itemize}
\end{defi}

Folgende Bemerkung liefert uns eine alternative Formulierung der Teile (a) und (b) von \thref{asmdef2}.

\begin{bem}	\thlabel{aspbm1}	Seien $\alpha \in \R$ und $n \in \No$. \dgfa
\begin{itemize}
  \item [\rm{(a)}] Es gilt $\sjn \in \Kqna$ bzw. $\sjn \in \Kpqna$ genau dann, wenn 
  $\Hn$ und im Fall $n \in \N$ auch $-\alpha\Hnm$ $+\Knm$ nichtnegativ bzw. positiv hermitesch sind.
  
  \item [\rm{(b)}] Es gilt $\sjne \in \Kqnea$ bzw. $\sjne \in \Kpqnea$ genau dann, wenn 
  $\Hn$ und $-\alpha\Hn+\Kn$ nichtnegativ bzw. positiv hermitesch sind.
\end{itemize}
\end{bem}


Wir können nun folgende Beobachtungen für die neu eingeführten Mengen aus \thref{asmdef2} machen.

\begin{satz}	\thlabel{asmsa1} \egfa
	\begin{itemize}
		\item [\rm{(a)}] Es gelten $\Kpqma\subset\Keqma\subset\Kqma$ für alle $m\in\No$ und $\Kpqa\subset\Kqa$.
		\item [\rm{(b)}] Seien $m\in\N$ und $\sjm\in\Kpqma$ bzw. $\sjm\in\Kqma$ bzw. $\sjm\in\Keqma$. Dann gilt $\sjl\in\Kpqla$ bzw. $\sjl\in\Kqla$ bzw. $\sjl\in\Keqla$ für alle $l\in\Zomm$.
		\item [\rm{(c)}] Seien $l,m\in\Noa$ mit $l<m$ und $\sjl\in\Kpqla$. Dann existiert eine Folge $(s_j)^{m}_{j=l+1}$ aus $\Cqq$ mit $\sjm\in\Kpqma$.
	\end{itemize}
\end{satz}

\bwanf Zu (a): Der endliche Fall folgt aus \cite[Remark 4.4]{12} und \cite[Remark 4.5]{12}. Hieraus folgt der unendliche Fall sogleich aus den Teilen (a) und (b) von \thref{asmdef2} (vergleiche auch \cite[Remark 2.13]{ot226}). 

Zu (b): Dies folgt aus \thref{asmdef2} und \thref{asmbm1} (vergleiche auch \cite[Remark 1.4]{ot226} und \cite[Remark 2.12]{ot226}).

Zu (c): Siehe \cite[Proposition 4.13]{ot226} und \cite[Proposition 4.14]{ot226}. \bwend

Wir kommen nun zum Hauptresultat dieses Abschnitts (vergleiche auch \cite[Theorem 1.6]{ot226} und \cite[Theorem 1.7]{ot226}).

\begin{theo}	\thlabel{asmth2} Sei $\alpha \in \R$. \dgfa
\begin{itemize}
  \item [\rm{(a)}] Seien $m \in \No$ und $\sjm$ eine Folge aus $\Cqq$. Dann ist
  $\Mqasmg$ genau dann nichtleer, wenn $\sjm \in \Keqma$ erfüllt ist.
  
  \item [\rm{(b)}] Sei $\sj$ eine Folge aus $\Cqq$. Dann ist
  $\Mqasg$ genau dann nichtleer, wenn $\sj \in \Kqa$ erfüllt ist.
  
  \item [\rm{(c)}] Seien $m \in \No$ und $\sjm$ eine Folge aus $\Cqq$. Dann ist
  $\Mqasmu$ genau dann nichtleer, wenn $\sjm \in \Kqma$ erfüllt ist.
\end{itemize}
\end{theo}

\bwanf Zu (a): Siehe \cite[Theorem 1.3]{12}.

Zu (b): Es gilt
\begin{align*}
  \Mqasg = \bigcap^{\infty}_{m=0} \Mqasmg.
\end{align*}
Hieraus folgt mithilfe von (a), den Teilen (a) und (c) von \thref{asmdef2} sowie der matriziellen Version des Helly-Prohorov-Theorems (siehe \cite[Satz 9]{14}) die Behauptung.

Zu (c): Siehe \cite[Theorem 1.4]{12}. \bwend

Mithilfe der $[\alpha,\infty)$-Stieltjes-Transformation (vergleiche \thref{asmth4} und \thref{asmdef4}) können wir nun beide Versionen des matriziellen rechtsseitigen $\alpha$-Stieltjes Momentenproblems wie folgt in ein äquivalentes Interpolationsproblem für holomorphe Matrixfunktionen der Klasse $\Soqa$ (vergleiche \thref{asmbz3}) umformulieren:
\label{Ss}
\begin{itemize}
  \item $\Saskg$: Seien $\alpha \in \R$, $\kappa \in \Noa$ und $\sjk$ eine Folge aus $\Cqq$. 
  Beschreibe die Menge $\Sqaskg$ aller $S \in \Soqa$, deren zugehöriges Stieltjes-Maß zu $\Mqaskg$ gehört.
  \item $\Sasmu$: Seien $\alpha \in \R$, $m \in \No$ und $\sjm$ eine Folge aus $\Cqq$. 
  Beschreibe die Menge $\Sqasmu$ aller $S \in \Soqa$, deren zugehöriges Stieltjes-Maß zu $\Mqasmu$ gehört.
\end{itemize}
In dieser so umformulierten Version werden wir uns im weiteren Verlauf der Arbeit dem ursprünglichen Momentenproblem zuwenden. Speziell werden wir uns in späteren Kapiteln verstärkt dem Problem $\Sasmu$ widmen. Es sei bemerkt, dass eine erste Beschreibung der Lösungsmenge $\Sqaskg$ in der noch nicht veröffentlichten Arbeit \cite[Part II, Chapter 6]{Sim} vorgenommen wird.

Wir nehmen nun eine erste Beobachtung über die Klasse $\Sqasmu$ für eine Folge $\sjm\in\Kqma$ vor.

\begin{satz}	\thlabel{asmsa3}
	Seien $\alpha \in \R$, $m \in \N$ und $\sjm\in\Kqma$. Weiterhin sei \linebreak $S\in\Sqasmu$. Dann gelten
	\begin{align*}
		{\cal R}(S(z)) = {\cal R}(s_0) \quad \text{und} \quad {\cal N}(S(z)) = {\cal N}(s_0)
	\end{align*}
	für alle $z\in\C\setminus[\alpha,\infty)$. Insbesondere ist im Fall $\det s_0 \neq 0$ sogar $\det S(z) \neq 0$ für alle $z\in\C\setminus[\alpha,\infty)$ erfüllt.
\end{satz}

\bwanf Sei $\mu$ das zu $S$ gehörige Stieltjes-Maß. Dann gilt $\mu\in\Mqasmu$. Hieraus folgt
\begin{align*}
	s_0 = s^{(\mu)}_0 = \mu([\alpha,\infty)).
\end{align*}
Hieraus folgt wegen \thref{mabm3} nun
\begin{align*}
		{\cal R}(S(z)) = {\cal R}(\mu([\alpha,\infty))) = {\cal R}(s_0)
\end{align*}
für alle $z\in\C\setminus[\alpha,\infty)$ und 
\begin{align} \label{asmsa3bw1}
		{\cal N}(S(z)) = {\cal N}(\mu([\alpha,\infty))) = {\cal N}(s_0)
\end{align}
für alle $z\in\C\setminus[\alpha,\infty)$.

Sei nun $\det s_0 \neq 0$ erfüllt. Wegen \fref{asmsa3bw1} gilt dann
\begin{align*}
	{\cal N}(S(z)) = {\cal N}(s_0) = \gklam{0_{q\times1}}
\end{align*}
für alle $z\in\C\setminus[\alpha,\infty)$, also folgt $\det S(z) \neq 0$ für alle $z\in\C\setminus[\alpha,\infty)$. \bwend

\subsection[Lösbarkeitsbedingungen für das matrizielle linksseitige \texorpdfstring{$\alpha$}{a}"=Stieltjes Momentenproblem]{Lösbarkeitsbedingungen für das matrizielle linksseitige \\ \texorpdfstring{$\alpha$}{a}"=Stieltjes Momentenproblem}	\label{chaplsm}

Wir wenden uns nun Lösbarkeitsbedingungen für das matrizielle linksseitige $\alpha$-Stieltjes Momentenproblem zu. Dieses Thema wurde in \cite{ot226} behandelt und wir werden die dortige Vorgehensweise rekapitulieren und sogar erweitern.
Zunächst führen wir den Begriff der linksseitig $\alpha$-Stieltjes-nichtnegativ bzw. -positiv Definitheit ein. Der Fall einer unendlichen Folge steht erneut vor dem Hintergrund von \thref{asmbm1}.

\begin{defi}	\thlabel{asmdef3} Sei $\alpha\in \R$.
\begin{itemize}
  \item [\rm{(a)}] Wir setzen $\Lqoa := \Hqo$ sowie für $n \in \N$
  \begin{align*}
    \Lqna &:= \bgklam{\sjn \in \Hqn \;\big|\; (\alpha s_{j} - s_{j+1})^{2(n-1)}_{j=0} \in \Hqnm} \quad \text{und} \\
    \Lqnea &:= \bgklam{\sjne \text{ Folge aus } \Cqq \;\big|\; \sjn,(\alpha s_{j} - s_{j+1})^{2n}_{j=0} \in \Hqn}.
  \end{align*}
  Weiterhin bezeichnet $\Lqa$ die Menge aller Folgen $\sj$ aus $\Cqq$, für die $\sjm \in \Lqma$ für alle $m \in \No$ erfüllt ist. Sei $\kappa \in \No$. Dann heißt $\sjk \in \Lqka$ \textbf{linksseitig $\alpha$-Stieltjes-nichtnegativ definit}.
  
  \item [\rm{(b)}] Wir setzen $\Lpqoa := \Hpqo$ sowie für $n \in \N$
  \begin{align*}
    \Lpqna &:= \bgklam{\sjn \in \Hpqn \;\big|\; (\alpha s_{j} - s_{j+1})^{2(n-1)}_{j=0} \in \Hpqnm} \quad \text{und} \\
    \Lpqnea &:= \bgklam{\sjne \text{ Folge aus } \Cqq \;\big|\; \sjn,(\alpha s_{j} - s_{j+1})^{2n}_{j=0} \in \Hpqn}.
  \end{align*}
  Weiterhin bezeichnet $\Lpqa$ die Menge aller Folgen $\sj$ aus $\Cqq$, für die $\sjm \in \Lpqma$ für alle $m \in \No$ erfüllt ist. Sei $\kappa \in \Noa$. Dann heißt $\sjk \in \Lpqka$ \textbf{linksseitig $\alpha$-Stieltjes-positiv definit}.
  
  \item [\rm{(c)}] Seien $m \in \No$ und $\sjm$ eine Folge aus $\Cqq$. Dann heißt $\sjm$ \textbf{linksseitig $\alpha$-Stieltjes-nichtnegativ definit fortsetzbar}, falls ein $s_{m+1} \in \Cqq$ existiert mit $\sjme \in \Lqmea$. Mit $\Leqma$ bezeichnen wir die Menge aller linksseitig $\alpha$-Stieltjes-nichtnegativ definit fortsetzbaren Folgen $\sjm$ aus $\Cqq$.
\end{itemize}
\end{defi}

Folgende Bemerkung liefert uns eine alternative Formulierung der Teile (a) und (b) von \thref{asmdef3}.

\begin{bem}	\thlabel{asplbm1}	Seien $\alpha \in \R$ und $n \in \No$. \dgfa
\begin{itemize}
  \item [\rm{(a)}] Es gilt $\sjn \in \Lqna$ bzw. $\sjn \in \Lpqna$ genau dann, wenn 
  $\Hn$ und im Fall $n \in \N$ auch $\alpha\Hnm-\Knm$ nichtnegativ bzw. positiv hermitesch sind.
  
  \item [\rm{(b)}] Es gilt $\sjne \in \Lqnea$ bzw. $\sjne \in \Lpqnea$ genau dann, wenn 
  $\Hn$ und $\alpha\Hn-\Kn$ nichtnegativ bzw. positiv hermitesch sind.
\end{itemize}
\end{bem}


Wir wollen nun für die beteiligten Größen einen ersten Zusammenhang zwischen rechtsseitigem und linksseitigem Fall liefern (vergleiche Teil (b) mit \cite[Lemma 4.11]{Trans}). Dafür benötigen wir folgende Bezeichnung.

\begin{bez}	\thlabel{asmbz2}
	Für $n\in\No$ sei
	\begin{align*}
		\Vn := \diag\rklam{(-1)^j\Iq)}^n_{j=0}.
	\end{align*}
\end{bez}

\begin{lemma}	\thlabel{asmlm1}
	Seien $\kappa \in \Noa$ und $\sjk$ eine Folge aus $\Cqq$. Weiterhin sei \linebreak $t_j := (-1)^j s_j$ für alle $j \in \Zok$. Dann gelten
	\begin{itemize}
		\item [\rm{(a)}] Sei $n\in\No$. Dann ist $\Vn$ unitär und es gilt $\Vn = \Vna$.
		\item [\rm{(b)}] Sei $n \in \Zofk$. Dann gilt
		\begin{align*}
			\Htn = \Vn\Hsn\Vna.
		\end{align*}
		\item [\rm{(c)}] Seien $\kappa\geq1$ und $n \in \Zofkm$. Dann gilt
		\begin{align*}
			\Ktn = -\Vn\Ksn\Vna.
		\end{align*}
		\item [\rm{(d)}] Seien $j,k \in \Zok$ mit $j \leq k$. Dann gelten	
		\begin{align*}
			y^{\sklam{t}}_{k,l} = (-1)^kV_{l-k}y^{\sklam{s}}_{k,l}
			\quad \text{und} \quad
			z^{\sklam{t}}_{k,l} = (-1)^kz^{\sklam{s}}_{k,l}V^{\ast}_{l-k}.
		\end{align*}		
		\item [\rm{(e)}] Sei $n\in\Zofk$. Dann gilt	
		\begin{align*}
			\dHn^{\sklam{t}} = \dHsn.
		\end{align*}
	\end{itemize}
\end{lemma}

\bwanf Einen Beweis findet man auch in \cite[Lemma 4.4]{FR}. Wir geben den kurzen Beweis zum besseren Verständnis selbst an.

Zu (a): Es gilt offensichtlich, dass $V_n$ regulär und $\Vn^{-1} = \Vn = \Vna$ erfüllt ist. Hieraus folgt weiterhin, dass $\Vn$ unitär ist.

Zu (b): Unter Beachtung von (a) gilt
\begin{align*}
	\Vn\Hsn\Vna = \rklam{(-1)^{j+k}s_{j+k}}^{n}_{j,k=0} = \rklam{t_{j+k}}^{n}_{j,k=0} = \Htn.
\end{align*}

Zu (c): Unter Beachtung von (a) gilt
\begin{align*}
	\Vn\Ksn\Vna = \rklam{(-1)^{j+k}s_{j+k+1}}^{n}_{j,k=0} = -\rklam{t_{j+k+1}}^{n}_{j,k=0} = -\Ktn.
\end{align*}

Zu (d): Sei $k$ gerade. Dann gelten
\begin{align*}
	V_{l-k}y^{\sklam{s}}_{k,l} = \begin{pmatrix} (-1)^0s_k \\ \vdots \\ (-1)^{l-k}s_l \end{pmatrix} = \begin{pmatrix} t_k \\ \vdots \\ t_l \end{pmatrix} = y^{\sklam{t}}_{k,l}
\end{align*}
und
\begin{align*}
	z^{\sklam{s}}_{k,l}V^{\ast}_{l-k} = \begin{pmatrix} (-1)^0s_k & \ldots & (-1)^{l-k}s_l \end{pmatrix} = \begin{pmatrix} t_k & \ldots & t_l \end{pmatrix} = z^{\sklam{t}}_{k,l}.
\end{align*}
Sei nun $k$ ungerade. Dann gelten
\begin{align*}
	V_{l-k}y^{\sklam{s}}_{k,l} = \begin{pmatrix} (-1)^0s_k \\ \vdots \\ (-1)^{l-k}s_l \end{pmatrix} = -\begin{pmatrix} t_k \\ \vdots \\ t_l \end{pmatrix} = -y^{\sklam{t}}_{k,l}
\end{align*}
und
\begin{align*}
	z^{\sklam{s}}_{k,l}V^{\ast}_{l-k} = \begin{pmatrix} (-1)^0s_k & \ldots & (-1)^{l-k}s_l \end{pmatrix} = -\begin{pmatrix} t_k & \ldots & t_l \end{pmatrix} = -z^{\sklam{t}}_{k,l}.
\end{align*}

Zu (e): Es gilt 
\begin{align*}
	\dH^{\sklam{t}}_0 = t_0 = s_0 = \dHsn. 
\end{align*}
Seien nun $\kappa\geq2$ und $n\geq1$. Wegen (a), (b) und (d) gilt dann
\begin{align*}
	\dHn^{\sklam{t}} &= t_{2n}-z^{\sklam{t}}_{n,2n-1}\big(\Hnm^{\sklam{t}}\big)^{+}y^{\sklam{t}}_{n,2n-1} \\
	&= s_{2n}-z^{\sklam{s}}_{n,2n-1}\Vnm\big(\Vnm\Hsnm\Vnma\big)^{+}\Vnm y^{\sklam{s}}_{n,2n-1} \\
	&= s_{2n}-z^{\sklam{s}}_{n,2n-1}\big(\Hsnm\big)^{+} y^{\sklam{s}}_{n,2n-1} = \dHsn. \tag*{$\Box$}
\end{align*}
	
Folgende Bemerkung liefert uns einen elementaren Zusammenhang zwischen linksseitig und rechtsseitig $\alpha$-Stieltjes-nichtnegativ definiten, -positiv definiten und -nichtnegativ definit fortsetzbaren Folgen.

\begin{bem}	\thlabel{asmbm3} Sei $\alpha \in \R$. \dgfa
\begin{itemize}
  \item [\rm{(a)}] Seien $\kappa \in \Noa$ und $\sjk$ eine Folge aus $\Cqq$. Weiterhin sei $t_j := (-1)^j s_j$ für alle $j \in \Zok$.
  Dann gilt $\sjk \in \Kqka$ bzw. $\sjk \in \Kpqka$ genau dann, wenn $(t_j)^{\kappa}_{j=0} \in {\cal L}^{\geq}_{q,\kappa,-\alpha}$ bzw. $(t_j)^{\kappa}_{j=0} \in \Lpqkma$ erfüllt ist.
  \item [\rm{(b)}] Seien $m \in \No$ und $\sjm$ eine Folge aus $\Cqq$. Weiterhin sei $t_j := (-1)^j s_j$ für alle $j \in \Zom$.
  Dann gilt $\sjm \in \Keqma$ genau dann, wenn $(t_j)^{m}_{j=0} \in {\cal L}^{\geq,e}_{q,m,-\alpha}$ erfüllt ist.
\end{itemize}  
\end{bem}

\bwanf Der nichtnegativ definite Fall wurde in \cite[Lemma 1.5]{ot226} gezeigt. Wir werden den Beweis zur besseren Anschauung noch einmal rekapitulieren. 

Zu (a): Für $n \in \Zofk$ ist wegen Teil (b) von \thref{asmlm1} die Matrix $\Htn$ genau dann nichtnegativ bzw. positiv hermitesch, wenn $\Hsn$ eine nichtnegativ bzw. positiv hermitesche Matrix ist. Wegen der Teile (a) und (b) von \thref{asmlm1} gilt im Fall $\kappa\geq1$ weiterhin
\begin{align*}
	(-\alpha)\Htn-\Ktn = \Vn\rklam{-\alpha\Hsn+\Ksn}\Vna
\end{align*}
für alle $n \in \Zofkm$. Für $\kappa\geq1$ und $n \in \Zofkm$ ist somit die Matrix $(-\alpha)\Htn-\Ktn$ genau dann nichtnegativ bzw. positiv hermitesch, wenn $-\alpha\Hsn+\Ksn$ eine nichtnegativ bzw. positiv hermitesche Matrix ist. Wegen \thref{aspbm1} folgt dann die Behauptung für $\kappa\in\No$. Unter Beachtung von Teil (a) bzw. (b) von \thref{asmdef2} und Teil (a) bzw. (b) von \thref{asmdef3} folgt hieraus auch die Behauptung für $\kappa=\infty$.

Zu (b): Dies folgt unter Beachtung von Teil (c) von \thref{asmdef2} und Teil (c) von \thref{asmdef3} aus (a). \bwend

Mithilfe von \thref{asmbm3} können wir eine Version von \thref{asmsa1} für den linksseitigen Fall formulieren.

\begin{satz}	\thlabel{asmsa2} \egfa
	\begin{itemize}
		\item [\rm{(a)}] Es gelten $\Lpqma\subset\Leqma\subset\Lqma$ für alle $m\in\No$ und $\Lpqa\subset\Lqa$.
		\item [\rm{(b)}] Seien $m\in\N$ und $\sjm\in\Lpqma$ bzw. $\sjm\in\Lqma$ bzw. $\sjm\in\Leqma$. Dann gilt $\sjl\in\Lpqla$ bzw. $\sjl\in\Lqla$ bzw. $\sjl\in\Leqla$ für alle $l\in\Zomm$.
		\item [\rm{(c)}] Seien $l,m\in\Na$ mit $l<m$ und $\sjl\in\Lpqla$. Dann existiert eine Folge $(s_j)^{m}_{j=l+1}$ aus $\Cqq$ mit $\sjm\in\Lpqma$.
	\end{itemize}
\end{satz}


\thref{asmbm3} führt uns weiterhin auf folgende Lösbarkeitsbedingungen für das matrizielle linksseitige $\alpha$-Stieltjes Momentenproblem.

\begin{theo}	\thlabel{asmth3} Sei $\alpha \in \R$. \dgfa
\begin{itemize}
  \item [\rm{(a)}] Seien $m \in \No$ und $\sjm$ eine Folge aus $\Cqq$. Dann ist
  $\Mqmasmg$ genau dann nichtleer, wenn $\sjm \in \Leqma$ erfüllt ist.
  
  \item [\rm{(b)}] Sei $\sj$ eine Folge aus $\Cqq$. Dann ist
  $\Mqmasg$ genau dann nichtleer, wenn $\sj \in \Lqa$ erfüllt ist.
  
  \item [\rm{(c)}] Seien $n \in \No$ und $\sjn$ eine Folge aus $\Cqq$. Dann ist
  $\Mqmasnu$ genau dann nichtleer, wenn $\sjn \in \Lqna$ erfüllt ist.
  
  \item [\rm{(d)}] Seien $n \in \No$ und $\sjne$ eine Folge aus $\Cqq$. Dann ist \linebreak
  $\Mqmasneuu$ genau dann nichtleer, wenn $\sjne \in \Lqnea$ erfüllt ist.
\end{itemize}
\end{theo}

\bwanf (a) und (b) wurde in \cite[Theorem 1.8]{ot226} gezeigt. (c) findet man auch in \cite[Theorem 1.9]{ot226}. Wir werden die Beweise zu besseren Anschauung noch einmal rekapitulieren.

Zu (a): Sei zunächst $\Mqmasmg \neq \emptyset$. Weiterhin seien $\mu \in \M^q_{\geq}\big[(-\infty,\alpha],\linebreak\sjm,=\negthickspace\big]$ und $t_j = (-1)^js_j$ für alle $j\in\Zom$. 
Wegen \thref{hmlm1} ist dann \linebreak $\widecheck{\mu} \in \M^q_{\geq}\eklam{[-\alpha,\infty),(t_j)^{m}_{j=0},=}$.
Hieraus folgt wegen Teil (a) von \thref{asmth2} dann $(t_j)^{m}_{j=0} \in {\cal K}^{\geq,e}_{q,m,-\alpha}$.
Hieraus folgt wiederum wegen Teil (b) von \thref{asmbm3} dann $\sjm \in \Leqma$. 
Da die verwendeten Aussagen alles Äquivalenzen sind, lässt sich umgekehrt zeigen, dass aus $\sjm \in \Leqma$ dann $\Mqmasmg \neq \emptyset$ folgt. 

Zu (b): Es gilt
\begin{align*}
  \Mqmasg = \bigcap^{\infty}_{m=0} \Mqmasmg.
\end{align*}
Hieraus folgt mithilfe von (a), den Teilen (a) und (c) von \thref{asmdef3} sowie der matriziellen Version des Helly-Prohorov-Theorems (siehe \cite[Satz 9]{14}) die Behauptung.

Zu (c): Sei zunächst $\Mqmasnu \neq \emptyset$. Weiterhin sei $\mu \in \M^q_{\geq}\big[(-\infty,\alpha],$ $\sjn,\leq\negmedspace\big]$. Wegen Teil (b) von \thref{hmlm1} gilt
\begin{align*}
  s^{(\widecheck{\mu})}_{2n} = (-1)^{2n}s^{(\mu)}_{2n} = s^{(\mu)}_{2n} \leq s_{2n}.
\end{align*}
Hieraus folgt wegen Teil (a) von \thref{hmlm1} dann $\widecheck{\mu} \in \M^q_{\geq}\eklam{[-\alpha,\infty),((-1)^js_j)^{2n}_{j=0},\leq}$.
Hieraus folgt wegen Teil (c) von \thref{asmth2} dann $((-1)^js_j)^{2n}_{j=0} \in {\cal K}^{\geq}_{q,2n,-\alpha}$.
Hieraus folgt wiederum wegen Teil (a) von \thref{asmbm3} dann $\sjn \in \Lqna$.
Da die verwendeten Aussagen alles Äquivalenzen sind, lässt sich umgekehrt zeigen, dass aus $\sjn \in \Lqna$ dann $\Mqmasnu \neq \emptyset$ folgt. 

Zu (d): Sei zunächst $\Mqmasneuu \neq \emptyset$. Weiterhin sei $\mu \in \M^q_{\geq}\big[(-\infty,\alpha],$ $\sjne,\geq\negmedspace\big]$. Wegen Teil (b) von \thref{hmlm1} gilt
\begin{align*}
  s^{(\widecheck{\mu})}_{2n+1} = (-1)^{2n+1}s^{(\mu)}_{2n} = -s^{(\mu)}_{2n+1} \leq -s_{2n+1}
\end{align*}
Hieraus folgt wegen Teil (a) von \thref{hmlm1} dann $\widecheck{\mu} \in \M^q_{\geq}\eklam{[-\alpha,\infty),((-1)^js_j)^{2n+1}_{j=0},\leq}$.
Hieraus folgt wegen Teil (c) von \thref{asmth2} dann $((-1)^js_j)^{2n+1}_{j=0} \in {\cal K}^{\geq}_{q,2n+1,-\alpha}$.
Hieraus folgt wiederum wegen Teil (a) von \thref{asmbm3} dann $\sjne \in \Lqnea$.
Da die verwendeten Aussagen alles Äquivalenzen sind, lässt sich umgekehrt zeigen, dass aus $\sjne \in \Lqnea$ dann $\Mqmasneuu \neq \emptyset$ folgt. \bwend

Mithilfe der $(-\infty,\alpha]$-Stieltjes-Transformation (vergleiche \thref{asmth5} und \thref{asmdef5}) können wir nun beide Versionen des matriziellen linksseitigen $\alpha$-Stieltjes Momentenproblems wie folgt in ein äquivalentes Interpolationsproblem für holomorphe Matrixfunktionen der Klasse $\Soqma$ (vergleiche \thref{asmbz4}) umformulieren:
\label{Ssl}
\begin{itemize}
  \item $\Smaskg$: Seien $\alpha \in \R$, $\kappa \in \Noa$ und $\sjk$ eine Folge aus $\Cqq$. 
  Beschreibe die Menge $\Sqmaskg$ aller $S \in \Soqma$, deren zugehöriges Stieltjes-Maß zu $\Mqmaskg$ gehört.
  \item $\Smasmu$: Seien $\alpha \in \R$, $m \in \No$ und $\sjm$ eine Folge aus $\Cqq$. 
  Beschreibe die Menge $\Sqmasmu$ aller $S \in \Soqma$, deren zugehöriges Stieltjes-Maß zu $\Mqmasmu$, falls $m$ gerade ist, bzw. zu \linebreak $\Mqmasmuu$, falls $m$ ungerade ist, gehört.
\end{itemize}
In dieser so umformulierten Version werden wir uns im weiteren Verlauf der Arbeit dem ursprünglichen Momentenproblem zuwenden. Speziell werden wir uns in späteren Kapiteln verstärkt dem Problem $\Smasmu$ widmen.

Wir nehmen nun eine erste Beobachtung über die Klasse $\Sqmasmu$ für eine Folge $\sjk\in\Lqma$ vor.

\begin{satz}	\thlabel{asmsa4}
	Seien $\alpha \in \R$, $m \in \N$ und $\sjm\in\Lqma$. Weiterhin sei \linebreak $S\in\Sqmasmu$. Dann gelten
	\begin{align*}
		{\cal R}(S(z)) = {\cal R}(s_0) \quad \text{und} \quad {\cal N}(S(z)) = {\cal N}(s_0)
	\end{align*}
	für alle $z\in\C\setminus(-\infty,\alpha]$. Insbesondere ist im Fall $\det s_0 \neq 0$ sogar $\det S(z) \neq 0$ für alle $z\in\C\setminus(-\infty,\alpha]$ erfüllt.
\end{satz}

\bwanf Sei $\mu$ das zu $S$ gehörige Stieltjes-Maß. Dann gilt $\mu\in\Mqasmu$. Hieraus folgt
\begin{align*}
	s_0 = s^{(\mu)}_0 = \mu((-\infty,\alpha]).
\end{align*}
Hieraus folgt wegen \thref{mabm5} nun
\begin{align*}
		{\cal R}(S(z)) = {\cal R}(\mu((-\infty,\alpha])) = {\cal R}(s_0)
\end{align*}
für alle $z\in\C\setminus(-\infty,\alpha]$ und 
\begin{align} \label{asmsa4bw1}
		{\cal N}(S(z)) = {\cal N}(\mu((-\infty,\alpha])) = {\cal N}(s_0)
\end{align}
für alle $z\in\C\setminus(-\infty,\alpha]$.

Sei nun $\det s_0 \neq 0$ erfüllt. Wegen \fref{asmsa4bw1} gilt dann
\begin{align*}
	{\cal N}(S(z)) = {\cal N}(s_0) = \gklam{0_{q\times1}}
\end{align*}
für alle $z\in\C\setminus(-\infty,\alpha]$, also folgt $\det S(z) \neq 0$ für alle $z\in\C\setminus(-\infty,\alpha]$. \bwend

Wir wollen abschließend für dieses Kapitel einige Zusammenhänge zwischen den Lösungsmengen im rechtsseitigen und linksseitigen Fall erwähnen.

\begin{bem}	\thlabel{asmbm5}
	Seien $\alpha \in \R$, $\kappa \in \Noa$ und $\sjk$ eine Folge aus $\Cqq$. Weiterhin sei $t_j := (-1)^j s_j$ für alle $j \in \Zok$. \dgfa
	\begin{itemize}
		\item [\rm{(a)}] \esfaa
		\begin{itemize}
			\item [\rm{(i)}] Es ist $\Sqaskg$ nichtleer.
			\item [\rm{(ii)}] Es ist $\Soqmma[\tjk$,""$=]$ nichtleer.
		\end{itemize}
		\item [\rm{(b)}] Seien {\rm (i)} erfüllt, $S: \C\setminus[\alpha,\infty) \rightarrow \Cqq$ und $\widecheck{S}: \C\setminus(-\infty,-\alpha] \rightarrow \Cqq$ definiert gemäß $\widecheck{S}(z) := -S(-z)$. Dann gilt $S\in\Sqaskg$ genau dann, wenn $\widecheck{S}\in\Soqmma[\tjk$,""$=]$ erfüllt ist.
		\item [\rm{(c)}] Sei $\kappa \in \No$. \dsfaa
		\begin{itemize}
			\item [\rm{(iii)}] Es ist $\Soqa[\sjk$,""$\leq]$ nichtleer.
			\item [\rm{(iv)}] Es ist $\Soqmma[\tjk$,""$\leq]$ nichtleer.
		\end{itemize}
		\item [\rm{(d)}] Seien $\kappa \in \No$, {\rm (iii)} erfüllt, $S: \C\setminus[\alpha,\infty) \rightarrow \Cqq$ und $\widecheck{S}: \C\setminus(-\infty,-\alpha] \rightarrow \Cqq$ definiert gemäß $\widecheck{S}(z) := -S(-z)$. Dann gilt $S\in\Soqa[\sjk$,""$\leq]$ genau dann, wenn $\widecheck{S}\in\Soqmma[\tjk$,""$\leq]$ erfüllt ist.	
	\end{itemize}
\end{bem}

\bwanf Zu (a): Dies folgt im Fall $\kappa\in\No$ wegen Teil (a) von \thref{asmth2} und Teil (a) von \thref{asmth3} aus Teil (b) von \thref{asmbm3}. Im Fall $\kappa=\infty$ folgt dies wegen Teil (b) von \thref{asmth2} und Teil (b) von \thref{asmth3} aus Teil (a) von \thref{asmdef2}, Teil (a) von \thref{asmdef3} und Teil (a) von \thref{asmbm3}.

Zu (b): Sei zunächst $S\in\Sqaskg$. Weiterhin sei $\mu$ das zu $S$ gehörige Stieltjes-Maß. Dann gilt $\mu\in\Mqaskg$. Wegen \thref{hmlm1} ist dann $\widecheck{\mu} \in \M^q_{\geq}\eklam{(-\infty,-\alpha],\tjk,=}$. Wegen \thref{asmdef4} und Teil (c) von \thref{hmsa1} gilt
\begin{align*}
	\widecheck{S}(z) &= -S(-z) 
	= -\int_{[\alpha,\infty)} \frac{1}{t-(-z)} \;\mu(dt) \\
	&= \int_{[\alpha,\infty)} \frac{1}{-t-z} \;\mu(dt)
	= \int_{(-\infty,-\alpha]} \frac{1}{t-z} \;\widecheck{\mu}(dt)
\end{align*}
für alle $z\in\C\setminus(-\infty,-\alpha]$. Somit ist wegen \thref{asmdef5} dann $\widecheck{\mu}$ das zu $\widecheck{S}$ gehörige Stieltjes-Maß. Hieraus folgt dann $\widecheck{S}\in\Soqmma[\tjk$,""$=]$.

Sei nun $\widecheck{S}\in\Soqmma[\tjk$,""$=]$. Weiterhin sei $\tau$ das zu $\widecheck{S}$ gehörige Stieltjes-Maß. Dann gilt $\tau\in\M^q_{\geq}\eklam{(-\infty,-\alpha],\tjk,=}$. Wegen \thref{hmlm1} ist dann $\widecheck{\tau}\in\M^q_{\geq}\big[[\alpha,\infty),$ $\sjk,=\negmedspace\big]$. Wegen \thref{asmdef5} und Teil (c) von \thref{hmsa1} gilt
\begin{align*}
	S(z) &= -\widecheck{S}(-z)
	= \int_{(-\infty,-\alpha]} \frac{1}{t-(-z)} \;\tau(dt) \\
	&= \int_{(-\infty,-\alpha]} \frac{1}{-t-z} \;\tau(dt)
	= \int_{[\alpha,\infty)} \frac{1}{t-z} \;\widecheck{\tau}(dt)
\end{align*}
für alle $z\in\C\setminus[\alpha,\infty)$. Somit ist wegen \thref{asmdef4} dann $\widecheck{\tau}$ das zu $S$ gehörige Stieltjes-Maß. Hieraus folgt dann $S\in\Sqaskg$.

Zu (c): Dies folgt wegen Teil (c) von \thref{asmth2} und der Teile (c) und (d) von \thref{asmth3} aus Teil (a) von \thref{asmbm3}. 

Zu (d): Sei zunächst $S\in\Soqa[\sjk$,""$\leq]$. Weiterhin sei $\mu$ das zu $S$ gehörige Stieltjes-Maß. Dann gilt $\mu\in\M^q_{\geq}\eklam{[\alpha,\infty),\sjk,\leq}$. Wegen Teil (b) von \thref{hmlm1} gilt
\begin{align*}
	s^{(\widecheck{\mu})}_{\kappa} = (-1)^{\kappa}s^{(\mu)}_{\kappa}
	= \begin{cases} s^{(\mu)}_{\kappa} \leq s_{\kappa} = t_{\kappa} & \text{falls } \kappa \text{ gerade} \\ -s^{(\mu)}_{\kappa} \geq -s_{\kappa} = t_{\kappa} & \text{falls } \kappa \text{ ungerade.} \end{cases}
\end{align*}
Hieraus folgt wegen Teil (a) von \thref{hmlm1} dann $\widecheck{\mu}\in\M^q_{\geq}\eklam{(-\infty,-\alpha],\tjk,\leq}$, falls $\kappa$ gerade ist, bzw. $\widecheck{\mu}\in\M^q_{\geq}\eklam{(-\infty,-\alpha],\tjk,\geq}$, falls $\kappa$ ungerade ist. Wegen \thref{asmdef4} und Teil (c) von \thref{hmsa1} gilt
\begin{align*}
	\widecheck{S}(z) &= -S(-z) 
	= -\int_{[\alpha,\infty)} \frac{1}{t-(-z)} \;\mu(dt) \\
	&= \int_{[\alpha,\infty)} \frac{1}{-t-z} \;\mu(dt)
	= \int_{(-\infty,-\alpha]} \frac{1}{t-z} \;\widecheck{\mu}(dt)
\end{align*}
für alle $z\in\C\setminus(-\infty,-\alpha]$. Somit ist wegen \thref{asmdef5} dann $\widecheck{\mu}$ das zu $\widecheck{S}$ gehörige Stieltjes-Maß. Hieraus folgt dann $\widecheck{S}\in\Soqmma[\tjk$,""$\leq]$.

Sei nun $\widecheck{S}\in\Soqmma[\tjk$,""$\leq]$. Weiterhin sei $\tau$ das zu $\widecheck{S}$ gehörige Stieltjes-Maß. Dann gilt $\tau\in\M^q_{\geq}\eklam{(-\infty,-\alpha],\tjk,\leq}$, falls $\kappa$ gerade ist, bzw. $\tau\in\M^q_{\geq}\big[(-\infty,-\alpha],$ $\tjk,\geq\negmedspace\big]$, falls $\kappa$ ungerade ist. Wegen Teil (b) von \thref{hmlm1} gilt
\begin{align*}
	s^{(\widecheck{\tau})}_{\kappa} = (-1)^{\kappa}s^{(\tau)}_{\kappa}
	\leq (-1)^{\kappa}t_{\kappa} = s_{\kappa}.
\end{align*}
Hieraus folgt wegen Teil (a) von \thref{hmlm1} dann $\widecheck{\tau}\in\M^q_{\geq}\eklam{[\alpha,\infty),\sjk,\leq}$. Wegen \thref{asmdef5} und Teil (c) von \thref{hmsa1} gilt
\begin{align*}
	S(z) &= -\widecheck{S}(-z)
	= \int_{(-\infty,-\alpha]} \frac{1}{t-(-z)} \;\tau(dt) \\
	&= \int_{(-\infty,-\alpha]} \frac{1}{-t-z} \;\tau(dt)
	= \int_{[\alpha,\infty)} \frac{1}{t-z} \;\widecheck{\tau}(dt)
\end{align*}
für alle $z\in\C\setminus[\alpha,\infty)$. Somit ist wegen \thref{asmdef4} dann $\widecheck{\tau}$ das zu $S$ gehörige Stieltjes-Maß. Hieraus folgt dann $S\in\Soqa[\sjk$,""$\leq]$. \bwend

\subsection{Die Potapovschen Fundamentalmatrizen für das matrizielle rechtssseitige \texorpdfstring{$\alpha$}{a}-Stieltjes Momentenproblem}	\label{chapFMR}

Wir kommen nun zu den Potapovschen Fundamentalmatrizen für das matrizielle rechtsseitige $\alpha$-Stieltjes Momentenproblem. Sie stehen im engen Zusammenhang zu dem System von Potapovschen fundamentalen Matrixungleichungen, dessen Lösung auch eine Lösung für das rechtsseitige $\alpha$-Stieltjes Momentenproblem in Form der Stieltjes-Transformierten darstellt. Da wir hier nicht näher auf diese Matrixungleichungen eingehen, sei darauf aufmerksam gemacht, dass sie in \cite[Kapitel 5-7]{Maka}, \cite[Kapitel 5]{Sh} und \cite[Kapitel 3]{MP} näher behandelt wurden. Wir werden die in diesem Abschnitt behandelten Resultate später in Abschnitt \ref{chapdrr} für den vollständig nichtdegenerierten Fall anwenden.
Wir führen zunächst einige Bezeichnungen ein, die aus \cite{dyu}, wo der Fall $\alpha=0$ des rechtsseitigen $\alpha$-Stieltjes Momentenproblems behandelt wurde, hervorgehen.

\begin{bez}	\thlabel{drbz1}
	Seien $L_0 := \Iq$, $\dL_0 := \Iq$, $T_0 := 0_{\qq}$ und $v_0 := \Iq$ sowie 
	\begin{align*}
  		L_n := \begin{pmatrix} \Oqn \\ \Inq \end{pmatrix}, \quad \dL_n := \begin{pmatrix} \Inq \\ \Oqn \end{pmatrix}, \quad T_n := \begin{pmatrix} \Oqn & \Oq \\ \Inq & \Onq \end{pmatrix} \quad \text{und} \quad v_n := \begin{pmatrix} \Iq \\ \Onq \end{pmatrix}
	\end{align*}
	für alle $n \in \N$. Unter Beachtung von $\det (I_{(n+1)q}-zT_n) = 1$ für alle $z \in \C$ und $n \in \No$ definieren wir die Funktion $R_n: \C \rightarrow \C^{(n+1)q \times (n+1)q}$ gemäß
	\begin{align*}
  		R_n(z) := (I_{(n+1)q}-zT_n)^{-1} = \begin{pmatrix} \Iq & \Oq & \ldots & \Oq \\ z\Iq & \Iq & \ddots & \vdots \\ \vdots & \ddots & \ddots & \Oq \\ z^n\Iq & \ldots & z\Iq & \Iq \end{pmatrix}
	\end{align*}
	für alle $n \in \No$.
\end{bez}

\begin{bez}	\thlabel{drbz3}
	Seien $\kappa \in \Na$, $\alpha \in \R$ und $\sjk$ eine Folge aus $\Cqq$. Dann seien weiterhin $\uso := \Oq$ und 
	\begin{align*}
  		\usn := \begin{pmatrix} \Oq \\ y^{\sklam{s}}_{0,n-1} \end{pmatrix}
	\end{align*}
	für alle $n \in \Zekp$ sowie $\usaro := s_0$ und 
	\begin{align*}
  		\usarn := \begin{pmatrix} s_0 \\ y^{\sklam{s}}_{\aro,n-1} \end{pmatrix} = \yn^{\sklam{s}}-\alpha\usn
	\end{align*}
	für alle $n \in \Zek$. Falls klar ist, von welchem $\sjk$ die Rede ist, lassen wir $\sklam{s}$ als oberen Index weg.
\end{bez}

Es sei bemerkt, dass die in \thref{drbz3} eingeführten Größen von der in anderen Arbeiten verwendeten Version abweichen kann. In jenen Arbeiten wurde im Allgemeinen noch ein Minus hinzugefügt.

Wir führen nun die Potapovschen Fundamentalmatrizen im rechtsseitigen Fall ein.

\begin{defi}[Potapovsche Fundamentalmatrizen]	\thlabel{drdef3}
  Seien $\kappa \in \Na$, $\alpha \in \R$ und $\sjk$ eine Folge aus $\Cqq$. 
  Weiterhin seien $\G$ eine Teilmenge von $\C$ mit $\G\setminus\R \neq \emptyset$ und $f: \G \rightarrow \Cqq$.
  Für alle $n \in \Zofk$ sei dann $\Ffns:\G\setminus\R \rightarrow \C^{(n+2)q\times (n+2)q}$ definiert gemäß
  \begin{align*}
    \Ffns(z) := \begin{pmatrix} \Hsn & \Rn(z)\rklam{v_nf(z)+\usn} \\ \eklam{\Rn(z)\rklam{v_nf(z)+\usn}}^{\ast} & \frac{f(z)-f^{\ast}(z)}{z-\za} \end{pmatrix}.
  \end{align*}
  Sei nun $f_{\ar}: \G \rightarrow \Cqq$ definiert gemäß $f_{\ar}(z) := (z-\alpha)f(z)$. 
  Für alle $n \in \Zofkm$ sei dann $\Ffarns:\G\setminus\R \rightarrow \C^{(n+2)q\times (n+2)q}$ definiert gemäß
  \begin{align*}
    \Ffarns(z) := \begin{pmatrix} \Hsarn & \Rn(z)\rklam{v_nf_{\ar}(z)+\usarn} \\ \eklam{\Rn(z)\rklam{v_nf_{\ar}(z)+\usarn}}^{\ast} & \frac{f_{\ar}(z)-f_{\ar}^{\ast}(z)}{z-\za} \end{pmatrix}.
  \end{align*}
  Für alle $z \in \G\setminus\R$ und $n \in \Zofk$ bezeichnen wir mit $\dFfns(z)$ das linke Schur-Komplement von $\Ffns(z)$, d.\,h.
  \begin{align*}
    \dFfns(z) = \frac{f(z)-f^{\ast}(z)}{z-\za} - \eklam{\Rn(z)\rklam{v_nf(z)+\usn}}^{\ast}\rklam{\Hsn}^{+}\Rn(z)\rklam{v_nf(z)+\usn}.
  \end{align*}
  Für alle $z \in \G\setminus\R$ und $n \in \Zofkm$ bezeichnen wir mit $\dFfarns(z)$ das linke Schur-Komplement von $\Ffarns(z)$, d.\,h.
  \begin{align*}
    \dFfarns(z) &= \frac{f_{\ar}(z)-f_{\ar}^{\ast}(z)}{z-\za} - \eklam{\Rn(z)\rklam{v_nf_{\ar}(z)+\usarn}}^{\ast} \\
    &\quad \cdot\rklam{\Hsarn}^{+}\Rn(z)\rklam{v_nf_{\ar}(z)+\usarn}.
  \end{align*}
  Falls klar ist, von welchem $\sjk$ die Rede ist, lassen wir das \anf{$s$} im unteren Index weg.
\end{defi}

Folgende Resultate zeigen nun in welchem Zusammenhang die in \thref{drdef3} eingeführten Matrizen zu dem rechtsseitigen $\alpha$-Stieltjes Momentenproblem stehen. Hierfür verwenden wir die am Ende von Abschnitt \ref{chaprsm} eingeführte Umformulierung des Momentenproblems.

\begin{satz}	\thlabel{drsa5}
  Seien $m \in \N$, $\alpha \in \R$ und $\sjm$ eine Folge aus $\Cqq$. Weiterhin seien $\D$ eine diskrete Teilmenge von $\Pp$ und $S:\Pp\setminus\D \rightarrow \Cqq$ eine holomorphe Funktion derart, dass $\FSfm(z)$ und $\FSarfm(z)$ für alle $z \in \Pp\setminus\D$ nichtnegativ hermitesche Matrizen sind.
  Dann gibt es genau ein $\widetilde{S} \in \Sqasmu$ mit $\Rstr_{\Pp\setminus\D}\widetilde{S} = S$.
\end{satz}

\bwanf Siehe \cite[Satz 7.14]{Maka} oder \cite[Satz 5.16]{Sh}. \bwend

\begin{satz}	\thlabel{drsa6}
  Seien $m \in \N$, $\alpha \in \R$ und $\sjm \in \Kqma$. Weiterhin sei \linebreak $S\in\Sqasmu$. 
  Dann sind $\FSfm(z)$ und $\FSarfm(z)$ für alle $z \in \C\setminus\R$ nichtnegativ hermitesche Matrizen.
\end{satz}

\bwanf Dies folgt wegen Teil (c) von \thref{asmth2} aus \cite[Teil (e) von Satz 5.27]{Maka} oder \cite[Teil (e) von Satz 3.22]{MP}. \bwend

\begin{theo}	\thlabel{drth2}
  Seien $m \in \N$, $\alpha \in \R$ und $\sjm$ eine Folge aus $\Cqq$. Weiterhin seien $\D$ eine diskrete Teilmenge von $\Pp$ und
  $S:\C\setminus[\alpha,\infty) \rightarrow \Cqq$ eine holomorphe Funktion. \dsfaa
  \begin{itemize}
    \item [\rm{(i)}] Es gilt $S \in \Sqasmu$.
    \item [\rm{(ii)}] Es sind $\FSfm(z)$ und $\FSarfm(z)$ für alle $z \in \Pp\setminus\D$ nichtnegativ hermitesche Matrizen.
  \end{itemize} 
\end{theo}

\bwanf Siehe \cite[Theorem 7.15]{Maka} oder \cite[Theorem 5.17]{Sh}. \bwend

\subsection{Die Potapovschen Fundamentalmatrizen für das matrizielle linksseitige \texorpdfstring{$\alpha$}{a}-Stieltjes Momentenproblem}	\label{chapFML}

Wir kommen nun zu den Potapovschen Fundamentalmatrizen für das matrizielle linksseitige $\alpha$-Stieltjes Momentenproblem. Sie stehen im engen Zusammenhang zu dem System von Potapovschen fundamentalen Matrixungleichungen, dessen Lösung auch eine Lösung für das linksseitige $\alpha$-Stieltjes Momentenproblem in Form der Stieltjes-Transformierten darstellt. Für die Beweisführung werden wir die Resultate des rechtsseitigen Falles von Abschnitt \ref{chapFMR} heranziehen. Wir werden die in diesem Abschnitt behandelten Resultate später in Abschnitt \ref{chapdrl} für den vollständig nichtdegenerierten Fall anwenden.
Zusätzlich zu den Bezeichnungen für den rechtsseitigen Fall benötigen wir im linksseitigen Fall noch folgende weitere Bezeichnung.

\begin{bez}	\thlabel{drlbz1}
	Seien $\kappa\in\Na$, $\alpha\in\R$ und $\sjk$ eine Folge aus $\Cqq$. Dann seien weiterhin $\usalo:=-s_0$ und 
	\begin{align*}
		\usaln := \begin{pmatrix} -s_0 \\ y^{\sklam{s}}_{\alo,n-1} \end{pmatrix} = -\yn^{\sklam{s}}+\alpha\usn
	\end{align*}
	für alle $n\in\Zek$. Falls klar ist, von welchem $\sjk$ die Rede ist, lassen wir das \anf{$\sklam{s}$} als oberen Index weg.
\end{bez}

Folgendes Lemma liefert uns einige Zusammenhänge zwischen dem rechtsseitigen und linksseitigen Fall und wird und als Grundlage für spätere Beweise dienen. In Teil (d) gehen wir auch kurz auf die in \thref{spbsp1} eingeführte Signaturmatrix $\tJq$ ein.

\begin{lemma}	\thlabel{drllm1}
	Seien $\alpha \in \R$, $\kappa \in \Na$, $\sjk$ eine Folge aus $\Cqq$ und $t_j:=(-1)^js_j$ für alle $j\in\Zok$. \dgfa
	\begin{itemize}
		\item [\rm{(a)}] Sei $n\in\No$. Dann gelten
		\begin{align*}
			\Tn = - \Vn\Tn\Vna
		\end{align*}
		und
		\begin{align*}
			\vn = \Vn\vn.
		\end{align*}
		\item [\rm{(b)}] Sei $n\in\No$ und $z\in\C$. Dann gilt
		\begin{align*}
			\Rn(-z) = \Vn\Rn(z)\Vna.
		\end{align*}
		\item [\rm{(c)}] Es gelten
		\begin{align*}
			\un^{\sklam{t}} = -\Vn\usn
		\end{align*}
		für alle $n\in\Zokp$ und
		\begin{align*}
			\umaln^{\sklam{t}} = -\Vn\usarn
		\end{align*}
		für alle $n\in\Zok$.
		\item [\rm{(d)}] Es gilt
		\begin{align*}
			-\tJq = \Ve\tJq\Vea.
		\end{align*}
	\end{itemize}
\end{lemma}

\bwanf Zu (a): Dies folgt sogleich aus der Definition der beteiligten Größen.

Zu (b): Wegen Teil (a) von \thref{asmlm1} und (a) gilt
\begin{align*}
	\Vn\Rn(z)\Vna &= \Vn(\Inpq-z\Tn)^{-1}\Vna = (\Vn\Vna-z\Vn\Tn\Vna)^{-1} \\
	&= (\Inpq-(-z)\Tn)^{-1} = \Rn(-z).
\end{align*}

Zu (c): Es gilt 
\begin{align*}
	\uo^{\sklam{t}} = \Oq = -\uso = -\Vo\uso.
\end{align*}
Wegen Teil (d) von \thref{asmlm1} gilt weiterhin
\begin{align}	\label{drllm1bw1}
	\un^{\sklam{t}} = \begin{pmatrix} \Oq \\ \ynm^{\sklam{t}} \end{pmatrix} = \begin{pmatrix} \Oq \\ \Vnm\ynm^{\sklam{s}} \end{pmatrix} = -\Vn \begin{pmatrix} \Oq \\ \ynm^{\sklam{s}} \end{pmatrix} = -\Vn\usn
\end{align}
für alle $n\in\Zekp$. Es gilt 
\begin{align*}
	\umalo^{\sklam{t}} = -t_0 = -s_0 = -\Vo\usaro.
\end{align*}
Wegen Teil (d) von \thref{asmlm1} und \fref{drllm1bw1} gilt weiterhin
\begin{align*}
	\umaln^{\sklam{t}} = -\yn^{\sklam{t}}-\alpha\un^{\sklam{t}} = -\Vn\rklam{\yn^{\sklam{s}}-\alpha\usn} = -\Vn\usarn
\end{align*}
für alle $n\in\Zek$. 

Zu (d): Wegen \thref{spbsp1} gilt
\begin{align*}
	\Ve\tJq\Vea = \begin{pmatrix} \Iq & \Oq \\ \Oq & -\Iq \end{pmatrix} \begin{pmatrix} \Oq & -i\Iq \\ i\Iq & \Oq \end{pmatrix} \begin{pmatrix} \Iq & \Oq \\ \Oq & -\Iq \end{pmatrix} = \begin{pmatrix} \Oq & i\Iq \\ -i\Iq & \Oq \end{pmatrix} = -\tJq. \tag*{$\Box$}
\end{align*}

Wir führen nun die Potapovschen Fundamentalmatrizen im linksseitigen Fall ein.

\begin{defi}[Potapovsche Fundamentalmatrizen]	\thlabel{drldef1}
	Seien $\kappa \in \Na$, $\alpha \in \R$ und $\sjk$ eine Folge aus $\Cqq$. 
  	Weiterhin seien $\G$ eine Teilmenge von $\C$ mit $\G\setminus\R \neq \emptyset$ und $f: \G \rightarrow \Cqq$.
  	Für alle $n \in \Zofk$ sei dann $\Ffns:\G\setminus\R \rightarrow \C^{(n+2)q\times (n+2)q}$ definiert gemäß
  	\begin{align*}
    	\Ffns(z) := \begin{pmatrix} \Hsn & \Rn(z)\rklam{v_nf(z)+\usn} \\ \eklam{\Rn(z)\rklam{v_nf(z)+\usn}}^{\ast} & \frac{f(z)-f^{\ast}(z)}{z-\za} \end{pmatrix}.
  	\end{align*}
  	Sei nun $f_{\al}:\G\rightarrow\Cqq$ definiert gemäß $f_{\al}(z) := (\alpha-z)f(z)$. 
  	Für alle $n \in \Zofkm$ sei dann $\Ffalns:\G\setminus\R \rightarrow \C^{(n+2)q\times (n+2)q}$ definiert gemäß
  	\begin{align*}
    	\Ffalns(z) := \begin{pmatrix} \Hsaln & \Rn(z)\rklam{v_nf_{\al}(z)+\usaln} \\ \eklam{\Rn(z)\rklam{v_nf_{\al}(z)+\usaln}}^{\ast} & \frac{f_{\al}(z)-f_{\al}^{\ast}(z)}{z-\za} \end{pmatrix}.
  	\end{align*}
  	Für alle $z \in \G\setminus\R$ und $n \in \Zofk$ bezeichnen wir mit $\dFfns(z)$ das linke Schur-Komplement von $\Ffns(z)$, d.\,h.
  	\begin{align*}
    	\dFfns(z) = \frac{f(z)-f^{\ast}(z)}{z-\za} - \eklam{\Rn(z)\rklam{v_nf(z)+\usn}}^{\ast}\brklam{\Hsn}^{+}\Rn(z)\rklam{v_nf(z)+\usn}.
  	\end{align*}
  	Für alle $z \in \G\setminus\R$ und $n \in \Zofkm$ bezeichnen wir mit $\dFfalns(z)$ das linke Schur-Komplement von $\Ffalns(z)$, d.\,h.
  	\begin{align*}
    	\dFfalns(z) &= \frac{f_{\al}(z)-f_{\al}^{\ast}(z)}{z-\za} - \eklam{\Rn(z)\rklam{v_nf_{\al}(z)+\usaln}}^{\ast} \\
    	&\quad \cdot\brklam{\Hsaln}^{+}\Rn(z)\rklam{v_nf_{\al}(z)+\usaln}.
  	\end{align*}
  	Falls klar ist, von welchem $\sjk$ die Rede ist, lassen wir das \anf{$s$} im unteren Index weg.
\end{defi}

Folgendes Lemma erlaubt uns eine Darstellung der Potapovschen Fundamentalmatrizen und deren linken Schurkomplemente für den linksseitigen Fall mithilfe der Potapovschen Fundamentalmatrizen und deren linken Schurkomplemente für den rechtsseitigen Fall. Dies gibt uns die Möglichkeit, die für den rechtsseitigen Fall erzielten Resultate für die Potapovschen fundamentalen Matrixungleichungen auf den linksseitigen Fall zu übertragen.

\begin{lemma}	\thlabel{drllm2}
	Seien $\kappa \in \Na$, $\alpha \in \R$, $\sjk$ eine Folge aus $\Cqq$ und $t_j:=(-1)^js_j$ für alle $j\in\Zok$. Weiterhin seien $\G$ eine Teilmenge von $\C$ mit $\G\setminus\R \neq \emptyset$, $f: \G \rightarrow \Cqq$, $\widecheck{\G} := \gklam{-z \;|\; z \in \G}$ und $g: \widecheck{\G} \rightarrow \Cqq$ definiert gemäß $g(z) := -f(-z)$. \dgfa
	\begin{itemize}
		\item [\rm{(a)}] Sei $z\in\G\setminus\R$. Dann gelten
		\begin{align*}
			\mathbf{F}^{[g]}_{n,t}(-z) = \begin{pmatrix} \Vn & -\Iq \end{pmatrix} \Ffns(z) \begin{pmatrix} \Vna \\ -\Iq \end{pmatrix}
		\end{align*}
		für alle $n \in \Zofk$ und
		\begin{align*}
			\mathbf{F}^{[g]}_{-\aln,t}(-z) = \begin{pmatrix} \Vn & -\Iq \end{pmatrix} \Ffarns(z) \begin{pmatrix} \Vna \\ -\Iq \end{pmatrix}
		\end{align*}
		für alle $n \in \Zofkm$.
		\item [\rm{(b)}] Sei $z\in\G\setminus\R$. Dann gelten
		\begin{align*}
			\widehat{\mathbf{F}}^{[g]}_{n,t}(-z) = \dFfns(z)
		\end{align*}
		für alle $n \in \Zofk$ und
		\begin{align*}
			\widehat{\mathbf{F}}^{[g]}_{-\aln,t}(-z) = \dFfarns(z)
		\end{align*}
		für alle $n \in \Zofkm$.
	\end{itemize}
\end{lemma}

\bwanf Zu (a): Wegen \thref{drldef1}, der Teile (a) und (b) von \thref{asmlm1}, der Teile (a)-(c) von \thref{drllm1} und \thref{drdef3} gilt
\begin{align*}
	&\ \mathbf{F}^{[g]}_{n,t}(-z) \\
	&= \begin{pmatrix} \Hn^{\sklam{t}} & \Rn(-z)\rklam{v_ng(-z)+\un^{\sklam{t}}} \\ \eklam{\Rn(-z)\rklam{v_ng(-z)+\un^{\sklam{t}}}}^{\ast} & \frac{g(-z)-g^{\ast}(-z)}{-z-(-\za)} \end{pmatrix} \\
	&= \begin{pmatrix} \Vn\Hsn\Vna & \Vn\Rn(z)\Vna\rklam{-\Vn\vn f(z)-\Vn\usn} \\ \eklam{\Vn\Rn(z)\Vna\rklam{-\Vn\vn f(z)-\Vn\usn}}^{\ast} & \frac{f(z)-f^{\ast}(z)}{z-\za} \end{pmatrix} \\
	&= \begin{pmatrix} \Vn & -\Iq \end{pmatrix} \begin{pmatrix} \Hsn & \Rn(z)\rklam{v_nf(z)+\usn} \\ \eklam{\Rn(z)\rklam{v_nf(z)+\usn}}^{\ast} & \frac{f(z)-f^{\ast}(z)}{z-\za} \end{pmatrix} \begin{pmatrix} \Vna \\ -\Iq \end{pmatrix} \\
	&= \begin{pmatrix} \Vn & -\Iq \end{pmatrix} \Ffns(z) \begin{pmatrix} \Vna \\ -\Iq \end{pmatrix}
\end{align*}
für alle $n \in \Zofk$. Es gilt weiterhin
\begin{align}	\label{drllm2bw1}
	g_{-\al}(-z) = (-\alpha-(-z))g(-z) = -(z-\alpha)f(z) = -f_{\ar}(z)
\end{align}
Hieraus folgt wegen \thref{drldef1}, Teil (a) von \thref{asmlm1}, Teil (b) von \thref{asplm1}, der Teile (a)-(c) von \thref{drllm1} und \thref{drdef3} dann
\begin{align*}
	&\ \mathbf{F}^{[g]}_{-\aln,t}(-z) \\
	&= \begin{pmatrix} \Hmaln^{\sklam{t}} & \Rn(-z)\rklam{\vn g_{-\al}(-z)+u^{\sklam{t}}_{-\aln}} \\ \eklam{\Rn(-z)\rklam{\vn g_{-\al}(-z)+u^{\sklam{t}}_{-\aln}}}^{\ast} & \frac{g_{-\al}(-z)-g_{-\al}^{\ast}(-z)}{-z-(-\za)} \end{pmatrix} \\
	&= \begin{pmatrix} \Vn\Hsaln\Vna & \Vn\Rn(z)\Vna\rklam{-\Vn\vn f_{\ar}(z)-\Vn\usarn} \\ \eklam{\Vn\Rn(z)\Vna\rklam{-\Vn\vn f_{\ar}(z)-\Vn\usarn}}^{\ast} & \frac{f_{\ar}(z)-f_{\ar}^{\ast}(z)}{z-\za} \end{pmatrix} \\
	&= \begin{pmatrix} \Vn & -\Iq \end{pmatrix} \begin{pmatrix} \Hsarn & \Rn(z)\rklam{\vn f_{\ar}(z)+\usarn} \\ \eklam{\Rn(z)\rklam{\vn f_{\ar}(z)+\usarn}}^{\ast} & \frac{f_{\ar}(z)-f^{\ast}_{\ar}(z)}{z-\za} \end{pmatrix} \begin{pmatrix} \Vna \\ -\Iq \end{pmatrix} \\
	&= \begin{pmatrix} \Vn & -\Iq \end{pmatrix} \Ffarns(z) \begin{pmatrix} \Vna \\ -\Iq \end{pmatrix}
\end{align*}
für alle $n \in \Zofkm$.

Zu (b): Wegen \thref{drldef1}, der Teile (a) und (b) von \thref{asmlm1}, der Teile (a)-(c) von \thref{drllm1} und \thref{drdef3} gilt
\begin{align*}
	&\ \widehat{\mathbf{F}}^{[g]}_{n,t}(-z) \\
	&= \frac{g(-z)-g^{\ast}(-z)}{-z-(-\za)} - \eklam{\Rn(-z)\rklam{v_ng(-z)+\un^{\sklam{t}}}}^{\ast} \rklam{\Hn^{\sklam{t}}}^{+} \Rn(-z)\rklam{v_ng(-z)+\un^{\sklam{t}}} \\
	&= \frac{f(z)-f^{\ast}(z)}{z-\za} - \eklam{\Vn\Rn(z)\Vna\rklam{-\Vn\vn f(z)-\Vn\usn}}^{\ast} \Vn\brklam{\Hsn}^{+}\Vna \\
	&\quad \cdot\Vn\Rn(z)\Vna\rklam{-\Vn\vn f(z)-\Vn\usn} \\
	&= \frac{f(z)-f^{\ast}(z)}{z-\za} - \eklam{\Rn(z)\rklam{v_nf(z)+\usn}}^{\ast}\brklam{\Hsn}^{+}\Rn(z)\rklam{v_nf(z)+\usn} \\
	&= \dFfns(z)
\end{align*}
für alle $n \in \Zofk$. Wegen \thref{drldef1}, \fref{drllm2bw1}, Teil (a) von \thref{asmlm1}, Teil (b) von \thref{asplm1}, der Teile (a)-(c) von \thref{drllm1} und \thref{drdef3} gilt weiterhin
\begin{align*}
	&\ \widehat{\mathbf{F}}^{[g]}_{-\aln,t}(-z) \\
	&= \frac{g_{-\al}(-z)-g^{\ast}_{-\al}(-z)}{-z-(-\za)} - \eklam{\Rn(-z)\rklam{v_ng_{-\al}(-z)+u^{\sklam{t}}_{-\aln}}}^{\ast} \rklam{\Hmaln^{\sklam{t}}}^{+} \\
	&\quad \cdot\Rn(-z)\rklam{v_ng_{-\al}(-z)+u^{\sklam{t}}_{-\aln}} \\
	&= \frac{f_{\ar}(z)-f^{\ast}_{\ar}(z)}{z-\za} - \eklam{\Vn\Rn(z)\Vna\rklam{-\Vn\vn f_{\ar}(z)-\Vn\usarn}}^{\ast} \Vn\brklam{\Hsarn}^{+}\Vna \\
	&\quad \cdot\Vn\Rn(z)\Vna\rklam{-\Vn\vn f_{\ar}(z)-\Vn\usarn} \\
	&= \frac{f_{\ar}(z)-f^{\ast}_{\ar}(z)}{z-\za} - \eklam{\Rn(z)\rklam{v_nf_{\ar}(z)+\usarn}}^{\ast}\brklam{\Hsarn}^{+}\Rn(z)\rklam{v_nf_{\ar}(z)+\usarn} \\
	&= \dFfarns(z)
\end{align*}
für alle $n \in \Zofkm$. \bwend

Folgende Resultate zeigen nun in welchem Zusammenhang die in \thref{drldef1} eingeführten Matrizen zu dem linksseitigen $\alpha$-Stieltjes Momentenproblem stehen. Hierfür verwenden wir die am Ende von Abschnitt \ref{chaplsm} eingeführte Umformulierung des Momentenproblems.

\begin{satz}	\thlabel{drlsa1}
  Seien $m \in \N$, $\alpha \in \R$ und $\sjm$ eine Folge aus $\Cqq$. Weiterhin seien $\D$ eine diskrete Teilmenge von $\Pm$ und
  $S:\Pm\setminus\D \rightarrow \Cqq$ eine holomorphe Funktion derart, dass $\FSfm(z)$ und $\FSalfm(z)$ für alle $z \in \Pm\setminus\D$ nichtnegativ hermitesche Matrizen sind.
  Dann gibt es genau ein $\widetilde{S} \in \Sqmasmu$ mit $\Rstr_{\Pm\setminus\D}\widetilde{S} = S$.
\end{satz}

\bwanf Seien $t_j = (-1)^js_j$ für alle $j\in\Zom$ und $T(-z) := -S(z)$ für alle $z\in\Pm\setminus\D$. Dann ist $\widecheck{\D}:=\gklam{-z|z\in\D}$ eine diskrete Teilmenge von $\Pp$ und es ist $T:\Pp\setminus\widecheck{\D} \rightarrow \Cqq$ eine holomorphe Funktion. Da unter Beachtung von Teil (a) von \thref{asmlm1} die Matrix $\begin{pmatrix} \Vn & -\Iq \end{pmatrix}$ für alle $n\in\No$ unitär ist und $\mathbf{F}^{[S]}_{\fklam{m},s}(z)$ und $\mathbf{F}^{[S]}_{\al\fklam{m-1},s}(z)$ für alle $z \in \Pm\setminus\D$ nichtnegativ hermitesche Matrizen sind, sind wegen Teil (a) von \thref{drllm2} auch  $\mathbf{F}^{[T]}_{\fklam{m},t}(z)$ und $\mathbf{F}^{[T]}_{-\ar\fklam{m-1},t}(z)$ für alle $z \in \Pp\setminus\widecheck{\D}$ nichtnegativ hermitesche Matrizen. Wegen \thref{drsa5} gibt es dann genau ein $\widetilde{T} \in {\cal S}_{0,q,[-\alpha,\infty)}[\tjm$,""$\leq]$ mit $\Rstr_{\Pp\setminus\widecheck{\D}}\widetilde{T} = T$. Sei nun $\widetilde{S}:\C\setminus(-\infty,\alpha]\rightarrow\Cqq$ definiert gemäß $\widetilde{S}(z) := -\widetilde{T}(-z)$. Hieraus folgen wegen Teil (d) von \thref{asmbm5} nun $\widetilde{S} \in \Sqmasmu$ und $\Rstr_{\Pm\setminus\D}\widetilde{S} = S$. Die Einzigartigkeit eines solchen $\widetilde{S}$ folgt unmittelbar aus der Einzigartigkeit von $\widetilde{T}$ und Teil (d) von \thref{asmbm5}. \bwend

\begin{satz}	\thlabel{drlsa2}
  Seien $m \in \N$, $\alpha \in \R$ und $\sjm \in \Lqma$. Weiterhin sei \linebreak $S\in\Sqmasmu$. 
  Dann sind $\FSfm(z)$ und $\FSalfm(z)$ für alle $z \in \C\setminus\R$ nichtnegativ hermitesche Matrizen.
\end{satz}

\bwanf Seien $t_j = (-1)^js_j$ für alle $j\in\Zom$ und $T:\C\setminus[-\alpha,\infty)$ definiert gemäß \linebreak $T(z) := -S(-z)$. Wegen Teil (a) von \thref{asmbm3} gilt dann $\tjm\in{\cal K}^{\geq}_{q,m,-\alpha}$. Weiterhin gilt wegen Teil (d) von \thref{asmbm5} dann $T \in {\cal S}_{0,q,[-\alpha,\infty)}[\tjm$,""$\leq]$. Hieraus folgt wegen \thref{drsa6} nun, dass $\mathbf{F}^{[T]}_{\fklam{m},t}(z)$ und $\mathbf{F}^{[T]}_{-\ar\fklam{m-1},t}(z)$ für alle $z \in \C\setminus\R$ nichtnegativ hermitesche Matrizen sind. Da unter Beachtung von Teil (a) von \thref{asmlm1} die Matrix $\begin{pmatrix} \Vn & -\Iq \end{pmatrix}$ für alle $n\in\No$ unitär ist, folgt hieraus wegen Teil (a) von \thref{drllm2} dann, dass $\mathbf{F}^{[S]}_{\fklam{m},s}(z)$ und $\mathbf{F}^{[S]}_{\al\fklam{m-1},s}(z)$ für alle $z \in \C\setminus\R$ nichtnegativ hermitesche Matrizen sind. \bwend

\begin{theo}	\thlabel{drlth1}
  Seien $m \in \N$, $\alpha \in \R$ und $\sjm$ eine Folge aus $\Cqq$. Weiterhin seien $\D$ eine diskrete Teilmenge von $\Pm$ und
  $S:\C\setminus(-\infty,\alpha] \rightarrow \Cqq$ eine holomorphe Funktion. \dsfaa
  \begin{itemize}
    \item [\rm{(i)}] Es gilt $S \in \Sqmasmu$.
    \item [\rm{(ii)}] Es sind $\FSfm(z)$ und $\FSalfm(z)$ für alle $z \in \Pm\setminus\D$ nichtnegativ hermitesche Matrizen.
  \end{itemize} 
\end{theo}

\bwanf Seien $t_j = (-1)^js_j$ für alle $j\in\Zom$ und $T:\C\setminus[-\alpha,\infty)\rightarrow\Cqq$ definiert gemäß $T(z) := -S(-z)$. Dann ist $\widecheck{\D}:=\gklam{-z \;|\; z\in\D}$ eine diskrete Teilmenge von $\Pp$ und es ist $T$ eine holomorphe Funktion. Wegen Teil (d) von \thref{asmbm5} ist (i) äquivalent zu
\begin{itemize}
	\item [\rm(iii)] Es gilt $T \in {\cal S}_{0,q,[-\alpha,\infty)}[\tjm$,""$\leq]$.
\end{itemize}
Wegen \thref{drth2} ist (iii) äquivalent zu
\begin{itemize}
	\item [\rm(iv)] Es sind $\mathbf{F}^{[T]}_{\fklam{m},t}(z)$ und $\mathbf{F}^{[T]}_{-\ar\fklam{m-1},t}(z)$ für alle $z \in \Pp\setminus\widecheck{\D}$ nichtnegativ hermitesche Matrizen.
\end{itemize}
Da unter Beachtung von Teil (a) von \thref{asmlm1} die Matrix $\begin{pmatrix} \Vn & -\Iq \end{pmatrix}$ für alle $n\in\No$ unitär ist, ist wegen Teil (a) von \thref{drllm2} dann (iv) äquivalent zu (ii). \bwend

%% file: sm1.tex
\newpage
\section{Über einige zu Matrizenfolgen gehörige Parametrisierungen und Matrixpolynome} \label{chappmp}

Bevor wir uns speziell den $\alpha$-Stieltjes-positiv definiten Matrizenfolgen widmen, gehen wir in diesem Kapitel von allgemeinen Matrizenfolgen aus, wobei wir auch schon Matrizenfolgen betrachten werden, die bis zu einem gewissen Folgenglied Hankel-positiv definit sind. Wir betrachten einige zugehörige Parametrisierungen und Matrixpolynome. Hierfür spielen die in Kapitel \ref{chapasm} eingeführten Block-Hankel-Matrizen und deren linke Schur-Komplemente eine tragende Rolle. Wir werden zeigen, dass einige dieser Parametrisierungen und Matrixpolynome für die gegebene Matrizenfolge einzigartig sind und aus ihnen die Matrizenfolge wiedergewonnen werden kann. Außerdem stellen wir einige Zusammenhänge zwischen den Parametrisierungen und Matrixpolynome her.

\subsection{Die \texorpdfstring{$\alpha$}{a}-Stieltjes-Parametrisierung von Matrizenfolgen}

In diesem Abschnitt betrachten wir die $\alpha$-Stieltjes-Parametrisierung von Matrizenfolgen, die, wie der Name vermuten lässt, im Zusammenhang zur $\alpha$-Stieltjes-nichtnegativ bzw. -positiv Definitheit stehen. Hierfür geben wir einige Resultate aus \cite[Chapter 4]{ot226} wieder.
Zunächst führen wir den Begriff der durch rechtsseitige bzw. linksseitige $\alpha$-Verschiebung generierten Folge ein, die wir schon in \thref{asmdef2} und \thref{asmdef3}, ohne sie beim Namen zu nennen, verwendet haben.

\begin{defi}	\thlabel{aspdef1} Seien $\alpha \in \C$, $\kappa \in \Na$ und $\sjk$ eine Folge aus $\Cpq$.
\begin{itemize}
  \item [\rm{(a)}] $\sarjk$ definiert gemäß $s_{\arj} := -\alpha s_{j}+s_{j+1}$ für alle $j \in \Zokm$ heißt
  die aus $\sjk$ \textbf{durch rechtsseitige $\alpha$-Verschiebung generierte Folge}. Sei \linebreak $(r_{j})^{\kappa-1}_{j=0} := \sarjk$.
  Dann führen wir folgende Bezeichnungen ein:
  \begin{align*}
    \Hsarn &:= H^{\sklam{r}}_{n} \quad \text{und} \quad \dHsarn := \widehat{H}^{\sklam{r}}_{n} 
  \end{align*}
  für alle $n \in \Zofkm$ sowie
  \begin{align*}
    \ysarjk &:= y^{\sklam{r}}_{j,k} \quad \text{und} \quad \zsarjk := z^{\sklam{r}}_{j,k} 
  \end{align*}
  für alle $j,k \in \Zokm$ mit $j \leq k$.

  \item [\rm{(b)}] $\saljk$ definiert gemäß $s_{\alj} := \alpha s_{j}-s_{j+1}$ für alle $j \in \Zokm$ heißt
  die aus $\sjk$ \textbf{durch linksseitige $\alpha$-Verschiebung generierte Folge}. Sei \linebreak $(r_{j})^{\kappa-1}_{j=0} := \saljk$.
  Dann führen wir folgende Bezeichnungen ein:
  \begin{align*}
    \Hsaln &:= H^{\sklam{r}}_{n} \quad \text{und} \quad \dHsaln := \widehat{H}^{\sklam{r}}_{n} 
  \end{align*}
  für alle $n \in \Zofkm$ sowie
  \begin{align*}
    \ysaljk &:= y^{\sklam{r}}_{j,k} \quad \text{und} \quad \zsaljk := z^{\sklam{r}}_{j,k} 
  \end{align*}
  für alle $j,k \in \Zokm$ mit $j \leq k$.
\end{itemize}
  Falls klar ist, von welchem $\sjk$ die Rede ist, lassen wir das \anf{$\sklam{s}$} als oberen Index weg.
\end{defi}

Es sei bemerkt, dass für $\alpha \in \C$, $\kappa \in \Na$ und eine Folge $\sjk$ aus $\Cqq$ für die in \thref{aspdef1} eingeführten Bezeichnungen $\Harn=-\alpha\Hn+\Kn$ bzw. $\Haln=\alpha\Hn-\Kn$ für alle $n\in\Zofkm$ erfüllt ist. 

Wir kommen nun zum zentralen Begriff dieses Abschnitts.

\begin{defi}	\thlabel{aspdef2} Seien $\alpha \in \C$, $\kappa \in \Noa$ und $\sjk$ eine Folge aus $\Cqq$.
\begin{itemize}
  \item [\rm{(a)}] $\Qsarjk$ definiert gemäß $Q^{\sklam{s}}_{\ar 2n} := \dHsn$ für alle $n \in \Zofk$ und im Fall $\kappa\geq1$ weiterhin $Q^{\sklam{s}}_{\ar 2n+1} := \dHsarn$ für alle $n \in \Zofkm$ heißt \textbf{rechtsseitige $\alpha$-Stieltjes-Parametrisierung} von $\sjk$.
  
  \item [\rm{(b)}] $\Qsaljk$ definiert gemäß $Q^{\sklam{s}}_{\al 2n} := \dHsn$ für alle $n \in \Zofk$ und im Fall $\kappa\geq1$ weiterhin $Q^{\sklam{s}}_{\al 2n+1} := \dHsaln$ für alle $n \in \Zofkm$ heißt \textbf{linksseitige $\alpha$-Stieltjes-Parametrisierung} von $\sjk$.
\end{itemize}
	Falls klar ist, von welchem $\sjk$ die Rede ist, lassen wir das \anf{$\sklam{s}$} als oberen Index weg.
\end{defi}

Wir können nun rekursiv die einzelnen Folgenglieder einer Matrizenfolge mithilfe ihrer $\alpha$-Stieltjes-Parametrisierung ausdrücken (vergleiche \cite[Remark 4.1]{ot226} und \cite[Remark 4.4]{ot226}).

\begin{bem}	\thlabel{aspbm5}
	Seien $\alpha \in \C$, $\kappa \in \Noa$ und $\sjk$ eine Folge aus $\Cqq$. \dgfa
	\begin{itemize}
		\item [\rm{(a)}] Sei $\Qarjk$ die rechtsseitige $\alpha$-Stieltjes-Parametrisierung von $\sjk$. Dann gelten $s_0 = Q_{\aro}$, im Fall $\kappa\geq1$ $s_1 = \alpha s_0 + Q_{\ar1}$, im Fall $\kappa\geq2$ 
		\begin{align*}
			s_{2n} = Q_{\ar2n}+z_{n,2n-1}\Hnm^{+}y_{n,2n-1}
		\end{align*}
		für alle $n\in\Zefk$ und im Fall $\kappa\geq3$
		\begin{align*}
			s_{2n+1} = \alpha s_{2n} + Q_{\ar2n+1}+z_{\arn,2n-1}\Harnm^{+}y_{\arn,2n-1}
		\end{align*}
		für alle $n\in\Zefkm$.
		\item [\rm{(b)}] Sei $\Qaljk$ die linksseitige $\alpha$-Stieltjes-Parametrisierung von $\sjk$. Dann gelten $s_0 = Q_{\alo}$, im Fall $\kappa\geq1$ $s_1 = \alpha s_0 - Q_{\al1}$, im Fall $\kappa\geq2$ 
		\begin{align*}
			s_{2n} = Q_{\al2n}+z_{n,2n-1}\Hnm^{+}y_{n,2n-1}
		\end{align*}
		für alle $n\in\Zefk$ und im Fall $\kappa\geq3$
		\begin{align*}
			s_{2n+1} = \alpha s_{2n} - Q_{\al2n+1}-z_{\aln,2n-1}\Halnm^{+}y_{\aln,2n-1}
		\end{align*}
		für alle $n\in\Zefkm$.
	\end{itemize}
\end{bem}

\bwanf Dies folgt unmittelbar aus \thref{aspdef2} und der Definition von $\dHn$ für alle $n\in\Zofk$ und $\dHarn$ bzw. $\dHaln$ für alle $n\in\Zofkm$. \bwend

\thref{aspdef2} führt uns zu folgender Beobachtung, die uns einen ersten Zusammenhang zwischen der rechtsseitigen und linksseitigen $\alpha$-Stieltjes-Parametrisierung liefert (vergleiche \cite[Remark 4.7]{ot226}).

\begin{bem}	\thlabel{aspbm2} Seien $\alpha \in \C$, $\kappa \in \Noa$ und $\sjk$ eine Folge aus $\Cqq$.
  Weiterhin sei $\Qarjk$ bzw. $\Qaljk$ die rechtsseitige bzw. linksseitige $\alpha$-Stieltjes-Parametrisierung von $\sjk$.
  Dann gilt $Q_{\alj} = (-1)^j Q_{\arj}$ für alle $j \in \Zok$. 
\end{bem}

\bwanf Sei $n\in\Zofkm$. Dann folgt aus der Definition von $\Harn$ und $\Haln$ sogleich $\Harn=-\Haln$. Weiterhin gelten $s_{\ar2n}=-s_{\al2n}$, $z_{\arn,2n-1}=-z_{\aln,2n-1}$ und \linebreak $y_{\arn,2n-1}=-y_{\aln,2n-1}$. Hieraus folgt nun
\begin{align*}
	\dHarn &= s_{\ar2n}-z_{\arn,2n-1}\Harnm^{+}y_{\arn,2n-1} \\
	& = -s_{\al2n}+z_{\aln,2n-1}\Halnm^{+}y_{\aln,2n-1} = -\dHaln.
\end{align*}
Hieraus folgt mithilfe von \thref{aspdef2} dann die Behauptung. \bwend

Für einen zweiten Zusammenhang zwischen der rechtsseitigen und linksseitigen $\alpha$-Stieltjes-Parametrisierung benötigen wir zunächst folgendes Lemma, das eine Weiterführung von \thref{asmlm1} für die aus der gegebenen Folge durch rechtsseitige bzw. linksseitige $\alpha$-Verschiebung generierte Folge darstellt.

\begin{lemma}	\thlabel{asplm1}
	Seien $\alpha \in \C$, $\kappa \in \Na$ und $\sjk$ eine Folge aus $\Cqq$. Weiterhin sei $t_j := (-1)^js_j$ für alle $j\in\Zok$. Dann gelten
	\begin{itemize}
		\item [\rm{(a)}] Sei $j\in\Zokm$. Dann gilt 
		\begin{align*}
			t_{-\alj} = (-1)^js_{\arj}.
		\end{align*}
		
		\item [\rm{(b)}] Sei $n\in\Zofkm$. Dann gilt
		\begin{align*}
			\Hmaln^{\sklam{t}} = \Vn\Hsarn\Vna.
		\end{align*}
		
		\item [\rm{(c)}] Seien $j,k \in \Zokm$ mit $j \leq k$. Dann gelten	
		\begin{align*}
			y^{\sklam{t}}_{-\alk,l} = (-1)^kV_{l-k}y^{\sklam{s}}_{\ark,l}
			\quad \text{und} \quad
			z^{\sklam{t}}_{-\alk,l} = (-1)^kz^{\sklam{s}}_{\ark,l}V^{\ast}_{l-k}.
		\end{align*}		
		
		\item [\rm{(d)}] Sei $n\in\Zofkm$. Dann gilt
		\begin{align*}
			\dHmaln^{\sklam{t}} = \dHsarn.
		\end{align*}
	\end{itemize}	
\end{lemma}

\bwanf Zu (a): Es gilt
\begin{align*}
	(-1)^js_{\arj} = -\alpha(-1)^js_j+(-1)^js_{j+1} = (-\alpha)t_j-t_{j+1} = t_{-\alj}.
\end{align*}

Zu (b)-(d): Dies folgt unter Beachtung von \thref{aspdef1} wegen (a) aus den Teilen (b), (d) und (e) von \thref{asmlm1}. \bwend

\begin{bem}	\thlabel{aspbm4} 
	Seien $\alpha \in \C$, $\kappa \in \Noa$, $\sjk$ eine Folge aus $\Cqq$ und $\Qsarjk$ die rechtsseitige $\alpha$-Stieltjes-Parametrisierung von $\sjk$. 
	Weiterhin sei $t_j := (-1)^js_j$ für alle $j\in\Zok$ und $\Qtmaljk$ die linksseitige $-\alpha$-Stieltjes-Parametrisierung von $\tjk$.
  	Dann gilt $Q^{\sklam{t}}_{-\alj} = Q^{\sklam{s}}_{\arj}$ für alle $j \in \Zok$. 
\end{bem}

\bwanf Unter Beachtung von \thref{aspdef1} folgt dies aus Teil (e) von \thref{asmlm1} und Teil (d) von \thref{asplm1}. \bwend

Nun kommen wir zur Einzigartigkeit der $\alpha$-Stieltjes-Parametrisierung einer Matrizenfolge (vergleiche \cite[Remark 4.3]{ot226} und \cite[Remark 4.6]{ot226}).

\begin{bem}	\thlabel{aspbm3}
	Seien $\alpha\in\C$, $\kappa \in \Noa$ und $\Qarjk$ bzw. $\Qaljk$ eine Folge aus $\Cqq$. Dann gibt es genau eine Folge $\sjk$ aus $\Cqq$, so dass $\Qarjk$ bzw. $\Qaljk$ die rechtsseitige bzw. linksseitige $\alpha$-Stieltjes-Parametrisierung von $\sjk$ ist.
\end{bem}

\bwanf Unter Beachtung von \thref{aspbm2} können wir uns auf den rechtsseitigen Fall beschränken. Seien durch rekursive Konstruktion $s_0 := Q_{\aro}$, im Fall $\kappa\geq1$ \linebreak $s_1 := \alpha s_0 + Q_{\ar1}$, im Fall $\kappa\geq2$
\begin{align*}
	s_{2n} := Q_{\ar2n}+z_{n,2n-1}\Hnm^{+}y_{n,2n-1}
\end{align*}
für alle $n\in\Zefk$ und im Fall $\kappa\geq3$
\begin{align*}
	s_{2n+1} := \alpha s_{2n} + Q_{\ar2n+1}+z_{\arn,2n-1}\Harnm^{+}y_{\arn,2n-1}
\end{align*}
für alle $n\in\Zefkm$. Dann folgt aus \thref{aspdef2} und der Definition von $\dHn$ für alle $n\in\Zofk$ bzw. $\dHarn$ für alle $n\in\Zofkm$, dass $\Qarjk$ die rechtsseitige $\alpha$-Stieltjes-Parametrisierung von $\sjk$ ist. Angenommen, es existiert eine von $\sjk$ verschiedene Folge $(t_j)^{\kappa}_{j=0}$ aus $\Cqq$, so dass $\Qarjk$ die rechtsseitige $\alpha$-Stieltjes-Parametrisierung von $(t_j)^{\kappa}_{j=0}$ ist. Sei $k := \min\gklam{j\in\Zok \;|\; s_j \neq t_j}$. Falls $k=0$ ist, gilt wegen \thref{aspdef2} dann 
\begin{align*}
	t_0 = \dH^{\sklam{t}}_0 = Q_{\aro} = s_0
\end{align*}
im Widerspruch zur Annahme. Falls $\kappa\geq1$ und $k=1$ erfüllt sind, gilt wegen \thref{aspdef2} und $s_0 = t_0$ dann
\begin{align*}
	t_1 = \alpha t_0 + \dH^{\sklam{t}}_{\aro} = \alpha s_0 + Q_{\ar1} = s_1
\end{align*}
im Widerspruch zur Annahme. Seien nun $\kappa\geq2$ und $k>1$ derart, dass ein $n\in\N$ mit $k=2n$ existiert. Wegen \thref{aspdef2} und $t_j = s_j$ für alle $j\in\Z{0}{k-1}$ gilt dann
\begin{align*}
	t_{2n} = \dH^{\sklam{t}}_{n} + z^{\sklam{t}}_{n,2n-1}\big(\Hnm^{\sklam{t}}\big)^{+}y^{\sklam{t}}_{n,2n-1} = Q_{\arn} + z^{\sklam{s}}_{n,2n-1}\big(\Hnm^{\sklam{s}}\big)^{+}y^{\sklam{s}}_{n,2n-1} = s_{2n}
\end{align*}
im Widerspruch zur Annahme. Seien nun $\kappa\geq3$ und $k>2$ derart, dass ein $n\in\N$ mit $k=2n+1$ existiert. Wegen \thref{aspdef2} und $t_j = s_j$ für alle $j\in\Z{0}{k-1}$ gilt dann
\begin{align*}
	t_{2n+1} &= \alpha t_{2n} + \dH^{\sklam{t}}_{\arn} + z^{\sklam{t}}_{\arn,2n-1}\big(\Harnm^{\sklam{t}}\big)^{+}y^{\sklam{t}}_{\arn,2n-1} \\
	&= \alpha s_{2n} + Q_{\ar2n+1}+z^{\sklam{s}}_{\arn,2n-1}\big(\Harnm^{\sklam{s}}\big)^{+}y^{\sklam{s}}_{\arn,2n-1} = s_{2n+1}
\end{align*}
im Widerspruch zur Annahme. Somit gilt $\sjk$ = $(t_j)^{\kappa}_{j=0}$. \bwend

Abschließend für dieses Kapitel stellen wir nun den Zusammenhang zwischen der $\alpha$-Stieltjes-Parametrisierung und der $\alpha$-Stieltjes-positiv bzw. -nichtnegativ Definitheit her. Zunächst betrachten wir den rechtsseitigen Fall.

\begin{satz}	\thlabel{aspsa1} Sei $\alpha \in \R$. \dgfa
\begin{itemize}
  \item [\rm{(a)}] Seien $\kappa \in \Noa$, $\sjk$ eine Folge aus $\Cqq$ und $\Qarjk$ die rechtsseitige $\alpha$-Stieltjes-Parametrisierung von $\sjk$. \dsfaa
  \begin{itemize}
    \item [\rm{(i)}] Es gilt $\sjk \in \Kqka$.
    \item [\rm{(ii)}] Es ist $\Qarjk$ eine Folge aus $\Cqq_{\geq}$ und im Fall $\kappa \geq 2$ gilt weiterhin ${\cal N}(Q_{\arj}) 
    \subseteq {\cal N}(Q_{\arj+1})$ für alle $j \in \Z{0}{\kappa-2}$.
  \end{itemize}
  
  \item [\rm{(b)}] Seien $m \in \No$, $\sjm$ eine Folge aus $\Cqq$ und $\Qarjm$ die rechtsseitige $\alpha$-Stieltjes-Parametrisierung von $\sjm$. \dsfaa
  \begin{itemize}
    \item [\rm{(i)}] Es gilt $\sjm \in \Keqma$.
    \item [\rm{(ii)}] Es ist $\Qarjm$ eine Folge aus $\Cqq_{\geq}$ und im Fall $m \geq 1$ gilt weiterhin ${\cal N}(Q_{\arj}) 
    \subseteq {\cal N}(Q_{\arj+1})$ für alle $j \in \Z{0}{m-1}$.
  \end{itemize}
  
  \item [\rm{(c)}] Seien $\kappa \in \Noa$, $\sjk$ eine Folge aus $\Cqq$ und $\Qarjk$ die rechtsseitige $\alpha$-Stieltjes-Parametrisierung von $\sjk$. \dsfaa
  \begin{itemize}
    \item [\rm{(i)}] Es gilt $\sjk \in \Kpqka$.
    \item [\rm{(ii)}] Es ist $\Qarjk$ eine Folge aus $\Cqq_{>}$.
  \end{itemize}

\end{itemize}  
\end{satz}

\bwanf Siehe \cite[Theorem 4.12]{ot226}. \bwend

Aus \thref{aspsa1} können wir sogleich die entsprechenden Aussagen für den linksseitigen Fall folgern.

\begin{satz}	\thlabel{aspsa2} Sei $\alpha \in \R$. \dgfa
\begin{itemize}
  \item [\rm{(a)}] Seien $\kappa \in \Noa$, $\sjk$ eine Folge aus $\Cqq$ und $\Qaljk$ die rechtsseitige $\alpha$-Stieltjes-Parametrisierung von $\sjk$. \dsfaa
  \begin{itemize}
    \item [\rm{(i)}] Es gilt $\sjk \in \Lqka$.
    \item [\rm{(ii)}] Es ist $\Qaljk$ eine Folge aus $\Cqq_{\geq}$ und im Fall $\kappa \geq 2$ gilt weiterhin ${\cal N}(Q_{\alj}) 
    \subseteq {\cal N}(Q_{\alj+1})$ für alle $j \in \Z{0}{\kappa-2}$.
  \end{itemize}
  
  \item [\rm{(b)}] Seien $m \in \No$, $\sjm$ eine Folge aus $\Cqq$ und $\Qaljm$ die rechtsseitige $\alpha$-Stieltjes-Parametrisierung von $\sjm$. \dsfaa
  \begin{itemize}
    \item [\rm{(i)}] Es gilt $\sjm \in \Leqma$.
    \item [\rm{(ii)}] Es ist $\Qaljm$ eine Folge aus $\Cqq_{\geq}$ und im Fall $m \geq 1$ gilt weiterhin ${\cal N}(Q_{\alj}) 
    \subseteq {\cal N}(Q_{\alj+1})$ für alle $j \in \Z{0}{m-1}$.
  \end{itemize}
  
  \item [\rm{(c)}] Seien $\kappa \in \Noa$, $\sjk$ eine Folge aus $\Cqq$ und $\Qaljk$ die rechtsseitige $\alpha$-Stieltjes-Parametrisierung von $\sjk$. \dsfaa
  \begin{itemize}
    \item [\rm{(i)}] Es gilt $\sjk \in \Lpqka$.
    \item [\rm{(ii)}] Es ist $\Qaljk$ eine Folge aus $\Cqq_{>}$.
  \end{itemize}

\end{itemize}  
\end{satz}

\bwanf Dies folgt wegen \thref{asmbm3} und \thref{aspbm4} aus \thref{aspsa1}. \bwend

%
%
%
%

\subsection[Die kanonische Hankel-Parametrisierung von Matrizenfolgen]{Die kanonische Hankel-Parametrisierung von \\ Matrizenfolgen}

In diesem Abschnitt geben wir einen kurzen Überblick zu der kanonischen Hankel"=Parametrisierung von Matrizenfolgen. Wir benutzen eine allgemeinere Definition jener Parametrisierung, wie sie auch schon in \cite[Definition 2.1]{Pert} formuliert wurde. Der Begriff der kanonischen Hankel-Parametrisierung wurde bereits in \cite[Definition 2.28]{13} eingeführt, wo man für $\kappa \in \Noa$ eine Matrizenfolge $\sjzk$ aus $\Cqq$ voraussetzte. Die Resultate aus \cite[Chapter 2]{13} lassen sich aber auch auf unseren allgemeineren Fall übertragen. Weitere grundlegende Ausführungen zu jener Parametrisierung findet man in \cite[Chapter 2]{15} und \cite[Chapter 3]{ot226}. 
Bevor wir uns nun der Definition der kanonischen Hankel-Parametrisierung widmen, führen wir eine weitere Bezeichnung ein. 

\begin{bez}	\thlabel{mpbz1}
	Seien $\kappa \in \Noa$ und $\sjk$ eine Folge aus $\Cpq$. Weiterhin seien $\Lambda^{\sklam{s}}_{0} := 0_{\pq}$ und
	\begin{align*}
  		\Lambda^{\sklam{s}}_{n} &:= z^{\sklam{s}}_{n,2n-1}\big(\Hsnm\big)^{+}y^{\sklam{s}}_{n+1,2n}+z^{\sklam{s}}_{n+1,2n}\big(\Hsnm\big)^{+}y^{\sklam{s}}_{n,2n-1} \\
  		&\;\quad -z^{\sklam{s}}_{n,2n-1}\big(\Hsnm\big)^{+}\Ksnm\big(\Hsnm\big)^{+}y^{\sklam{s}}_{n,2n-1}.
	\end{align*}
	für alle $n \in \Zefk$. Falls klar ist, von welchem $\sjk$ die Rede ist, lassen wir das \anf{$\sklam{s}$} als oberen Index weg. 
\end{bez}

\begin{defi}	\thlabel{mpdef1}
  	Seien $\kappa \in \Noa$ und $\sjk$ eine Folge aus $\Cpq$. Weiterhin seien 
	\begin{align*}
    	D^{\sklam{s}}_{n} := \dHsn
  	\end{align*}
  	für alle $n\in\Zofk$ und im Fall $\kappa\geq1$
   	\begin{align*}
    	C^{\sklam{s}}_{n} := s_{2n-1}-\Lambda^{\sklam{s}}_{n-1}
  	\end{align*}
  	für alle $n \in \Zefkp$.  
  	Dann heißt im Fall $\kappa=0$ die Folge $(D^{\sklam{s}}_{n})^{0}_{n=0}$ bzw. im Fall $\kappa \geq 1$ das Paar $\CDsfk$ die \textbf{kanonische Hankel-Parametrisierung} von $\sjk$. Falls klar ist, von welchem $\sjk$ die Rede ist, lassen wir das \anf{$\sklam{s}$} als oberen Index weg.
\end{defi}

Seien $\kappa \in \Noa$, $\sjk$ eine Folge aus $\Cpq$ und $(D_{n})^{0}_{n=0}$ bzw. im Fall $\kappa \geq 1$ das Paar $\CDfk$ die kanonische Hankel-Parametrisierung von $\sjk$. Anhand \thref{mpdef1} erkennt man, dass für $n\in\Zofk$ die Matrix $D_n$ nur von $\sjn$ abhängt und im Fall $\kappa\geq1$ für $n \in \Zefkp$ die Matrix $C_n$ nur von $\sjnm$ abhängt. Dies führt uns auf folgende Behauptung.

\begin{bem}	\thlabel{mpbm3}
	Seien $\kappa \in \Na$ und $\sjk$ eine Folge aus $\Cpq$. Weiterhin sei \linebreak $\CDfk$ die kanonische Hankel-Parametrisierung von $\sjk$. Dann ist \linebreak $(D_{n})^{0}_{n=0}$ die kanonische Hankel-Parametrisierung von $(s_j)^{0}_{j=0}$. Sei im Fall $\kappa\geq2$ nun $m\in\Zekm$. Dann ist $[(C_{n})^{\fklam{m+1}}_{n=1},$ $(D_{n})^{\fklam{m}}_{n=0}]$ die kanonische Hankel-Parametrisierung von $\sjm$.
\end{bem}

Wir können nun rekursiv die einzelnen Folgenglieder einer Matrizenfolge mithilfe ihrer Hankel-Parametrisierung ausdrücken (vergleiche \cite[Remark 2.4]{15}). 

\begin{bem}	\thlabel{mpbm4}
	Seien $\kappa \in \Noa$ und $\sjk$ eine Folge aus $\Cpq$. Weiterhin sei \linebreak $\CDfk$ die kanonische Hankel-Parametrisierung von $\sjk$. Dann gelten $s_0 = D_0$, im Fall $\kappa\geq1$
	\begin{align*}
		s_{2n-1} = C_n + \Lambda_{n-1}
	\end{align*}
	für alle $n\in\N$ mit $2n-1\leq\kappa$ und im Fall $\kappa\geq2$
	\begin{align*}
		s_{2n} = D_n + z_{n,2n-1}\Hnm^{+}y_{n,2n-1}
	\end{align*}
	für alle $n\in\Zefk$. 
\end{bem}

\bwanf Dies folgt unmittelbar aus \thref{mpdef1} und der Definition von $\dHn$ für alle $n\in\Zofk$. \bwend

Folgendes Resultat liefert uns die Einzigartigkeit der kanonischen Hankel"=Parametrisierung einer Matrizenfolge (vergleiche \cite[Remark 2.29]{13}).

\begin{bem}	\thlabel{mpbm2}
	Seien $\kappa \in \Noa$ und $(D_{n})^{\fklam{\kappa}}_{n=0}$ sowie im Fall $\kappa\geq1$ auch $(C_{n})^{\fklam{\kappa+1}}_{n=1}$ eine Folge aus $\Cpq$. Dann existiert genau eine Folge $\sjk$ aus $\Cpq$, so dass die Folge $(D_{n})^{0}_{n=0}$ bzw. im Fall $\kappa \geq 1$ das Paar $\CDfk$ die kanonische Hankel-Parametrisierung von $\sjk$ ist.
\end{bem}

\bwanf Seien durch rekursive Konstruktion $s_0 := D_0$, im Fall $\kappa\geq1$
\begin{align*}
	s_{2n-1} := C_n + \Lambda_{n-1}
\end{align*}
für alle $n\in\N$ mit $2n-1\leq\kappa$ und im Fall $\kappa\geq2$
\begin{align*}
	s_{2n} := D_n + z_{n,2n-1}\Hnm^{+}y_{n,2n-1}
\end{align*}
für alle $n\in\Zefk$. Dann folgt aus \thref{mpdef1} und der Definition von $\dHn$ für alle $n\in\Zofk$, dass die Folge $(D_{n})^{0}_{n=0}$ bzw. im Fall $\kappa \geq 1$ das Paar $\CDfk$ die kanonische Hankel-Parametrisierung von $\sjk$ ist. Angenommen, es existiert eine von $\sjk$ verschiedene Folge $(t_j)^{\kappa}_{j=0}$ aus $\Cpq$, so dass die Folge $(D_{n})^{0}_{n=0}$ bzw. im Fall $\kappa \geq 1$ das Paar $\CDfk$ die kanonische Hankel-Parametrisierung von $(t_j)^{\kappa}_{j=0}$ ist. Sei $k := \min\gklam{j\in\Zok \;|\; s_j \neq t_j}$. Falls $k=0$ ist, gilt wegen \thref{mpdef1} dann 
\begin{align*}
	t_0 = \dH^{\sklam{t}}_0 = D_0 = s_0
\end{align*}
im Widerspruch zur Annahme. Seien nun $\kappa\geq1$ und $k$ derart, dass ein $n\in\N$ mit $k=2n-1$ existiert. Wegen \thref{mpdef1} und $t_j = s_j$ für alle $j\in\Z{0}{k-1}$ gilt dann
\begin{align*}
	t_{2n-1} = C_n-\Lambda^{\sklam{t}}_{n-1} = C_n-\Lambda^{\sklam{s}}_{n-1} = s_{2n-1}
\end{align*}
im Widerspruch zur Annahme. Seien nun $\kappa\geq2$ und $k>0$ derart, dass ein $n\in\N$ mit $k=2n$ existiert. Wegen \thref{mpdef1} und $t_j = s_j$ für alle $j\in\Z{0}{k-1}$ gilt dann
\begin{align*}
	t_{2n} = D_n + z^{\sklam{t}}_{n,2n-1}\big(\Hnm^{\sklam{t}}\big)^{+}y^{\sklam{t}}_{n,2n-1} = D_n + z^{\sklam{s}}_{n,2n-1}\rklam{\Hnm^{\sklam{s}}}^{+}y^{\sklam{s}}_{n,2n-1} = s_{2n}
\end{align*}
im Widerspruch zur Annahme. Somit gilt $\sjk$ = $(t_j)^{\kappa}_{j=0}$. \bwend

Abschließend für dieses Kapitel stellen wir einen Zusammenhang zwischen der kanonischen Hankel-Parametrisierung und der Hankel-positiv Definitheit her.

\begin{satz}	\thlabel{mpsa1}
  Seien $\kappa \in \Noa$ und $\sjk$ eine Folge aus $\Cqq$. Weiterhin sei $(D_{n})^{0}_{n=0}$ bzw. im Fall $\kappa \geq 1$ $\CDfk$ die kanonische Hankel-Parametrisierung von $\sjk$. \dsfaa
  \begin{itemize}
    \item [\rm{(i)}] Es gilt $\sjfk \in \Hpqfk$.
    \item [\rm{(ii)}] Es ist $(D_{n})^{\fklam{\kappa}}_{n=0}$ eine Folge aus $\Cqq_{>}$ und für $\kappa \geq 1$ ist $(C_{n})^{\fklam{\kappa}}_{n=1}$ eine Folge aus $\Cqq_{H}$.
  \end{itemize}
\end{satz}

\bwanf Im Fall $\kappa=0$ folgt dies direkt aus $D_{0} = s_{0} = H_{0}$. 
Der Fall $\kappa \in \N$ folgt unter Beachtung von \thref{mpbm3} aus \cite[Proposition 2.10(d)]{15}. Der Fall $\kappa=\infty$ wurde in \cite[Proposition 2.15(c)]{15} gezeigt. \bwend 

Es sei bemerkt, dass man \thref{mpsa1} unter Beachtung von \thref{aspdef1} auch auf die kanonische Hankel-Parametrisierung von der durch rechtsseitige bzw. linksseitige $\alpha$-Verschiebung generierten Folge anwenden kann. Dies wird uns später unter Beachtung von \thref{asmdef2} und \thref{asmdef3} eine Verbindung zu den rechtsseitig bzw. linksseitig $\alpha$-Stieltjes-positiv definiten Folgen liefern.

%
%
%
%

\subsection{Das Favard-Paar bezüglich Hankel-positiv definiter Folgen} \label{chapfp}

Im engen Zusammenhang zur kanonischen Hankel-Parametrisierung einer Matrizenfolge, die bis zu einem gewissen Folgenglied Hankel-positiv definit ist, steht das Favard-Paar bezüglich jener Folge. Dieses Paar wurde in \cite[Chapter 8]{CR1} für eine unendliche Matrizenfolge eingeführt und behandelt. Wir betrachten hingegen den allgemeineren Fall einer Matrizenfolge mit beliebig vielen Folgengliedern. Die Resultate aus \cite[Chapter 8]{CR1} lassen sich leicht auf den Fall einer endlichen Matrizenfolge übertragen.

Bevor wir zum zentralen Begriff dieses Abschnitts kommen, benötigen wir noch folgende Bemerkung. Sie erlaubt uns die Betrachtung der linken Schur-Komplemente der zu einer Hankel-positiv definiten Folge gehörigen Block-Hankel-Matrizen. Wir werden im weiteren Verlauf nicht mehr explizit dieses Resultat aufrufen.

\begin{bem}	\thlabel{fpbm5}
	Seien $\kappa\in\Noa$ und $\sjk\in\Hpqfk$. Dann ist $\dHn$ für alle $n\in\Zofk$ positiv hermitesch und insbesondere regulär.
\end{bem}

\bwanf Wegen Teil (a) von \thref{asmdef1} und \thref{asmbm1} ist $\Hn$ für alle $n\in\Zofk$ positiv hermitesch. Hieraus folgt wegen \thref{amlm4} dann, dass $\dHn$ für alle $n\in\Zofk$ positiv hermitesch und somit insbesondere regulär ist. \bwend

\begin{defi}	\thlabel{fpdef1}
	Seien $\kappa\in\Noa$ und $\sjk$ eine Folge aus $\Cqq$ mit \linebreak $\sjfkm\in\Hpqfkm$, falls $\kappa\geq1$. Weiterhin seien $B^{\sklam{s}}_0 := s_0$, im Fall $\kappa\geq1$ \linebreak $A^{\sklam{s}}_0 := s_1s^{-1}_0$, im Fall $\kappa\geq2$
	\begin{align*}
		B^{\sklam{s}}_n := \brklam{\dHsnm}^{-1}\dHsn
	\end{align*}
	für alle $n\in\Zefk$ und im Fall $\kappa\geq3$
	\begin{align*}
		A^{\sklam{s}}_n := \begin{pmatrix} -z^{\sklam{s}}_{n,2n-1}\brklam{\Hsnm}^{-1} & \Iq \end{pmatrix} \Ksn \begin{pmatrix} -\brklam{\Hsnm}^{-1}y^{\sklam{s}}_{n,2n-1} \\ \Iq \end{pmatrix} \brklam{\dHn^{\sklam{s}}}^{-1}
	\end{align*}
	für alle $n\in\Zefkm$. Dann heißt im Fall $\kappa=0$ die Folge $(B^{\sklam{s}}_n)^{0}_{n=0}$ bzw. im Fall $\kappa\geq1$ das Paar $\ABsfk$ das \textbf{Favard-Paar} bezüglich $\sjk$. Falls klar ist, von welchem $\sjk$ die Rede ist, lassen wir $\sklam{s}$ als oberen Index weg.
\end{defi}

Seien $\kappa\in\Noa$, $\sjk$ eine Folge aus $\Cqq$ mit $\sjfkm\in\Hpqfkm$, falls $\kappa\geq1$, und $(B_n)^{0}_{n=0}$ bzw. im Fall $\kappa\geq1$ das Paar $\ABfk$ das Favard-Paar bezüglich $\sjk$. Anhand \thref{fpdef1} erkennt man, dass für $n\in\Zofk$ die Matrix $B_n$ nur von $\sjn$ abhängt und im Fall $\kappa\geq1$ für $n \in \Zofkm$ die Matrix $A_n$ nur von $\sjne$ abhängt. Dies führt uns auf folgende Behauptung.

\begin{bem}	\thlabel{fpbm2}
	Seien $\kappa\in\Na$ und $\sjk$ eine Folge aus $\Cqq$ mit \linebreak $\sjfkm\in\Hpqfkm$. Weiterhin sei $\ABfk$ das Favard-Paar bezüglich $\sjk$. Dann ist $(B_{n})^{0}_{n=0}$ das Favard-Paar bezüglich $(s_j)^{0}_{j=0}$. Sei im Fall $\kappa\geq2$ nun $m\in\Zekm$. Dann ist $[(A_{n})^{\fklam{m-1}}_{n=0},$ $(B_{n})^{\fklam{m}}_{n=0}]$ das Favard-Paar bezüglich $\sjm$.
\end{bem}

Wir befassen uns nun mit einem direkten Zusammenhang zwischen der kanonischen Hankel-Parametrisierung einer Matrizenfolge, die bis zu einem gewissen Folgenglied Hankel-positiv ist, und dem Favard-Paar bezüglich jener Folge (im Fall $\kappa=\infty$ vergleiche \cite[Proposition 8.3]{CR1}). Wir werden im späteren Verlauf der Arbeit bevorzugt das Favard-Paar anstelle der kanonischen Hankel"=Parametrisierung verwenden.

\begin{bem}	\thlabel{fpbm4}
	Sei $s_0\in\Cqq$. Dann stimmt die kanonische Hankel"=Parametrisierung von $(s_j)^0_{j=0}$ mit dem Favard-Paar bezüglich $(s_j)^0_{j=0}$ überein.
\end{bem}

\bwanf Unter Beachtung von $\dH_0 = s_0$ folgt dies unmittelbar aus \thref{mpdef1} und \thref{fpdef1}. \bwend

\begin{satz}	\thlabel{fpsa1}
	Seien $\kappa\in\Na$ und $\sjk$ eine Folge aus $\Cqq$ mit $\sjfkm\in\Hpqfkm$. Weiterhin seien $\ABfk$ das Favard-Paar bezüglich $\sjk$ und $\CDfk$ die kanonische Hankel-Parametrisierung von $\sjk$. \dgfa
	\begin{itemize}
		\item [\rm{(a)}] Es gelten $B_0 = D_0$,
		\begin{align*}
			A_n = C_{n+1}D^{-1}_n
		\end{align*}
		für alle $n\in\Zofkm$ und im Fall $\kappa\geq2$
		\begin{align*}
			B_n = D^{-1}_{n-1}D_n
		\end{align*}
		für alle $n\in\Zefk$.
		\item [\rm{(b)}] Es gelten
		\begin{align*}
			D_n = \prodr^n_{j=0} B_j
		\end{align*}
		für alle $n\in\Zofk$ und
		\begin{align*}
			C_n = A_{n-1} \prodr^{n-1}_{j=0} B_j
		\end{align*}
		für alle $n\in\Zefkp$.
	\end{itemize}
\end{satz}

\bwanf Zu (a): Wegen \thref{mpdef1} gilt $C_1 = s_1$ und im Fall $\kappa\geq3$
\begin{align*}
	&\ C_{n+1} = s_{2n+1} - \Lambda_n \\
	&= s_{2n+1} - z_{n,2n-1}\Hnm^{-1}y_{n+1,2n} - z_{n+1,2n}\Hnm^{-1}y_{n,2n-1} + z_{n,2n-1}\Hnm^{-1}\Knm\Hnm^{-1}y_{n,2n-1} \\
	&= \begin{pmatrix} -z_{n,2n-1}\Hnm^{-1} & \Iq \end{pmatrix} \begin{pmatrix} K_{n-1} & y_{n+1,2n} \\ z_{n+1,2n} & s_{2n+1} \end{pmatrix} \begin{pmatrix} -\Hnm^{-1}y_{n,2n-1} \\ \Iq \end{pmatrix} \\
	&= \begin{pmatrix} -z_{n,2n-1}\Hnm^{-1} & \Iq \end{pmatrix} K_n \begin{pmatrix} -\Hnm^{-1}y_{n,2n-1} \\ \Iq \end{pmatrix}
\end{align*}
für alle $n\in\Zefkm$ sowie $D_n = \dHn$ für alle $n\in\Zofk$. Hieraus folgt wegen \thref{fpdef1} dann die Behauptung.

Zu (b): Wegen \thref{mpdef1} und (a) gelten $B_0 = D_0$ und im Fall $\kappa\geq2$
\begin{align*}
	\prodr^n_{j=0} B_j = D_0 \prodr^n_{j=1} \rklam{D^{-1}_{j-1}D_j} = D_n
\end{align*}
für alle $n\in\Zefk$. Hieraus folgt wegen (a) dann
\begin{align*}
	A_{n-1} \prodr^{n-1}_{j=0} B_j = C_nD^{-1}_{n-1}D_{n-1} = C_n
\end{align*}
für alle $n\in\Zefkp$. \bwend

Mithilfe der Einzigartigkeit der kanonischen Hankel-Parametrisierung einer beliebigen Matrizenfolge können wir nun die Einzigartigkeit des Favard-Paar bezüglich einer Matrizenfolge, die bis zu einem gewissen Folgenglied Hankel-positiv ist, zeigen (im Fall $\kappa=\infty$ vergleiche \cite[Proposition 8.5]{CR1}).

\begin{bem}	\thlabel{fpbm1}
	Seien $\kappa\in\Na$, $\Dnfk$ eine Folge aus $\Cqq$ derart, dass $D_n\in\Cqq_>$ für alle $n\in\Zofkm$ erfüllt ist, und $\Cnfk$ eine Folge aus $\Cqq$ derart, dass im Fall $\kappa\geq2$ $\;C_n\in\Cqq_H$ für alle $n\in\Zefkm$ erfüllt ist. Weiterhin seien $B_0 := D_0$,
	\begin{align*}
		A_n := C_{n+1}D^{-1}_n
	\end{align*}
	für alle $n\in\Zofkm$ und im Fall $\kappa\geq2$
	\begin{align*}
		B_n := D^{-1}_{n-1}D_n
	\end{align*}
	für alle $n\in\Zefk$. Dann existiert genau eine Folge $\sjk$, so dass $\CDfk$ die kanonische Hankel-Parametrisierung von $\sjk$ ist. Weiterhin gilt \linebreak $\sjfkm\in\Hpqfkm$ und $\ABfk$ ist das Favard-Paar bezüglich $\sjk$.
\end{bem}

\bwanf Wegen \thref{mpbm2} existiert genau eine Folge $\sjk$, so dass $\CDfk$ die kanonische Hankel-Parametrisierung von $\sjk$ ist. Hieraus folgt wegen \thref{mpbm3} dann, dass im Fall $\kappa=1$ die Folge $(D_n)^0_{n=0}$ bzw. im Fall $\kappa\geq2$ das Paar $[(C_{n})^{\fklam{\kappa}}_{n=1},$ $(D_{n})^{\fklam{\kappa-1}}_{n=0}]$ die kanonische Hankel-Parametrisierung von $\sjfkm$ ist. Weiterhin ist $(D_{n})^{\fklam{\kappa-1}}_{n=0}$ eine Folge aus $\Cqq_{>}$ und für $\kappa \geq 2$ ist $(C_{n})^{\fklam{\kappa-1}}_{n=1}$ eine Folge aus $\Cqq_{H}$. Wegen \thref{mpsa1} gilt dann $\sjfkm\in\Hpqfkm$. Hieraus folgt wegen \thref{fpsa1} nun, dass $\ABfk$ das Favard-Paar bezüglich $\sjk$ ist. \bwend

Es sei bemerkt, dass für eine rechtsseitig bzw. linksseitig $\alpha$-Stieltjes-positive Folge unter Beachtung von \thref{asmdef2} bzw. \thref{asmdef3} und \thref{aspdef1} jene Folge und die aus jener Folge durch rechtsseitige bzw. linksseitige $\alpha$-Verschiebung generierte Folge jeweils Matrizenfolgen sind, die bis zu einem gewissen Folgenglied Hankel-positiv sind. Somit können wir später jeweils das Favard-Paar bezüglich jener Folge und der aus jener Folge durch rechtsseitige bzw. linksseitige $\alpha$-Verschiebung generierten Folge betrachten.

%
%
%
%

\subsection{Über einige zu Matrizenfolgen gehörige Matrixpolynome} \label{chapmp}

Wir befassen uns nun mit zwei Folgen von Matrixpolynomen bezüglich einer gegebenen Matrizenfolge, wobei wir speziell auf bis zu einem gewissen Folgenglied Hankel-positiv definite Matrizenfolgen eingehen werden. In diesem Spezialfall können wir dann die betrachteten Folgen von Matrixpolynomen mithilfe von Hankelmatrizen direkt darstellen sowie einen ersten Zusammenhang zu der kanonischen Hankel-Parametrisierung und dem Favard-Paar finden. Zuerst befassen wir uns mit dem monischen links-orthogonalen System von Matrixpolynomen und verwenden die Herangehensweise von \cite[Chapter 3]{Pert} und \cite[Chapter 5]{15}. Dort wurde für die gegebene Matrizenfolge $\sjk$ aus $\Cqq$, wobei $\kappa\in\Noa$, vorausgesetzt, dass $\kappa$ gerade oder gleich unendlich ist. Wir werden den Begriff auf beliebige $\kappa$ erweitern. Danach befassen wir uns mit dem linken System von Matrixpolynomen zweiter Art, das für gerade $\kappa$ oder $\kappa=\infty$ in \cite[Chapter 4]{Pert} behandelt wurde. Es sei bemerkt, dass man auch vom monischen rechts-orthogonalen System von Matrixpolynomen ausgehen kann. Dies wurde schon in \cite[Chapter 5]{15} behandelt und das rechte System von Matrixpolynomen kann dann auf die gleiche Weise eingeführt werden, wie im linken Fall. Der Übersicht halber beschränken wir uns aber auf den linken Fall. Zunächst stellen wir einige Bezeichnungen bezüglich Matrixpolynomen bereit.

\begin{bez}	\thlabel{mpbz2}
	Sei $P$ ein $\pq$-Matrixpolynom. Dann gibt es eine eindeutige Folge $(P^{[j]})^{\infty}_{j=0}$ aus $\Cpq$, die $P(z)=\sum^{\infty}_{j=0}z^{j}P^{[j]}$ für alle $z\in\C$ erfüllt.
	Wir bezeichnen mit
	\begin{align*}
  	\deg P := \sup \gklam{j \in \No \;\big|\; P^{[j]} \neq 0_{\pq}}
	\end{align*}
	den Grad von P. Im Fall $\{j \in \No \;|\; P^{[j]} \neq 0_{\pq}\} = \emptyset$, also $P \equiv 0_{\pq}$, ist dann \linebreak $\deg P = -\infty$.
	Im Fall $\deg P \geq 0$ heißt $P^{[\deg P]}$ der Leitkoeffizient von $P$.
\end{bez}

Nun kommen wir zum ersten zentralen Begriff dieses Abschnitts.

\begin{defi}	\thlabel{mpdef2}
  Seien $\kappa \in \Noa$ und $\sjk$ eine Folge aus $\Cqq$. Die Folge $\Pnfk$ von \textit{q}$\times$\textit{q}-Matrixpolynomen heißt ein
  \textbf{monisches links-orthogonales System von Matrixpolynomen} bezüglich $\sjk$, falls
  \begin{itemize}
    \item [\rm{(i)}] $\deg P_{n} = n$ für alle $n \in \Zofkp$,
    \item [\rm{(ii)}] $P_{n}$ hat den Leitkoeffizienten $\Iq$ für alle $n \in \Zofkp$ 
  \end{itemize}
  und im Fall $\kappa\in\Na$ weiterhin
  \begin{itemize}
    \item [\rm{(iii)}] $\sum_{l \in \Z{0}{j}, m \in \Z{0}{k}} P^{[l]}_{j} s_{l+m} \brklam{P^{[m]}_{k}}^{\ast} = 0_{\qq}$
      für alle $j,k \in \Zofkp$ mit $j \neq k$
  \end{itemize}
  erfüllt sind.
\end{defi}

Falls $\kappa \in \Na$ gerade ist, lässt sich \rm{(iii)} von \thref{mpdef2} auch schreiben als
\begin{itemize}
  \item [\rm{(iii')}] $\begin{pmatrix} P^{[0]}_{j} & \ldots & P^{[n]}_{j} \end{pmatrix} H_{n} \begin{pmatrix} P^{[0]}_{k} & \ldots & P^{[n]}_{k} \end{pmatrix}^{\ast} = 0_{\qq}$
    für alle $j,k \in \Zofk$ mit $j \neq k$, wobei $n:=\max\gklam{j,k}$.
\end{itemize}

Wir zeigen nun die Existenz eines monischen links-orthogonalen Systems von Matrixpolynomen bezüglich einer Matrizenfolge, die bis zu einem gewissen Folgenglied Hankel-positiv ist. Die Aussage wurde in \cite[Proposition 3.4]{Pert} für den speziellen Fall einer gegebenen Folge aus $\Hpqk$ mit $\kappa\in\Na$ formuliert. Die Eindeutigkeit wurde in \cite[Proposition 5.6(c)]{15} dargestellt und die Identität für den dualen Fall eines monischen rechts-orthogonalen Systems von Matrixpolynomen wurde im Beweis von \cite[Proposition 5.6(b)]{15} gezeigt. Wir zeigen analog die folgende Identität für den linken Fall.

\begin{satz}	\thlabel{mpsa2}
  Seien $\kappa \in \Noa$ und $\sjk$ eine Folge aus $\Cqq$ mit $\sjfkm \in \Hpqfkm$, falls $\kappa \geq 1$. 
  Dann existiert genau ein monisches links-orthogonales System von Matrixpolynomen $\Pnfk$ bezüglich $\sjk$ und es gilt
  \begin{align*}
    \begin{pmatrix} P^{[0]}_{n} & \ldots & P^{[n]}_{n} \end{pmatrix} = \begin{cases} \Iq & \text{falls } n=0 \\
    \begin{pmatrix} -z_{n,2n-1}H^{-1}_{n-1} & \Iq \end{pmatrix} & \text{falls } n \in \Zefkp. \end{cases}
  \end{align*}
\end{satz}

\bwanf Wir konstruieren eine Folge $\Pnfk$ von Matrixpolynomen, die den Bedingungen aus \thref{mpdef2} gerecht wird. Wegen Teil (i) von \thref{mpdef2} gilt
\begin{align*}
	P_n(z)=\sum^{n}_{j=0}z^jP^{[j]}_n
\end{align*}
für alle $z\in\C$ und $n\in\Zofkp$. Es genügt also für jedes $n\in\Zofkp$ Matrizenfolgen $(P^{[j]}_n)^n_{j=0}$ zu finden, die die Teile (ii) und im Fall $\kappa \geq 1$ auch (iii) von \thref{mpdef2} erfüllen, um die Existenz des monischen links-orthogonalen Systems von Matrixpolynomen bezüglich $\sjk$ zu zeigen. Gibt es zudem für jedes $n\in\Zofkp$ nur genau eine Matrizenfolge $(P^{[j]}_n)^n_{j=0}$, die die Teile (ii) und im Fall $\kappa \geq 1$ auch (iii) von \thref{mpdef2} erfüllt, so haben wir auch die Eindeutigkeit des monischen links-orthogonalen Systems von Matrixpolynomen bezüglich $\sjk$ gezeigt.

Es folgt aus Teil (ii) von \thref{mpdef2} direkt $P^{[0]}_0 = \Iq$ und somit ist für den Fall $\kappa=0$ alles gezeigt.
Sei nun $\kappa \geq 1$. Weiterhin sei
\begin{align*}
  Z^{[P_n]}_n := \begin{pmatrix} P^{[0]}_n & \ldots & P^{[n]}_n \end{pmatrix}
\end{align*}
für alle $n \in \Zofkp$. Wegen Teil (ii) von \thref{mpdef2} existiert eine Matrix $t_n \in \C^{q \times nq}$ mit
\begin{align} \label{mpsa2bw1}
  Z^{[P_n]}_n = \begin{pmatrix} t_n & \Iq \end{pmatrix}
\end{align}
für alle $n \in \Zefkp$. Sei
\begin{align*}
  Z^{[P_k]}_n := \begin{pmatrix} Z^{[P_k]}_k & 0_{q \times (n-k)q} \end{pmatrix}
\end{align*}
für alle $n,k \in \Zofkp$ mit $k < n$. Weiterhin sei
\begin{align*}
  U_n := \begin{pmatrix} Z^{[P_0]}_n \\ \vdots \\ Z^{[P_n]}_n \end{pmatrix}
\end{align*}
für alle $n \in \Zefkp$. Wegen \fref{mpsa2bw1} ist $U_n$ für alle $n \in \Zefkp$ eine Matrix
von unterer Dreiecksgestalt und Einsen auf der Hauptdiagonale, somit also insbesondere regulär. Es existiert eine Matrix 
$V_n \in \C^{nq \times nq}$ mit
\begin{align} \label{mpsa2bw2}
  U^{-\ast}_n \begin{pmatrix} I_{nq} \\ 0_{q \times nq} \end{pmatrix} = \begin{pmatrix} V_n \\ 0_{q \times nq} \end{pmatrix}
\end{align}
für alle $n \in \Zefkp$. Sei 
\begin{align*}
  \tHn := \begin{pmatrix} H_{n-1} & y_{n,2n-1} \\ z_{n,2n-1} & 0_{\qq} \end{pmatrix}
\end{align*}
für alle $n \in \Zefkp$. Dann gilt wegen Teil (iii) von \thref{mpdef2}
\begin{align} \label{mpsa2bw4}
  Z^{[P_n]}_n \tHn \brklam{Z^{[P_k]}_n}^{\ast} = \sum_{l \in \Zon, m \in {\mathbb Z}_{0,k}} P^{[l]}_{n} s_{l+m} \brklam{P^{[m]}_{k}}^{\ast} = 0_{\qq}
\end{align}
für alle $n,k \in \Zefkp$ mit $k < n$. Wegen \fref{mpsa2bw1} - \fref{mpsa2bw4} gilt dann
\begin{align*}
  t_n H_{n-1} & = \begin{pmatrix} t_n & \Iq \end{pmatrix} \begin{pmatrix} H_{n-1} & y_{n,2n-1} \\ z_{n,2n-1} & 0_{\qq} \end{pmatrix} \begin{pmatrix} I_{nq} \\ 0_{q \times nq} \end{pmatrix} - z_{n,2n-1}\\
  & = Z^{[P_n]}_n \tHn U^{\ast}_n U^{-\ast}_n \begin{pmatrix} I_{nq} \\ 0_{q \times nq} \end{pmatrix} - z_{n,2n-1} \\
  & = \begin{pmatrix} Z^{[P_n]}_n \tHn \brklam{Z^{[P_0]}_n}^{\ast} & \ldots & Z^{[P_n]}_n \tHn \brklam{Z^{[P_n]}_n}^{\ast} \end{pmatrix} \begin{pmatrix} V_n \\ 0_{q \times nq} \end{pmatrix} - z_{n,2n-1} \\
  & = \begin{pmatrix} 0_{q \times nq} & Z^{[P_n]}_n \tHn \brklam{Z^{[P_n]}_n}^{\ast} \end{pmatrix} \begin{pmatrix} V_n \\ 0_{q \times nq} \end{pmatrix} - z_{n,2n-1} = -z_{n,2n-1}
\end{align*}
für alle $n \in \Zefkp$. Hieraus folgt unter Beachtung von \thref{asmbm1} wegen \linebreak $\sjfkm \in \Hpqfkm$ dann
\begin{align*}
  t_n = -z_{n,2n-1}H^{-1}_{n-1}
\end{align*}
für alle $n \in \Zefkp$. Hieraus folgt wegen \fref{mpsa2bw1} nun
\begin{align} \label{mpsa2bw3}
  \begin{pmatrix} P^{[0]}_{n} & \ldots & P^{[n]}_{n} \end{pmatrix} = \begin{pmatrix} -z_{n,2n-1}H^{-1}_{n-1} & \Iq \end{pmatrix}
\end{align} 
für alle $n \in \Zefkp$. Wegen \fref{mpsa2bw3} gilt weiterhin
\begin{align*}
  &\ \sum_{l \in {\mathbb Z}_{0,k}, m \in \Zon} P^{[l]}_{k} s_{l+m} \brklam{P^{[m]}_{n}}^{\ast} = Z^{[P_k]}_n \tHn \brklam{Z^{[P_n]}_n}^{\ast} \\
  &= Z^{[P_k]}_n \begin{pmatrix} H_{n-1} & y_{n,2n-1} \\ z_{n,2n-1} & 0_{\qq} \end{pmatrix} \begin{pmatrix} -H^{-1}_{n-1}y_{n,2n-1} \\ \Iq \end{pmatrix} \\
  &= \begin{pmatrix} Z^{[P_k]}_k & 0_{q \times (n-k)q} \end{pmatrix} \begin{pmatrix} 0_{nq\times q} \\ -z_{n,2n-1}H^{-1}_{n-1}y_{n,2n-1} \end{pmatrix}
  = 0_{\qq}
\end{align*}
für alle $n,k \in \Zefkp$ mit $k < n$ und somit ist in Verbindung mit \fref{mpsa2bw4} Teil (iii) von \thref{mpdef2} erfüllt. Teil (ii) von \thref{mpdef2} und die Eindeutigkeit des monischen links-orthogonalen Systems von Matrixpolynomen bezüglich $\sjk$ folgt direkt aus \fref{mpsa2bw3}. \bwend

Man kann das monische links-orthogonale System von Matrixpolynomen bezüglich einer Matrizenfolge, die bis zu einem gewissen Folgenglied Hankel-positiv definit ist, auch mithilfe der kanonischen Hankel-Parametrisierung jener Folge darstellen, wie wir im folgenden Satz zeigen werden. Für eine gegebene Folge aus $\Hqk$ mit $\kappa\in\Na$ wurde die Aussage in \cite[Theorem 5.5(b)]{15} dargestellt, aber nur für den dualen Fall des monischen rechts-orthogonale System von Matrixpolynomen detailliert bewiesen. Wir zeigen analog die Aussage im linken Fall. Hierfür benötigen wir zunächst noch ein Hilfsmittel, das im rechten Fall für eine gegebene Folge aus ${\cal H}^{\geq,e}_{2n-1}$ mit $n\in\N\setminus\gklam{1}$ in \cite[Proposition 3.21]{13} gezeigt wurde. Wir verwenden eine analoge Version des Beweises für den linken Fall.

\begin{lemma}	\thlabel{mplm1}
  Seien $\kappa \in \Na$ und $\sjk$ eine Folge aus $\Cqq$ mit $\sjfkm \in \Hpqfkm$. Weiterhin seien $t_1 := -s_1s^{-1}_0$ und im Fall $\kappa\geq3$
  \begin{align*}
    t_n := \begin{pmatrix} 0_{\qq} & t_{n-1} \end{pmatrix} - (s_{2n-1} - \Lambda_{n-1})\dH^{-1}_{n-1} Z_{n-1} - \dH_{n-1}\dH^{-1}_{n-2} \begin{pmatrix} Z_{n-2} & 0_{\qq} \end{pmatrix}
  \end{align*}
  für alle $n\in\Zzfkp$,
  wobei $Z_0 := \Iq$ und $Z_k := \begin{pmatrix} t_k & \Iq \end{pmatrix}$ für alle $k \in \Zen$. Dann gilt $t_n = -z_{n,2n-1}\Hnm^{-1}$ für alle $n\in\Zefkp$.
\end{lemma}

\bwanf 
Wir zeigen die Behauptung mithilfe vollständiger Induktion. Es gelten
\begin{align*}
	t_1 = -s_1s^{-1}_0 = -z_{1,1}H^{-1}_0
\end{align*}
und im Fall $\kappa\geq3$
\begin{align}	\label{mplm1bw1}
	t_2 = \begin{pmatrix} \Oq & t_1 \end{pmatrix} - (s_3-\Lambda_1)\dH^{-1}_{1} \begin{pmatrix} t_1 & \Iq \end{pmatrix} - \dH_1\dH^{-1}_0 \begin{pmatrix} \Iq & \Oq \end{pmatrix}.
\end{align}
Weiterhin gelten
\begin{align*}
	\begin{pmatrix} \Oq & t_1 \end{pmatrix} H_1 
	& = \begin{pmatrix} \Oq & -s_1s^{-1}_0 \end{pmatrix} \begin{pmatrix} s_0 & s_1 \\ s_1 & s_2 \end{pmatrix} 
	= \begin{pmatrix} -s_1s^{-1}_0s_1 & -s_1s^{-1}_0s_2 \end{pmatrix} \\
	& = \begin{pmatrix} \dH_1-s_2 & -s_1s^{-1}_0s_2 \end{pmatrix}, \\
	\begin{pmatrix} t_1 & \Iq \end{pmatrix} H_1
	& = \begin{pmatrix} -s_1s^{-1}_0 & \Iq \end{pmatrix} \begin{pmatrix} s_0 & s_1 \\ s_1 & s_2 \end{pmatrix} 
	= \begin{pmatrix} -s_1+s_1 & -s_1s^{-1}_0s_1+s_2 \end{pmatrix}
	= \begin{pmatrix} \Oq & \dH_1 \end{pmatrix}, \\	
	\begin{pmatrix} \Iq & \Oq \end{pmatrix} H_1
	& = \begin{pmatrix} \Iq & \Oq \end{pmatrix} \begin{pmatrix} s_0 & s_1 \\ s_1 & s_2 \end{pmatrix}
	= \begin{pmatrix} s_0 & s_1 \end{pmatrix}
	= \begin{pmatrix} \dH_0 & s_1 \end{pmatrix}
\end{align*}
und
\begin{align*}
	& -s_1s^{-1}_0s_2-(s_3-\Lambda_1)-\dH_1\dH^{-1}_0s_1 \\
	&= -s_1s^{-1}_0s_2-(s_3-s_1s^{-1}_0s_2-s_2s^{-1}_0s_1+s_1s^{-1}_0s_1s^{-1}_0s_1)-(s_2-s_1s^{-1}_0s_1)s^{-1}_0s_1 \\
	&= -s_3.
\end{align*}
Hieraus folgt wegen \fref{mplm1bw1} dann
\begin{align*}
	t_2H_1 &= \begin{pmatrix} \dH_1-s_2 & -s_1s^{-1}_0s_2 \end{pmatrix} - \begin{pmatrix} \Oq & s_3-\Lambda_1 \end{pmatrix} - \begin{pmatrix} \dH_1 & \dH_1\dH^{-1}_0s_1 \end{pmatrix} = - \begin{pmatrix} s_2 & s_3 \end{pmatrix} \\
	&= -z_{2,3},
\end{align*}
also $t_2 = -z_{2,3}H^{-1}_1$.
Seien im Fall $\kappa\ge5$ nun $n\in\Z{3}{\fklam{\kappa+1}}$ und $t_k = -z_{k,2k-1}H^{-1}_{k-1}$ für alle $k \in \Zenm$ erfüllt. Dann gelten
\begin{align}	\label{mplm1bw2}
	Z_{n-1}\Hnm &= \begin{pmatrix} t_{n-1} & \Iq \end{pmatrix} \begin{pmatrix} H_{n-2} & y_{n-1,2n-3} \\ z_{n-1,2n-3} & s_{2n-2} \end{pmatrix} \notag \\
	&= \begin{pmatrix} t_{n-1}H_{n-2}+z_{n-1,2n-3} & t_{n-1}y_{n-1,2n-3}+s_{2n-2} \end{pmatrix}
	= \begin{pmatrix} 0_{q\times (n-1)q} & \dH_{n-1} \end{pmatrix}
\end{align}
und
\begin{align*}
	\Tnma\Hnm-\Hnm\Tnm &= \begin{pmatrix} K_{n-2} & y_{n,2n-2} \\ 0_{(n-1)q \times q} & \Oq \end{pmatrix}-\begin{pmatrix} K_{n-2} & 0_{q \times (n-1)q} \\ z_{n,2n-2} & \Oq \end{pmatrix} \\
	&= \begin{pmatrix} 0_{(n-1)q \times (n-1)q} & y_{n,2n-2} \\ -z_{n,2n-2} & \Oq \end{pmatrix}.
\end{align*}
Hieraus folgt
\begin{align}	\label{mplm1bw3}
	& \begin{pmatrix} \Oq & t_{n-1} \end{pmatrix} \Hnm = Z_{n-1}\Tnma\Hnm \notag \\
	&= Z_{n-1}\Hnm\Tnm + \begin{pmatrix} t_{n-1} & \Iq \end{pmatrix}\begin{pmatrix} 0_{(n-1)q \times (n-1)q} & y_{n,2n-2} \\ -z_{n,2n-2} & \Oq \end{pmatrix} \notag \\
	&= \begin{pmatrix} 0_{q\times (n-1)q} & \dH_{n-1} \end{pmatrix}\Tnm + \begin{pmatrix} -z_{n,2n-2} & t_{n-1}y_{n,2n-2} \end{pmatrix} \notag \\
	&= \begin{pmatrix} 0_{q\times (n-2)q} & \dH_{n-1} & \Oq \end{pmatrix} - z_{n,2n-1} + \begin{pmatrix} 0_{q\times(n-1)q} & s_{2n-1}-z_{n-1,2n-3}H^{-1}_{n-2}y_{n,2n-2} \end{pmatrix} \notag\\
	&= \begin{pmatrix} 0_{q\times (n-2)q} & \dH_{n-1} & s_{2n-1}-z_{n-1,2n-3}H^{-1}_{n-2}y_{n,2n-2} \end{pmatrix} - z_{n,2n-1}.
\end{align}
Weiterhin gelten
\begin{align*}
	Z_{n-2}H_{n-2} &= \begin{pmatrix} t_{n-2} & \Iq \end{pmatrix} \begin{pmatrix} H_{n-3} & y_{n-2,2n-5} \\ z_{n-2,2n-5} & s_{2n-4} \end{pmatrix} \\
	&= \begin{pmatrix} t_{n-2}H_{n-3}+z_{n-2,2n-5} & t_{n-2}y_{n-2,2n-5}+s_{2n-4} \end{pmatrix}
	= \begin{pmatrix} 0_{q\times (n-2)q} & \dH_{n-2} \end{pmatrix}
\end{align*}
und
\begin{align*}
	Z_{n-2}y_{n-1,2n-3} &= \begin{pmatrix} t_{n-2} & \Iq \end{pmatrix} \begin{pmatrix} y_{n-1,2n-4} & s_{2n-3} \end{pmatrix}
	= -z_{n-2,2n-5}H^{-1}_{n-3}y_{n-1,2n-4}+s_{2n-3}.
\end{align*}
Hieraus folgt
\begin{align}	\label{mplm1bw4}
	\begin{pmatrix} Z_{n-2} & \Oq \end{pmatrix} H_{n-1} 
	&= \begin{pmatrix} Z_{n-2} & \Oq \end{pmatrix} \begin{pmatrix} H_{n-2} & y_{n-1,2n-3} \\ z_{n-1,2n-3} & s_{2n-2} \end{pmatrix} \notag \\
	&= \begin{pmatrix} Z_{n-2}H_{n-2} & Z_{n-2}y_{n-1,2n-3} \end{pmatrix} \notag \\
	&= \begin{pmatrix} 0_{q\times(n-2)q} & \dH_{n-2} & s_{2n-3}-z_{n-2,2n-5}H^{-1}_{n-3}y_{n-1,2n-4} \end{pmatrix}.
\end{align}
Wegen $(s_j)^{2n-2}_{j=0}\in{\cal H}^{>}_{q,2n-2}$ (siehe \thref{asmbm1}) und ${\cal H}^{>}_{q,2n-2} \subseteq {\cal H}^{\geq,e}_{q,2n-2}$ (siehe \thref{asmbm6}) gilt $(s_j)^{2n-2}_{j=0}\in{\cal H}^{\geq,e}_{q,2n-2}$. Hieraus folgt wegen \cite[Lemma 3.14]{13} dann
\begin{align*}
	\Lambda_{n-1} = z_{n-1,2n-3}H^{-1}_{n-2}y_{n,2n-2}+\dH_{n-1}\dH^{-1}_{n-2}(s_{2n-3}-z_{n-2,2n-5}H^{-1}_{n-3}y_{n-1,2n-4})
\end{align*}
(diese Aussage kann auch unter Beachtung von \thref{asmbm1} mithilfe von \cite[Lemma 5.5]{Pert} gewonnen werden).
Hieraus folgt wegen \fref{mplm1bw2}-\fref{mplm1bw4} nun
\begin{align*}
	t_n\Hnm	&= \begin{pmatrix} 0_{\qq} & t_{n-1} \end{pmatrix}\Hnm - (s_{2n-1} - \Lambda_{n-1})\dH^{-1}_{n-1} Z_{n-1}\Hnm \\
	&\quad - \dH_{n-1}\dH^{-1}_{n-2} \begin{pmatrix} Z_{n-2} & 0_{\qq} \end{pmatrix}\Hnm \\
	&= \begin{pmatrix} 0_{q\times (n-2)q} & \dH_{n-1} & s_{2n-1}-z_{n-1,2n-3}H^{-1}_{n-2}y_{n,2n-2} \end{pmatrix} \\
	&\quad - z_{n,2n-1} - \begin{pmatrix} 0_{q\times (n-1)q} & s_{2n-1}-\Lambda_{n-1} \end{pmatrix} \\
	&\quad - \begin{pmatrix} 0_{q\times(n-2)q} & \dH_{n-1} & \dH_{n-1}\dH^{-1}_{n-2}(s_{2n-3}-z_{n-2,2n-5}H^{-1}_{n-3}y_{n-1,2n-4}) \end{pmatrix} \\
	&= - z_{n,2n-1},
\end{align*}
also $t_n = -z_{n,2n-1}\Hnm^{-1}$. \bwend

\begin{satz}	\thlabel{mpsa3}
  Seien $\kappa \in \Noa$ und $\sjk$ eine Folge aus $\Cqq$ mit $\sjfkm \in \Hpqfkm$, falls $\kappa \geq 1$.
  Weiterhin sei $\Pnfk$ eine Folge von \textit{q}$\times$\textit{q}-Matrixpolynomen und im Fall $\kappa\geq1$ $\CDfk$ die kanonische Hankel-Parametrisierung von $\sjk$. \dsfaa
  \begin{itemize}
    \item [\rm{(i)}] Es ist $\Pnfk$ das monische links-orthogonale System von Matrixpolynomen bezüglich $\sjk$.
    \item [\rm{(ii)}] Es gelten $P_{0} \equiv \Iq$, im Fall $\kappa \geq 1$
      \begin{align*}
	P_{1}(z) = (z\Iq-C_{1}D^{-1}_{0})P_{0}(z)
      \end{align*}
      für alle $z \in \C$ und im Fall $\kappa \geq 3$
      \begin{align*}
	P_{n}(z) = (z\Iq-C_{n}D^{-1}_{n-1})P_{n-1}(z)-D_{n-1}D^{-1}_{n-2}P_{n-2}(z)
      \end{align*}
      für alle $z \in \C$ und $n \in \Zzfkp$.
  \end{itemize}
\end{satz}


\bwanf Im Fall $\kappa=0$ folgt dies aus Teil (ii) von \thref{mpdef2}.

Im Fall $\kappa \in \gklam{1,2}$ gilt unter Beachtung von \thref{mpdef1} dann
\begin{align}	\label{mpsa3bw1}
  P_1 = (z\Iq-C_{1}D^{-1}_{0})P_{0}(z) = (z\Iq-s_1s^{-1}_0)\Iq = z\Iq-z_{1,1}H^{-1}_0
\end{align}
für alle $z \in \C$. Wegen \thref{mpsa2} ist dann $\Pnfk$ genau das eine monische links-orthogonale System von Matrixpolynomen bezüglich $\sjk$.

Sei nun $\kappa \geq 3$. Weiterhin seien
\begin{align*}
  Z^{[P_n]}_n := \begin{pmatrix} P^{[0]}_n & \ldots & P^{[n]}_n \end{pmatrix}
\end{align*}
für alle $n \in \Zofkp$ und
\begin{align*}
  Z^{[P_k]}_n := \begin{pmatrix} Z^{[P_k]}_k & 0_{q \times (n-k)q} \end{pmatrix}
\end{align*}
für alle $n,k \in \Zofkp$ mit $k<n$. Anhand der Rekursionsvorschrift erkennt man leicht, dass für alle $n \in \Zefkp$ das Matrixpolynom $P_n$ vom Grad $n$ mit Leitkoeffizient $\Iq$ ist.
Es existiert also eine Matrix $t_n \in \C^{q \times nq}$ mit
\begin{align}	\label{mpsa3bw2}
  Z^{[P_n]}_n = \begin{pmatrix} t_n & \Iq \end{pmatrix}
\end{align}
für alle $n \in \Zefkp$. Wegen \fref{mpsa3bw1} gilt $t_1 = -s_1s^{-1}_0$. Unter Beachtung von \thref{mpdef1} gilt weiterhin
\begin{align*}
  Z^{[P_n]}_n & = \begin{pmatrix} 0_{\qq} & Z^{[P_{n-1}]}_{n-1} \end{pmatrix}-C_{n}D^{-1}_{n-1}Z^{[P_{n-1}]}_n-D_{n-1}D^{-1}_{n-2}Z^{[P_{n-2}]}_n \\
  & = \begin{pmatrix} 0_{\qq} & t_{n-1} & \Iq \end{pmatrix}-(s_{2n-1}-\Lambda_{n-1})\widehat{H}^{-1}_{n-1}\begin{pmatrix} t_{n-1} & \Iq & 0_{\qq} \end{pmatrix}
  -\widehat{H}_{n-1}\widehat{H}^{-1}_{n-2}Z^{[P_{n-2}]}_n
\end{align*}
für alle $n \in \Zzfkp$, wobei
\begin{align*}
  Z^{[P_{n-2}]}_n = \begin{cases} \begin{pmatrix} \Iq & 0_{q \times 2q} \end{pmatrix} & \text{falls } n=2 \\ \begin{pmatrix} t_{n-2} & \Iq & 0_{q \times 2q} \end{pmatrix} & \text{falls } n > 2. \end{cases}
\end{align*}
Hieraus folgt wegen \fref{mpsa3bw2} dann 
\begin{align*}
  t_n = \begin{pmatrix} 0_{\qq} & t_{n-1} \end{pmatrix}-(s_{2n-1}-\Lambda_{n-1})\widehat{H}^{-1}_{n-1}\begin{pmatrix} t_{n-1} & \Iq \end{pmatrix}
  -\widehat{H}_{n-1}\widehat{H}^{-1}_{n-2}Z^{[P_{n-2}]}_{n-1}
\end{align*}
für alle $n \in \Zzfkp$, wobei
\begin{align*}
  Z^{[P_{n-2}]}_{n-1} = \begin{cases} \begin{pmatrix} \Iq & 0_{q \times q} \end{pmatrix} & \text{falls } n=2 \\ \begin{pmatrix} t_{n-2} & \Iq & 0_{q \times q} \end{pmatrix} & \text{falls } n > 2. \end{cases}
\end{align*}
Hieraus folgt wegen \thref{mplm1} dann $t_n=-z_{n,2n-1}\Hnm^{-1}$ für alle $n \in \Zefkp$, das heißt
\begin{align*}
  \begin{pmatrix} P^{[0]}_{n} & \ldots & P^{[n]}_{n} \end{pmatrix} = \begin{cases} \Iq & \text{falls } n=0 \\
  \begin{pmatrix} -z_{n,2n-1}H^{-1}_{n-1} & \Iq \end{pmatrix} & \text{falls } n \in \Zefkp. \end{cases}
\end{align*}
Somit stimmt wegen \thref{mpsa2} die Folge $\Pnfk$ mit dem monischen links-orthogonalen System von Matrixpolynomen bezüglich $\sjk$ überein. \bwend

Mithilfe von \thref{fpsa1} und  \thref{mpsa3} können wir nun das monische links-orthogonale System von Matrixpolynomen bezüglich einer Matrizenfolge, die bis zu einem gewissen Folgenglied Hankel-positiv definit ist, unter Verwendung des Favard-Paar bezüglich jener Folge darstellen (vergleiche \cite[Proposition 8.11(a)]{CR1} im Fall $\kappa=\infty$).

\begin{folg}	\thlabel{mpfo2}
	Seien $\kappa \in \Noa$ und $\sjk$ eine Folge aus $\Cqq$ mit \linebreak $\sjfkm \in \Hpqfkm$, falls $\kappa\geq1$.
  	Weiterhin seien $\Pnfk$ eine Folge von \textit{q}$\times$\textit{q}-Matrixpolynomen und im Fall $\kappa\geq1$ $\ABfk$ das Favard-Paar bezüglich $\sjk$. \dsfaa
  	\begin{itemize}
    \item [\rm{(i)}] Es ist $\Pnfk$ das monische links-orthogonale System von Matrixpolynomen bezüglich $\sjk$.
    \item [\rm{(ii)}] Es gelten $P_{0} \equiv \Iq$, im Fall $\kappa\geq1$
    \begin{align*}
		P_{1}(z) = (z\Iq-A_0)P_{0}(z)
    \end{align*}
    für alle $z \in \C$ und im Fall $\kappa \geq 3$
    \begin{align*}
		P_{n}(z) = (z\Iq-A_{n-1})P_{n-1}(z)-B^{\ast}_{n-1}P_{n-2}(z)
    \end{align*}
    für alle $z \in \C$ und $n \in \Zzfkp$.
  	\end{itemize}
\end{folg}

\bwanf Der Fall $\kappa=0$ folgt unmittelbar aus \thref{mpsa3}. Seien nun $\kappa\geq1$ und $\CDfk$ die kanonische Hankel-Parametrisierung von $\sjk$. Unter Beachtung von \thref{mpbm3} gilt wegen $\sjfkm \in \Hpqfkm$ und \thref{mpsa1} dann $D_n\in\Cpqq$ für alle $n\in\Zofkm$. Hieraus folgt wegen Teil (a) von \thref{fpsa1} nun
\begin{align*}
	A_n = C_{n+1}D^{-1}_n
\end{align*}
für alle $n\in\Zofkm$ und im Fall $\kappa\geq3$
\begin{align*}
	B^{\ast}_n = \rklam{D^{-1}_{n-1}D_n}^{\ast} = D^{\ast}_nD^{-\ast}_{n-1} = D_nD^{-1}_{n-1}
\end{align*}
für alle $n\in\Zefkm$. Hieraus folgt wegen \thref{mpsa3} dann die Behauptung. \bwend

Es sei bemerkt, das unter Beachtung von Teil (b) von \thref{asmdef2} bzw. \thref{asmdef3} eine rechtsseitig bzw. linksseitig $\alpha$-Stieltjes-positiv definite Folge auch bis zu einem gewissen Folgenglied Hankel-positiv definit ist und weiterhin unter Beachtung von \thref{aspdef1} die durch rechtsseitige bzw. linksseitige $\alpha$-Verschiebung generierte Folge ebenso bis zu einem gewissen Folgenglied Hankel-positiv definit ist. Somit können wir die zuvor gewonnen Resultate auf $\alpha$-Stieltjes-positiv definite Folgen erweitern. Dies werden wir später in Kapitel \ref{chapsq} noch ausführlicher behandeln und zunächst nur folgendes Hilfsresultat bereitstellen, welches für gerade $\kappa$ oder $\kappa=\infty$ bereits in \cite[Proposition 7.2]{Pert} gezeigt wurde. Zur besseren Anschauung geben wir den dortigen Beweis noch einmal wieder.

\begin{lemma}	\thlabel{mplm2}
	Seien $\kappa\in\Na$, $\alpha\in\R$ und $\sjk\in\Kpqka$. Weiterhin seien $\Pnfk$ bzw. $\Parnfk$ das monische links-orthogonale System von Matrixpolynomen bezüglich $\sjk$ bzw. $\sarjk$. Dann gilt
	\begin{align*}
		(z-\alpha)P_{\arn}(z) = P_{n+1}(z)+\dHarn\dHn^{-1}P_n(z)
	\end{align*}
	für alle $z\in\C$ und $n\in\Zofkm$.
\end{lemma}

\bwanf 
Wegen \thref{mpsa2} gilt
\begin{align*}
	&\ P_1(z)+\dH_{\aro}\dH^{-1}_0P_0(z) = z\Iq-z_{1,1}H^{-1}_0+s_{\aro}s^{-1}_0\Iq \\ 
	&= z\Iq-s_1s^{-1}_0+(-\alpha s_0+s_1)s^{-1}_0 = (z-\alpha)\Iq = (z-\alpha)P_{\aro}(z)
\end{align*}
für alle $z\in\C$. Seien nun $\kappa\geq3$ und $n\in\Zefkm$. Weiterhin sei
\begin{align*}
	Q_n(z) := (z-\alpha)P_{\arn}(z)-\dHarn\dHn^{-1}P_n(z)
\end{align*}
für alle $z\in\C$. Wegen \thref{mpdef2} gilt dann $\deg Q_n \leq n+1$. Sei
\begin{align*}
	X_n := \begin{pmatrix} \Oq & -z_{\arn,2n-1}\Harnm^{-1} \end{pmatrix} - \alpha \begin{pmatrix} -z_{\arn,2n-1}\Harnm^{-1} & \Iq \end{pmatrix}.
\end{align*}
Es gilt
\begin{align*}
	\eklam{\begin{pmatrix} \Inpq \\ 0_{q \times (n+1)q} \end{pmatrix} - \alpha \begin{pmatrix} 0_{q \times (n+1)q} \\ \Inpq \end{pmatrix}}\Hn 
	&= \begin{pmatrix} \Hn \\ 0_{q \times (n+1)q} \end{pmatrix} - \alpha \begin{pmatrix} 0_{q \times (n+1)q} \\ \Hn \end{pmatrix} \\
	&= \begin{pmatrix} z_{0,n-1} & s_n \\ \Knm & y_{n+1,2n} \\ \Oqn & \Oq \end{pmatrix} - \alpha \begin{pmatrix} \Oqn & \Oq \\ \Hnm & y_{n,2n-1} \\ z_{n,2n-1} & s_{2n} \end{pmatrix} \\
	&= \begin{pmatrix} z_{0,n-1} & s_n \\ \Harnm & y_{\arn,2n-1} \\ -\alpha z_{n,2n-1} & -\alpha s_{2n} \end{pmatrix}.
\end{align*}
Hieraus folgt
\begin{align} \label{mplm2bw1}
	X_n\Hn &= \begin{pmatrix} \Oq & -z_{\arn,2n-1}\Harnm^{-1} & \Iq \end{pmatrix}\eklam{\begin{pmatrix} \Inpq \\ 0_{q \times (n+1)q} \end{pmatrix} - \alpha \begin{pmatrix} 0_{q \times (n+1)q} \\ \Inpq \end{pmatrix}}\Hn \notag \\
	&= \begin{pmatrix} \Oq & -z_{\arn,2n-1}\Harnm^{-1} & \Iq \end{pmatrix} \begin{pmatrix} z_{0,n-1} & s_n \\ \Harnm & y_{\arn,2n-1} \\ -\alpha z_{n,2n-1} & -\alpha s_{2n} \end{pmatrix} \notag \\
	&= \begin{pmatrix} -z_{\arn,2n-1}-\alpha z_{n,2n-1} & -z_{\arn,2n-1}\Harnm^{-1}y_{\arn,2n-1}-\alpha s_{2n} \end{pmatrix} \notag \\
	&= \begin{pmatrix} -z_{n+1,2n} & \dHarn-s_{2n+1} \end{pmatrix} \notag \\
	&= \begin{pmatrix} \Oqn & \dHarn \end{pmatrix} - z_{n+1,2n+1}.
\end{align}
Weiterhin sei
\begin{align*}
	Y_n := \dHarn\dHn^{-1} \begin{pmatrix} -z_{n,2n-1}\Hn^{-1} & \Iq \end{pmatrix}.
\end{align*}
Dann gilt
\begin{align*}
	Y_n\Hn &= \dHarn\dHn^{-1} \begin{pmatrix} -z_{n,2n-1}\Hn^{-1} & \Iq \end{pmatrix} \begin{pmatrix} \Hnm & y_{n,2n-1} \\ z_{n,2n-1} & s_{2n} \end{pmatrix} \\
	&= \dHarn\dHn^{-1} \begin{pmatrix} \Oqn & \dHn \end{pmatrix} \\
	&= \begin{pmatrix} \Oqn & \dHarn \end{pmatrix}.
\end{align*}
Hieraus folgt wegen \fref{mplm2bw1} dann
\begin{align*}
	(X_n-Y_n)\Hn = -z_{n+1,2n-1},
\end{align*}
also 
\begin{align*}
	X_n-Y_n = -z_{n+1,2n-1}\Hn^{-1}.
\end{align*}
Hieraus folgt wegen \thref{mpsa2} nun
\begin{align*}
	\begin{pmatrix} Q^{[0]}_n & \ldots & Q^{[n+1]}_n \end{pmatrix} 
	&= \begin{pmatrix} \Oq & P^{[0]}_{\arn} & \ldots & P^{[n]}_{\arn} \end{pmatrix} - \alpha \begin{pmatrix} P^{[0]}_{\arn} & \ldots & P^{[n]}_{\arn} & \Oq \end{pmatrix} \\
	&\quad -\dHarn\dHn^{-1}\begin{pmatrix} P^{[0]}_n & \ldots & P^{[n]}_n & \Oq \end{pmatrix} \\
	&= \begin{pmatrix} \Oq & -z_{\arn,2n-1}\Harnm^{-1} & \Iq \end{pmatrix} - \alpha \begin{pmatrix} -z_{\arn,2n-1}\Harnm^{-1} & \Iq & \Oq \end{pmatrix} \\
	&\quad -\dHarn\dHn^{-1} \begin{pmatrix} -z_{n,2n-1}\Hn^{-1} & \Iq & \Oq \end{pmatrix} \\
	&= \begin{pmatrix} X_n & \Iq \end{pmatrix} - \begin{pmatrix} Y_n & \Oq \end{pmatrix} \\
	&= \begin{pmatrix} -z_{n+1,2n-1}\Hn^{-1} & \Iq \end{pmatrix} \\
	&= \begin{pmatrix} P^{[0]}_{n+1} & \ldots & P^{[n+1]}_{n+1} \end{pmatrix}.
\end{align*}
Somit gilt also $Q_n(z) = P_{n+1}(z)$ für alle $z\in\C$. \bwend

Wir kommen nun zum zweiten zentralen Begriff dieses Abschnitts, dem linken System von Matrixpolynomen zweiter Art bezüglich einer Matrizenfolge, die bis zu einem gewissen Folgenglied Hankel-positiv definit ist. Hierfür benötigen wir folgende Bezeichnungen. 

\begin{bez}	\thlabel{mpbz3}
	Seien $\kappa \in \Noa$ und $\sjk$ aus $\Cqr$. Für jedes $n \in \Zok$ heißt dann
	\begin{align*}
  		S_{n} := \begin{cases} s_0 & \text{falls } n=0 \\
  		\begin{pmatrix} s_{0} & 0_{\qr} & \ldots & 0_{\qr} \\
  		s_{1} & s_{0} & \ddots & \vdots \\
  		\vdots & \ddots & \ddots & 0_{\qr} \\
  		s_{n} & \ldots & s_{1} & s_{0} \end{pmatrix} & \text{falls } n>0 \end{cases}
	\end{align*}
	\textbf{untere n-te Blockdreiecksmatrix} von $\sjk$. Sei $P$ ein \textit{p}$\times$\textit{q}-Matrixpolynom mit $k:=\deg P \leq \kappa+1$. Dann heißt $\Ps:\C \rightarrow \Cpr$ definiert gemäß
	\begin{align*}
  		\Ps(z) := \begin{cases} 0_{\pr} & \text{falls } k \leq 0 \\
  		\begin{pmatrix} P^{[0]} & \ldots & P^{[k]} \end{pmatrix} \begin{pmatrix} 0_{q \times kr} \\ S_{k-1} \end{pmatrix} \begin{pmatrix} \Ir \\ z\Ir \\ \vdots \\ z^{k-1}\Ir \end{pmatrix} & \text{falls } k \geq 1\end{cases}
	\end{align*}
	das \textbf{zu} $\sjk$ \textbf{gehörige Matrixpolynom} bezüglich $P$.
\end{bez}

\begin{defi}	\thlabel{mpdef3}
  Seien $\kappa \in \Noa$ und $\sjk$ eine Folge aus $\Cqq$ mit \linebreak $\sjfkm \in \Hpqfkm$, falls $\kappa \geq 1$. Weiterhin sei $\Pnfk$ das monische links"=orthogonale System von Matrixpolynomen bezüglich $\sjk$.
  Dann heißt $\Psnfk$ \textbf{linkes System von Matrixpolynomen zweiter Art} bezüglich $\sjk$.
\end{defi}

Folgendes Resultat liefert uns eine genaue Beschreibung des linken System von Matrixpolynomen zweiter Art bezüglich einer Matrizenfolge, die bis zu einem gewissen Folgenglied Hankel-positiv definit ist (vergleiche \cite[Remark 4.7]{Pert} für eine gegebene Folge aus $\Hpqk$ mit $\kappa\in\Na$).

\begin{bem}	\thlabel{mpbem1}
  Seien $\kappa \in \Noa$ und $\sjk$ eine Folge aus $\Cqq$ mit \linebreak $\sjfkm \in \Hpqfkm$, falls $\kappa \geq 1$. 
  Weiterhin sei $\Psnfk$ das linke System von Matrixpolynomen zweiter Art bezüglich $\sjk$.
  Dann gelten $\Ps_{0} \equiv 0_{\qq}$ und für alle $n\in \Zefkp$ weiterhin $\deg \Ps_{n} = n-1$ sowie
  \begin{align*}
    \begin{pmatrix} (\Ps_{n})^{\eklam{0}} & \ldots & (\Ps_{n})^{\eklam{n-1}} \end{pmatrix}
    = \begin{pmatrix} -z_{n,2n-1}H^{-1}_{n-1} & \Iq \end{pmatrix} \begin{pmatrix} 0_{q\times nq} \\ S_{n-1}\end{pmatrix}.
  \end{align*}
\end{bem}

\bwanf Dies ergibt sich sofort aus \thref{mpsa2} und \thref{mpdef3}. \bwend

Auch das linke System von Matrixpolynomen zweiter Art bezüglich einer Matrizenfolge, die bis zu einem gewissen Folgenglied Hankel-positiv definit ist, lässt sich wie schon für das links-orthogonale System von Matrixpolynomen bezüglich jener Folge mithilfe der kanonischen Hankel-Parametrisierung jener Folge darstellen, wie wir im folgenden Satz zeigen werden (vergleiche \cite[Proposition 4.9]{Pert} für eine gegebene Folge aus $\Hpqk$ mit $\kappa\in\Na$).

\begin{satz}	\thlabel{mpsa5}
  Seien $\kappa \in \Noa$ und $\sjk$ eine Folge aus $\Cqq$ mit $\sjfkm \in \Hpqfkm$, falls $\kappa \geq 1$.
  Weiterhin seien $(Q_{n})^{\fklam{\kappa+1}}_{n=0}$ eine Folge von \textit{q}$\times$\textit{q}-Matrixpolynomen und im Fall $\kappa\geq1$ $\CDfk$ die kanonische Hankel-Parametrisierung von $\sjk$. \dsfaa
  \begin{itemize}
    \item [\rm{(i)}] Es ist $(Q_{n})^{\fklam{\kappa+1}}_{n=0}$ das linke System von Matrixpolynomen zweiter Art bezüglich $\sjk$.
    \item [\rm{(ii)}] Es gelten $Q_{0} \equiv 0_{\qq}$, im Fall $\kappa \geq 1$ $\;Q_{1} \equiv D_{0}$ und im Fall $\kappa \geq 3$
      \begin{align*}
	Q_{n}(z) = (z\Iq-C_{n}D^{-1}_{n-1})Q_{n-1}(z) - D_{n-1}D^{-1}_{n-2}Q_{n-2}(z)
      \end{align*}
      für alle $z \in \C$, $n \in \Zzfkp$.  
  \end{itemize}
\end{satz}


\bwanf Sei $\Pnfk$ das monische links-orthogonale System von Matrixpolynomen bezüglich $\sjk$. 
Wegen \thref{mpsa3} gelten dann $P_{0} \equiv \Iq$ und für $\kappa \geq 1$
\begin{align}	\label{mpsa5bw1}
  P_{1}(z) = (z\Iq-C_{1}D^{-1}_{0})P_{0}(z)
\end{align}
für alle $z \in \C$ und für $\kappa \geq 3$
\begin{align}	\label{mpsa5bw2}
  P_{n}(z) = (z\Iq-C_{n}D^{-1}_{n-1})P_{n-1}(z)-D_{n-1}D^{-1}_{n-2}P_{n-2}(z)
\end{align}
für alle $z \in \C$ und $n \in \Zzfkp$.
Sei $\Psnfk$ das linke System von Matrixpolynomen zweiter Art bezüglich $\sjk$. 
Dann gilt $\Ps_{0} \equiv 0_{\qq}$ und wegen \fref{mpsa5bw1} und \thref{mpdef1} für $\kappa \geq 1$
\begin{align*}
  \Ps_1 = \begin{pmatrix} P^{[0]}_1 & P^{[1]}_1 \end{pmatrix} \begin{pmatrix} 0_{\qq} \\ S_0 \end{pmatrix} \Iq
  = \begin{pmatrix} -s_1s^{-1}_0 & \Iq \end{pmatrix} \begin{pmatrix} 0_{\qq} \\ s_0 \end{pmatrix}
  = s_0 = D_0.
\end{align*}
Sei nun $\kappa \geq 3$. Dann gilt
\begin{align}	\label{mpsa5bw3}
  \begin{pmatrix} (\Ps_n)^{[0]} & \ldots & (\Ps_n)^{[n-1]} \end{pmatrix} = \begin{pmatrix} P^{[0]}_n & \ldots & P^{[n]}_n \end{pmatrix} \begin{pmatrix} 0_{q \times nq} \\ S_{n-1} \end{pmatrix}
\end{align}
für alle $n \in \Zzfkp$. Wegen \fref{mpsa5bw2} gilt
\begin{align}	\label{mpsa5bw4}
  \begin{pmatrix} P^{[0]}_n & \ldots & P^{[n]}_n \end{pmatrix} & = \begin{pmatrix} 0_{\qq} & P^{[0]}_{n-1} & \ldots & P^{[n-1]}_{n-1} \end{pmatrix} \notag \\
  &\quad - C_nD^{-1}_{n-1} \begin{pmatrix} P^{[0]}_{n-1} & \ldots & P^{[n-1]}_{n-1} & 0_{\qq} \end{pmatrix} \notag \\ 
  &\quad - D_{n-1}D^{-1}_{n-2} \begin{pmatrix} P^{[0]}_{n-2} & \ldots & P^{[n-2]}_{n-2} & 0_{q \times 2q} \end{pmatrix}
\end{align}
für alle $n \in \Zzfkp$. Wegen Teil (iii) von \thref{mpdef2} gilt
\begin{align}	\label{mpsa5bw5}
  & \begin{pmatrix} 0_{\qq} & P^{[0]}_{n-1} & \ldots & P^{[n-1]}_{n-1} \end{pmatrix} \begin{pmatrix} 0_{q \times nq} \\ S_{n-1} \end{pmatrix} \notag \\
  & = \begin{pmatrix} 0_{\qq} & P^{[0]}_{n-1} & \ldots & P^{[n-1]}_{n-1} \end{pmatrix} \begin{pmatrix} 0_{\qq} & 0_{q \times (n-1)q} \\ s_0 & 0_{q \times (n-1)q} \\ y_{1,n-1} & S_{n-2} \end{pmatrix} \notag \\
  & = \begin{pmatrix} 0_{\qq} & (\Ps_{n-1})^{[0]} & \ldots & (\Ps_{n-1})^{[n-1]} \end{pmatrix}
\end{align}
für alle $n \in \Zzfkp$. Weiterhin gelten
\begin{align}	\label{mpsa5bw6}
  & \begin{pmatrix} P^{[0]}_{n-1} & \ldots & P^{[n-1]}_{n-1} & 0_{\qq} \end{pmatrix} \begin{pmatrix} 0_{q \times nq} \\ S_{n-1} \end{pmatrix} \notag \\
  & = \begin{pmatrix} P^{[0]}_{n-1} & \ldots & P^{[n-1]}_{n-1} & 0_{\qq} \end{pmatrix} \begin{pmatrix} 0_{q \times (n-1)q} & 0_{\qq} \\ S_{n-2} & 0_{q \times (n-1)q} \\ s_{n-1} \, \ldots \, s_1 & s_0 \end{pmatrix} \notag \\
  & = \begin{pmatrix} (\Ps_{n-1})^{[0]} & \ldots & (\Ps_{n-1})^{[n-1]} & 0_{\qq} \end{pmatrix}
\end{align}
für alle $n \in \Zzfkp$,
\begin{align}	\label{mpsa5bw7}
  \begin{pmatrix} P^{[0]}_{0} & 0_{q \times 2q} \end{pmatrix} \begin{pmatrix} 0_{q \times 2q} \\ S_1 \end{pmatrix} 
  = \begin{pmatrix} (\Ps_0)^{[0]} & 0_{q \times 2q} \end{pmatrix}
\end{align}
und für $\kappa \geq 5$
\begin{align}	\label{mpsa5bw8}
  & \begin{pmatrix} P^{[0]}_{n-2} & \ldots & P^{[n-2]}_{n-2} & 0_{q \times 2q} \end{pmatrix} \begin{pmatrix} 0_{q \times nq} \\ S_{n-1} \end{pmatrix} \notag \\
  & = \begin{pmatrix} P^{[0]}_{n-2} & \ldots & P^{[n-2]}_{n-2} & 0_{q \times 2q} \end{pmatrix} \begin{pmatrix} 0_{q \times (n-2)q} & 0_{q \times q} & 0_{q \times q} \\ S_{n-3} & 0_{(n-2)q \times q} & 0_{(n-2)q \times q} \\ s_{n-2} \, \ldots \, s_1 & s_0 & 0_{\qq} \\ s_{n-1} \, \ldots \, s_2 & s_1 & s_0 \end{pmatrix} \notag \\
  & = \begin{pmatrix} (\Ps_{n-2})^{[0]} & \ldots & (\Ps_{n-2})^{[n-2]} & 0_{q \times 2q} \end{pmatrix}
\end{align}
für alle $n \in \Z{3}{\fklam{\kappa+1}}$. Wegen \fref{mpsa5bw3}, \fref{mpsa5bw4} und \fref{mpsa5bw5} - \fref{mpsa5bw8} gilt dann
\begin{align*}
  \begin{pmatrix} (\Ps_n)^{[0]} & \ldots & (\Ps_n)^{[n-1]} \end{pmatrix} & = \begin{pmatrix} P^{[0]}_n & \ldots & P^{[n]}_n \end{pmatrix} \begin{pmatrix} 0_{q \times nq} \\ S_{n-1} \end{pmatrix} \\
  & = \begin{pmatrix} 0_{\qq} & P^{[0]}_{n-1} & \ldots & P^{[n-1]}_{n-1} \end{pmatrix} \begin{pmatrix} 0_{q \times nq} \\ S_{n-1} \end{pmatrix} \\
  &\quad - C_nD^{-1}_{n-1} \begin{pmatrix} P^{[0]}_{n-1} & \ldots & P^{[n-1]}_{n-1} & 0_{\qq} \end{pmatrix} \begin{pmatrix} 0_{q \times nq} \\ S_{n-1} \end{pmatrix} \\ 
  &\quad - D_{n-1}D^{-1}_{n-2} \begin{pmatrix} P^{[0]}_{n-2} & \ldots & P^{[n-2]}_{n-2} & 0_{q \times 2q} \end{pmatrix} \begin{pmatrix} 0_{q \times nq} \\ S_{n-1} \end{pmatrix} \\
  & = \begin{pmatrix} 0_{\qq} & (\Ps_{n-1})^{[0]} & \ldots & (\Ps_{n-1})^{[n-1]} \end{pmatrix} \\
  &\quad - C_nD^{-1}_{n-1} \begin{pmatrix} (\Ps_{n-1})^{[0]} & \ldots & (\Ps_{n-1})^{[n-1]} & 0_{\qq} \end{pmatrix} \\ 
  &\quad - D_{n-1}D^{-1}_{n-2} \begin{pmatrix} (\Ps_{n-2})^{[0]} & \ldots & (\Ps_{n-2})^{[n-2]} & 0_{q \times 2q} \end{pmatrix}
\end{align*}
für alle $n \in \Zzfkp$.
Somit stimmt also die Folge $(Q_{n})^{\fklam{\kappa+1}}_{n=0}$ mit dem linken System von Matrixpolynomen zweiter Art $\Psnfk$ bezüglich $\sjk$ überein. \bwend

Mithilfe von \thref{fpsa1} und  \thref{mpsa5} können wir nun das linke System von Matrixpolynomen zweiter Art bezüglich einer Matrizenfolge, die bis zu einem gewissen Folgenglied Hankel-positiv definit ist, unter Verwendung des Favard-Paars bezüglich jener Folge darstellen (vergleiche \cite[Proposition 8.11(b)]{CR1} im Fall $\kappa=\infty$).

\begin{folg}	\thlabel{mpfo3}
  Seien $\kappa \in \Noa$ und $\sjk$ eine Folge aus $\Cqq$ mit \linebreak $\sjfkm \in \Hpqfkm$, falls $\kappa \geq 1$.
  Weiterhin seien $(Q_{n})^{\fklam{\kappa+1}}_{n=0}$ eine Folge von \textit{q}$\times$\textit{q}-Matrixpolynomen und im Fall $\kappa\geq1$ $\ABfk$ das Favard-Paar bezüglich $\sjk$. \dsfaa
  \begin{itemize}
    \item [\rm{(i)}] Es ist $(Q_{n})^{\fklam{\kappa+1}}_{n=0}$ das linke System von Matrixpolynomen zweiter Art bezüglich $\sjk$.
    \item [\rm{(ii)}] Es gelten $Q_{0} \equiv 0_{\qq}$, im Fall $\kappa\geq1$ $\;Q_{1} \equiv B_{0}$ und im Fall $\kappa \geq 3$
      \begin{align*}
	Q_{n}(z) = (z\Iq-A_{n-1})Q_{n-1}(z) - B^{\ast}_{n-1}Q_{n-2}(z)
      \end{align*}
      für alle $z \in \C$ und $n \in \Zzfkp$.  
  \end{itemize}
\end{folg}

\bwanf Der Fall $\kappa=0$ folgt unmittelbar aus \thref{mpsa5}. Seien nun $\kappa\geq1$ und $\CDfk$ die kanonische Hankel-Parametrisierung von $\sjk$. Unter Beachtung von \thref{mpbm3} gilt wegen $\sjfkm \in \Hpqfkm$ und \thref{mpsa1} dann $D_n\in\Cpqq$ für alle $n\in\Zofkm$. Hieraus folgt wegen Teil (a) von \thref{fpsa1} nun
\begin{align*}
	A_n = C_{n+1}D^{-1}_n
\end{align*}
für alle $n\in\Zofkm$, $B_0 = D_0$ und im Fall $\kappa\geq3$
\begin{align*}
	B^{\ast}_n = \rklam{D^{-1}_{n-1}D_n}^{\ast} = D^{\ast}_nD^{-\ast}_{n-1} = D_nD^{-1}_{n-1}
\end{align*}
für alle $n\in\Zefkm$. Hieraus folgt wegen \thref{mpsa5} dann die Behauptung. \bwend

%% file: sm2.tex
\newpage
\section[Die \texorpdfstring{$\alpha$}{a}-Dyukarev-Stieltjes-Parametrisierung von \texorpdfstring{$\alpha$}{a}-Stieltjes-positiv definiten Folgen]{Die \texorpdfstring{$\alpha$}{a}-Dyukarev-Stieltjes-Parametrisierung von \\ \texorpdfstring{$\alpha$}{a}-Stieltjes-positiv definiten Folgen} \label{chapadp}

In diesem Kapitel beschäftigen wir uns mit speziellen Aspekten der Struktur $\alpha$-Stieltjes-positiv definiter Folgen. Insbesondere werden wir die $\alpha$-Dyukarev-Stieltjes"=Parametrisierung einer $\alpha$-Stieltjes-positiv definiten Folge einführen und schauen uns einige Zusammenhänge zur $\alpha$-Stieltjes-Parametrisierung (vergleiche \thref{aspdef2}) an. 

Wir werden beide Fälle separat behandeln und beginnen mit dem rechtsseitigen Fall, also wenn die gegebene Folge rechtsseitig $\alpha$-Stieltjes-positiv definit ist. Wir bearbeiten dann den linksseitigen Fall durch geeignete Zurückführung auf den rechtsseitigen Fall.

\subsection{Der rechtsseitige Fall}

Im Mittelpunkt dieses Abschnitts steht die Diskussion einer speziellen inneren Parametrisierung rechtsseitig $\alpha$-Stieltjes-positiv definiter Folgen von komplexen \textit{q}$\times$\textit{q}-Matrizen, welche im Fall $q=1$ und $\alpha=0$ bereits auf T.-J. Stieltjes \cite{Sur} zurückgeht und von M.\,G. Krein in \cite[Anhang 2]{GK} mit einer mechanischen Interpretation versehen wurde. Die Behandlung der matriziellen Situation erfolgte im Fall $\alpha = 0$ erstmals bei Yu.\,M. Dyukarev \cite[Chapter 6]{dyu}.

Unsere Vorgehensweise orientiert sich an \cite[Chapter 8]{Trans}, wo der Fall $\sj\in{\cal K}^{>}_{q,\infty,0}$ als gegebene Folge behandelt wird. Dort findet man auch schon erste Ansätze für den Fall einer rechtsseitig $\alpha$-Stieltjes-positiv definiten Folge (vergleiche z.\,B. \cite[Lemma 8.14]{Trans} und \cite[Lemma 8.19]{Trans}), welche als Ausgangspunkt dieser Arbeit dienten.

Folgende Bemerkung ist von fundamentaler Bedeutung für unser weiteres Vorgehen. Sie erlaubt uns die Betrachtung der Inversen der für uns relevanten Hankel-Matrizen. Wir werden im weiteren Verlauf nicht mehr explizit dieses Resultat aufrufen.

\begin{bem}	\thlabel{adpbm1}
  Seien $\alpha \in \R$, $\kappa \in \Noa$ und $\sjk \in \Kpqka$. Dann sind $\Hn$ für alle $n \in \Zofk$ und 
  im Fall $\kappa \geq 1$ auch $\Harn$ für alle $n \in \Zofkm$ positiv hermitesch und insbesondere regulär.
\end{bem}

\bwanf  Wegen der Definition von $\Harn$ gilt
\begin{align}	\label{adpbm1bw1}
  \Harn=-\alpha\Hn+\Kn
\end{align}
für alle $n \in \Zofkm$. Im Fall $\kappa\in\gklam{0,1}$ folgt aus \thref{aspbm1} dann die Behauptung.

Im Fall $\kappa = 2n$ für ein $n \in \N$ gilt wegen Teil (a) von \thref{aspbm1} und \fref{adpbm1bw1} dann, dass $\Hn$ und $\Harnm$ positiv hermitesch und insbesondere regulär sind.
Seien nun $k \in \Zonm$ bzw. $l \in \Z{0}{n-2}$. Dann liefert ein wohlbekanntes Resultat über positiv hermitesche Matrizen,
dass $\Hk$ bzw. $\Harl$ als Hauptuntermatrix von $\Hn$ bzw. $\Harn$ ebenfalls positiv hermitesch und somit insbesondere regulär ist.

Im Fall $\kappa = 2n+1$ für ein $n \in \N$ gilt wegen Teil (b) von \thref{aspbm1} und \fref{adpbm1bw1} dann, dass $\Hn$ und $\Harn$ positiv hermitesch und insbesondere regulär sind.
Sei nun $k \in \Zonm$. Dann liefert ein wohlbekanntes Resultat über positiv hermitesche Matrizen,
dass $\Hk$ bzw. $\Hark$ als Hauptuntermatrix von $\Hn$ bzw. $\Harn$ ebenfalls positiv hermitesch und somit insbesondere regulär ist.

Im Fall $\kappa = \infty$ folgt wegen Teil (b) von \thref{asmdef2}, Teil (a) von \thref{aspbm1} und \fref{adpbm1bw1} dann die Behauptung.\bwend

Folgendes Lemma erlaubt uns eine Darstellung der Inversen der für uns relevanten Hankel-Matrizen (vergleiche \cite[Remark 2.2(d)]{CR1} für den Fall $\kappa=\infty$). Die allgemeine Formel für die Inverse einer Blockmatrix in Termen von Schurkomplementen wurde in \cite{Bana} behandelt.

\begin{lemma}	\thlabel{adplm2}
  Seien $\alpha \in \R$, $\kappa \in \Noa$ und $\sjk \in \Kpqka$. Dann gelten $H^{-1}_0 = s^{-1}_0$ und im Fall $\kappa \geq 1$ auch $\Haro^{-1} = s^{-1}_{\aro}$ sowie im Fall $\kappa \geq 2$
  \begin{align}	\label{adplm2bw1}
    \Hn^{-1} = \begin{pmatrix} \Hnm^{-1} & \Onq \\ \Oqn & \Oq \end{pmatrix} + \begin{pmatrix} -\Hnm^{-1}y_{n,2n-1} \\ \Iq \end{pmatrix} \dHn^{-1} \begin{pmatrix} -z_{n,2n-1}\Hnm^{-1} & \Iq \end{pmatrix}
  \end{align}
  für alle $n \in \Zefk$ und im Fall $\kappa \geq 3$
  \begin{align}	\label{adplm2bw2}
    \Harn^{-1} = \begin{pmatrix} \Harnm^{-1} & \Onq \\ \Oqn & \Oq \end{pmatrix} + \begin{pmatrix} -\Harnm^{-1}y_{\arn,2n-1} \\ \Iq \end{pmatrix} \dHarn^{-1} \begin{pmatrix} -z_{\arn,2n-1}\Harnm^{-1} & \Iq \end{pmatrix}
  \end{align}
  für alle $n \in \Zefkm$.  
\end{lemma}

\bwanf Aus der Definition von $H_0$ bzw. im Fall $\kappa\geq1$ von $H_{\aro}$ folgt sogleich $H^{-1}_0 = s^{-1}_0$ bzw. $\Haro^{-1} = s^{-1}_{\aro}$. Seien nun $\kappa \geq 2$ und $n \in \Zefk$. Dann gilt
\begin{align*}
	&\ \eklam{\begin{pmatrix} \Hnm^{-1} & \Onq \\ \Oqn & \Oq \end{pmatrix} + \begin{pmatrix} -\Hnm^{-1}y_{n,2n-1} \\ \Iq \end{pmatrix} \dHn^{-1} \begin{pmatrix} -z_{n,2n-1}\Hnm^{-1} & \Iq \end{pmatrix}} \begin{pmatrix} \Hnm & y_{n,2n-1} \\ z_{n,2n-1} & s_{2n} \end{pmatrix} \\
	&= \begin{pmatrix} \Inq & \Hnm^{-1} y_{n,2n-1} \\ \Oqn & \Oq \end{pmatrix} + \begin{pmatrix} -\Hnm^{-1}y_{n,2n-1} \\ \Iq \end{pmatrix} \dHn^{-1} \begin{pmatrix} \Oqn & \dHn \end{pmatrix} \\
	&=  \begin{pmatrix} \Inq & \Hnm^{-1} y_{n,2n-1} \\ \Oqn & \Oq \end{pmatrix} + \begin{pmatrix} -\Hnm^{-1}y_{n,2n-1} \\ \Iq \end{pmatrix} \begin{pmatrix} \Oqn & \Iq \end{pmatrix} \\
	&=  \begin{pmatrix} \Inq & \Hnm^{-1} y_{n,2n-1} \\ \Oqn & \Oq \end{pmatrix} + \begin{pmatrix} \Onqn & -\Hnm^{-1}y_{n,2n-1} \\ \Oqn & \Iq \end{pmatrix} = \Inpq.
\end{align*}
Hieraus folgt \fref{adplm2bw1}. Im Fall $\kappa\geq3$ kann man unter Beachtung von Teil (a) von \thref{aspdef1} dann \fref{adplm2bw2} analog beweisen. \bwend

Wir kommen nun zur zentralen Begriffsdefinition dieses Abschnitts. Hierfür benötigen wir noch folgende Bezeichnung.

\begin{bez}	\thlabel{adpbz1}
	Für $n\in\No$ sei $E_n:\C\rightarrow\C^{q\times(n+1)q}$ definiert gemäß $E_0(z):=\Iq$ bzw. im Fall $n\in\N$
	\begin{align*}
		\En(z) := \begin{pmatrix} \Iq \\ z\Iq \\ \vdots \\ z^{n}\Iq \end{pmatrix}.
	\end{align*}
\end{bez}

Es sei bemerkt, dass \thref{adpbz1} in Verbindung mit \thref{mpbz2} eine alternative Darstellung eines \textit{p}$\times$\textit{q}-Matrixpolynoms $P$ mit $\deg P=n$ für ein $n\in\No$ liefert: Es gilt $P(z) = \begin{pmatrix} P^{[0]} & \ldots & P^{[n]} \end{pmatrix} E_n(z)$ für alle $z\in\C$.

\begin{defi}	\thlabel{adpdef1}
  Seien $\alpha \in \R$, $\kappa \in \Na$ und $\sjk \in \Kpqka$. Weiterhin seien \linebreak $\Maro^{\sklam{s}} := s^{-1}_{0}$, $\Laro^{\sklam{s}} := s_{0}s^{-1}_{\ar 0}s_{0}$, im Fall $\kappa \geq 2$
  \begin{align*}
    \Marn^{\sklam{s}} := \En^{\ast}(\alpha)\big(\Hsn\big)^{-1}\En(\alpha)-\Enm^{\ast}(\alpha)\big(\Hsnm\big)^{-1}\Enm(\alpha)
  \end{align*}
  für alle $n \in \Zefk$ und im Fall $\kappa \geq 3$
  \begin{align*}
    \Larn^{\sklam{s}} := z^{\sklam{s}}_{0,n}\big(\Hsarn\big)^{-1}y^{\sklam{s}}_{0,n}-z^{\sklam{s}}_{0,n-1}\big(\Hsarnm\big)^{-1}y^{\sklam{s}}_{0,n-1}
  \end{align*}
  für alle $n \in \Zefkm$. 
  Dann heißt $\LMskarn$ die \textbf{rechtsseitige $\alpha$"=Dyukarev"=Stieltjes"=Parametrisierung} von $\sjk$. Falls klar ist, von welchem $\sjk$ die Rede ist, lassen wir das \anf{$\sklam{s}$} als oberen Index weg.
\end{defi}

Der Fall $\kappa = 0$ ist für uns wegen $\Kpqoa = \Hpqo$ nicht von Bedeutung, da hier $\alpha$ keine Rolle mehr spielt. Es sei bemerkt, dass die in \thref{adpdef1} eingeführten Größen mit denen in \cite[Definition 8.2]{Trans} im Fall $\alpha=0$ und $\kappa=\infty$ übereinstimmen.

Nun können wir die rechtsseitige $\alpha$-Dyukarev-Stieltjes-Parametrisierung einer rechtsseitig $\alpha$-Stieltjes-positiv definiten Folge mithilfe der rechtsseitigen $\alpha$"=Stieltjes"=Parametrisierung jener Folge darstellen (vergleiche \cite[Theorem 8.22]{Trans} im Fall $\alpha=0$ und $\kappa=\infty$).
Dafür benötigen wir einige Hilfsmittel, die wir im folgenden Lemma bereitstellen, das insbesondere eine Verbindung zu den in \thref{mpdef2} und \thref{mpdef3} eingeführten \textit{q}$\times$\textit{q}-Matrixpolynomen herstellt (vergleiche \cite[Lemma 8.14]{Trans}, \cite[Lemma 8.15]{Trans}, \cite[Lemma 8.19]{Trans} und \cite[Lemma 8.20]{Trans} für den Fall $\kappa=\infty$).

\begin{lemma}	\thlabel{adplm1}
  Seien $\alpha \in \R$, $\kappa \in \Na$ , $\sjk \in \Kpqka$, $\Qarjk$ die rechtsseitige $\alpha$-Stieltjes-Parametrisierung von $\sjk$ und  $\LMkarn$ die rechtsseitige $\alpha$-Dyukarev-Stieltjes-Parametrisierung von $\sjk$. 
  Weiterhin seien $\Pnfk$ das monische links-orthogonale System von Matrixpolynomen bezüglich $\sjk$ und $\Psnfk$ das linke System von Matrixpolynomen zweiter Art bezüglich $\sjk$. \dgfa
  \begin{itemize}
    \item [\rm{(a)}] Sei $n \in \Zofk$. Dann gilt
      \begin{align*}
	\Marn = P^{\ast}_{n}(\alpha)Q^{-1}_{\ar 2n} P_{n}(\alpha).
      \end{align*}
    \item [\rm{(b)}] Sei $n \in \Zefkp$. Dann gilt
      \begin{align*}
	P_{n}(\alpha) = (-1)^{n}\prodl^{n-1}_{j=0}Q_{\ar 2j+1}Q^{-1}_{\ar 2j}.
      \end{align*}
    \item [\rm{(c)}] Sei $n \in \Zofkm$. Dann gilt
      \begin{align*}
	\Larn & = \rklam{\Ps_{n+1}(\alpha)+Q_{\ar 2n+1}Q^{-1}_{\ar 2n}\Ps_{n}(\alpha)}^{\ast}Q^{-1}_{\ar 2n+1} \\
	&\quad \cdot \rklam{\Ps_{n+1}(\alpha)+Q_{\ar 2n+1}Q^{-1}_{\ar 2n}\Ps_{n}(\alpha)}.
      \end{align*}
    \item [\rm{(d)}] Seien $\kappa\geq3$ und $n \in \Zefkm$. Dann gilt
      \begin{align*}
	\Ps_{n+1}(\alpha)+Q_{\ar 2n+1}Q^{-1}_{\ar 2n}\Ps_{n}(\alpha) = (-1)^{n}Q_{\ar 2n}\prodl^{n-1}_{j=0}Q^{-1}_{\ar 2j+1}Q_{\ar 2j} .
      \end{align*}
  \end{itemize}
\end{lemma}


\bwanf Wegen Teil (c) von \thref{aspsa1} ist $\Qarjk$ eine Folge von regulären Matrizen.

Zu (a): Wegen \thref{mpsa3}, Teil (a) von \thref{aspdef2} und \thref{adpdef1} gilt
\begin{align*}
  P^{\ast}_{0}(\alpha)Q^{-1}_{\ar 0} P_{0}(\alpha) = \Iq^{\ast} s^{-1}_{0} \Iq = s^{-1}_{0} = \Maro.
\end{align*}
Seien nun $\kappa \geq 2$ und $n \in \Zefk$. Wegen \thref{adpdef1}, \thref{adplm2}, \thref{mpsa2} und Teil (a) von \thref{aspdef2} gilt dann
\begin{align*}
  \Marn & = \En^{\ast}(\alpha)\Hn^{-1}\En(\alpha)-\Enm^{\ast}(\alpha)\Hnm^{-1}\Enm(\alpha) \\
  & = \Ena(\alpha)\rklam{\Hn^{-1}-\diag(\Hnm^{-1},\Oq)}\En(\alpha) \\
  & = \Ena(\alpha)\begin{pmatrix} -z_{n,2n-1}\Hnm^{-1} & \Iq \end{pmatrix}^{\ast}\dHn^{-1}\begin{pmatrix} -z_{n,2n-1}\Hnm^{-1} & \Iq \end{pmatrix}\En(\alpha) \\
  & = P^{\ast}_{n}(\alpha)Q^{-1}_{\ar 2n} P_{n}(\alpha).
\end{align*}

Zu (b): 
Sei $\CDfk$ die kanonische Hankel-Parametrisierung von $\sjk$. 
Wegen \thref{mpsa3} gelten dann $P_{0} \equiv \Iq$,
\begin{align}	\label{adplm1bw1}
	P_{1}(z) = (z\Iq-C_{1}D^{-1}_{0})P_{0}(z)
\end{align}
für alle $z \in \C$ und im Fall $\kappa \geq 3$
\begin{align}	\label{adplm1bw2}
  	P_{n}(z) = (z\Iq-C_{n}D^{-1}_{n-1})P_{n-1}(z)-D_{n-1}D^{-1}_{n-2}P_{n-2}(z)
\end{align}
für alle $z \in \C$ und $n \in \Zzfkp$.
Wir zeigen die Behauptung mithilfe vollständiger Induktion. Wegen \fref{adplm1bw1}, \thref{mpdef1} und Teil (a) von \thref{aspdef2} gilt
\begin{align*}
	P_{1}(\alpha) = (\alpha\Iq-C_{1}D^{-1}_{0})P_{0}(\alpha) = \alpha\Iq-s_1s^{-1}_0 = -(-\alpha s_0+s_1)s^{-1}_0 = -Q_{\ar1}Q^{-1}_{\aro}.
\end{align*}
Seien nun $\kappa\geq3$, $n \in \Zzfkp$ und
\begin{align}	\label{adplm1bw3}
	P_k(\alpha) = (-1)^k\prodl^{k-1}_{j=0}Q_{\ar2j+1}Q^{-1}_{\ar2j}
\end{align}
für alle $k\in\Zenm$ erfüllt. Hieraus folgt wegen $P_{0} \equiv \Iq$ dann
\begin{align}	\label{adplm1bw3b}
	P_k(\alpha) = -Q_{\ar2k-1}Q^{-1}_{\ar2k-2}P_{k-1}(\alpha)
\end{align}
für alle $k\in\Zenm$. Wegen \cite[Lemma 6.9]{ot226} und Teil (a) von \thref{aspdef2} gilt
\begin{align*}
	Q_{\ar2n-1} &= s_{2n-1}-\Lambda_{n-1}-(\alpha\Iq+Q_{\ar2n-2}Q^{-1}_{\ar2n-3})Q_{\ar2n-2},
\end{align*}
also wegen \thref{mpdef1} und Teil (a) von \thref{aspdef2} gilt dann
\begin{align*}
	Q_{\ar2n-1}Q^{-1}_{\ar2n-2} &= C_nD^{-1}_{n-1}-\alpha\Iq-Q_{\ar2n-2}Q^{-1}_{\ar2n-3}.
\end{align*}
Hieraus folgt dann
\begin{align*}
	\alpha\Iq-C_nD^{-1}_{n-1} &= -Q_{\ar2n-1}Q^{-1}_{\ar2n-2}-Q_{\ar2n-2}Q^{-1}_{\ar2n-3}.
\end{align*}
Hieraus folgt wegen \fref{adplm1bw2}, \thref{mpdef1}, Teil (a) von \thref{aspdef2}, \fref{adplm1bw3b} und \fref{adplm1bw3} nun
\begin{align*}
	P_{n}(\alpha) &= (\alpha\Iq-C_{n}D^{-1}_{n-1})P_{n-1}(\alpha)-D_{n-1}D^{-1}_{n-2}P_{n-2}(\alpha) \\
	&= (-Q_{\ar2n-1}Q^{-1}_{\ar2n-2}-Q_{\ar2n-2}Q^{-1}_{\ar2n-3})P_{n-1}(\alpha)-Q_{\ar2n-2}Q^{-1}_{\ar2n-4}P_{n-2}(\alpha) \\
	&= -Q_{\ar2n-1}Q^{-1}_{\ar2n-2}P_{n-1}(\alpha)+Q_{\ar2n-2}Q^{-1}_{\ar2n-3}Q_{\ar2n-3}Q^{-1}_{\ar2n-4}P_{n-2}(\alpha) \\
	&\quad -Q_{\ar2n-2}Q^{-1}_{\ar2n-4}P_{n-2}(\alpha) \\
	&= -Q_{\ar2n-1}Q^{-1}_{\ar2n-2}P_{n-1}(\alpha) \\
	&= (-1)^{n}\prodl^{n-1}_{j=0}Q_{\ar 2j+1}Q^{-1}_{\ar 2j}.
\end{align*}

Zu (c): Wegen \thref{mpbem1}, Teil (a) von \thref{aspdef2} und \thref{adpdef1} gilt
\begin{align*}
  &\ \eklam{\Ps_{1}(\alpha)+Q_{\ar 1}Q^{-1}_{\ar 0}\Ps_{0}(\alpha)}^{\ast}Q^{-1}_{\ar 1}\eklam{\Ps_{1}(\alpha)+Q_{\ar 1}Q^{-1}_{\ar 0}\Ps_{0}(\alpha)} \\
  & = \rklam{s_{0}+s_{\ar 0}s^{-1}_{0}0_{\qq}}^{\ast}s^{-1}_{\ar 0}\rklam{s_{0}+s_{\ar 0}s^{-1}_{0}0_{\qq}} = s_{0}s^{-1}_{\ar 0}s_{0} = \Laro.
\end{align*}
Seien nun $\kappa \geq 2$ und $n \in \Zefkm$. Wegen (b) gilt dann
\begin{align}	\label{adplm1bw3c}
	P_{n+1}(\alpha) = -Q_{\ar2n+1}Q^{-1}_{\ar2n}P_{n}(\alpha).
\end{align}
Sei $\Parnfk$ das monische links-orthogonale System von Matrixpolynomen bezüglich $\sarjk$. Wegen \thref{mplm2}, Teil (a) von \thref{aspdef2} und \fref{adplm1bw3c} gilt dann
\begin{align*}
	(z-\alpha)P_{\arn}(z) = \eklam{P_{n+1}(z)+Q_{\ar2n+1}Q^{-1}_{\ar2n}P_{n}(z)} 
	- \eklam{P_{n+1}(\alpha)+Q_{\ar2n+1}Q^{-1}_{\ar2n}P_{n}(\alpha)}.
\end{align*}
Hieraus folgt wegen \cite[Remark 4.2]{Pert}, \cite[Lemma 4.3]{Pert} und \thref{mpsa2} dann
\begin{align*}
	\Ps_{n+1}(\alpha)+Q_{\ar2n+1}Q^{-1}_{\ar2n}\Ps_{n}(\alpha) &= \eklam{P_{n+1}(\alpha)+Q_{\ar2n+1}Q^{-1}_{\ar2n}P_{n}(\alpha)}^{\sklam{s}} \\
	&= \begin{pmatrix} P^{[0]}_{\arn} & \ldots & P^{[n]}_{\arn} \end{pmatrix} \yn \\
	&= \begin{pmatrix} -z_{\arn,2n-1}\Harnm^{-1} & \Iq \end{pmatrix} \yn.
\end{align*}
Hieraus folgt wegen \thref{adpdef1}, \thref{adplm2} und Teil (a) von \thref{aspdef2} nun
\begin{align*}
	\Larn &= z_{0,n}\Harn^{-1}y_{0,n}-z_{0,n-1}\Harnm^{-1}y_{0,n-1} \\
	&= z_{0,n}\rklam{\Harn^{-1}-\diag(\Harnm^{-1}, \Oq)}y_{0,n} \\
	&= z_{0,n}\begin{pmatrix} -\Harnm^{-1}y_{\arn,2n-1} \\ \Iq \end{pmatrix} \dHarn^{-1} \begin{pmatrix} -z_{\arn,2n-1}\Harnm^{-1} & \Iq \end{pmatrix}y_{0,n} \\
	&= \eklam{\Ps_{n+1}(\alpha)+Q_{\ar2n+1}Q^{-1}_{\ar2n}\Ps_{n}(\alpha)}^{\ast} Q^{-1}_{\ar2n+1} \eklam{\Ps_{n+1}(\alpha)+Q_{\ar2n+1}Q^{-1}_{\ar2n}\Ps_{n}(\alpha)}.
\end{align*}

Zu (d): 
Sei $\CDfk$ die kanonische Hankel-Parametrisierung von $\sjk$.
Wegen \thref{mpsa5} gelten dann $\Ps_{0} \equiv\Oq$, $\Ps_{1} \equiv D_0$ und
\begin{align}	\label{adplm1bw4}
  	\Ps_{n}(z) = (z\Iq-C_{n}D^{-1}_{n-1})\Ps_{n-1}(z)-D_{n-1}D^{-1}_{n-2}\Ps_{n-2}(z)
\end{align}
für alle $z \in \C$ und $n \in \Zzfkp$.
Wir zeigen die Behauptung mithilfe vollständiger Induktion. Wegen \cite[Lemma 6.9]{ot226} und Teil (a) von \thref{aspdef2} gilt
\begin{align*}
	Q_{\ar2n+1} &= s_{2n+1}-\Lambda_{n}-(\alpha\Iq+Q_{\ar2n}Q^{-1}_{\ar2n-1})Q_{\ar2n}
\end{align*}
für alle $n \in \Zefkm$, also wegen \thref{mpdef1} und Teil (a) von \thref{aspdef2} gilt dann
\begin{align*}
	Q_{\ar2n+1}Q^{-1}_{\ar2n} &= C_{n+1}D^{-1}_n-\alpha\Iq+Q_{\ar2n}Q^{-1}_{\ar2n-1}
\end{align*}
für alle $n \in \Zefkm$. Hieraus folgt dann
\begin{align}	\label{adplm1bw5}
	\alpha\Iq-C_{n+1}D^{-1}_n &= -Q_{\ar2n+1}Q^{-1}_{\ar2n}-Q_{\ar2n}Q^{-1}_{\ar2n-1}
\end{align}
für alle $n \in \Zefkm$. Hieraus folgt wegen \fref{adplm1bw4}, \thref{mpdef1} und Teil (a) von \thref{aspdef2} nun
\begin{align*}
  	\Ps_{2}(\alpha) &= (\alpha\Iq-C_{2}D^{-1}_{1})\Ps_{1}(\alpha)-D_{1}D^{-1}_{0}\Ps_{0}(\alpha) \\
  	&= (-Q_{\ar3}Q^{-1}_{\ar2}-Q_{\ar2}Q^{-1}_{\ar1})\Ps_{1}(\alpha)-Q_{\ar2}Q^{-1}_{\aro}\Ps_{0}(\alpha).
\end{align*}
Hieraus folgt wegen $\Ps_{0} \equiv\Oq$ und $\Ps_{1} \equiv Q_{\aro}$ dann
\begin{align*}
	\Ps_{2}(\alpha) + Q_{\ar3}Q^{-1}_{\ar2}\Ps_{1}(\alpha) = -Q_{\ar2}Q^{-1}_{\ar1}Q_{\aro}.
\end{align*}
Seien nun $\kappa \geq 5$, $n \in \Zzfkm$ und
\begin{align}	\label{adplm1bw6}
	\Ps_{k+1}(\alpha)+Q_{\ar2k+1}Q^{-1}_{\ar2k}\Ps_{k}(\alpha) = (-1)^{k}Q_{\ar2k}\prodl^{k-1}_{j=0}Q^{-1}_{\ar2j+1}Q_{\ar2j}
\end{align}
für alle $k\in\Zenm$ erfüllt. Hieraus folgt wegen $\Ps_{0} \equiv\Oq$ und $\Ps_{1} \equiv Q_{\aro}$ dann
\begin{align*}
	&\ \Ps_{k+1}(\alpha)+Q_{\ar2k+1}Q^{-1}_{\ar2k}\Ps_{k}(\alpha) \\
	&= -Q_{\ar2k}Q^{-1}_{\ar2k-1}\rklam{\Ps_{k}(\alpha)+Q_{\ar2k-1}Q^{-1}_{\ar2k-2}\Ps_{k-1}(\alpha)}
\end{align*}
für alle $k\in\Zenm$.
Hieraus folgt wegen \fref{adplm1bw4}, \fref{adplm1bw5}, \thref{mpdef1}, Teil (a) von \thref{aspdef2} und \fref{adplm1bw6} nun
\begin{align*}
	&\ \Ps_{n+1}(\alpha)+Q_{\ar2n+1}Q^{-1}_{\ar2n}\Ps_{n}(\alpha) \\
	&= (\alpha\Iq-C_{n+1}D^{-1}_{n})\Ps_{n}(\alpha)-D_{n}D^{-1}_{n-1}\Ps_{n-1}(\alpha)+Q_{\ar2n+1}Q^{-1}_{\ar2n}\Ps_{n}(\alpha) \\
	&= (-Q_{\ar2n+1}Q^{-1}_{\ar2n}-Q_{\ar2n}Q^{-1}_{\ar2n-1})\Ps_{n}(\alpha) \\
	&\quad -Q_{\ar2n}Q^{-1}_{\ar2n-2}\Ps_{n-1}(\alpha)+Q_{\ar2n+1}Q^{-1}_{\ar2n}\Ps_{n}(\alpha) \\
	&= -Q_{\ar2n}Q^{-1}_{\ar2n-1}\Ps_{n}(\alpha)-Q_{\ar2n}Q^{-1}_{\ar2n-1}Q_{\ar2n-1}Q^{-1}_{\ar2n-2}\Ps_{n-1}(\alpha) \\
	&= -Q_{\ar2n}Q^{-1}_{\ar2n-1}\eklam{\Ps_{n}(\alpha)+Q_{\ar2n-1}Q^{-1}_{\ar2n-2}\Ps_{n-1}(\alpha)} \\
	&= (-1)^{n}Q_{\ar 2n}\prodl^{n-1}_{j=0}Q^{-1}_{\ar 2j+1}Q_{\ar 2j}. \tag*{$\Box$}
\end{align*}

\begin{satz}	\thlabel{adpsa1}
  Seien $\alpha \in \R$, $\kappa \in \Na$, $\sjk \in \Kpqka$ und $\Qarjk$ die rechtsseitige $\alpha$-Stieltjes-Parametrisierung von $\sjk$. 
  Weiterhin sei $\LMkarn$ die rechtsseitige $\alpha$-Dyukarev-Stieltjes-Parametrisierung von $\sjk$.
  Dann gelten $\Maro = Q^{-1}_{\ar 0}$, 
  \begin{align*}
    \Larn = \rklam{\prodr^{n}_{j=0}Q_{\ar 2j}Q^{-1}_{\ar 2j+1}}Q_{\ar 2n+1}
    \rklam{\prodr^{n}_{j=0}Q_{\ar 2j}Q^{-1}_{\ar 2j+1}}^{\ast}
  \end{align*}
  für alle $n \in \Zofkm$ und im Fall $\kappa\geq2$
  \begin{align*}
    \Marn = \rklam{\prodr^{n-1}_{j=0}Q^{-1}_{\ar 2j}Q_{\ar 2j+1}} Q^{-1}_{\ar 2n}
    \rklam{\prodr^{n-1}_{j=0}Q^{-1}_{\ar 2j}Q_{\ar 2j+1}}^{\ast}
  \end{align*}
  für alle $n \in \Zefk$.
\end{satz}

\bwanf Wegen Teil (c) von \thref{aspsa1} ist $\Qarjk$ eine Folge von positiv hermiteschen und insbesondere regulären Matrizen. Wegen Teil (a) von \thref{aspdef2}, der Definition von $\widehat{H}_{0}$ und $\widehat{H}_{\ar 0}$ sowie \thref{adpdef1} gelten
\begin{align*}
  Q^{-1}_{\ar 0} = \widehat{H}^{-1}_{0} = s^{-1}_{0} = \Maro
\end{align*}
und
\begin{align*}
  &\ (Q_{\ar 0}Q^{-1}_{\ar 1})Q_{\ar 1}(Q_{\ar 0}Q^{-1}_{\ar 1})^{\ast} 
  = Q_{\ar 0}Q^{-1}_{\ar 1}Q_{\ar 1}Q^{-\ast}_{\ar 1}Q^{\ast}_{\ar 0}\\
  & = Q_{\ar 0}Q^{-1}_{\ar 1}Q_{\ar 0}
  = \widehat{H}_{0}\widehat{H}^{-1}_{\ar 0}\widehat{H}_{0} = s_{0}s^{-1}_{\ar 0}s_{0} = \Laro.
\end{align*}
Seien nun $\kappa \geq 2$ und $\Pnfk$ das monische links-orthogonale System von Matrixpolynomen bezüglich $\sjk$. Wegen der Teile (a) und (b) von \thref{adplm1} gilt dann
\begin{align*}
  \Marn & = P^{\ast}_{n}(\alpha)Q_{\ar 2n} P_{n}(\alpha) \\
  & = (-1)^{n}\rklam{\prodl^{n-1}_{j=0}Q_{\ar 2j+1}Q^{-1}_{\ar 2j}}^{\ast} Q_{\ar 2n} (-1)^{n}\prodl^{n-1}_{j=0}Q_{\ar 2j+1}Q^{-1}_{\ar 2j} \\
  & = \rklam{\prodr^{n-1}_{j=0}Q^{-1}_{\ar 2j}Q_{\ar 2j+1}} Q_{\ar 2n} \rklam{\prodr^{n-1}_{j=0}Q^{-1}_{\ar 2j}Q_{\ar 2j+1}}^{\ast}
\end{align*}
für alle $n \in \Zefk$. Seien nun $\kappa \geq 3$ und $\Psnfk$ das linke System von Matrixpolynomen zweiter Art bezüglich $\sjk$. Wegen Teil (d) von \thref{adplm1} gelten
\begin{align*}
  &\ \Ps_{n+1}(\alpha)+Q_{\ar 2n+1}Q^{-1}_{\ar 2n}\Ps_{n}(\alpha) = (-1)^{n}Q_{\ar 2n}\prodl^{n-1}_{j=0}Q^{-1}_{\ar 2j+1}Q_{\ar 2j} \\
  & = (-1)^{n}Q_{\ar 2n+1}Q^{-1}_{\ar 2n+1}Q_{\ar 2n}\prodl^{n-1}_{j=0}Q^{-1}_{\ar 2j+1}Q_{\ar 2j} \\
  & = (-1)^{n}Q_{\ar 2n+1}\prodl^{n}_{j=0}Q^{-1}_{\ar 2j+1}Q_{\ar 2j}
\end{align*}
und
\begin{align*}
  \rklam{\Ps_{n+1}(\alpha)+Q_{\ar 2n+1}Q^{-1}_{\ar 2n}\Ps_{n}(\alpha)}^{\ast} &= (-1)^{n}\rklam{Q_{\ar 2n+1}\prodl^{n}_{j=0}Q^{-1}_{\ar 2j+1}Q_{\ar 2j}}^{\ast} \\
  & = (-1)^{n}\rklam{\prodr^{n}_{j=0}Q_{\ar 2j}Q^{-1}_{\ar 2j+1}} Q_{\ar 2n+1}
\end{align*}
für alle $n \in \Zefkm$. Hieraus und aus Teil (c) von \thref{adplm1} folgt dann
\begin{align*}
  \Larn & = \eklam{\Ps_{n+1}(\alpha)+Q_{\ar 2n+1}Q^{-1}_{\ar 2n}\Ps_{n}(\alpha)}^{\ast}Q^{-1}_{\ar 2n+1}\eklam{\Ps_{n+1}(\alpha)+Q_{\ar 2n+1}Q^{-1}_{\ar 2n}\Ps_{n}(\alpha)} \\
  & = (-1)^{n}\rklam{\prodr^{n}_{j=0}Q_{\ar 2j}Q^{-1}_{\ar 2j+1}} Q_{\ar 2n+1}Q^{-1}_{\ar 2n+1}(-1)^{n}Q_{\ar 2n+1}\prodl^{n}_{j=0}Q^{-1}_{\ar 2j+1}Q_{\ar 2j} \\
  & = \rklam{\prodr^{n}_{j=0}Q_{\ar 2j}Q^{-1}_{\ar 2j+1}} Q_{\ar 2n+1} \rklam{\prodr^{n}_{j=0}Q_{\ar 2j}Q^{-1}_{\ar 2j+1}}^{\ast}
\end{align*}
für alle $n \in \Zefkm$. \bwend

Folgendes Resultat zeigt nun, dass die einzelnen Matrizen der rechtsseitigen $\alpha$"=Dyukarev"=Stieltjes"=Parametrisierung einer rechtsseitig $\alpha$-Stieltjes-positiv definiten Folge jeweils positiv hermitesch und somit regulär sind (vergleiche \cite[Remark 8.23]{Trans} für den Fall $\alpha=0$ und $\kappa=\infty$).

\begin{bem}	\thlabel{adpbm3}
  Seien $\alpha \in \R$, $\kappa \in \Na$,  $\sjk \in \Kpqka$ und $\LMkarn$ die rechtsseitige $\alpha$"=Dyukarev"=Stieltjes"=Parametrisierung von $\sjk$. 
  Dann sind $\Larn$ für alle $n \in \Zofkm$ und $\Marn$ für alle $n \in \Zofk$ positiv hermitesch und insbesondere regulär.
\end{bem}

\bwanf Dies folgt wegen Teil (c) von \thref{aspsa1} und \thref{adpsa1} aus \thref{amlm1}. \bwend

Umgekehrt können wir nun die rechtsseitige $\alpha$-Stieltjes-Parametrisierung einer rechtsseitig $\alpha$-Stieltjes-positiv definiten Folge mithilfe der rechtsseitigen $\alpha$-Dyukarev-Stieltjes-Parametrisierung jener Folge beschreiben (vergleiche \cite[Theorem 8.24]{Trans} für den Fall $\alpha=0$ und $\kappa=\infty$).

\begin{satz}	\thlabel{adpsa2}
  Seien $\alpha \in \R$, $\kappa \in \Na$,  $\sjk \in \Kpqka$ und $\Qarjk$ die rechtsseitige $\alpha$-Stieltjes-Parametrisierung von $\sjk$. Weiterhin sei $\LMkarn$ die rechtsseitige $\alpha$-Dyukarev-Stieltjes-Parametrisierung von $\sjk$.
  Dann gelten $Q_{\aro} = \Maro^{-1}$,
  \begin{align*}
    Q_{\ar 2n+1} = \rklam{\prodr^{n}_{j=0}\Marj\Larj}^{-\ast} \Larn
    \rklam{\prodr^{n}_{j=0}\Marj\Larj}^{-1}
  \end{align*}
  für alle $n\in \Zofkm$ und im Fall $\kappa\geq2$  
  \begin{align*}
    Q_{\ar 2n} = \rklam{\prodr^{n-1}_{j=0}\Marj\Larj}^{-\ast} \Marn^{-1}
    \rklam{\prodr^{n-1}_{j=0}\Marj\Larj}^{-1}
  \end{align*}
  für alle $n\in \Zofk$.  
\end{satz}

\bwanf Dies folgt wegen Teil (c) von \thref{aspsa1} und \thref{adpsa1} aus \thref{amlm1}. \bwend

Wir können nun rekursiv die einzelnen Folgenglieder einer rechtsseitig $\alpha$-Stieltjes-positiv definiten Folge mithilfe ihrer $\alpha$-Dyukarev-Stieltjes-Parametrisierung ausdrücken (vergleiche \cite[Proposition 8.26]{Trans} für den Fall $\alpha=0$ und $\kappa=\infty$).

\begin{satz}	\thlabel{adpbm4}
	Seien $\alpha \in \R$, $\kappa \in \Na$, $\sjk \in \Kpqka$ und $\LMkarn$ die rechtsseitige $\alpha$-Dyukarev-Stieltjes-Parametrisierung von $\sjk$. Dann gelten $s_0 = \Maro^{-1}$, $s_1 = \alpha s_0 + (\Maro\Laro)^{-\ast}\Laro(\Maro\Laro)^{-1}$, im Fall $\kappa\geq2$
	\begin{align*}
		s_{2n} = \rklam{\prodr^{n-1}_{j=0}\Marj\Larj}^{-\ast} \Marn^{-1}
    \rklam{\prodr^{n-1}_{j=0}\Marj\Larj}^{-1}+z_{n,2n-1}\Hnm^{-1}y_{n,2n-1}
	\end{align*}
	für alle $n\in \Zefk$ und im Fall $\kappa\geq3$
	\begin{align*}
		s_{2n+1} = \alpha s_{2n} + \rklam{\prodr^{n}_{j=0}\Marj\Larj}^{-\ast} \Larn
    \rklam{\prodr^{n}_{j=0}\Marj\Larj}^{-1}+z_{\arn,2n-1}\Harnm^{-1}y_{\arn,2n-1}
	\end{align*}
	für alle $n\in\Zefkm$.
\end{satz}

\bwanf Dies folgt aus Teil (a) von \thref{aspbm5} und \thref{adpsa2}. \bwend

Wir können nun bei beliebig vorgegebenem $\alpha\in\R$ mithilfe zweier Folgen von positiv hermiteschen Matrizen eine rechtsseitige $\alpha$-Stieltjes-positiv definite Folge konstruieren, sodass jene zwei Folgen die rechtsseitige $\alpha$-Dyukarev-Stieltjes-Parametrisierung der konstruierten Folge bilden (vergleiche \cite[Proposition 8.27]{Trans} für den Fall $\alpha=0$ und $\kappa=\infty$).

\begin{satz}	\thlabel{adpbm5}
	Seien $\alpha \in \R$, $\kappa \in \Na$ sowie $(\Larn)^{\fklam{\kappa-1}}_{n=0}$ und $(\Marn)^{\fklam{\kappa}}_{n=0}$ Folgen von positiv hermiteschen Matrizen aus $\Cqq$. Weiterhin seien durch rekursive Konstruktion \linebreak $s_0 := \Maro^{-1}$, $s_1 := \alpha s_0 + (\Maro\Laro)^{-\ast}\Laro(\Maro\Laro)^{-1}$, im Fall $\kappa\geq2$
	\begin{align*}
		s_{2n} := \rklam{\prodr^{n-1}_{j=0}\Marj\Larj}^{-\ast} \Marn^{-1}
    \rklam{\prodr^{n-1}_{j=0}\Marj\Larj}^{-1}+z_{n,2n-1}\Hnm^{-1}y_{n,2n-1}
	\end{align*}
	für alle $n\in \Zefk$ und im Fall $\kappa\geq3$
	\begin{align*}
		s_{2n+1} := \alpha s_{2n} + \rklam{\prodr^{n}_{j=0}\Marj\Larj}^{-\ast} \Larn
    \rklam{\prodr^{n}_{j=0}\Marj\Larj}^{-1}+z_{\arn,2n-1}\Harnm^{-1}y_{\arn,2n-1}
	\end{align*}
	für alle $n\in\Zefkm$. Dann gilt $\sjk\in\Kpqka$ und $\LMkarn$ ist die rechtsseitige $\alpha$-Dyukarev-Stieltjes-Parametrisierung von $\sjk$.	
\end{satz}

\bwanf Sei $\Qarjk$ die rechtsseitige $\alpha$-Stieltjes-Parametrisierung von $\sjk$. Wegen Teil (a) von \thref{aspbm5} gelten dann $s_0 = Q_{\aro}$, $s_1 = \alpha s_0 + Q_{\ar1}$, im Fall $\kappa\geq2$ 
\begin{align*}
	s_{2n} = Q_{\ar2n}+z_{n,2n-1}\Hnm^{+}y_{n,2n-1}
\end{align*}
für alle $n\in\Zefk$ und im Fall $\kappa\geq3$
\begin{align*}
	s_{2n+1} = \alpha s_{2n} + Q_{\ar2n+1}+z_{\arn,2n-1}\Harnm^{+}y_{\arn,2n-1}
\end{align*}
für alle $n\in\Zefkm$. Hieraus folgt aus der Definition der Folge $\sjk$ dann \linebreak $Q_{\aro} = \Maro^{-1}$,
\begin{align}	\label{adpbm4bw1}
  	Q_{\ar 2n+1} = \rklam{\prodr^{n}_{j=0}\Marj\Larj}^{-\ast} \Larn
    \rklam{\prodr^{n}_{j=0}\Marj\Larj}^{-1}
\end{align}
für alle $n\in \Zofkm$ und im Fall $\kappa\geq2$  
\begin{align}	\label{adpbm4bw2}
    Q_{\ar 2n} = \rklam{\prodr^{n-1}_{j=0}\Marj\Larj}^{-\ast} \Marn^{-1}
    \rklam{\prodr^{n-1}_{j=0}\Marj\Larj}^{-1}
\end{align}
für alle $n\in \Zefk$. Hieraus folgt aus der Tatsache, dass $(\Larn)^{\fklam{\kappa-1}}_{n=0}$ und $(\Marn)^{\fklam{\kappa}}_{n=0}$ Folgen von positiv hermiteschen Matrizen aus $\Cqq$ sind, dann, dass auch $\Qarjk$ eine Folge von positiv hermiteschen Matrizen aus $\Cqq$ ist. Somit gilt wegen Teil (c) von \thref{aspsa1} dann $\sjk\in\Kpqka$. Wegen $Q_{\aro} = \Maro^{-1}$, \fref{adpbm4bw1}, \fref{adpbm4bw2} und \thref{amlm3} gelten $\Maro = Q^{-1}_{\ar 0}$, 
\begin{align*}
    \Larn = \rklam{\prodr^{n}_{j=0}Q_{\ar 2j}Q^{-1}_{\ar 2j+1}}Q_{\ar 2n+1}
    \rklam{\prodr^{n}_{j=0}Q_{\ar 2j}Q^{-1}_{\ar 2j+1}}^{\ast}
\end{align*}
für alle $n \in \Zofkm$ und im Fall $\kappa\geq2$
\begin{align*}
    \Marn = \rklam{\prodr^{n-1}_{j=0}Q^{-1}_{\ar 2j}Q_{\ar 2j+1}} Q^{-1}_{\ar 2n}
    \rklam{\prodr^{n-1}_{j=0}Q^{-1}_{\ar 2j}Q_{\ar 2j+1}}^{\ast}
\end{align*}
für alle $n \in \Zefk$. Hieraus folgt wegen \thref{adpsa1} dann, dass $\LMkarn$ die rechtsseitige $\alpha$-Dyukarev-Stieltjes-Parametrisierung von $\sjk$ ist. \bwend

\subsection{Der linksseitige Fall} \label{chapadpl}

Folgende Bemerkung ist von fundamentaler Bedeutung für unser weiteres Vorgehen. Sie erlaubt uns die Betrachtung der Inversen der für uns relevanten Hankel-Matrizen. Wir werden im weiteren Verlauf nicht mehr explizit dieses Resultat aufrufen.

\begin{bem}	\thlabel{adplbm1}
  Seien $\alpha \in \R$, $\kappa \in \Noa$ und $\sjk \in \Lpqka$. Dann sind $\Hn$ für alle $n \in \Zofk$ und 
  im Fall $\kappa \geq 1$ auch $\Haln$ für alle $n \in \Zofkm$ positiv hermitesch und insbesondere regulär.
\end{bem}

\bwanf Sei $t_j := (-1)^js_j$ für alle $j\in\eklam{\kappa}_0$. Wegen Teil (a) von \thref{asmbm3} gilt dann $\tjk\in\Kpqkma$. Hieraus folgt unter Beachtung von Teil (b) von \thref{asmlm1} und Teil (b) von \thref{asplm1} wegen \thref{adpbm1} dann die Behauptung. 
\bwend

Folgendes Lemma erlaubt uns eine Darstellung der Inversen der für uns relevanten Hankel-Matrizen.

\begin{lemma}	\thlabel{adpllm2}
  Seien $\alpha \in \R$, $\kappa \in \Noa$ und $\sjk \in \Lpqka$. Dann gelten $H^{-1}_0 = s^{-1}_0$ und im Fall $\kappa\geq1$ auch $\Haro^{-1} = s^{-1}_{\aro}$ sowie im Fall $\kappa \geq 2$
  \begin{align}	\label{adpllm2bw1}
    \Hn^{-1} = \begin{pmatrix} \Hnm^{-1} & \Onq \\ \Oqn & \Oq \end{pmatrix} + \begin{pmatrix} -\Hnm^{-1}y_{n,2n-1} \\ \Iq \end{pmatrix} \dHn^{-1} \begin{pmatrix} -z_{n,2n-1}\Hnm^{-1} & \Iq \end{pmatrix}
  \end{align}
  für alle $n \in \Zefk$ und im Fall $\kappa \geq 3$
  \begin{align}	\label{adpllm2bw2}
    \Haln^{-1} = \begin{pmatrix} \Halnm^{-1} & \Onq \\ \Oqn & \Oq \end{pmatrix} + \begin{pmatrix} -\Halnm^{-1}y_{\aln,2n-1} \\ \Iq \end{pmatrix} \dHarn^{-1} \begin{pmatrix} -z_{\aln,2n-1}\Halnm^{-1} & \Iq \end{pmatrix}
  \end{align}
  für alle $n \in \Zefkm$.  
\end{lemma}

\bwanf  Aus der Definition von $H_0$ bzw. im Fall $\kappa\geq1$ von $H_{\alo}$ folgt sogleich $H^{-1}_0 = s^{-1}_0$ bzw. $\Halo^{-1} = s^{-1}_{\aro}$. Der Beweis von \fref{adpllm2bw1} im Fall $\kappa\geq2$ wurde schon im Beweis von \thref{adplm2} ausgeführt. Im Fall $\kappa\geq3$ kann man unter Beachtung von Teil (b) von \thref{aspdef1} dann \fref{adpllm2bw2} analog beweisen.  \bwend

Wir kommen nun zur zentralen Begriffsdefinition dieses Abschnitts.

\begin{defi}	\thlabel{adpldef1}
  Seien $\alpha \in \R$, $\kappa \in \Na$ und $\sjk \in \Lpqka$. Weiterhin seien \linebreak $\Malo^{\sklam{s}} := s^{-1}_{0}$, $\Lalo^{\sklam{s}} := s_{0}s^{-1}_{\al 0}s_{0}$, im Fall $\kappa \geq 2$
  \begin{align*}
    \Maln^{\sklam{s}} := \En^{\ast}(\alpha)\brklam{\Hsn}^{-1}\En(\alpha)-\Enm^{\ast}(\alpha)\brklam{\Hsnm}^{-1}\Enm(\alpha)
  \end{align*}
  für alle $n \in \Zefk$ und im Fall $\kappa \geq 3$
  \begin{align*}
    \Laln^{\sklam{s}} := z^{\sklam{s}}_{0,n}\brklam{\Hsaln}^{-1}y^{\sklam{s}}_{0,n}-z^{\sklam{s}}_{0,n-1}\brklam{\Hsalnm}^{-1}y^{\sklam{s}}_{0,n-1}
  \end{align*}
  für alle $n \in \Zefkm$.
  Dann heißt $\LMskaln$ die \textbf{linksseitige $\alpha$"=Dyukarev"=Stieltjes"=Parametrisierung} von $\sjk$. Falls klar ist, von welchem $\sjk$ die Rede ist, lassen wir das \anf{$\sklam{s}$} als oberen Index weg.
\end{defi}

Der Fall $\kappa = 0$ ist für uns wegen $\Lpqoa = \Hpqo$ nicht von Bedeutung, da hier $\alpha$ keine Rolle mehr spielt.

Wir werden die folgenden Resultate hauptsächlich mithilfe den entsprechenden Resultaten für den rechtsseitigen Fall beweisen. Hierfür wird uns das nächste Lemma den benötigten Zusammenhang liefern.

\begin{lemma}	\thlabel{adpllm1}
	Seien $\alpha \in \R$, $\kappa \in \Na$, $\sjk$ eine Folge aus $\Cqq$ und $t_j:=(-1)^js_j$ für alle $j\in\Zok$. \dgfa
	\begin{itemize}
		\item [\rm{(a)}] \esfaa
		\begin{itemize}
			\item [\rm{(i)}] Es gilt $\sjk\in\Kpqka$.
			\item [\rm{(ii)}] Es gilt $\tjk\in\Lpqkma$.
		\end{itemize}
		\item [\rm{(b)}] Seien $n\in\No$ und $z\in\C$. Dann gilt
		\begin{align*}
			\En(-z) = \Vn\En(z).
		\end{align*}
		\item [\rm{(c)}] Sei nun {\rm (i)} erfüllt. Weiterhin seien $\LMskarn$ die rechtsseitige $\alpha$-Dyukarev-Stieltjes-Parametrisierung von $\sjk$ und $\LMtkmaln$ die linksseitige $-\alpha$-Dyukarev-Stieltjes-Parametrisierung von $\tjk$. Dann gelten
		\begin{align*}
			\Lmaln^{\sklam{t}} = \Larn^{\sklam{s}}
		\end{align*}
		für alle $n \in \Zofkm$ und
		\begin{align*}
			\Mmaln^{\sklam{t}} = \Marn^{\sklam{s}}
		\end{align*}
		für alle $n \in \Zofk$.
	\end{itemize}
\end{lemma}

\bwanf Zu (a): Dies folgt aus Teil (a) von \thref{asmbm3}.

Zu (b): Es gelten $E_0(z)=\Iq=E_0(-z)$ und
\begin{align*}
	\Vn\En(z) = \begin{pmatrix} \Iq \\ -z\Iq \\ \vdots \\ (-1)^nz^n\Iq \end{pmatrix} = \En(-z).
\end{align*}

Zu (c): Es gilt
\begin{align*}
	\Mmalo^{\sklam{t}} = t^{-1}_0 = s^{-1}_0 = \Maro^{\sklam{s}}.
\end{align*}
Wegen Teil (a) von \thref{asplm1} gilt weiterhin
\begin{align*}
	\Lmalo^{\sklam{t}} = t_0t^{-1}_{-\alo}t_0 = s_0s^{-1}_{\alo}s_0 = \Laro^{\sklam{s}}.
\end{align*}
Sei nun $\kappa \geq 2$ und $n \in \Zefk$. Wegen (b) und der Teile (a) und (b) von \thref{asmlm1} gilt dann
\begin{align*}
	&\ \Mmaln^{\sklam{t}} = \En^{\ast}(-\alpha)\brklam{\Hn^{\sklam{t}}}^{-1}\En(-\alpha)-\Enm^{\ast}(-\alpha)\brklam{\Hnm^{\sklam{t}}}^{-1}\Enm(-\alpha) \\
	&= \En^{\ast}(\alpha)\Vna\Vn\brklam{\Hsn}^{-1}\Vna\Vn\En(\alpha)-\Enm^{\ast}(\alpha)\Vnma\Vnm\brklam{\Hsn}^{-1}\Vnma\Vnm\Enm(\alpha) \\
	&= \En^{\ast}(\alpha)\brklam{\Hsn}^{-1}\En(\alpha)-\Enm^{\ast}(\alpha)\brklam{\Hsn}^{-1}\Enm(\alpha)
	= \Marn^{\sklam{s}}.
\end{align*}
Seien nun $\kappa \geq 3$ und $n \in \Zefkm$. Wegen der Teile (a) und (d) von \thref{asmlm1} sowie Teil (b) von \thref{asplm1} gilt dann
\begin{align*}
	\Lmaln^{\sklam{t}} &= z^{\sklam{t}}_{0,n}\brklam{\Haln^{\sklam{t}}}^{-1}y^{\sklam{t}}_{0,n}-z^{\sklam{t}}_{0,n-1}\brklam{\Halnm^{\sklam{t}}}^{-1}y^{\sklam{t}}_{0,n-1} \\
	&= z^{\sklam{s}}_{0,n}\Vna\Vn\brklam{\Hsarn}^{-1}\Vna\Vn y^{\sklam{s}}_{0,n}-z^{\sklam{s}}_{0,n-1}\Vnma\Vnm\brklam{\Hsarnm}^{-1}\Vnma\Vnm y^{\sklam{s}}_{0,n-1} \\
	&= z^{\sklam{s}}_{0,n}\brklam{\Hsarn}^{-1}y^{\sklam{s}}_{0,n}-z^{\sklam{s}}_{0,n-1}\brklam{\Hsarnm}^{-1}y^{\sklam{s}}_{0,n-1}
	= \Larn^{\sklam{s}}. \tag*{$\Box$}
\end{align*}

Nun können wir die linksseitige $\alpha$-Dyukarev-Stieltjes-Parametrisierung einer linksseitig $\alpha$-Stieltjes-positiv definiten Folge mithilfe der linksseitigen $\alpha$"=Stieltjes"=Parametrisierung jener Folge darstellen.

\begin{satz}	\thlabel{adplsa1}
  Seien $\alpha \in \R$, $\kappa \in \Na$, $\sjk \in \Lpqka$ und $\Qaljk$ die linksseitige $\alpha$-Stieltjes-Parametrisierung von $\sjk$. 
  Weiterhin sei $\LMkaln$ die linksseitige $\alpha$-Dyukarev-Stieltjes-Parametrisierung von $\sjk$.
  Dann gelten $\Malo = Q^{-1}_{\al 0}$,
  \begin{align*}
    \Laln = \rklam{\prodr^{n}_{j=0}Q_{\al 2j}Q^{-1}_{\al 2j+1}}Q_{\al 2n+1}
    \rklam{\prodr^{n}_{j=0}Q_{\al 2j}Q^{-1}_{\al 2j+1}}^{\ast}
  \end{align*}
  für alle $n \in \Zofkm$ und im Fall $\kappa\geq2$
  \begin{align*}
    \Maln = \rklam{\prodr^{n-1}_{j=0}Q^{-1}_{\al 2j}Q_{\ar 2j+1}} Q^{-1}_{\al 2n}
    \rklam{\prodr^{n-1}_{j=0}Q^{-1}_{\al 2j}Q_{\al 2j+1}}^{\ast}
  \end{align*}
  für alle $n \in \Zefk$.
\end{satz}

\bwanf Dies folgt wegen \thref{aspbm4} und der Teile (a) und (c) von \thref{adpllm1} aus \thref{adpsa1}. \bwend

Folgendes Resultat zeigt nun, dass die einzelnen Matrizen der linksseitigen $\alpha$-Dyukarev-Stieltjes-Parametrisierung einer linksseitig $\alpha$-Stieltjes-positiv definiten Folge jeweils positiv hermitesch und somit regulär sind.

\begin{bem}	\thlabel{adplbm2}
  Seien $\alpha \in \R$, $\kappa \in \Na$,  $\sjk \in \Lpqka$ und $\LMkaln$ die linksseitige $\alpha$-Dyukarev-Stieltjes-Parametrisierung von $\sjk$. 
  Dann sind $\Laln$ für alle $n \in \Zofkm$ und $\Maln$ für alle $n \in \Zofk$ positiv hermitesch und insbesondere regulär.
\end{bem}

\bwanf Dies folgt wegen der Teile (a) und (c) von \thref{adpllm1} aus \thref{adpbm3}. \bwend

Umgekehrt können wir nun die linksseitige $\alpha$-Stieltjes-Parametrisierung einer linksseitig $\alpha$-Stieltjes-positiv definiten Folge mithilfe der linksseitigen $\alpha$-Dyukarev-Stieltjes-Parametrisierung jener Folge beschreiben.

\begin{satz}	\thlabel{adplsa2}
  Seien $\alpha \in \R$, $\kappa \in \Na$,  $\sjk \in \Lpqka$ und $\Qaljk$ die linksseitige $\alpha$-Stieltjes-Parametrisierung von $\sjk$. Weiterhin sei $\LMkaln$ die linksseitige $\alpha$-Dyukarev-Stieltjes-Parametrisierung von $\sjk$.
  Dann gelten $Q_{\alo} = \Malo^{-1}$,
  \begin{align*}
    Q_{\al 2n+1} = \rklam{\prodr^{n}_{j=0}\Malj\Lalj}^{-\ast} \Laln
    \rklam{\prodr^{n}_{j=0}\Malj\Lalj}^{-1}
  \end{align*}
  für alle $n\in \Zofkm$ und im Fall $\kappa\geq2$
  \begin{align*}
    Q_{\al 2n} = \rklam{\prodr^{n-1}_{j=0}\Malj\Lalj}^{-\ast} \Maln^{-1}
    \rklam{\prodr^{n-1}_{j=0}\Malj\Lalj}^{-1}
  \end{align*}
  für alle $n\in \Zofk$.
\end{satz}

\bwanf Dies folgt wegen \thref{aspbm4} und der Teile (a) und (c) von \thref{adpllm1} aus \thref{adpsa2}. \bwend

Wir können nun rekursiv die einzelnen Folgenglieder einer linksseitig $\alpha$-Stieltjes-positiv definiten Folge mithilfe ihrer $\alpha$-Dyukarev-Stieltjes-Parametrisierung ausdrücken.

\begin{satz}	\thlabel{adplbm3}
	Seien $\alpha \in \R$, $\kappa \in \Na$, $\sjk \in \Lpqka$ und $\LMkaln$ die linksseitige $\alpha$-Dyukarev-Stieltjes-Parametrisierung von $\sjk$. Dann gelten $s_0 = \Malo^{-1}$, $s_1 = \alpha s_0 - (\Malo\Lalo)^{-\ast}\Lalo(\Malo\Lalo)^{-1}$, im Fall $\kappa\geq2$
	\begin{align*}
		s_{2n} = \rklam{\prodr^{n-1}_{j=0}\Malj\Lalj}^{-\ast} \Maln^{-1}
    \rklam{\prodr^{n-1}_{j=0}\Malj\Lalj}^{-1}+z_{n,2n-1}\Hnm^{+}y_{n,2n-1}
	\end{align*}
	für alle $n\in \Zefk$ und im Fall $\kappa\geq3$
	\begin{align*}
		s_{2n+1} = \alpha s_{2n} - \rklam{\prodr^{n}_{j=0}\Malj\Lalj}^{-\ast} \Laln
    \rklam{\prodr^{n}_{j=0}\Malj\Lalj}^{-1}-z_{\aln,2n-1}\Halnm^{+}y_{\aln,2n-1}
	\end{align*}
	für alle $n\in\Zefkm$.
\end{satz}

\bwanf Dies folgt aus Teil (b) von \thref{aspbm5} und \thref{adplsa2}. \bwend

Wir können nun bei beliebig vorgegebenen $\alpha\in\R$ mithilfe zweier Folgen von positiv hermiteschen Matrizen eine linksseitige $\alpha$-Stieltjes-positiv definite Folge konstruieren, sodass jene zwei Folgen die linksseitige $\alpha$-Dyukarev-Stieltjes-Parametrisierung der konstruierten Folge bilden.

\begin{satz}	\thlabel{adplbm4}
	Seien $\alpha \in \R$, $\kappa \in \Na$ sowie $(\Laln)^{\fklam{\kappa-1}}_{n=0}$ und $(\Maln)^{\fklam{\kappa}}_{n=0}$ Folgen von positiv hermiteschen Matrizen aus $\Cqq$. Weiterhin seien durch rekursive Konstruktion \linebreak $s_0 := \Malo^{-1}$, $s_1 := \alpha s_0 - (\Malo\Lalo)^{-\ast}\Lalo(\Malo\Lalo)^{-1}$, im Fall $\kappa\geq2$
	\begin{align*}
		s_{2n} := \rklam{\prodr^{n-1}_{j=0}\Malj\Lalj}^{-\ast} \Maln^{-1}
    \rklam{\prodr^{n-1}_{j=0}\Malj\Lalj}^{-1}+z_{n,2n-1}\Hnm^{+}y_{n,2n-1}
	\end{align*}
	für alle $n\in \Zefk$ und im Fall $\kappa\geq3$
	\begin{align*}
		s_{2n+1} := \alpha s_{2n} - \rklam{\prodr^{n}_{j=0}\Malj\Lalj}^{-\ast} \Laln
    \rklam{\prodr^{n}_{j=0}\Malj\Lalj}^{-1}-z_{\aln,2n-1}\Halnm^{+}y_{\aln,2n-1}
	\end{align*}
	für alle $n\in\Zefkm$. Dann gilt $\sjk\in\Lpqka$ und $\LMkaln$ ist die linksseitige $\alpha$-Dyukarev-Stieltjes-Parametrisierung von $\sjk$.
\end{satz}

\bwanf Sei $\Qaljk$ die linksseitige $\alpha$-Stieltjes-Parametrisierung von $\sjk$. Wegen Teil (b) von \thref{aspbm5} gelten dann $s_0 = Q_{\alo}$, $s_1 = \alpha s_0 - Q_{\al1}$, im Fall $\kappa\geq2$ 
\begin{align*}
	s_{2n} = Q_{\al2n}+z_{n,2n-1}\Hnm^{+}y_{n,2n-1}
\end{align*}
für alle $n\in\Zefk$ und im Fall $\kappa\geq3$
\begin{align*}
	s_{2n+1} = \alpha s_{2n} - Q_{\al2n+1}-z_{\aln,2n-1}\Halnm^{+}y_{\aln,2n-1}
\end{align*}
für alle $n\in\Zefkm$. Hieraus folgt aus der Definition der Folge $\sjk$ dann \linebreak $Q_{\alo} = \Malo^{-1}$,
\begin{align}	\label{adplbm4bw1}
  	Q_{\al 2n+1} = \rklam{\prodr^{n}_{j=0}\Malj\Lalj}^{-\ast} \Laln
    \rklam{\prodr^{n}_{j=0}\Malj\Lalj}^{-1}
\end{align}
für alle $n\in \Zofkm$ und im Fall $\kappa\geq2$  
\begin{align}	\label{adplbm4bw2}
    Q_{\al 2n} = \rklam{\prodr^{n-1}_{j=0}\Malj\Lalj}^{-\ast} \Maln^{-1}
    \rklam{\prodr^{n-1}_{j=0}\Malj\Lalj}^{-1}
\end{align}
für alle $n\in \Zefk$. Hieraus folgt aus der Tatsache, dass $(\Laln)^{\fklam{\kappa-1}}_{n=0}$ und $(\Maln)^{\fklam{\kappa}}_{n=0}$ Folgen von positiv hermiteschen Matrizen aus $\Cqq$ sind, dann, dass auch $\Qaljk$ eine Folge von positiv hermiteschen Matrizen aus $\Cqq$ ist. Somit gilt wegen Teil (c) von \thref{aspsa2} dann $\sjk\in\Lpqka$. Wegen $Q_{\alo} = \Malo^{-1}$, \fref{adplbm4bw1}, \fref{adplbm4bw2} und \thref{amlm3} gelten $\Malo = Q^{-1}_{\al 0}$, 
\begin{align*}
    \Laln = \rklam{\prodr^{n}_{j=0}Q_{\al 2j}Q^{-1}_{\al 2j+1}}Q_{\al 2n+1}
    \rklam{\prodr^{n}_{j=0}Q_{\al 2j}Q^{-1}_{\al 2j+1}}^{\ast}
\end{align*}
für alle $n \in \Zofkm$ und im Fall $\kappa\geq2$
\begin{align*}
    \Maln = \rklam{\prodr^{n-1}_{j=0}Q^{-1}_{\al 2j}Q_{\al 2j+1}} Q^{-1}_{\al 2n}
    \rklam{\prodr^{n-1}_{j=0}Q^{-1}_{\al 2j}Q_{\al 2j+1}}^{\ast}
\end{align*}
für alle $n \in \Zefk$. Hieraus folgt wegen \thref{adplsa1} dann, dass $\LMkaln$ die linksseitige $\alpha$-Dyukarev-Stieltjes-Parametrisierung von $\sjk$ ist. \bwend

%% file: sm3.tex
\newpage
\section[Konstruktion einer Resolventenmatrix für vollständig nichtdegenerierte matrizielle \texorpdfstring{$\alpha$}{a}-Stieltjes Momentenprobleme]{Konstruktion einer Resolventenmatrix für voll-\\ständig nichtdegenerierte matrizielle \texorpdfstring{$\alpha$}{a}-Stieltjes \\ Momentenprobleme} \label{chapdr}

In den 1980er Jahren erfolgte in der Schule von V.\,P. Potapov ein intensives Studium von Matrixversionen klassischer Interpolations- und Momentenprobleme (siehe z.\,B. I.\,V. Kovalishina \cite{Ko2}, V.\,K. Dubovoj \cite{Du1}, \cite{Du2} und V.\,E. Katsnelson \cite{Kat1}, \cite{Ka1}, \cite{Kat2}, \cite{Kat3}). Aus diesen Arbeiten wurde ein wichtiges gemeinsames Merkmal dieser Aufgabenstellungen deutlich. Dieses besteht darin, dass sich im sogenannten vollständig nichtdegenerierten Fall die jeweilige Lösungsmenge durch eine gebrochen lineare Transformation von Matrizen parametrisieren lässt. Die erzeugende Matrixfunktion dieser gebrochen linearen Transformationen wird hierbei aus den Ausgangsdaten des ursprünglichen Problems konstruiert. Als Parametermenge fungiert eine von der jeweiligen Aufgabenstellung abhängige Klasse von in einem gewissen Gebiet der komplexen Ebene meromorphen Matrixfunktionen oder auch geordneten Paaren von meromorphen Matrixfunktionen. 

Im Hintergrund einer jeden in der Schule von V.\,P. Potapov betrachteten Aufgabe steht eine spezielle Signaturmatrix (vergleiche \thref{spdef1}). Diejenigen Matrixfunktionen, welche die gebrochen lineare Transformation, die die Lösungsmenge parametrisiert, erzeugen, werden auch als Resolventenmatrizen des Problems bezeichnet. Diese sind auf besondere Weise mit der zugrundeliegenden Signaturmatrix verknüpft. Dies trifft auch auf die als Parametermenge fungierende Klasse von meromorphen Matrixfunktionen bzw. geordneten Paaren von meromorphen Matrixfunktionen zu.

In diesem Kapitel werden wir nun den entsprechenden Apparat für die hier behandelten matriziellen Momentenprobleme vom $\alpha$-Stieltjes-Typ bereitstellen. In unserem Fall werden wir es mit der $2$\textit{q}$\times2$\textit{q}-Signaturmatrix
\begin{align*}
	\tJq := \begin{pmatrix} \Oq & -i\Iq \\ i\Iq & \Oq \end{pmatrix}
\end{align*}
zu tun bekommen (vergleiche \thref{spbsp1}). Im Anhang \ref{chapsp} stellen wir die für die Parametrisierung der Lösungsmenge der via Stieltjes-Transformation äquivalent umformulierten Momentenprobleme benötigten Klassen von Paaren meromorpher Matrixfunktionen, den Stieltjes-Paaren (vergleiche \thref{spdef4} und \thref{spdef5}), bereit.
Vor dem Hintergrund der vorangehenden Ausführungen prägen wir nun folgende Begriffsbildung, welche in diesem Kapitel eine zentrale Rolle einnehmen wird.

\begin{defi}	\thlabel{drdef0}
	Seien $m\in\N$, $\alpha\in\R$ und $U:\C\rightarrow\C^{2q\times2q}$ ein $2$\textit{q}$\times2$\textit{q}-Matrixpolynom. Weiterhin bezeichne
	\begin{align*}
		U = \begin{pmatrix} U^{(1,1)} & U^{(1,2)} \\ U^{(2,1)} & U^{(2,2)} \end{pmatrix}
	\end{align*}
	die \textit{q}$\times$\textit{q}-Blockzerlegung von $U$.
	\begin{itemize}
		\item [\rm{(a)}] Sei $\sjm\in\Kpqma$. Dann heißt $U$ eine \textbf{Resolventenmatrix des Momentenproblems} $\Masmu$, falls
		\begin{align*}
			\rklam{U^{(1,1)}\phi+U^{(1,2)}\psi}\rklam{U^{(2,1)}\phi+U^{(2,2)}\psi}^{-1}
		\end{align*}
		eine Bijektion zwischen der Menge der Äquivalenzklassen von \textit{q}$\times$\textit{q}"=Stieltjes"=Paaren $\phipsi$ in $\C\setminus[\alpha,\infty)$ und der Menge $\Sqasmu$ erzeugt.		
		\item [\rm{(b)}] Sei $\sjm\in\Lpqma$. Dann heißt $U$ eine \textbf{Resolventenmatrix des Momentenproblems} $\Mmasmu$ (im Fall, dass $m$ gerade ist) bzw. \linebreak $\Mmasmuu$ (im Fall, dass $m$ ungerade ist), falls
		\begin{align*}
			\rklam{U^{(1,1)}\phi+U^{(1,2)}\psi}\rklam{U^{(2,1)}\phi+U^{(2,2)}\psi}^{-1}
		\end{align*}
		eine Bijektion zwischen der Menge der Äquivalenzklassen von \textit{q}$\times$\textit{q}"=Stieltjes"=Paaren $\phipsi$ in $\C\setminus(-\infty,\alpha]$ und der Menge $\Sqmasmu$ erzeugt.
	\end{itemize}
\end{defi}

Das Hauptziel dieses Abschnitts besteht nun in der Konstruktion von Resolventenmatrizen für die beiden hier betrachteten nichtdegenerierten matriziellen $\alpha$-Stieltjes Momentenprobleme. Hierbei wenden wir uns zunächst dem rechtsseitigen Fall zu. Am Ausgangspunkt unserer Betrachtungen stehen die Untersuchungen von Yu.\,M. Dyukarev \cite{dyu} im Fall $\alpha=0$. Er behandelte das Problem mithilfe der Methode der fundamentalen Matrixungleichungen von V.\,P. Potapov und konstruierte eine konkrete Resolventenmatrix (siehe \cite[Theorem 2]{dyu}), welche ein $2$\textit{q}$\times2$\textit{q}-Matrixpolynom mit speziellen Eigenschaften bezüglich der $2$\textit{q}$\times2$\textit{q}-Signaturmatrix $\tJq$ ist. Eine weitergehende Analyse der Struktur der Dyukarevschen Resolventenmatrix erfolgte in der kürzlich erschienenen Arbeit von A.\,E. Choque Rivero \cite{CR1}.

In ihrer Dissertation \cite{Maka} zeigte T. Makarevich bereits, dass via Stieltjes-Transformation die Lösungsmenge des ursprünglichen Momentenproblems mit der Lösungsmenge des Systems der beiden Potapovschen fundamentalen Matrixungleichungen im Fall einer vorgegebenen Momentenfolge aus ${\cal K}^{\geq,e}_{q,2n+1,\alpha}$ für $n\in\No$ übereinstimmt, und formulierte eine entsprechende Resolventenmatrix. 

Wir wenden uns hier dem sogenannten vollständig nichtdegenerierten Fall zu, der genau dann vorliegt, wenn die Folge der vorgegebenen Momente $\alpha$-Stieltjes-positiv definit ist. In diesem Fall kann jede der beiden Potapovschen fundamentalen Matrixungleichungen mithilfe einer Modifikation der in der Schule von V.\,P. Potapov ausgearbeiteten Faktorisierungsmethode gelöst werden. Das Problem besteht nun darin, eine geeignete Kopplung zwischen den beiden Potapovschen fundamentalen Matrixungleichungen herzustellen. Eine ähnliche Situation lag in den Untersuchungen von A.\,E. Choque Rivero, Yu.\,M. Dyukarev, B. Fritzsche und B. Kirstein \cite{C06} zum finiten matriziellen Hausdorffschen Momentenproblem vor. Der dort entwickelte Formalismus zur Herstellung einer Kopplung zwischen beiden fundamentalen Matrixungleichungen lieferte den Ausgangspunkt für unsere Vorgehensweise. Zur Behandlung der hier vorliegenden Situation werden wir eine entsprechende Modifikation der Konstruktion in \cite[Chapter 6]{C06} vornehmen.

Zudem führen wir eine weitere Begriffsbildung ein, das sogenannte $\alpha$-Dyukarev"=Quadrupel, welches die vier \textit{q}$\times$\textit{q}-Einträge der von uns betrachteten Resolventenmatrix umfasst. Dieses Quadrupel von Folgen von \textit{q}$\times$\textit{q}-Matrixpolynomen wurde für den Fall $\sj\in{\cal K}^{>}_{q,\infty,0}$ schon in \cite[Chapter 3]{CR1} eingeführt.

Weiterhin widmen wir uns zwei ausgezeichneten Elementen der Lösungsmenge des via Stieltjes"=Transformation umformulierten vollständig nichtdegenerierten $\alpha$-Stieltjes Momentenproblems. Es wird sich herausstellen, dass diese rationalen \textit{q}$\times$\textit{q}"=Matrixfunktionen eine gewisse extremale Stellung einnehmen. Wir werden diese beiden Funktionen explizit darstellen können und zeigen, dass ihre zugehörigen Stieltjes-Maße molekular sind.

Für die Beweise im linksseitigen Fall werden wir auf die Resultate des rechtsseitigen Falles zurückgreifen.

\subsection{Der rechtsseitige Fall} \label{chapdrr}

Wir knüpfen nun an die Ausführungen von Abschnitt \ref{chapFMR} an. Haben wir dort noch allgemeine rechtsseitige $\alpha$-Stieltjes Momentenprobleme betrachtet, wollen wir uns nun auf den vollständig nichtdegenerierten Fall beschränken, das heißt unsere gegebene Folge ist rechtsseitig $\alpha$-Stieltjes-positiv definit anstatt nur -nichtnegativ definit.

\begin{bem}	\thlabel{drbm3}
  Seien $\kappa \in \Na$, $\alpha \in \R$ und $\sjk \in \Kpqka$. Weiterhin sei $S:$ \linebreak $\C\setminus[\alpha,\infty) \rightarrow \Cqq$.\dgfa
  \begin{itemize}
    \item [\rm{(a)}] Seien $z \in \C\setminus\R$ und $n \in \Zofk$. \dsfaa
    \begin{itemize}
      \item [\rm{(i)}] Es gilt $\FSn(z) \in \C^{(n+2)q\times(n+2)q}_{\geq}$.
      \item [\rm{(ii)}] Es gilt $\dFSn(z) \in \Cqq_{\geq}$.
    \end{itemize}
    \item [\rm{(b)}] Seien $z \in \C\setminus\R$ und $n \in \Zofkm$. \dsfaa
    \begin{itemize}
      \item [\rm{(iii)}] Es gilt $\FSarn(z) \in \C^{(n+2)q\times(n+2)q}_{\geq}$.
      \item [\rm{(iv)}] Es gilt $\dFSarn(z) \in \Cqq_{\geq}$.
    \end{itemize}
  \end{itemize}
\end{bem}

\bwanf Zu (a): Wegen \thref{adpbm1} gelten $\Hn \in \C^{(n+1)q\times(n+1)q}_{\geq}$ und $\Hn^{+} = \Hn^{-1}$.
Hieraus folgt mithilfe von \thref{drdef3} und \thref{amlm2} dann die Äquivalenz von (i) und (ii).

Zu (b): Wegen \thref{adpbm1} gelten $\Harn \in \C^{(n+1)q\times(n+1)q}_{\geq}$ und $\Harn^{+} = \Harn^{-1}$.
Hieraus folgt mithilfe von \thref{drdef3} und \thref{amlm2} dann die Äquivalenz von (iii) und (iv). \bwend

Wir stellen nun einige für die in \thref{drbz1} und \thref{drbz3} eingeführten Matrizen gültige Identitäten auf, die für unsere weitere Vorgehensweise nützlich werden. Einige der folgenden Identitäten oder Identitäten in ähnlicher Form wurden in \cite{dyu} für den Fall $\alpha=0$ auch ohne Beweis angegeben. Man findet in \cite[Kapitel 2.2 und Kapitel 5]{Sch} oder \cite[Kapitel 2]{MP} weitere Ausführungen zu jenen Matrizen. Wir verwenden hier aber speziell in \thref{drbm1} stark abgewandelte Identitäten, die unserer später folgenden Beweisform angepasst sind.

\begin{lemma}	\thlabel{drlm1}
  \egfa
  \begin{itemize}
    \item [\rm{(a)}] Für alle $z \in \C$ und $n \in \No$ gilt
      \begin{align*}
	R_n(z)T_n = T_nR_n(z)
      \end{align*}
    \item [\rm{(b)}] Für alle $z, \omega \in \C$ und $n \in \No$ gilt
      \begin{align*}
	R_n(z)R_n(\omega) = R_n(\omega)R_n(z)
      \end{align*}
    \item [\rm{(c)}] Für alle $z, \omega \in \C$ und $n \in \No$ gilt
      \begin{align*}
	R_n(z)R^{-1}_n(\omega) = R^{-1}_n(\omega)R_n(z)
      \end{align*}
    \item [\rm{(d)}] Für alle $z, \omega \in \C$ und $n \in \No$ gilt
      \begin{align*}
	R_n(z)-R_n(\omega) = (z-\omega)R_n(z)T_nR_n(\omega)
      \end{align*}
  \end{itemize}
\end{lemma}

\bwanf Die Identitäten sind leicht selbst nachzuweisen. Einen detaillierten Beweis findet man z.\,B. unter \cite[Lemma 5.12]{Sch}. \bwend

\begin{lemma}	\thlabel{drlm2}
  \egfa
  \begin{itemize}
    \item [\rm{(a)}] Für alle $n \in \No$ gelten
      \begin{align*}
	T_n = L_n\dL^{\ast}_n, \quad T_n\dL_n = L_n \quad \text{und} \quad T^{\ast}_nL_n = \dL_n.
      \end{align*}
    \item [\rm{(b)}] Für alle $n \in \N$ gelten
      \begin{align*}
	\dL_nv_{n-1} = v_n \quad \text{und} \quad \dL^{\ast}_nv_n = v_{n-1}.
      \end{align*}
    \item [\rm{(c)}] Für alle $n \in \N$ gelten
      \begin{align*}
	L_nT_{n-1} = T_nL_n \quad \text{und} \quad \dL^{\ast}_nT_n = T_{n-1}\dL^{\ast}_n
      \end{align*}
    \item [\rm{(d)}] Für alle $z \in \C$ und $n \in \N$ gelten
      \begin{align*}
	L_nR_{n-1}(z) = R_n(z)L_n \quad \text{und} \quad \dL^{\ast}_nR_n(z) = R_{n-1}(z)\dL^{\ast}_n.
      \end{align*}
  \end{itemize}
\end{lemma}

\bwanf Zu (a)-(b): Dies folgt direkt aus der Definition der beteiligten Größen (vergleiche \thref{drbz1}).

Zu (c): Es gelten
\begin{align*}
  L_nT_{n-1} = \begin{pmatrix} 0_{2q \times (n-1)q} & 0_{2q \times q} \\ I_{(n-1)q} & 0_{(n-1)q \times q} \end{pmatrix} = T_nL_n
\end{align*}
und
\begin{align*}
  \dL^{\ast}_nT_n = \begin{pmatrix} 0_{q \times (n-1)q} & 0_{q \times 2q} \\ I_{(n-1)q} & 0_{(n-1)q \times 2q} \end{pmatrix} = T_{n-1}\dL^{\ast}_n.
\end{align*}

Zu (d): Es gelten
\begin{align*}
  L_nR_{n-1}(z) = \begin{pmatrix} \Oqn \\ R_{n-1}(z) \end{pmatrix} = R_n(z)L_n
\end{align*}
und
\begin{align*}
  \dL^{\ast}_nR_n(z) = \begin{pmatrix} R_{n-1}(z) & \Onq \end{pmatrix} = R_{n-1}(z)\dL^{\ast}_n. \tag*{$\Box$}
\end{align*}

\begin{lemma}	\thlabel{drlm3}
  Seien $\kappa \in \Na$, $\alpha \in \R$ und $\sjk$ eine Folge aus $\Cqq$. \dgfa
  \begin{itemize}
    \item [\rm{(a)}] Für alle $n \in \Zok$ gilt
      \begin{align*}
	T_ny_{0,n} = u_n.
      \end{align*}
    \item [\rm{(b)}] Für alle $n \in \Zokp$ gilt
      \begin{align*}
	L_ny_{0,n-1} = u_n
      \end{align*}
      und für alle $n \in \Zok$ gilt
      \begin{align*}
	\dL^{\ast}_ny_{0,n}=y_{0,n-1}.
      \end{align*}
    \item [\rm{(c)}] Für alle $n \in \Zokp$ gilt
      \begin{align*}
	\dL^{\ast}_nu_n = u_{n-1}.
      \end{align*}
    \item [\rm{(d)}] Für alle $n \in \Zek$ gilt
      \begin{align*}
	\dL^{\ast}_nu_{\arn} = u_{\arn-1}.
      \end{align*}
    \item [\rm{(e)}] Für alle $n \in \Zok$ gilt
      \begin{align*}
	R_n(\alpha)u_{\arn} = y_{0,n}.
      \end{align*}
  \end{itemize}
\end{lemma}

\bwanf Zu (a)-(d): Dies folgt direkt aus der Definition der beteiligten Größen (vergleiche \thref{drbz1}, \thref{asmbz1} und \thref{drbz3}).

Zu (e): Aus der Definition von $u_{\arn}$ und $\Rn(\alpha)$ (vergleiche \thref{drbz3} und \thref{drbz1}) sowie (a) folgt
\begin{align*}
  u_{\arn} = y_{0,n} - \alpha u_n = R^{-1}_n(\alpha)y_{0,n}.
\end{align*}
Hieraus folgt dann die Behauptung. \bwend

\begin{bem}[Kopplungsidentitäten]	\thlabel{drbm1}
  Seien $\kappa \in \Na$, $\alpha \in \R$ und $\sjk$ eine Folge aus $\Cqq$. \dgfa
  \begin{itemize}
    \item [\rm{(a)}] Für alle $n \in \Zofkm$ gilt
      \begin{align*}
	v_nz_{0,n} = H_n - T_nK_n.
      \end{align*}
    \item [\rm{(b)}] Für alle $n \in \Zefk$ gilt
      \begin{align*}
	v_nz_{0,n-1} = H_n\dL_n - L_nK_{n-1}.
      \end{align*}
    \item [\rm{(c)}] Für alle $n \in \Zofkm$ gilt
      \begin{align*}
	v_nz_{0,n} = R^{-1}_n(\alpha)H_n - T_n\Harn.
      \end{align*}
    \item [\rm{(d)}] Für alle $n \in \Zefk$ gilt
      \begin{align*}
	v_nz_{0,n-1} = R^{-1}_n(\alpha)H_n\dL_n - L_n\Harnm.
      \end{align*}
    \item [\rm{(e)}] Für alle $n \in \Zofk$ gilt
      \begin{align*}
	v_nz_{0,n} = R^{-1}_n(\alpha)H_n - T_n\tHarn,
      \end{align*}
      wobei
      \begin{align*}
	\tHarn := \begin{cases} \Harn	& \text{falls } n\in\Zofkm \\ 
	\begin{pmatrix} \Harnm & y_{\arn,2n-1} \\ z_{\arn,2n-1} & \Oq \end{pmatrix} & \text{falls } 2n=\kappa. \end{cases}
      \end{align*}
    \item [\rm{(f)}] Für alle $n \in \Zefkp$ gilt
      \begin{align*}
	v_nz_{0,n-1} = R^{-1}_n(\alpha)\tHn\dL_n - L_n\Harnm,
      \end{align*}
      wobei
      \begin{align*}
	\tHn := \begin{cases} \Hn	& \text{falls } n\in\Zofk \\ 
	\begin{pmatrix} \Hnm & y_{n,2n-1} \\ z_{n,2n-1} & \Oq \end{pmatrix} & \text{falls } 2n-1=\kappa. \end{cases}
      \end{align*}
  \end{itemize}
\end{bem}

\bwanf Teil (c) findet man in ähnlicher Form auch unter \cite[Lemma 5.17]{Sch}.

Zu (a): Es gilt $H_0-T_0K_0=s_0=v_0z_{0,0}$ und für $n \in \Zefkm$ gilt weiterhin
\begin{align*}
  H_n - T_nK_n = \begin{pmatrix} s_0 & \ldots & \ s_n \\ \vdots & & \vdots \\ s_n & \ldots & s_{2n} \end{pmatrix}
  - \begin{pmatrix} \Oq & \ldots & \Oq \\ s_1 & \ldots & \ s_{n+1} \\ \vdots & & \vdots \\ s_n & \ldots & s_{2n} \end{pmatrix}
  = \begin{pmatrix} s_0 \ldots \ s_n \\ 0_{nq \times (n+1)q} \end{pmatrix} = v_nz_{0,n}.
\end{align*}

Zu (b): Es gilt
\begin{align*}
  H_n\dL_n - L_nK_{n-1} & = \begin{pmatrix} s_0 & \ldots & \ s_{n-1} \\ \vdots & & \vdots \\ s_n & \ldots & s_{2n-1} \end{pmatrix}
  - \begin{pmatrix} \Oq & \ldots & \Oq \\ s_1 & \ldots & \ s_n \\ \vdots & & \vdots \\ s_n & \ldots & s_{2n-1} \end{pmatrix} \\
  & = \begin{pmatrix} s_0 \ldots \ s_{n-1} \\ 0_{nq \times nq} \end{pmatrix} = v_nz_{0,n-1}.
\end{align*}

Zu (c): Wegen der Definition von $\Harn$ gilt
\begin{align*}
  \Harn = -\alpha\Hn + \Kn.
\end{align*}
Hieraus folgt wegen (a) dann
\begin{align*}
  v_nz_{0,n} = H_n - T_n(\Harn + \alpha H_n) = R^{-1}_n(\alpha)H_n - T_n\Harn.
\end{align*}

Zu (d): Wegen der Definition von $\Harnm$ gilt
\begin{align}	\label{drbm1bw1}
  \Harnm = -\alpha\Hnm + \Knm.
\end{align}
Weiterhin gilt
\begin{align*}
  T_nH_n\dL_n = \begin{pmatrix} \Oqn & \Oq \\ \Hnm & y_{n,2n-1} \end{pmatrix} \begin{pmatrix} \Inq \\ \Oqn \end{pmatrix}
  = \begin{pmatrix} \Oqn \\ \Hnm \end{pmatrix} = L_n\Hnm.
\end{align*}
Hieraus folgt wegen \fref{drbm1bw1} und (b) dann
\begin{align*}
  v_nz_{0,n-1} = H_n\dL_n - L_n(\Harnm + \alpha \Hnm) = R^{-1}_n(\alpha)H_n\dL_n - L_n\Harnm.
\end{align*}

Zu (e): Dies folgt aus (c). Man beachte hier, dass $\Tn\Harn$ unabhängig von $s_{\ar2n}$ ist.

Zu (f): Dies folgt aus (d). Man beachte hier, dass $\Hn\dLn$ unabhängig von $s_{2n}$ ist. \bwend

Unsere nächsten Betrachtungen führen uns vor dem Hintergrund von \thref{drbm3} von den Potapovschen Fundamentalmatrizen auf eine Resolventenmatrix für das rechtsseitige $\alpha$-Stieltjes Momentenproblem. Hierfür behandeln wir zunächst einige Matrixfunktionen, die sich als spezielle $\tJq$-Potapov-Funktionen herausstellen werden (vergleiche \thref{spdef3} und \thref{spbsp1}). Sie werden uns helfen, die linken Schur-Komplemente der Potapovschen Fundamentalmatrizen jeweils in einer für unsere Zwecke besonders geeigneten Form darzustellen. Wir verwenden hier eine ähnliche Vorgehensweise wie das in \cite[Chapter 6]{C06} für das Hausdorffsche Momentenproblem vorgenommen wurde, allerdings mit auf unser Problem zugeschnittenen Matrixfunktionen. Es sei bemerkt, dass diese Vorgehensweise sich mit der von \cite[Chapter 8]{Maka} für den allgemeineren Fall einer gegebenen Folge aus $\Keqka$ unterscheidet. Die dort behandelten Resultate sind speziell auf das finite rechtsseitige $\alpha$-Stieltjes Momentenproblem mit einer geraden Anzahl von vorgegebenen Momenten angepasst. Da wir hier im vollständig nichtdegenerten Fall auch eine ungerade Anzahl von vorgegebenen Momenten behandeln wollen, geht unsere Vorgehensweise in eine etwas andere Richtung. Wir werden aber nach jedem Resultat, falls möglich, Zusammenhänge zu \cite{Maka} erläutern.

\begin{bez}	\thlabel{drbz4}
	Seien $\kappa \in \Na$, $\alpha \in \R$ und $\sjk \in \Kpqka$. Für $n \in \Zofk$ sei $\Varn: \C \rightarrow \C^{2q\times 2q}$ definiert gemäß
	\begin{align*}
		\Varn^{\sklam{s}}(z) := \Izq + (z-\alpha)\begin{pmatrix} \big(\un^{\sklam{s}}\big)^{\ast} \\ -\vna \end{pmatrix}\Rna(\za)\brklam{\Hsn}^{-1}\Rn(\alpha)\begin{pmatrix} \vn & \un^{\sklam{s}}\end{pmatrix}.
	\end{align*}
	Für $n \in \Zofkm$ sei weiterhin $\Varnp: \C \rightarrow \C^{2q\times 2q}$ definiert gemäß
	\begin{align*}
		\Varnp^{\sklam{s}}(z) := \Izq + (z-\alpha)\begin{pmatrix} \big(\uarn^{\sklam{s}}\big)^{\ast} \\ -\vna \end{pmatrix}\Rna(\za)\brklam{\Hsarn}^{-1}\Rn(\alpha)\begin{pmatrix} \vn & \uarn^{\sklam{s}} \end{pmatrix}.
	\end{align*}
	Falls klar ist, von welchem $\sjk$ die Rede ist, lassen wir das \anf{$\sklam{s}$} als oberen Index weg.
\end{bez}

\begin{satz}	\thlabel{drsa2}
  Seien $\kappa \in \Na$, $\alpha \in \R$ und $\sjk \in \Kpqka$. 
  \dgfa
  \begin{itemize}
    \item [\rm{(a)}] Es gelten
    \begin{align*}
      \tJq-\Varn(z)\tJq\Varn^{\ast}(\omega) = -i(z-\overline{\omega})\begin{pmatrix} u^{\ast}_n \\ -v^{\ast}_n \end{pmatrix}\Rna(\za)\Hn^{-1}\Rn(\overline{\omega})\begin{pmatrix}u_n & -\vn\end{pmatrix}
    \end{align*}
    für alle $z, \omega \in \C$ und $n \in \Zofk$ sowie
    \begin{align*}
      \tJq-\Varnp(z)\tJq\Varnp^{\ast}(\omega) = -i(z-\overline{\omega})\begin{pmatrix} \uarn^{\ast} \\ -v^{\ast}_n \end{pmatrix}\Rna(\za)\Harn^{-1}\Rn(\overline{\omega})\begin{pmatrix}\uarn & -\vn\end{pmatrix}
    \end{align*}
    für alle $z, \omega \in \C$ und $n \in \Zofkm$.
    \item [\rm{(b)}] Es gilt $\Varm \in \tbtJqPp$ für alle $m \in \Zok$.
    \item [\rm{(c)}] Für alle $z \in \C$ und $n \in \Zofk$ ist $\Varn(z)$ regulär und es gilt
    \begin{align*}
      \Varn^{-1}(z) = \Izq - (z-\alpha)\begin{pmatrix} u^{\ast}_n \\ -v^{\ast}_n \end{pmatrix}\Rna(\alpha)\Hn^{-1}\Rn(z)\begin{pmatrix}v_n & u_n\end{pmatrix}.
    \end{align*}
    Für alle $z \in \C$ und $n \in \Zofkm$ ist $\Varnp(z)$ regulär und es gilt
    \begin{align*}
      \Varnp^{-1}(z) = \Izq - (z-\alpha)\begin{pmatrix} \uarn^{\ast} \\ -v^{\ast}_n \end{pmatrix}\Rna(\alpha)\Harn^{-1}\Rn(z)\begin{pmatrix}v_n & \uarn\end{pmatrix}.
    \end{align*}
    \item [\rm{(d)}] Unter Beachtung von (c) gelten
    \begin{align*}
      \tJq-\Varn^{-\ast}(z)\tJq\Varn^{-1}(\omega) = i(\omega-\za)\tJq\begin{pmatrix} u^{\ast}_n \\ -v^{\ast}_n \end{pmatrix}\Rna(z)\Hn^{-1}\Rn(\omega)\begin{pmatrix}u_n & -\vn\end{pmatrix}\tJq
    \end{align*}
    für alle $z, \omega \in \C$ und $n \in \Zofk$ sowie
    \begin{align*}
      &\ \tJq-\Varnp^{-\ast}(z)\tJq\Varnp^{-1}(\omega) \\
      &= i(\omega-\za)\tJq\begin{pmatrix} \uarn^{\ast} \\ -v^{\ast}_n \end{pmatrix}\Rna(z)\Harn^{-1}\Rn(\omega)\begin{pmatrix}\uarn & -\vn\end{pmatrix}\tJq
    \end{align*}
    für alle $z, \omega \in \C$ und $n \in \Zofkm$.
    \item [\rm{(e)}] Seien $f: \C\setminus\R \rightarrow \Cqq$ und $f_{\ar}(z) := (z-\alpha)f(z)$ für alle $z \in \C\setminus\R$. 
    Unter Beachtung von (c) gelten dann
    \begin{align*}
      \dFfn(z) = \frac{1}{i(z-\za)}\binom{f(z)}{\Iq}^{\ast}\binom{f(z)}{\Iq}^{\ast}\Varn^{-\ast}(z)\tJq\Varn^{-1}(z)\binom{f(z)}{\Iq}
    \end{align*}
    für alle $z \in \C\setminus\R$ und $n \in \Zofk$ sowie
    \begin{align*}
      \dFfarn(z) = \frac{1}{i(z-\za)}\binom{f_{\ar}(z)}{\Iq}^{\ast}\Varnp^{-\ast}(z)\tJq\Varnp^{-1}(z)\binom{f_{\ar}(z)}{\Iq}
    \end{align*}
    für alle $z \in \C\setminus\R$ und $n \in \Zofkm$.
  \end{itemize}
\end{satz}

\bwanf Zu (a): Seien $z, \omega \in \C$ und zunächst $n \in \Zofk$. Dann gilt
\begin{align}	\label{drsa2bw1}
  \begin{pmatrix} \vn & \un \end{pmatrix} \tJq \begin{pmatrix} \vna \\ \una \end{pmatrix} 
  = \begin{pmatrix} i\un & -i\vn \end{pmatrix} \begin{pmatrix} \vna \\ \una \end{pmatrix} = i(\un\vna-\vn\una).
\end{align}
Weiterhin gilt
\begin{align}	\label{drsa2bw2}
  \Hn\Tna-\Tn\Hn & = \begin{pmatrix} \Oq & s_0 & \ldots & s_{n-1} \\ \vdots & \vdots & & \vdots \\ \Oq & s_n & \ldots & s_{2n-1} \end{pmatrix}
  - \begin{pmatrix} \Oq & \ldots & \Oq \\ s_0 & \ldots & s_n \\ \vdots & & \vdots \\ s_{n-1} & \ldots & s_{2n-1} \end{pmatrix} \notag \\
  & = \begin{pmatrix} \Oq & \ynma \\ -\ynm & \Onqn \end{pmatrix} \notag \\
  & = -\begin{pmatrix} \Oq & \Oqn \\ \ynm & \Onqn \end{pmatrix} + \begin{pmatrix} \Oq & \ynma \\ \Oqn & \Onqn \end{pmatrix} \notag \\
  & = -(\un\vna-\vn\una).
\end{align}
Es gilt
\begin{align}	\label{drsa2bw3}
  \Rn^{-1}(\alpha)\Hn\Tna-\Tn\Hn\Rn^{-\ast}(\alpha) & = (\Inpq-\alpha\Tn)\Hn\Tna-\Tn\Hn(\Inpq-\alpha\Tna) \notag \\
  & = \Hn\Tna-\Tn\Hn.
\end{align}
Weiterhin gelten
\begin{align*}
  (\omegaa-\alpha)\Tn = (\Inpq - \alpha\Tn) - (\Inpq - \omegaa\Tn) = \Rn^{-1}(\alpha) - \Rn^{-1}(\omegaa)
\end{align*}
und
\begin{align*}
  (z-\alpha)\Tna = (\Inpq - \alpha\Tn)^{\ast} - (\Inpq - \za\Tn)^{\ast} = \Rn^{-\ast}(\alpha) - \Rn^{-\ast}(\za).
\end{align*}
Hieraus folgen dann
\begin{align}	\label{drsa2bw4}
  (\omegaa-\alpha)\Rn(\alpha)\Tn = \Inpq-\Rn(\alpha)\Rn^{-1}(\omegaa)
\end{align}
und
\begin{align}	\label{drsa2bw5}
  (z-\alpha)\Tna\Rna(\alpha) = \Inpq-\Rn^{-\ast}(\za)\Rna(\alpha).
\end{align}
Wegen \fref{drsa2bw1} - \fref{drsa2bw3} sowie \fref{drsa2bw4} und \fref{drsa2bw5} gilt nun
\begin{align*}
  &\ (z-\alpha)(\omegaa-\alpha)\Rn(\alpha) \begin{pmatrix} \vn & \un \end{pmatrix} \tJq \begin{pmatrix} \vna \\ \una \end{pmatrix} \Rna(\alpha) \\
  & = -i(z-\alpha)(\omegaa-\alpha)\Rn(\alpha) \eklam{\Rn^{-1}(\alpha)\Hn\Tna-\Tn\Hn\Rn^{-\ast}(\alpha)} \Rna(\alpha) \\
  & = -i(z-\alpha)(\omegaa-\alpha) \eklam{\Hn\Tna\Rna(\alpha)-\Rn(\alpha)\Tn\Hn} \\
  & = -i(\omegaa-\alpha)\Hn\eklam{\Inpq-\Rn^{-\ast}(\za)\Rna(\alpha)} + i(z-\alpha)\eklam{\Inpq-\Rn(\alpha)\Rn^{-1}(\omegaa)}\Hn \\
  & = i\eklam{(z-\omegaa)\Hn + (\omegaa-\alpha)\Hn\Rn^{-\ast}(\za)\Rna(\alpha) - (z-\alpha)\Rn(\alpha)\Rn^{-1}(\omegaa)\Hn}.
\end{align*}
Hieraus folgt
\begin{align*}
  &\ (z-\alpha)(\omegaa-\alpha)\begin{pmatrix} \una \\ -\vna \end{pmatrix}\Rna(\za)\Hn^{-1}\Rn(\alpha)\begin{pmatrix} \vn & \un \end{pmatrix}\tJq\begin{pmatrix} \vna \\ \una \end{pmatrix}\Rna(\alpha)\Hn^{-1}\Rn(\omegaa)\begin{pmatrix} \un & -\vn \end{pmatrix} \\
  & = i\begin{pmatrix} \una \\ -\vna \end{pmatrix}\Rna(\za)\Hn^{-1}\big[(z-\omegaa)\Hn + (\omegaa-\alpha)\Hn\Rn^{-\ast}(\za)\Rna(\alpha) \\
  &\quad - (z-\alpha)\Rn(\alpha)\Rn^{-1}(\omegaa)\Hn\big]\Hn^{-1}\Rn(\omegaa)\begin{pmatrix} \un & -\vn \end{pmatrix} \\
  & = i(z-\omegaa)\begin{pmatrix} \una \\ -\vna \end{pmatrix}\Rna(\za)\Hn^{-1}\Rn(\omegaa)\begin{pmatrix} \un & -\vn \end{pmatrix} \\
  &\quad +i(\omegaa-\alpha)\begin{pmatrix} \una \\ -\vna \end{pmatrix}\Rna(\alpha)\Hn^{-1}\Rn(\omegaa)\begin{pmatrix} \un & -\vn \end{pmatrix} \\
  &\quad -i(z-\alpha)\begin{pmatrix} \una \\ -\vna \end{pmatrix}\Rna(\za)\Hn^{-1}\Rn(\alpha)\begin{pmatrix} \un & -\vn \end{pmatrix}.
\end{align*}
Hieraus folgt wiederum
\begin{align*}
  &\ \tJq-\Varn(z)\tJq\Varn^{\ast}(\omega) \\
  & = \tJq-\eklam{\Izq + (z-\alpha)\begin{pmatrix} \una \\ -\vna \end{pmatrix}\Rna(\za)\Hn^{-1}\Rn(\alpha)\begin{pmatrix} \vn & \un \end{pmatrix}}\tJq \\
  &\quad \cdot \eklam{\Izq + (\omegaa-\alpha)\begin{pmatrix} \vna \\ \una \end{pmatrix}\Rna(\alpha)\Hn^{-1}\Rn(\omegaa)\begin{pmatrix} \un & -\vn \end{pmatrix}} \\
  & = -i(z-\alpha)\begin{pmatrix} \una \\ -\vna \end{pmatrix}\Rna(\za)\Hn^{-1}\Rn(\alpha)\begin{pmatrix} \un & -\vn \end{pmatrix} \\
  &\quad +i(\omegaa-\alpha)\begin{pmatrix} \una \\ -\vna \end{pmatrix}\Rna(\alpha)\Hn^{-1}\Rn(\omegaa)\begin{pmatrix} \un & -\vn \end{pmatrix} \\
  &\quad -(z-\alpha)(\omegaa-\alpha)\begin{pmatrix} \una \\ -\vna \end{pmatrix}\Rna(\za)\Hn^{-1}\Rn(\alpha)\begin{pmatrix} \vn & \un \end{pmatrix}\tJq\begin{pmatrix} \vna \\ \una \end{pmatrix} \\
  &\quad \cdot \Rna(\alpha)\Hn^{-1}\Rn(\omegaa)\begin{pmatrix} \un & -\vn \end{pmatrix} \\
  & = -i(z-\omegaa)\begin{pmatrix} \una \\ -\vna \end{pmatrix}\Rna(\za)\Hn^{-1}\Rn(\omegaa)\begin{pmatrix} \un & -\vn \end{pmatrix}.
\end{align*}

Sei nun $n \in \Zofkm$. Dann gilt
\begin{align}	\label{drsa2bw6}
  \begin{pmatrix} \vn & \uarn \end{pmatrix} \tJq \begin{pmatrix} \vna \\ \uarna \end{pmatrix} 
  = \begin{pmatrix} i\uarn & -i\vn \end{pmatrix} \begin{pmatrix} \vna \\ \uarna \end{pmatrix} = i(\uarn\vna-\vn\uarna).
\end{align}
Weiterhin gilt
\begin{align}	\label{drsa2bw7}
  \Harn\Tna-\Tn\Harn & = \begin{pmatrix} \Oq & s_{\aro} & \ldots & s_{\arn-1} \\ \vdots & \vdots & & \vdots \\ \Oq & s_{\arn} & \ldots & s_{\ar2n-1} \end{pmatrix}
  - \begin{pmatrix} \Oq & \ldots & \Oq \\ s_{\aro} & \ldots & s_{\arn} \\ \vdots & & \vdots \\ s_{\arn-1} & \ldots & s_{\ar2n-1} \end{pmatrix} \notag \\
  & = \begin{pmatrix} \Oq & \yarnma \\ -\yarnm & \Onqn \end{pmatrix} \notag \\
  & = -\begin{pmatrix} s_0 & \Oqn \\ \yarnm & \Onqn \end{pmatrix} + \begin{pmatrix} s_0 & \yarnma \\ \Oqn & \Onqn \end{pmatrix} \notag \\
  & = -(\uarn\vna-\vn\uarna).
\end{align}
Es gilt
\begin{align}	\label{drsa2bw8}
  \Rn^{-1}(\alpha)\Harn\Tna-\Tn\Harn\Rn^{-\ast}(\alpha) & = (\Inpq-\alpha\Tn)\Harn\Tna-\Tn\Harn(\Inpq-\alpha\Tna) \notag \\
  & = \Harn\Tna-\Tn\Harn.
\end{align}
Wegen \fref{drsa2bw6} - \fref{drsa2bw8} sowie \fref{drsa2bw4} und \fref{drsa2bw5} gilt nun
\begin{align*}
  &\ (z-\alpha)(\omegaa-\alpha)\Rn(\alpha) \begin{pmatrix} \vn & \uarn \end{pmatrix} \tJq \begin{pmatrix} \vna \\ \uarna \end{pmatrix} \Rna(\alpha) \\
  & = -i(z-\alpha)(\omegaa-\alpha)\Rn(\alpha) \eklam{\Rn^{-1}(\alpha)\Harn\Tna-\Tn\Harn\Rn^{-\ast}(\alpha)} \Rna(\alpha) \\
  & = -i(z-\alpha)(\omegaa-\alpha) \eklam{\Harn\Tna\Rna(\alpha)-\Rn(\alpha)\Tn\Harn} \\
  & = -i(\omegaa-\alpha)\Harn\eklam{\Inpq-\Rn^{-\ast}(\za)\Rna(\alpha)} + i(z-\alpha)\eklam{\Inpq-\Rn(\alpha)\Rn^{-1}(\omegaa)}\Harn \\
  & = i\eklam{(z-\omegaa)\Harn + (\omegaa-\alpha)\Harn\Rn^{-\ast}(\za)\Rna(\alpha) - (z-\alpha)\Rn(\alpha)\Rn^{-1}(\omegaa)\Harn}.
\end{align*}
Hieraus folgt
\begin{align*}
  &\ (z-\alpha)(\omegaa-\alpha)\begin{pmatrix} \uarna \\ -\vna \end{pmatrix}\Rna(\za)\Harn^{-1}\Rn(\alpha)\begin{pmatrix} \vn & \uarna \end{pmatrix}\tJq \\
  &\ \cdot\begin{pmatrix} \vna \\ \uarna \end{pmatrix}\Rna(\alpha)\Harn^{-1}\Rn(\omegaa)\begin{pmatrix} \uarn & -\vn \end{pmatrix} \\
  & = i\begin{pmatrix} \uarna \\ -\vna \end{pmatrix}\Rna(\za)\Harn^{-1}\big[(z-\omegaa)\Harn + (\omegaa-\alpha)\Harn\Rn^{-\ast}(\za)\Rna(\alpha) \\
  &\quad - (z-\alpha)\Rn(\alpha)\Rn^{-1}(\omegaa)\Harn\big]\Harn^{-1}\Rn(\omegaa)\begin{pmatrix} \uarn & -\vn \end{pmatrix} \\
  & = i(z-\omegaa)\begin{pmatrix} \uarna \\ -\vna \end{pmatrix}\Rna(\za)\Harn^{-1}\Rn(\omegaa)\begin{pmatrix} \uarn & -\vn \end{pmatrix} \\
  &\quad +i(\omegaa-\alpha)\begin{pmatrix} \uarna \\ -\vna \end{pmatrix}\Rna(\alpha)\Harn^{-1}\Rn(\omegaa)\begin{pmatrix} \uarn & -\vn \end{pmatrix} \\
  &\quad -i(z-\alpha)\begin{pmatrix} \uarna \\ -\vna \end{pmatrix}\Rna(\za)\Harn^{-1}\Rn(\alpha)\begin{pmatrix} \uarn & -\vn \end{pmatrix}.
\end{align*}
Hieraus folgt wiederum
\begin{align*}
  &\ \tJq-\Varnp(z)\tJq\Varnp^{\ast}(\omega) \\
  & = \tJq-\eklam{\Izq + (z-\alpha)\begin{pmatrix} \uarna \\ -\vna \end{pmatrix}\Rna(\za)\Harn^{-1}\Rn(\alpha)\begin{pmatrix} \vn & \uarn \end{pmatrix}}\tJq \\
  &\quad \cdot \eklam{\Izq + (\omegaa-\alpha)\begin{pmatrix} \vna \\ \uarna \end{pmatrix}\Rna(\alpha)\Harn^{-1}\Rn(\omegaa)\begin{pmatrix} \uarn & -\vn \end{pmatrix}} \\
  & = -i(z-\alpha)\begin{pmatrix} \uarna \\ -\vna \end{pmatrix}\Rna(\za)\Harn^{-1}\Rn(\alpha)\begin{pmatrix} \uarn & -\vn \end{pmatrix} \\
  &\quad +i(\omegaa-\alpha)\begin{pmatrix} \uarna \\ -\vna \end{pmatrix}\Rna(\alpha)\Harn^{-1}\Rn(\omegaa)\begin{pmatrix} \uarn & -\vn \end{pmatrix} \\
  &\quad -(z-\alpha)(\omegaa-\alpha)\begin{pmatrix} \uarna \\ -\vna \end{pmatrix}\Rna(\za)\Harn^{-1}\Rn(\alpha)\begin{pmatrix} \vn & \uarn \end{pmatrix}\tJq\begin{pmatrix} \vna \\ \uarna \end{pmatrix} \\
  &\quad \cdot \Rna(\alpha)\Harn^{-1}\Rn(\omegaa)\begin{pmatrix} \uarn & -\vn \end{pmatrix} \\
  & = -i(z-\omegaa)\begin{pmatrix} \uarna \\ -\vna \end{pmatrix}\Rna(\za)\Harn^{-1}\Rn(\omegaa)\begin{pmatrix} \uarn & -\vn \end{pmatrix}.
\end{align*}

Zu (b): Aus der Definition der Funktionen $\Varm$ für alle $m \in \Zok$ und $\Rn$ für alle $n \in \Zofk$ ergibt sich,
dass $\Varm$ für alle $m \in \Zok$ ein Matrixpolynom vom Grad nicht größer als $\fklamo{\kappa}+1$ ist und somit eine in $\C$ holomorphe $2$\textit{q}$\times2$\textit{q}-Matrixfunktion.

Sei zunächst $n \in \Zofk$. Wegen (a) gilt dann
\begin{align*}
  \tJq-\Varn(z)\tJq\Varn^{\ast}(z) = 2\im z\rklam{\Rn(\za)\begin{pmatrix} \un & -\vn \end{pmatrix}}^{\ast}\Hn^{-1}\Rn(\za)\begin{pmatrix} \un & -\vn \end{pmatrix}
\end{align*}
für alle $z \in \C$. Wegen \thref{adpbm1} gilt $\Hn^{-1} \in \C^{(n+1)q\times(n+1)q}_{>}$, 
also ist unter Beachtung von \thref{spbsp1}, \thref{spdef2} sowie der Teile (a) und (c) von \thref{spbm1} dann $\Varn(z)$ für alle $z \in \Pp$ eine $\tJq$-kontraktive Matrix und insbesondere ist $\Varn(x)$ für alle $x \in \R$ eine $\tJq$-unitäre Matrix.
Wegen \thref{spbsp1} und \thref{spdef3} gilt somit $\Varn \in \btJqPp$.

Sei nun $n \in \Zofkm$. Wegen (a) gilt dann
\begin{align*}
  \tJq-\Varnp(z)\tJq\Varnp^{\ast}(z) = 2\im z\eklam{\Rn(\za)\begin{pmatrix} \uarn & -\vn \end{pmatrix}}^{\ast}\Harn^{-1}\Rn(\za)\begin{pmatrix} \uarn & -\vn \end{pmatrix}
\end{align*}
für alle $z \in \C$. Wegen \thref{adpbm1} gilt $\Harn^{-1} \in \C^{(n+1)q\times(n+1)q}_{>}$, 
also ist unter Beachtung von \thref{spbsp1}, \thref{spdef2} sowie der Teile (a) und (c) von \thref{spbm1} dann $\Varnp(z)$ für alle $z \in \Pp$ eine $\tJq$-kontraktive Matrix und insbesondere ist $\Varnp(x)$ für alle $x \in \R$ eine $\tJq$-unitäre Matrix.
Wegen \thref{spbsp1} und \thref{spdef3} gilt somit $\Varnp \in \btJqPp$.

Zu (c): Für alle $z \in \C$ und $n \in \Zofk$ ist unter Beachtung von (b) und \thref{spbsp1} wegen Teil (a) von \thref{splm1} $\Varn(z)$ regulär und es gilt
\begin{align*}
  \Varn^{-1}(z) & = \tJq \Varn^{\ast}(\za) \tJq \\
  & = \tJq^2 + (z-\alpha)\tJq\begin{pmatrix} \vna \\ \un \end{pmatrix}\Rna(\alpha)\Hn^{-1}\Rn(z)\begin{pmatrix} \un & -\vn \end{pmatrix}\tJq \\  
  & = \Izq - (z-\alpha)\begin{pmatrix} u^{\ast}_n \\ -v^{\ast}_n \end{pmatrix}\Rna(\alpha)\Hn^{-1}\Rn(z)\begin{pmatrix}v_n & u_n\end{pmatrix}.
\end{align*}

Für alle $z \in \C$ und $n \in \Zofkm$ ist unter Beachtung von (b) und \thref{spbsp1} wegen Teil (a) von \thref{splm1} $\Varnp(z)$ regulär und es gilt
\begin{align*}
  \Varnp^{-1}(z) & = \tJq \Varnp^{\ast}(\za) \tJq \\
  & = \tJq^2 + (z-\alpha)\tJq\begin{pmatrix} \vna \\ \uarn \end{pmatrix}\Rna(\alpha)\Harn^{-1}\Rn(z)\begin{pmatrix} \uarn & -\vn \end{pmatrix}\tJq \\  
  & = \Izq - (z-\alpha)\begin{pmatrix} \uarna \\ -v^{\ast}_n \end{pmatrix}\Rna(\alpha)\Harn^{-1}\Rn(z)\begin{pmatrix}v_n & \uarn\end{pmatrix}.
\end{align*}

Zu (d): Seien $z, \omega \in \C$. Unter Beachtung von (b) und \thref{spbsp1} gelten wegen Teil (b) von \thref{splm1} und (a) dann
\begin{align*}
  \tJq-\Varn^{-\ast}(z)\tJq\Varn^{-1}(\omega) & = \tJq\beklam{\tJq-\Varn(\za)\tJq\Varn^{\ast}(\omegaa)}\tJq \\
  & = i(\omega-\za)\tJq\begin{pmatrix} u^{\ast}_n \\ -v^{\ast}_n \end{pmatrix}\Rna(z)\Hn^{-1}\Rn(\omega)\begin{pmatrix}u_n & -\vn\end{pmatrix}\tJq
\end{align*}
für alle $n \in \Zofk$ sowie
\begin{align*}
  \tJq-\Varnp^{-\ast}(z)\tJq\Varnp^{-1}(\omega) & = \tJq\beklam{\tJq-\Varnp(\za)\tJq\Varnp^{\ast}(\omegaa)}\tJq \\
  & = i(\omega-\za)\tJq\begin{pmatrix} \uarna \\ -v^{\ast}_n \end{pmatrix}\Rna(z)\Hn^{-1}\Rn(\omega)\begin{pmatrix}\uarn & -\vn\end{pmatrix}\tJq
\end{align*}
für alle $n \in \Zofkm$.

Zu (e): Sei $z \in \C\setminus\R$. Wegen \thref{spbsp1} gelten dann
\begin{align}	\label{drsa2bw9}
  \frac{1}{i(z-\za)}\binom{f(z)}{\Iq}^{\ast}\tJq\binom{f(z)}{\Iq} = \frac{f(z)-f^{\ast}(z)}{z-\za}
\end{align}
und
\begin{align}	\label{drsa2bw10}
  \frac{1}{i(z-\za)}\binom{f_{\ar}(z)}{\Iq}^{\ast}\tJq\binom{f_{\ar}(z)}{\Iq} = \frac{f_{\ar}(z)-f^{\ast}_{\ar}(z)}{z-\za}.
\end{align}
Weiterhin gelten
\begin{align}	\label{drsa2bw11}
  i\begin{pmatrix} \un & -\vn \end{pmatrix}\tJq\binom{f(z)}{\Iq} = i\begin{pmatrix} -i\vn & -i\un \end{pmatrix}\binom{f(z)}{\Iq} = \vn f(z)+\un
\end{align}
für alle $n \in \Zofk$ und
\begin{align}	\label{drsa2bw12}
  i\begin{pmatrix} \uarn & -\vn \end{pmatrix}\tJq\binom{f_{\ar}(z)}{\Iq} = i\begin{pmatrix} -i\vn & -i\uarn \end{pmatrix}\binom{f_{\ar}(z)}{\Iq} = \vn f_{\ar}(z)+\uarn
\end{align}
für alle $n \in \Zofkm$. Wegen (d), \fref{drsa2bw9}, \fref{drsa2bw11} und \thref{drdef3} gilt
\begin{align*}
  &\ \frac{1}{i(z-\za)}\binom{f(z)}{\Iq}^{\ast}\Varn^{-\ast}(z)\tJq\Varn^{-1}(z)\binom{f(z)}{\Iq} \\
  & = \frac{1}{i(z-\za)}\binom{f(z)}{\Iq}^{\ast}\eklam{\tJq-i(z-\za)\tJq\begin{pmatrix} u^{\ast}_n \\ -v^{\ast}_n \end{pmatrix}\Rna(z)\Hn^{-1}\Rn(z)\begin{pmatrix}u_n & -\vn\end{pmatrix}\tJq}\binom{f(z)}{\Iq} \\
  & = \frac{f(z)-f^{\ast}(z)}{z-\za} - \eklam{i\begin{pmatrix}u_n & -\vn\end{pmatrix}\tJq\binom{f(z)}{\Iq}}^{\ast}\Rna(z)\Hn^{-1}\Rn(z)i\begin{pmatrix}u_n & -\vn\end{pmatrix}\tJq\binom{f(z)}{\Iq} \\
  & = \frac{f(z)-f^{\ast}(z)}{z-\za} - \eklam{\Rn(z)(\vn f(z)+\un)}^{\ast}\Hn^{-1}\Rn(z)(\vn f(z)+\un) \\
  & = \dFfn(z)
\end{align*}
für alle $n \in \Zofk$. Wegen (d), \fref{drsa2bw10}, \fref{drsa2bw12} und \thref{drdef3} gilt weiterhin
\begin{align*}
  &\ \frac{1}{i(z-\za)}\binom{f_{\ar}(z)}{\Iq}^{\ast}\Varnp^{-\ast}(z)\tJq\Varnp^{-1}(z)\binom{f_{\ar}(z)}{\Iq} \\
  & = \frac{1}{i(z-\za)}\binom{f_{\ar}(z)}{\Iq}^{\ast}\eklam{\tJq-i(z-\za)\tJq\begin{pmatrix} \uarna \\ -v^{\ast}_n \end{pmatrix}\Rna(z)\Harn^{-1}\Rn(z)\begin{pmatrix}\uarn & -\vn\end{pmatrix}\tJq} \\
  &\quad \cdot\binom{f_{\ar}(z)}{\Iq} \\
  & = \frac{f_{\ar}(z)-f^{\ast}_{\ar}(z)}{z-\za} - \eklam{i\begin{pmatrix}\uarn & -\vn\end{pmatrix}\tJq\binom{f_{\ar}(z)}{\Iq}}^{\ast}\Rna(z)\Harn^{-1}\Rn(z) \\
  &\quad \cdot i\begin{pmatrix}\uarn & -\vn\end{pmatrix}\tJq\binom{f_{\ar}(z)}{\Iq} \\
  & = \frac{f_{\ar}(z)-f^{\ast}_{\ar}(z)}{z-\za} - \eklam{\Rn(z)(\vn f_{\ar}(z)+\uarn)}^{\ast}\Harn^{-1}\Rn(z)(\vn f_{\ar}(z)+\uarn) \\
  & = \dFfarn(z)
\end{align*}
für alle $n \in \Zofkm$. \bwend

Es sei bemerkt, dass die in \thref{drbz4} eingeführten Matrizen $\Varn$ bzw. $\Varnp$ mit den in \cite[Bemerkung 8.6]{Maka} bzw. \cite[Bemerkung 8.7]{Maka} für den allgemeineren Fall einer gegebenen Folge aus $\Keqka$ eingeführten Matrizen $U_{n,\alpha}$ bzw. $\widetilde{U}_{n,\alpha}$ übereinstimmen, wie einige einfache Rechnungen zeigen können. Man kann dann Teil (a) von \thref{drsa2} mithilfe einiger Rechnungen alternativ aus \cite[Lemma 8.21]{Maka} und \cite[Lemma 8.22]{Maka} gewinnen.

Teil (e) von \thref{drsa2} verdeutlicht in Kombination mit \thref{drbm3} die Relevanz der in \thref{drbz4} eingeführten Matrixfunktionen für unsere Zwecke. Mit ihrer Hilfe gelingt nämlich eine für unsere weiteren Überlegungen sehr hilfreiche Darstellung der linken Schurkomplemente der beiden Potapovschen Fundamentalmatrizen.
Wir werden nun mithilfe speziell gewählter $2$\textit{q}$\times2$\textit{q}-Matrizen die beiden Folgen von Matrixpolynomen aus \thref{drbz4} miteinander koppeln, ohne dabei die Zugehörigkeit zu $\tbtJqPp$ zu verlieren. Diese neu gewonnenen Folgen von Matrixpolynomen können wir dann weiterhin verwenden, um die linken Schur-Komplemente der beiden Potapovschen Fundamentalmatrizen darzustellen. Wir verwenden zwei Arten solcher Folgen von Matrixpolynomen: Die einen stellen später die geraden Glieder unserer Resolventenmatrix dar und die anderen die ungeraden Glieder, wobei wir zweite zuerst behandeln.

\begin{bez}	\thlabel{drbz5}
	Seien $\kappa \in \Na$, $\alpha \in \R$ und $\sjk \in \Kpqka$. Für $n \in \Zofkm$ sei
	\begin{align*}
    	M^{\sklam{s}}_{\arn} := \begin{pmatrix} \Iq & \ysna\big(\Hsarn\big)^{-1}\ysn \\ \Oq & \Iq \end{pmatrix}.
	\end{align*}
	Für $n \in \Zofk$ sei weiterhin
	\begin{align*}
		\widetilde{M}^{\sklam{s}}_{\arn} := \begin{pmatrix} \Iq & \Oq \\ -\vna\Rna(\alpha)\big(\Hsn\big)^{-1}\Rn(\alpha)\vn & \Iq \end{pmatrix}.
	\end{align*}
	Falls klar ist, von welchem $\sjk$ die Rede ist, lassen wir das \anf{$\sklam{s}$} als oberen Index weg.
\end{bez}

\begin{bem}	\thlabel{drbm4}
  	Seien $\kappa \in \Na$, $\alpha \in \R$ und $\sjk \in \Kpqka$. 
  	Dann sind $M_{\arn}$ für alle $n \in \Zofkm$ und $\widetilde{M}_{\arn}$ für alle $n \in \Zofk$ jeweils $\tJq$-unitär.
\end{bem}

\bwanf Wegen \thref{adpbm1} gelten
\begin{align*}
  \eklam{y^{\ast}_{0,n}\Harn^{-1}y_{0,n}}^{\ast} = y^{\ast}_{0,n}\Harn^{-1}y_{0,n}
\end{align*}
für alle $n \in \Zofkm$ und
\begin{align*}
  \eklam{-v^{\ast}_n\Rna(\alpha)\Hn^{-1}\Rn(\alpha)v_n}^{\ast} = -v^{\ast}_n\Rna(\alpha)\Hn^{-1}\Rn(\alpha)v_n
\end{align*}
für alle $n \in \Zofk$.
Hieraus folgt mithilfe von Teil (c) von \thref{spbm2} dann die Behauptung. \bwend

Es sei bemerkt, dass die in \thref{drbz5} eingeführten Matrizen $M_{\arn}$ bzw. $\widetilde{M}_{\arn}$ mit den in \cite[Lemma 8.10]{Maka} für den allgemeineren Fall einer gegebenen Folge aus $\Keqka$ eingeführten Matrizen $B_{n,\alpha}$ bzw. $\widetilde{B}_{n,\alpha}$ übereinstimmen. Somit kann \thref{drbm4} alternativ aus \cite[Lemma 8.10]{Maka} gewonnen werden.

\begin{bez}	\thlabel{drbz6}
	Seien $\kappa \in \Na$, $\alpha \in \R$ und $\sjk \in \Kpqka$. Für $z\in\C$ und \linebreak $n \in \Zofkm$ seien 
	\begin{align*}
    	\Tarn^{\sklam{s}}(z) := \Varn^{\sklam{s}}(z) M^{\sklam{s}}_{\arn} \quad \text{und} \quad \tTarn^{\sklam{s}}(z) := \Varnp^{\sklam{s}}(z) \widetilde{M}^{\sklam{s}}_{\arn}.
	\end{align*}
	Falls klar ist, von welchem $\sjk$ die Rede ist, lassen wir das \anf{$\sklam{s}$} als oberen Index weg.
\end{bez}

Wie eine nähere Analyse des Beweises von \thref{drsa2} zeigt, wurden die in \thref{drbz4} eingeführten $2$\textit{q}$\times2$\textit{q}-Matrixpolynome jeweils aus einer der beiden Potapovschen Fundamentalmatrizen gewonnen. Es kommt uns nun darauf an, eine entsprechende Kopplung zwischen jenen Matrixpolynomen herzustellen. Hierzu werden die beiden in \thref{drbz5} eingeführten $2$\textit{q}$\times2$\textit{q}-Matrizen herangezogen. Da beide Matrizen $\tJq$-unitär sind (siehe \thref{drbm4}), werden die in \thref{drbz6} eingeführten $2$\textit{q}$\times2$\textit{q}-Matrixpolynome dieselben Eigenschaften bezüglich der Signaturmatrix $\tJq$ besitzen. Die wesentliche neue Eigenschaft jener Matrixpolynome ist aber, dass diese durch eine spezielle Kopplung miteinander verbunden sind, die in Teil (b) des nachfolgenden Satzes explizit beschrieben wird.

\begin{satz}	\thlabel{drsa3}
  Seien $\kappa \in \Na$, $\alpha \in \R$ und $\sjk \in \Kpqka$. 
  \dgfa
  \begin{itemize}
    \item [\rm{(a)}] Es gelten
    \begin{align*}
      \Tarn(z) = \begin{pmatrix} \Iq+(z-\alpha)u^{\ast}_n\Rna(\za)\Hn^{-1}\Rn(\alpha)v_n & \uarn^{\ast}\Rna(\za)\Harn^{-1}y_{0,n} \\
      -(z-\alpha)v^{\ast}_n\Rna(\za)\Hn^{-1}\Rn(\alpha)v_n & \Iq-(z-\alpha)v^{\ast}_n\Rna(\za)\Harn^{-1}y_{0,n} \end{pmatrix}
    \end{align*}
    und
    \begin{align*}
      \tTarn(z) = \begin{pmatrix} \Iq+(z-\alpha)u^{\ast}_n\Rna(\za)\Hn^{-1}\Rn(\alpha)v_n & (z-\alpha)\uarn^{\ast}\Rna(\za)\Harn^{-1}y_{0,n} \\
      -v^{\ast}_n\Rna(\za)\Hn^{-1}\Rn(\alpha)v_n & \Iq-(z-\alpha)v^{\ast}_n\Rna(\za)\Harn^{-1}y_{0,n} \end{pmatrix}
    \end{align*}
    für alle $z \in \C$ und $n \in \Zofkm$.
    \item [\rm{(b)}] Es gilt
    \begin{align*}
      \tTarn(z) = \begin{pmatrix} (z-\alpha)\Iq & \Oq \\ \Oq & \Iq \end{pmatrix} \Tarn(z) \begin{pmatrix} (z-\alpha)^{-1}\Iq & \Oq \\ \Oq & \Iq \end{pmatrix}
    \end{align*}
    für alle $z \in \C\setminus\gklam{\alpha}$ und $n \in \Zofkm$.
    \item [\rm{(c)}] Es gelten
    \begin{align*}
      \Tarn(z)\tJq\Tarn^{\ast}(\omega) = \Varn(z)\tJq\Varn^{\ast}(\omega)
    \end{align*}
    und
    \begin{align*}
      \tTarn(z)\tJq\tTarn^{\ast}(\omega) = \Varnp(z)\tJq\Varnp^{\ast}(\omega)
    \end{align*}
    für alle $z, \omega \in \C$ und $n \in \Zofkm$.
    \item [\rm{(d)}] Es gelten $\Tarn,\tTarn \in \tbtJqPp$ für alle $n \in \Zofkm$.
    \item [\rm{(e)}] Für alle $z, \omega \in \C$ und $n \in \Zofkm$ sind $\Tarn(z)$ und $\tTarn(z)$ regulär und es gelten
    \begin{align*}
      \Tarn^{-\ast}(z)\tJq\Tarn^{-1}(\omega) = \Varn^{-\ast}(z)\tJq\Varn^{-1}(\omega)
    \end{align*}
    sowie
    \begin{align*}
      \tTarn^{-\ast}(z)\tJq\tTarn^{-1}(\omega) = \Varnp^{-\ast}(z)\tJq\Varnp^{-1}(\omega).
    \end{align*}
    \item [\rm{(f)}] Seien $f: \C\setminus\R \rightarrow \Cqq$ und $f_{\ar}(z) := (z-\alpha)f(z)$ für alle $z \in \C\setminus\R$. 
    Unter Beachtung von (e) gelten dann
    \begin{align*}
      \dFfn(z) = \frac{1}{i(z-\za)}\binom{f(z)}{\Iq}^{\ast}\Tarn^{-\ast}(z)\tJq\Tarn^{-1}(z)\binom{f(z)}{\Iq}
    \end{align*}
    und
    \begin{align*}
      \dFfarn(z) = \frac{1}{i(z-\za)}\binom{f_{\ar}(z)}{\Iq}^{\ast}\tTarn^{-\ast}(z)\tJq\tTarn^{-1}(z)\binom{f_{\ar}(z)}{\Iq}
    \end{align*}
    für alle $z \in \C\setminus\R$ und $n \in \Zofkm$.
  \end{itemize}
\end{satz}

\bwanf Zu (a): Seien $z \in \C$ und $n \in \Zofkm$. Bezeichne
\begin{align*}
  \Tarn = \begin{pmatrix} \Tarn^{(1,1)} & \Tarn^{(1,2)} \\ \Tarn^{(2,1)} & \Tarn^{(2,2)} \end{pmatrix} \quad \text{bzw.} \quad
  \Varn = \begin{pmatrix} \Varn^{(1,1)} & \Varn^{(1,2)} \\ \Varn^{(2,1)} & \Varn^{(2,2)} \end{pmatrix}
\end{align*}
die \textit{q}$\times$\textit{q}-Blockzerlegung von $\Tarn$ bzw. $\Varn$. Dann gelten
\begin{align*}
  \Tarn^{(1,1)} = \Varn^{(1,1)} = \Iq+(z-\alpha)\una\Rna(\za)\Hn^{-1}\Rn(\alpha)\vn
\end{align*}
und
\begin{align*}
  \Tarn^{(2,1)} = \Varn^{(2,1)} = -(z-\alpha)\vna\Rna(\za)\Hn^{-1}\Rn(\alpha)\vn.
\end{align*}
Wegen Teil (c) von \thref{drbm1} gilt
\begin{align}	\label{drsa3bw1}
  \Rn(\alpha)\vn\yna = \Hn - \Rn(\alpha)\Tn\Harn.
\end{align}
Wegen der Teile (e) und (a) von \thref{drlm3} sowie Teil (d) von \thref{drlm1} gilt weiterhin
\begin{align}	\label{drsa3bw2}
  \uarna\Rna(\za)-\yna & = \uarna\eklam{\Rna(\za)-\Rna(\alpha)} \notag \\
  & = (z-\alpha)\uarna\Rna(\alpha)\Tna\Rna(\za) \notag \\
  & = (z-\alpha)\yna\Tna\Rna(\za) \notag \\
  & = (z-\alpha)\una\Rna(\za).
\end{align}
Wegen \fref{drsa3bw1}, \fref{drsa3bw2} und Teil (a) von \thref{drlm3} gilt nun
\begin{align*}
  \Tarn^{(1,2)} & = \Varn^{(1,1)}\yna\Harn^{-1}\yn + \Varn^{(1,2)} \\
  & = \yna\Harn^{-1}\yn + (z-\alpha)\una\Rna(\za)\Hn^{-1}\Rn(\alpha)\vn\yna\Harn^{-1}\yn \\
  &\quad + (z-\alpha)\una\Rna(\za)\Hn^{-1}\Rn(\alpha)\un \\
  & = \yna\Harn^{-1}\yn + (z-\alpha)\una\Rna(\za)\Harn^{-1}\yn \\
  &\quad - (z-\alpha)\una\Rna(\za)\Hn^{-1}\Rn(\alpha)\Tn\yn + (z-\alpha)\una\Rna(\za)\Hn^{-1}\Rn(\alpha)\un \\
  & = \yna\Harn^{-1}\yn + (\uarna\Rna(\za)-\yna)\Harn^{-1}\yn \\
  &\quad - (z-\alpha)\una\Rna(\za)\Hn^{-1}\Rn(\alpha)\un + (z-\alpha)\una\Rna(\za)\Hn^{-1}\Rn(\alpha)\un \\
  & = \uarna\Rna(\za)\Harn^{-1}\yn.
\end{align*}
Wegen \fref{drsa3bw1} und Teil (a) von \thref{drlm3} gilt weiterhin
\begin{align*}
  \Tarn^{(2,2)} & = \Varn^{(2,1)}\yna\Harn^{-1}\yn + \Varn^{(2,2)} \\
  & = -(z-\alpha)\vna\Rna(\za)\Hn^{-1}\Rn(\alpha)\vn\yna\Harn^{-1}\yn \\
  &\quad + \Iq-(z-\alpha)\vna\Rna(\za)\Hn^{-1}\Rn(\alpha)\un \\
  & = -(z-\alpha)\vna\Rna(\za)\Harn^{-1}\yn + (z-\alpha)\vna\Rna(\za)\Hn^{-1}\Rn(\alpha)\Tn\yn \\
  &\quad + \Iq-(z-\alpha)\vna\Rna(\za)\Hn^{-1}\Rn(\alpha)\un \\
  & = -(z-\alpha)\vna\Rna(\za)\Harn^{-1}\yn + (z-\alpha)\vna\Rna(\za)\Hn^{-1}\Rn(\alpha)\un \\
  &\quad + \Iq-(z-\alpha)\vna\Rna(\za)\Hn^{-1}\Rn(\alpha)\un \\
  & = \Iq-(z-\alpha)\vna\Rna(\za)\Harn^{-1}\yn.  
\end{align*}

Bezeichne
\begin{align*}
  \tTarn = \begin{pmatrix} \tTarn^{(1,1)} & \tTarn^{(1,2)} \\ \tTarn^{(2,1)} & \tTarn^{(2,2)} \end{pmatrix} \quad \text{bzw.} \quad
  \Varnp = \begin{pmatrix} \Varnp^{(1,1)} & \Varnp^{(1,2)} \\ \Varnp^{(2,1)} & \Varnp^{(2,2)} \end{pmatrix}
\end{align*}
die \textit{q}$\times$\textit{q}-Blockzerlegung von $\tTarn$ bzw. $\Varnp$. Wegen Teil (e) von \thref{drlm3} gelten dann
\begin{align*}
  \tTarn^{(1,2)} & = \Varnp^{(1,2)} = (z-\alpha)\uarna\Rna(\za)\Harn^{-1}\Rn(\alpha)\uarn \\
  & = (z-\alpha)\uarna\Rna(\za)\Harn^{-1}\yn
\end{align*}
und
\begin{align*}
  \tTarn^{(2,2)} & = \Varnp^{(2,2)} = \Iq-(z-\alpha)\vna\Rna(\za)\Harn^{-1}\Rn(\alpha)\uarn \\
  & = \Iq-(z-\alpha)\vna\Rna(\za)\Harn^{-1}\yn.
\end{align*}
Wegen Teil (e) von \thref{drlm3} sowie Teil (c) von \thref{drbm1} gilt
\begin{align}	\label{drsa3bw3}
  \Rn(\alpha)\uarn\vna\Rna(\alpha) = \yn\vna\Rna(\alpha) = \Hn-\Harn\Tna\Rna(\alpha).
\end{align}
Wegen der Teile (a) und (b) von \thref{drlm1} sowie der Teile (e) und (a) von \thref{drlm3} gilt weiterhin
\begin{align}	\label{drsa3bw4}
  \uarna\Rna(\za)\Tna\Rna(\alpha) = \uarna\Rna(\alpha)\Tna\Rna(\za) = \yna\Tna\Rna(\za) = \una\Rna(\za).
\end{align}
Wegen \fref{drsa3bw3} und \fref{drsa3bw4} gilt nun
\begin{align*}
  \tTarn^{(1,1)} & = \Varnp^{(1,1)} - \Varnp^{(1,2)}\vna\Rna(\alpha)\Hn^{-1}\Rn(\alpha)\vn \\
  & = \Iq+(z-\alpha)\uarna\Rna(\za)\Harn^{-1}\Rn(\alpha)\vn \\
  &\quad - (z-\alpha)\uarna\Rna(\za)\Harn^{-1}\Rn(\alpha)\uarn\vna\Rna(\alpha)\Hn^{-1}\Rn(\alpha)\vn \\
  & = \Iq+(z-\alpha)\uarna\Rna(\za)\Harn^{-1}\Rn(\alpha)\vn \\
  &\quad - (z-\alpha)\uarna\Rna(\za)\Harn^{-1}\Rn(\alpha)\vn + (z-\alpha)\uarna\Rna(\za)\Tna\Rna(\alpha)\Hn^{-1}\Rn(\alpha)\vn \\
  & = \Iq+(z-\alpha)\una\Rna(\za)\Hn^{-1}\Rn(\alpha)\vn.
\end{align*}
Wegen \fref{drsa3bw3} und Teil (d) von \thref{drlm1} gilt weiterhin 
\begin{align*}
  \tTarn^{(2,1)} & = \Varnp^{(2,1)} - \Varnp^{(2,2)}\vna\Rna(\alpha)\Hn^{-1}\Rn(\alpha)\vn \\
  & = -(z-\alpha)\vna\Rna(\za)\Harn^{-1}\Rn(\alpha)\vn - \vna\Rna(\alpha)\Hn^{-1}\Rn(\alpha)\vn \\
  &\quad + (z-\alpha)\vna\Rna(\za)\Harn^{-1}\Rn(\alpha)\uarn\vna\Rna(\alpha)\Hn^{-1}\Rn(\alpha)\vn \\
  & = -(z-\alpha)\vna\Rna(\za)\Harn^{-1}\Rn(\alpha)\vn - \vna\Rna(\alpha)\Hn^{-1}\Rn(\alpha)\vn \\
  &\quad + (z-\alpha)\vna\Rna(\za)\Harn^{-1}\Rn(\alpha)\vn - (z-\alpha)\vna\Rna(\za)\Tna\Rna(\alpha)\Hn^{-1}\Rn(\alpha)\vn \\
  & = - \vna\Rna(\alpha)\Hn^{-1}\Rn(\alpha)\vn + \vna(\Rna(\alpha)-\Rna(\za))\Hn^{-1}\Rn(\alpha)\vn \\
  & = - \vna\Rna(\za)\Hn^{-1}\Rn(\alpha)\vn.  
\end{align*}

Zu (b): Dies ist eine direkte Konsequenz aus (a).

Zu (c): Seien $z, \omega \in \C$ und $n \in \Zofkm$. Wegen \thref{drbm4} gelten dann
\begin{align*}
  \Tarn(z)\tJq\Tarn^{\ast}(\omega) = \Varn(z)M_{\arn}\tJq M^{\ast}_{\arn}\Varn^{\ast}(\omega) = \Varn(z)\tJq\Varn^{\ast}(\omega)
\end{align*}
und
\begin{align*}
  \tTarn(z)\tJq\tTarn^{\ast}(\omega) = \Varnp(z)\widetilde{M}_{\arn}\tJq\widetilde{M}^{\ast}_{\arn}\Varnp^{\ast}(\omega) = \Varnp(z)\tJq\Varnp^{\ast}(\omega).
\end{align*}

Zu (d): Dies folgt wegen (c), \thref{spbsp1} und Teil (b) von \thref{drsa2} aus \thref{spdef3}.

Zu (e): Dies folgt wegen (d), Teil (b) von \thref{drsa2} und \thref{spbsp1} aus (c) sowie den Teilen (a) und (b) von \thref{splm1}.

Zu (f): Dies folgt wegen (e) aus Teil (e) von \thref{drsa2}. \bwend

Es sei bemerkt, dass Teil (a) von \thref{drsa3} und \cite[Bemerkung 8.15]{Maka} zeigen, dass die in \thref{drbz6} eingeführten Matrizen $\Tarn$ bzw. $\tTarn$ mit den in \cite[Lemma 8.14]{Maka} für den allgemeineren Fall einer gegebenen Folge aus $\Keqka$ eingeführten Matrizen $\Theta_{n,\alpha}$ bzw. $\widetilde{\Theta}_{n,\alpha}$ übereinstimmen, wie eine einfache Rechnung zeigen kann. Man kann dann Teil (b) bzw. (c) bzw. (d) bzw. (e) bzw. (f) von \thref{drsa3} mithilfe einiger Rechnungen alternativ aus \cite[Lemma 8.16]{Maka} bzw. \cite[Lemma 8.17]{Maka} in Verbindung mit \cite[Lemma 8.8]{Maka} und \cite[Lemma 8.9]{Maka} bzw. \cite[Lemma 8.19]{Maka} in Verbindung mit \thref{spdef3} bzw. \cite[Lemma 8.20]{Maka} in Verbindung mit Teil (d) von \thref{drsa2} bzw. \cite[Lemma 9.1]{Maka} und \cite[Lemma 9.2]{Maka} gewinnen.

Die in \thref{drbz6} eingeführten $2$\textit{q}$\times2$\textit{q}-Matrixpolynome werden uns später im Falle ungeradzahliger Indizes nützliche Dienste erweisen. Nachfolgend führen wir nun analoge Betrachtungen durch, welche im Fall geradzahliger Indizes Anwendung finden werden.

\begin{bez}	\thlabel{drbz7}
	Seien $\kappa \in \Na\setminus\gklam{1}$, $\alpha \in \R$ und $\sjk \in \Kpqka$. Für $z\in\C$ und $n \in \Zefk$ seien 
	\begin{align*}
    	\Yarn^{\sklam{s}}(z) := \Varn^{\sklam{s}}(z) M^{\sklam{s}}_{\arn-1} \quad \text{und} \quad \tYarn^{\sklam{s}}(z) := \Varnm^{\sklam{s}}(z) \widetilde{M}^{\sklam{s}}_{\arn}.
	\end{align*}
	Falls klar ist, von welchem $\sjk$ die Rede ist, lassen wir das \anf{$\sklam{s}$} als oberen Index weg.
\end{bez}

\begin{satz}	\thlabel{drsa4}
  Seien $\kappa \in \Na\setminus\gklam{1}$, $\alpha \in \R$ und $\sjk \in \Kpqka$. 
  \dgfa
  \begin{itemize}
    \item [\rm{(a)}] Es gelten
    \begin{align*}
      \Yarn(z) &= \left( \begin{matrix} \Iq+(z-\alpha)u^{\ast}_n\Rna(\za)\Hn^{-1}\Rn(\alpha)v_n \\ -(z-\alpha)v^{\ast}_n\Rna(\za)\Hn^{-1}\Rn(\alpha)v_n \end{matrix} \right. \\
      &\left. \hspace{6cm} \begin{matrix} \uarnm^{\ast}\Rnma(\za)\Harnm^{-1}y_{0,n-1} \\ \Iq-(z-\alpha)v^{\ast}_{n-1}\Rnma(\za)\Harnm^{-1}y_{0,n-1} \end{matrix} \right)
    \end{align*}
    und
    \begin{align*}
      \tYarn(z) &= \left( \begin{matrix} \Iq+(z-\alpha)u^{\ast}_n\Rna(\za)\Hn^{-1}\Rn(\alpha)v_n \\ -v^{\ast}_n\Rna(\za)\Hn^{-1}\Rn(\alpha)v_n \end{matrix} \right. \\
      &\left. \hspace{6cm} \begin{matrix} (z-\alpha)\uarnm^{\ast}\Rnma(\za)\Harnm^{-1}y_{0,n-1} \\ \Iq-(z-\alpha)v^{\ast}_{n-1}\Rnma(\za)\Harnm^{-1}y_{0,n-1} \end{matrix} \right)
    \end{align*}
    für alle $z \in \C$ und $n \in \Zefk$.
    \item [\rm{(b)}] Es gilt
    \begin{align*}
      \tYarn(z) = \begin{pmatrix} (z-\alpha)\Iq & \Oq \\ \Oq & \Iq \end{pmatrix} \Yarn(z) \begin{pmatrix} (z-\alpha)^{-1}\Iq & \Oq \\ \Oq & \Iq \end{pmatrix}
    \end{align*}
    für alle $z \in \C\setminus\gklam{\alpha}$ und $n \in \Zefk$.
    \item [\rm{(c)}] Es gelten
    \begin{align*}
      \Yarn(z)\tJq\Yarn^{\ast}(\omega) = \Varn(z)\tJq\Varn^{\ast}(\omega)
    \end{align*}
    und
    \begin{align*}
      \tYarn(z)\tJq\tYarn^{\ast}(\omega) = \Varnm(z)\tJq\Varnm^{\ast}(\omega)
    \end{align*}
    für alle $z, \omega \in \C$ und $n \in \Zefk$.
    \item [\rm{(d)}] Es gelten $\Yarn,\tYarn \in \tbtJqPp$ für alle $n \in \Zefk$.
    \item [\rm{(e)}] Für alle $z, \omega \in \C$ und $n \in \Zefk$ sind $\Yarn(z)$ und $\tYarn(z)$ regulär und es gelten
    \begin{align*}
      \Yarn^{-\ast}(z)\tJq\Yarn^{-1}(\omega) = \Varn^{-\ast}(z)\tJq\Varn^{-1}(\omega)
    \end{align*}
    sowie
    \begin{align*}
      \tYarn^{-\ast}(z)\tJq\tYarn^{-1}(\omega) = \Varnm^{-\ast}(z)\tJq\Varnm^{-1}(\omega).
    \end{align*}
    \item [\rm{(f)}] Seien $f: \C\setminus\R \rightarrow \Cqq$ und $f_{\ar}(z) := (z-\alpha)f(z)$ für alle $z \in \C\setminus\R$. 
    Unter Beachtung von (e) gelten dann
    \begin{align*}
      \dFfn(z) = \frac{1}{i(z-\za)}\binom{f(z)}{\Iq}^{\ast}\Yarn^{-\ast}(z)\tJq\Yarn^{-1}(z)\binom{f(z)}{\Iq}
    \end{align*}
    und
    \begin{align*}
      \dFfarnm(z) = \frac{1}{i(z-\za)}\binom{f_{\ar}(z)}{\Iq}^{\ast}\tYarn^{-\ast}(z)\tJq\tYarn^{-1}(z)\binom{f_{\ar}(z)}{\Iq}
    \end{align*}
    für alle $z \in \C\setminus\R$ und $n \in \Zefk$.
  \end{itemize}
\end{satz}

\bwanf Zu (a): Seien $z \in \C$ und $n \in \Zefk$. Bezeichne
\begin{align*}
  \Yarn = \begin{pmatrix} \Yarn^{(1,1)} & \Yarn^{(1,2)} \\ \Yarn^{(2,1)} & \Yarn^{(2,2)} \end{pmatrix} \quad \text{bzw.} \quad
  \Varn = \begin{pmatrix} \Varn^{(1,1)} & \Varn^{(1,2)} \\ \Varn^{(2,1)} & \Varn^{(2,2)} \end{pmatrix}
\end{align*}
die \textit{q}$\times$\textit{q}-Blockzerlegung von $\Yarn$ bzw. $\Varn$. Dann gelten
\begin{align*}
  \Yarn^{(1,1)} = \Varn^{(1,1)} = \Iq+(z-\alpha)\una\Rna(\za)\Hn^{-1}\Rn(\alpha)\vn
\end{align*}
und
\begin{align*}
  \Yarn^{(2,1)} = \Varn^{(2,1)} = -(z-\alpha)\vna\Rna(\za)\Hn^{-1}\Rn(\alpha)\vn.
\end{align*}
Wegen Teil (d) von \thref{drbm1} gilt
\begin{align}	\label{drsa4bw1}
  \Rn(\alpha)\vn\ynma = \Hn\dLn-\Rn(\alpha)\Ln\Harnm.
\end{align}
Wegen der Teile (e) und (a) von \thref{drlm3} sowie Teil (d) von \thref{drlm1} gilt weiterhin
\begin{align}	\label{drsa4bw2}
  \uarnma\Rnma(\za)-\ynma & = \uarnma(\Rnma(\za)-\Rnma(\alpha)) \notag \\
  & = (z-\alpha)\uarnma\Rnma(\alpha)\Tnma\Rnma(\za) \notag \\
  & = (z-\alpha)\ynma\Tnma\Rnma(\za) \notag \\
  & = (z-\alpha)\unma\Rnma(\za).
\end{align}
Wegen \fref{drsa4bw1}, Teil (d) von \thref{drlm2}, der Teile (c) und (b) von \thref{drlm3} sowie \fref{drsa4bw2} gilt nun
\begin{align*}
  \Yarn^{(1,2)} & = \Varn^{(1,1)}\ynma\Harnm^{-1}\ynm + \Varn^{(1,2)} \\
  & = \ynma\Harnm^{-1}\ynm + (z-\alpha)\una\Rna(\za)\Hn^{-1}\Rn(\alpha)\vn\ynma\Harnm^{-1}\ynm \\
  &\quad + (z-\alpha)\una\Rna(\za)\Hn^{-1}\Rn(\alpha)\un \\
  & = \ynma\Harnm^{-1}\ynm + (z-\alpha)\una\Rna(\za)\dLn\Harnm^{-1}\ynm \\
  &\quad - (z-\alpha)\una\Rna(\za)\Hn^{-1}\Rn(\alpha)\Ln\ynm + (z-\alpha)\una\Rna(\za)\Hn^{-1}\Rn(\alpha)\un \\
  & = \ynma\Harnm^{-1}\ynm + (z-\alpha)\unma\Rnma(\za)\Harnm^{-1}\ynm \\
  &\quad - (z-\alpha)\una\Rna(\za)\Hn^{-1}\Rn(\alpha)\un + (z-\alpha)\una\Rna(\za)\Hn^{-1}\Rn(\alpha)\un \\
  & = \ynma\Harnm^{-1}\ynm + (\uarnma\Rnma(\za)-\ynma)\Harnm^{-1}\ynm \\
  & = \uarnma\Rnma(\za)\Harnm^{-1}\ynm.
\end{align*}
Wegen \fref{drsa4bw1}, der Teile (d) und (b) von \thref{drlm2} sowie Teil (b) von \thref{drlm3} gilt weiterhin
\begin{align*}
  \Yarn^{(2,2)} & = \Varn^{(2,1)}\ynma\Harnm^{-1}\ynm + \Varn^{(2,2)} \\
  & = -(z-\alpha)\vna\Rna(\za)\Hn^{-1}\Rn(\alpha)\vn\ynma\Harnm^{-1}\ynm \\
  &\quad + \Iq-(z-\alpha)\vna\Rna(\za)\Hn^{-1}\Rn(\alpha)\un \\
  & = -(z-\alpha)\vna\Rna(\za)\dLn\Harnm^{-1}\ynm + (z-\alpha)\vna\Rna(\za)\Hn^{-1}\Rn(\alpha)\Ln\ynm \\
  &\quad + \Iq-(z-\alpha)\vna\Rna(\za)\Hn^{-1}\Rn(\alpha)\un \\
  & = -(z-\alpha)\vnma\Rnma(\za)\Harnm^{-1}\ynm + (z-\alpha)\vna\Rna(\za)\Hn^{-1}\Rn(\alpha)\un \\
  &\quad + \Iq-(z-\alpha)\vna\Rna(\za)\Hn^{-1}\Rn(\alpha)\un \\
  & = \Iq-(z-\alpha)\vnma\Rnma(\za)\Harnm^{-1}\ynm.
\end{align*}
Bezeichne
\begin{align*}
  \tYarn = \begin{pmatrix} \Yarn^{(1,1)} & \Yarn^{(1,2)} \\ \Yarn^{(2,1)} & \Yarn^{(2,2)} \end{pmatrix} \quad \text{bzw.} \quad
  \Varnm = \begin{pmatrix} \Varnm^{(1,1)} & \Varnm^{(1,2)} \\ \Varnm^{(2,1)} & \Varnm^{(2,2)} \end{pmatrix}
\end{align*}
die \textit{q}$\times$\textit{q}-Blockzerlegung von $\tYarn$ bzw. $\Varnm$. Wegen Teil (e) von \thref{drlm3} gelten dann
\begin{align*}
  \tYarn^{(1,2)} & = \Varnm^{(1,2)} = (z-\alpha)\uarnma\Rnma(\za)\Harnm^{-1}\Rnm(\alpha)\uarnm \\
  & = (z-\alpha)\uarnma\Rnma(\za)\Harnm^{-1}\ynm
\end{align*}
und
\begin{align*}
  \tYarn^{(2,2)} & = \Varnm^{(2,2)} = \Iq-(z-\alpha)\vnma\Rnma(\za)\Harnm^{-1}\Rnm(\alpha)\uarnm \\
  & = \Iq-(z-\alpha)\vnma\Rnma(\za)\Harnm^{-1}\ynm.
\end{align*}
Wegen Teil (e) von \thref{drlm3} sowie Teil (d) von \thref{drbm1} gilt
\begin{align}	\label{drsa4bw3}
  \Rnm(\alpha)\uarnm\vna\Rna(\alpha) = \ynm\vna\Rna(\alpha) = \dLna\Hn-\Harnm\Lna\Rna(\alpha).
\end{align}
Wegen Teil (d) von \thref{drlm2}, Teil (b) von \thref{drlm1} sowie der Teile (e) und (b) von \thref{drlm3} gilt weiterhin
\begin{align}	\label{drsa4bw4}
  \uarnma\Rnma(\za)\Lna\Rna(\alpha) = \uarnma\Rnma(\alpha)\Lna\Rna(\za) = \ynma\Lna\Rna(\za) = \una\Rna(\za).
\end{align}
Wegen der Teile (b), (d) und (a) von \thref{drlm2} sowie der Teile (a), (b) und (d) von \thref{drlm1} gilt weiterhin
\begin{align}	\label{drsa4bw5}
  (z-\alpha)\vnma\Rnma(\za)\Lna\Rna(\alpha) & = (z-\alpha)\vna\dLn\Rnma(\za)\Lna\Rna(\alpha) \notag \\
  & = (z-\alpha)\vna\Rna(\za)\Tna\Rna(\alpha) \notag \\
  & = (z-\alpha)\vna\Rna(\alpha)\Tna\Rna(\za) \notag \\
  & = \vna\eklam{\Rna(\za)-\Rna(\alpha)}.
\end{align}
Wegen \fref{drsa4bw3}, der Teile (d) und (b) von \thref{drlm2} sowie \fref{drsa4bw4} gilt nun
\begin{align*}
  \tYarn^{(1,1)} & = \Varnm^{(1,1)} - \Varnm^{(1,2)}\vna\Rna(\alpha)\Hn^{-1}\Rn(\alpha)\vn \\
  & = \Iq+(z-\alpha)\uarnma\Rnma(\za)\Harnm^{-1}\Rnm(\alpha)\vnm \\
  &\quad - (z-\alpha)\uarnma\Rnma(\za)\Harnm^{-1}\Rnm(\alpha)\uarnm\vna\Rna(\alpha)\Hn^{-1}\Rn(\alpha)\vn \\
  & = \Iq+(z-\alpha)\uarnma\Rnma(\za)\Harnm^{-1}\Rnm(\alpha)\vnm \\
  &\quad - (z-\alpha)\uarnma\Rnma(\za)\Harnm^{-1}\dLna\Rn(\alpha)\vn \\
  &\quad + (z-\alpha)\uarnma\Rnma(\za)\Lna\Rna(\alpha)\Hn^{-1}\Rn(\alpha)\vn \\
  & = \Iq+(z-\alpha)\uarnma\Rnma(\za)\Harnm^{-1}\Rnm(\alpha)\vnm \\
  &\quad - (z-\alpha)\uarnma\Rnma(\za)\Harnm^{-1}\Rnm(\alpha)\vnm + (z-\alpha)\una\Rna(\za)\Hn^{-1}\Rn(\alpha)\vn \\
  & = \Iq+(z-\alpha)\una\Rna(\za)\Hn^{-1}\Rn(\alpha)\vn.
\end{align*}
Wegen \fref{drsa4bw3}, der Teile (d) und (b) von \thref{drlm2} sowie \fref{drsa4bw5} gilt weiterhin
\begin{align*}
  \tYarn^{(2,1)} & = \Varnm^{(2,1)} - \Varnm^{(2,2)}\vna\Rna(\alpha)\Hn^{-1}\Rn(\alpha)\vn \\
  & = -(z-\alpha)\vnma\Rnma(\za)\Harnm^{-1}\Rnm(\alpha)\vnm - \vna\Rna(\alpha)\Hn^{-1}\Rn(\alpha)\vn \\
  &\quad + (z-\alpha)\vnma\Rnma(\za)\Harnm^{-1}\Rnm(\alpha)\uarnm\vna\Rna(\alpha)\Hn^{-1}\Rn(\alpha)\vn \\
  & = -(z-\alpha)\vnma\Rnma(\za)\Harnm^{-1}\Rnm(\alpha)\vnm - \vna\Rna(\alpha)\Hn^{-1}\Rn(\alpha)\vn \\
  &\quad + (z-\alpha)\vnma\Rnma(\za)\Harnm^{-1}\dLna\Rn(\alpha)\vn \\
  &\quad - (z-\alpha)\vnma\Rnma(\za)\Lna\Rna(\alpha)\Hn^{-1}\Rn(\alpha)\vn \\
  & = -(z-\alpha)\vnma\Rnma(\za)\Harnm^{-1}\Rnm(\alpha)\vnm - \vna\Rna(\alpha)\Hn^{-1}\Rn(\alpha)\vn \\
  &\quad + (z-\alpha)\vnma\Rnma(\za)\Harnm^{-1}\Rnm(\alpha)\vnm - \vna(\Rna(\za)-\Rna(\alpha))\Hn^{-1}\Rn(\alpha)\vn \\
  & = - \vna\Rna(\za)\Hn^{-1}\Rn(\alpha)\vn.
\end{align*}

Zu (b): Dies ist eine direkte Konsequenz aus (a).

Zu (c): Seien $z, \omega \in \C$ und $n \in \Zefk$. Wegen \thref{drbm4} gelten dann
\begin{align*}
  \Yarn(z)\tJq\Yarn^{\ast}(\omega) = \Varn(z)M_{\arn-1}\tJq M^{\ast}_{\arn-1}\Varn^{\ast}(\omega) = \Varn(z)\tJq\Varn^{\ast}(\omega)
\end{align*}
und
\begin{align*}
  \tYarn(z)\tJq\tYarn^{\ast}(\omega) = \Varnm(z)\widetilde{M}_{\arn}\tJq\widetilde{M}^{\ast}_{\arn}\Varnm^{\ast}(\omega) = \Varnm(z)\tJq\Varnm^{\ast}(\omega).
\end{align*}

Zu (d): Dies folgt wegen (c), \thref{spbsp1} und Teil (b) von \thref{drsa2} aus \thref{spdef3}.

Zu (e): Dies folgt wegen (d), Teil (b) von \thref{drsa2} und \thref{spbsp1} aus (c) sowie den Teilen (a) und (b) von \thref{splm1}.

Zu (f): Dies folgt wegen (e) aus Teil (e) von \thref{drsa2}. \bwend

Es sei bemerkt, dass die in \thref{drbz7} eingeführten Matrizen $\Yarn$ und $\tYarn$ \textit{nicht} mit den in \cite[Formel (8.57)]{Maka} für den allgemeineren Fall einer gegebenen Folge aus $\Keqka$ eingeführten Matrizen $\Phi_{n,\alpha}$ bzw. $\widetilde{\Phi}_{n,\alpha}$ übereinstimmen. Die dortige Definition zeigt in Wahrheit, dass $\Phi_{n,\alpha}$ bzw. $\widetilde{\Phi}_{n,\alpha}$ mit den in \thref{drbz6} eingeführten Matrizen $\Tarn$ bzw. $\tTarn$ übereinstimmen.

Die Sätze \ref{drsa3} und \ref{drsa4} führen uns auf ein spezielles Quadrupel von Folgen von \textit{q}$\times$\textit{q}-Matrixpolynomen, welches in unseren weiteren Betrachtungen eine zentrale Rolle spielen wird.

\begin{defi}	\thlabel{drdef1}
  Seien $\kappa \in \Na$, $\alpha \in \R$ und $\sjk \in \Kpqka$. Weiterhin seien für alle $z \in \C$
  \begin{align*}
    \quad \Bsaro(z) := \Oq, \quad \quad \Dsaro(z) := \Iq
  \end{align*}
  und
  \begin{align*}
    \Asarn(z) & := \Iq+(z-\alpha)\usna\Rna(\za)\big(\Hsn\big)^{-1}\Rn(\alpha)\vn, \\
    \Csarn(z) & := -(z-\alpha)\vna\Rna(\za)\big(\Hsn\big)^{-1}\Rn(\alpha)\vn
  \end{align*}
  für alle $n \in \Zofk$ sowie
  \begin{align*}
    \Bsarn(z) & := \usarnma\Rnma(\za)\big(\Hsarnm\big)^{-1}\ysnm, \\
    \Dsarn(z) & := \Iq-(z-\alpha)\vnma\Rnma(\za)\big(\Hsarnm\big)^{-1}\ysnm
  \end{align*}
  für alle $n \in \Zefkp$. Dann heißt $\ABCDsark$ das \textbf{rechtsseitige $\alpha$-Dyukarev-Quadrupel} bezüglich $\sjk$. Falls klar ist, von welchem $\sjk$ die Rede ist, lassen wir das \anf{$\sklam{s}$} als oberen Index weg.
\end{defi}

Es sei bemerkt, dass die in \thref{drdef1} eingeführten Größen mit denen in \cite[Chapter 3]{CR1} im Fall $\alpha=0$ und $\kappa=\infty$ übereinstimmen, wobei dort in der Formel (3.9) ein $v_j$ statt $u_{1,j}$ stehen müsste.

%

In der Situation von \thref{drdef1} bilden wir nun aus den dort eingeführten Folgen von \textit{q}$\times$\textit{q}-Matrixpolynomen eine spezielle Folge von $2$\textit{q}$\times2$\textit{q}-Matrixpolynomen.

\begin{defi}	\thlabel{drdef2}
  Seien $\kappa \in \Na$, $\alpha \in \R$, $\sjk \in \Kpqka$ und $\ABCDsark$ das rechtsseitige $\alpha$-Dyukarev-Quadrupel bezüglich $\sjk$.
  Weiterhin sei
  \begin{align*}
    \Usarm := \begin{pmatrix} \Asarfm & \Bsarfm \\ \Csarfm & \Dsarfm \end{pmatrix}
  \end{align*}
  für alle $m \in \Zok$. Dann heißt $\Usarmk$ die \textbf{Folge von rechtsseitigen }$\mathbf{2}$\textbf{q}$\mathbf{\times2}$\textbf{q-$\alpha$-Dyukarev-Matrixpolynomen} bezüglich $\sjk$.
  Im Fall $\kappa \in \N$ heißt $\Usark$ das \textbf{rechtsseitige }$\mathbf{2}$\textbf{q}$\mathbf{\times2}$\textbf{q-$\alpha$-Dyukarev-Matrixpolynom} bezüglich $\sjk$. Falls klar ist, von welchem $\sjk$ die Rede ist, lassen wir das \anf{$\sklam{s}$} als oberen Index weg.
\end{defi}

Es sei wieder bemerkt, dass die in \thref{drdef2} eingeführten Größen mit denen in \cite[Chapter 3]{CR1} im Fall $\alpha=0$ und $\kappa=\infty$ übereinstimmen.

Unsere bereits erzielten Resultate erlauben uns eine alternative Darstellung der Folge von rechtsseitigen $2$\textit{q}$\times2$\textit{q}-$\alpha$-Dyukarev-Matrixpolynomen bezüglich einer rechtsseitig $\alpha$-Stieltjes-positiv definiten Folge (vergleiche die Teile (c) und (d) mit \cite[Proposition 3.1]{CR1} im Fall $\alpha=0$ und $\kappa=\infty$).

\begin{satz}	\thlabel{drbm5}
  Seien $\kappa \in \Na$, $\alpha \in \R$, $\sjk \in \Kpqka$ und $\Uarmk$ die Folge von rechtsseitigen $2$\textit{q}$\times2$\textit{q}-$\alpha$-Dyukarev-Matrixpolynomen bezüglich $\sjk$. Weiterhin sei
    \begin{align*}
      \tUarm(z) := \begin{pmatrix} (z-\alpha)\Iq & \Oq \\ \Oq & \Iq \end{pmatrix} \Uarm(z) \begin{pmatrix} (z-\alpha)^{-1}\Iq & \Oq \\ \Oq & \Iq \end{pmatrix}
    \end{align*}
    für alle $z\in\C\setminus\gklam{\alpha}$ und $m \in \Zok$.
  	\dgfa
  	\begin{itemize}
    	\item [\rm{(a)}] Es gelten $U_{\ar2n+1} = \Tarn$ für alle $n \in \Zofkm$ und $U_{\ar2n} = \Yarn$ für alle \linebreak $n \in \Zefk$.
    	\item [\rm{(b)}] Es gelten $\widetilde{U}_{\ar2n+1} = \tTarn$ für alle $n \in \Zofkm$ und $\widetilde{U}_{\ar2n} = \tYarn$ für alle \linebreak $n \in \Zefk$.
    	\item [\rm{(c)}] Es gelten $\Uarm,\tUarm\in\tbtJqPp$ für alle $m \in \Zok$.
    	\item [\rm{(d)}] Für alle $m\in\Zok$ sind $\det\Uarm$ und $\det\tUarm$ konstante Funktionen auf $\C$ und deren Wert jeweils von Null verschieden. Insbesondere gelten
    	\begin{align*}
    		\Uarm^{-1}(z)=\tJq\Uarm^{\ast}(\za)\tJq \quad \text{und} \quad
    		\tUarm^{-1}(z)=\tJq\tUarm^{\ast}(\za)\tJq
    	\end{align*}
    	für alle $z\in\C$ und $m \in \Zok$.
  	\end{itemize}
\end{satz}

\bwanf Zu (a): Dies folgt unter Beachtung von \thref{drdef2} und \thref{drdef1} aus Teil (a) von \thref{drsa3} bzw. Teil (a) von \thref{drsa4}.

Zu (b): Dies folgt unter Beachtung von Teil (b) von \thref{drsa3} bzw. Teil (b) von \thref{drsa4} aus (a). 

Zu (c): Wegen \thref{drdef2} und \thref{drdef1} gelten
\begin{align*}
	\Uaro(z) = \begin{pmatrix} \Iq & \Oq \\ -(z-\alpha)s^{-1}_0 & \Iq \end{pmatrix} \quad \text{und} \quad \tUaro(z) = \begin{pmatrix} \Iq & \Oq \\ -s^{-1}_0 & \Iq \end{pmatrix}
\end{align*}
für alle $z\in\C$. Hieraus folgen
\begin{align*}
	&\ \tJq-\Uaro^{\ast}(z)\tJq\Uaro(z) \\
	&= \begin{pmatrix} \Oq & -i\Iq \\ i\Iq & \Oq \end{pmatrix} - \begin{pmatrix} \Iq & -(\za-\alpha)s^{-1}_0 \\ \Oq & \Iq \end{pmatrix} \begin{pmatrix} \Oq & -i\Iq \\ i\Iq & \Oq \end{pmatrix} \begin{pmatrix} \Iq & \Oq \\ -(z-\alpha)s^{-1}_0 & \Iq \end{pmatrix} \\
	&= \begin{pmatrix} \Oq & -i\Iq \\ i\Iq & \Oq \end{pmatrix} - \begin{pmatrix} -i(\za-\alpha)s^{-1}_0 & -i\Iq \\ i\Iq & \Oq \end{pmatrix} \begin{pmatrix} \Iq & \Oq \\ -(z-\alpha)s^{-1}_0 & \Iq \end{pmatrix} \\
	&= \begin{pmatrix} \Oq & -i\Iq \\ i\Iq & \Oq \end{pmatrix} - \begin{pmatrix} -i(\za-z)s^{-1}_0 & -i\Iq \\ i\Iq & \Oq \end{pmatrix} = 2\im z\begin{pmatrix} s^{-1}_0 & \Oq \\ \Oq & \Oq \end{pmatrix}
\end{align*}
und
\begin{align*}
	&\ \tJq-\tUaro^{\ast}(z)\tJq\tUaro(z) \\
	&= \begin{pmatrix} \Oq & -i\Iq \\ i\Iq & \Oq \end{pmatrix} - \begin{pmatrix} \Iq & -s^{-1}_0 \\ \Oq & \Iq \end{pmatrix} \begin{pmatrix} \Oq & -i\Iq \\ i\Iq & \Oq \end{pmatrix} \begin{pmatrix} \Iq & \Oq \\ -s^{-1}_0 & \Iq \end{pmatrix} \\
	&= \begin{pmatrix} \Oq & -i\Iq \\ i\Iq & \Oq \end{pmatrix} - \begin{pmatrix} -is^{-1}_0 & -i\Iq \\ i\Iq & \Oq \end{pmatrix} \begin{pmatrix} \Iq & \Oq \\ -s^{-1}_0 & \Iq \end{pmatrix} \\
	&= \begin{pmatrix} \Oq & -i\Iq \\ i\Iq & \Oq \end{pmatrix} - \begin{pmatrix} \Oq & -i\Iq \\ i\Iq & \Oq \end{pmatrix} = 0_{2q\times 2q}
\end{align*}
für alle $z\in\C$. Hieraus folgt wegen \thref{spbsp1} und \thref{spdef3} nun \linebreak $\Uaro,\tUaro\in\tbtJqPp$. Unter Beachtung von Teil (d) von \thref{drsa3} bzw. Teil (d) von \thref{drsa4} folgt aus (a) und (b) weiterhin $\Uarm,\tUarm\in\tbtJqPp$ für alle $m \in \Zek$. 

Zu (d): Sei $m\in\Zok$. Wegen \thref{drdef2} und \thref{drdef1} sind dann $\Uarm$ und $\tUarm$ jeweils $2$\textit{q}$\times2$\textit{q}-Matrixpolynome. Hieraus folgt wegen \thref{spbsp1} und \thref{spfo1} dann, dass $\det\Uarm$ und $\det\tUarm$ konstante Funktionen auf $\C$ sind und deren Wert jeweils von Null verschieden ist. Wegen (c), \thref{spbsp1} und Teil (a) von \thref{splm1} gelten
\begin{align*}
	\Uarm^{-1}(z)=\tJq\Uarm^{\ast}(\za)\tJq \quad \text{und} \quad
	\tUarm^{-1}(z)=\tJq\tUarm^{\ast}(\za)\tJq
\end{align*}
für alle $z\in\C$. \bwend

Wir wollen uns nun die Determinanten der einzelnen \textit{q}$\times$\textit{q}-Matrixpolynome des rechtsseitigen $\alpha$-Dyukarev-Quadrupels bezüglich einer rechtsseitig $\alpha$-Stieltjes-positiv definiten Folge anschauen (eine entsprechende Aussage wurde in \cite{dyu} nach Formel (8) im Fall $\alpha=0$ ohne Beweis angegeben; einen Beweis, an dem wir uns auch orientieren werden, für diesen Fall formulierte Yu.\,M. Dyukarev hingegen in einer E-Mail vom 23. Oktober 2016 an B. Kirstein). Hierfür benötigen wir zunächst noch folgendes Lemma.

\begin{lemma}	\thlabel{drlm6}
	Seien $\kappa\in\Na$, $\alpha\in\R$ und $\sjk\in\Kpqka$. Weiterhin sei $\tHarn$ für alle $n\in\Zofk$ definiert wie in Teil (e) von \thref{drbm1}. \dgfa
	\begin{itemize}
		\item [\rm{(a)}] Es gelten
		\begin{align*}
			\det\eklam{H_{\arn}-(z-\alpha)H_{n}} \neq 0
		\end{align*}
		für alle $z\in\C\setminus[\alpha,\infty)$ und $n\in\Zofkm$ sowie
		\begin{align*}
			\det\beklam{T_{n}\tHarn T^{\ast}_{n}-(z-\alpha)^{-1}R^{-1}_{n}(\alpha)H_{n}R^{-\ast}_{n}(\alpha)} \neq 0
		\end{align*}
		für alle $z\in\C\setminus[\alpha,\infty)$ und $n\in\Zofk$.
		\item [\rm{(b)}] Es gelten
		\begin{align*}
			\det\eklam{\yna\eklam{H_{\arn}-(z-\alpha)H_{n}}^{-1}\yn} \neq 0
		\end{align*}
		für alle $z\in\C\setminus[\alpha,\infty)$ und $n\in\Zofkm$ sowie
		\begin{align*}
			\det\eklam{\vna\beklam{T_{n}\tHarn T^{\ast}_{n}-(z-\alpha)^{-1}R^{-1}_{n}(\alpha)H_{n}R^{-\ast}_{n}(\alpha)}^{-1}\vn} \neq 0
		\end{align*}
		für alle $z\in\C\setminus[\alpha,\infty)$ und $n\in\Zofk$.
	\end{itemize}
\end{lemma}

\bwanf Wir zeigen zunächst folgende Aussage:
\begin{itemize}
	\item [\rm{(I)}] Sei $A\in\Cpp$ mit $\im A > 0_{\pp}$ oder $-\im  A > 0_{\pp}$. Dann gilt $\det A \neq 0$.
\end{itemize}
Sei hierfür $x\in{\cal N}(A)$, d.\,h. es gilt $Ax=0_{p\times1}$. Dann gilt
\begin{align*}
	x^{\ast}[\im A]x = \im[x^{\ast}Ax] =\im[x^{\ast}0_{p\times1}] = 0.
\end{align*}
Hieraus folgt wegen $\im A > 0_{\pp}$ oder $-\im  A > 0_{\pp}$ dann $x=0_{p\times1}$, also $\det A \neq 0$.

Zu (a): Sei zunächst $n\in\Zofkm$. Es gilt $-(x-\alpha)\Hn > 0_{(n+1)q\times(n+1)q}$ für alle \linebreak $x\in(-\infty,\alpha)$. Hieraus folgt dann
\begin{align}	\label{drlm6bw1}
	\Harn-(x-\alpha)\Hn > 0_{(n+1)q\times(n+1)q}
\end{align}
für alle $x\in(-\infty,\alpha)$. Unter Beachtung von
\begin{align}	\label{drlm6bw0}
	\im \eklam{zA} = \frac{1}{2i}(zA-\za A^{\ast}) = \frac{1}{2i}(z-\za)A =[\im z] A
	\quad \text{und} \quad
	\im A = 0_{\pp}
\end{align}
für alle $A\in\Cpp_H$ gilt weiterhin
\begin{align*}
	\frac{-1}{\im z}\im\eklam{\Harn-(z-\alpha)\Hn}
	= \frac{\im(z-\alpha)}{\im z}\Hn 
	= \Hn > 0_{(n+1)q\times(n+1)q}
\end{align*}
für alle $z\in\C\setminus\R$. Hieraus folgen dann
\begin{align}	\label{drlm6bw2}
	-\im\eklam{\Harn-(z-\alpha)\Hn} > 0_{(n+1)q\times(n+1)q}
\end{align}
für alle $z\in\Pp$ und
\begin{align}	\label{drlm6bw3}
	\im\eklam{\Harn-(z-\alpha)\Hn} > 0_{(n+1)q\times(n+1)q}
\end{align}
für alle $z\in\Pm$. Wegen \fref{drlm6bw1}, \fref{drlm6bw2}, \fref{drlm6bw3} und (I) gilt dann
\begin{align*}
	\det\eklam{H_{\arn}-(z-\alpha)H_{n}} \neq 0
\end{align*}
für alle $z\in\C\setminus[\alpha,\infty)$.

Sei nun $n\in\Zofk$. Dann gelten 
\begin{align*}
	\Tn\tHarn\Tna = \begin{pmatrix} \Oqn & \Oq \\ \Harnm & \yarnm \end{pmatrix} \begin{pmatrix} \Onq & \Inq \\ \Oq & \Oqn \end{pmatrix} = \begin{pmatrix} \Oq & \Oqn \\ \Onq & \Harnm \end{pmatrix} \geq 0_{(n+1)q\times(n+1)q}
\end{align*}
und
\begin{align*}
	-(x-\alpha)^{-1}\Rn^{-1}(\alpha)\Hn\Rn^{-\ast}(\alpha) > 0_{(n+1)q\times(n+1)q}
\end{align*}
für alle $x\in(-\infty,\alpha)$. Hieraus folgt dann
\begin{align}	\label{drlm6bw4}
	\Tn\tHarn\Tna-(x-\alpha)^{-1}\Rn^{-1}(\alpha)\Hn\Rn^{-\ast}(\alpha) > 0_{(n+1)q\times(n+1)q}
\end{align}
für alle $x\in(-\infty,\alpha)$. Unter Beachtung von
\begin{align*}
	\im z^{-1} = \frac{1}{2i}\brklam{z^{-1}-\overline{z^{-1}}}
	= \frac{1}{2i}\rklam{z^{-1}\za[\za]^{-1}-z^{-1}z[\za]^{-1}}
	= z^{-1}[\za]^{-1}\frac{1}{2i}\rklam{\za-z} = \frac{-1}{\abs{z}^2}\im z
\end{align*}
für alle $z\in\C\setminus\gklam{0}$ gilt wegen $\Tn\tHarn\Tna\in\C^{(n+1)q\times(n+1)q}_H$ und \fref{drlm6bw0} weiterhin
\begin{align*}
	&\ \frac{1}{\im z}\im\beklam{\Tn\tHarn\Tna-(z-\alpha)^{-1}\Rn^{-1}(\alpha)\Hn\Rn^{-\ast}(\alpha)} \\
	&= -\frac{\im(z-\alpha)^{-1}}{\im z}\Rn^{-1}(\alpha)\Hn\Rn^{-\ast}(\alpha) \\
	&= \frac{1}{\abs{z-\alpha}^2}\Rn^{-1}(\alpha)\Hn\Rn^{-\ast}(\alpha) > 0_{(n+1)q\times(n+1)q}
\end{align*}
für alle $z\in\C\setminus\R$. Hieraus folgen dann
\begin{align}	\label{drlm6bw5}
	\im\beklam{\Tn\tHarn\Tna-(z-\alpha)^{-1}\Rn^{-1}(\alpha)\Hn\Rn^{-\ast}(\alpha)} > 0_{(n+1)q\times(n+1)q}
\end{align}
für alle $z\in\Pp$ und
\begin{align}	\label{drlm6bw6}
	-\im\beklam{\Tn\tHarn\Tna-(z-\alpha)^{-1}\Rn^{-1}(\alpha)\Hn\Rn^{-\ast}(\alpha)} > 0_{(n+1)q\times(n+1)q}
\end{align}
für alle $z\in\Pm$. Wegen \fref{drlm6bw4}-\fref{drlm6bw6} und (I) gilt dann
\begin{align*}
	\det\beklam{\Tn\tHarn\Tna-(z-\alpha)^{-1}\Rn^{-1}(\alpha)\Hn\Rn^{-\ast}(\alpha)} \neq 0
\end{align*}
für alle $z\in\C\setminus[\alpha,\infty)$. 

Zu (b): Sei zunächst $n\in\Zofkm$. Wegen $\sjk\in\Kpqma$ ist $H_0 = s_0$ positiv hermitesch. Hieraus folgt wegen \thref{asmbz1} dann
\begin{align}	\label{drlm6bw7}
	\rank \yn = q.
\end{align}
Wegen \fref{drlm6bw7} und \fref{drlm6bw1} gilt
\begin{align}	\label{drlm6bw8}
	\yna\eklam{\Harn-(x-\alpha)\Hn}^{-1}\yn > \Oq
\end{align}
für alle $x\in(-\infty,\alpha)$. Es gilt
\begin{align}	\label{drlm6bw15}
	\im A^{-1} &= \frac{1}{2i}\rklam{A^{-1}-A^{-\ast}}
	= \frac{1}{2i}\rklam{A^{-\ast}A^{\ast}A^{-1}-A^{-\ast}AA^{-1}} \notag \\
	&= A^{-\ast}\frac{1}{2i}\rklam{A^{\ast}-A} A^{-1} = A^{-\ast}(-\im A)A^{-1}
\end{align} 
für alle reguläre $A\in\Cpp$. Unter Beachtung von (a) folgt hieraus dann
\begin{align*}
	&\ \im\eklam{\yna\eklam{\Harn-(z-\alpha)\Hn}^{-1}\yn}
	=\yna\im\eklam{\Harn-(z-\alpha)\Hn}^{-1}\yn \\
	&=\yna\eklam{\Harn-(z-\alpha)\Hn}^{-\ast}\brklam{{-\im\eklam{\Harn-(z-\alpha)\Hn}}}\eklam{\Harn-(z-\alpha)\Hn}^{-1}\yn
\end{align*}
für alle $z\in\C\setminus[\alpha,\infty)$. Hieraus folgen wegen \fref{drlm6bw7}, (a) und \fref{drlm6bw2} bzw. \fref{drlm6bw3} nun
\begin{align}	\label{drlm6bw9}
	\im\eklam{\yna\eklam{\Harn-(z-\alpha)\Hn}^{-1}\yn} > \Oq
\end{align} 
für alle $z\in\Pp$ bzw.
\begin{align}	\label{drlm6bw10}
	-\im\eklam{\yna\eklam{\Harn-(z-\alpha)\Hn}^{-1}\yn} > \Oq
\end{align} 
für alle $z\in\Pm$. Wegen \fref{drlm6bw8}, \fref{drlm6bw9}, \fref{drlm6bw10} und (I) gilt dann
\begin{align*}
	\det\eklam{\yna\eklam{\Harn-(z-\alpha)\Hn}^{-1}\yn} \neq 0
\end{align*}
für alle $z\in\C\setminus[\alpha,\infty)$. 

Sei nun $n\in\Zofk$. Wegen \thref{drbz1} gilt dann
\begin{align}	\label{drlm6bw11}
	\rank \vn = q. 
\end{align}
Wegen \fref{drlm6bw11} und \fref{drlm6bw4} gilt weiterhin
\begin{align}	\label{drlm6bw12}
	\vna\beklam{\Tn\tHarn\Tna-(x-\alpha)^{-1}\Rn^{-1}(\alpha)\Hn\Rn^{-\ast}(\alpha)}^{-1}\vn > \Oq
\end{align}
für alle $x\in(-\infty,\alpha)$. Unter Beachtung von (a) gilt wegen \fref{drlm6bw15} dann
\begin{align*}
	&\ \im\eklam{\vna\beklam{\Tn\tHarn\Tna-(z-\alpha)^{-1}\Rn^{-1}(\alpha)\Hn\Rn^{-\ast}(\alpha)}^{-1}\vn} \\
	&=\vna\im\beklam{\Tn\tHarn\Tna-(z-\alpha)^{-1}\Rn^{-1}(\alpha)\Hn\Rn^{-\ast}(\alpha)}^{-1}\vn \\
	&=\vna\beklam{\Tn\tHarn\Tna-(z-\alpha)^{-1}\Rn^{-1}(\alpha)\Hn\Rn^{-\ast}(\alpha)}^{-\ast} \\
	&\quad \cdot\brklam{{-\im\beklam{\Tn\tHarn\Tna-(z-\alpha)^{-1}\Rn^{-1}(\alpha)\Hn\Rn^{-\ast}(\alpha)}}} \\
	&\quad \cdot\beklam{\Tn\tHarn\Tna-(z-\alpha)^{-1}\Rn^{-1}(\alpha)\Hn\Rn^{-\ast}(\alpha)}^{-1}\vn
\end{align*}
für alle $z\in\C\setminus[\alpha,\infty)$. Hieraus folgen wegen \fref{drlm6bw11}, (a) und \fref{drlm6bw5} bzw. \fref{drlm6bw6} nun
\begin{align}	\label{drlm6bw13}
	-\im\eklam{\vna\beklam{\Tn\tHarn\Tna-(z-\alpha)^{-1}\Rn^{-1}(\alpha)\Hn\Rn^{-\ast}(\alpha)}^{-1}\vn} > \Oq
\end{align} 
für alle $z\in\Pp$ bzw.
\begin{align}	\label{drlm6bw14}
	\im\eklam{\vna\beklam{\Tn\tHarn\Tna-(z-\alpha)^{-1}\Rn^{-1}(\alpha)\Hn\Rn^{-\ast}(\alpha)}^{-1}\vn} > \Oq
\end{align} 
für alle $z\in\Pm$. Wegen \fref{drlm6bw12}-\fref{drlm6bw14} und (I) gilt dann
\begin{align*}
	\det\eklam{\vna\beklam{\Tn\tHarn\Tna-(z-\alpha)^{-1}\Rn^{-1}(\alpha)\Hn\Rn^{-\ast}(\alpha)}^{-1}\vn} \neq 0
\end{align*}
für alle $z\in\C\setminus[\alpha,\infty)$. \bwend

\begin{satz}	\thlabel{drsa12}
	Seien $\kappa \in \Na$, $\alpha \in \R$ und $\sjm \in \Kpqka$. Weiterhin sei $\ABCDark$ das rechtsseitige $\alpha$"=Dyukarev"=Quadrupel bezüglich \linebreak $\sjk$ und $\tHarn$ für alle $n\in\Zofk$ definiert wie in Teil (e) von \thref{drbm1}. \dgfa
	\begin{itemize}
		\item [\rm{(a)}] Unter Beachtung von Teil (a) von \thref{drlm6} gelten
		\begin{align*}
			\yna\beklam{\Harn-(z-\alpha)\Hn}^{-1}\yn\Darnp(z) = \Barnp(z)
		\end{align*}
		für alle $z\in\C\setminus[\alpha,\infty)$ und $n\in\Zofkm$ sowie
		\begin{align*}
			\vna\beklam{\Tn\tHarn\Tna-(z-\alpha)^{-1}\Rn^{-1}(\alpha)\Hn\Rn^{-\ast}(\alpha)}^{-1}\vn\Aarn(z) = \Carn(z)
		\end{align*}	
		für alle $z\in\C\setminus[\alpha,\infty)$ und $n\in\Zofk$.
		
		\item [\rm{(b)}] Es gelten $\det \Aarn(z) \neq 0$ für alle $z\in\C\setminus(\alpha,\infty)$ und $n\in\Zofk$, $\det \Barn(z) \neq 0$ für alle $z\in\C\setminus(\alpha,\infty)$ und $n\in\Zefkp$, $\det \Carn(z) \neq 0$ für alle $z\in\C\setminus[\alpha,\infty)$ und $n\in\Zofk$ sowie $\det \Darn(z) \neq 0$ für alle $z\in\C\setminus(\alpha,\infty)$ und $n\in\Zofkp$.
	\end{itemize}
\end{satz}

\bwanf Zu (a): Sei zunächst $n\in\Zofkm$. Wegen Teil (c) von \thref{drbm1} gilt
\begin{align*}
	\Harn\Rn^{-\ast}(\za)-(z-\alpha)\yn\vna
	&= \Harn\Rn^{-\ast}(\za)-(z-\alpha)\beklam{\Hn\Rn^{-\ast}(\alpha)-\Harn\Tna} \\
	&= \Harn(I_{(n+1)q}-z\Tna)-(z-\alpha)\beklam{\Hn\Rn^{-\ast}(\alpha)-\Harn\Tna} \\
	&= \Harn(I_{(n+1)q}-\alpha\Tna)-(z-\alpha)\Hn\Rn^{-\ast}(\alpha) \\
	&= \beklam{\Harn-(z-\alpha)\Hn}\Rn^{-\ast}(\alpha)
\end{align*}
für alle $z\in\C$. Hieraus folgt wegen \thref{drdef1} und Teil (c) von \thref{drlm3} dann
\begin{align*}
	&\ \yna\beklam{\Harn-(z-\alpha)\Hn}^{-1}\yn\Darnp(z) \\
	&= \yna\beklam{\Harn-(z-\alpha)\Hn}^{-1}\yn\beklam{\Iq-(z-\alpha)\vna\Rna(\za)\Harn^{-1}\yn} \\
	&= \uarna\Rn^{\ast}(\alpha)\beklam{\Harn-(z-\alpha)\Hn}^{-1}\beklam{\Harn\Rn^{-\ast}(\za)-(z-\alpha)\yn\vna}\Rna(\za)\Harn^{-1}\yn \\
	&= \uarna\Rna(\za)\Harn^{-1}\yn = \Barnp(z)
\end{align*}
für alle $z\in\C\setminus[\alpha,\infty)$.

Sei nun $n\in\Zofk$. Wegen Teil (a) von \thref{drlm3} und Teil (e) von \thref{drbm1} gilt
\begin{align*}
	&\ \Rn^{-1}(\alpha)\Hn\Rn^{-\ast}(\za)+(z-\alpha)\vn\una
	= \Rn^{-1}(\alpha)\Hn\Rn^{-\ast}(\za)+(z-\alpha)\vn\yna\Tna \\
	&= \Rn^{-1}(\alpha)\Hn(I_{(n+1)q}-z\Tna)+(z-\alpha)\beklam{\Rn^{-1}(\alpha)\Hn-\Tn\tHarn}\Tna \\
	&= \Rn^{-1}(\alpha)\Hn(I_{(n+1)q}-\alpha\Tna)-(z-\alpha)\Tn\tHarn\Tna \\
	&= \Rn^{-1}(\alpha)\Hn\Rn^{-\ast}(\alpha)-(z-\alpha)\Tn\tHarn\Tna
\end{align*}
für alle $z\in\C$. Hieraus folgt wegen \thref{drdef1} dann
\begin{align*}
	&\ \vna\beklam{\Tn\tHarn\Tna-(z-\alpha)^{-1}\Rn^{-1}(\alpha)\Hn\Rn^{-\ast}(\alpha)}^{-1}\vn\Aarn(z) \\
	&= -(z-\alpha)\vna\beklam{\Rn^{-1}(\alpha)\Hn\Rn^{-\ast}(\alpha)-(z-\alpha)\Tn\tHarn\Tna}^{-1}\vn \\
	&\quad \cdot \beklam{\Iq+(z-\alpha)\una\Rna(\za)\Hn^{-1}\Rn(\alpha)\vn} \\
	&= -(z-\alpha)\vna\beklam{\Rn^{-1}(\alpha)\Hn\Rn^{-\ast}(\alpha)-(z-\alpha)\Tn\tHarn\Tna}^{-1} \\
	&\quad \cdot \beklam{\Rn^{-1}(\alpha)\Hn\Rn^{-\ast}(\za)+(z-\alpha)\vn\una}\Rna(\za)\Hn^{-1}\Rn(\alpha)\vn \\
	&= -(z-\alpha)\vna\Rna(\za)\Hn^{-1}\Rn(\alpha)\vn = \Carn(z)
\end{align*}
für alle $z\in\C\setminus[\alpha,\infty)$.

Zu (b): Seien $z\in\C\setminus[\alpha,\infty)$ und zunächst $n\in\Zofk$. Wegen (a) und Teil (b) von \thref{drlm6} gilt dann
\begin{align}	\label{drsa12bw1}
	{\cal N}\brklam{\Aarn(z)} = {\cal N}\brklam{\Carn(z)}.
\end{align}
Sei nun $x\in{\cal N}\brklam{\Aarn(z)}$. Wegen \thref{drdef2} und \fref{drsa12bw1} gilt dann
\begin{align}	\label{drsa12bw2}
	U_{\ar2n}\binom{x}{0_{q\times1}} = \begin{pmatrix} \Aarn(z) & \Barn(z) \\ \Carn(z) & \Darn(z) \end{pmatrix}\binom{x}{0_{q\times1}} = 0_{2q\times 1}.
\end{align}
Wegen Teil (d) von \thref{drbm5} gilt ${\cal N}\brklam{U_{\ar2n}(z)} = \gklam{0_{2q\times1}}$. Hieraus folgt wegen \fref{drsa12bw2} nun $x=0_{q\times1}$, also wegen \fref{drsa12bw1} gelten dann $\det\Aarn(z)\neq0$ und $\det\Carn(z)\neq0$.

Wegen \thref{drdef1} gilt $\Daro \equiv \Iq$ und somit $\det\Daro(z)\neq0$. Sei nun $n\in\Zofkm$. Wegen (a) und Teil (b) von \thref{drlm6} gilt dann
\begin{align}	\label{drsa12bw3}
	{\cal N}\brklam{\Barnp(z)} = {\cal N}\brklam{\Darnp(z)}.
\end{align}
Sei nun $x\in{\cal N}\brklam{\Barnp(z)}$. Wegen \thref{drdef2} und \fref{drsa12bw3} gilt dann
\begin{align}	\label{drsa12bw4}
	U_{\ar2n+1}\binom{0_{q\times1}}{x} = \begin{pmatrix} \Aarn(z) & \Barnp(z) \\ \Carn(z) & \Darnp(z) \end{pmatrix}\binom{0_{q\times1}}{x} = 0_{2q\times 1}.
\end{align}
Wegen Teil (d) von \thref{drbm5} gilt ${\cal N}\brklam{U_{\ar2n+1}(z)} = \gklam{0_{2q\times1}}$. Hieraus folgt wegen \fref{drsa12bw3} nun $x=0_{q\times1}$, also wegen \fref{drsa12bw3} gelten dann $\det\Barnp(z)\neq0$ und $\det\Darnp(z)\neq0$. 

Wegen \thref{drdef1} gelten $\Aarn(\alpha) = \Iq$ für alle $n\in\Zofk$ und $\Darn(\alpha)=\Iq$ für alle $n\in\Zofkp$. Hieraus folgen $\det\Aarn(\alpha) \neq 0$ für alle $n\in\Zofk$ und $\det\Darn(\alpha) \neq 0$ für alle $n\in\Zofkp$. Wegen \thref{drdef1} und Teil (e) von \thref{drlm3} gilt weiterhin
\begin{align}	\label{drsa12bw5}
	\Barnp(\alpha) = \uarna\Rna(\alpha)\Harn^{-1}\yn = \yna\Harn^{-1}\yn
\end{align}
für alle $n\in\Zofkm$. Wegen $\sjk\in\Kpqka$ sind $H_0 = s_0$ und $\Harn$ für alle $n\in\Zofkm$ positiv hermitesche Matrizen. Hieraus folgen wegen \thref{asmbz1} dann $\rank\yn=q$ für alle $n\in\Zofkm$ und wegen \fref{drsa12bw5} weiterhin, dass $\Barn(\alpha)$ für alle $n\in\Zefkp$ eine positiv hermitesche und somit reguläre Matrix ist. \bwend

\thref{drbm5} erlaubt uns, das rechtsseitige $2$\textit{q}$\times2$\textit{q}-$\alpha$-Dyukarev-Matrixpolynom bezüglich einer endlichen rechtsseitig $\alpha$-Stieltjes-positiv definiten Folge mit den Potapovschen Fundamentalmatrizen in Verbindung zu bringen. Es wird sich herausstellen, dass das rechtsseitige $2$\textit{q}$\times2$\textit{q}-$\alpha$-Dyukarev-Matrixpolynom bezüglich jener Folge auch eine Resolventenmatrix für das zugehörige rechtsseitige $\alpha$-Stieltjes Momentenproblem ist, wie wir im Folgenden zeigen werden (vergleiche \cite[Theorem 3.2]{CR1} im Fall $\alpha=0$). Wir orientieren uns dafür weiterhin an der Vorgehensweise von \cite[Chapter 6]{C06} für das Hausdorffsche Momentenproblem.

\begin{lemma}	\thlabel{drlm4}
  Seien $m \in \N$, $\alpha \in \R$, $\sjm \in \Kpqma$ und $\Uarm$ das rechtsseitige $2$\textit{q}$\times2$\textit{q}-$\alpha$-Dyukarev-Matrixpolynom bezüglich $\sjm$.
  Es bezeichne
  \begin{align*}
    \Uarm = \begin{pmatrix} \Uarm^{(1,1)} & \Uarm^{(1,2)} \\ \Uarm^{(2,1)} & \Uarm^{(2,2)} \end{pmatrix}
  \end{align*}
  die \textit{q}$\times$\textit{q}-Blockzerlegung von $\Uarm$. Weiterhin seien $\phipsi \in \PtJqCa$ sowie
  \begin{align*}
    \dphi := \Uarm^{(1,1)}\phi+\Uarm^{(1,2)}\psi \quad \text{und} \quad \dpsi := \Uarm^{(2,1)}\phi+\Uarm^{(2,2)}\psi.
  \end{align*}
  \dgfa
  \begin{itemize}
    \item [\rm{(a)}] Es gilt $\dphipsi \in \PtJqCa$.
    \item [\rm{(b)}] Es sind $\det\dphi$ und $\det\dpsi$ jeweils nicht die Nullfunktion.
  \end{itemize}
\end{lemma}

\bwanf Zu (a): Es gilt
\begin{align}	\label{drlm4bw1}
  \Uarm\phipsi = \begin{pmatrix}\dphi\\\dpsi\end{pmatrix}.
\end{align}
Wegen Teil (a) von \thref{spdef4} existiert eine diskrete Teilmenge $\D$ von $\C\setminus[\alpha,\infty)$ mit
\begin{itemize}
  \item [\rm{(i)}] $\phi$ und $\psi$ sind in $\C\setminus\rklam{[\alpha,\infty)\cup\D}$ holomorph.
  \item [\rm{(ii)}] Es gilt $\rank \phipsiz = q$ für alle $z \in \C\setminus\rklam{[\alpha,\infty)\cup\D}$.
  \item [\rm{(iii)}] Es gelten
  \begin{align*}
    \binom{(z-\alpha)\phi(z)}{\psi(z)}^{\ast}\bbrklam{\frac{-\tJq}{2\im z}}\binom{(z-\alpha)\phi(z)}{\psi(z)} \geq \Oq
  \end{align*}
  und
  \begin{align*}
    \phipsiz^{\ast}\bbrklam{\frac{-\tJq}{2\im z}}\phipsiz \geq \Oq
  \end{align*}
  für alle $z \in \C\setminus\rklam{\R\cup\D}$.
\end{itemize}
Wegen \thref{spbsp1}, Teil (c) von \thref{drbm5} und \thref{spdef3} gilt, dass $\Uarm$ in $\C$ holomorph ist.
Hieraus folgt wegen \fref{drlm4bw1} und (i) dann
\begin{itemize}
  \item [\rm{(iv)}] $\dphi$ und $\dpsi$ sind in $\C\setminus\rklam{[\alpha,\infty)\cup\D}$ holomorph.
\end{itemize}
Wegen Teil (d) von \thref{drbm5}, \fref{drlm4bw1} und (ii) gilt
\begin{itemize}
 \item [\rm{(v)}] Es gilt $\rank \bbinom{\dphi(z)}{\dpsi(z)} = q$ für alle $z \in \C\setminus\rklam{[\alpha,\infty)\cup\D}$.
\end{itemize}
Sei $\tUarm$ definiert wie in \thref{drbm5}. Wegen \fref{drlm4bw1} gilt dann
\begin{align}	\label{drlm4bw2}
  \tUarm(z)\binom{(z-\alpha)\phi(z)}{\psi(z)} = \begin{pmatrix}(z-\alpha)\dphi(z)\\\dpsi(z)\end{pmatrix}
\end{align}
für alle $z\in\C\setminus[\alpha,\infty)$. 
Unter Beachtung von \thref{spbsp1} und Teil (c) von \thref{drbm5} gilt wegen \fref{drlm4bw2}, Teil (d) von \thref{splm1} und (iii) dann
\begin{align}	\label{drlm4bw3}
  &\ \begin{pmatrix}(z-\alpha)\dphi(z)\\\dpsi(z)\end{pmatrix}^{\ast}\bbrklam{\frac{-\tJq}{2\im z}}\begin{pmatrix}(z-\alpha)\dphi(z)\\\dpsi(z)\end{pmatrix} \notag \\
  & = \binom{(z-\alpha)\phi(z)}{\psi(z)}^{\ast}\tUarm^{\ast}\bbrklam{\frac{-\tJq}{2\im z}}\tUarm\binom{(z-\alpha)\phi(z)}{\psi(z)} \notag \\
  & \geq \binom{(z-\alpha)\phi(z)}{\psi(z)}^{\ast}\bbrklam{\frac{-\tJq}{2\im z}}\binom{(z-\alpha)\phi(z)}{\psi(z)} \geq \Oq
\end{align}
für alle $z \in \C\setminus\rklam{\R\cup\D}$.
Unter Beachtung von \thref{spbsp1} und Teil (c) von \thref{drbm5} gilt wegen \fref{drlm4bw1}, Teil (d) von \thref{splm1} und (iii) weiterhin
\begin{align}	\label{drlm4bw4}
  &\ \dphipsiz^{\ast}\bbrklam{\frac{-\tJq}{2\im z}}\dphipsiz \notag \\
  & = \phipsiz^{\ast}\Uarm^{\ast}\bbrklam{\frac{-\tJq}{2\im z}}\Uarm\phipsiz \notag \\
  & \geq \phipsiz^{\ast}\bbrklam{\frac{-\tJq}{2\im z}}\phipsiz \geq \Oq
\end{align}
für alle $z \in \C\setminus\rklam{\R\cup\D}$.
Wegen (iv), (v), \fref{drlm4bw3}, \fref{drlm4bw4} und Teil (a) von \thref{spdef4} gilt nun $\dphipsi \in \PtJqCa$.

Zu (b): Sei $z \in \C\setminus\rklam{\R\cup\D}$. Wegen Teil (d) von \thref{drbm5}, \fref{drlm4bw1} und (iii) gilt
\begin{align}	\label{drlm4bw5}
  \dphipsiz^{\ast}\Uarm^{-\ast}\bbrklam{\frac{-\tJq}{2\im z}}\Uarm^{-1}\dphipsiz = \phipsiz^{\ast}\bbrklam{\frac{-\tJq}{2\im z}}\phipsiz \geq \Oq.
\end{align}
Sei zunächst $y \in {\cal N}\brklam{\dphi(z)}$. Wegen \fref{drlm4bw5} gilt dann
\begin{align}	\label{drlm4bw6}
  \binom{\Oq}{\dpsi(z)y}^{\ast}\Uarm^{-\ast}\bbrklam{\frac{-\tJq}{2\im z}}\Uarm^{-1}\binom{\Oq}{\dpsi(z)y} \geq \Oq.
\end{align}
Es gilt
\begin{align*}
   \binom{\Oq}{\dpsi(z)y}^{\ast}\tJq\binom{\Oq}{\dpsi(z)y} = \binom{i\dpsi(z)y}{\Oq}^{\ast}\binom{\Oq}{\dpsi(z)y} = \Oq.
\end{align*}
Hieraus folgt wegen $2\im z = -i(z-\za)$ und \fref{drlm4bw6} dann
\begin{align}	\label{drlm4bw7}
  \binom{\Oq}{\dpsi(z)y}^{\ast}\frac{\tJq-\Uarm^{-\ast}\tJq\Uarm^{-1}}{i(z-\za)}\binom{\Oq}{\dpsi(z)y} \leq \Oq.
\end{align}
Unter Beachtung von Teil (a) von \thref{drbm5} gilt wegen Teil (e) von \thref{drsa3} bzw. Teil (e) von \thref{drsa4} und Teil (d) von \thref{drsa2} nun
\begin{align*}
  \frac{\tJq-\Uarm^{-\ast}\tJq\Uarm^{-1}}{i(z-\za)} = \tJq\binom{u^{\ast}_{\fklam{m}}}{-v^{\ast}_{\fklam{m}}}R^{\ast}_{\fklam{m}}(z)H^{-1}_{\fklam{m}}R_{\fklam{m}}(z)\begin{pmatrix} u_{\fklam{m}} & -v_{\fklam{m}} \end{pmatrix}\tJq.
\end{align*}
Hieraus folgt wegen $H^{-1}_{\fklam{m}} \in \C^{(n+1)q\times(n+1)q}_>$ und \fref{drlm4bw7} dann
\begin{align*}
  &\ \binom{\Oq}{\dpsi(z)y}^{\ast}\tJq\binom{u^{\ast}_{\fklam{m}}}{-v^{\ast}_{\fklam{m}}}R^{\ast}_{\fklam{m}}(z)H^{-1}_{\fklam{m}}R_{\fklam{m}}(z)\begin{pmatrix} u_{\fklam{m}} & -v_{\fklam{m}} \end{pmatrix}\tJq\binom{\Oq}{\dpsi(z)y}  \\
  &= \binom{\Oq}{\dpsi(z)y}^{\ast}\begin{pmatrix} iv_{\fklam{m}} & iu_{\fklam{m}} \end{pmatrix}^{\ast}R^{\ast}_{\fklam{m}}(z)H^{-1}_{\fklam{m}}R_{\fklam{m}}(z)\begin{pmatrix} iv_{\fklam{m}} & iu_{\fklam{m}} \end{pmatrix}\binom{\Oq}{\dpsi(z)y}
  = \Oq.
\end{align*}
Hieraus folgt wegen $H^{-1}_{\fklam{m}} \in \C^{(n+1)q\times(n+1)q}_>$ und $\det R_{\fklam{m}}(z) \neq 0$ nun
\begin{align*}
  0_{(\fklam{m}+1)q\times q} = \begin{pmatrix} iv_{\fklam{m}} & iu_{\fklam{m}} \end{pmatrix}\binom{\Oq}{\dpsi(z)y} = iu_{\fklam{m}}\dpsi(z)y.
\end{align*}
Hieraus folgt $is_0\dpsi(z)y = \Oq$. Wegen \thref{adpbm1} ist $s_0 = H_0$ regulär, also ist $y \in {\cal N}\brklam{\dpsi(z)}$, 
das heißt $\bbinom{\dphi(z)}{\dpsi(z)} y = 0_{2q\times q}$ und somit wegen (v) dann $y = 0_{q\times1}$. 
Da $z \in \C\setminus\rklam{\R\cup\D}$ beliebig gewählt ist, ist $\det\dphi$ nicht die Nullfunktion.

Sei nun $y \in {\cal N}\brklam{\dpsi(z)}$. Wegen \fref{drlm4bw5} gilt dann
\begin{align}	\label{drlm4bw8}
  \binom{\dphi(z)y}{\Oq}^{\ast}\Uarm^{-\ast}\bbrklam{\frac{-\tJq}{2\im z}}\Uarm^{-1}\binom{\dphi(z)y}{\Oq} \geq \Oq.
\end{align}
Es gilt
\begin{align*}
   \binom{\dphi(z)y}{\Oq}^{\ast}\tJq\binom{\dphi(z)y}{\Oq} = \binom{\Oq}{-i\dphi(z)}^{\ast}\binom{\dphi(z)y}{\Oq} = \Oq.
\end{align*}
Hieraus folgt wegen $2\im z = -i(z-\za)$ und \fref{drlm4bw8} dann
\begin{align}	\label{drlm4bw9}
  \binom{\dphi(z)y}{\Oq}^{\ast}\frac{\tJq-\Uarm^{-\ast}\tJq\Uarm^{-1}}{i(z-\za)}\binom{\dphi(z)y}{\Oq} \leq \Oq.
\end{align}
Unter Beachtung von Teil (a) von \thref{drbm5} gilt wegen Teil (e) von \thref{drsa3} bzw. Teil (e) von \thref{drsa4} und Teil (d) von \thref{drsa2}
\begin{align*}
  \frac{\tJq-\Uarm^{-\ast}\tJq\Uarm^{-1}}{i(z-\za)} = \tJq\binom{u^{\ast}_{\fklam{m}}}{-v^{\ast}_{\fklam{m}}}R^{\ast}_{\fklam{m}}(z)H^{-1}_{\fklam{m}}R_{\fklam{m}}(z)\begin{pmatrix} u_{\fklam{m}} & -v_{\fklam{m}} \end{pmatrix}\tJq.
\end{align*}
Hieraus folgt wegen $H^{-1}_{\fklam{m}} \in \C^{(n+1)q\times(n+1)q}_>$ und \fref{drlm4bw9} dann
\begin{align*}
  &\ \binom{\dphi(z)y}{\Oq}^{\ast}\tJq\binom{u^{\ast}_{\fklam{m}}}{-v^{\ast}_{\fklam{m}}}R^{\ast}_{\fklam{m}}(z)H^{-1}_{\fklam{m}}R_{\fklam{m}}(z)\begin{pmatrix} u_{\fklam{m}} & -v_{\fklam{m}} \end{pmatrix}\tJq\binom{\dphi(z)y}{\Oq} \\
  &= \binom{\dphi(z)y}{\Oq}^{\ast}\begin{pmatrix} iv_{\fklam{m}} & iu_{\fklam{m}} \end{pmatrix}^{\ast}R^{\ast}_{\fklam{m}}(z)H^{-1}_{\fklam{m}}R_{\fklam{m}}(z)\begin{pmatrix} iv_{\fklam{m}} & iu_{\fklam{m}} \end{pmatrix}\binom{\dphi(z)y}{\Oq} = \Oq.
\end{align*}
Wegen $\det\Rm(z) \neq 0$ gilt nun
\begin{align*}
  0_{(\fklam{m}+1)q\times q} = \begin{pmatrix} iv_{\fklam{m}} & iu_{\fklam{m}} \end{pmatrix}\binom{\dphi(z)y}{\Oq} = iv_{\fklam{m}}\dphi(z)y.
\end{align*}
Hieraus folgt $i\dphi(z)y = \Oq$, also ist $y \in {\cal N}\brklam{\dpsi(z)}$, das heißt $\bbinom{\dphi(z)}{\dpsi(z)} y = 0_{2q\times q}$ und somit wegen (v) dann $y = 0_{q\times1}$. 
Da $z \in \C\setminus\rklam{\R\cup\D}$ beliebig gewählt ist, ist $\det\dphi$ nicht die Nullfunktion. \bwend

Es sei bemerkt, dass man Teil (b) von \thref{drlm4} in ähnlicher und allgemeinerer Form für ungerade $m$ in \cite[Lemma 10.22]{Maka} findet.

\begin{theo}	\thlabel{drth1}
  Seien $m \in \N$, $\alpha \in \R$, $\sjm \in \Kpqma$ und $\Uarm$ das rechtsseitige $2$\textit{q}$\times2$\textit{q}-$\alpha$-Dyukarev-Matrixpolynom bezüglich $\sjm$.
  Es bezeichne
  \begin{align*}
    \Uarm = \begin{pmatrix} \Uarm^{(1,1)} & \Uarm^{(1,2)} \\ \Uarm^{(2,1)} & \Uarm^{(2,2)} \end{pmatrix}
  \end{align*}
  die \textit{q}$\times$\textit{q}-Blockzerlegung von $\Uarm$. \dgfa
  \begin{itemize}
    \item [\rm{(a)}] Sei $\phipsi \in \PtJqCa$. Dann ist $\det\big(\Uarm^{(2,1)}\phi+\Uarm^{(2,2)}\psi\big)$ nicht die Nullfunktion und
    \begin{align*}
      S := \rklam{\Uarm^{(1,1)}\phi+\Uarm^{(1,2)}\psi}\rklam{\Uarm^{(2,1)}\phi+\Uarm^{(2,2)}\psi}^{-1}
    \end{align*}
    gehört zu $\Sqasmu$.
    \item [\rm{(b)}] Sei $S\in\Sqasmu$. Dann existiert ein $\phipsi \in \dPtJqCa$ derart, dass
    \begin{align*}
      \det\eklam{\Uarm^{(2,1)}(z)\phi(z)+\Uarm^{(2,2)}(z)\psi(z)} \neq 0
    \end{align*}
    und
    \begin{align*}
      S(z) = \eklam{\Uarm^{(1,1)}(z)\phi(z)+\Uarm^{(1,2)}(z)\psi(z)}\eklam{\Uarm^{(2,1)}(z)\phi(z)+\Uarm^{(2,2)}(z)\psi(z)}^{-1}
    \end{align*}
    für alle $z \in \C\setminus[\alpha,\infty)$ erfüllt sind.
    \item [\rm{(c)}] Seien $\binom{\phi_1}{\psi_1}, \binom{\phi_2}{\psi_2} \in \PtJqCa$. 
    Dann sind $\det\big(\Uarm^{(2,1)}\phi_1+\Uarm^{(2,2)}\psi_1\big)$ und \linebreak $\det\big(\Uarm^{(2,1)}\phi_2+\Uarm^{(2,2)}\psi_2\big)$ jeweils nicht die Nullfunktion und es sind folgende Aussagen äquivalent:
    \begin{itemize}
      \item [\rm{(i)}] Es gilt
      \begin{align*}
	&\ \rklam{\Uarm^{(1,1)}\phi_1+\Uarm^{(1,2)}\psi_1}\rklam{\Uarm^{(2,1)}\phi_1+\Uarm^{(2,2)}\psi_1}^{-1} \\
	&= \rklam{\Uarm^{(1,1)}\phi_2+\Uarm^{(1,2)}\psi_2}\rklam{\Uarm^{(2,1)}\phi_2+\Uarm^{(2,2)}\psi_2}^{-1}.
      \end{align*}
      \item [\rm{(ii)}] Es gilt $\sklam{\binom{\phi_1}{\psi_1}} = \sklam{\binom{\phi_2}{\psi_2}}$.
    \end{itemize}
  \end{itemize}
\end{theo}

\bwanf Zu (a): Seien
\begin{align*}
  \dphi := \Uarm^{(1,1)}\phi+\Uarm^{(1,2)}\psi \quad \text{und} \quad \dpsi := \Uarm^{(2,1)}\phi+\Uarm^{(2,2)}\psi.
\end{align*}
Wegen \thref{drlm4} gilt 
\begin{itemize}
  \item [\rm{(iii)}] Es gilt $\dphipsi \in \PtJqCa$ und $\det\dpsi$ ist nicht die Nullfunktion. 
\end{itemize}  
Weiterhin gelten
\begin{align}	\label{drth1bw1}
  \Uarm\phipsi = \begin{pmatrix}\dphi\\\dpsi\end{pmatrix}
\end{align}
und
\begin{align}	\label{drth1bw2}
  S = \dphi\cdot\dpsi^{-1}.
\end{align}
Wegen Teil (a) von \thref{spdef4}, (iii) und \fref{drth1bw2} existiert dann eine diskrete Teilmenge $\D$ von $\C\setminus[\alpha,\infty)$ 
(als Vereinigung endlich vieler diskreter Teilmengen von $\C\setminus[\alpha,\infty)$) mit
\begin{itemize}
  \item [\rm{(iv)}] Es sind $S, \phi, \psi, \dphi$ und $\dpsi$ in $\C\setminus([\alpha,\infty)\cup\D)$ holomorph. 
  \item [\rm{(v)}] Es gilt $\det\dpsi(z) \neq 0$ für alle $z \in \C\setminus([\alpha,\infty)\cup\D)$.
  \item [\rm{(vi)}] Es gelten
  \begin{align*}
    \binom{(z-\alpha)\phi(z)}{\psi(z)}^{\ast}\bbrklam{\frac{-\tJq}{2\im z}}\binom{(z-\alpha)\phi(z)}{\psi(z)} \geq \Oq
  \end{align*}
  und
  \begin{align*}
    \phipsiz^{\ast}\bbrklam{\frac{-\tJq}{2\im z}}\phipsiz \geq \Oq
  \end{align*}
  für alle $z \in \C\setminus\rklam{\R\cup\D}$.  
\end{itemize}  
Wegen Teil (d) von \thref{drbm5} gilt, dass $\Uarm(z)$ für alle $z \in \C$ regulär ist.
Hieraus folgt wegen (v), \fref{drth1bw1} und \fref{drth1bw2} dann
\begin{align}	\label{drth1bw3}
  \phipsiz\dpsi^{-1}(z) = \Uarm^{-1}\dphipsiz\dpsi^{-1}(z) = \Uarm^{-1}\binom{S(z)}{\Iq}
\end{align}
für alle $z \in \C\setminus([\alpha,\infty)\cup\D)$.
Sei $\tUarm$ definiert wie in \thref{drbm5}. Dann gilt
\begin{align}	\label{drth1bw7}
  \tUarm(z)\binom{(z-\alpha)\phi(z)}{\psi(z)} = \begin{pmatrix}(z-\alpha)\dphi(z)\\\dpsi(z)\end{pmatrix}
\end{align}
für alle $z \in \C\setminus[\alpha,\infty)$.
Wegen Teil (d) von \thref{drbm5} gilt, dass $\tUarm(z)$ für alle $z \in \C$ regulär ist.
Hieraus folgt wegen (v), \fref{drth1bw7} und \fref{drth1bw2} dann
\begin{align}	\label{drth1bw4}
  \binom{(z-\alpha)\phi(z)}{\psi(z)}\dpsi^{-1}(z) = \tUarm^{-1}(z)\begin{pmatrix}(z-\alpha)\dphi(z)\\\dpsi(z)\end{pmatrix}\dpsi^{-1}(z) = \tUarm^{-1}(z)\binom{(z-\alpha)S(z)}{\Iq}
\end{align}
für alle $z \in \C\setminus([\alpha,\infty)\cup\D)$.
Wegen \fref{drth1bw3} und (vi) gilt
\begin{align}	\label{drth1bw5}
  &\ \binom{S(z)}{\Iq}^{\ast}\Uarm^{-\ast}\bbrklam{\frac{-\tJq}{2\im z}}\Uarm^{-1}\binom{S(z)}{\Iq} \notag \\
  & = \dpsi^{-\ast}(z)\phipsiz^{\ast}\bbrklam{\frac{-\tJq}{2\im z}}\phipsiz\dpsi^{-1}(z) \geq \Oq
\end{align}
für alle $z \in \C\setminus\rklam{\R\cup\D}$.
Wegen \fref{drth1bw4} und (vi) gilt weiterhin
\begin{align}	\label{drth1bw6}
  &\ \binom{(z-\alpha)S(z)}{\Iq}^{\ast}\tUarm^{-\ast}\bbrklam{\frac{-\tJq}{2\im z}}\tUarm^{-1}\binom{(z-\alpha)S(z)}{\Iq} \notag \\
  & = \dpsi^{-\ast}(z)\binom{(z-\alpha)\phi(z)}{\psi(z)}^{\ast}\bbrklam{\frac{-\tJq}{2\im z}}\binom{(z-\alpha)\phi(z)}{\psi(z)}\dpsi^{-1}(z) \geq \Oq
\end{align}
für alle $z \in \C\setminus\rklam{\R\cup\D}$.
Unter Beachtung von den Teilen (a) und (b) von \thref{drbm5} gilt wegen \fref{drth1bw5}, \fref{drth1bw6}, \thref{drbm3} und Teil (f) von \thref{drsa3} bzw. Teil (f) von \thref{drsa4} dann 
$\FSfm(z) \in \C^{(\fklam{m}+2)q\times (\fklam{m}+2)q}_{\geq}$ und $\FSarfm(z) \in \C^{(\fklam{m-1}+2)q\times (\fklam{m-1}+2)q}_{\geq}$ für alle $z \in \Pp\setminus\D$.
Wegen \thref{drsa5} gibt es genau ein $\widetilde{S} \in \Sqasmu$ mit $\Rstr_{\Pp\setminus\D}\widetilde{S} = \Rstr_{\Pp\setminus\D}S$.
Da $\widetilde{S}$ in $\C\setminus[\alpha,\infty)$ holomorph ist, gilt wegen (iv) und des Identitätssatzes für meromorphe Funktionen (vergleiche z.\,B. im skalaren Fall \cite[Satz 10.3.2]{Funk}; im matriziellen Fall betrachtet man die einzelnen Einträge der Matrixfunktion) dann 
$\widetilde{S} = S$ und somit $S \in \Sqasmu$.

Zu (b): Wegen Teil (d) von \thref{drbm5} gilt, dass $\Uarm(z)$ für alle $z \in \C$ regulär ist.
Sei
\begin{align*}
  \phipsi := \Uarm^{-1}\binom{S}{\Iq}.
\end{align*}
Wegen $S \in \Soqa$ (vergleiche \thref{asmbz3}) gilt dann
\begin{itemize}
  \item [\rm{(vii)}] Es sind $\phi$ und $\psi$ in $\C\setminus[\alpha,\infty)$ holomorph.
\end{itemize}
Wegen $\rank\Iq = q$ gilt weiterhin
\begin{align}	\label{drth1bw8}
  \rank\phipsiz = \rank\binom{S(z)}{\Iq} = q
\end{align}
für alle $z \in \C\setminus[\alpha,\infty)$.
Unter Beachtung von Teil (a) von \thref{drbm5} gilt wegen Teil (f) von \thref{drsa3} bzw. Teil (f) von \thref{drsa4}, Teil (a) von \thref{drbm3} und \thref{drsa6} dann
\begin{align}	\label{drth1bw9}
  &\ \phipsiz^{\ast}\bbrklam{\frac{-\tJq}{2\im z}}\phipsiz \notag \\
  & = \binom{S}{\Iq}^{\ast}\Uarm^{-\ast}\bbrklam{\frac{-\tJq}{2\im z}}\Uarm^{-1}\binom{S}{\Iq} \geq \Oq
\end{align}
für alle $z \in \C\setminus\R$.
Sei $\tUarm$ definiert wie in \thref{drbm5}. 
Wegen Teil (d) von \thref{drbm5} gilt, dass $\tUarm(z)$ für alle $z \in \C$ regulär ist.
Dann gilt
\begin{align*}
  \binom{(z-\alpha)\phi(z)}{\psi(z)} = \tUarm^{-1}(z)\binom{(z-\alpha)S(z)}{\Iq}
\end{align*}
für alle $z \in \C\setminus[\alpha,\infty)$. 
Hieraus folgt unter Beachtung von Teil (b) von \thref{drbm5} wegen Teil (f) von \thref{drsa3} bzw. Teil (f) von \thref{drsa4}, Teil (b) von \thref{drbm3} und \thref{drsa6} weiterhin
\begin{align}	\label{drth1bw10}
  &\ \binom{(z-\alpha)\phi(z)}{\psi(z)}^{\ast}\bbrklam{\frac{-\tJq}{2\im z}}\binom{(z-\alpha)\phi(z)}{\psi(z)} \notag \\
  & = \binom{(z-\alpha)S(z)}{\Iq}^{\ast}\tUarm^{-\ast}\bbrklam{\frac{-\tJq}{2\im z}}\tUarm^{-1}\binom{(z-\alpha)S(z)}{\Iq} \geq \Oq
\end{align}
für alle $z \in \C\setminus\R$. 
Wegen (vii), \fref{drth1bw8}, \fref{drth1bw9}, \fref{drth1bw10} und Teil (a) von \thref{spdef4} gilt nun $\phipsi \in \dPtJqCa$.
Insbesondere gilt
\begin{align*}
  \Uarm\phipsi = \binom{S}{\Iq}.
\end{align*}
Hieraus folgen
\begin{align*}
  \det\eklam{\Uarm^{(2,1)}(z)\phi(z)+\Uarm^{(2,2)}(z)\psi(z)} = \det\Iq \neq 0
\end{align*}
und
\begin{align*}
  S = S\cdot\Iq^{-1} = \eklam{\Uarm^{(1,1)}(z)\phi(z)+\Uarm^{(1,2)}(z)\psi(z)}\eklam{\Uarm^{(2,1)}(z)\phi(z)+\Uarm^{(2,2)}(z)\psi(z)}^{-1}
\end{align*}
für alle $z \in \C\setminus[\alpha,\infty)$.

Zu (c): Wegen (a) sind $\det\big(\Uarm^{(2,1)}\phi_1+\Uarm^{(2,2)}\psi_1\big)$ und $\det\big(\Uarm^{(2,1)}\phi_2+\Uarm^{(2,2)}\psi_2\big)$ jeweils nicht die Nullfunktion.

Sei zunächst (i) erfüllt. Wegen Teil (d) von \thref{drbm5} gilt, dass $\Uarm(z)$ für alle $z \in \C$ regulär ist.
Hieraus folgt nun
\begin{align*}
  \binom{\phi_j}{\psi_j} & = \Uarm^{-1}\Uarm\binom{\phi_j}{\psi_j} \\
  & = \Uarm^{-1}\binom{\Uarm^{(1,1)}\phi_j+\Uarm^{(1,2)}\psi_j}{\Uarm^{(2,1)}\phi_j+\Uarm^{(2,2)}\psi_j} \\
  & = \Uarm^{-1}\binom{\big(\Uarm^{(1,1)}\phi_j+\Uarm^{(1,2)}\psi_j\big)\big(\Uarm^{(2,1)}\phi_j+\Uarm^{(2,2)}\psi_j\big)^{-1}}{\Iq}\rklam{\Uarm^{(2,1)}\phi_j+\Uarm^{(2,2)}\psi_j}
\end{align*}
für alle $j \in \gklam{1,2}$. Hieraus folgt wegen (i) dann
\begin{align*}
  \binom{\phi_2}{\psi_2} & = \Uarm^{-1}\binom{\big(\Uarm^{(1,1)}\phi_2+\Uarm^{(1,2)}\psi_2\big)\big(\Uarm^{(2,1)}\phi_2+\Uarm^{(2,2)}\psi_2\big)^{-1}}{\Iq}\rklam{\Uarm^{(2,1)}\phi_2+\Uarm^{(2,2)}\psi_2} \\
  & = \Uarm^{-1}\binom{\big(\Uarm^{(1,1)}\phi_1+\Uarm^{(1,2)}\psi_1\big)\big(\Uarm^{(2,1)}\phi_1+\Uarm^{(2,2)}\psi_1)^{-1}}{\Iq}\rklam{\Uarm^{(2,1)}\phi_2+\Uarm^{(2,2)}\psi_2} \\
  & = \Uarm^{-1}\binom{\Uarm^{(1,1)}\phi_1+\Uarm^{(1,2)}\psi_1}{\Uarm^{(2,1)}\phi_1+\Uarm^{(2,2)}\psi_1}\rklam{\Uarm^{(2,1)}\phi_1+\Uarm^{(2,2)}\psi_1}^{-1}\rklam{\Uarm^{(2,1)}\phi_2+\Uarm^{(2,2)}\psi_2} \\ 
  & = \binom{\phi_1}{\psi_1}(\Uarm^{(2,1)}\phi_1+\Uarm^{(2,2)}\psi_1)^{-1}\rklam{\Uarm^{(2,1)}\phi_2+\Uarm^{(2,2)}\psi_2}.
\end{align*}
Sei $g := \big(\Uarm^{(2,1)}\phi_1+\Uarm^{(2,2)}\psi_1\big)^{-1}\big(\Uarm^{(2,1)}\phi_2+\Uarm^{(2,2)}\psi_2\big)$. 
Dann existiert eine Teilmenge $\D$ von $\C\setminus[\alpha,\infty)$ (als Vereinigung endlich vieler diskreter Teilmengen von $\C\setminus[\alpha,\infty)$), so dass
$\phi_1, \psi_1, \phi_2, \psi_2$ und $g$ in $\C\setminus([\alpha,\infty)\cup\D$ holomorph sind 
und $\det g(z) \neq 0$ sowie 
\begin{align*}
	\binom{\phi_2(z)}{\psi_2(z)} = \binom{\phi_1(z)}{\psi_1(z)}g(z)
\end{align*}
für alle $z \in \C\setminus([\alpha,\infty)\cup\D$ erfüllt sind. 
Wegen \thref{spdef4b} gilt dann $\sklam{\binom{\phi_1}{\psi_1}} = \sklam{\binom{\phi_2}{\psi_2}}$.

Sei nun (ii) erfüllt, das heißt wegen \thref{spdef4b} existieren eine in $\C\setminus[\alpha,\infty)$ meromorphe \textit{q}$\times$\textit{q}"=Matrixfunktion $g$ und diskrete Teilmenge $\D$ von $\C\setminus[\alpha,\infty)$,
so dass $\phi_1, \psi_1, \phi_2, \psi_2$ und $g$ in $\C\setminus([\alpha,\infty)\cup\D)$ holomorph sind 
sowie $\det g(z) \neq 0$ und
\begin{align*}
  \binom{\phi_2(z)}{\psi_2(z)} = \binom{\phi_1(z)}{\psi_2(z)}g(z)
\end{align*}
für alle $z \in \C\setminus([\alpha,\infty)\cup\D)$ erfüllt sind. Hieraus folgt dann
\begin{align*}
  &\ \eklam{\Uarm^{(1,1)}(z)\phi_2(z)+\Uarm^{(1,2)}(z)\psi_2(z)}\eklam{\Uarm^{(2,1)}(z)\phi_2(z)+\Uarm^{(2,2)}(z)\psi_2(z)}^{-1} \\
  & = \eklam{\Uarm^{(1,1)}(z)\phi_1(z)g(z)+\Uarm^{(1,2)}(z)\psi_1(z)g(z)}\eklam{\Uarm^{(2,1)}\phi_1(z)g(z)+\Uarm^{(2,2)}(z)\psi_1(z)g(z)}^{-1} \\
  & = \eklam{\Uarm^{(1,1)}(z)\phi_1(z)+\Uarm^{(1,2)}(z)\psi_1(z)}g(z)g^{-1}(z)\eklam{\Uarm^{(2,1)}(z)\phi_1(z)+\Uarm^{(2,2)}(z)\psi_1(z)}^{-1} \\
  & = \eklam{\Uarm^{(1,1)}(z)\phi_1(z)+\Uarm^{(1,2)}(z)\psi_1(z)}\eklam{\Uarm^{(2,1)}(z)\phi_1(z)+\Uarm^{(2,2)}(z)\psi_1(z)}^{-1}
\end{align*}
für alle $z \in \C\setminus([\alpha,\infty)\cup\D)$. Hieraus folgt wegen des Identitätssatzes für meromorphe Funktionen (vergleiche z.\,B. im skalaren Fall \cite[Satz 10.3.2]{Funk}; im matriziellen Fall betrachtet man die einzelnen Einträge der Matrixfunktion) dann (i). \bwend

Es sei bemerkt, dass \thref{drth1} für ungerade $m$ im allgemeineren Fall einer gegebenen Folge aus $\Keqma$ in \cite[Kapitel 11]{Maka} behandelt wurde. Insbesondere ist dort auch der nichtdegnerierte Fall in \cite[Theorem 11.2]{Maka} aufgeführt.

\thref{drth1} zeigt nun, dass über die in \thref{drdef2} eingeführte Begriffsbildung tatsächlich eine Resolventenmatrix für das rechtsseitige $\alpha$-Stieltjes Momentenproblem gewonnen wird (vergleiche Teil (a) von \thref{drdef0}).

Abschließend können wir folgende Beobachtung für Funktionen aus $\Sqasmu$ für eine Folge $\sjm \in \Kpqma$ mit $m \in \N$ vornehmen.

\begin{satz}	\thlabel{drsa11}
	Seien $m \in \N$, $\alpha \in \R$ und $\sjm \in \Kpqma$. Weiterhin sei \linebreak $S\in\Sqasmu$. Dann gilt $\det S(z) \neq 0$ für alle $z\in\C\setminus[\alpha,\infty)$. Insbesondere ist $S(x)$ für alle $x\in(-\infty,\alpha)$ eine positiv hermitesche Matrix.
\end{satz}

\bwanf Wegen \thref{adpbm1} ist $s_0 = H_0$ eine reguläre Matrix. Hieraus folgt wegen \thref{asmsa3} dann $\det S(z) \neq 0$ für alle $z\in\C\setminus[\alpha,\infty)$. Hieraus folgt wegen $S\in\Soqa$ (vergleiche \thref{asmbz3}) dann $S(x) > \Oq$ für alle $x\in(-\infty,\alpha)$. \bwend

\subsection{Zwei extremale Elemente der Menge \texorpdfstring{$\Sqasmu$}{Sqasmu}} \label{EE}

Seien $m \in \N$, $\alpha \in \R$ und $\sjm \in \Kpqma$ gegeben. Der vorliegende Teilabschnitt ist dann der Diskussion zweier ausgezeichneter Elemente der Menge $\Sqasmu$ gewidmet. Hierbei handelt es sich um rationale \textit{q}$\times$\textit{q}-Matrixfunktionen, die durch spezielle Extremaleigenschaften auf dem Intervall $(-\infty,\alpha)$ bezüglich der Löwner"=Halbordnung für hermitesche \textit{q}$\times$\textit{q}-Matrizen gekennzeichnet werden. Ausgangspunkt unserer Betrachtungen ist die folgende Beobachtung.

\begin{bem}	\thlabel{drbsp1}
	Seien $\kappa \in \Na$, $\alpha \in \R$, $\sjk \in \Kpqka$ und $\ABCDark$ das rechtsseitige $\alpha$-Dyukarev-Quadrupel bezüglich $\sjk$. Unter Beachtung von Teil (b) von \thref{drsa12} seien weiterhin $\Sarmmin, \Sarmmax: \C\setminus[\alpha,\infty)\rightarrow\Cqq$ für alle $m\in\Zek$ definiert gemäß
	\begin{align*}
		\Sarmmin(z) := \Barfm(z)\Darfm^{-1}(z) \quad \text{bzw.} \quad \Sarmmax(z) := \Aarfm(z)\Carfm^{-1}(z).
	\end{align*}
	Dann gelten $\Sarmmin, \Sarmmax\in\Sqasmu$ für alle $m\in\Zek$.	
\end{bem}

\bwanf Sei $m\in\Zek$. Wegen Teil (b) von \thref{asmsa1} bzw. Teil (b) von \thref{asmdef2} gilt dann $\sjm \in \Kpqma$. Bezeichne ${\cal I}$ bzw. ${\cal O}$ die in $\C\setminus[\alpha,\infty)$ konstante Matrixfunktion mit dem Wert $\Iq$ bzw. $\Oq$. Wegen \thref{spbsp2} und Teil (a) von \thref{spdef4} gilt dann $\binom{\cal I}{\cal O}, \binom{\cal O}{\cal I} \in \PtJqCa$. Unter Beachtung von \thref{drdef2} folgt dann die Behauptung aus Teil (a) von \thref{drth1}, wobei $\Sarmmin$ bzw. $\Sarmmax$ die zum \textit{q}$\times$\textit{q}-Stieltjes-Paar $\binom{\cal O}{\cal I}$ bzw. $\binom{\cal I}{\cal O}$  in $\C\setminus[\alpha,\infty)$ zugehörige Funktion aus $\Sqasmu$ darstellt. \bwend

\begin{defi}	\thlabel{drdef4}
	Seien $m \in \N$, $\alpha \in \R$, $\sjm \in \Kpqma$ und $\ABCDsarm$ das rechtsseitige $\alpha$-Dyukarev-Quadrupel bezüglich $\sjm$. Weiterhin seien $\Sarmmins, \Sarmmaxs: \C\setminus[\alpha,\infty)\rightarrow\Cqq$ definiert gemäß 
	\begin{align*}
		\Sarmmins(z) := \Bsarfm(z)\beklam{\Dsarfm(z)}^{-1} \text{ bzw. } \Sarmmaxs(z) := \Asarfm(z)\beklam{\Csarfm(z)}^{-1}.
	\end{align*}	
	Dann heißt $\Sarmmins$ bzw. $\Sarmmaxs$ das \textbf{untere} bzw. \textbf{obere Extremalelement} von \linebreak $\Sqasmu$.
	Falls klar ist, von welchem $\sjm$ die Rede ist, lassen wir das \anf{$s$} im oberen Index weg.
\end{defi}

Es sei bemerkt, dass die in \thref{drdef4} eingeführten Größen mit denen in \cite[Chapter 3]{CR1} für den Fall $\alpha=0$ übereinstimmen, wobei dort die zwei Extremalelemente fälschlicherweise vertauscht wurden.

Wir können nun folgende Beobachtung über die zugehörigen Stieltjes-Maße der in \thref{drdef4} eingeführten Funktionen machen.

\begin{bem}	\thlabel{drbm7}
	Seien $m \in \N$, $\alpha \in \R$, $\sjm \in \Kpqma$ und $\Sarmmin$ bzw. $\Sarmmax$ das untere bzw. obere Extremalelement von $\Sqasmu$. Weiterhin sei $\mu^{(\arm)}_{min}$ bzw. $\mu^{(\arm)}_{max}$ das zu $\Sarmmin$ bzw. $\Sarmmax$ gehörige Stieltjes-Maß. \dgfa 
	\begin{itemize} 
		\item [\rm{(a)}] Es sind $\mu^{(\arm)}_{min}$ und $\mu^{(\arm)}_{max}$ molekulare Maße aus $\M^q_{\geq,\infty}([\alpha,\infty))$.
		\item [\rm{(b)}] Sei $\ABCDark$ das rechtsseitige $\alpha$-Dyukarev-Quadrupel bezüglich $\sjk$. Es bezeichne ${\cal N}^{(\arm)}_{min}$ bzw. ${\cal N}^{(\arm)}_{max}$ die Nullstellenmenge von $\det\Darfm$ bzw. $\det\Carfm$. Dann gelten
		\begin{align*}
			\mu^{(min)}_{\arm}\Brklam{[\alpha,\infty)\setminus{\cal N}^{(\arm)}_{min}} = \Oq
			\quad \text{und} \quad
			\mu^{(max)}_{\arm}\Brklam{[\alpha,\infty)\setminus{\cal N}^{(\arm)}_{max}} = \Oq.
		\end{align*}
	\end{itemize}
\end{bem}

\bwanf Dies folgt unter Beachtung von \thref{asmth4} aus \cite[Lemma B.4]{13} in Verbindung mit \thref{drdef4} und Teil (b) von \thref{drsa12}. \bwend

Wir wollen nun erläutern, warum die in \thref{drdef4} definierten Funktionen gerade das untere und obere Extremalelement von $\Sqasmu$ genannt werden. Hierfür verwenden wir eine ähnliche Vorgehensweise wie in \cite[Chapter 3]{dyu}, wo der Fall $\alpha=0$ behandelt wurde. Zunächst benötigen wir aber noch folgendes Hilfsresultat.

\begin{lemma}	\thlabel{drlm5}
	Seien $\kappa \in \Na$, $\alpha \in \R$, $\sjk \in \Kpqka$ und $\ABCDark$ das rechtsseitige $\alpha$-Dyukarev-Quadrupel bezüglich $\sjk$. \dgfa
	\begin{itemize}
		\item [\rm{(a)}] Es gelten
		\begin{align*}
			\Aarfm^{\ast}(x)\Carfm(x)-\Carfm^{\ast}(x)\Aarfm(x) &= \Oq \\
			\Aarfm^{\ast}(x)\Darfm(x)-\Carfm^{\ast}(x)\Barfm(x) &= \Iq \\
			\Darfm^{\ast}(x)\Aarfm(x)-\Barfm^{\ast}(x)\Carfm(x) &= \Iq \\
			\Barfm^{\ast}(x)\Darfm(x)-\Darfm^{\ast}(x)\Barfm(x) &= \Oq
		\end{align*}
		für alle $x\in\R$ und $m\in\Zok$.
		\item [\rm{(b)}] Es gelten
		\begin{align*}
			\Aarfm(x)\Barfm^{\ast}(x)+\Barfm(x)\Aarfm^{\ast}(x) \geq \Oq
		\end{align*}
		und
		\begin{align*}
			\Carfm(x)\Darfm^{\ast}(x)+\Darfm(x)\Carfm^{\ast}(x) \geq \Oq
		\end{align*}
		für alle $x\in(-\infty,\alpha)$ und $m\in\Zok$.
	\end{itemize}
\end{lemma}

\bwanf Zu (a): Sei $\Uarmk$ die Folge von rechtsseitigen $2$\textit{q}$\times2$\textit{q}"=$\alpha$"=Dyukarev"=Matrixpolynomen bezüglich $\sjk$. Wegen \thref{spbsp1}, Teil (c) von \thref{drbm5} und Teil (b) von \thref{spdef3} ist $\Uarm(x)$ für alle $x\in\R$ und $m\in\Zok$ eine $\tJq$-unitäre Matrix. Hieraus folgt wegen Teil (b) von \thref{spdef2} und \thref{drdef2} nun
\begin{align*}
	0_{2q\times2q} &= \tJq - \Uarm^{\ast}(x) \tJq \Uarm(x) \\
	&= \begin{pmatrix} \Oq & -i\Iq \\ i\Iq & \Oq \end{pmatrix} - \begin{pmatrix} \Aarfm(x) & \Barfm(x) \\ \Carfm(x) & \Darfm(x) \end{pmatrix}^{\ast} \begin{pmatrix} \Oq & -i\Iq \\ i\Iq & \Oq \end{pmatrix} \\
	&\quad \cdot\begin{pmatrix} \Aarfm(x) & \Barfm(x) \\ \Carfm(x) & \Darfm(x) \end{pmatrix} \\
	&= \begin{pmatrix} \Oq & -i\Iq \\ i\Iq & \Oq \end{pmatrix} - \begin{pmatrix} \Aarfm^{\ast}(x) & \Carfm^{\ast}(x) \\ \Barfm^{\ast}(x) & \Darfm^{\ast}(x) \end{pmatrix} \\
	&\quad \cdot\begin{pmatrix} -i\Carfm(x) & -i\Darfm(x) \\ i\Aarfm(x) & i\Barfm(x) \end{pmatrix} \\
	&= i \left( \begin{matrix} \Aarfm^{\ast}(x)\Carfm(x)-\Carfm^{\ast}(x)\Aarfm(x) \\ \Iq+\Barfm^{\ast}(x)\Carfm(x)-\Darfm^{\ast}(x)\Aarfm(x) \end{matrix} \right. \\
	&\quad \left. \hspace{3cm} \begin{matrix} -\Iq+\Aarfm^{\ast}(x)\Darfm(x)-\Carfm^{\ast}(x)\Barfm(x) \\ \Barfm^{\ast}(x)\Darfm(x)-\Darfm^{\ast}(x)\Barfm(x) \end{matrix} \right)
\end{align*}
für alle $x\in\R$ und $m\in\Zok$. Hieraus folgen dann alle Behauptungen.

Zu (b): Wegen \thref{drdef1} gelten
\begin{align*}
	\Aaro(x)\Baro^{\ast}(x)+\Baro(x)\Aaro^{\ast}(x) = \Oq
\end{align*}
und
\begin{align*}
	\Caro(x)\Daro^{\ast}(x)+\Daro(x)\Caro^{\ast}(x) 
	&= -(x-\alpha)s^{-1}_0\Iq^{\ast}-\Iq[(x-\alpha)s^{-1}_0]^{\ast} \\
	&= -2(x-\alpha)s^{-1}_0 \geq \Oq
\end{align*}
für alle $x\in(-\infty,\alpha)$. 

Sei nun $m\in\Zek$ ungerade, das heißt es existiert ein $n\in\Zofk$ mit $m=2n+1$. Wegen der Teile (a) und (b) von \thref{drlm1} und der Teile (a) und (e) von \thref{drlm3} gilt
\begin{align}	\label{drlm5bw3}
	\Rn(\alpha)\Tn\Rn(x)\uarn = \Rn(x)\Tn\Rn(\alpha)\uarn = \Rn(x)\Tn\yn = \Rn(x)\un
\end{align}
für alle $x\in\R$. Wegen der Teile (a) und (e) von \thref{drlm3} gilt weiterhin
\begin{align}	\label{drlm5bw4}
	\yna+(x-\alpha)\una\Rna(x) &= \yna\eklam{\Rn^{-\ast}(x)+(x-\alpha)\Tna}\Rna(x) \notag \\
	&= \yna\Rn^{-\ast}(\alpha)\Rna(x) = \uarna\Rna(x)
\end{align}
für alle $x\in\R$. Wegen \thref{drdef1}, Teil (c) von \thref{drbm1}, \fref{drlm5bw4} und \fref{drlm5bw3} gilt nun
\begin{align*}
	&\ \Aarfm(x)\Barfm^{\ast}(x) = \Aarn(x)\Barnp^{\ast}(x) \\
	&= \eklam{\Iq+(x-\alpha)\una\Rna(x)\Hn^{-1}\Rn(\alpha)\vn}\yna\Harn^{-1}\Rn(x)\uarn \\
	&= \yna\Harn^{-1}\Rn(x)\uarn \\
	&\quad + (x-\alpha)\una\Rna(x)\Hn^{-1}\Rn(\alpha)\eklam{\Rn^{-1}(\alpha)\Hn-\Tn\Harn}\Harn^{-1}\Rn(x)\uarn \\
	&= \yna\Harn^{-1}\Rn(x)\uarn + (x-\alpha)\una\Rna(x)\Harn^{-1}\Rn(x)\uarn \\
	&\quad -(x-\alpha)\una\Rna(x)\Hn^{-1}\Rn(\alpha)\Tn\Rn(x)\uarn \\
	&= \uarna\Rna(x)\Harn^{-1}\Rn(x)\uarn - (x-\alpha)\una\Rna(x)\Hn^{-1}\Rn(x)\un
\end{align*}
für alle $x\in\R$. Hieraus folgt dann
\begin{align*}
	&\ \Aarfm(x)\Barfm^{\ast}(x)+\Barfm(x)\Aarfm^{\ast}(x) \\
	&= 2\uarna\Rna(x)\Harn^{-1}\Rn(x)\uarn - 2(x-\alpha)\una\Rna(x)\Hn^{-1}\Rn(x)\un \geq \Oq
\end{align*}
für alle $x\in(-\infty,\alpha)$. 
Weiterhin gilt
\begin{align*}
	\Rn(\alpha)\vn+(x-\alpha)\Rn(\alpha)\Tn\Rn(x)\vn
	&= \Rn(\alpha)\eklam{\Rn^{-1}(x)+(x-\alpha)\Tn}\Rn(x)\vn \\
	&= \Rn(\alpha)\Rn^{-1}(\alpha)\Rn(x)\vn = \Rn(x)\vn
\end{align*}
für alle $x\in\R$. Hieraus folgt wegen \thref{drdef1} und Teil (c) von \thref{drbm1} nun
\begin{align*}
	&\ \Carfm(x)\Darfm^{\ast}(x) = \Carn(x)\Darnp^{\ast}(x) \\
	&= -(x-\alpha)\vna\Rna(x)\Hn^{-1}\Rn(\alpha)\vn\eklam{\Iq-(x-\alpha)\vna\Rna(x)\Harn^{-1}\yn}^{\ast} \\
	&= -(x-\alpha)\vna\Rna(x)\Hn^{-1}\Rn(\alpha)\vn \\
	&\quad +(x-\alpha)^2\vna\Rna(x)\Hn^{-1}\Rn(\alpha)\eklam{\Rn^{-1}(\alpha)\Hn-\Tn\Harn}\Harn^{-1}\Rn(x)\vn \\
	&= -(x-\alpha)\vna\Rna(x)\Hn^{-1}\Rn(\alpha)\vn + (x-\alpha)^2\vna\Rna(x)\Harn^{-1}\Rn(x)\vn \\
	&\quad -(x-\alpha)^2\vna\Rna(x)\Hn^{-1}\Rn(\alpha)\Tn\Rn(x)\vn \\
	&= -(x-\alpha)\vna\Rna(x)\Hn^{-1}\Rn(x)\vn + (x-\alpha)^2\vna\Rna(x)\Harn^{-1}\Rn(x)\vn
\end{align*}
für alle $x\in\R$. Hieraus folgt dann
\begin{align*}
	&\ \Carfm(x)\Darfm^{\ast}(x)+\Carfm(x)\Darfm^{\ast}(x) \\
	&= -2(x-\alpha)\vna\Rna(x)\Hn^{-1}\Rn(x)\vn + 2(x-\alpha)^2\vna\Rna(x)\Harn^{-1}\Rn(x)\vn \geq \Oq
\end{align*}
für alle $x\in(-\infty,\alpha)$.

Seien nun $\kappa\geq2$ und $m\in\Zzk$ gerade, das heißt es existiert ein $n\in\Zefk$ mit $m=2n$. Wegen Teil (d) von \thref{drlm2}, Teil (b) von \thref{drlm1} und der Teile (b) und (e) von \thref{drlm3} gilt
\begin{align}	\label{drlm5bw1}
	\Rn(\alpha)\Ln\Rnm(x)\uarnm = \Rn(x)\Ln\Rnm(\alpha)\uarnm = \Rn(x)\Ln\ynm = \Rn(x)\un
\end{align}
für alle $x\in\R$. Wegen der Teile (a) und (e) von \thref{drlm3} gilt weiterhin
\begin{align}	\label{drlm5bw2}
	\ynma+(x-\alpha)\unma\Rnma(x) &= \ynma\eklam{\Rnm^{-\ast}(x)+(x-\alpha)\Tnma}\Rnma(x) \notag \\
	&= \ynma\Rnm^{-\ast}(\alpha)\Rnma(x) = \uarnma\Rnma(x)
\end{align}
für alle $x\in\R$. Wegen \thref{drdef1}, Teil (d) von \thref{drbm1}, Teil (d) von \thref{drlm2}, Teil (c) von \thref{drlm3}, \fref{drlm5bw2} und \fref{drlm5bw1} gilt nun
\begin{align*}
	&\ \Aarfm(x)\Barfm^{\ast}(x) = \Aarn(x)\Barn^{\ast}(x) \\
	&= \eklam{\Iq+(x-\alpha)\una\Rna(x)\Hn^{-1}\Rn(\alpha)\vn}\ynma\Harnm^{-1}\Rnm(x)\uarnm \\
	&= \ynma\Harnm^{-1}\Rnm(x)\uarnm \\
	&\quad + (x-\alpha)\una\Rna(x)\Hn^{-1}\Rn(\alpha)\beklam{\Rn^{-1}(\alpha)\Hn\dLn-\Ln\Harnm}\Harnm^{-1}\Rnm(x)\uarnm \\
	&= \ynma\Harnm^{-1}\Rnm(x)\uarnm + (x-\alpha)\unma\Rnma(x)\Harnm^{-1}\Rnm(x)\uarnm \\
	&\quad -(x-\alpha)\una\Rna(x)\Hn^{-1}\Rn(\alpha)\Ln\Rnm(x)\uarnm \\
	&= \uarnma\Rnma(x)\Harnm^{-1}\Rnm(x)\uarnm - (x-\alpha)\una\Rna(x)\Hn^{-1}\Rn(x)\un
\end{align*}
für alle $x\in\R$. Hieraus folgt dann
\begin{align*}
	&\ \Aarfm(x)\Barfm^{\ast}(x)+\Barfm(x)\Aarfm^{\ast}(x) \\
	&= 2\uarnma\Rnma(x)\Harnm^{-1}\Rnm(x)\uarnm - 2(x-\alpha)\una\Rna(x)\Hn^{-1}\Rn(x)\un \geq \Oq
\end{align*}
für alle $x\in(-\infty,\alpha)$.
Wegen der Teile (a), (b) und (d) von \thref{drlm2} und Teil (a) von \thref{drlm1} gilt
\begin{align*}
	\Rn(\alpha)\vn+(x-\alpha)\Rn(\alpha)\Ln\Rnm(x)\vnm
	&= \Rn(\alpha)\vn+(x-\alpha)\Rn(\alpha)\Rn(x)\Tn\vn \\
	&= \Rn(\alpha)\eklam{\Rn^{-1}(x)+(x-\alpha)\Tn}\Rn(x)\vn \\
	&= \Rn(\alpha)\Rn^{-1}(\alpha)\Rn(x)\vn = \Rn(x)\vn
\end{align*}
für alle $x\in\R$. Hieraus folgt wegen \thref{drdef1}, Teil (d) von \thref{drbm1} und der Teile (b) und (d) von \thref{drlm2} nun
\begin{align*}
	&\ \Carfm(x)\Darfm^{\ast}(x) = \Carn(x)\Darn^{\ast}(x) \\
	&= -(x-\alpha)\vna\Rna(x)\Hn^{-1}\Rn(\alpha)\vn\eklam{\Iq-(x-\alpha)\vnma\Rnma(x)\Harnm^{-1}\ynm}^{\ast} \\
	&= -(x-\alpha)\vna\Rna(x)\Hn^{-1}\Rn(\alpha)\vn \\
	&\quad +(x-\alpha)^2\vna\Rna(x)\Hn^{-1}\Rn(\alpha)\beklam{\Rn^{-1}(\alpha)\Hn\dLn-\Ln\Harnm}\Harnm^{-1}\Rnm(x)\vnm \\
	&= -(x-\alpha)\vna\Rna(x)\Hn^{-1}\Rn(\alpha)\vn + (x-\alpha)^2\vnma\Rnma(x)\Harnm^{-1}\Rnm(x)\vnm \\
	&\quad -(x-\alpha)^2\vna\Rna(x)\Hn^{-1}\Rn(\alpha)\Ln\Rnm(x)\vnm \\
	&= -(x-\alpha)\vna\Rna(x)\Hn^{-1}\Rn(x)\vn + (x-\alpha)^2\vnma\Rnma(x)\Harnm^{-1}\Rnm(x)\vnm
\end{align*}
für alle $x\in\R$. Hieraus folgt dann
\begin{align*}
	&\ \Carfm(x)\Darfm^{\ast}(x)+\Carfm(x)\Darfm^{\ast}(x) \\
	&= -2(x-\alpha)\vna\Rna(x)\Hn^{-1}\Rn(x)\vn + 2(x-\alpha)^2\vnma\Rnma(x)\Harnm^{-1}\Rnm(x)\vnm \geq \Oq
\end{align*}
für alle $x\in(-\infty,\alpha)$. \bwend

Da wegen \thref{asmbz3} für jedes $S\in\Soqa$ die Matrix $S(x)$ für alle $x\in(-\infty,\alpha)$ (nichtnegativ) hermitesch ist, können wir nun folgende Ungleichungen bezüglich der Löwner Halbordnung auf dem Intervall $(-\infty,\alpha)$ für Funktionen aus $\Sqasmu$ für eine Folge $\sjm \in \Kpqma$ mit $m\in\N$ betrachten (vergleiche \cite[Theorem 4]{dyu} für den Fall $\alpha=0$).

\begin{satz}	\thlabel{drsa8}
	Seien $m \in \N$, $\alpha \in \R$, $\sjm \in \Kpqma$ und $\Sarmmin$ bzw. $\Sarmmax$ das untere bzw. obere Extremalelement von $\Sqasmu$. \dgfa
	\begin{itemize}
		\item [\rm{(a)}] Es gilt
		\begin{align*}
			\Sarmmin(x) < \Sarmmax(x)
		\end{align*}
		für alle $x\in(-\infty,\alpha)$.
		\item [\rm{(b)}] Sei $S\in\Sqasmu$. Dann gilt
		\begin{align*}
			\Sarmmin(x) \leq S(x) \leq \Sarmmax(x)
		\end{align*}
		für alle $x\in(-\infty,\alpha)$.
		\item [\rm{(c)}] Es ist $\beklam{\Sarmmax-\Sarmmin}^{-1}$ holomorph in $\C\setminus[\alpha,\infty)$ und es gilt
		\begin{align*}
			\beklam{\Sarmmax(z)-\Sarmmin(z)}^{-1} &= -(z-\alpha)\vna\Rna(\za)\Hn^{-1}\Rn(z)\vn \\
			&\quad + (z-\alpha)^2\vnma\Rnma(\za)\Harnm^{-1}\Rnm(z)\vnm
		\end{align*}
		für alle $z\in\C\setminus[\alpha,\infty)$.
	\end{itemize}
\end{satz}

\bwanf Sei $\ABCDarm$ das rechtsseitige $\alpha$"=Dyukarev"=Quadrupel bezüglich $\sjm$. 

Zu (a): Wegen $\Sarmmax\in\Sqasmu$ (vergleiche \thref{asmbz3}) ist $\Sarmmax(x)$ für alle $x\in(-\infty,\alpha)$ (nichtnegativ) hermitesch. Hieraus folgt wegen \thref{drdef4} und Teil (a) von \thref{drlm5} dann
\begin{align}	\label{drsa8bw8}
	&\ \Sarmmax(x)-\Sarmmin(x) = \eklam{\Sarmmax(x)}^{\ast}-\Sarmmin(x) \notag \\
	&= \Carfm^{-\ast}(x)\Aarfm^{\ast}(x)-\Barfm(x)\Darfm^{-1}(x) \notag \\
	&= \Carfm^{-\ast}(x)\eklam{\Aarfm^{\ast}(x)\Darfm(x)-\Carfm^{\ast}(x)\Barfm(x)}\Darfm^{-1}(x) \notag \\
	&= \Carfm^{-\ast}(x)\Darfm^{-1}(x)
\end{align}
für alle $x\in(-\infty,\alpha)$. 

Sei nun $m$ ungerade, das heißt es existiert ein $n\in\No$ mit $m=2n+1$. Dann gilt
\begin{align*}
	\vna\Rna(\alpha)+(x-\alpha)\vna\Rna(x)\Tna\Rna(\alpha) 
	&= \vna\Rna(x)\eklam{\Rn^{-\ast}(x)+(x-\alpha)\Tna}\Rna(\alpha) \\
	&= \vna\Rna(x)\eklam{I_{(n+1)q}-x\Tn+(x-\alpha)\Tn}^{\ast}\Rna(\alpha) \\
	&= \vna\Rna(x)\Rn^{-\ast}(\alpha)\Rna(\alpha) = \vna\Rna(x)
\end{align*}
für alle $x\in\R$. Hieraus folgt wegen \fref{drsa8bw8}, \thref{drdef1} und Teil (c) von \thref{drbm1} dann
\begin{align*}
	&\ \beklam{\Sarmmax(x)-\Sarmmin(x)}^{-1} = \Darfm(x)\Carfm^{\ast}(x) = \Darnp(x)\Carn(x) \\
	&= \eklam{\Iq-(x-\alpha)\vna\Rna(x)\Harn^{-1}\yn}\eklam{-(x-\alpha)\vna\Rna(x)\Hn^{-1}\Rn(\alpha)\vn}^{\ast} \\
	&= -(x-\alpha)\vna\Rna(\alpha)\Hn^{-1}\Rn(x)\vn+(x-\alpha)^2\vna\Rna(x)\Harn^{-1}\yn\vna\Rna(\alpha)\Hn^{-1}\Rn(x)\vn \\
	&= -(x-\alpha)\vna\Rna(\alpha)\Hn^{-1}\Rn(x)\vn \\
	&\quad +(x-\alpha)^2\vna\Rna(x)\Harn^{-1}\eklam{\Hn\Rn^{-\ast}(\alpha)-\Harn\Tna}\Rna(\alpha)\Hn^{-1}\Rn(x)\vn \\
	&= -(x-\alpha)\vna\Rna(\alpha)\Hn^{-1}\Rn(x)\vn + (x-\alpha)^2\vna\Rna(x)\Harn^{-1}\Rn(x)\vn \\
	&\quad -(x-\alpha)^2\vna\Rna(x)\Tna\Rna(\alpha)\Hn^{-1}\Rn(x)\vn  \\
	&= -(x-\alpha)\vna\Rna(x)\Hn^{-1}\Rn(x)\vn + (x-\alpha)^2\vna\Rna(x)\Harn^{-1}\Rn(x)\vn	
	> \Oq
\end{align*}
für alle $x\in(-\infty,\alpha)$. Hieraus folgt nun $\Sarmmax(x)-\Sarmmin(x) > \Oq$, also gilt \linebreak $\Sarmmin(x) < \Sarmmax(x)$ für alle $x\in(-\infty,\alpha)$. 

Sei nun $m$ gerade, das heißt es existiert ein $n\in\N$ mit $m=2n$. Wegen der Teile (a), (b) und (d) von \thref{drlm2} und Teil (a) von \thref{drlm1} gilt
\begin{align*}
	\vna\Rna(\alpha)+(x-\alpha)\vnma\Rnma(x)\Lna\Rna(\alpha)
	&= \vna\Rna(\alpha)+(x-\alpha)\vna\Tna\Rna(x)\Rna(\alpha) \\
	&= \vna\Rna(x)\eklam{\Rn^{-\ast}(x)+(x-\alpha)\Tna}\Rna(\alpha) \\
	&= \vna\Rna(x)\eklam{I_{(n+1)q}-x\Tn+(x-\alpha)\Tn}^{\ast}\Rna(\alpha) \\
	&= \vna\Rna(x)\Rn^{-\ast}(\alpha)\Rna(\alpha) = \vna\Rna(x)
\end{align*}
für alle $x\in\R$. Hieraus folgt wegen \fref{drsa8bw8}, \thref{drdef1}, Teil (d) von \thref{drbm1} und der Teile (b) und (d) von \thref{drlm2} nun
\begin{align}	\label{drsa8bw9}
	&\ \beklam{\Sarmmax(x)-\Sarmmin(x)}^{-1} = \Darfm(x)\Carfm^{\ast}(x) = \Darn(x)\Carn(x) \notag \\
	&= \eklam{\Iq-(x-\alpha)\vnma\Rnma(x)\Harnm^{-1}\ynm}\eklam{-(x-\alpha)\vna\Rna(x)\Hn^{-1}\Rn(\alpha)\vn}^{\ast} \notag \\
	&= -(x-\alpha)\vna\Rna(\alpha)\Hn^{-1}\Rn(x)\vn \notag \\
	&\quad +(x-\alpha)^2\vnma\Rnma(x)\Harnm^{-1}\ynm\vna\Rna(\alpha)\Hn^{-1}\Rn(x)\vn \notag \\
	&= -(x-\alpha)\vna\Rna(\alpha)\Hn^{-1}\Rn(x)\vn \notag \\
	&\quad +(x-\alpha)^2\vnma\Rnma(x)\Harnm^{-1}\beklam{\dLna\Hn\Rn^{-\ast}(\alpha)-\Harnm\Lna}\Rna(\alpha)\Hn^{-1}\Rn(x)\vn \notag \\
	&= -(x-\alpha)\vna\Rna(\alpha)\Hn^{-1}\Rn(x)\vn + (x-\alpha)^2\vnma\Rnma(x)\Harnm^{-1}\Rnm(x)\vnm \notag \\
	&\quad +(x-\alpha)^2\vnma\Rnma(x)\Lna\Rna(\alpha)\Hn^{-1}\Rn(x)\vn \notag \\
	&= -(x-\alpha)\vna\Rna(x)\Hn^{-1}\Rn(x)\vn + (x-\alpha)^2\vnma\Rnma(x)\Harnm^{-1}\Rnm(x)\vnm
	> \Oq
\end{align}
für alle $x\in(-\infty,\alpha)$. Hieraus folgt nun $\Sarmmax(x)-\Sarmmin(x) > \Oq$ und somit $\Sarmmin(x) < \Sarmmax(x)$ für alle $x\in(-\infty,\alpha)$.

Zu (b): Wegen Teil (b) von \thref{drth1} und \thref{drdef2} existiert ein \linebreak $\phipsi \in \dPtJqCa$ derart, dass
\begin{align*}
	\det\eklam{\Carfm(z)\phi(z)+\Darfm(z)\psi(z)} \neq 0
\end{align*}
und
\begin{align}	\label{drsa8bw1}
	S(z) = \eklam{\Aarfm(z)\phi(z)+\Barfm(z)\psi(z)}\eklam{\Carfm(z)\phi(z)+\Darfm(z)\psi(z)}^{-1}
\end{align}
für alle $z \in \C\setminus[\alpha,\infty)$ erfüllt sind. Wegen $\Sarmmax\in\Sqasmu$ (vergleiche \thref{asmbz3}) ist $\Sarmmax(x)$ für alle $x\in(-\infty,\alpha)$ (nichtnegativ) hermitesch. Hieraus folgt wegen \thref{drdef4}, \fref{drsa8bw1} und Teil (a) von \thref{drlm5} dann
\begin{align}	\label{drsa8bw2}
	\Sarmmax(x)-S(x) &= \eklam{\Sarmmax(x)}^{\ast}-S(x) \notag \\
	&= \Carfm^{-\ast}(x)\Aarfm^{\ast}(x)-\eklam{\Aarfm(x)\phi(x)+\Barfm(x)\psi(x)} \notag \\ 
	&\quad \cdot\eklam{\Carfm(x)\phi(x)+\Darfm(x)\psi(x)}^{-1} \notag \\
	&= \Carfm^{-\ast}(x)\Big(\Aarfm^{\ast}(x)\eklam{\Carfm(x)\phi(x)+\Darfm(x)\psi(x)} \notag \\
	&\quad -\Carfm^{\ast}(x)\eklam{\Aarfm(x)\phi(x)+\Barfm(z)\psi(x)}\Big) \notag \\
	&\quad \cdot\eklam{\Carfm(x)\phi(x)+\Darfm(x)\psi(x)}^{-1} \notag \\
	&= \Carfm^{-\ast}(x)\Big(\eklam{\Aarfm^{\ast}(x)\Carfm(x)-\Carfm^{\ast}(x)\Aarfm(x)}\phi(x) \notag \\
	&\quad +\eklam{\Aarfm^{\ast}(x)\Darfm(x)-\Carfm^{\ast}(x)\Barfm(x)}\psi(x)\Big) \notag \\
	&\quad \cdot\eklam{\Carfm(x)\phi(x)+\Darfm(x)\psi(x)}^{-1} \notag \\
	&= \Carfm^{-\ast}(x)\psi(x)\eklam{\Carfm(x)\phi(x)+\Darfm(x)\psi(x)}^{-1}
\end{align}
für alle $x\in(-\infty,\alpha)$. Sei
\begin{align*}
	W^{(1)}_{\arm} := \begin{pmatrix} \Oq & \Carfm^{-\ast} \\ \Carfm & \Darfm \end{pmatrix}.
\end{align*}
Wegen \fref{drsa8bw2} gilt dann
\begin{align}	\label{drsa8bw3}
	&\ W^{(1)}_{\arm}(x)\binom{\phi(x)}{\psi(x)}\eklam{\Carfm(x)\phi(x)+\Darfm(x)\psi(x)}^{-1} \notag \\
	&= \binom{\Carfm^{-\ast}(x)\psi(x)}{\Carfm(x)\phi(x)+\Darfm(x)\psi(x)}\eklam{\Carfm(x)\phi(x)+\Darfm(x)\psi(x)}^{-1} \notag \\
	&= \binom{\Sarmmax(x)-S(x)}{\Iq}
\end{align}
für alle $x\in(-\infty,\alpha)$. Sei $\Jq$ definiert wie in \thref{spbsp4}. Wegen Teil (b) von \thref{drlm5} gilt weiterhin
\begin{align*}
	&\ \eklam{W^{(1)}_{\arm}(x)}^{\ast}(-\Jq)W^{(1)}_{\arm}(x)+\Jq  \\
	&= \begin{pmatrix} \Oq & \Carfm^{\ast}(x) \\ \Carfm^{-1}(x) & \Darfm^{\ast}(x) \end{pmatrix} \begin{pmatrix} \Carfm(x) & \Darfm(x) \\ \Oq & \Carfm^{-\ast}(x) \end{pmatrix} + \begin{pmatrix} \Oq & -\Iq \\ -\Iq & \Oq \end{pmatrix} \\
	&= \begin{pmatrix} \Oq & \Oq \\ \Oq & \Carfm^{-1}(x)\Darfm(x)+\Darfm^{\ast}(x)\Carfm^{-\ast}(x) \end{pmatrix} \\
	&= \begin{pmatrix} \Oq & \Oq \\ \Oq & \Carfm^{-1}(x)\beklam{\Darfm(x)\Carfm^{\ast}(x)+\Carfm(x)\Darfm^{\ast}(x)}\Carfm^{-\ast}(x) \end{pmatrix} \\
	&\geq \Oq
\end{align*}
für alle $x\in(-\infty,\alpha)$. Hieraus folgt dann 
\begin{align*}
	\eklam{W^{(1)}_{\arm}(x)}^{\ast}(-\Jq)W^{(1)}_{\arm}(x)\geq -\Jq
\end{align*}
für alle $x\in(-\infty,\alpha)$. Hieraus folgt wegen \thref{splm5} nun
\begin{align*}
	\binom{\phi(x)}{\psi(x)}^{\ast}[W^{(1)}_{\arm}(x)]^{\ast}(-\Jq)W^{(1)}_{\arm}(x)\binom{\phi(x)}{\psi(x)}
	\geq \binom{\phi(x)}{\psi(x)}^{\ast}(-\Jq)\binom{\phi(x)}{\psi(x)} \geq \Oq
\end{align*}
für alle $x\in(-\infty,\alpha)$. Hieraus folgt wegen \fref{drsa8bw3} dann
\begin{align}	\label{drsa8bw4}
	& \eklam{\Sarmmax(x)}^{\ast}-S^{\ast}(x)+\Sarmmax(x)-S(x) \notag \\
	&= \binom{\Sarmmax(x)-S(x)}{\Iq}^{\ast}(-\Jq)\binom{\Sarmmax(x)-S(x)}{\Iq} \notag \\
	&= \eklam{\Carfm(x)\phi(x)+\Darfm(x)\psi(x)}^{-\ast}\binom{\phi(x)}{\psi(x)}^{\ast}[W^{(1)}_{\arm}(x)]^{\ast}(-\Jq) \notag \\
	&\quad \cdot W^{(1)}_{\arm}(x)\binom{\phi(x)}{\psi(x)}\eklam{\Carfm(x)\phi(x)+\Darfm(x)\psi(x)}^{-1} \geq \Oq
\end{align}
für alle $x\in(-\infty,\alpha)$. Wegen $\Sarmmax,S\in\Sqasmu$ (vergleiche \thref{asmbz3}) sind $\Sarmmax(x)$ und $S(x)$ für alle $x\in(-\infty,\alpha)$ jeweils (nichtnegativ) hermitesch. Hieraus folgt wegen \fref{drsa8bw4} nun $\Sarmmax(x)\geq S(x)$ für alle $x\in(-\infty,\alpha)$.

Wegen $\Sarmmin\in\Sqasmu$ (vergleiche \thref{asmbz3}) ist $\Sarmmin(x)$ für alle $x\in(-\infty,\alpha)$ (nichtnegativ) hermitesch. Hieraus folgt wegen \thref{drdef4}, \fref{drsa8bw1} und Teil (a) von \thref{drlm5} dann
\begin{align}	\label{drsa8bw5}
	&\ S(x)-\Sarmmin(x) = S(x)-\beklam{\Sarmmin(x)}^{\ast} \notag \\
	&= \eklam{\Aarfm\phi(x)+\Barfm(x)\psi(x)}\eklam{\Carfm(x)\phi(x)+\Darfm(x)\psi(x)}^{-1} \notag \\ 
	&\quad -\Darfm^{-\ast}(x)\Barfm^{\ast}(x) \notag \\
	&= \Darfm^{-\ast}(x)\Big(\Darfm^{\ast}(x)\eklam{\Aarfm(x)\phi(x)+\Barfm(x)\psi(x)} \notag \\
	&\quad -\Barfm^{\ast}(x)\eklam{\Carfm(x)\phi(x)+\Darfm(z)\psi(x)}\Big) \notag \\
	&\quad \cdot\eklam{\Carfm(x)\phi(x)+\Darfm(x)\psi(x)}^{-1} \notag \\
	&= \Darfm^{-\ast}(x)\Big(\eklam{\Darfm^{\ast}(x)\Aarfm(x)-\Barfm^{\ast}(x)\Carfm(x)}\phi(x) \notag \\
	&\quad +\eklam{\Darfm^{\ast}(x)\Barfm(x)-\Barfm^{\ast}(x)\Darfm(x)}\psi(x)\Big) \notag \\
	&\quad \cdot\eklam{\Carfm(x)\phi(x)+\Darfm(x)\psi(x)}^{-1} \notag \\
	&= \Darfm^{-\ast}(x)\phi(x)\eklam{\Carfm(x)\phi(x)+\Darfm(x)\psi(x)}^{-1}
\end{align}
für alle $x\in(-\infty,\alpha)$. Sei
\begin{align*}
	W^{(2)}_{\arm} := \begin{pmatrix} \Darfm^{-\ast} & \Oq \\ \Carfm & \Darfm \end{pmatrix}.
\end{align*}
Wegen \fref{drsa8bw5} gilt dann
\begin{align}	\label{drsa8bw6}
	&\ W^{(2)}_{\arm}(x)\binom{\phi(x)}{\psi(x)}\eklam{\Carfm(x)\phi(x)+\Darfm(x)\psi(x)}^{-1} \notag \\
	&= \binom{\Darfm^{-\ast}(x)\phi(x)}{\Carfm(x)\phi(x)+\Darfm(x)\psi(x)}\eklam{\Carfm(x)\phi(x)+\Darfm(x)\psi(x)}^{-1} \notag \\
	&= \binom{S(x)-\Sarmmin(x)}{\Iq}
\end{align}
für alle $x\in(-\infty,\alpha)$. Wegen Teil (b) von \thref{drlm5} gilt weiterhin
\begin{align*}
	&\ \eklam{W^{(2)}_{\arm}(x)}^{\ast}(-\Jq)W^{(2)}_{\arm}(x)+\Jq  \\
	&= \begin{pmatrix} \Darfm^{-1}(x) & \Carfm^{\ast}(x) \\ \Oq & \Darfm^{\ast}(x) \end{pmatrix} \begin{pmatrix} \Carfm(x) & \Darfm(x) \\ \Darfm^{-\ast}(x) & \Oq \end{pmatrix} + \begin{pmatrix} \Oq & -\Iq \\ -\Iq & \Oq \end{pmatrix} \\
	&= \begin{pmatrix} \Darfm^{-1}(x)\Carfm(x)+\Carfm^{\ast}(x)\Darfm^{-\ast}(x) & \Oq \\ \Oq & \Oq \end{pmatrix} \\
	&= \begin{pmatrix} \Darfm^{-1}(x)\beklam{\Carfm(x)\Darfm^{\ast}(x)+\Darfm(x)\Carfm^{\ast}(x)}\Darfm^{-\ast}(x) & \Oq \\ \Oq & \Oq \end{pmatrix} \\
	&\geq \Oq
\end{align*}
für alle $x\in(-\infty,\alpha)$. Hieraus folgt dann 
\begin{align*}
	\eklam{W^{(2)}_{\arm}(x)}^{\ast}(-\Jq)W^{(2)}_{\arm}(x)\geq -\Jq
\end{align*}
für alle $x\in(-\infty,\alpha)$. Hieraus folgt wegen \thref{splm5} nun
\begin{align*}
	\binom{\phi(x)}{\psi(x)}^{\ast}[W^{(2)}_{\arm}(x)]^{\ast}(-\Jq)W^{(2)}_{\arm}(x)\binom{\phi(x)}{\psi(x)}
	\geq \binom{\phi(x)}{\psi(x)}^{\ast}(-\Jq)\binom{\phi(x)}{\psi(x)} \geq \Oq
\end{align*}
für alle $x\in(-\infty,\alpha)$. Hieraus folgt wegen \fref{drsa8bw6} dann
\begin{align}	\label{drsa8bw7}
	&\ S^{\ast}(x)-\beklam{\Sarmmin(x)}^{\ast}+S(x)-\Sarmmin(x) \notag \\
	&= \binom{S(x)-\Sarmmin(x))}{\Iq}^{\ast}(-\Jq)\binom{S(x)-\Sarmmin(x))}{\Iq} \notag \\
	&= \eklam{\Carfm(x)\phi(x)+\Darfm(x)\psi(x)}^{-\ast}\binom{\phi(x)}{\psi(x)}^{\ast}[W^{(2)}_{\arm}(x)]^{\ast}(-\Jq) \notag \\
	&\quad \cdot W^{(2)}_{\arm}(x)\binom{\phi(x)}{\psi(x)}\eklam{\Carfm(x)\phi(x)+\Darfm(x)\psi(x)}^{-1} \geq \Oq
\end{align}
für alle $x\in(-\infty,\alpha)$. Wegen $S,\Sarmmin\in\Sqasmu$ (vergleiche \thref{asmbz3}) sind $S(x)$ und $\Sarmmin(x)$ für alle $x\in(-\infty,\alpha)$ jeweils (nichtnegativ) hermitesch. Hieraus folgt wegen \fref{drsa8bw7} nun $S(x)\geq \Sarmmin(x)$ für alle $x\in(-\infty,\alpha)$. 

Zu (c): Unter Beachtung von Teil (b) von \thref{drsa12} sei $F:\C\setminus[\alpha,\infty)\rightarrow\Cqq$ definiert gemäß $F(z) = \Carfm^{-\ast}(\za)\Darfm^{-1}(z)$. Da $\Carfm$ und $\Darfm$ als Matrixpolynome holomorph in $\C$ sind, folgt wegen Teil (b) von \thref{drsa12} dann, dass $F$ holomorph in $\C\setminus[\alpha,\infty)$ ist und $\det F(z) \neq 0$ für alle $z\in\C\setminus[\alpha,\infty)$ erfüllt ist. Wegen \fref{drsa8bw8} gilt
\begin{align*}
	\Sarmmax(x)-\Sarmmin(x) = F(x)
\end{align*}
für alle $x\in(-\infty,\alpha)$. Hieraus folgt wegen des Identitätssatzes für meromorphe Funktionen (vergleiche z.\,B. im skalaren Fall \cite[Satz 10.3.2]{Funk}; im matriziellen Fall betrachtet man die einzelnen Einträge der Matrixfunktion) dann
\begin{align*}
	\Sarmmax(z)-\Sarmmin(z) = F(z)
\end{align*}
für alle $z\in\C\setminus[\alpha,\infty)$. Somit ist $\Sarmmax-\Sarmmin$ holomorph in $\C\setminus[\alpha,\infty)$ und $\det\beklam{\Sarmmax(z)-\Sarmmin(z)} \neq 0$ für alle $z\in\C\setminus[\alpha,\infty)$ erfüllt. Hieraus folgt dann, dass $\beklam{\Sarmmax-\Sarmmin}^{-1}$ holomorph in $\C\setminus[\alpha,\infty)$ ist. Sei $G:\C\rightarrow\Cqq$ definiert gemäß
\begin{align*}
	G(z) = -(z-\alpha)\vna\Rna(\za)\Hn^{-1}\Rn(z)\vn + (z-\alpha)^2\vnma\Rnma(\za)\Harnm^{-1}\Rnm(z)\vnm.
\end{align*}
Offensichtlich ist $G$ ein \textit{q}$\times$\textit{q}-Matrixpolynom und somit holomorph in $\C$. Wegen \fref{drsa8bw9} gilt
\begin{align*}
	\beklam{\Sarmmax(x)-\Sarmmin(x)}^{-1} = G(x)
\end{align*}
für alle $x\in(-\infty,\alpha)$. Hieraus folgt wegen des Identitätssatzes für meromorphe Funktionen (vergleiche z.\,B. im skalaren Fall \cite[Satz 10.3.2]{Funk}; im matriziellen Fall betrachtet man die einzelnen Einträge der Matrixfunktion) dann
\begin{align*}
	\beklam{\Sarmmax(z)-\Sarmmin(z)}^{-1} = G(z)
\end{align*}
für alle $z\in\C\setminus[\alpha,\infty)$. \bwend

Teil (b) von \thref{drsa8} dokumentiert nun, in welcher Weise die beiden in \thref{drdef4} eingeführten rationalen \textit{q}$\times$\textit{q}-Matrixfunktionen eine extremale Position der Menge $\Sqasmu$ für eine Folge $\sjm \in \Kpqma$ mit $m\in\N$ einnehmen.

Der folgende Satz liefert uns eine alternative Darstellung des unteren bzw. oberen Extremalelements von $\Sqasmu$ für eine Folge $\sjm \in \Kpqma$ mit $m\in\N$, welche uns eine bessere Vorstellung von der Struktur dieser \textit{q}$\times$\textit{q}-Matrixfunktionen liefert (vergleiche \cite[Theorem 3]{dyu} für den Fall $\alpha=0$).

\begin{satz}	\thlabel{drsa9}
	Seien $m \in \N$, $\alpha \in \R$ und $\sjm \in \Kpqma$. Weiterhin seien $\Sarmmin$ bzw. $\Sarmmax$ das untere bzw. obere Extremalelement von $\Sqasmu$ und $\tHarn$ für alle $n\in\Zofm$ definiert wie in Teil (e) von \thref{drbm1} (mit $m$ statt $\kappa$). Unter Beachtung von \thref{drlm6} gelten dann
	\begin{align*}
		\Sarmmin(z) = y^{\ast}_{0,\fklam{m-1}}\beklam{H_{\ar\fklam{m-1}}-(z-\alpha)H_{\fklam{m-1}}}^{-1}y_{0,\fklam{m-1}}
	\end{align*}
	für alle $z\in\C\setminus[\alpha,\infty)$ und
	\begin{align*}
		\Sarmmax(z) = \eklam{v^{\ast}_{\fklam{m}}\beklam{T_{\fklam{m}}\widetilde{H}_{\ar\fklam{m}}T^{\ast}_{\fklam{m}}-(z-\alpha)^{-1}R^{-1}_{\fklam{m}}(\alpha)H_{\fklam{m}}R^{-\ast}_{\fklam{m}}(\alpha)}^{-1}v_{\fklam{m}}}^{-1}
	\end{align*}
	für alle $z\in\C\setminus[\alpha,\infty)$.
\end{satz}

\bwanf Zu (a): Seien $\ABCDarm$ das rechtsseitige $\alpha$"=Dyukarev"=Quadrupel bezüglich $\sjm$ und zunächst $n:=\fklamo{m-1}$. Wegen \thref{drdef4} und Teil (a) von \thref{drsa12} gilt dann
\begin{align*}
	\Sarmmin(z) &= \Barfm(z)\Darfm^{-1}(z) = \Barnp(z)\Darnp^{-1}(z) \\
	&= \yna\beklam{\Harn-(z-\alpha)\Hn}^{-1}\yn \\
	&= y^{\ast}_{0,\fklam{m-1}}\beklam{H_{\ar\fklam{m-1}}-(z-\alpha)H_{\fklam{m-1}}}^{-1}y_{0,\fklam{m-1}}
\end{align*}
für alle $z\in\C\setminus[\alpha,\infty)$. 

Sei nun $n:=\fklamo{m}$. Wegen \thref{drdef4} und Teil (a) von \thref{drsa12} gilt dann
\begin{align*}
	\Sarmmax(z) &= \Aarfm(z)\Carfm^{-1}(z) = \Aarn(z)\Carn^{-1}(z) \\
	&= \eklam{\vna\beklam{\Tn\tHarn\Tna-(z-\alpha)^{-1}\Rn^{-1}(\alpha)\Hn\Rn^{-\ast}(\alpha)}^{-1}\vn}^{-1} \\
	&= \eklam{v^{\ast}_{\fklam{m}}\beklam{T_{\fklam{m}}\widetilde{H}_{\ar\fklam{m}}T^{\ast}_{\fklam{m}}-(z-\alpha)R^{-1}_{\fklam{m}}(\alpha)H_{\fklam{m}}R^{-\ast}_{\fklam{m}}(\alpha)}^{-1}v_{\fklam{m}}}^{-1}
\end{align*}
für alle $z\in\C\setminus[\alpha,\infty)$. \bwend

Teil (a) von \thref{drsa8} und \thref{drsa11} führen uns auf folgende Definition.

\begin{defi}	\thlabel{drdef5}
	Seien $m \in \N$, $\alpha \in \R$, $\sjm \in \Kpqma$ und $\Sarmmin$ bzw. $\Sarmmax$ das untere bzw. obere Extremalelement von $\Sqasmu$. Weiterhin sei $x\in(-\infty,\alpha)$. Dann heißt $\beklam{\Sarmmin(x),\Sarmmax(x)}$ das zu $\sjm$ und $x$ zugehörige \textbf{Weylsche Intervall}.
\end{defi}

Das so eben eingeführte zu einer rechtsseitig $\alpha$-Stieltjes-positiv definiten Folge und einem Punkt $x\in(-\infty,\alpha)$ gehörige Weylsche Intervall ist ein nichtdegneriertes Matrixintervall bezüglich der Löwner-Halbordnung. Der folgende Satz zeigt nun, dass dieses Intervall gerade mit der Menge der Funktionswerte $S(x)$ aller Lösungen \linebreak $S\in\Sqasmu$ übereinstimmt (der Fall $\alpha=0$ wurde in \cite{dyu} am Ende von Kapitel 3 ohne Beweis angegeben; einen Beweis, an dem wir uns größtenteils orientieren werden, für diesen Fall formulierte Yu.\,M. Dyukarev freundlicherweise während seines Aufenthaltes in Leipzig Ende Januar 2017).

\begin{satz}	\thlabel{drbm6}
	Seien $m \in \N$, $\alpha \in \R$ und $\sjm \in \Kpqma$. Weiterhin sei $\big[\Sarmmin(x)$, $\Sarmmax(x)\big]$ für alle $x\in(-\infty,\alpha)$ das zu $\sjm$ und $x$ zugehörige Weylsche Intervall. Dann gilt
	\begin{align*}
		 \gklam{S(x) \;\big|\; S\in\Soqa[\sjm,\leq]} = \beklam{\Sarmmin(x),\Sarmmax(x)}
	\end{align*}
	für alle $x\in(-\infty,\alpha)$.
\end{satz}

\bwanf Wegen Teil (b) von \thref{drsa8} gilt
\begin{align*}
	\gklam{S(x) \;\big|\; S\in\Soqa[\sjm,\leq]} \subseteq \beklam{\Sarmmin(x),\Sarmmax(x)}. 
\end{align*}
Sei nun $T\in\beklam{\Sarmmin(x),\Sarmmax(x)}$. Unter Beachtung von Teil (a) von \thref{drsa8} existiert wegen Teil (b) von \thref{ambm1} ein $K\in\Cqq_H$ mit $\Oq\leq K \leq\Iq$ und
\begin{align}	\label{drbm6bw1}
	T = \Sarmmax(x)-\sqrt{\Sarmmax(x)-\Sarmmin(x)}K\sqrt{\Sarmmax(x)-\Sarmmin(x)}.
\end{align}
Weiterhin sei $\ABCDarm$ das rechtsseitige $\alpha$"=Dyukarev"=Quadrupel bezüglich $\sjm$.

Sei zunächst $\Oq< K \leq\Iq$ erfüllt. Dann ist insbesondere $K$ regulär. Unter Beachtung von Teil (a) von \thref{drsa8} und Teil (b) von \thref{drsa12} sei
\begin{align*}
	W &:= \Carfm^{-1}(x)\beklam{\Sarmmax(x)-\Sarmmin(x)}^{-\frac{1}{2}}\rklam{K^{-1}-\Iq} \\
	&\quad\ \cdot\beklam{\Sarmmax(x)-\Sarmmin(x)}^{-\frac{1}{2}}\Carfm^{-\ast}(x).
\end{align*}
Es gilt $K^{-1}\geq\Iq$ und somit ist $W\in\Cqq_{\geq}$. Hieraus folgt $\im W=\Oq$. Bezeichne ${\cal W}$ bzw. ${\cal I}$ die in $\C\setminus[\alpha,\infty)$ konstante Matrixfunktion mit dem Wert $W$ bzw. $\Iq$. Wegen \thref{asmbz3} ist dann ${\cal W}\in\Soqa$ und somit ist wegen Teil (a) von \thref{spbm11} nun $\binom{\cal W}{\cal I}$ ein eigentliches (und konstantes) Paar aus $\dPtJqCa$. Bezeichne $S$ die aus Teil (a) von \thref{drth1} zu $\binom{\cal W}{\cal I}$ gehörige Funktion aus $\Sqasmu$. Wegen Teil (a) von \thref{drth1} und \thref{drdef2} sind dann 
\begin{align*}
	\det\eklam{\Carfm(x)W+\Darfm(x)}\neq0
\end{align*}
und
\begin{align}	\label{drbm6bw2}
	S(x) = \eklam{\Aarfm(x)W+\Barfm(x)}\eklam{\Carfm(x)W+\Darfm(x)}^{-1}
\end{align}
erfüllt. Wegen $\Sarmmax\in\Sqasmu$ (vergleiche \thref{asmbz3}) ist $\Sarmmax(x)$ (nichtnegativ) hermitesch. Hieraus folgen wegen \thref{drdef4} und Teil (a) von \thref{drlm5} dann
\begin{align*}
	&\ \Sarmmax(x)-\Sarmmin(x) = \eklam{\Sarmmax(x)}^{\ast}-\Sarmmin(x) \\
	&= \Carfm^{-\ast}(x)\Aarfm^{\ast}(x)-\Barfm(x)\Darfm^{-1}(x) \\
	&= \Carfm^{-\ast}(x)\eklam{\Aarfm^{\ast}(x)\Darfm(x)-\Carfm^{\ast}(x)\Barfm(x)}\Darfm^{-1}(x) \\
	&= \Carfm^{-\ast}(x)\Darfm^{-1}(x)
\end{align*}
und unter Beachtung von \fref{drbm6bw2} weiterhin
\begin{align*}
	&\ \Sarmmax(x)-S(x) = \eklam{\Sarmmax(x)}^{\ast}-S(x) \\
	&= \Carfm^{-\ast}(x)\Aarfm^{\ast}(x)-\eklam{\Aarfm(x)W+\Barfm(x)}\eklam{\Carfm(x)W+\Darfm(x)}^{-1} \notag \\
	&= \Carfm^{-\ast}(x)\Big(\Aarfm^{\ast}(x)\eklam{\Carfm(x)W+\Darfm(x)} \\
	&\quad -\Carfm^{\ast}(x)\eklam{\Aarfm(x)W+\Barfm(z)}\Big)\eklam{\Carfm(x)W+\Darfm(x)}^{-1} \\
	&= \Carfm^{-\ast}(x)\Big(\eklam{\Aarfm^{\ast}(x)\Carfm(x)-\Carfm^{\ast}(x)\Aarfm(x)}W \\
	&\quad +\eklam{\Aarfm^{\ast}(x)\Darfm(x)-\Carfm^{\ast}(x)\Barfm(x)}\Big)\eklam{\Carfm(x)W+\Darfm(x)}^{-1} \\
	&= \Carfm^{-\ast}(x)\eklam{\Carfm(x)W+\Darfm(x)}^{-1}.
\end{align*}
Hieraus folgt unter Beachtung von Teil (a) von \thref{drsa8} und \fref{drbm6bw1} nun
\begin{align*}
	S(x) &= \Sarmmax(x)-\Carfm^{-\ast}(x)\eklam{\Carfm(x)W+\Darfm(x)}^{-1} \\
	&= \Sarmmax(x)-\eklam{\Carfm(x)W\Carfm^{\ast}(x)+\Darfm(x)\Carfm^{\ast}(x)}^{-1} \\
	&= \Sarmmax(x)-\left(\beklam{\Sarmmax(x)-\Sarmmin(x)}^{-\frac{1}{2}}\rklam{K^{-1}-\Iq}\right. \\
	&\quad \left.\cdot\beklam{\Sarmmax(x)-\Sarmmin(x)}^{-\frac{1}{2}}+\beklam{\Sarmmax(x)-\Sarmmin(x)}^{-1}\right)^{-1} \\
	&= \Sarmmax(x)-\rklam{\beklam{\Sarmmax(x)-\Sarmmin(x)}^{-\frac{1}{2}}K^{-1}\beklam{\Sarmmax(x)-\Sarmmin(x)}^{-\frac{1}{2}}}^{-1} \\
	&= T,
\end{align*}
also gilt $T\in\{S(x) \;|\; S\in\Soqa[\sjm,\leq]\}$.

Seien nun $\Oq\leq K \leq\Iq$ und nicht $\Oq< K$ erfüllt. Im Fall $K=\Oq$ folgt wegen \fref{drbm6bw1} sogleich $T=\Sarmmax\in\Soqa[\sjm,\leq]$. Sei also nun $K\neq\Oq$ erfüllt. Weiterhin seien $r\in\Z{2}{q}$ die Anzahl der verschiedenen Eigenwerte von $K$ und $\lambda_1, \ldots, \lambda_r$ die verschiedenen Eigenwerte in aufsteigender Reihenfolge von $K$. Dann gilt
\begin{align*}
	0=\lambda_1<\lambda_2<\ldots<\lambda_r\leq1.
\end{align*}
Nach dem Spektralsatz ergibt sich
\begin{align*}
	K=0P_1+\lambda_2P_2+\ldots+\lambda_rP_r,
\end{align*}
wobei $P_j$ für jedes $j\in\Z{1}{r}$ die Orthoprojektionsmatrix von $\Cq$ auf den zu $\lambda_j$ gehörigen Eigenraum bezeichnet. Sei $(\eps_k)^{\infty}_{k=0}$ eine gegen Null strebende Folge aus $(0,\lambda_2)$. Weiterhin sei
\begin{align*}
	K_k:=\eps_kP_1+\lambda_2P_2+\ldots+\lambda_rP_r
\end{align*}
für alle $k\in\No$. Für jedes $k\in\No$ gilt dann $\Oq<K_k\leq\Iq$ und, wie wir schon gezeigt haben, gehört
\begin{align*}
	T_k = \Sarmmax(x)-\sqrt{\Sarmmax(x)-\Sarmmin(x)}K_k\sqrt{\Sarmmax(x)-\Sarmmin(x)}
\end{align*}
zu $\{S(x) \;|\; S\in\Soqa[\sjm,\leq]\}$. Hieraus folgt wegen $\lim_{k\rightarrow\infty}K_k=K$ und \fref{drbm6bw1} dann $\lim_{k\rightarrow\infty}T_k=T$. Wie wir schon gezeigt haben, existiert für jedes $k\in\No$ ein konstantes Paar $\binom{\phi_k}{\psi_k}\in\dPtJqCa$, sodass 
\begin{align*}
	\det\eklam{\Carfm(x)\phi_k+\Darfm(x)\psi_k}\neq0
\end{align*}
und
\begin{align*}
	T_k = \eklam{\Aarfm(x)\phi_k+\Barfm(x)\psi_k}\eklam{\Carfm(x)\phi_k+\Darfm(x)\psi_k}^{-1}
\end{align*}
erfüllt sind. Wegen Teil (a) von \thref{spbm17} existiert für jedes $k\in\No$ ein zu $\binom{\phi_k}{\psi_k}$ äquivalentes (und konstantes) Paar $\bbinom{\widetilde{\phi}_k}{\widetilde{\psi}_k}\in\dPtJqCa$ derart, dass $\tpsi_k-i\tphi_k=2{\cal I}$ erfüllt ist. Wegen Teil (c) von \thref{drth1} und \thref{drdef2} sind dann
\begin{align*}
	\det\beklam{\Carfm(x)\widetilde{\phi}_k+\Darfm(x)\widetilde{\psi}_k}\neq0
\end{align*}
und
\begin{align*}
	T_k = \beklam{\Aarfm(x)\widetilde{\phi}_k+\Barfm(x)\widetilde{\psi}_k}\beklam{\Carfm(x)\widetilde{\phi}_k+\Darfm(x)\widetilde{\psi}_k}^{-1}
\end{align*}
für alle $k\in\No$ erfüllt. Wegen Teil (b) von \thref{spbm17} ist die Folge $\Brklam{\bbinom{\widetilde{\phi}_k}{\widetilde{\psi}_k}}^{\infty}_{k=0}$ beschränkt. Wegen des Satzes von Bolzano-Weierstraß existiert dann eine Teilfolge dieser Folge, die gegen ein $\bbinom{\widetilde{\phi}}{\widetilde{\psi}}$ konvergiert. Unter Beachtung von Teil (a) von \thref{spdef4} und Teil (b) von \thref{spbm17} gilt dann $\bbinom{\widetilde{\phi}}{\widetilde{\psi}}\in\dPtJqCa$. Weiterhin gelten
\begin{align*}
	\det\beklam{\Carfm(x)\widetilde{\phi}+\Darfm(x)\widetilde{\psi}}\neq0
\end{align*}
und
\begin{align*}
	T = \beklam{\Aarfm(x)\widetilde{\phi}+\Barfm(x)\widetilde{\psi}}\beklam{\Carfm(x)\widetilde{\phi}+\Darfm(x)\widetilde{\psi}}^{-1}.
\end{align*}
Hieraus folgt wegen Teil (a) von \thref{drth1} nun $T\in\{S(x) \;|\; S\in\Soqa[\sjm,$ $\leq]\}$. \bwend

\subsection{Der linksseitige Fall} \label{chapdrl}

Wir knüpfen nun an die Ausführungen von Abschnitt \ref{chapFML} an. Haben wir dort noch allgemeine linksseitige $\alpha$-Stieltjes Momentenprobleme betrachtet, wollen wir uns nun auf den vollständig nichtdegenerierten Fall beschränken, das heißt unsere gegebene Folge ist linksseitig $\alpha$-Stieltjes-positiv definit anstatt nur -nichtnegativ definit.

\begin{bem}	\thlabel{drlbm1}
  Seien $\kappa \in \Na$, $\alpha \in \R$ und $\sjk \in \Lpqka$. Weiterhin sei $S:$ \linebreak $\C\setminus(-\infty,\alpha] \rightarrow \Cqq$.\dgfa
  \begin{itemize}
    \item [\rm{(a)}] Seien $z \in \C\setminus\R$ und $n \in \Zofk$. \dsfaa
    \begin{itemize}
      \item [\rm{(i)}] Es gilt $\FSn(z) \in \C^{(n+2)q\times(n+2)q}_{\geq}$.
      \item [\rm{(ii)}] Es gilt $\dFSn(z) \in \Cqq_{\geq}$.
    \end{itemize}
    \item [\rm{(b)}] Seien $z \in \C\setminus\R$ und $n \in \Zofkm$. \dsfaa
    \begin{itemize}
      \item [\rm{(iii)}] Es gilt $\FSaln(z) \in \C^{(n+2)q\times(n+2)q}_{\geq}$.
      \item [\rm{(iv)}] Es gilt $\dFSaln(z) \in \Cqq_{\geq}$.
    \end{itemize}
  \end{itemize}
\end{bem}

\bwanf Zu (a): Wegen \thref{adplbm1} gilt $\Hn \in \C^{(n+1)q\times(n+1)q}_{\geq}$ und $\Hn^{+} = \Hn^{-1}$.
Hieraus folgt mithilfe von \thref{drldef1} und \thref{amlm2} dann die Äquivalenz von (i) und (ii).

Zu (b): Wegen \thref{adplbm1} gilt $\Haln \in \C^{(n+1)q\times(n+1)q}_{\geq}$ und $\Haln^{+} = \Haln^{-1}$.
Hieraus folgt mithilfe von \thref{drldef1} und \thref{amlm2} dann die Äquivalenz von (iii) und (iv). \bwend

Unsere nächsten Betrachtungen führen uns vor dem Hintergrund von \thref{drlbm1} von den Potapovschen Fundamentalmatrizen auf eine Resolventenmatrix für das linksseitige $\alpha$-Stieltjes Momentenproblem. Hierfür behandeln wir zunächst einige Matrixfunktionen, die sich als spezielle $\tJq$-Potapov-Funktionen herausstellen werden (vergleiche \thref{spdef3} und \thref{spbsp1}). Sie werden uns helfen die linken Schur-Komplemente der Potapovschen Fundamentalmatrizen anders darzustellen.

\begin{bez}	\thlabel{drlbz2}
	Seien $\kappa \in \Na$, $\alpha \in \R$ und $\sjk \in \Lpqka$. Für $n \in \Zofk$ sei $\Valn: \C \rightarrow \C^{2q\times 2q}$ definiert gemäß
	\begin{align*}
		\Valn^{\sklam{s}}(z) := \Izq - (\alpha-z)\begin{pmatrix} \usna \\ -\vna \end{pmatrix}\Rna(\za)\brklam{\Hsn}^{-1}\Rn(\alpha)\begin{pmatrix} \vn & \usn \end{pmatrix}.
	\end{align*}
	Für $n \in \Zofkm$ sei weiterhin $\Valnp: \C \rightarrow \C^{2q\times 2q}$ definiert gemäß
	\begin{align*}
		\Valnp^{\sklam{s}}(z) := \Izq - (\alpha-z)\begin{pmatrix} \usalna \\ -\vna \end{pmatrix}\Rna(\za)\brklam{\Hsaln}^{-1}\Rn(\alpha)\begin{pmatrix} \vn & \usaln \end{pmatrix}.
	\end{align*}
	Für $n \in \Zofkm$ sei
	\begin{align*}
    	M^{\sklam{s}}_{\aln} := \begin{pmatrix} \Iq & -\ysna\big(\Hsaln\big)^{-1}\ysn \\ \Oq & \Iq \end{pmatrix}.
	\end{align*}
	Für $n \in \Zofk$ sei weiterhin
	\begin{align*}
		\widetilde{M}^{\sklam{s}}_{\aln} := \begin{pmatrix} \Iq & \Oq \\ \vna\Rna(\alpha)\big(\Hsn\big)^{-1}\Rn(\alpha)\vn & \Iq \end{pmatrix}.
	\end{align*}
	Für $z\in\C$ und $n \in \Zofkm$ seien 
	\begin{align*}
    	\Taln^{\sklam{s}}(z) := \Valn^{\sklam{s}}(z) M^{\sklam{s}}_{\aln} \quad \text{und} \quad \tTaln^{\sklam{s}}(z) := \Valnp^{\sklam{s}}(z) \widetilde{M}^{\sklam{s}}_{\aln}.
	\end{align*}
	Im Fall $\kappa\geq2$ seien für $z\in\C$ und $n \in \Zefk$ weiterhin
	\begin{align*}
    	\Yaln^{\sklam{s}}(z) := \Valn^{\sklam{s}}(z) M^{\sklam{s}}_{\aln-1} \quad \text{und} \quad \tYaln^{\sklam{s}}(z) := \Valnm^{\sklam{s}}(z) \widetilde{M}^{\sklam{s}}_{\aln}.
	\end{align*}
	Falls klar ist, von welchem $\sjk$ die Rede ist, lassen wir das \anf{$\sklam{s}$} als oberen Index weg.
\end{bez}

Wir werden nun die in \thref{drlbz2} eingeführten Funktionen und Matrizen für den linksseitigen Fall mithilfe der entsprechenden Funktionen und Matrizen für den rechtsseitigen Fall darstellen. Somit können wir dann die nachfolgenden Resultate mithilfe der entsprechenden Aussagen für den rechtsseitigen Fall beweisen.

\begin{lemma}	\thlabel{drllm3}
	Seien $\kappa \in \Na$, $\alpha \in \R$, $\sjk$ eine Folge aus $\Cqq$ und $t_j:=(-1)^js_j$ für alle $j\in\Zok$. Weiterhin sei $\tjk\in\Lpqkma$ oder $\sjk\in\Kpqka$ erfüllt. Dann gelten $\sjk\in\Kpqka$ bzw. $\tjk\in\Lpqkma$ sowie folgende Aussagen:
	\begin{itemize}		
		\item [\rm{(a)}] Seien $z\in\C$ und $m\in\Zok$. Dann gilt
		\begin{align*}
			\mathbf{V}^{\sklam{t}}_{-\alm}(-z) = \Ve\Varm^{\sklam{s}}(z)\Vea.
		\end{align*}
		\item [\rm{(b)}] Es gelten
		\begin{align*}
			M^{\sklam{t}}_{-\aln} = \Ve M^{\sklam{s}}_{\arn} \Vea
		\end{align*}
		für alle $n \in \Zofkm$ und
		\begin{align*}
			\widetilde{M}^{\sklam{t}}_{-\aln} = \Ve \widetilde{M}^{\sklam{s}}_{\arn} \Vea
		\end{align*}
		für alle $n \in \Zofk$.
		\item [\rm{(c)}] Seien $z\in\C$ und $n \in \Zofkm$. Dann gelten
		\begin{align*}
			\Theta^{\sklam{t}}_{-\aln}(-z) = \Ve\Tarn^{\sklam{s}}(z)\Vea \quad \text{und} \quad
			\widetilde{\Theta}^{\sklam{t}}_{-\aln}(-z) = \Ve\tTarn^{\sklam{s}}(z)\Vea.
		\end{align*}
		\item [\rm{(d)}] Seien $\kappa\geq2$, $z\in\C$ und $n \in \Zofk$. Dann gelten
		\begin{align*}
			\Phi^{\sklam{t}}_{-\aln}(-z) = \Ve\Yarn^{\sklam{s}}(z)\Vea \quad \text{und} \quad
			\widetilde{\Phi}^{\sklam{t}}_{-\aln}(-z) = \Ve\tYarn^{\sklam{s}}(z)\Vea.
		\end{align*}		
	\end{itemize}
\end{lemma}

\bwanf Wegen Teil (a) von \thref{asmbm3} gilt $\sjk\in\Kpqka$ bzw. $\tjk\in\Lpqkma$.

Zu (a): Wegen der Teile (a) und (b) von \thref{asmlm1} und der Teile (a)-(c) von \thref{drllm1} gilt
\begin{align*}
	&\ \mathbf{V}^{\sklam{t}}_{-\al2n}(-z) \\
	&= \Izq - (-\alpha-(-z))\begin{pmatrix} \big(\un^{\sklam{t}}\big)^{\ast} \\ -\vna \end{pmatrix}\Rna(-\za)\rklam{\Hn^{\sklam{t}}}^{-1}\Rn(-\alpha)\begin{pmatrix} \vn & \un^{\sklam{t}}\end{pmatrix} \\
	&= \Izq - (z-\alpha)\begin{pmatrix} -\big(\un^{\sklam{s}}\big)^{\ast}\Vna \\ -\vna\Vna \end{pmatrix}\Vn\Rna(\za)\Vna\Vn\brklam{\Hsn}^{-1}\Vna\Vn\Rn(\alpha)\Vna\begin{pmatrix} \Vn\vn & -\Vn\un^{\sklam{s}}\end{pmatrix} \\
	&= \Izq + (z-\alpha)\begin{pmatrix} \big(\un^{\sklam{s}}\big)^{\ast} \\ \vna \end{pmatrix}\Rna(\za)\brklam{\Hsn}^{-1}\Rn(\alpha)\begin{pmatrix} \vn & -\un^{\sklam{s}}\end{pmatrix} \\
	&= \Ve\Izq\Vea + (z-\alpha)\Ve\begin{pmatrix} \big(\un^{\sklam{s}}\big)^{\ast} \\ -\vna \end{pmatrix}\Rna(\za)\brklam{\Hsn}^{-1}\Rn(\alpha)\begin{pmatrix} \vn & \un^{\sklam{s}}\end{pmatrix}\Vea \\
	&= \Ve\Varn^{\sklam{s}}(z)\Vea
\end{align*}
für alle $n \in \Zofk$. Wegen Teil (a) von \thref{asmlm1}, Teil (b) von \thref{asplm1} und der Teile (a)-(c) von \thref{drllm1} gilt
\begin{align*}
	&\ \mathbf{V}^{\sklam{t}}_{-\al2n+1}(-z) \\
	&= \Izq - (-\alpha-(-z))\begin{pmatrix} \big(\umaln^{\sklam{t}}\big)^{\ast} \\ -\vna \end{pmatrix}\Rna(-\za)\rklam{\Hmaln^{\sklam{t}}}^{-1}\Rn(-\alpha)\begin{pmatrix} \vn & \umaln^{\sklam{t}}\end{pmatrix} \\
	&= \Izq - (z-\alpha)\begin{pmatrix} -\big(\uarn^{\sklam{s}}\big)^{\ast}\Vna \\ -\vna\Vna \end{pmatrix}\Vn\Rna(\za)\Vna\Vn\brklam{\Hsarn}^{-1}\Vna \\
	&\quad \cdot\Vn\Rn(\alpha)\Vna\begin{pmatrix} \Vn\vn & -\Vn\uarn^{\sklam{s}}\end{pmatrix} \\
	&= \Izq + (z-\alpha)\begin{pmatrix} \big(\uarn^{\sklam{s}}\big)^{\ast} \\ \vna \end{pmatrix}\Rna(\za)\brklam{\Hsarn}^{-1}\Rn(\alpha)\begin{pmatrix} \vn & -\uarn^{\sklam{s}}\end{pmatrix} \\
	&= \Ve\Izq\Vea + (z-\alpha)\Ve\begin{pmatrix} \big(\uarn^{\sklam{s}}\big)^{\ast} \\ -\vna \end{pmatrix}\Rna(\za)\brklam{\Hsarn}^{-1}\Rn(\alpha)\begin{pmatrix} \vn & \uarn^{\sklam{s}}\end{pmatrix}\Vea \\
	&= \Ve\Varnp^{\sklam{s}}(z)\Vea
\end{align*}
für alle $n \in \Zofkm$.

Zu (b): Wegen der Teile (a) und (d) von \thref{asmlm1} und Teil (b) von \thref{asplm1} gilt
\begin{align*}
	M^{\sklam{t}}_{-\aln} &= \begin{pmatrix} \Iq & -\big(\yn^{\sklam{t}}\big)^{\ast}\big(\Hmaln^{\sklam{t}}\big)^{-1}\yn^{\sklam{t}} \\ \Oq & \Iq \end{pmatrix}	\\
	&= \begin{pmatrix} \Iq & -\big(\yn^{\sklam{s}}\big)^{\ast}\Vna\Vn\big(\Hsarn\big)^{-1}\Vna\Vn\yn^{\sklam{s}} \\ \Oq & \Iq \end{pmatrix} \\
	&= \Ve\begin{pmatrix} \Iq & \big(\yn^{\sklam{s}}\big)^{\ast}\big(\Hsarn\big)^{-1}\yn^{\sklam{s}} \\ \Oq & \Iq \end{pmatrix}\Vea = \Ve M^{\sklam{s}}_{\arn} \Vea
\end{align*}
für alle $n \in \Zofkm$. Wegen der Teile (a) und (b) von \thref{asmlm1} und der Teile (a) und (b) von \thref{drllm1} gilt
\begin{align*}
	\widetilde{M}^{\sklam{t}}_{-\aln} &= \begin{pmatrix} \Iq & \Oq \\ \vna\Rna(-\alpha)\big(\Hn^{\sklam{t}}\big)^{-1}\Rn(-\alpha)\vn & \Iq \end{pmatrix} \\
	&= \begin{pmatrix} \Iq & \Oq \\ \vna\Vna\Vn\Rna(\alpha)\Vna\Vn\big(\Hn^{\sklam{t}}\big)^{-1}\Vna\Vn\Rn(\alpha)\Vna\Vn\vn & \Iq \end{pmatrix} \\
	&= \Ve\begin{pmatrix} \Iq & \Oq \\ -\vna\Rna(\alpha)\big(\Hn^{\sklam{t}}\big)^{-1}\Rn(\alpha)\vn & \Iq \end{pmatrix}\Vea = \Ve \widetilde{M}^{\sklam{s}}_{\arn} \Vea
\end{align*}
für alle $n \in \Zofk$.

Zu (c): Wegen (a) und (b) sowie Teil (a) von \thref{asmlm1} gelten
\begin{align*}
	\Theta^{\sklam{t}}_{-\aln}(-z) &= \mathbf{V}^{\sklam{t}}_{-\al2n}(-z)M^{\sklam{t}}_{-\aln} \\
	&= \Ve\Varn^{\sklam{s}}(z)\Vea\Ve M^{\sklam{s}}_{\arn} \Vea \\
	&= \Ve\Varn^{\sklam{s}}(z)M^{\sklam{s}}_{\arn} \Vea
	= \Ve\Tarn^{\sklam{s}}(z)\Vea
\end{align*}
und
\begin{align*}
	\widetilde{\Theta}^{\sklam{t}}_{-\aln}(-z) &= \mathbf{V}^{\sklam{t}}_{-\al2n+1}(-z)\widetilde{M}^{\sklam{t}}_{-\aln} \\
	&= \Ve\Varnp^{\sklam{s}}(z)\Vea\Ve \widetilde{M}^{\sklam{s}}_{\arn} \Vea \\
	&= \Ve\Varnp^{\sklam{s}}(z)\widetilde{M}^{\sklam{s}}_{\arn} \Vea
	= \Ve\tTarn^{\sklam{s}}(z)\Vea.
\end{align*}

Zu (d): Wegen (a) und (b) sowie Teil (a) von \thref{asmlm1} gelten
\begin{align*}
	\Phi^{\sklam{t}}_{-\aln}(-z) &= \mathbf{V}^{\sklam{t}}_{-\al2n}(-z)M^{\sklam{t}}_{-\aln-1} \\
	&= \Ve\Varn^{\sklam{s}}(z)\Vea\Ve M^{\sklam{s}}_{\arn-1} \Vea \\
	&= \Ve\Varn^{\sklam{s}}(z)M^{\sklam{s}}_{\arn-1} \Vea
	= \Ve\Yarn^{\sklam{s}}(z)\Vea
\end{align*}
und
\begin{align*}
	\widetilde{\Phi}^{\sklam{t}}_{-\aln}(-z) &= \mathbf{V}^{\sklam{t}}_{-\al2n-1}(-z)\widetilde{M}^{\sklam{t}}_{-\aln} \\
	&= \Ve\Varnm^{\sklam{s}}(z)\Vea\Ve \widetilde{M}^{\sklam{s}}_{\arn} \Vea \\
	&= \Ve\Varnm^{\sklam{s}}(z)\widetilde{M}^{\sklam{s}}_{\arn} \Vea
	= \Ve\tYarn^{\sklam{s}}(z)\Vea. \tag*{$\Box$}
\end{align*}

\begin{satz}	\thlabel{drlsa3}
	Seien $\kappa \in \Na$, $\alpha \in \R$ und $\sjk \in \Lpqka$. \dgfa
	\begin{itemize}
		\item [\rm{(a)}] Es gelten
		\begin{align*}
			\Taln(z) = \begin{pmatrix} \Iq-(\alpha-z)\una\Rna(\za)\Hn^{-1}\Rn(\alpha)\vn &
			\ualna\Rna(\za)\Haln^{-1}\yn \\
			(\alpha-z)\vna\Rna(\za)\Hn^{-1}\Rn(\alpha)\vn &
			\Iq-(\alpha-z)\vna\Rna(\za)\Haln^{-1}\yn \end{pmatrix}
		\end{align*}
		und
		\begin{align*}
			\tTaln(z) = \begin{pmatrix} \Iq-(\alpha-z)\una\Rna(\za)\Hn^{-1}\Rn(\alpha)\vn &
			(\alpha-z)\ualna\Rna(\za)\Haln^{-1}\yn \\
			\vna\Rna(\za)\Hn^{-1}\Rn(\alpha)\vn &
			\Iq-(\alpha-z)\vna\Rna(\za)\Haln^{-1}\yn \end{pmatrix}
		\end{align*}
		für alle $z\in\C$ und $n\in\Zofkm$.
		\item [\rm{(b)}] Es gilt
    	\begin{align*}
      		\tTaln(z) = \begin{pmatrix} (\alpha-z)\Iq & \Oq \\ \Oq & \Iq \end{pmatrix} \Taln(z) \begin{pmatrix} (\alpha-z)^{-1}\Iq & \Oq \\ \Oq & \Iq \end{pmatrix}
    	\end{align*}
    	für alle $z \in \C\setminus\gklam{\alpha}$ und $n \in \Zofkm$.
    	\item [\rm{(c)}] Es gelten $\Taln,\tTaln \in \tbtJqPp$ für alle $n \in \Zofkm$.
    	\item [\rm{(d)}] Für alle $z \in \C$ und $n \in \Zofkm$ sind $\Taln(z)$ und $\tTaln(z)$ regulär.
    	\item [\rm{(e)}] Seien $f: \C\setminus\R \rightarrow \Cqq$ und $f_{\al}(z) := (\alpha-z)f(z)$ für alle $z \in \C\setminus\R$. 
    	Unter Beachtung von (d) gelten dann
    	\begin{align*}
     		\dFfn(z) = \frac{1}{i(z-\za)}\binom{f(z)}{\Iq}^{\ast}\Taln^{-\ast}(z)\tJq\Taln^{-1}(z)\binom{f(z)}{\Iq}
    	\end{align*}
    	und
    	\begin{align*}
      		\dFfaln(z) = \frac{1}{i(z-\za)}\binom{f_{\al}(z)}{\Iq}^{\ast}\tTaln^{-\ast}(z)\tJq\tTaln^{-1}(z)\binom{f_{\al}(z)}{\Iq}
    	\end{align*}
    	für alle $z \in \C\setminus\R$ und $n \in \Zofkm$.	
	\end{itemize}
\end{satz}

\bwanf Seien $t_j=(-1)^js_j$ für alle $j\in\Zok$. Wegen Teil (a) von \thref{asmbm3} gilt dann $\tjk\in\Kpqkma$.

Zu (a): Seien $z \in \C$, $n \in \Zofkm$. Es bezeichne
\begin{align*}
  \Taln^{\sklam{s}} = \begin{pmatrix} \big(\Taln^{\sklam{s}}\big)^{(1,1)} & \big(\Taln^{\sklam{s}}\big)^{(1,2)} \\ \big(\Taln^{\sklam{s}}\big)^{(2,1)} & \big(\Taln^{\sklam{s}}\big)^{(2,2)} \end{pmatrix} \quad \text{bzw.} \quad
  \Theta^{\sklam{t}}_{-\arn} = \begin{pmatrix} \big(\Theta^{\sklam{t}}_{-\arn}\big)^{(1,1)} & \big(\Theta^{\sklam{t}}_{-\arn}\big)^{(1,2)} \\ \big(\Theta^{\sklam{t}}_{-\arn}\big)^{(2,1)} & \big(\Theta^{\sklam{t}}_{-\arn}\big)^{(2,2)} \end{pmatrix}
\end{align*}
die \textit{q}$\times$\textit{q}-Blockzerlegung von $\Taln^{\sklam{s}}$ bzw. $\Theta^{\sklam{t}}_{-\arn}$. Wegen Teil (c) von \thref{drllm3}, Teil (a) von \thref{drsa3}, der Teile (a), (b) und (d) von \thref{asmlm1}, Teil (b) von \thref{asplm1} und der Teile (a)-(c) von \thref{drllm1} gelten
\begin{align*}
	\big(\Taln^{\sklam{s}}\big)^{(1,1)}(z) 
	&= \big(\Theta^{\sklam{t}}_{-\arn}\big)^{(1,1)}(-z) \\
	&= \Iq+(-z-(-\alpha))\big(\un^{\sklam{t}}\big)^{\ast}\Rna(-\za)\big(\Hn^{\sklam{t}}\big)^{-1}\Rn(-\alpha)\vn \\
	&= \Iq-(\alpha-z)\usna\Vna\Vn\Rna(\za)\Vna\Vn\big(\Hsn\big)^{-1}\Vna\Vn\Rn(\alpha)\Vna\Vn\vn \\
	&= \Iq-(\alpha-z)\usna\Rna(\za)\big(\Hsn\big)^{-1}\Rn(\alpha)\vn,
\end{align*}
\begin{align*}
	\big(\Taln^{\sklam{s}}\big)^{(1,2)}(z) 
	&= -\big(\Theta^{\sklam{t}}_{-\arn}\big)^{(1,2)}(-z) \\
	&= -\big(u^{\sklam{t}}_{-\arn}\big)^{\ast}\Rna(-\za)\big(H^{\sklam{t}}_{-\arn}\big)^{-1}\yn^{\sklam{t}} \\
	&= \usalna\Vna\Vn\Rna(\za)\Vna\Vn\big(\Hsn\big)^{-1}\Vna\Vn\ysn \\
	&= \usalna\Rna(\za)\big(\Hsn\big)^{-1}\ysn,
\end{align*}
\begin{align*}
	\big(\Taln^{\sklam{s}}\big)^{(2,1)}(z) 
	&= -\big(\Theta^{\sklam{t}}_{-\arn}\big)^{(2,1)}(-z) \\
	&= (-z-(-\alpha))\vna\Rna(-\za)\big(\Hn^{\sklam{t}}\big)^{-1}\Rn(-\alpha)\vn \\
	&= (\alpha-z)\vna\Vna\Vn\Rna(\za)\Vna\Vn\big(\Hsn\big)^{-1}\Vna\Vn\Rn(\alpha)\Vna\Vn\vn \\
	&= (\alpha-z)\vna\Rna(\za)\big(\Hsn\big)^{-1}\Rn(\alpha)\vn
\end{align*}
und
\begin{align*}
	\big(\Taln^{\sklam{s}}\big)^{(2,2)}(z) 
	&= \big(\Theta^{\sklam{t}}_{-\arn}\big)^{(2,2)}(-z) \\
	&= \Iq-(-z-(-\alpha))\vna\Rna(-\za)\big(H^{\sklam{t}}_{-\arn}\big)^{-1}\yn^{\sklam{t}} \\
	&= \Iq-(\alpha-z)\vna\Vna\Vn\Rna(\za)\Vna\Vn\big(\Hsaln\big)^{-1}\Vna\Vn\ysn \\
	&= \Iq-(\alpha-z)\vna\Rna(\za)\big(\Hsaln\big)^{-1}\ysn.
\end{align*}

Die Aussage für $\tTaln^{\sklam{s}}$ zeigt man analog. 

Zu (b): Dies ist eine direkte Konsequenz aus (a).

Zu (c): Sei $n \in \Zofkm$. Wegen Teil (d) von \thref{drsa3} gilt dann $\Theta^{\sklam{t}}_{-\arn} \in \tbtJqPp$. Wegen Teil (c) von \thref{drllm3}, Teil (d) von \thref{drllm1} und Teil (a) von \thref{asmlm1} gilt
\begin{align}	\label{drlsa3bw1}
	\tJq-\big[\Taln^{\sklam{s}}(z)\big]^{\ast}\tJq\Taln^{\sklam{s}}(z)
	&= \Ve\big[\Theta^{\sklam{t}}_{-\arn}(-z)\big]^{\ast}\Vea\Ve\tJq\Vea\Ve\Theta^{\sklam{t}}_{-\arn}(-z)\Vea-\Ve\tJq\Vea \notag \\
	&= \Ve\rklam{\big[\Theta^{\sklam{t}}_{-\arn}(-z)\big]^{\ast}\tJq\Theta^{\sklam{t}}_{-\arn}(-z)-\tJq}\Vea
\end{align}
für alle $z\in\C$. Wegen Teil (c) von \thref{splm1} ist $\Theta^{\sklam{t}}_{-\arn}(z)$ für alle $z\in\Pm$ eine $\tJq$-expansive Matrix. Hieraus folgt wegen \fref{drlsa3bw1}, Teil (a) von \thref{spdef2} und Teil (a) von \thref{asmlm1} dann, dass $\Taln^{\sklam{s}}(z)$ für alle $z\in\Pp$ eine $\tJq$-kontraktive Matrix ist. Wegen Teil (c) von \thref{drllm3}, Teil (d) von \thref{drllm1}, Teil (a) von \thref{asmlm1}, Teil (b) von \thref{spdef3} und Teil (b) von \thref{spdef2} gilt
\begin{align*}
	\tJq-\big[\Taln^{\sklam{s}}(x)\big]^{\ast}\tJq\Taln^{\sklam{s}}(x)
	&= -\Ve\tJq\Vea+\Ve\big[\Theta^{\sklam{t}}_{-\arn}(-x)\big]^{\ast}\Vea\Ve\tJq\Vea\Ve\Theta^{\sklam{t}}_{-\arn}(-x)\Vea \\
	&= -\Ve\rklam{\tJq-\big[\Theta^{\sklam{t}}_{-\arn}(-x)\big]^{\ast}\tJq\Theta^{\sklam{t}}_{-\arn}(-x)}\Vea = 0_{2q\times2q}
\end{align*}
für alle $x\in\R$, also ist $\Taln^{\sklam{s}}(x)$ für alle $x\in\R$ eine $\tJq$-unitäre Matrix. Somit ist wegen \thref{spdef3} dann $\Taln^{\sklam{s}}\in\tbtJqPp$. 

Die Aussage für $\tTaln^{\sklam{s}}$ zeigt man analog. 

Zu (d): Dies folgt aus (c) und Teil (a) von \thref{splm1}.

Zu (e): Seien $n \in \Zofkm$ und $g:\C\setminus\R\rightarrow\Cqq$ definiert gemäß $g(z) := -f(-z)$. Wegen Teil (b) von \thref{drllm2}, Teil (f) von \thref{drsa3}, Teil (c) von \thref{drllm3}, Teil (d) von \thref{drllm1} und Teil (a) von \thref{asmlm1} gilt dann
\begin{align*}
	\dFfns(z) &= \widehat{\mathbf{F}}^{[g]}_{n,t}(-z) \\
	&= \frac{1}{i(-z-(-\za))}\binom{g(-z)}{\Iq}^{\ast}\big[\Theta^{\sklam{t}}_{-\arn}(-z)\big]^{-\ast}\tJq\big[\Theta^{\sklam{t}}_{-\arn}(-z)\big]^{-1}\binom{g(-z)}{\Iq} \\
	&= \frac{-1}{i(z-\za)}\binom{f(z)}{\Iq}^{\ast}(-\Vea)\Ve\big[\Taln^{\sklam{s}}(z)\big]^{-\ast}\Vea\Ve\rklam{-\tJq}\Vea \\
	&\quad \cdot\Ve\big[\Taln^{\sklam{s}}(z)\big]^{-1}\Vea(-\Ve)\binom{f(z)}{\Iq} \\
	&= \frac{1}{i(z-\za)}\binom{f(z)}{\Iq}^{\ast}\big[\Taln^{\sklam{s}}(z)\big]^{-\ast}\tJq\big[\Taln^{\sklam{s}}(z)\big]^{-1}\binom{f(z)}{\Iq}
\end{align*}
für alle $z\in\C\setminus\R$. Sei $g_{-\ar}(z):=(z-(-\alpha))g(z)$ für alle $z\in\C\setminus\R$. Dann gilt
\begin{align*}
	g_{-\ar}(-z) = (-z-(-\alpha))g(-z) = -(\alpha-z)f(z) = -f_{\al}(z)
\end{align*}
für alle $z\in\C\setminus\R$. Hieraus folgt wegen Teil (b) von \thref{drllm2}, Teil (f) von \thref{drsa3}, Teil (c) von \thref{drllm3}, Teil (d) von \thref{drllm1} und Teil (a) von \thref{asmlm1} dann
\begin{align*}
	&\ \dFfalns(z) = \widehat{\mathbf{F}}^{[g]}_{-\arn,t}(-z) \\
	&= \frac{1}{i(-z-(-\za))}\binom{g_{-\ar}(-z)}{\Iq}^{\ast}\beklam{\widetilde{\Theta}^{\sklam{t}}_{-\arn}(-z)}^{-\ast}\tJq\beklam{\widetilde{\Theta}^{\sklam{t}}_{-\arn}(-z)}^{-1}\binom{g_{-\ar}(-z)}{\Iq} \\
	&= \frac{-1}{i(z-\za)}\binom{f_{\al}(z)}{\Iq}^{\ast}(-\Vea)\Ve\beklam{\tTaln^{\sklam{s}}(z)}^{-\ast}\Vea\Ve\rklam{-\tJq}\Vea \\
	&\quad \cdot\Ve\beklam{\tTaln^{\sklam{s}}(z)}^{-1}\Vea(-\Ve)\binom{f_{\al}(z)}{\Iq} \\
	&= \frac{1}{i(z-\za)}\binom{f_{\al}(z)}{\Iq}^{\ast}\beklam{\tTaln^{\sklam{s}}(z)}^{-\ast}\tJq\beklam{\tTaln^{\sklam{s}}(z)}^{-1}\binom{f_{\al}(z)}{\Iq}
\end{align*}
für alle $z\in\C\setminus\R$. \bwend

\begin{satz}	\thlabel{drlsa4}
	Seien $\kappa \in \Na\setminus\gklam{1}$, $\alpha \in \R$ und $\sjk \in \Lpqka$. \dgfa
	\begin{itemize}
		\item [\rm{(a)}] Es gelten
		\begin{align*}
			\Yaln(z) &= \left( \begin{matrix} \Iq-(\alpha-z)\una\Rna(\za)\Hn^{-1}\Rn(\alpha)\vn \\
			(\alpha-z)\vna\Rna(\za)\Hn^{-1}\Rn(\alpha)\vn \end{matrix} \right. \\
			& \hspace{6cm} \left. \begin{matrix} \ualnma\Rnma(\za)\Halnm^{-1}\ynm \\
			\Iq-(\alpha-z)\vnma\Rnma(\za)\Halnm^{-1}\ynm \end{matrix} \right)
		\end{align*}
		und
		\begin{align*}
			\tYaln(z) &= \left( \begin{matrix} \Iq-(\alpha-z)\una\Rna(\za)\Hn^{-1}\Rn(\alpha)\vn \\			
			\vna\Rna(\za)\Hn^{-1}\Rn(\alpha)\vn \end{matrix} \right. \\
			& \hspace{6cm} \left. \begin{matrix} (\alpha-z)\ualnma\Rnma(\za)\Halnm^{-1}\ynm \\			
			\Iq-(\alpha-z)\vnma\Rnma(\za)\Halnm^{-1}\ynm \end{matrix} \right)
		\end{align*}
		für alle $z\in\C$ und $n\in\Zofk$.
		\item [\rm{(b)}] Es gilt
    	\begin{align*}
      		\tYaln(z) = \begin{pmatrix} (\alpha-z)\Iq & \Oq \\ \Oq & \Iq \end{pmatrix} \Yaln(z) \begin{pmatrix} (\alpha-z)^{-1}\Iq & \Oq \\ \Oq & \Iq \end{pmatrix}
    	\end{align*}
    	für alle $z \in \C\setminus\gklam{\alpha}$ und $n \in \Zofk$.
    	\item [\rm{(c)}] Es gelten $\Yaln,\tYaln \in \tbtJqPp$ für alle $n \in \Zofk$.
    	\item [\rm{(d)}] Für alle $z \in \C$ und $n \in \Zofk$ sind $\Yaln(z)$ und $\tYaln(z)$ regulär.
    	\item [\rm{(e)}] Seien $f: \C\setminus\R \rightarrow \Cqq$ und $f_{\al}(z) := (\alpha-z)f(z)$ für alle $z \in \C\setminus\R$. 
    	Unter Beachtung von (d) gelten dann
    	\begin{align*}
     		\dFfn(z) = \frac{1}{i(z-\za)}\binom{f(z)}{\Iq}^{\ast}\Yaln^{-\ast}(z)\tJq\Yaln^{-1}(z)\binom{f(z)}{\Iq}
    	\end{align*}
    	und
    	\begin{align*}
      		\dFfalnm(z) = \frac{1}{i(z-\za)}\binom{f_{\al}(z)}{\Iq}^{\ast}\tYaln^{-\ast}(z)\tJq\tYaln^{-1}(z)\binom{f_{\al}(z)}{\Iq}
    	\end{align*}
    	für alle $z \in \C\setminus\R$ und $n \in \Zofk$.	
	\end{itemize}
\end{satz}

\bwanf Seien $t_j=(-1)^js_j$ für alle $j\in\Zok$. Wegen Teil (a) von \thref{asmbm3} gilt dann $\tjk\in\Kpqkma$.

Zu (a): Seien $z \in \C$, $n \in \Zofk$. Es bezeichne
\begin{align*}
  \Yaln^{\sklam{s}} = \begin{pmatrix} \big(\Yaln^{\sklam{s}}\big)^{(1,1)} & \big(\Yaln^{\sklam{s}}\big)^{(1,2)} \\ \big(\Yaln^{\sklam{s}}\big)^{(2,1)} & \big(\Yaln^{\sklam{s}}\big)^{(2,2)} \end{pmatrix} \quad \text{bzw.} \quad
  \Phi^{\sklam{t}}_{-\arn} = \begin{pmatrix} \big(\Phi^{\sklam{t}}_{-\arn}\big)^{(1,1)} & \big(\Phi^{\sklam{t}}_{-\arn}\big)^{(1,2)} \\ \big(\Phi^{\sklam{t}}_{-\arn}\big)^{(2,1)} & \big(\Phi^{\sklam{t}}_{-\arn}\big)^{(2,2)} \end{pmatrix}
\end{align*}
die \textit{q}$\times$\textit{q}-Blockzerlegung von $\Yaln^{\sklam{s}}$ bzw. $\Phi^{\sklam{t}}_{-\arn}$. Wegen Teil (d) von \thref{drllm3}, Teil (a) von \thref{drsa4}, der Teile (a), (b) und (d) von \thref{asmlm1}, Teil (b) von \thref{asplm1} und der Teile (a)-(c) von \thref{drllm1} gelten
\begin{align*}
	\big(\Yaln^{\sklam{s}}\big)^{(1,1)}(z) 
	&= \big(\Phi^{\sklam{t}}_{-\arn}\big)^{(1,1)}(-z) \\
	&= \Iq+(-z-(-\alpha))\big(\un^{\sklam{t}}\big)^{\ast}\Rna(-\za)\big(\Hn^{\sklam{t}}\big)^{-1}\Rn(-\alpha)\vn \\
	&= \Iq-(\alpha-z)\usna\Vna\Vn\Rna(\za)\Vna\Vn\big(\Hsn\big)^{-1}\Vna\Vn\Rn(\alpha)\Vna\Vn\vn \\
	&= \Iq-(\alpha-z)\usna\Rna(\za)\big(\Hsn\big)^{-1}\Rn(\alpha)\vn,
\end{align*}
\begin{align*}
	\big(\Yaln^{\sklam{s}}\big)^{(1,2)}(z) 
	&= -\big(\Phi^{\sklam{t}}_{-\arn}\big)^{(1,2)}(-z) \\
	&= -\big(u^{\sklam{t}}_{-\arn-1}\big)^{\ast}\Rnma(-\za)\big(H^{\sklam{t}}_{-\arn-1}\big)^{-1}\ynm^{\sklam{t}} \\
	&= \usalnma\Vnma\Vnm\Rnma(\za)\Vnma\Vnm\big(\Hsnm\big)^{-1}\Vnma\Vnm\ysnm \\
	&= \usalnma\Rnma(\za)\big(\Hsnm\big)^{-1}\ysnm,
\end{align*}
\begin{align*}
	\big(\Yaln^{\sklam{s}}\big)^{(2,1)}(z) 
	&= -\big(\Phi^{\sklam{t}}_{-\arn}\big)^{(2,1)}(-z) \\
	&= (-z-(-\alpha))\vna\Rna(-\za)\big(\Hn^{\sklam{t}}\big)^{-1}\Rn(-\alpha)\vn \\
	&= (\alpha-z)\vna\Vna\Vn\Rna(\za)\Vna\Vn\big(\Hsn\big)^{-1}\Vna\Vn\Rn(\alpha)\Vna\Vn\vn \\
	&= (\alpha-z)\vna\Rna(\za)\big(\Hsn\big)^{-1}\Rn(\alpha)\vn
\end{align*}
und
\begin{align*}
	\big(\Yaln^{\sklam{s}}\big)^{(2,2)}(z) 
	&= \big(\Phi^{\sklam{t}}_{-\arn}\big)^{(2,2)}(-z) \\
	&= \Iq-(-z-(-\alpha))\vnma\Rnma(-\za)\big(H^{\sklam{t}}_{-\arn-1}\big)^{-1}\ynm^{\sklam{t}} \\
	&= \Iq-(\alpha-z)\vnma\Vnma\Vnm\Rnma(\za)\Vnma\Vnm\big(\Hsalnm\big)^{-1}\Vnma\Vnm\ysnm \\
	&= \Iq-(\alpha-z)\vnma\Rnma(\za)\big(\Hsalnm\big)^{-1}\ysnm.
\end{align*}

Die Aussage für $\tYaln^{\sklam{s}}$ zeigt man analog. 

Zu (b): Dies ist eine direkte Konsequenz aus (a).

Zu (c): Sei $n \in \Zofk$. Wegen Teil (d) von \thref{drsa4} gelten dann $\Phi^{\sklam{t}}_{-\arn} \in \tbtJqPp$. Wegen Teil (d) von \thref{drllm3}, Teil (d) von \thref{drllm1} und Teil (a) von \thref{asmlm1} gilt
\begin{align}	\label{drlsa4bw1}
	\tJq-\big[\Yaln^{\sklam{s}}(z)\big]^{\ast}\tJq\Yaln^{\sklam{s}}(z)
	&= \Ve\big[\Phi^{\sklam{t}}_{-\arn}(-z)\big]^{\ast}\Vea\Ve\tJq\Vea\Ve\Phi^{\sklam{t}}_{-\arn}(-z)\Vea-\Ve\tJq\Vea \notag \\
	&= \Ve\rklam{\big[\Phi^{\sklam{t}}_{-\arn}(-z)\big]^{\ast}\tJq\Phi^{\sklam{t}}_{-\arn}(-z)-\tJq}\Vea
\end{align}
für alle $z\in\C$. Wegen Teil (c) von \thref{splm1} ist $\Phi^{\sklam{t}}_{-\arn}(z)$ für alle $z\in\Pm$ eine $\tJq$-expansive Matrix. Hieraus folgt wegen \fref{drlsa4bw1}, Teil (a) von \thref{spdef2} und Teil (a) von \thref{asmlm1} dann, dass $\Yaln^{\sklam{s}}(z)$ für alle $z\in\Pp$ eine $\tJq$-kontraktive Matrix ist. Wegen Teil (d) von \thref{drllm3}, Teil (d) von \thref{drllm1}, Teil (a) von \thref{asmlm1}, Teil (b) von \thref{spdef3} und Teil (b) von \thref{spdef2} gilt
\begin{align*}
	\tJq-\big[\Yaln^{\sklam{s}}(x)\big]^{\ast}\tJq\Yaln^{\sklam{s}}(x)
	&= -\Ve\tJq\Vea+\Ve\big[\Phi^{\sklam{t}}_{-\arn}(-x)\big]^{\ast}\Vea\Ve\tJq\Vea\Ve\Phi^{\sklam{t}}_{-\arn}(-x)\Vea \\
	&= -\Ve\rklam{\tJq-\big[\Phi^{\sklam{t}}_{-\arn}(-x)\big]^{\ast}\tJq\Phi^{\sklam{t}}_{-\arn}(-x)}\Vea = 0_{2q\times2q}
\end{align*}
für alle $x\in\R$, also ist $\Yaln^{\sklam{s}}(x)$ für alle $x\in\R$ eine $\tJq$-unitäre Matrix. Somit ist wegen \thref{spdef3} dann $\Yaln^{\sklam{s}}\in\tbtJqPp$. 

Die Aussage für $\tYaln^{\sklam{s}}$ zeigt man analog. 

Zu (d): Dies folgt aus (c) und Teil (a) von \thref{splm1}.

Zu (e): Seien $n \in \Zofk$ und $g:\C\setminus\R\rightarrow\Cqq$ definiert gemäß $g(z) := -f(-z)$. Wegen Teil (b) von \thref{drllm2}, Teil (f) von \thref{drsa4}, Teil (d) von \thref{drllm3}, Teil (d) von \thref{drllm1} und Teil (a) von \thref{asmlm1} gilt dann
\begin{align*}
	\dFfns(z) &= \widehat{\mathbf{F}}^{[g]}_{n,t}(-z) \\
	&= \frac{1}{i(-z-(-\za))}\binom{g(-z)}{\Iq}^{\ast}\big[\Phi^{\sklam{t}}_{-\arn}(-z)\big]^{-\ast}\tJq\big[\Phi^{\sklam{t}}_{-\arn}(-z)\big]^{-1}\binom{g(-z)}{\Iq} \\
	&= \frac{-1}{i(z-\za)}\binom{f(z)}{\Iq}^{\ast}(-\Vea)\Ve\big[\Yaln^{\sklam{s}}(z)\big]^{-\ast}\Vea\Ve\rklam{-\tJq}\Vea \\
	&\quad \cdot\Ve\big[\Yaln^{\sklam{s}}(z)\big]^{-1}\Vea(-\Ve)\binom{f(z)}{\Iq} \\
	&= \frac{1}{i(z-\za)}\binom{f(z)}{\Iq}^{\ast}\big[\Yaln^{\sklam{s}}(z)\big]^{-\ast}\tJq\big[\Yaln^{\sklam{s}}(z)\big]^{-1}\binom{f(z)}{\Iq}
\end{align*}
für alle $z\in\C\setminus\R$. Sei $g_{-\ar}(z):=(z-(-\alpha))g(z)$ für alle $z\in\C\setminus\R$. Dann gilt
\begin{align*}
	g_{-\ar}(-z) = (-z-(-\alpha))g(-z) = -(\alpha-z)f(z) = -f_{\al}(z)
\end{align*}
für alle $z\in\C\setminus\R$. Hieraus folgt wegen Teil (b) von \thref{drllm2}, Teil (f) von \thref{drsa4}, Teil (d) von \thref{drllm3}, Teil (d) von \thref{drllm1} und Teil (a) von \thref{asmlm1} dann
\begin{align*}
	&\ \dFfalnms(z) = \widehat{\mathbf{F}}^{[g]}_{-\arn-1,t}(-z) \\
	&= \frac{1}{i(-z-(-\za))}\binom{g_{-\ar}(-z)}{\Iq}^{\ast}\beklam{\widetilde{\Phi}^{\sklam{t}}_{-\arn}(-z)}^{-\ast}\tJq\beklam{\widetilde{\Phi}^{\sklam{t}}_{-\arn}(-z)}^{-1}\binom{g_{-\ar}(-z)}{\Iq} \\
	&= \frac{-1}{i(z-\za)}\binom{f_{\al}(z)}{\Iq}^{\ast}(-\Vea)\Ve\beklam{\tYaln^{\sklam{s}}(z)}^{-\ast}\Vea\Ve\rklam{-\tJq}\Vea \\
	&\quad \cdot\Ve\beklam{\tYaln^{\sklam{s}}(z)}^{-1}\Vea(-\Ve)\binom{f_{\al}(z)}{\Iq} \\
	&= \frac{1}{i(z-\za)}\binom{f_{\al}(z)}{\Iq}^{\ast}\beklam{\tYaln^{\sklam{s}}(z)}^{-\ast}\tJq\beklam{\tYaln^{\sklam{s}}(z)}^{-1}\binom{f_{\al}(z)}{\Iq}
\end{align*}
für alle $z\in\C\setminus\R$. \bwend

Die Sätze \ref{drlsa3} und \ref{drlsa4} führen uns auf die Betrachtung spezieller Quadrupel von Folgen von \textit{q}$\times$\textit{q}-Matrixpolynomen, welche in unseren weiteren Betrachtungen eine zentrale Rolle spielen werden.

\begin{defi}	\thlabel{drldef2}
  Seien $\kappa \in \Na$, $\alpha \in \R$ und $\sjk \in \Lpqka$. Weiterhin seien für alle $z \in \C$
  \begin{align*}
    \quad \Bsalo(z) := \Oq, \quad \quad \Dsalo(z) := \Iq
  \end{align*}
  und
  \begin{align*}
    \Asaln(z) & := \Iq-(\alpha-z)\usna\Rna(\za)\big(\Hsn\big)^{-1}\Rn(\alpha)\vn, \\
    \Csaln(z) & := (\alpha-z)\vna\Rna(\za)\big(\Hsn\big)^{-1}\Rn(\alpha)\vn
  \end{align*}
  für alle $n \in \Zofk$ sowie
  \begin{align*}
    \Bsaln(z) & := \rklam{\usalnm}^{\ast}\Rnma(\za)\big(\Hsalnm\big)^{-1}\ysnm, \\
    \Dsaln(z) & := \Iq-(\alpha-z)\vnma\Rnma(\za)\big(\Hsalnm\big)^{-1}\ysnm
  \end{align*}
  für alle $n \in \Zefkp$. Dann heißt $\ABCDsalk$ das \textbf{linksseitige $\alpha$-Dyukarev-Quadrupel} bezüglich $\sjk$. Falls klar ist, von welchem $\sjk$ die Rede ist, lassen wir das \anf{$\sklam{s}$} als oberen Index weg.
\end{defi}

In der Situation von \thref{drldef2} bilden wir nun aus den dort eingeführten Folgen von \textit{q}$\times$\textit{q}-Matrixpolynomen eine spezielle Folge von $2$\textit{q}$\times2$\textit{q}-Matrixpolynomen.

\begin{defi}	\thlabel{drldef3}
  Seien $\kappa \in \Na$, $\alpha \in \R$, $\sjk \in \Lpqka$ und $\ABCDsalk$ das linksseitige $\alpha$-Dyukarev-Quadrupel bezüglich $\sjk$.
  Weiterhin sei
  \begin{align*}
    \Usalm := \begin{pmatrix} \Asalfm & \Bsalfm \\ \Csalfm & \Dsalfm \end{pmatrix}
  \end{align*}
  für alle $m \in \Zok$. Dann heißt $\Usalmk$ die \textbf{Folge von linksseitigen }$\mathbf{2}$\textbf{q}$\mathbf{\times2}$\textbf{q-$\alpha$-Dyukarev-Matrixpolynomen} bezüglich $\sjk$.
  Im Fall $\kappa \in \N$ heißt $\Usalk$ das \textbf{linksseitige }$\mathbf{2}$\textbf{q}$\mathbf{\times2}$\textbf{q-$\alpha$-Dyukarev-Matrixpolynom} bezüglich $\sjk$. Falls klar ist, von welchem $\sjk$ die Rede ist, lassen wir das \anf{$\sklam{s}$} als oberen Index weg.
\end{defi}

Unsere bereits erzielten Resultate erlauben uns eine alternative Darstellung der Folge von linksseitigen $2$\textit{q}$\times2$\textit{q}-$\alpha$-Dyukarev-Matrixpolynomen bezüglich einer linksseitig $\alpha$-Stieltjes-positiv definiten Folge.

\begin{satz}	\thlabel{drlbm2}
  Seien $\kappa \in \Na$, $\alpha \in \R$, $\sjk \in \Lpqka$ und $\Ualmk$ die Folge von linksseitigen $2$\textit{q}$\times2$\textit{q}-$\alpha$-Dyukarev-Matrixpolynomen bezüglich $\sjk$. Weiterhin sei
    \begin{align*}
      \tUalm(z) := \begin{pmatrix} (\alpha-z)\Iq & \Oq \\ \Oq & \Iq \end{pmatrix} \Ualm(z) \begin{pmatrix} (\alpha-z)^{-1}\Iq & \Oq \\ \Oq & \Iq \end{pmatrix}
    \end{align*}
    für alle $z\in\C\setminus\gklam{\alpha}$ und $m \in \Zok$.
  	\dgfa
  	\begin{itemize}
   		\item [\rm{(a)}] Es gelten $U_{\al2n+1} = \Taln$ für alle $n \in \Zofkm$ und $U_{\al2n} = \Yaln$ für alle \linebreak $n \in \Zefk$.
    	\item [\rm{(b)}] Es gelten $\widetilde{U}_{\al2n+1} = \tTaln$ für alle $n \in \Zofkm$ und $\widetilde{U}_{\al2n} = \tYaln$ für alle \linebreak $n \in \Zefk$.
    	\item [\rm{(c)}] Es gelten $\Ualm,\tUalm\in\tbtJqPp$ für alle $m \in \Zok$.
    	\item [\rm{(d)}] Für alle $m\in\Zok$ sind $\det\Ualm$ und $\det\tUalm$ konstante Funktionen auf $\C$ und deren Wert jeweils von Null verschieden. Insbesondere gelten
    	\begin{align*}
    		\Ualm^{-1}(z)=\tJq\Ualm^{\ast}(\za)\tJq \quad \text{und} \quad
    		\tUalm^{-1}(z)=\tJq\tUalm^{\ast}(\za)\tJq
    	\end{align*}
    	für alle $z\in\C$ und $m \in \Zok$.
  \end{itemize}
\end{satz}

\bwanf Zu (a): Dies folgt unter Beachtung von \thref{drldef3} und \thref{drldef2} aus Teil (a) von \thref{drlsa3} bzw. Teil (a) von \thref{drlsa4}.

Zu (b): Dies folgt unter Beachtung von Teil (b) von \thref{drlsa3} bzw. Teil (b) von \thref{drlsa4} aus (a). 

Zu (c): Wegen \thref{drldef3} und \thref{drldef2} gelten
\begin{align*}
	\Ualo(z) = \begin{pmatrix} \Iq & \Oq \\ (\alpha-z)s^{-1}_0 & \Iq \end{pmatrix} \quad \text{und} \quad \tUalo(z) = \begin{pmatrix} \Iq & \Oq \\ s^{-1}_0 & \Iq \end{pmatrix}
\end{align*}
für alle $z\in\C$. Hieraus folgen
\begin{align*}
	&\ \tJq-\Ualo^{\ast}(z)\tJq\Ualo(z) \\
	&= \begin{pmatrix} \Oq & -i\Iq \\ i\Iq & \Oq \end{pmatrix} - \begin{pmatrix} \Iq & (\alpha-\za)s^{-1}_0 \\ \Oq & \Iq \end{pmatrix} \begin{pmatrix} \Oq & -i\Iq \\ i\Iq & \Oq \end{pmatrix} \begin{pmatrix} \Iq & \Oq \\ (\alpha-z)s^{-1}_0 & \Iq \end{pmatrix} \\
	&= \begin{pmatrix} \Oq & -i\Iq \\ i\Iq & \Oq \end{pmatrix} - \begin{pmatrix} i(\alpha-\za)s^{-1}_0 & -i\Iq \\ i\Iq & \Oq \end{pmatrix} \begin{pmatrix} \Iq & \Oq \\ (\alpha-z)s^{-1}_0 & \Iq \end{pmatrix} \\
	&= \begin{pmatrix} \Oq & -i\Iq \\ i\Iq & \Oq \end{pmatrix} - \begin{pmatrix} i(z-\za)s^{-1}_0 & -i\Iq \\ i\Iq & \Oq \end{pmatrix} = 2\im z\begin{pmatrix} s^{-1}_0 & \Oq \\ \Oq & \Oq \end{pmatrix}
\end{align*}
und
\begin{align*}
	&\ \tJq-\tUalo^{\ast}(z)\tJq\tUalo(z) \\
	&= \begin{pmatrix} \Oq & -i\Iq \\ i\Iq & \Oq \end{pmatrix} - \begin{pmatrix} \Iq & s^{-1}_0 \\ \Oq & \Iq \end{pmatrix} \begin{pmatrix} \Oq & -i\Iq \\ i\Iq & \Oq \end{pmatrix} \begin{pmatrix} \Iq & \Oq \\ s^{-1}_0 & \Iq \end{pmatrix} \\
	&= \begin{pmatrix} \Oq & -i\Iq \\ i\Iq & \Oq \end{pmatrix} - \begin{pmatrix} is^{-1}_0 & -i\Iq \\ i\Iq & \Oq \end{pmatrix} \begin{pmatrix} \Iq & \Oq \\ s^{-1}_0 & \Iq \end{pmatrix} \\
	&= \begin{pmatrix} \Oq & -i\Iq \\ i\Iq & \Oq \end{pmatrix} - \begin{pmatrix} \Oq & -i\Iq \\ i\Iq & \Oq \end{pmatrix} = 0_{2q\times 2q}
\end{align*}
für alle $z\in\C$. Hieraus folgt wegen \thref{spbsp1} und \thref{spdef3} nun \linebreak $\Ualo,\tUalo\in\tbtJqPp$. Unter Beachtung von Teil (d) von \thref{drlsa3} bzw. Teil (b) von \thref{drlsa4} aus (a) und (b) weiterhin $\Ualm,\tUalm\in\tbtJqPp$ für alle $m \in \Zek$.. 

Zu (d): Sei $m\in\Zok$. Wegen \thref{drldef3} und \thref{drldef2} sind dann $\Ualm$ und $\tUalm$ jeweils $2$\textit{q}$\times2$\textit{q}-Matrixpolynome. Hieraus folgt wegen \thref{spbsp1} und \thref{spfo1} dann, dass $\det\Ualm$ und $\det\tUalm$ konstante Funktionen auf $\C$ sind und deren Wert jeweils von Null verschieden ist. Wegen (c), \thref{spbsp1} und Teil (a) von \thref{splm1} gelten
\begin{align*}
	\Ualm^{-1}(z)=\tJq\Ualm^{\ast}(\za)\tJq \quad \text{und} \quad
	\tUalm^{-1}(z)=\tJq\tUalm^{\ast}(\za)\tJq
\end{align*}
für alle $z\in\C$. \bwend

Wir werden nun das rechtsseitige $2$\textit{q}$\times2$\textit{q}-$\alpha$-Dyukarev-Matrixpolynom bezüglich einer rechtsseitig $\alpha$-Stieltjes-positiv definiten Folge verwenden, um das linksseitige $2$\textit{q}$\times2$\textit{q}-$\alpha$-Dyukarev-Matrixpolynom bezüglich einer zugehörigen linksseitig $-\alpha$-Stieltjes-positiv definiten Folge anders darzustellen.

\begin{lemma}	\thlabel{drllm4}
	Seien $\kappa \in \Na$, $\alpha \in \R$, $\sjk$ eine Folge aus $\Cqq$ und $t_j:=(-1)^js_j$ für alle $j\in\Zok$. Weiterhin sei $\tjk\in\Lpqkma$ oder $\sjk\in\Kpqka$ erfüllt. Dann gelten $\sjk\in\Kpqka$ bzw. $\tjk\in\Lpqkma$ sowie folgende Aussagen:
	\begin{itemize}
		\item [\rm{(a)}] Seien $m\in\Zok$, $\Uarm^{\sklam{s}}$ das rechtsseitige $2$\textit{q}$\times2$\textit{q}-$\alpha$-Dyukarev-Matrixpolynom bezüglich $\sjm$ und $U^{\sklam{t}}_{-\alm}$ das linksseitige $2$\textit{q}$\times2$\textit{q}-$\alpha$-Dyukarev-Matrixpolynom bezüglich $\tjm$. Dann gilt
		\begin{align*}
			U^{\sklam{t}}_{-\alm}(-z) = \Ve\Uarm^{\sklam{s}}(z)\Vea
		\end{align*}
		für alle $z\in\C$.
		\item [\rm{(b)}] Seien $m\in\Zok$ und $\tUarm^{\sklam{s}}$ bzw. $\widetilde{U}^{\sklam{t}}_{-\alm}$ definiert wie in \thref{drbm5} bzw. \thref{drlbm2}. Dann gilt
		\begin{align*}
			\widetilde{U}^{\sklam{t}}_{-\alm}(-z) = \Ve\tUarm^{\sklam{s}}(z)\Vea
		\end{align*}
		für alle $z\in\C$.	
	\end{itemize}
\end{lemma}

\bwanf Zu (a): Dies folgt wegen Teil (a) von \thref{drbm5} und Teil (a) von \thref{drlbm2} aus den Teilen (c) und (d) von \thref{drllm3}.

Zu (b): Dies folgt wegen Teil (b) von \thref{drbm5} und Teil (b) von \thref{drlbm2} aus den Teilen (c) und (d) von \thref{drllm3}. \bwend

Mithilfe von \thref{drllm4} können wir das rechtsseitige $\alpha$-Dyukarev-Quadrupel bezüglich einer rechtsseitig $\alpha$-Stieltjes-positiv definiten Folge verwenden, um das linksseitige $-\alpha$-Dyukarev-Quadrupel bezüglich einer zugehörigen linksseitig $-\alpha$-Stieltjes-positiv definiten Folge anders darzustellen.

\begin{lemma}	\thlabel{drllm5}
	Seien $\kappa \in \Na$, $\alpha \in \R$, $\sjk$ eine Folge aus $\Cqq$ und $t_j:=(-1)^js_j$ für alle $j\in\Zok$. Weiterhin sei $\tjk\in\Lpqkma$ oder $\sjk\in\Kpqka$ erfüllt. Dann gelten $\sjk\in\Kpqka$ bzw. $\tjk\in\Lpqkma$. Seien nun $\ABCDsark$ das rechtsseitige $\alpha$-Dyukarev-Quadrupel bezüglich $\sjk$ und $[(\mathbf{A}^{\sklam{t}}_{-\aln})^{\fklam{\kappa}}_{n=0}$,""$(\mathbf{B}^{\sklam{t}}_{-\aln})^{\fklam{\kappa+1}}_{n=0}$,""$(\mathbf{C}^{\sklam{t}}_{-\aln})^{\fklam{\kappa}}_{n=0}$,""$(\mathbf{D}^{\sklam{t}}_{-\aln})^{\fklam{\kappa+1}}_{n=0}]$ das linksseitige $-\alpha$-Dyukarev-Quadrupel bezüglich $\tjk$. Dann gelten
	\begin{align*}
		\mathbf{A}^{\sklam{t}}_{-\aln}(-z) = \Asarn(z) \quad \text{und} \quad
		\mathbf{C}^{\sklam{t}}_{-\aln}(-z) = -\Csarn(z)
	\end{align*}
	für alle $z\in\C$ und $n\in\Zofk$ sowie
	\begin{align*}
		\mathbf{B}^{\sklam{t}}_{-\aln}(-z) = -\Bsarn(z) \quad \text{und} \quad
		\mathbf{D}^{\sklam{t}}_{-\aln}(-z) = \Dsarn(z)
	\end{align*}
	für alle $z\in\C$ und $n\in\Zofkp$.
\end{lemma}

\bwanf Dies folgt aus \thref{drldef3} und Teil (a) von \thref{drllm4}. \bwend

Wir wollen uns nun die Determinanten der einzelnen \textit{q}$\times$\textit{q}-Matrixpolynome des linksseitigen $\alpha$-Dyukarev-Quadrupels bezüglich einer linksseitig $\alpha$-Stieltjes-positiv definiten Folge anschauen.

\begin{satz}	\thlabel{drlsa10}
	Seien $\kappa \in \Na$, $\alpha \in \R$ und $\sjk \in \Lpqka$. Weiterhin sei $\ABCDalk$ das linksseitige $\alpha$"=Dyukarev"=Quadrupel bezüglich \linebreak $\sjk$. Dann gelten $\det \Aaln(z) \neq 0$ für alle $z\in\C\setminus(-\infty,\alpha)$ und $n\in\Zofk$, \linebreak $\det \Baln(z) \neq 0$ für alle $z\in\C\setminus(-\infty,\alpha)$ und $n\in\Zefkp$, $\det \Caln(z) \neq 0$ für alle $z\in\C\setminus(-\infty,\alpha]$ und $n\in\Zofk$ sowie $\det \Daln(z) \neq 0$ für alle $z\in\C\setminus(-\infty,\alpha)$ und $n\in\Zofkp$.
\end{satz}

\bwanf Dies folgt sogleich aus \thref{drllm5} und Teil (b) von \thref{drsa12}. \bwend

\thref{drlbm2} erlaubt uns, das linksseitige $2$\textit{q}$\times2$\textit{q}-$\alpha$-Dyukarev-Matrixpolynom bezüglich einer endlichen linksseitig $\alpha$-Stieltjes-positiv definiten Folge mit den Potapovschen Fundamentalmatrizen in Verbindung zu bringen. Es wird sich herausstellen, dass das linksseitige $2$\textit{q}$\times2$\textit{q}-$\alpha$-Dyukarev-Matrixpolynom bezüglich jener Folge auch eine Resolventenmatrix für das zugehörige linksseitige $\alpha$-Stieltjes Momentenproblem ist, wie wir im Folgenden zeigen werden.

\begin{theo}	\thlabel{drlth2}
  Seien $m \in \N$, $\alpha \in \R$, $\sjm \in \Lpqma$ und $\Ualm$ das linksseitige $2$\textit{q}$\times2$\textit{q}-$\alpha$-Dyukarev-Matrixpolynom bezüglich $\sjm$.
  Es bezeichne
  \begin{align*}
    \Ualm = \begin{pmatrix} \Ualm^{(1,1)} & \Ualm^{(1,2)} \\ \Ualm^{(2,1)} & \Ualm^{(2,2)} \end{pmatrix}
  \end{align*}
  die \textit{q}$\times$\textit{q}-Blockzerlegung von $\Ualm$. \dgfa
  \begin{itemize}
    \item [\rm{(a)}] Sei $\phipsi \in \PtJqCma$. Dann ist $\det\big(\Ualm^{(2,1)}\phi+\Ualm^{(2,2)}\psi\big)$ nicht die Nullfunktion und
    \begin{align*}
      S := \rklam{\Ualm^{(1,1)}\phi+\Ualm^{(1,2)}\psi}\rklam{\Ualm^{(2,1)}\phi+\Ualm^{(2,2)}\psi}^{-1}
    \end{align*}
    gehört zu $\Sqmasmu$.
    \item [\rm{(b)}] Sei $S\in\Sqmasmu$. Dann existiert ein $\phipsi \in \dPtJqCma$ derart, dass
    \begin{align*}
      \det\eklam{\Ualm^{(2,1)}(z)\phi(z)+\Ualm^{(2,2)}(z)\psi(z)} \neq 0
    \end{align*}
    und
    \begin{align*}
      S(z) = \eklam{\Ualm^{(1,1)}(z)\phi(z)+\Ualm^{(1,2)}(z)\psi(z)}\eklam{\Ualm^{(2,1)}(z)\phi(z)+\Ualm^{(2,2)}(z)\psi(z)}^{-1}
    \end{align*}
    für alle $z \in \C\setminus(-\infty,\alpha]$ erfüllt sind.
    \item [\rm{(c)}] Seien $\binom{\phi_1}{\psi_1}, \binom{\phi_2}{\psi_2} \in \PtJqCma$. 
    Dann sind $\det\big(\Ualm^{(2,1)}\phi_1+\Ualm^{(2,2)}\psi_1\big)$ und  $\det\big(\Ualm^{(2,1)}\phi_2+\Ualm^{(2,2)}\psi_2\big)$ jeweils nicht die Nullfunktion und es sind folgende Aussagen äquivalent:
    \begin{itemize}
      \item [\rm{(i)}] Es gilt
      \begin{align*}
	&\ \rklam{\Ualm^{(1,1)}\phi_1+\Ualm^{(1,2)}\psi_1}\rklam{\Ualm^{(2,1)}\phi_1+\Ualm^{(2,2)}\psi_1}^{-1} \\
	& = \rklam{\Ualm^{(1,1)}\phi_2+\Ualm^{(1,2)}\psi_2}\rklam{\Ualm^{(2,1)}\phi_2+\Ualm^{(2,2)}\psi_2}^{-1}.
      \end{align*}
      \item [\rm{(ii)}] Es gilt $\sklam{\binom{\phi_1}{\psi_1}} = \sklam{\binom{\phi_2}{\psi_2}}$.
    \end{itemize}
  \end{itemize}
\end{theo}

\bwanf Seien $t_j:=(-1)^js_j$ für alle $j\in\Zom$. Wegen Teil (a) von \thref{asmbm3} gilt dann $\tjm\in\Kpqmma$. Weiterhin sei $U^{\sklam{t}}_{-\arm}$ das rechtsseitige $2$\textit{q}$\times2$\textit{q}-$\alpha$-Dyukarev-Matrixpolynom bezüglich $\tjm$. Es bezeichne
\begin{align*}
  U^{\sklam{t}}_{-\arm} = \begin{pmatrix} \big(U^{\sklam{t}}_{-\arm}\big)^{(1,1)} & \big(U^{\sklam{t}}_{-\arm}\big)^{(1,2)} \\ \big(U^{\sklam{t}}_{-\arm}\big)^{(2,1)} & \big(U^{\sklam{t}}_{-\arm}\big)^{(2,2)} \end{pmatrix}
\end{align*}
die \textit{q}$\times$\textit{q}-Blockzerlegung von $U^{\sklam{t}}_{-\arm}$. 

Zu (a): Seien $\widecheck{\phi}, \widecheck{\psi}:\C\setminus[-\alpha,\infty)\rightarrow\Cqq$ definiert gemäß $\widecheck{\phi}(z) := -\phi(-z)$ bzw. \linebreak $\widecheck{\psi}(z) := \psi(-z)$. Wegen Teil (a) von \thref{spbm3} ist dann $\bbinom{\widecheck{\phi}}{\widecheck{\psi}}\in{\cal P}^{(q,q)}_{-\tJq,\geq}(\C$,""$[-\alpha,\infty))$. Hieraus folgt wegen Teil (a) von \thref{drth1} dann, dass \linebreak $\det\beklam{(U^{\sklam{t}}_{-\arn})^{(2,1)}\widecheck{\phi}+(U^{\sklam{t}}_{-\arn})^{(2,2)}\widecheck{\psi}}$ nicht die Nullfunktion ist und
\begin{align*}
	T := \eklam{\big(U^{\sklam{t}}_{-\arm}\big)^{(1,1)}\widecheck{\phi}+\big(U^{\sklam{t}}_{-\arm}\big)^{(1,2)}\widecheck{\psi}}\eklam{\big(U^{\sklam{t}}_{-\arm}\big)^{(2,1)}\widecheck{\phi}+\big(U^{\sklam{t}}_{-\arm}\big)^{(2,2)}\widecheck{\psi}}^{-1}
\end{align*}
zu ${\cal S}_{0,q,[-\alpha,\infty)}[\tjm$,""$\leq]$ gehört. Hieraus folgt wegen Teil (a) von \thref{drllm4} nun, dass wegen
\begin{align*}
	&\ \big(\Ualm^{\sklam{s}}\big)^{(2,1)}(z)\phi(z)+\big(\Ualm^{\sklam{s}}\big)^{(2,2)}(z)\psi(z) \\
	&= -\big(U^{\sklam{t}}_{-\arm}\big)^{(2,1)}(-z)\brklam{{-\widecheck{\phi}}}(-z)+\big(U^{\sklam{t}}_{-\arm}\big)^{(2,2)}(-z)\widecheck{\psi}(-z) \\
	&= \big(U^{\sklam{t}}_{-\arm}\big)^{(2,1)}(-z)\widecheck{\phi}(-z)+\big(U^{\sklam{t}}_{-\arm}\big)^{(2,2)}(-z)\widecheck{\psi}(-z)
\end{align*}
für alle $z\in\C\setminus(-\infty,\alpha]$ dann $\det\beklam{(\Ualm^{\sklam{s}})^{(2,1)}\phi+(\Ualm^{\sklam{s}})^{(2,2)}\psi}$ nicht die Nullfunktion ist und eine diskrete Teilmenge $\D$ von $\C\setminus(-\infty,\alpha]$ existiert, so dass
\begin{align*}
	 S(z) &= \eklam{\big(\Ualm^{\sklam{s}}\big)^{(1,1)}(z)\phi(z)+\big(\Ualm^{\sklam{s}}\big)^{(1,2)}(z)\psi(z)} \\
	 &\quad \cdot\eklam{\big(\Ualm^{\sklam{s}}\big)^{(2,1)}(z)\phi(z)+\big(\Ualm^{\sklam{s}}\big)^{(2,2)}(z)\psi(z)}^{-1} \\
	&= \eklam{\big(U^{\sklam{t}}_{-\arm}\big)^{(1,1)}(-z)\brklam{{-\widecheck{\phi}}}(-z)-\big(U^{\sklam{t}}_{-\arm}\big)^{(1,2)}(-z)\widecheck{\psi}(-z)} \\
	&\quad \cdot\eklam{-\big(U^{\sklam{t}}_{-\arm}\big)^{(2,1)}(-z)\brklam{{-\widecheck{\phi}}}(-z)+\big(U^{\sklam{t}}_{-\arm}\big)^{(2,2)}(-z)\widecheck{\psi}(-z)}^{-1} \\
	&= -\eklam{\big(U^{\sklam{t}}_{-\arm}\big)^{(1,1)}(-z)\widecheck{\phi}(-z)+\big(U^{\sklam{t}}_{-\arm}\big)^{(1,2)}(-z)\widecheck{\psi}(-z)} \\
	&\quad \cdot\eklam{\big(U^{\sklam{t}}_{-\arm}\big)^{(2,1)}(-z)\widecheck{\phi}(-z)+\big(U^{\sklam{t}}_{-\arm}\big)^{(2,2)}(-z)\widecheck{\psi}(-z)}^{-1} = -T(-z)
\end{align*}
für alle $z\in\C\setminus((-\infty,\alpha]\cup\D)$ erfüllt ist. Hieraus folgt wegen des Identitätssatzes für meromorphe Funktionen (vergleiche z.\,B. im skalaren Fall \cite[Satz 10.3.2]{Funk}; im matriziellen Fall betrachtet man die einzelnen Einträge der Matrixfunktion) und Teil (d) von \thref{asmbm5} dann $S\in\Sqmasmu$.

Zu (b): Sei $\widecheck{S}: \C\setminus[-\alpha,\infty) \rightarrow \Cqq$ definiert gemäß $\widecheck{S}(z) := -S(-z)$. Wegen Teil (d) von \thref{asmbm5} ist dann $\widecheck{S}\in{\cal S}_{0,q,[-\alpha,\infty)}[\tjm$,""$\leq]$. Hieraus folgt wegen Teil (b) von \thref{drth1} dann, dass ein $\bbinom{\widetilde{\phi}}{\widetilde{\psi}}\in\widetilde{\cal P}^{(q,q)}_{-\tJq,\geq}(\C$,""$[-\alpha,\infty))$ derart existiert, dass $\widetilde{\phi}$ und $\widetilde{\psi}$ in $\C\setminus[-\alpha,\infty)$ holomorph sowie
\begin{align}	\label{drlth2w1}
	\det\eklam{\big(U^{\sklam{t}}_{-\arm}\big)^{(2,1)}(z)\widetilde{\phi}(z)+\big(U^{\sklam{t}}_{-\arm}\big)^{(2,2)}(z)\widetilde{\psi}(z)} \neq 0
\end{align}
und
\begin{align}	\label{drlth2w2}
	\widecheck{S}(z) &= \eklam{\big(U^{\sklam{t}}_{-\arm}\big)^{(1,1)}(z)\widetilde{\phi}(z)+\big(U^{\sklam{t}}_{-\arm}\big)^{(1,2)}(z)\widetilde{\psi}(z)} \notag \\
	&\quad \cdot\eklam{\big(U^{\sklam{t}}_{-\arm}\big)^{(2,1)}(z)\widetilde{\phi}(z)+\big(U^{\sklam{t}}_{-\arm}\big)^{(2,2)}(z)\widetilde{\psi}(z)}^{-1}
\end{align}
für alle $z \in \C\setminus[-\alpha,\infty)$ erfüllt sind. Seien $\phi, \psi:\C\setminus(-\infty,\alpha]\rightarrow\Cqq$ definiert gemäß $\phi(z) := -\widetilde{\phi}(-z)$ bzw. $\psi(z) := \widetilde{\psi}(-z)$. Wegen Teil (b) von \thref{spbm3} ist dann $\phipsi \in \dPtJqCma$. Weiterhin gilt wegen Teil (a) von \thref{drllm4} und \fref{drlth2w1} bzw. \fref{drlth2w2} nun
\begin{align*}
	&\ \det\eklam{\big(\Ualm^{\sklam{s}}\big)^{(2,1)}(z)\phi(z)+\big(\Ualm^{\sklam{s}}\big)^{(2,2)}(z)\psi(z)} \\
	&= \det\eklam{-\big(U^{\sklam{t}}_{-\arm}\big)^{(2,1)}(-z)\big({-\widetilde{\phi}}\big)(-z)+\big(U^{\sklam{t}}_{-\arm}\big)^{(2,2)}(-z)\widetilde{\psi}(-z)} \\
	&= \det\eklam{\big(U^{\sklam{t}}_{-\arm}\big)^{(2,1)}(-z)\widetilde{\phi}(-z)+\big(U^{\sklam{t}}_{-\arm}\big)^{(2,2)}(-z)\widetilde{\psi}(-z)} \neq 0
\end{align*}
bzw.
\begin{align*}
	&\ S(z) = -\widecheck{S}(-z) \\
	&= -\eklam{\big(U^{\sklam{t}}_{-\arm}\big)^{(1,1)}(-z)\widetilde{\phi}(-z)+\big(U^{\sklam{t}}_{-\arm}\big)^{(1,2)}(-z)\widetilde{\psi}(-z)} \notag \\
	&\quad \cdot\eklam{\big(U^{\sklam{t}}_{-\arm}\big)^{(2,1)}(-z)\widetilde{\phi}(-z)+\big(U^{\sklam{t}}_{-\arm}\big)^{(2,2)}(-z)\widetilde{\psi}(-z)}^{-1} \\
	&= -\eklam{\big(\Ualm^{\sklam{s}}\big)^{(1,1)}(z)\rklam{-\phi}(z)-\big(\Ualm^{\sklam{s}}\big)^{(1,2)}(z)\psi(z)} \\
	&\quad \cdot\eklam{-\big(\Ualm^{\sklam{s}}\big)^{(2,1)}(z)\rklam{-\phi}(z)+\Ualm^{(2,2)}(z)\psi(z)}^{-1} \\
	&= \eklam{\big(\Ualm^{\sklam{s}}\big)^{(1,1)}(z)\phi(z)+\big(\Ualm^{\sklam{s}}\big)^{(1,2)}(z)\psi(z)}\eklam{\big(\Ualm^{\sklam{s}}\big)^{(2,1)}(z)\phi(z)+\Ualm^{(2,2)}(z)\psi(z)}^{-1}
\end{align*}
für alle $z \in \C\setminus(-\infty,\alpha]$.

Zu (c): Wegen (a) ist $\det\beklam{(\Ualm^{\sklam{s}})^{(2,1)}\phi_1+(\Ualm^{\sklam{s}})^{(2,2)}\psi_1}$ für alle $j\in\gklam{1,2}$ nicht die Nullfunktion. Seien $\widecheck{\phi}_1, \widecheck{\psi}_1, \widecheck{\phi}_2, \widecheck{\psi}_2:\C\setminus[-\alpha,\infty)\rightarrow\Cqq$ definiert gemäß \linebreak $\widecheck{\phi}_1(z) := -\phi_1(-z)$, $\widecheck{\psi}_1(z) := \psi_1(-z)$, $\widecheck{\phi}_2(z) := -\phi_2(-z)$ bzw. $\widecheck{\psi}_2(z) := \psi_2(-z)$. Wegen Teil (a) von \thref{spbm3} sind dann $\bbinom{\widecheck{\phi}_1}{\widecheck{\psi}_1}, \bbinom{\widecheck{\phi}_2}{\widecheck{\psi}_2}\in{\cal P}^{(q,q)}_{-\tJq,\geq}(\C$,""$[-\alpha,\infty))$. Hieraus folgt wegen Teil (c) von \thref{drth1} dann, dass $\det\beklam{\big(U^{\sklam{t}}_{-\arm}\big)^{(2,1)}\widecheck{\phi}_j+\big(U^{\sklam{t}}_{-\arm}\big)^{(2,2)}\widecheck{\psi}_j}$ für alle $j\in\gklam{1,2}$ nicht die Nullfunktion ist und folgende Aussagen äquivalent sind:
\begin{itemize}
	\item [\rm{(iii)}] Es gilt
	\begin{align*}
		&\ \eklam{\big(U^{\sklam{t}}_{-\arm}\big)^{(1,1)}\widecheck{\phi}_1+\big(U^{\sklam{t}}_{-\arm}\big)^{(1,2)}\widecheck{\psi}_1}\eklam{\big(U^{\sklam{t}}_{-\arm}\big)^{(2,1)}\widecheck{\phi}_1+\big(U^{\sklam{t}}_{-\arm}\big)^{(2,2)}\widecheck{\psi}_1}^{-1} \\
		& = \eklam{\big(U^{\sklam{t}}_{-\arm}\big)^{(1,1)}\widecheck{\phi}_2+\big(U^{\sklam{t}}_{-\arm}\big)^{(1,2)}\widecheck{\psi}_2}\eklam{\big(U^{\sklam{t}}_{-\arm}\big)^{(2,1)}\widecheck{\phi}_2+\big(U^{\sklam{t}}_{-\arm}\big)^{(2,2)}\widecheck{\psi}_2}^{-1}.
	\end{align*}
	\item [\rm{(iv)}] Es gilt $\sklam{\bbinom{\widecheck{\phi}_1}{\widecheck{\psi}_1}} = \sklam{\bbinom{\widecheck{\phi}_2}{\widecheck{\psi}_2}}$.
\end{itemize}
Wegen Teil (a) von \thref{drllm4} gilt nun
\begin{align*}
	&\ \eklam{\big(\Ualm^{\sklam{s}}\big)^{(1,1)}(z)\phi_j(z)+\big(\Ualm^{\sklam{s}}\big)^{(1,2)}(z)\psi_j(z)} \\
	&\ \cdot\eklam{\big(\Ualm^{\sklam{s}}\big)^{(2,1)}(z)\phi_j(z)+\big(\Ualm^{\sklam{s}}\big)^{(2,2)}(z)\psi_j(z)}^{-1} \\
	&= \eklam{\big(U^{\sklam{t}}_{-\arm}\big)^{(1,1)}(-z)\brklam{{-\widecheck{\phi}_j}}(-z)-\big(U^{\sklam{t}}_{-\arm}\big)^{(1,2)}(-z)\widecheck{\psi}_j(-z)} \\
	&\quad \cdot\eklam{-\big(U^{\sklam{t}}_{-\arm}\big)^{(2,1)}(-z)\brklam{{-\widecheck{\phi}_j}}(-z)+\big(U^{\sklam{t}}_{-\arm}\big)^{(2,2)}(-z)\widecheck{\psi}_j(-z)}^{-1} \\
	&= -\eklam{\big(U^{\sklam{t}}_{-\arm}\big)^{(1,1)}(-z)\widecheck{\phi}_j(-z)+\big(U^{\sklam{t}}_{-\arm}\big)^{(1,2)}(-z)\widecheck{\psi}_j(-z)} \\
	&\quad \cdot\eklam{\big(U^{\sklam{t}}_{-\arm}\big)^{(2,1)}(-z)\widecheck{\phi}_j(-z)+\big(U^{\sklam{t}}_{-\arm}\big)^{(2,2)}(-z)\widecheck{\psi}_j(-z)}^{-1}
\end{align*}
für alle $z\in\C\setminus(-\infty,\alpha]$ und $j\in\gklam{1,2}$. Hieraus folgt die Äquivalenz von (i) und (iii). Weiterhin gilt
\begin{align*}
	\binom{\phi_j(z)}{\psi_j(z)} = -\Ve\begin{pmatrix}\widecheck{\phi}_j(-z)\\\widecheck{\psi}_j(-z)\end{pmatrix}
\end{align*}
für alle $z\in\C\setminus(-\infty,\alpha]$ und $j\in\gklam{1,2}$. Hieraus folgt wegen Teil (a) von \thref{asmlm1}, \thref{spdef4b} und \thref{spdef5b} dann die Äquivalenz von (ii) und (iv). Somit sind die Aussagen (i) und (ii) äquivalent. \bwend

\thref{drlth2} zeigt nun, dass über die in \thref{drldef3} eingeführte Begriffsbildung tatsächlich eine Resolventenmatrix für das linksseitige $\alpha$-Stieltjes Momentenproblem gewonnen wird (vergleiche Teil (b) von \thref{drdef0}).

Abschließend können wir folgende Beobachtung für Funktionen aus \linebreak $\Sqmasmu$ für eine Folge $\sjm \in \Lpqma$ mit $m \in \N$ vornehmen.

\begin{satz}	\thlabel{drlsa9}
	Seien $m \in \N$, $\alpha \in \R$ und $\sjm \in \Lpqma$. Weiterhin sei \linebreak $S\in\Sqmasmu$. Dann gilt $\det S(z) \neq 0$ für alle $z\in\C\setminus(-\infty,\alpha]$. Insbesondere ist $S(x)$ für alle $x\in(\alpha,\infty)$ eine positiv hermitesche Matrix.
\end{satz}

\bwanf Wegen \thref{adpbm1} ist $s_0 = H_0$ eine reguläre Matrix. Hieraus folgt wegen \thref{asmsa4} dann $\det S(z) \neq 0$ für alle $z\in\C\setminus(-\infty,\alpha]$. Hieraus folgt wegen $S\in\Soqma$ (vergleiche \thref{asmbz4}) dann $-S(x) > \Oq$ für alle $x\in(\alpha,\infty)$. \bwend

\subsection{Zwei extremale Elemente der Menge \texorpdfstring{$\Sqmasmu$}{Sqmasmu}} \label{EEl}

Seien $m \in \N$, $\alpha \in \R$ und $\sjm \in \Lpqma$ gegeben. Der vorliegende Teilabschnitt ist dann der Diskussion zweier ausgezeichneter Elemente der Menge $\Sqmasmu$ gewidmet. Hierbei handelt es sich um rationale \textit{q}$\times$\textit{q}-Matrixfunktionen, die durch spezielle Extremaleigenschaften auf dem Intervall $(\alpha,\infty)$ bezüglich der Löwner"=Halbordnung für hermitesche \textit{q}$\times$\textit{q}-Matrizen gekennzeichnet werden. Ausgangspunkt unserer Betrachtungen ist die folgende Beobachtung.

\begin{bem}	\thlabel{drlbsp1}
	Seien $\kappa \in \Na$, $\alpha \in \R$, $\sjk \in \Lpqka$ und $\ABCDalk$ das linksseitige $\alpha$-Dyukarev-Quadrupel bezüglich $\sjk$. Unter Beachtung von \thref{drlsa10} seien weiterhin $\Salmmin,\Salmmax:\C\setminus(-\infty,\alpha]\rightarrow\Cqq$ für alle $m\in\Zek$ definiert gemäß
	\begin{align*}
		\Salmmin(z) := \Aalfm(z)\Calfm^{-1}(z) \quad \text{und} \quad \Salmmax(z) := \Balfm(z)\Dalfm^{-1}(z).
	\end{align*}
	Dann gelten $\Salmmin,\Salmmax\in\Sqmasmu$ für alle $m\in\Zek$.
\end{bem}

\bwanf Sei $m\in\Zek$. Wegen Teil (b) von \thref{asmsa2} bzw. Teil (b) von \thref{asmdef3} gilt dann $\sjm \in \Lpqma$. Bezeichne ${\cal I}$ bzw. ${\cal O}$ die in $\C\setminus(-\infty,\alpha]$ konstante Matrixfunktion mit dem Wert $\Iq$ bzw. $\Oq$. Wegen \thref{spbsp3} gilt dann $\binom{\cal I}{\cal O}, \binom{\cal O}{\cal I} \in \PtJqCma$. Unter Beachtung von \thref{drldef3} folgt dann die Behauptung aus Teil (a) von \thref{drlth1}, wobei $\Salmmin$ bzw. $\Salmmax$ die zum \textit{q}$\times$\textit{q}-Stieltjes-Paar $\binom{\cal I}{\cal O}$ bzw. $\binom{\cal O}{\cal I}$  in $\C\setminus(-\infty,\alpha]$ zugehörige Funktion aus $\Sqmasmu$ darstellt. \bwend

\begin{defi}	\thlabel{drldef4}
	Seien $m \in \N$, $\alpha \in \R$, $\sjm \in \Lpqma$ und $\ABCDsalm$ das linksseitige $\alpha$-Dyukarev-Quadrupel bezüglich $\sjm$. Weiterhin seien $\Salmmins,\Salmmaxs:\C\setminus(-\infty,\alpha]\rightarrow\Cqq$ definiert gemäß
	\begin{align*}
		\Salmmins(z) := \Asalfm(z)\beklam{\Csalfm(z)}^{-1} \text{ bzw. } \Salmmaxs(z) := \Bsalfm(z)\beklam{\Dsalfm(z)}^{-1}.
	\end{align*}
	Dann heißt $\Salmmins$ bzw. $\Salmmaxs$ das \textbf{untere} bzw. \textbf{obere Extremalelement} von \linebreak $\Sqmasmu$. 
	Falls klar ist, von welchem $\sjm$ die Rede ist, lassen wir das \anf{$s$} im oberen Index weg.
\end{defi}

Wir können nun folgende Beobachtung über die zugehörigen Stieltjes-Maße der in \thref{drldef4} eingeführten Funktionen machen.

\begin{bem}	\thlabel{drlbm7}
	Seien $m \in \N$, $\alpha \in \R$, $\sjm \in \Lpqma$ und $\Salmmin$ bzw. $\Salmmax$ das untere bzw. obere Extremalelement von $\Sqmasmu$. Weiterhin sei $\mu^{(\alm)}_{min}$ bzw. $\mu^{(\alm)}_{max}$ das zu $\Salmmin$ bzw. $\Salmmax$ gehörige Stieltjes-Maß. \dgfa 
	\begin{itemize} 
		\item [\rm{(a)}] Es sind $\mu^{(\alm)}_{min}$ und $\mu^{(\alm)}_{max}$ molekulare Maße aus $\M^q_{\geq,\infty}((-\infty,\alpha])$.
		\item [\rm{(b)}] Sei $\ABCDalk$ das linksseitige $\alpha$-Dyukarev- \linebreak Quadrupel bezüglich $\sjk$. Es bezeichne ${\cal N}^{(\alm)}_{min}$ bzw. ${\cal N}^{(\alm)}_{max}$ die Nullstellenmenge von $\det\Calfm$ bzw. $\det\Dalfm$. Dann gelten
		\begin{align*}
			\mu^{(min)}_{\alm}\Brklam{(-\infty,\alpha]\setminus{\cal N}^{(\alm)}_{min}} = \Oq
			\quad \text{und} \quad
			\mu^{(max)}_{\alm}\Brklam{(-\infty,\alpha]\setminus{\cal N}^{(\alm)}_{max}} = \Oq.
		\end{align*}
	\end{itemize}
\end{bem}

\bwanf Dies folgt unter Beachtung von \thref{asmth5} aus \cite[Lemma B.4]{13} in Verbindung mit \thref{drldef4} und \thref{drlsa10}. \bwend

Wir wollen nun erläutern, warum die in \thref{drldef4} definierten Funktionen gerade das untere und obere Extremalelement von $\Sqmasmu$ genannt werden. Hierfür benötigen wir zunächst folgendes Hilfsresultat, das eine Verbindung zum rechtsseitigen Fall schafft.

\begin{lemma}	\thlabel{drllm6}
	Seien $m \in \N$, $\alpha \in \R$, $\sjm$ eine Folge aus $\Cqq$ und $t_j:=(-1)^js_j$ für alle $j\in\Zom$. Weiterhin sei $\tjm\in\Lpqmma$ oder $\sjm\in\Kpqma$ erfüllt. Dann gelten $\sjm\in\Kpqma$ bzw. $\tjm\in\Lpqmma$. Seien nun $\Sarmmins$ bzw. $\Sarmmaxs$ das untere bzw. obere Extremalelement von $\Sqasmu$ und $\Smin^{(-\alm,t)}$ bzw. $\Smax^{(-\alm,t)}$ das untere bzw. obere Extremalelement von $\Soqmma$""$[\tjm$,""$\leq]$. Dann gelten
	\begin{align*}
		\Smin^{(-\alm,t)}(-z) = -\Sarmmaxs(z) \quad \text{und} \quad
		\Smax^{(-\alm,t)}(-z) = -\Sarmmins(z)
	\end{align*}
	für alle $z\in\C\setminus[\alpha,\infty)$.
\end{lemma}

\bwanf Wegen Teil (a) von \thref{asmbm3} gilt $\sjm\in\Kpqma$ bzw. \linebreak $\tjm\in\Lpqmma$. Seien nun $[(\Asarn)^{\fklam{m}}_{n=0}$,""$(\Bsarn)^{\fklam{m+1}}_{n=0}$,""$(\Csarn)^{\fklam{m}}_{n=0}$,""$(\Dsarn)^{\fklam{m+1}}_{n=0}]$ das rechtsseitige $\alpha$-Dyukarev-Quadrupel bezüglich $\sjm$ und $[(\mathbf{A}^{\sklam{t}}_{-\aln})^{\fklam{m}}_{n=0}$,""$(\mathbf{B}^{\sklam{t}}_{-\aln})^{\fklam{m+1}}_{n=0}$, \linebreak $(\mathbf{C}^{\sklam{t}}_{-\aln})^{\fklam{m}}_{n=0}$,""$(\mathbf{D}^{\sklam{t}}_{-\aln})^{\fklam{m+1}}_{n=0}]$ das linksseitige $-\alpha$-Dyukarev-Quadrupel bezüglich $\tjm$. Wegen \thref{drldef4}, \thref{drllm5} und \thref{drdef4} gelten für alle $z\in\C\setminus[\alpha,\infty)$ dann
\begin{align*}
	\Smin^{(-\alm,t)}(-z) &= \mathbf{A}^{\sklam{t}}_{-\al\fklam{m}}(-z)\beklam{\mathbf{C}^{\sklam{t}}_{-\al\fklam{m}}(-z)}^{-1} \\
	&= -\Asarfm(z)\beklam{\Csarfm(z)}^{-1} = -\Sarmmaxs(z)
\end{align*}
und
\begin{align*}
	\Smax^{(-\alm,t)}(-z) &= \mathbf{B}^{\sklam{t}}_{-\al\fklam{m+1}}(-z)\beklam{\mathbf{D}^{\sklam{t}}_{-\al\fklam{m+1}}(-z)}^{-1} \\
	&= -\Bsarfm(z)\beklam{\Dsarfm(z)}^{-1} = -\Sarmmins(z). \tag*{$\Box$}
\end{align*}

Da wegen \thref{asmbz4} für jedes $S\in\Soqa$ die Matrix $S(x)$ für alle $x\in(\alpha,\infty)$ (nichtnegativ) hermitesch ist, können wir nun folgende Ungleichungen bezüglich der Löwner Halbordnung auf dem Intervall $(\alpha,\infty)$ für Funktionen aus $\Sqmasmu$ für eine Folge $\sjm \in \Lpqma$ mit $m\in\N$ betrachten.

\begin{satz}	\thlabel{drlsa6}
	Seien $m \in \N$, $\alpha \in \R$, $\sjm \in \Lpqma$ und $\Salmmin$ bzw. $\Salmmax$ das untere bzw. obere Extremalelement von $\Sqmasmu$.  \dgfa
	\begin{itemize}
		\item [\rm{(a)}] Es gilt
		\begin{align*}
			\Salmmin(x) < \Salmmax(x)
		\end{align*}
		für alle $x\in(\alpha,\infty)$.
		\item [\rm{(b)}] Sei $S\in\Sqmasmu$. Dann gilt
		\begin{align*}
			\Salmmin(x) \leq S(x) \leq \Salmmax(x)
		\end{align*}
		für alle $x\in(\alpha,\infty)$.
		\item [\rm{(c)}] Es ist $\beklam{\Salmmax-\Salmmin}^{-1}$ holomorph in $\C\setminus(-\infty,\alpha]$ und es gilt
		\begin{align*}
			\beklam{\Salmmax(z)-\Salmmin(z)}^{-1} &= -(\alpha-z)\vna\Rna(\za)\Hn^{-1}\Rn(z)\vn \\
			&\quad + (\alpha-z)^2\vnma\Rnma(\za)\Halnm^{-1}\Rnm(z)\vnm
		\end{align*}
		für alle $z\in\C\setminus(-\infty,\alpha]$.
	\end{itemize}
\end{satz}

\bwanf Seien $t_j:=(-1)^js_j$ für alle $j\in\Zom$. Wegen Teil (a) von \thref{asmbm3} gilt dann $\tjm\in\Kpqmma$. Weiterhin sei $\Smin^{(-\arm,t)}$ bzw. $\Smax^{(-\arm,t)}$ das untere bzw. obere Extremalelement von ${\cal S}_{0,q,[-\alpha,\infty)}[\tjm,\leq]$. 

Zu (a): Wegen Teil (a) von \thref{drsa8} gilt
\begin{align*}
	\Smin^{(-\arm,t)}(x) < \Smax^{(-\arm,t)}(x)
\end{align*}
für alle $x\in(-\infty,-\alpha)$. Hieraus folgt wegen \thref{drllm6} dann
\begin{align*}
	-\Salmmaxs(-x) < -\Salmmins(-x)
\end{align*}
für alle $x\in(-\infty,-\alpha)$. Hieraus folgt dann die Behauptung.

Zu (b): Sei $\widecheck{S}: \C\setminus[-\alpha,\infty) \rightarrow \Cqq$ definiert gemäß $\widecheck{S}(z) := -S(-z)$. Wegen Teil (d) von \thref{asmbm5} gilt dann $\widecheck{S}\in{\cal S}_{0,q,[-\alpha,\infty)}[\tjm,\leq]$. Wegen Teil (b) von \thref{drsa8} gilt weiterhin
\begin{align*}
	\Smin^{(-\arm,t)}(x) \leq \widecheck{S}(x) \leq \Smax^{(-\arm,t)}(x)
\end{align*}
für alle $x\in(-\infty,-\alpha)$. Hieraus folgt wegen \thref{drllm6} dann
\begin{align*}
	-\Salmmaxs(-x) \leq -S(-x) \leq -\Salmmins(-x)
\end{align*}
für alle $x\in(-\infty,-\alpha)$. Hieraus folgt dann die Behauptung. 

Zu (c): Wegen \thref{drllm6} gilt
\begin{align}	\label{drlsa6bw1}
	\Salmmaxs(z)-\Salmmins(z) = \Smax^{(-\arm,t)}(-z)-\Smin^{(-\arm,t)}(-z)
\end{align}
für alle $z\in\C\setminus(-\infty,\alpha]$. Wegen Teil (c) von \thref{drsa8} ist $\beklam{\Smax^{(-\arm,t)}-\Smin^{(-\arm,t)}}^{-1}$ holomorph in $\C\setminus[-\alpha,\infty)$. Hieraus folgt wegen \fref{drlsa6bw1} dann, dass $\beklam{\Salmmaxs-\Salmmins}^{-1}$ holomorph in $\C\setminus(-\infty,\alpha]$ ist. Wegen \fref{drlsa6bw1}, Teil (c) von \thref{drsa8}, der Teile (a) und (b) von \thref{asmlm1}, Teil (b) von \thref{asplm1} sowie der Teile (a) und (b) von \thref{drllm1} gilt weiterhin
\begin{align*}
	&\ \beklam{\Salmmaxs(z)-\Salmmins(z)}^{-1} = \beklam{\Smax^{(-\arm,t)}(-z)-\Smin^{(-\arm,t)}(-z)}^{-1} \\
	&= -(-z+\alpha)\vna\Rna(-\za)\rklam{\Hn^{\sklam{t}}}^{-1}\Rn(-z)\vn \\
	&\quad + (-z+\alpha)^2\vnma\Rnma(-\za)\rklam{H^{\sklam{t}}_{-\arn-1}}^{-1}\Rnm(-z)\vnm \\
	&= -(\alpha-z)\vna\Vna\Vn\Rna(\za)\Vna\Vn\rklam{\Hsn}^{-1}\Vna\Vn\Rn(z)\Vna\Vn\vn \\
	&\quad + (\alpha-z)^2\vnma\Vnma\Vnm\Rnma(\za)\Vnma\Vnm\rklam{\Hsalnm}^{-1}\Vnma\Vnm\Rnm(z)\Vnma\Vnm\vnm \\
	&= -(\alpha-z)\vna\Rna(\za)\rklam{\Hsn}^{-1}\Rn(z)\vn + (\alpha-z)^2\vnma\Rnma(\za)\rklam{\Hsalnm}^{-1}\Rnm(z)\vnm
\end{align*}
für alle $z\in\C\setminus(-\infty,\alpha]$. \bwend

Teil (b) von \thref{drlsa6} dokumentiert nun, in welcher Weise die beiden in \thref{drldef4} eingeführten rationalen \textit{q}$\times$\textit{q}-Matrixfunktionen eine extremale Position der Menge $\Sqmasmu$ für eine Folge $\sjm\in\Lpqma$ mit $m\in\N$ einnehmen.

Der folgende Satz liefert uns eine alternative Darstellung des unteren bzw. oberen Extremalelements von $\Sqmasmu$ für eine Folge $\sjm\in\Lpqma$ mit $m\in\N$, welche uns eine bessere Vorstellung von der Struktur dieser \textit{q}$\times$\textit{q}-Matrixfunktionen liefert. Zuvor benötigen wir aber noch folgendes Lemma.

\begin{lemma}	\thlabel{drllm7}
	Seien $\kappa\in\Na$, $\alpha\in\R$ und $\sjk\in\Lpqka$. Weiterhin sei 
	\begin{align*}
		\widetilde{H}_{\aln} := \begin{cases} \Haln & \text{falls } n\in\Zofkm \\
		\begin{pmatrix} \Halnm & y_{\aln,2n-1} \\ z_{\aln,2n-1} & \Oq \end{pmatrix} & \text{falls } 2n=\kappa \end{cases}
	\end{align*}
	für alle $n\in\Zofk$ \dgfa
	\begin{itemize}
		\item [\rm{(a)}] Es gelten
		\begin{align*}
			\det\eklam{(\alpha-z)H_{n}-H_{\aln}} \neq 0	
		\end{align*}
		für alle $z\in\C\setminus(-\infty,\alpha]$ und $n\in\Zofkm$ sowie
		\begin{align*}
			\det\beklam{(\alpha-z)^{-1}R^{-1}_{n}(\alpha)H_{n}R^{-\ast}_{n}(\alpha)-T_{n}\widetilde{H}_{\aln}T^{\ast}_{n}} \neq 0
		\end{align*}
		für alle $z\in\C\setminus(-\infty,\alpha]$ und $n\in\Zofk$.
		
		\item [\rm{(b)}] Es gelten
		\begin{align*}
			\det\eklam{\yna\eklam{(\alpha-z)H_{n}-H_{\aln}}^{-1}\yn} \neq 0
		\end{align*}
		für alle $z\in\C\setminus(-\infty,\alpha]$ und $n\in\Zofkm$ sowie
		\begin{align*}
			\det\Beklam{\vna\beklam{(\alpha-z)^{-1}R^{-1}_{n}(\alpha)H_{n}R^{-\ast}_{n}(\alpha)-T_{n}\widetilde{H}_{\aln}T^{\ast}_{n}}^{-1}\vn} \neq 0
		\end{align*}
		für alle $z\in\C\setminus(-\infty,\alpha]$ und $n\in\Zofk$.		
	\end{itemize}
\end{lemma}

\bwanf Seien $t_j:=(-1)^js_j$ für alle $j\in\Zok$. Wegen Teil (a) von \thref{asmbm3} gilt dann $\tjk\in\Kpqkma$. 

Zu (a): Wegen der Teile (a) und (b) von \thref{asmlm1}, Teil (b) von \thref{asplm1} sowie Teil (a) von \thref{drlm6} gilt
\begin{align*}
	\det\eklam{(\alpha-z)\Hsn-\Hsaln} &= \det\eklam{(-z+\alpha)\Vn\Hn^{\sklam{t}}\Vna-\Vn H^{\sklam{t}}_{-\arn}\Vna} \\
	&= \det\eklam{\Vn\eklam{(-z+\alpha)\Hn^{\sklam{t}}-H^{\sklam{t}}_{-\arn}}\Vna} \\
	&= (-1)^{(n+1)q}\det\eklam{H^{\sklam{t}}_{-\arn}-(-z+\alpha)\Hn^{\sklam{t}}} \neq 0
\end{align*}
für alle $z\in\C\setminus(-\infty,\alpha]$ und $n\in\Zofkm$. Wegen der Teile (a) und (b) von \thref{asmlm1}, Teil (b) von \thref{asplm1}, der Teile (a) und (b) von \thref{drllm1} sowie Teil (a) von \thref{drlm6} gilt weiterhin
\begin{align*}
	&\ \det\beklam{(\alpha-z)^{-1}\Rn^{-1}(\alpha)\Hsn\Rn^{-\ast}(\alpha)-\Tn\widetilde{H}^{\sklam{s}}_{\arn}\Tna} \\
	&= \det\big[(-z+\alpha)^{-1}\Vn\Rn^{-1}(-\alpha)\Vna\Vn\Hn^{\sklam{t}}\Vna\Vn\Rn^{-\ast}(-\alpha)\Vna \\
	&\quad -\Vn\Tn\Vna\Vn \widetilde{H}^{\sklam{t}}_{-\arn}\Vna\Vn\Tna\Vna\big] \\
	&= \det\eklam{\Vn\beklam{(-z+\alpha)^{-1}\Rn^{-1}(-\alpha)\Hn^{\sklam{t}}\Rn^{-\ast}(-\alpha)-\Tn \widetilde{H}^{\sklam{t}}_{-\arn}\Tna}\Vna} \\
	&= (-1)^{(n+1)q}\det\beklam{\Tn \widetilde{H}^{\sklam{t}}_{-\arn}\Tna-(-z+\alpha)^{-1}\Rn^{-1}(-\alpha)\Hn^{\sklam{t}}\Rn^{-\ast}(-\alpha)} \neq 0
\end{align*}
für alle $z\in\C\setminus(-\infty,\alpha]$ und $n\in\Zofk$. 

Zu (b): Wegen der Teile (a), (b) und (d) von \thref{asmlm1}, Teil (b) von \thref{asplm1} sowie Teil (b) von \thref{drlm6} gilt
\begin{align*}
	&\ \det\eklam{(\ysn)^{\ast}\beklam{(\alpha-z)\Hsn-\Hsaln}^{-1}\ysn} \\
	&= \det\eklam{(\yn^{\sklam{t}})^{\ast}\Vna\beklam{(-z+\alpha)\Vn\Hn^{\sklam{t}}\Vna-\Vn H^{\sklam{t}}_{-\arn}\Vna}^{-1}\Vn\yn^{\sklam{t}}} \\
	&= \det\eklam{(\yn^{\sklam{t}})^{\ast}\beklam{(-z+\alpha)\Hn^{\sklam{t}}-H^{\sklam{t}}_{-\arn}}^{-1}\yn^{\sklam{t}}} \\
	&= (-1)^q\det\eklam{(\yn^{\sklam{t}})^{\ast}\beklam{H^{\sklam{t}}_{-\arn}-(-z+\alpha)\Hn^{\sklam{t}}}^{-1}\yn^{\sklam{t}}} \neq 0
\end{align*}
für alle $z\in\C\setminus(-\infty,\alpha]$ und $n\in\Zofkm$. Wegen der Teile (a) und (b) von \thref{asmlm1}, Teil (b) von \thref{asplm1}, der Teile (a) und (b) von \thref{drllm1} sowie Teil (b) von \thref{drlm6} gilt weiterhin
\begin{align*}
	&\ \det\Beklam{\vna\eklam{(\alpha-z)^{-1}\Rn^{-1}(\alpha)\Hsn\Rn^{-\ast}(\alpha)-\Tn\widetilde{H}^{\sklam{s}}_{\arn}\Tna}^{-1}\vn} \\
	&= \det\Big[\vna\Vna\big[(-z+\alpha)^{-1}\Vn\Rn^{-1}(-\alpha)\Vna\Vn\Hn^{\sklam{t}}\Vna\Vn\Rn^{-\ast}(-\alpha)\Vna \\
	&\quad -\Vn\Tn\Vna\Vn \widetilde{H}^{\sklam{t}}_{-\arn}\Vna\Vn\Tna\Vna\big]^{-1}\Vn\vn\Big]  \\
	&= \det\eklam{\vna\beklam{(-z+\alpha)^{-1}\Rn^{-1}(-\alpha)\Hn^{\sklam{t}}\Rn^{-\ast}(-\alpha)-\Tn \widetilde{H}^{\sklam{t}}_{-\arn}\Tna}^{-1}\vn} \\
	&= (-1)^{q}\det\Beklam{\vna\beklam{\Tn \widetilde{H}^{\sklam{t}}_{-\arn}\Tna-(-z+\alpha)^{-1}\Rn^{-1}(-\alpha)\Hn^{\sklam{t}}\Rn^{-\ast}(-\alpha)}^{-1}\vn} \neq 0
\end{align*}
für alle $z\in\C\setminus(-\infty,\alpha]$ und $n\in\Zofk$. \bwend

\begin{satz}	\thlabel{drlsa7}
	Seien $m \in \N$, $\alpha \in \R$ und $\sjm \in \Lpqma$. Weiterhin sei $\Salmmin$ bzw. $\Salmmax$ das untere bzw. obere Extremalelement von $\Sqmasmu$ und $\widetilde{H}_{aln}$ für alle $n\in\Zom$ definiert wie in \thref{drllm7} (mit $m$ statt $\kappa$). Unter Beachtung von \thref{drllm7} gelten dann
	\begin{align*}
		\Salmmin(z) = \eklam{v^{\ast}_{\fklam{m}}\beklam{(\alpha-z)^{-1}R^{-1}_{\fklam{m}}(\alpha)H_{\fklam{m}}R^{-\ast}_{\fklam{m}}(\alpha)-T_{\fklam{m}}\widetilde{H}_{\al\fklam{m}}T^{\ast}_{\fklam{m}}}^{-1}v_{\fklam{m}}}^{-1}
	\end{align*}
	für alle $z\in\C\setminus(-\infty,\alpha]$ und
	\begin{align*}
		\Salmmax(z) = y^{\ast}_{0,\fklam{m-1}}\beklam{(\alpha-z)H_{\fklam{m-1}}-H_{\al\fklam{m-1}}}^{-1}y_{0,\fklam{m-1}}
	\end{align*}
	für alle $z\in\C\setminus(-\infty,\alpha]$.
\end{satz}

\bwanf Seien $t_j:=(-1)^js_j$ für alle $j\in\Zom$. Wegen Teil (a) von \thref{asmbm3} gilt dann $\tjm\in\Kpqmma$. Weiterhin seien $\Smin^{(-\arm,t)}$ bzw. $\Smax^{(-\arm,t)}$ das untere bzw. obere Extremalelement von ${\cal S}_{0,q,[-\alpha,\infty)}[\tjm,\leq]$. Wegen \thref{drllm6}, \thref{drsa9}, der Teile (a) und (b) von \thref{asmlm1}, Teil (b) von \thref{asplm1} sowie der Teile (a) und (b) von \thref{drllm1} gilt dann
\begin{align*}
	&\ \Salmmins(z) = -\Smax^{(-\arm,t)}(-z) \\
	&= -\eklam{v^{\ast}_{\fklam{m}}\beklam{T_{\fklam{m}}\widetilde{H}^{\sklam{t}}_{-\ar\fklam{m}}T^{\ast}_{\fklam{m}}-(-z+\alpha)^{-1}R^{-1}_{\fklam{m}}(-\alpha)H^{\sklam{t}}_{\fklam{m}}R^{-\ast}_{\fklam{m}}(-\alpha)}^{-1}v_{\fklam{m}}}^{-1} \\
	&= -\Big[v^{\ast}_{\fklam{m}}V^{\ast}_{\fklam{m}}\beklam{V_{\fklam{m}}T_{\fklam{m}}V^{\ast}_{\fklam{m}}V_{\fklam{m}}\widetilde{H}^{\sklam{s}}_{\al\fklam{m}}V^{\ast}_{\fklam{m}}V_{\fklam{m}}T^{\ast}_{\fklam{m}}V^{\ast}_{\fklam{m}} \\
	&\quad -(\alpha-z)^{-1}V_{\fklam{m}}R^{-1}_{\fklam{m}}(\alpha)V^{\ast}_{\fklam{m}}V_{\fklam{m}}H^{\sklam{s}}_{\fklam{m}}V^{\ast}_{\fklam{m}}V_{\fklam{m}}R^{-\ast}_{\fklam{m}}(\alpha)V^{\ast}_{\fklam{m}}}^{-1}V_{\fklam{m}}v_{\fklam{m}}\Big]^{-1} \\
	&= \eklam{v^{\ast}_{\fklam{m}}\beklam{(\alpha-z)^{-1}R^{-1}_{\fklam{m}}(\alpha)H^{\sklam{s}}_{\fklam{m}}R^{-\ast}_{\fklam{m}}(\alpha)-T_{\fklam{m}}\widetilde{H}^{\sklam{s}}_{\al\fklam{m}}T^{\ast}_{\fklam{m}}}^{-1}v_{\fklam{m}}}^{-1}
\end{align*}
für alle $z\in\C\setminus(-\infty,\alpha]$. Wegen \thref{drllm6}, \thref{drsa9}, der Teile (a), (b) und (d) von \thref{asmlm1} sowie Teil (b) von \thref{asplm1} gilt weiterhin
\begin{align*}
	&\ \Salmmaxs(z) = -\Smin^{(-\arm,t)}(-z) \\
	&= -(y^{\sklam{t}}_{0,\fklam{m-1}})^{\ast}\beklam{H^{\sklam{t}}_{-\ar\fklam{m-1}}-(-z+\alpha)H^{\sklam{t}}_{\fklam{m-1}}}^{-1}y^{\sklam{t}}_{0,\fklam{m-1}} \\
	&= -(y^{\sklam{s}}_{0,\fklam{m-1}})^{\ast}V^{\ast}_{\fklam{m-1}}\big[V_{\fklam{m-1}}H^{\sklam{s}}_{\al\fklam{m-1}}V^{\ast}_{\fklam{m-1}} \\
	&\quad -(\alpha-z)V_{\fklam{m-1}}H^{\sklam{s}}_{\fklam{m-1}}V^{\ast}_{\fklam{m-1}}\big]^{-1}V_{\fklam{m-1}}y^{\sklam{s}}_{0,\fklam{m-1}} \\
	&= (y^{\sklam{s}}_{0,\fklam{m-1}})^{\ast}\beklam{(\alpha-z)H^{\sklam{s}}_{\fklam{m-1}}-H^{\sklam{s}}_{\al\fklam{m-1}}}^{-1}y^{\sklam{s}}_{0,\fklam{m-1}}
\end{align*}
für alle $z\in\C\setminus(-\infty,\alpha]$. \bwend

Teil (a) von \thref{drlsa6} und \thref{drlsa9} führt uns auf folgende Definition.

\begin{defi}	\thlabel{drldef5}
	Seien $m \in \N$, $\alpha \in \R$, $\sjm \in \Lpqma$ und $\Salmmin$ bzw. $\Salmmax$ das untere bzw. obere Extremalelement von $\Sqmasmu$. Weiterhin sei $x\in(\alpha,\infty)$. Dann heißt $\beklam{\Salmmin(x),\Salmmax(x)}$ das zu $\sjm$ und $x$ zugehörige \textbf{Weylsche Intervall}.
\end{defi}

Das so eben eingeführte zu einer linksseitig $\alpha$-Stieltjes-positiv definiten Folge und einem Punkt $x\in(\alpha,\infty)$ gehörige Weylsche Intervall ist ein nichtdegeneriertes Matrixintervall bezüglich der Löwner-Halbordnung. Der folgende Satz zeigt nun, dass dieses Intervall gerade mit der Menge der Funktionswerte $S(x)$ aller Lösungen $S\in\Soqma[\sjm,$ $\leq]$ übereinstimmt.

\begin{satz}	\thlabel{drlbm6}
	Seien $m \in \N$, $\alpha \in \R$ und $\sjm \in \Lpqma$. Weiterhin sei $\big[\Salmmin(x)$, $\Salmmax(x)\big]$ für alle $x\in(\alpha,\infty)$ das zu $\sjm$ und $x$ zugehörige Weylsche Intervall. Dann gilt
	\begin{align*}
		 \gklam{S(x) \;\big|\; S\in\Soqma[\sjm,\leq]} = \beklam{\Salmmin(x),\Salmmax(x)}
	\end{align*}
	für alle $x\in(\alpha,\infty)$.
\end{satz}

\bwanf Dies folgt wegen \thref{asmbm4} und \thref{drllm6} sogleich aus \thref{drbm6}. \bwend

\newpage
\section{Die multiplikative Struktur der Folge von \texorpdfstring{$2$\textit{q}$\times2$\textit{q}}{2qx2q}-\texorpdfstring{$\alpha$}{a}-Dyukarev"=Matrixpolynomen bezüglich \texorpdfstring{$\alpha$}{a}-Stieltjes-positiv definiter Folgen} \label{chapfa}

Am Ausgangspunkt dieses Kapitels steht die in Teil (c) von \thref{drbm5} bzw. \thref{drlbm2} enthaltene Aussage, dass die Folge von $2$\textit{q}$\times2$\textit{q}-$\alpha$"=Dyukarev"=Matrixpolynomen bezüglich einer $\alpha$-Stieltjes-positiv definiten Folge eine Folge von Funktionen aus $\tbtJqPp$ ist. Für Funktionen der die Klasse $\tbtJqPp$ umfassenden Klasse $\btJqPp$ (vergleiche \thref{spdef3}) existiert eine tiefliegende Faktorisierungstheorie, welche in ihren Grundzügen auf V.\,P. Potapov \cite{Po} zurückgeht.

In diesem Zusammenhang sei bemerkt, dass V.\,P. Potapov in \cite{Po} eine beliebige \textit{p}$\times$\textit{p}-Signaturmatrix $J$ (vergleiche \thref{spdef1}) und statt $\Pp$ den offenen Einheitskreis ${\mathbb D} := \gklam{z\in\C \;| \abs{z} < 1}$ sowie in ${\mathbb D}$ meromorphe \textit{p}$\times$\textit{p}-Matrixfunktionen betrachtete. Eine solche Funktion $W$ heißt eine $J$-Potapov-Funktion bezüglich ${\mathbb D}$, falls für jeden in ${\mathbb D}$ gelegenen Holomorphiepunkt $z$ von $W$ dann $W(z)$ eine $J$-kontraktive Matrix ist. Es bezeichne dann $\mathfrak{P}_J({\mathbb D})$ die Menge aller $J$-Potapov-Funktionen bezüglich ${\mathbb D}$. Im Spezialfall $J=\Ip$ wird man hierbei genau auf die Menge ${\cal S}_{\pp}({\mathbb D})$ der  \textit{p}$\times$\textit{p}-Schur-Funktionen auf ${\mathbb D}$ (vergleiche Teil (a) von \thref{spdef6}) geführt. Für die Klasse ${\cal S}_{1\times1}({\mathbb D})$ wurde am Anfang des 20. Jahrhunderts eine spezielle Faktorisierungstheorie ausgearbeitet. Diese Theorie stand am Ausgangspunkt der Arbeit \cite{Po} von V.\,P. Potapov. Sein ursprüngliches Ziel bestand darin, diese Faktorisierungstheorie auf die Klasse ${\cal S}_{\pp}({\mathbb D})$ zu erweitern. Während der Arbeit an dieser Aufgabenstellung bemerkte er, dass eine Verallgemeinerung des klassischen skalaren Resultats sogar für die weitaus allgemeinere Klasse $\mathfrak{P}_J({\mathbb D})$ erreicht werden kann und gelangte auf diese Weise zu einer der bedeutendsten Resultate der matriziellen komplexen Funktionentheorie des 20. Jahrhunderts. 

Dieses Resultat von V.\,P. Potapov lässt sich unter Verwendung einer konformen Abbildung von ${\mathbb D}$ auf die offene obere Halbebene $\Pp$ übertragen. Die Behandlung von Matrixversionen von Interpolations- und Momentenproblemen lenkte dann die Aufmerksamkeit auf eine spezielle Teilklasse von $J$-Potapov-Funktionen, nämlich die Klasse der $J$-inneren Funktionen. Im Falle einer polynomialen $J$-Potapov-Funktion besteht deren multiplikative Zerlegung in $J$-Elementarfunktionen lediglich aus linearen Polynomen, im Falle eines $J$-inneren Polynoms dann sogar aus $J$-inneren linearen Polynomen.

Eine solche Situation trifft dann also insbesondere für die Folge der $2$\textit{q}$\times2$\textit{q}-$\alpha$"=Dyukarev"=Matrixpolynomen zu. Im Mittelpunkt dieses Abschnitts steht eine elementare Herleitung der multiplikativen Zerlegung der $2$\textit{q}$\times2$\textit{q}-$\alpha$"=Dyukarev"=Matrixpolynome in ein Produkt von linearen Polynomen aus der Klasse $\tbtJqPp$. Unter elementarer Herleitung verstehen wir hierbei einen Zugang, der keinen Zugriff auf den allgemeinen Faktorisierungssatz von V.\,P. Potapov benötigt. Unser Zugang basiert auf der Herleitung entsprechender Rekursionen für die Polynome des $\alpha$-Dyukarev-Quadrupels von Folgen von \textit{q}$\times$\textit{q}-Matrixpolynomen. Das sind gerade diejenigen Matrixpolynome, welche aus der \textit{q}$\times$\textit{q}-Blockzerlegung der Folge von $2$\textit{q}$\times2$\textit{q}-$\alpha$"=Dyukarev"=Matrixpolynomen hervorgehen. Es wird sich zeigen, dass diese Rekursion gerade durch die $\alpha$-Dyukarev-Stieltjes-Parametrisierung der zugrundeliegenden $\alpha$-Stieltjes-positiv definiten Folge vermittelt werden.

Hierbei übernehmen wir im rechtsseitigen Fall die in \cite[Chapter 3]{CR1} für den Fall einer gegebenen Folge aus ${\cal K}^{>}_{q,\infty,0}$ verwendete Vorgehensweise und werfen auch einen Blick auf die Matrixpolynomstruktur des rechtsseitigen $\alpha$-Dyukarev-Quadrupels bezüglich einer rechtsseitig $\alpha$-Stieltjes-positiv definiten Folge. Im linksseitigen Fall greifen wir auf die Ergebnisse des rechtsseitigen Falles zurück.

\subsection{Der rechtsseitige Fall}

Unsere nachfolgende Überlegung ist darauf gerichtet, den Nachweis dafür zu erbringen, dass die Matrixpolynome des rechtsseitigen $\alpha$-Dyukarev-Quadrupels bezüglich einer rechtsseitig $\alpha$-Stieltjes-positiv definiten Folge durch spezielle Rekursionsbeziehungen miteinander verbunden sind, wobei die Koeffizienten in diesen Rekursionsformeln gerade über die rechtsseitige $\alpha$-Dyukarev-Stieltjes-Parametrisierung der zugrundeliegenden rechtsseitig $\alpha$-Stieltjes-positiv definiten Folge gewonnen werden. Dies ist der Inhalt des nachfolgenden Lemmas.

\begin{lemma}	\thlabel{drbm2}
  Seien $\kappa \in \Na$, $\alpha \in \R$, $\sjk \in \Kpqka$ und $[(\Aarn)^{\fklam{\kappa}}_{n=0}$,""$(\Barn)^{\fklam{\kappa+1}}_{n=0}$, \linebreak $(\Carn)^{\fklam{\kappa}}_{n=0}$,""$(\Darn)^{\fklam{\kappa+1}}_{n=0}]$ das rechtsseitige $\alpha$-Dyukarev-Quadrupel bezüglich $\sjk$.
  Weiterhin sei $\LMkarn$ die rechtsseitige $\alpha$-Dyukarev-Stieltjes-""Parametrisierung von  $\sjk$. Dann gelten
  \begin{align*}
    \Aarn(z)\Larn+\Barn(z) = \Barnp(z), \\
    \Carn(z)\Larn+\Darn(z) = \Darnp(z)
  \end{align*}
  für alle $z \in \C$ und $n \in \Zofkm$ und im Fall $\kappa\geq2$
  \begin{align*}
    \Aarnm(z)-(z-\alpha)\Barn(z)\Marn = \Aarn(z), \\
    \Carnm(z)-(z-\alpha)\Darn(z)\Marn = \Carn(z)
  \end{align*}
  für alle $z \in \C$ und $n \in \Zefk$.
\end{lemma}

\bwanf Sei $z \in \C$. Wegen \thref{drdef1} und \thref{adpdef1} gelten
\begin{align*}
  \Aaro(z)\Laro+\Baro(z) = \Iq s_0s^{-1}_{\aro}s_0 + \Oq = u^{\ast}_{\aro}R^{\ast}_0(\za)\Haro^{-1}y_{0,0} = \mathbf{B}_{\ar 1}
\end{align*}
und
\begin{align*}
  \Caro(z)\Laro+\Daro(z) & = -(z-\alpha)s^{-1}_0 s_0s^{-1}_{\aro}s_0 + \Iq \\
  & = \Iq - (z-\alpha)v^{\ast}_0R^{\ast}_0(\za)\Haro^{-1}y_{0,0} = \mathbf{D}_{\ar 1}.
\end{align*}

Sei nun $n \in \Zefkm$. Wegen Teil (c) von \thref{drbm1} gilt
\begin{align}	\label{drbm2bw1}
  R_n(\alpha)v_ny^{\ast}_{0,n} = H_n-R_n(\alpha)T_n\Harn.
\end{align}
Hieraus folgt wegen Teil (a) von \thref{drlm3} dann
\begin{align}	\label{drbm2bw2}
  &\ (z-\alpha)u^{\ast}_nR^{\ast}_n(\za)H^{-1}_nR_n(\alpha)v_ny^{\ast}_{0,n}\Harn^{-1}y_{0,n} \notag \\
  & = (z-\alpha)u^{\ast}_nR^{\ast}_n(\za)\Harn^{-1}y_{0,n}-(z-\alpha)u^{\ast}_nR^{\ast}_n(\za)H^{-1}_nR_n(\alpha)T_ny_{0,n} \notag \\
  & = (z-\alpha)u^{\ast}_nR^{\ast}_n(\za)\Harn^{-1}y_{0,n}-(z-\alpha)u^{\ast}_nR^{\ast}_n(\za)H^{-1}_nR_n(\alpha)u_n.
\end{align}
Wegen Teil (d) von \thref{drbm1} gilt
\begin{align}	\label{drbm2bw3}
  R_n(\alpha)v_ny^{\ast}_{0,n-1} = H_n\dL_n-R_n(\alpha)L_n\Harnm.
\end{align}
Hieraus folgt wegen Teil (d) von \thref{drlm2} sowie der Teile (c) und (b) von \thref{drlm3} dann
\begin{align}	\label{drbm2bw4}
  &\ (z-\alpha)u^{\ast}_nR^{\ast}_n(\za)H^{-1}_nR_n(\alpha)v_ny^{\ast}_{0,n-1}\Harnm^{-1}y_{0,n-1} \notag \\
  & = (z-\alpha)u^{\ast}_nR^{\ast}_n(\za)\dL_n\Harnm^{-1}y_{0,n-1}-(z-\alpha)u^{\ast}_nR^{\ast}_n(\za)H^{-1}_nR_n(\alpha)L_ny_{0,n-1} \notag \\
  & = (z-\alpha)u^{\ast}_{n-1}R^{\ast}_{n-1}(\za)\Harnm^{-1}y_{0,n-1}-(z-\alpha)u^{\ast}_nR^{\ast}_n(\za)H^{-1}_nR_n(\alpha)u_n.
\end{align}
Wegen Teil (d) von \thref{drlm1} sowie der Teile (e) und (a) von \thref{drlm3} gelten
\begin{align}	\label{drbm2bw5}
  u^{\ast}_{\arn}\eklam{R^{\ast}_n(\za)-R^{\ast}_n(\alpha)} = (z-\alpha)u^{\ast}_{\arn}R^{\ast}_n(\alpha)T^{\ast}_nR^{\ast}_n(\za) = (z-\alpha)u^{\ast}_nR^{\ast}_n(\za)
\end{align}
und
\begin{align}	\label{drbm2bw6}
  u^{\ast}_{\arn-1}\eklam{R^{\ast}_{n-1}(\za)-R^{\ast}_{n-1}(\alpha)} = (z-\alpha)u^{\ast}_{n-1}R^{\ast}_{n-1}(\za).
\end{align}
Wegen \thref{drdef1}, \thref{adpdef1}, Teil (e) von \thref{drlm3}, \fref{drbm2bw2}, \fref{drbm2bw4}, \fref{drbm2bw5} und \fref{drbm2bw6} gilt nun
\begin{align*}
  &\ \Aarn(z)\Larn+\Barn(z) \\
  & = \eklam{\Iq+(z-\alpha)u^{\ast}_nR^{\ast}_n(\za)H^{-1}_nR_n(\alpha)v_n}\eklam{y^{\ast}_{0,n}\Harn^{-1}y_{0,n} - y^{\ast}_{0,n-1}\Harnm^{-1}y_{0,n-1}} \\
  &\quad + u^{\ast}_{\arn-1}R^{\ast}_{n-1}(\za)\Harnm^{-1}y_{0,n-1} \\
  & = y^{\ast}_{0,n}\Harn^{-1}y_{0,n} - y^{\ast}_{0,n-1}\Harnm^{-1}y_{0,n-1} + (z-\alpha)u^{\ast}_nR^{\ast}_n(\za)H^{-1}_nR_n(\alpha)v_ny^{\ast}_{0,n}\Harn^{-1}y_{0,n} \\
  &\quad - (z-\alpha)u^{\ast}_nR^{\ast}_n(\za)H^{-1}_nR_n(\alpha)v_ny^{\ast}_{0,n-1}\Harnm^{-1}y_{0,n-1} + u^{\ast}_{\arn-1}R^{\ast}_{n-1}(\za)\Harnm^{-1}y_{0,n-1} \\
  & = u^{\ast}_{\arn}R^{\ast}_n(\alpha)\Harn^{-1}y_{0,n} - u^{\ast}_{\arn-1}R^{\ast}_{n-1}(\alpha)\Harnm^{-1}y_{0,n-1} \\
  &\quad + (z-\alpha)u^{\ast}_nR^{\ast}_n(\za)\Harn^{-1}y_{0,n} - (z-\alpha)u^{\ast}_{n-1}R^{\ast}_{n-1}(\za)\Harnm^{-1}y_{0,n-1} \\
  &\quad + u^{\ast}_{\arn-1}R^{\ast}_{n-1}(\za)\Harnm^{-1}y_{0,n-1} \\
  & = u^{\ast}_{\arn}R^{\ast}_n(\alpha)\Harn^{-1}y_{0,n} - u^{\ast}_{\arn-1}R^{\ast}_{n-1}(\alpha)\Harnm^{-1}y_{0,n-1} \\
  &\quad + u^{\ast}_{\arn}(R^{\ast}_n(\za)-R^{\ast}_n(\alpha))\Harn^{-1}y_{0,n} - u^{\ast}_{\arn-1}(R^{\ast}_{n-1}(\za)-R^{\ast}_{n-1}(\alpha))\Harnm^{-1}y_{0,n-1} \\
  &\quad + u^{\ast}_{\arn-1}R^{\ast}_{n-1}(\za)\Harnm^{-1}y_{0,n-1} \\
  & = u^{\ast}_{\arn}R^{\ast}_n(\za)\Harn^{-1}y_{0,n} = \Barnp(z).
\end{align*}
Wegen \fref{drbm2bw1} und Teil (a) von \thref{drlm3} gilt
\begin{align}	\label{drbm2bw7}
  &\ (z-\alpha)v^{\ast}_nR^{\ast}_n(\za)H^{-1}_nR_n(\alpha)v_ny^{\ast}_{0,n}\Harn^{-1}y_{0,n} \notag \\
  & = (z-\alpha)v^{\ast}_nR^{\ast}_n(\za)\Harn^{-1}y_{0,n} - (z-\alpha)v^{\ast}_nR^{\ast}_n(\za)H^{-1}_nR_n(\alpha)T_ny_{0,n} \notag \\
  & = (z-\alpha)v^{\ast}_nR^{\ast}_n(\za)\Harn^{-1}y_{0,n} - (z-\alpha)v^{\ast}_nR^{\ast}_n(\za)H^{-1}_nR_n(\alpha)u_n.
\end{align}
Wegen \fref{drbm2bw3}, der Teile (d) und (b) von \thref{drlm2} sowie Teil (b) von \thref{drlm3} gilt
\begin{align}	\label{drbm2bw8}
  &\ (z-\alpha)v^{\ast}_nR^{\ast}_n(\za)H^{-1}_nR_n(\alpha)v_ny^{\ast}_{0,n-1}\Harnm^{-1}y_{0,n-1} \notag \\
  & = (z-\alpha)v^{\ast}_nR^{\ast}_n(\za)\dL_n\Harnm^{-1}y_{0,n-1} - (z-\alpha)v^{\ast}_nR^{\ast}_n(\za)H^{-1}_nR_n(\alpha)L_ny_{0,n-1} \notag \\
  & = (z-\alpha)v^{\ast}_{n-1}R^{\ast}_{n-1}(\za)\Harnm^{-1}y_{0,n-1} - (z-\alpha)v^{\ast}_nR^{\ast}_n(\za)H^{-1}_nR_n(\alpha)u_n.
\end{align}
Wegen \thref{drdef1}, \thref{adpdef1}, \fref{drbm2bw7} und \fref{drbm2bw8} gilt
\begin{align*}
  &\ \Carn(z)\Larn+\Darn(z) \\
  & = -(z-\alpha)v^{\ast}_nR^{\ast}_n(\za)H^{-1}_nR_n(\alpha)v_n\eklam{y^{\ast}_{0,n}\Harn^{-1}y_{0,n} - y^{\ast}_{0,n-1}\Harnm^{-1}y_{0,n-1}} \\
  &\quad + \Iq-(z-\alpha)v^{\ast}_{n-1}R^{\ast}_{n-1}(\za)\Harnm^{-1}y_{0,n-1} \\
  & = \Iq-(z-\alpha)v^{\ast}_nR^{\ast}_n(\za)H^{-1}_nR_n(\alpha)v_ny^{\ast}_{0,n}\Harn^{-1}y_{0,n} \\
  &\quad + (z-\alpha)v^{\ast}_nR^{\ast}_n(\za)H^{-1}_nR_n(\alpha)v_ny^{\ast}_{0,n-1}\Harnm^{-1}y_{0,n-1} \\
  &\quad - (z-\alpha)v^{\ast}_{n-1}R^{\ast}_{n-1}(\za)\Harnm^{-1}y_{0,n-1} \\
  & = \Iq-(z-\alpha)v^{\ast}_nR^{\ast}_n(\za)\Harn^{-1}y_{0,n} + (z-\alpha)v^{\ast}_{n-1}R^{\ast}_{n-1}(\za)\Harnm^{-1}y_{0,n-1} \\
  &\quad - (z-\alpha)v^{\ast}_{n-1}R^{\ast}_{n-1}(\za)\Harnm^{-1}y_{0,n-1} \\
  & = \Iq-(z-\alpha)\vna\Rna(\za)\Harn^{-1}\yn = \Darnm(z).
\end{align*}

Seien nun $\kappa\geq2$ und $n \in \Zefk$. Wegen Teil (d) von \thref{drbm1} gilt
\begin{align}	\label{drbm2bw9}
  y_{0,n-1}v^{\ast}_nR^{\ast}_n(\alpha) = \dL^{\ast}_nH_n-\Harnm L^{\ast}_nR^{\ast}_n(\alpha).
\end{align}
Hieraus folgt wegen der Teile (d) und (b) von \thref{drlm2}, der Teile (b) und (e) von \thref{drlm3} sowie Teil (b) von \thref{drlm1} dann
\begin{align}	\label{drbm2bw10}
  &\ (z-\alpha)\uarnm^{\ast}R^{\ast}_{n-1}(\za)\Harnm^{-1}y_{0,n-1}v^{\ast}_nR^{\ast}_n(\alpha)H^{-1}_nR_n(\alpha)v_n \notag \\
  & = (z-\alpha)\uarnm^{\ast}R^{\ast}_{n-1}(\za)\Harnm^{-1}\dL^{\ast}_nR_n(\alpha)v_n \notag \\
  &\quad - (z-\alpha)\uarnm^{\ast}R^{\ast}_{n-1}(\za)L^{\ast}_nR^{\ast}_n(\alpha)H^{-1}_nR_n(\alpha)v_n \notag \\
  & = (z-\alpha)\uarnm^{\ast}R^{\ast}_{n-1}(\za)\Harnm^{-1}R_{n-1}(\alpha)v_{n-1} - (z-\alpha)u^{\ast}_nR^{\ast}_n(\za)H^{-1}_nR_n(\alpha)v_n.
\end{align}
Wegen Teil (c) von \thref{drbm1} gilt
\begin{align}	\label{drbm2bw11}
  y_{0,n-1}v^{\ast}_{n-1}R^{\ast}_{n-1}(\alpha) = H_{n-1}-\Harnm T^{\ast}_{n-1}R^{\ast}_{n-1}(\alpha).
\end{align}
Hieraus folgt wegen der Teile (a) und (b) von \thref{drlm1} sowie der Teile (e) und (a) von \thref{drlm3} dann
\begin{align}	\label{drbm2bw12}
  &\ (z-\alpha)\uarnm^{\ast}R^{\ast}_{n-1}(\za)\Harnm^{-1}y_{0,n-1}v^{\ast}_{n-1}R^{\ast}_{n-1}(\alpha)H^{-1}_{n-1}R_{n-1}(\alpha)v_{n-1} \notag \\
  & = (z-\alpha)\uarnm^{\ast}R^{\ast}_{n-1}(\za)\Harnm^{-1}R_{n-1}(\alpha)v_{n-1} \notag \\
  &\quad - (z-\alpha)\uarnm^{\ast}R^{\ast}_{n-1}(\za)T^{\ast}_{n-1}R^{\ast}_{n-1}(\alpha)H^{-1}_{n-1}R_{n-1}(\alpha)v_{n-1} \notag \\
  & = (z-\alpha)\uarnm^{\ast}R^{\ast}_{n-1}(\za)\Harnm^{-1}R_{n-1}(\alpha)v_{n-1} - (z-\alpha)u^{\ast}_{n-1}R^{\ast}_{n-1}(\za)H^{-1}_{n-1}R_{n-1}(\alpha)v_{n-1}.
\end{align}
Wegen \thref{drdef1}, \thref{adpdef1}, Teil (a) von \thref{sqlm1}, \fref{drbm2bw10} und \fref{drbm2bw12} gilt
\begin{align*}
  &\ \Aarnm(z)-(z-\alpha)\Barn(z)\Marn \\
  & = \Iq + (z-\alpha)u^{\ast}_{n-1}R^{\ast}_{n-1}(\za)H^{-1}_{n-1}R_{n-1}(\alpha)v_{n-1} - (z-\alpha)\uarnm^{\ast}R^{\ast}_{n-1}(\za)\Harnm^{-1}y_{0,n-1} \\
  &\quad \cdot \eklam{v^{\ast}_nR^{\ast}_n(\alpha)H^{-1}_nR_n(\alpha)v_n - v^{\ast}_{n-1}R^{\ast}_{n-1}(\alpha)H^{-1}_{n-1}R_{n-1}(\alpha)v_{n-1}} \\
  & = \Iq + (z-\alpha)u^{\ast}_{n-1}R^{\ast}_{n-1}(\za)H^{-1}_{n-1}R_{n-1}(\alpha)v_{n-1} \\
  &\quad - (z-\alpha)\uarnm^{\ast}R^{\ast}_{n-1}(\za)\Harnm^{-1}y_{0,n-1}v^{\ast}_nR^{\ast}_n(\alpha)H^{-1}_nR_n(\alpha)v_n \\
  &\quad + (z-\alpha)\uarnm^{\ast}R^{\ast}_{n-1}(\za)\Harnm^{-1}y_{0,n-1}v^{\ast}_{n-1}R^{\ast}_{n-1}(\alpha)H^{-1}_{n-1}R_{n-1}(\alpha)v_{n-1} \\
  & = \Iq + (z-\alpha)u^{\ast}_{n-1}R^{\ast}_{n-1}(\za)H^{-1}_{n-1}R_{n-1}(\alpha)v_{n-1} + (z-\alpha)u^{\ast}_nR^{\ast}_n(\za)H^{-1}_nR_n(\alpha)v_n \\
  &\quad - (z-\alpha)u^{\ast}_{n-1}R^{\ast}_{n-1}(\za)H^{-1}_{n-1}R_{n-1}(\alpha)v_{n-1} \\
  & = \Iq +(z-\alpha)\una\Rna(\za)\Hn^{-1}\Rn(\alpha)\vn = \Aarn(z).
\end{align*}
Wegen \fref{drbm2bw9}, der Teile (d), (b) und (a) von \thref{drlm2} sowie der Teile (a) und (d) von \thref{drlm1} gilt
\begin{align}	\label{drbm2bw13}
  &\ (z-\alpha)^2v^{\ast}_{n-1}R^{\ast}_{n-1}(\za)\Harnm^{-1}y_{0,n-1}v^{\ast}_nR^{\ast}_n(\alpha)H^{-1}_nR_n(\alpha)v_n \notag \\
  & = (z-\alpha)^2v^{\ast}_{n-1}R^{\ast}_{n-1}(\za)\Harnm^{-1}\dL^{\ast}_nR_n(\alpha)v_n - (z-\alpha)^2v^{\ast}_{n-1}R^{\ast}_{n-1}(\za)L^{\ast}_nR^{\ast}_n(\alpha)H^{-1}_nR_n(\alpha)v_n \notag \\
  & = (z-\alpha)^2v^{\ast}_{n-1}R^{\ast}_{n-1}(\za)\Harnm^{-1}R_{n-1}(\alpha)v_{n-1} - (z-\alpha)v^{\ast}_n\eklam{R^{\ast}_n(\za)-R^{\ast}_n(\alpha)}H^{-1}_nR_n(\alpha)v_n.
\end{align}
Wegen \fref{drbm2bw11} und der Teile (a) und (d) von \thref{drlm1} gilt
\begin{align}	\label{drbm2bw14}
  & (z-\alpha)^2v^{\ast}_{n-1}R^{\ast}_{n-1}(\za)\Harnm^{-1}y_{0,n-1}v^{\ast}_{n-1}R^{\ast}_{n-1}(\alpha)H^{-1}_{n-1}R_{n-1}(\alpha)v_{n-1} \notag \\
  & = (z-\alpha)^2v^{\ast}_{n-1}R^{\ast}_{n-1}(\za)\Harnm^{-1}R_{n-1}(\alpha)v_{n-1} \notag \\
  &\quad - (z-\alpha)^2v^{\ast}_{n-1}R^{\ast}_{n-1}(\za)T^{\ast}_{n-1}R^{\ast}_{n-1}(\alpha)H^{-1}_{n-1}R_{n-1}(\alpha)v_{n-1} \notag \\
  & = (z-\alpha)^2v^{\ast}_{n-1}R^{\ast}_{n-1}(\za)\Harnm^{-1}R_{n-1}(\alpha)v_{n-1} \notag \\
  &\quad - (z-\alpha)v^{\ast}_{n-1}\eklam{R^{\ast}_{n-1}(\za)-R^{\ast}_{n-1}(\alpha)}H^{-1}_{n-1}R_{n-1}(\alpha)v_{n-1}.
\end{align}
Wegen \thref{drdef1}, \thref{adpdef1}, Teil (a) von \thref{sqlm1}, \fref{drbm2bw13} und \fref{drbm2bw14} gilt
\begin{align*}
  &\ \Carnm(z)-(z-\alpha)\Darn(z)\Marn \\
  & = -(z-\alpha)v^{\ast}_{n-1}R^{\ast}_{n-1}(\za)H^{-1}_{n-1}\Rnm(\alpha)\vnm \\
  &\quad -(z-\alpha)\eklam{\Iq-(z-\alpha)v^{\ast}_{n-1}R^{\ast}_{n-1}(\za)\Harnm^{-1}y_{0,n-1}} \\
  &\quad \cdot \eklam{v^{\ast}_nR^{\ast}_n(\alpha)H^{-1}_nR_n(\alpha)v_n - v^{\ast}_{n-1}R^{\ast}_{n-1}(\alpha)H^{-1}_{n-1}R_{n-1}(\alpha)v_{n-1}} \\
  & = -(z-\alpha)v^{\ast}_{n-1}R^{\ast}_{n-1}(\za)H^{-1}_{n-1}\Rnm(\alpha)\vnm \\
  &\quad - (z-\alpha)v^{\ast}_nR^{\ast}_n(\alpha)H^{-1}_nR_n(\alpha)v_n + (z-\alpha)v^{\ast}_{n-1}R^{\ast}_{n-1}(\alpha)H^{-1}_{n-1}R_{n-1}(\alpha)v_{n-1} \\
  &\quad + (z-\alpha)^2v^{\ast}_{n-1}R^{\ast}_{n-1}(\za)\Harnm^{-1}y_{0,n-1}v^{\ast}_nR^{\ast}_n(\alpha)H^{-1}_nR_n(\alpha)v_n \\
  &\quad - (z-\alpha)^2v^{\ast}_{n-1}R^{\ast}_{n-1}(\za)\Harnm^{-1}y_{0,n-1}v^{\ast}_{n-1}R^{\ast}_{n-1}(\alpha)H^{-1}_{n-1}R_{n-1}(\alpha)v_{n-1} \\
  & = -(z-\alpha)v^{\ast}_{n-1}R^{\ast}_{n-1}(\za)H^{-1}_{n-1}\Rnm(\alpha)\vnm \\
  &\quad - (z-\alpha)v^{\ast}_nR^{\ast}_n(\alpha)H^{-1}_nR_n(\alpha)v_n + (z-\alpha)v^{\ast}_{n-1}R^{\ast}_{n-1}(\alpha)H^{-1}_{n-1}R_{n-1}(\alpha)v_{n-1} \\
  &\quad - (z-\alpha)v^{\ast}_n\eklam{R^{\ast}_n(\za)-R^{\ast}_n(\alpha)}H^{-1}_nR_n(\alpha)v_n \\
  &\quad + (z-\alpha)v^{\ast}_{n-1}\eklam{R^{\ast}_{n-1}(\za)-R^{\ast}_{n-1}(\alpha)}H^{-1}_{n-1}R_{n-1}(\alpha)v_{n-1} \\
  & = -(z-\alpha)\vna\Rna(\za)\Hn^{-1}\Rn(\alpha)\vn = \Carn(z). \tag*{$\Box$}
\end{align*}

Mithilfe von \thref{drbm2} können wir nun die Folge von rechtsseitigen $2$\textit{q}$\times2$\textit{q}-$\alpha$-Dyukarev-Matrixpolynomen bezüglich einer rechtsseitig $\alpha$-Stieltjes-positiv definiten Folge unter Verwendung der rechtsseitigen $\alpha$-Dyukarev-Stieltjes-Parametrisierung jener Folge wie folgt faktorisieren (vergleiche die Teile (a) und (b) mit \cite[Theorem 3.3]{CR1} für den Fall $\alpha=0$ und $\kappa=\infty$).

\begin{satz}	\thlabel{drsa1}
  Seien $\kappa \in \Na$, $\alpha \in \R$, $\sjk \in \Kpqka$ und $\LMkarn$ die rechtsseitige $\alpha$-Dyukarev-Stieltjes-Parametrisierung von $\sjk$.
  Weiterhin sei $\Uarmk$ die Folge von rechtsseitigen $2$\textit{q}$\times2$\textit{q}-$\alpha$-Dyukarev-Matrixpolynomen bezüglich $\sjk$. \dgfa
  \begin{itemize}
    \item [\rm{(a)}] Es gelten
      \begin{align*}
	\Uarznp(z) = \Uarzn(z) \begin{pmatrix} \Iq & \Larn \\ \Oq & \Iq \end{pmatrix}
      \end{align*}
      für alle $z \in \C$ und $n \in \Zofkm$ sowie im Fall $\kappa\geq2$
      \begin{align*}
	\Uarzn(z) = \Uarznm(z) \begin{pmatrix} \Iq & \Oq \\ -(z-\alpha)\Marn & \Iq \end{pmatrix}
      \end{align*}
      für alle $z \in \C$ und $n \in \Zefk$.
    \item [\rm{(b)}] Es gelten
      \begin{align*}
	\Uaro(z) = \begin{pmatrix} \Iq & \Oq \\ -(z-\alpha)\Maro & \Iq \end{pmatrix}
      \end{align*}
      für alle $z \in \C$,
      \begin{align*}
	\Uarznp(z) = \prodr^{n}_{j=0} \begin{pmatrix} \Iq & \Oq \\ -(z-\alpha)\Marj & \Iq \end{pmatrix} \begin{pmatrix} \Iq & \Larj \\ \Oq & \Iq \end{pmatrix}
      \end{align*}
      für alle $z \in \C$ und $n \in \Zofkm$ sowie im Fall $\kappa\geq2$
      \begin{align*}
	\Uarzn(z) = \eklam{\prodr^{n-1}_{j=0} \begin{pmatrix} \Iq & \Oq \\ -(z-\alpha)\Marj & \Iq \end{pmatrix} \begin{pmatrix} \Iq & \Larj \\ \Oq & \Iq \end{pmatrix}} \begin{pmatrix} \Iq & \Oq \\ -(z-\alpha)\Marn & \Iq \end{pmatrix}
      \end{align*}
      für alle $z \in \C$ und $n \in \Zefk$.   
      \item [\rm{(c)}] Seien
      	\begin{align*}
      		W_{\ar2n+1}(z) := \begin{pmatrix} \Iq & \Larn \\ \Oq & \Iq \end{pmatrix} \quad \text{bzw.} \quad
      		W_{\ar2n}(z) := \begin{pmatrix} \Iq & \Oq \\ -(z-\alpha)\Marn & \Iq \end{pmatrix}
      	\end{align*}
      	für alle $z\in\C$ und $n \in \Zofkm$ bzw. $n \in \Zofk$. Dann gelten $W_{\arm}\in\tbtJqPp$ für alle $m\in\Zok$ sowie
      	\begin{align*}
      		U_{\arm}(z) = \prodr^{m}_{j=0} W_{\arj}(z)
      	\end{align*}
      	für alle $z\in\C$ und $m\in\Zok$.
      \item [\rm{(d)}] Mit den Bezeichnungen von (c) gelten $W_{\ar2n+1}\in\bJqLam$ für alle $n \in \Zofkm$, $W_{\ar2n}\in\tbJqLam$ für alle $n \in \Zofk$ und $\Uarm\in\bJqLam$ für alle $m\in\Zok$.
  \end{itemize}
\end{satz}

\bwanf Zu (a): Dies folgt aus \thref{drbm2}.

Zu (b): Sei $z \in \C$. Wegen \thref{drdef2}, \thref{drdef1} und \thref{adpdef1} gilt dann
\begin{align*}
  \Uaro(z) = \begin{pmatrix} \Aaro(z) & \Baro(z) \\ \Caro(z) & \Daro(z) \end{pmatrix} = \begin{pmatrix} \Iq & \Oq \\ -(z-\alpha)s^{-1}_0 & \Iq \end{pmatrix} = \begin{pmatrix} \Iq & \Oq \\ -(z-\alpha)\Maro & \Iq \end{pmatrix}.
\end{align*}
Die weiteren Behauptungen folgen dann iterativ durch Anwendung von (a). 

Zu (c): Offensichtlich ist $W_{\arm}$ für alle $m\in\Zok$ eine in $\C$ holomorphe Funktion. Wegen \thref{adpbm3} gilt $\Larn^{\ast}=\Larn$ für alle $n \in \Zofkm$. Hieraus folgt wegen Teil (c) von \thref{spbm2} dann, dass $W_{\ar2n+1}(z)$ für alle $z\in\C$ und $n \in \Zofkm$ eine $\tJq$-unitäre und somit auch $\tJq$-kontraktive Matrix ist. Hieraus folgt wegen \thref{spbsp1} und \thref{spdef3} nun $W_{\ar2n+1}\in\tbtJqPp$ für alle $n \in \Zofkm$. Wegen \thref{adpbm3} gilt weiterhin $\Marn\in\Cqq_{\geq}$ für alle $n \in \Zofk$. Hieraus folgt wegen Teil (a) von \thref{spbm2} dann
\begin{align*}
	W^{\ast}_{\ar2n}(z)\tJq W_{\ar2n}(z) = \tJq + \diag(i(z-\za)\Marn,\Oq)
\end{align*}
für alle $z\in\C$ und $n \in \Zofk$. Hieraus folgt für alle $n \in \Zofk$ dann, dass $W_{\ar2n}(z)$ für alle $z\in\Pp$ eine $\tJq$-kontraktive sowie $W_{\ar2n}(x)$ für alle $x\in\R$ eine $\tJq$-unitäre Matrix ist. Hieraus folgt wegen \thref{spbsp1} und \thref{spdef3} nun $W_{\ar2n}\in\tbtJqPp$ für alle $n \in \Zofk$. Wegen (b) gilt weiterhin
\begin{align*}
	U_{\arm}(z) = \prodr^{m}_{j=0} W_{\arj}(z)
\end{align*}
für alle $z\in\C$ und $m\in\Zok$. 

Zu (d): Wegen (c) ist $W_{\arm}$ für alle $m\in\Zok$ eine in $\C$ holomorphe Funktion. Wegen \thref{adpbm3} gilt $\Larn\in\Cqq_{\geq}$ für alle $n \in \Zofkm$. Hieraus folgt wegen Teil (c) von \thref{spbm5} dann, dass $W_{\ar2n+1}(z)$ für alle $z\in\C$ und $n \in \Zofkm$ eine $\Jq$-kontraktive Matrix ist. Hieraus folgt wegen \thref{spbsp4} und Teil (a) von \thref{spdef8} nun $W_{\ar2n+1}\in\bJqLam$ für alle $n \in \Zofkm$. 
Wegen \thref{adpbm3} gilt weiterhin $\Marn\in\Cqq_{\geq}$ für alle $n \in \Zofk$. Hieraus folgt wegen Teil (a) von \thref{spbm5} dann 
\begin{align*}
	\Jq-W^{\ast}_{\ar2n}(z)\Jq W_{\ar2n}(z) &= \diag(-(z-\alpha)\Marn-(\za-\alpha)\Marn^{\ast},\Oq) \\
	&= 2(\alpha-\re z)\diag(\Marn,\Oq)
\end{align*}
für alle $z\in\C$ und $n \in \Zofk$. Hieraus folgt für alle $n \in \Zofk$ dann, dass $W_{\ar2n}(z)$ für alle $z\in\Lam$ eine $\Jq$-kontraktive Matrix sowie für alle $z\in\C$ mit $\re z=\alpha$ eine $\Jq$-unitäre Matrix ist. Hieraus folgt wegen \thref{spbsp4} und \thref{spdef8} nun $W_{\ar2n}\in\tbJqLam$ für alle $n \in \Zofk$.
Wegen (c) gilt
\begin{align*}
	U_{\arm} = \prodr^{m}_{j=0} W_{\arj}
\end{align*}
für alle $m\in\Zok$. Hieraus folgt für alle $m\in\Zok$ wegen $W_{\arj}\in\bJqLam$ für alle $j\in\Zom$, \thref{spbsp4} und \thref{spbm8} dann $\Uarm\in\bJqLam$. \bwend

Da unter Beachtung von \thref{spbsp1} und \thref{spbm6} jedes Produkt von Funktionen aus $\tbtJqPp$ selbst eine Funktion aus $\tbtJqPp$ ist, liefert Teil (c) von \thref{drsa1} sogleich einen alternativen Beweis zu Teil (c) von \thref{drbm5}.

Abschließend schauen wir uns noch die Leitkoeffizienten und die Absolutglieder der einzelnen Matrixpolynome des rechtsseitigen $\alpha$-Dyukarev-Quadrupels bezüglich einer $\alpha$-Stieltjes-positiv definiten Folge an (vergleiche \cite[Proposition 3.5]{CR1} für den Fall $\alpha=0$ und $\kappa=\infty$).

\begin{satz}	\thlabel{drsa7}
	Seien $\kappa\in\Na$, $\alpha\in\R$, $\sjk\in\Kpqka$ und $\ABCDark$ das rechtsseitige $\alpha$-Dyukarev-Quadrupel bezüglich $\sjk$. Weiterhin sei $\LMkarn$ die rechtsseitige $\alpha$-Dyukarev-Stieltjes-Parametrisierung von \linebreak $\sjk$. Dann gelten
	\begin{align}
		\Barn(z)&=(-1)^{n-1}(z-\alpha)^{n-1}\Laro\prodr^{n-1}_{j=1}\rklam{\Marj\Larj} + \ldots + \sum^{n-1}_{j=0}\Larj, \label{drsa7bw1}\\
		\Darn(z)&=(-1)^n(z-\alpha)^n\prodr^{n-1}_{j=0}\rklam{\Marj\Larj} + \ldots + \Iq \label{drsa7bw2}
	\end{align}
	für alle $z\in\C$ und $n\in\Zefkp$ und im Fall $\kappa\geq2$
	\begin{align}
		\Aarn(z)&=(-1)^n(z-\alpha)^n\prodr^{n-1}_{j=0}\rklam{\Larj\Marjp} + \ldots + \Iq, \label{drsa7bw3}\\
		\Carn(z)&=(-1)^{n+1}(z-\alpha)^{n+1}\Maro\prodr^{n-1}_{j=0}\rklam{\Larj\Marjp} + \ldots - (z-\alpha)\sum^n_{j=0}\Marj \label{drsa7bw4}
	\end{align}
	für alle $z\in\C$ und $n\in\Zefk$.
\end{satz}

\bwanf Wir zeigen die Behauptungen mithilfe vollständiger Induktion. Sei $z\in\C$. Wegen \thref{drdef1} gelten dann
\begin{align*}
	\Baro(z) = \Oq, \quad \Daro(z) = \Iq, \quad \Aaro(z) = \Iq \quad \text{und} \quad \Caro(z) = -(z-\alpha)s^{-1}_0.
\end{align*}
Hieraus folgt wegen \thref{drbm2} und \thref{adpdef1} nun
\begin{align*}
	\Bare(z) &= \Aaro(z)\Laro+\Baro(z) = \Laro \quad \text{und} \\
	\Dare(z) &= \Caro(z)\Laro+\Daro(z) = -(z-\alpha)\Maro\Laro + \Iq
\end{align*}
sowie im Fall $\kappa\geq2$
\begin{align*}
	\Aare(z) &= \Aaro(z)-(z-\alpha)\Bare(z)\Mare = -(z-\alpha)\Laro\Mare + \Iq \quad \text{und} \\
	\Care(z) &= \Caro(z)-(z-\alpha)\Dare(z)\Mare \\
	&= (z-\alpha)^2\Maro\Laro\Mare-(z-\alpha)\rklam{\Maro+\Mare}.
\end{align*}

Seien nun $\kappa\geq3$, $m\in\Zzfkp$ und \fref{drsa7bw1}-\fref{drsa7bw4} für alle $n\in\Zemm$ erfüllt. Wegen \thref{drbm2} gelten dann
\begin{align*}
	&\ \mathbf{B}_{\arm}(z) = \mathbf{A}_{\arm-1}(z)\mathbf{L}_{\arm-1}+\mathbf{B}_{\arm-1}(z) \\
	&= \eklam{(-1)^{m-1}(z-\alpha)^{m-1}\prodr^{m-2}_{j=0}\rklam{\Larj\Marjp} + \ldots + \Iq}\mathbf{L}_{\arm-1} \\
	&\quad + (-1)^{m-2}(z-\alpha)^{m-2}\Laro\prodr^{m-2}_{j=1}\rklam{\Marj\Larj} + \ldots + \sum^{m-2}_{j=0}\Larj \\
	&= (-1)^{m-1}(z-\alpha)^{m-1}\Laro\prodr^{m-1}_{j=1}\rklam{\Marj\Larj} + \ldots + \sum^{m-1}_{j=0}\Larj
\end{align*}
und
\begin{align*}
	&\ \mathbf{D}_{\arm}(z) = \mathbf{C}_{\arm-1}(z)\mathbf{L}_{\arm-1}+\mathbf{D}_{\arm-1}(z) \\
	&= \eklam{(-1)^m(z-\alpha)^m\Maro\prodr^{m-2}_{j=0}\rklam{\Larj\Marjp} + \ldots - (z-\alpha)\sum^{m-1}_{j=0}\Marj}\mathbf{L}_{\arm-1} \\
	&\quad + (-1)^{m-1}(z-\alpha)^{m-1}\prodr^{m-2}_{j=1}\rklam{\Marj\Larj} + \ldots + \Iq \\
	&= (-1)^m(z-\alpha)^m\prodr^{m-1}_{j=1}\rklam{\Marj\Larj} + \ldots + \Iq.
\end{align*}

Seien nun $\kappa\geq4$, $m\in\Zzfk$ und \fref{drsa7bw1}-\fref{drsa7bw2} für alle $n\in\Zem$ sowie \fref{drsa7bw3}-\fref{drsa7bw4} für alle $n\in\Zemm$ erfüllt. Wegen \thref{drbm2} gelten dann
\begin{align*}
	&\ \mathbf{A}_{\arm} = \mathbf{A}_{\arm-1}-(z-\alpha)\mathbf{B}_{\arm}\mathbf{M}_{\arm} \\
	&= (-1)^{m-1}(z-\alpha)^{m-1}\prodr^{m-2}_{j=0}\rklam{\Larj\Marjp} + \ldots + \Iq \\
	&\quad -(z-\alpha)\eklam{(-1)^{m-1}(z-\alpha)^{m-1}\Laro\prodr^{m-1}_{j=1}\rklam{\Marj\Larj} + \ldots + \sum^{m-1}_{j=0}\Larj}\mathbf{M}_{\arm} \\
	&= (-1)^m(z-\alpha)^m\prodr^{m-1}_{j=0}\rklam{\Larj\Marjp} + \ldots + \Iq
\end{align*}
und
\begin{align*}
	&\ \mathbf{C}_{\arm} = \mathbf{C}_{\arm-1}-(z-\alpha)\mathbf{D}_{\arm}\mathbf{M}_{\arm} \\
	&= (-1)^m(z-\alpha)^m\Maro\prodr^{m-2}_{j=0}\rklam{\Larj\Marjp} + \ldots -(z-\alpha)\sum^{m-1}_{j=0}\Marj \\
	&\quad -(z-\alpha)\eklam{(-1)^m(z-\alpha)^m\prodr^{m-1}_{j=1}\rklam{\Marj\Larj} + \ldots + \Iq}\mathbf{M}_{\arm} \\
	&= (-1)^{m+1}(z-\alpha)^{m+1}\Maro\prodr^{m-1}_{j=0}\rklam{\Larj\Marjp} + \ldots -(z-\alpha)\sum^m_{j=0}\Marj. \tag*{$\Box$}
\end{align*}

\subsection{Der linksseitige Fall}

Wir behandeln nun eine Faktorisierung der Folge von linksseitigen $2$\textit{q}$\times2$\textit{q}-$\alpha$-Dyukarev"=Matrixpolynomen einer linksseitig $\alpha$-Stieltjes-positiv definiten Folge unter Verwendung der linksseitigen $\alpha$-Dyukarev-Stieltjes-Parametrisierung jener Folge.

\begin{satz}	\thlabel{drlsa5}
 	Seien $\kappa \in \Na$, $\alpha \in \R$, $\sjk \in \Lpqka$ und $\LMkaln$ die linksseitige $\alpha$-Dyukarev-Stieltjes-Parametrisierung von $\sjk$.
  	Weiterhin sei $\Ualmk$ die Folge von linksseitigen $2$\textit{q}$\times2$\textit{q}-$\alpha$-Dyukarev-Matrixpolynomen bezüglich $\sjk$. \dgfa
  	\begin{itemize}
    	\item [\rm{(a)}] Es gelten
      	\begin{align*}
			\Ualznp(z) = \Ualzn(z) \begin{pmatrix} \Iq & -\Laln \\ \Oq & \Iq \end{pmatrix}
      	\end{align*}
     	für alle $z \in \C$ und $n \in \Zofkm$ sowie im Fall $\kappa\geq2$
      	\begin{align*}
			\Ualzn(z) = \Ualznm(z) \begin{pmatrix} \Iq & \Oq \\ (\alpha-z)\Maln & \Iq \end{pmatrix}
      	\end{align*}
      	für alle $z \in \C$ und $n \in \Zefk$.
    	\item [\rm{(b)}] Es gelten
      	\begin{align*}
			\Ualo(z) = \begin{pmatrix} \Iq & \Oq \\ (\alpha-z)\Malo & \Iq \end{pmatrix}
      	\end{align*}
      	für alle $z \in \C$,
      	\begin{align*}
			\Ualznp(z) = \prodr^{n}_{j=0} \begin{pmatrix} \Iq & \Oq \\ (\alpha-z)\Malj & \Iq \end{pmatrix} \begin{pmatrix} \Iq & -\Lalj \\ \Oq & \Iq \end{pmatrix}
      	\end{align*}
      	für alle $z \in \C$ und $n \in \Zofkm$ sowie im Fall $\kappa\geq2$
      	\begin{align*}
			\Ualzn(z) = \eklam{\prodr^{n-1}_{j=0} \begin{pmatrix} \Iq & \Oq \\ (\alpha-z)\Malj & \Iq \end{pmatrix} \begin{pmatrix} \Iq & -\Lalj \\ \Oq & \Iq \end{pmatrix}} \begin{pmatrix} \Iq & \Oq \\ (\alpha-z)\Maln & \Iq \end{pmatrix}
      	\end{align*}
      	für alle $z \in \C$ und $n \in \Zefk$.
      	\item [\rm{(c)}] Seien
      	\begin{align*}
      		W_{\al2n+1}(z) := \begin{pmatrix} \Iq & -\Laln \\ \Oq & \Iq \end{pmatrix} \quad \text{bzw.} \quad
      		W_{\al2n}(z) := \begin{pmatrix} \Iq & \Oq \\ (\alpha-z)\Maln & \Iq \end{pmatrix}
      	\end{align*}
      	für alle $z\in\C$ und $n \in \Zofkm$ bzw. $n \in \Zofk$. Dann gelten $W_{\alm}\in\tbtJqPp$ für alle $m\in\Zok$ sowie
      	\begin{align*}
      		U_{\alm}(z) = \prodr^{m}_{j=0} W_{\alj}(z)
      	\end{align*}
      	für alle $z\in\C$ und $m\in\Zok$.   
  	\end{itemize}
\end{satz}

\bwanf Zu (a): Seien $t_j:=(-1)^js_j$ für alle $j\in\Zok$. Wegen Teil (a) von \thref{asmbm3} gilt dann $\tjk\in\Kpqkma$. Weiterhin seien $[(\mathbf{L}^{\sklam{t}}_{-\arn})^{\fklam{\kappa-1}}_{n=0}$,""$(\mathbf{M}^{\sklam{t}}_{-\arn})^{\fklam{\kappa}}_{n=0}]$ die rechtsseitige $-\alpha$-Dyukarev-Stieltjes-Parametrisierung von $\tjk$ und $(U^{\sklam{t}}_{-\arm})^{\kappa}_{m=0}$ die Folge von rechtsseitigen $2$\textit{q}$\times2$\textit{q}-$-\alpha$-Dyukarev-Matrixpolynomen bezüglich $\tjk$. Wegen Teil (a) von \thref{drllm4}, Teil (a) von \thref{drsa1}, Teil (a) von \thref{asmlm1} und Teil (c) von \thref{adpllm1} gelten dann
\begin{align*}
	\Ualznp^{\sklam{s}}(z) &= \Ve U^{\sklam{t}}_{-\ar2n+1}(-z) \Vea \\
	&= \Ve U^{\sklam{t}}_{-\ar2n}(-z) \begin{pmatrix} \Iq & \mathbf{L}^{\sklam{t}}_{-\arn} \\ \Oq & \Iq \end{pmatrix} \Vea \\
	&= \Ve U^{\sklam{t}}_{-\ar2n}(-z) \Vea\Ve \begin{pmatrix} \Iq & \Laln^{\sklam{s}} \\ \Oq & \Iq \end{pmatrix} \Vea \\
	&= \Ualzn^{\sklam{s}}(z) \begin{pmatrix} \Iq & -\Laln^{\sklam{s}} \\ \Oq & \Iq \end{pmatrix}
\end{align*}
für alle $z \in \C$ und $n \in \Zofkm$ sowie im Fall $\kappa\geq2$
\begin{align*}
	\Ualzn^{\sklam{s}}(z) &= \Ve U^{\sklam{t}}_{-\ar2n}(-z) \Vea \\
	&= \Ve U^{\sklam{t}}_{-\ar2n-1}(-z) \begin{pmatrix} \Iq & \Oq \\ -(-z-(-\alpha))\mathbf{M}^{\sklam{t}}_{-\arn} & \Iq \end{pmatrix} \Vea \\
	&= \Ve U^{\sklam{t}}_{-\ar2n-1}(-z) \Vea\Ve \begin{pmatrix} \Iq & \Oq \\ -(\alpha-z)\Maln^{\sklam{s}} & \Iq \end{pmatrix} \Vea \\
	&= \Ualznm^{\sklam{s}}(z) \begin{pmatrix} \Iq & \Oq \\ (\alpha-z)\Maln^{\sklam{s}} & \Iq \end{pmatrix}
\end{align*}
für alle $z \in \C$ und $n \in \Zefk$.

Zu (b): Sei $z \in \C$. Wegen \thref{drldef3}, \thref{drldef2} und \thref{adpldef1} gilt dann
\begin{align*}
  \Ualo(z) = \begin{pmatrix} \Aalo(z) & \Balo(z) \\ \Calo(z) & \Dalo(z) \end{pmatrix} = \begin{pmatrix} \Iq & \Oq \\ (\alpha-z)s^{-1}_0 & \Iq \end{pmatrix} = \begin{pmatrix} \Iq & \Oq \\ (\alpha-z)\Malo & \Iq \end{pmatrix}.
\end{align*}
Die weiteren Behauptungen folgen dann iterativ durch Anwendung von (a). 

Zu (c): Offensichtlich ist $W_{\alm}$ für alle $m\in\Zok$ eine in $\C$ holomorphe Funktion. Wegen \thref{adplbm2} gilt $\Laln^{\ast}=\Laln$ für alle $n \in \Zofkm$. Hieraus folgt wegen Teil (c) von \thref{spbm2} dann, dass $W_{\al2n+1}(z)$ für alle $z\in\C$ und $n \in \Zofkm$ eine $\tJq$-unitäre und somit auch $\tJq$-kontraktive Matrix ist. Hieraus folgt wegen \thref{spbsp1} und \thref{spdef3} nun $W_{\al2n+1}\in\tbtJqPp$ für alle $n \in \Zofkm$. Wegen \thref{adpbm3} gilt weiterhin $\Maln^{\ast}=\Maln$ für alle $n \in \Zofk$. Hieraus folgt wegen Teil (a) von \thref{spbm2} dann
\begin{align*}
	W^{\ast}_{\al2n}(z)\tJq W_{\al2n}(z) = \tJq + \diag(i(z-\za)\Maln,\Oq)
\end{align*}
für alle $z\in\C$ und $n \in \Zofk$. Hieraus folgt für alle $n \in \Zofk$ dann, dass $W_{\al2n}(z)$ für alle $z\in\Pp$ eine $\tJq$-kontraktive sowie $W_{\al2n}(x)$ für alle $x\in\R$ eine $\tJq$-unitäre Matrix ist. Hieraus folgt wegen \thref{spbsp1} und \thref{spdef3} nun $W_{\al2n}\in\tbtJqPp$ für alle $n \in \Zofk$. Wegen (b) gilt weiterhin
\begin{align*}
	U_{\alm}(z) = \prodr^{m}_{j=0} W_{\alj}(z)
\end{align*}
für alle $z\in\C$ und $m\in\Zok$. \bwend

Da unter Beachtung von \thref{spbsp1} und \thref{spbm6} jedes Produkt von Funktionen aus $\tbtJqPp$ selbst eine Funktion aus $\tbtJqPp$ ist, liefert Teil (c) von \thref{drlsa5} sogleich einen alternativen Beweis zu Teil (c) von \thref{drlbm2}.

%% file: sm4.tex
\newpage
\section{Eine alternative Beschreibung der Lösungsmenge für vollständig nichtdegenerierte matrizielle \texorpdfstring{$\alpha$}{a}"=Stieltjes Momentenprobleme} \label{chapar}

Im Mittelpunkt dieses Kapitels steht die Herleitung einer im Vergleich zu Kapitel \ref{chapdr} alternativen Beschreibung der Lösungsmenge für mithilfe der $[\alpha,\infty)$-Stieltjes"=Transformation umformulierte vollständig nichtdegenerierte matrizielle $\alpha$-Stieltjes Momentenprobleme. Statt einer speziellen Klasse von Stieltjes-Paaren als Parametermenge verwenden wir nun allerdings eine Teilklasse der Schurfunktionen. Auch jetzt erfolgt die Beschreibung mithilfe einer gebrochen linearen Transformation von Matrizen. 

Der rechtsseitige Fall wurde für eine gegebene Folge aus ${\cal K}^{\geq,e}_{q,2n+1,\alpha}$ in \cite{Wb} behandelt und wir werden die dortige Herangehensweise auf unser Problem übertragen können. Im linksseitigen Fall greifen wir dann auf die Resultate des rechtsseitigen Falles zurück.

\subsection{Der rechtsseitige Fall}

Am Ausgangspunkt dieses Abschnitts steht die in \thref{drth1} erhaltene Beschreibung der Lösungsmenge des zu einer endlichen rechtsseitig $\alpha$-Stieltjes-positiv definiten Folge zugehörigen mithilfe der $[\alpha,\infty)$-Stieltjes"=Transformation umformulierten rechtsseitigen $\alpha$-Stieltjes Momentenproblems. Diese Beschreibung erfolgte dort durch Parametrisierung mithilfe einer erzeugenden Matrixfunktion, welche ein $2$\textit{q}$\times2$\textit{q}-Matrixpolynom der Klasse $\tbtJqPp$ ist. Als Parametermenge tritt hierbei eine spezielle Klasse von geordneten Paaren von in $\C\setminus[\alpha,\infty)$ meromorphen \textit{q}$\times$\textit{q}-Matrixfunktionen auf. 

Aufgrund der in Anhang \ref{chapsm} entwickelten Theorie ausgewählter Resultate für gebrochen lineare Transformationen von Matrizen können wir nun die gebrochen lineare Transformation etwas einfacher gestalten, indem wie die erzeugende Matrixfunktion durch ein Matrixpolynom ersetzen, welches entsprechende Kontraktivitätseigenschaften bezüglich der Signaturmatrix
$\jqq := \diag(\Iq,-\Iq)$
und $\tJq$ besitzt. Der große Vorteil ist dann, dass die Menge der Parameterfunktionen der gebrochen linearen Transformation eine Teilklasse der Menge der \textit{q}$\times$\textit{q}-Schurfunktionen auf $\Pp$ (vergleiche \thref{spdef6}) ist.

Diese Beschreibung der Lösungsmenge des rechtsseitigen $\alpha$-Stieltjes Momentenproblems wurde für eine gegebene Folge aus ${\cal K}^{\geq,e}_{q,2n+1,\alpha}$ in \cite[Kapitel 9]{Wb} behandelt und wir verwenden für unseren Fall eine ähnliche Vorgehensweise. Bevor wir dies tun, wollen wir uns zunächst oben erwähnten Kontraktivitätseigenschaften der erzeugenden Matrixfunktion widmen. 

\begin{satz}	\thlabel{drsa10}
	Seien $\kappa\in\Na$, $\alpha\in\R$, $\sjk\in\Kpqka$ und $\Uarmk$ die Folge von rechtsseitigen $2$\textit{q}$\times2$\textit{q}-$\alpha$-Dyukarev-Matrixpolynomen bezüglich $\sjk$. Weiterhin seien
	\begin{align*}
  		E := \frac{1}{\sqrt{2}}\begin{pmatrix} -i\Iq & i\Iq \\ \Iq & \Iq \end{pmatrix} \quad \text{und} \quad \Sigma_{\arm} := \Uarm E
  	\end{align*}
  	für alle $m\in\Zok$. Dann gilt $\Sarm\in\tbtJJPp$ für alle $m\in\Zok$.
\end{satz}

\bwanf Dies folgt wegen Teil (c) von \thref{drbm5} sogleich aus Teil (a2) von \thref{spbm9}. \bwend

Wir kommen nun zur Beschreibung der Lösungsmenge des zu einer endlichen rechtsseitig $\alpha$-Stieltjes-positiv definiten Folge zugehörigen mithilfe der $[\alpha,\infty)$-Stieltjes"=Transformation umformulierten rechtsseitigen $\alpha$-Stieltjes Momentenproblems mit $\Dqa$ (vergleiche Teil (b) von \thref{spdef6}) als Parametermenge.

\begin{theo}	\thlabel{drth3}
  Seien $m \in \N$, $\alpha \in \R$, $\sjm \in \Kpqma$, $\Uarm$ das rechtsseitige $2$\textit{q}$\times2$\textit{q}-$\alpha$-Dyukarev-Matrixpolynom bezüglich $\sjm$,
  \begin{align*}
  	E := \frac{1}{\sqrt{2}}\begin{pmatrix} -i\Iq & i\Iq \\ \Iq & \Iq \end{pmatrix} \quad \text{und} \quad \Sigma_{\arm} := \Uarm E.
  \end{align*}
  Weiterhin bezeichne
  \begin{align*}
    \Sarm = \begin{pmatrix} \Sarm^{(1,1)} & \Sarm^{(1,2)} \\ \Sarm^{(2,1)} & \Sarm^{(2,2)} \end{pmatrix}
  \end{align*}
  die \textit{q}$\times$\textit{q}-Blockzerlegung von $\Sarm$. \dgfa
  \begin{itemize}
    \item [\rm{(a)}]  Sei $F\in\Dqa$. Weiterhin seien
	\begin{align*}
		\MF := \gklam{z\in\Pm \;|\; \det F(\za) = 0}
	\end{align*}
	und $\tF: \Pp\cup(-\infty,\alpha)\cup(\Pm\setminus\MF) \rightarrow \Cqq$ definiert gemäß
	\begin{align*}
		\tF(z) := \begin{cases} F(z) & \text{falls } z\in\Pp \\
		\lim_{\omega\rightarrow z}F(\omega) & \text{falls } z\in(-\infty,\alpha) \\
		F^{-\ast}(\za) & \text{falls } z\in\Pm\setminus\MF \end{cases}.
	\end{align*}
	Dann ist $\det\big(\Sarm^{(2,1)}\tF+\Sarm^{(2,2)}\big)$ nicht die Nullfunktion und
	\begin{align*}
      	S := \brklam{\Sarm^{(1,1)}\tF+\Sarm^{(1,2)}}\brklam{\Sarm^{(2,1)}\tF+\Sarm^{(2,2)}}^{-1}
    \end{align*}
    gehört zu $\Sqasmu$.
    \item [\rm{(b)}] Sei $S\in\Sqasmu$. Dann existiert genau eine in $\C\setminus[\alpha,\infty)$ meromorphe \textit{q}$\times$\textit{q}-Matrixfunktion $\tF$ mit $\Rstr_{\Pp}\tF\in\Dqa$ derart, dass
  	$\det\big(\Sarm^{(2,1)}\tF+\Sarm^{(2,2)}\big)$ nicht die Nullfunktion ist und
    \begin{align*}
      	S = \brklam{\Sarm^{(1,1)}\tF+\Sarm^{(1,2)}}\brklam{\Sarm^{(2,1)}\tF+\Sarm^{(2,2)}}^{-1}
    \end{align*}
    erfüllt ist.
  \end{itemize}
\end{theo}

\bwanf Es bezeichne
\begin{align*}
	\Uarm = \begin{pmatrix} \Uarm^{(1,1)} & \Uarm^{(1,2)} \\ \Uarm^{(2,1)} & \Uarm^{(2,2)} \end{pmatrix}
\end{align*}
die \textit{q}$\times$\textit{q}-Blockzerlegung von $\Uarm$. Dann gilt
\begin{align}	\label{drth3bw1}
	\Sarm = \frac{1}{\sqrt{2}} \begin{pmatrix} -i\Uarm^{(1,1)}+\Uarm^{(1,2)} & i\Uarm^{(1,1)}+\Uarm^{(1,2)} \\ -i\Uarm^{(2,1)}+\Uarm^{(2,2)} & i\Uarm^{(2,1)}+\Uarm^{(2,2)} \end{pmatrix}.
\end{align}

Zu (a): Seien $\phi := i\big(\Iq-\tF\big)$ und $\psi := \Iq+\tF$. Wegen Teil (b) von \thref{spsa1} gilt dann $\phipsi\in\PtJqCa$. Wegen \fref{drth3bw1} gelten weiterhin
\begin{align}	\label{drth3bw2}
	\Uarm^{(1,1)}\phi+\Uarm^{(1,2)}\psi 
	&= \rklam{-i\Uarm^{(1,1)}+\Uarm^{(1,2)}}\tF+i\Uarm^{(1,1)}+\Uarm^{(1,2)} \notag \\
	&= \sqrt{2}\brklam{\Sarm^{(1,1)}\tF+\Sarm^{(1,2)}}
\end{align}
und
\begin{align}	\label{drth3bw3}
	\Uarm^{(2,1)}\phi+\Uarm^{(2,2)}\psi 
	&= \rklam{-i\Uarm^{(2,1)}+\Uarm^{(2,2)}}\tF+i\Uarm^{(2,1)}+\Uarm^{(2,2)} \notag \\
	&= \sqrt{2}\brklam{\Sarm^{(2,1)}\tF+\Sarm^{(2,2)}}.
\end{align}
Wegen Teil (a) von \thref{drth1} ist $\det\big(\Uarm^{(2,1)}\phi+\Uarm^{(2,2)}\psi\big)$ nicht die Nullfunktion, also ist wegen \fref{drth3bw3} auch $\det\big(\Sarm^{(2,1)}\tF+\Sarm^{(2,2)}\big)$ nicht die Nullfunktion. Wegen \fref{drth3bw2} und \fref{drth3bw3} gilt nun
\begin{align*}
	S &= \brklam{\Sarm^{(1,1)}\tF+\Sarm^{(1,2)}}\brklam{\Sarm^{(2,1)}\tF+\Sarm^{(2,2)}}^{-1} \\
	&= \rklam{\Uarm^{(1,1)}\phi+\Uarm^{(1,2)}\psi}\rklam{\Uarm^{(2,1)}\phi+\Uarm^{(2,2)}\psi}^{-1}.
\end{align*}
Hieraus folgt wegen Teil (a) von \thref{drth1} dann $S\in\Sqasmu$.

Zu (b): Wegen Teil (b) von \thref{drth1} existiert ein $\phipsi \in \dPtJqCa$ derart, dass
\begin{align}	\label{drth3bw4}
	\det\eklam{\Uarm^{(2,1)}(z)\phi(z)+\Uarm^{(2,2)}(z)\psi(z)} \neq 0
\end{align}
und
\begin{align}	\label{drth3bw5}
	S(z) = \eklam{\Uarm^{(1,1)}(z)\phi(z)+\Uarm^{(1,2)}(z)\psi(z)}\eklam{\Uarm^{(2,1)}(z)\phi(z)+\Uarm^{(2,2)}(z)\psi(z)}^{-1}
\end{align}
für alle $z \in \C\setminus[\alpha,\infty)$ erfüllt sind. 
Wegen Teil (a) von \thref{spsa1} ist $\det(\psi-i\phi)$ nicht die Nullfunktion. Sei $\tF := (\psi+i\phi)(\psi-i\phi)^{-1}$. Wegen Teil (a) von \thref{spsa1} gelten dann weiterhin $\Rstr_{\Pp}\tF\in\Dqa$ sowie
\begin{align}	\label{drth3bw9}
	\phi = \frac{i}{2}\brklam{\Iq-\tF}\rklam{\psi-i\phi} \quad \text{und} \quad
	\psi = \frac{1}{2}\brklam{\Iq+\tF}\rklam{\psi-i\phi}.
\end{align}
Seien $\tphi := i\big(\Iq-\tF\big)$, $\tpsi := \Iq+\tF$ und $g := \frac{1}{2}(\psi-i\phi)$.
Wegen \fref{drth3bw9} gelten dann
\begin{align}	\label{drth3bw6}
	\phi = \tphi g \quad \text{und} \quad \psi = \tpsi g.
\end{align}
Wegen \fref{drth3bw1} gelten weiterhin
\begin{align}	\label{drth3bw7}
	\Uarm^{(1,1)}\tphi+\Uarm^{(1,2)}\tpsi 
	&= \rklam{-i\Uarm^{(1,1)}+\Uarm^{(1,2)}}\tF+i\Uarm^{(1,1)}+\Uarm^{(1,2)} \notag \\
	&= \sqrt{2}\brklam{\Sarm^{(1,1)}\tF+\Sarm^{(1,2)}}
\end{align}
und
\begin{align}	\label{drth3bw8}
	\Uarm^{(2,1)}\tphi+\Uarm^{(2,2)}\tpsi 
	&= \rklam{-i\Uarm^{(2,1)}+\Uarm^{(2,2)}}\tF+i\Uarm^{(2,1)}+\Uarm^{(2,2)} \notag \\
	&= \sqrt{2}\brklam{\Sarm^{(2,1)}\tF+\Sarm^{(2,2)}}.
\end{align}
Wegen \fref{drth3bw6} und \fref{drth3bw8} gilt nun
\begin{align*}
	\det\rklam{\Uarm^{(2,1)}\phi+\Uarm^{(2,2)}\psi}
	&= \det\brklam{\Uarm^{(2,1)}\tphi+\Uarm^{(2,2)}\tpsi}\det g \\
	&= 2^{\frac{q}{2}}\det\brklam{\Sarm^{(2,1)}\tF+\Sarm^{(2,2)}}\det g.
\end{align*}
Unter Beachtung, dass $\det(\psi-i\phi)$ und somit $\det g$ nicht die Nullfunktion ist, folgt hieraus wegen \fref{drth3bw4} dann, dass $\det\big(\Sarm^{(2,1)}\tF+\Sarm^{(2,2)}\big)$ nicht die Nullfunktion ist.
Unter Beachtung, dass $\det g$ nicht die Nullfunktion ist, folgt hieraus wegen \fref{drth3bw5}, \fref{drth3bw6}, \fref{drth3bw7} und \fref{drth3bw8} nun
\begin{align*}
	S &= \rklam{\Uarm^{(1,1)}\phi+\Uarm^{(1,2)}\psi}\rklam{\Uarm^{(2,1)}\phi+\Uarm^{(2,2)}\psi}^{-1} \\
	&= \brklam{\Uarm^{(1,1)}\tphi+\Uarm^{(1,2)}\tpsi}gg^{-1}\brklam{\Uarm^{(2,1)}\tphi+\Uarm^{(2,2)}\tpsi}^{-1} \\
	&= \brklam{\Sarm^{(1,1)}\tF+\Sarm^{(1,2)}}\brklam{\Sarm^{(2,1)}\tF+\Sarm^{(2,2)}}^{-1}.
\end{align*}

Bleibt noch die Eindeutigkeitsaussage zu zeigen. Seien hierzu $\tF_1$ und $\tF_2$ in $\C\setminus[\alpha,\infty)$ meromorphe \textit{q}$\times$\textit{q}-Matrixfunktionen mit
\begin{itemize}
	\item [\rm{(i)}] Es gilt $F_k := \Rstr_{\Pp}\tF_k \in \Dqa$.
	\item [\rm{(ii)}] $\det\big(\Sarm^{(2,1)}\tF_k+\Sarm^{(2,2)}\big)$ ist nicht die Nullfunktion.
	\item [\rm{(iii)}] Es gilt
	\begin{align*}
		S = \brklam{\Sarm^{(1,1)}\tF_k+\Sarm^{(1,2)}}\brklam{\Sarm^{(2,1)}\tF_k+\Sarm^{(2,2)}}^{-1}.
	\end{align*}
\end{itemize}
für alle $k\in\gklam{1,2}$. Sei nun $k\in\gklam{1,2}$.
Weiterhin seien ${\cal M}_{F_k} := \gklam{z\in\Pp \;|\; \det F_k(\za)=0}$ und 
$\widehat{F}_k: \Pp\cup(-\infty,\alpha)\cup(\Pm\setminus{\cal M}_{F_k}) \rightarrow \Cqq$ definiert gemäß
\begin{align*}
	\widehat{F}_k(z) := \begin{cases} F_k(z) & \text{falls } z\in\Pp \\
	\lim_{\omega\rightarrow z}F_k(\omega) & \text{falls } z\in(-\infty,\alpha) \\
	F_k^{-\ast}(\za) & \text{falls } z\in\Pm\setminus{\cal M}_{F_k} \end{cases}.
\end{align*}
Wegen (i) und Teil (b) von \thref{spsa1} ist dann $\widehat{F}_k$ eine in $\C\setminus[\alpha,\infty)$ meromorphe \textit{q}$\times$\textit{q}-Matrixfunktion und es gilt
\begin{align*}
	\Rstr_{\Pp}\widehat{F}_k = F_k = \Rstr_{\Pp}\tF_k.
\end{align*}
Hieraus folgt wegen des Identitätssatzes für meromorphe Funktionen (vergleiche z.\,B. im skalaren Fall \cite[Satz 10.3.2]{Funk}; im matriziellen Fall betrachtet man die einzelnen Einträge der Matrixfunktion) dann 
\begin{align}	\label{drth3bw14}
	\widehat{F}_k = \tF_k.
\end{align}
Seien $\phi_k := i\big(\Iq-\tF_k\big)$ sowie $\psi_k := \Iq+\tF_k$. Wegen (i), \fref{drth3bw14} und Teil (b) von \thref{spsa1} gilt dann weiterhin $\binom{\phi_k}{\psi_k}\in\PtJqCa$. Wegen \fref{drth3bw1} gelten
\begin{align}	\label{drth3bw11}
	\Uarm^{(1,1)}\phi_k+\Uarm^{(1,2)}\psi_k
	&= \rklam{-i\Uarm^{(1,1)}+\Uarm^{(1,2)}}\tF_k+i\Uarm^{(1,1)}+\Uarm^{(1,2)} \notag \\
	&= \sqrt{2}\brklam{\Sarm^{(1,1)}\tF_k+\Sarm^{(1,2)}}
\end{align}
und
\begin{align}	\label{drth3bw12}
	\Uarm^{(2,1)}\phi_k+\Uarm^{(2,2)}\psi_k 
	&= \rklam{-i\Uarm^{(2,1)}+\Uarm^{(2,2)}}\tF_k+i\Uarm^{(2,1)}+\Uarm^{(2,2)} \notag \\
	&= \sqrt{2}\brklam{\Sarm^{(2,1)}\tF_k+\Sarm^{(2,2)}}.
\end{align}
Wegen Teil (b) von \thref{drth1} ist $\det\big(\Uarm^{(2,1)}\phi_k+\Uarm^{(2,2)}\psi_k\big)$ nicht die Nullfunktion. Hieraus folgt wegen (ii), (iii), \fref{drth3bw11} und \fref{drth3bw12} dann
\begin{align*}
	& \rklam{\Uarm^{(1,1)}\phi_1+\Uarm^{(1,2)}\psi_1}\rklam{\Uarm^{(2,1)}\phi_1+\Uarm^{(2,2)}\psi_1}^{-1} \\
	&= \brklam{\Sarm^{(1,1)}\tF_1+\Sarm^{(1,2)}}\brklam{\Sarm^{(2,1)}\tF_1+\Sarm^{(2,2)}}^{-1} \\
	&= S \\
	&= \brklam{\Sarm^{(1,1)}\tF_2+\Sarm^{(1,2)}}\brklam{\Sarm^{(2,1)}\tF_2+\Sarm^{(2,2)}}^{-1} \\
	&= \rklam{\Uarm^{(1,1)}\phi_2+\Uarm^{(1,2)}\psi_2}\rklam{\Uarm^{(2,1)}\phi_2+\Uarm^{(2,2)}\psi_2}^{-1}
\end{align*}
und somit wegen Teil (c) von \thref{drth1} dann $\sklam{\binom{\phi_1}{\psi_1}} = \sklam{\binom{\phi_2}{\psi_2}}$, das heißt wegen \thref{spdef4b} existiert eine in $\C\setminus[\alpha,\infty)$ meromorphe \textit{q}$\times$\textit{q}-Matrixfunktion $h$, so dass $\det h$ nicht die Nullfunktion ist und
\begin{align}	\label{drth3bw13}
	\binom{\phi_2}{\psi_2} = \binom{\phi_1}{\psi_1}h
\end{align}
erfüllt ist. Sei $k\in\gklam{1,2}$. Wegen (i), \fref{drth3bw14} und Teil (b) von \thref{spsa1} gilt dann weiterhin, dass $\det(\psi_k-i\phi_k)$ nicht die Nullfunktion ist und
\begin{align*}
	\tF_k = \rklam{\psi_k+i\phi_k}\rklam{\psi_k-i\phi_k}^{-1}
\end{align*}
erfüllt ist. Unter Beachtung, dass $\det h$ nicht die Nullfunktion ist, folgt hieraus wegen \fref{drth3bw13} nun
\begin{align*}
	\tF_2 &= \rklam{\psi_2+i\phi_2}\rklam{\psi_2-i\phi_2}^{-1} \\
	&= \rklam{\psi_1+i\phi_1}hh^{-1}\rklam{\psi_1-i\phi_1}^{-1} \\
	&= \rklam{\psi_1+i\phi_1}\rklam{\psi_1-i\phi_1}^{-1} = \tF_1. \tag*{$\Box$}
\end{align*}

Seien $m\in\N$, $\alpha\in\R$ und $\sjm\in\Kpqma$. Weiterhin sei $\Sarmmin$ bzw. $\Sarmmax$ das untere bzw. obere Extremalelement von $\Sqasmu$. Aus dem Beweis von \thref{drbsp1} geht hervor, dass $\Sarmmin$ bzw. $\Sarmmax$ die in Teil (a) von \thref{drth1} zum \linebreak \textit{q}$\times$\textit{q}-Stieltjes-Paar $\binom{\cal O}{\cal I}$ bzw. $\binom{\cal I}{\cal O}$ in $\C\setminus[\alpha,\infty)$ zugehörige Funktion aus $\Soqa[\sjm,$ $\leq]$ darstellt, wobei ${\cal I}$ bzw. ${\cal O}$ die in $\C\setminus[\alpha,\infty)$ konstante Matrixfunktion mit dem Wert $\Iq$ bzw. $\Oq$ bezeichnet. Es folgt aus dem Beweis von \thref{drth3} nun, dass $\Sarmmin$ bzw. $\Sarmmax$ auch in Teil (a) von \thref{drth3} mithilfe der Funktion ${\cal I}$ bzw. $-{\cal I}$ gewonnen werden kann. Mit den Bezeichnungen aus \thref{drth3} und \thref{drdef4} gelten wegen \fref{drth3bw1} und \thref{drdef2} nämlich
\begin{align*}
	&\ \brklam{\Sarm^{(1,1)}{\cal I}+\Sarm^{(1,2)}}\brklam{\Sarm^{(2,1)}{\cal I}+\Sarm^{(2,2)}}^{-1} \\
	&= \brklam{2\Uarm^{(1,2)}}\brklam{2\Uarm^{(2,2)}}^{-1}
	= \Barfm\Darfm^{-1} = \Sarmmin
\end{align*}
und
\begin{align*}
	&\ \brklam{\Sarm^{(1,1)}[-{\cal I}]+\Sarm^{(1,2)}}\brklam{\Sarm^{(2,1)}[-{\cal I}]+\Sarm^{(2,2)}}^{-1} \\
	&= \brklam{2i\Uarm^{(1,1)}}\brklam{2i\Uarm^{(2,1)}}^{-1}
	= \Aarfm\Carfm^{-1} = \Sarmmax.
\end{align*}

\subsection{Der linksseitige Fall}

In diesem Abschnitt werden wir die Lösungsmenge des zu einer endlichen linksseitig $\alpha$-Stieltjes-positiv definiten Folge zugehörigen mithilfe der $(-\infty,\alpha]$-Stieltjes"=Transformation umformulierten linksseitigen $\alpha$-Stieltjes Momentenproblems unter Verwendung einer gebrochen linearen Transformation so beschreiben, dass wir als erzeugende Matrixfunktion ein Matrixpolynom verwenden, welches gewisse Kontraktivitätseigenschaften bezüglich 
$\jqq := \diag(\Iq,-\Iq)$
und $\tJq$ besitzt. Weiterhin verwenden wir als Parametermenge nun eine Teilklasse der Menge von \textit{q}$\times$\textit{q}-Schurfunktionen auf $\Pp$ (vergleiche \thref{spdef6}). Für die Beweisführung werden wir hauptsächlich auf die Resultate des rechtsseitigen Falles zurückgreifen.

Bevor wir uns aber dieser alternativen Beschreibung der Lösungsmenge widmen, wollen wir uns zunächst oben erwähnten Kontraktivitätseigenschaften der erzeugenden Matrixfunktion widmen.

\begin{satz}	\thlabel{drlsa8}
	Seien $\kappa\in\Na$, $\alpha\in\R$, $\sjk\in\Lpqka$ und $\Ualmk$ die Folge von linksseitigen $2$\textit{q}$\times2$\textit{q}-$\alpha$-Dyukarev-Matrixpolynomen bezüglich $\sjk$. Weiterhin seien
	\begin{align*}
  		\widetilde{E} := \frac{1}{\sqrt{2}}\begin{pmatrix} -i\Iq & -i\Iq \\ -\Iq & \Iq \end{pmatrix} \quad \text{und} \quad \Salm := \Ualm \widetilde{E}
  	\end{align*}
  	für alle $m\in\Zok$. Dann gilt $\Salm\in\tbtJJPp$ für alle $m\in\Zok$.
\end{satz}

\bwanf Dies folgt wegen Teil (c) von \thref{drlbm2} sogleich aus Teil (b2) von \thref{spbm9}. \bwend

Wir kommen nun zur Beschreibung der Lösungsmenge des zu einer endlichen linksseitig $\alpha$-Stieltjes-positiv definiten Folge zugehörigen mithilfe der $(-\infty,\alpha]$-Stieltjes"=Transformation umformulierten linksseitigen $\alpha$-Stieltjes Momentenproblems mit $\Eqa$ (vergleiche Teil (c) von \thref{spdef6}) als Parametermenge.

\begin{theo}	\thlabel{drlth3}
  Seien $m \in \N$, $\alpha \in \R$, $\sjm \in \Lpqma$, $\Ualm$ das linksseitige $2$\textit{q}$\times2$\textit{q}-$\alpha$-Dyukarev-Matrixpolynom bezüglich $\sjm$,
  \begin{align*}
  	\widetilde{E} := \frac{1}{\sqrt{2}}\begin{pmatrix} -i\Iq & -i\Iq \\ -\Iq & \Iq \end{pmatrix} \quad \text{und} \quad \Salm := \Ualm \widetilde{E}.
  \end{align*}
  Weiterhin bezeichne
  \begin{align*}
    \Salm = \begin{pmatrix} \Salm^{(1,1)} & \Salm^{(1,2)} \\ \Salm^{(2,1)} & \Salm^{(2,2)} \end{pmatrix}
  \end{align*}
  die \textit{q}$\times$\textit{q}-Blockzerlegung von $\Salm$. \dgfa
  \begin{itemize}
    \item [\rm{(a)}]  Sei $F\in\Eqa$. Weiterhin seien
	\begin{align*}
		\MF := \gklam{z\in\Pp \;|\; \det F(\za) = 0}
	\end{align*}
	und $\tF: \Pm\cup(\alpha,\infty)\cup(\Pp\setminus\MF) \rightarrow \Cqq$ definiert gemäß
	\begin{align*}
		\tF(z) := \begin{cases} F(z) & \text{falls } z\in\Pm \\
		\lim_{\omega\rightarrow z}F(\omega) & \text{falls } z\in(\alpha,\infty) \\
		F^{-\ast}(\za) & \text{falls } z\in\Pp\setminus\MF \end{cases}.
	\end{align*}
	Dann ist $\det\big(\Salm^{(2,1)}\tF+\Salm^{(2,2)}\big)$ nicht die Nullfunktion und
	\begin{align*}
      	S := \brklam{\Salm^{(1,1)}\tF+\Salm^{(1,2)}}\brklam{\Salm^{(2,1)}\tF+\Salm^{(2,2)}}^{-1}
    \end{align*}
    gehört zu $\Sqmasmu$.
    \item [\rm{(b)}] Sei $S\in\Sqmasmu$. Dann existiert genau eine in $\C\setminus(-\infty,\alpha]$ meromorphe \textit{q}$\times$\textit{q}-Matrixfunktion $\tF$ mit $\Rstr_{\Pm}\tF\in\Eqa$ derart, dass
  	$\det\big(\Salm^{(2,1)}\tF+\Salm^{(2,2)}\big)$ nicht die Nullfunktion ist und
    \begin{align*}
      	S = \brklam{\Salm^{(1,1)}\tF+\Salm^{(1,2)}}\brklam{\Salm^{(2,1)}\tF+\Salm^{(2,2)}}^{-1}
    \end{align*}
    erfüllt ist.
  \end{itemize}
\end{theo}

\bwanf Seien $t_j:=(-1)^js_j$ für alle $j\in\Zom$. Wegen Teil (a) von \thref{asmbm3} gilt dann $\tjm\in\Kpqmma$. Weiterhin seien $U^{\sklam{t}}_{-\arm}$ das rechtsseitige $2$\textit{q}$\times2$\textit{q}-$-\alpha$-Dyukarev-Matrixpolynom bezüglich $\tjm$,
\begin{align*}
  	E := \frac{1}{\sqrt{2}}\begin{pmatrix} -i\Iq & i\Iq \\ \Iq & \Iq \end{pmatrix} \quad \text{und} \quad \Sigma^{\sklam{t}}_{-\arm} := U^{\sklam{t}}_{-\arm} E.
\end{align*}
Es bezeichne
\begin{align*}
 	\Sigma^{\sklam{t}}_{-\arm} = \begin{pmatrix} \big(\Sigma^{\sklam{t}}_{-\arm}\big)^{(1,1)} & \big(\Sigma^{\sklam{t}}_{-\arm}\big)^{(1,2)} \\ \big(\Sigma^{\sklam{t}}_{-\arm}\big)^{(2,1)} & \big(\Sigma^{\sklam{t}}_{-\arm}\big)^{(2,2)} \end{pmatrix}
\end{align*}
die \textit{q}$\times$\textit{q}-Blockzerlegung von $\Sigma^{\sklam{t}}_{-\arm}$. Aus der Definition von $\widetilde{E}$ und $E$ folgt
\begin{align*}
	\widetilde{E} = \Ve E \Vea.
\end{align*} 
Hieraus folgt wegen Teil (a) von \thref{drllm4} und Teil (a) von \thref{asmlm1} dann
\begin{align}	\label{drlth3bw1}
	\Salm^{\sklam{s}}(z) 
	= \Ualm^{\sklam{s}}(z) \widetilde{E} 
	= \Ve U^{\sklam{t}}_{-\arm}(-z) \Vea \Ve E \Vea
	= \Ve \Sigma^{\sklam{t}}_{-\arm}(-z) \Vea
\end{align}
für alle $z\in\C$.

Zu (a): Sei $G:\Pp\rightarrow\Cqq$ definiert gemäß $G(z) = -F(-z)$. Wegen Teil (b) von \thref{spbm4} gilt dann $G\in\Dqma$. Weiterhin seien
\begin{align*}
	{\cal M}_G := \gklam{z\in\Pm \;|\; \det F(\za) = 0}
\end{align*}
und $\widetilde{G}: \Pp\cup(-\infty,-\alpha)\cup(\Pm\setminus{\cal M}_G) \rightarrow \Cqq$ definiert gemäß
\begin{align*}
	\widetilde{G}(z) := \begin{cases} G(z) & \text{falls } z\in\Pp \\
	\lim_{\omega\rightarrow z}G(\omega) & \text{falls } z\in(-\infty,-\alpha) \\
	G^{-\ast}(\za) & \text{falls } z\in\Pm\setminus{\cal M}_G \end{cases}.
\end{align*}
Aus der Definition von $\tF$ und $\widetilde{G}$ folgt
\begin{align}	\label{drlth3bw2}
	\widetilde{F}(z) = -\widetilde{G}(-z)
\end{align}
für alle $z\in\Pm\cup(\alpha,\infty)\cup(\Pp\setminus\MF)$.
Unter Beachtung von $G\in\Dqma$ gilt wegen Teil (a) von \thref{drth3} dann, dass $\det\beklam{(\Sigma^{\sklam{t}}_{-\arm})^{(2,1)}\widetilde{G}+(\Sigma^{\sklam{t}}_{-\arm})^{(2,2)}}$ nicht die Nullfunktion ist und
\begin{align*}
    T := \eklam{\big(\Sigma^{\sklam{t}}_{-\arm}\big)^{(1,1)}\widetilde{G}+\big(\Sigma^{\sklam{t}}_{-\arm}\big)^{(1,2)}}\eklam{\big(\Sigma^{\sklam{t}}_{-\arm}\big)^{(2,1)}\widetilde{G}+\big(\Sigma^{\sklam{t}}_{-\arm}\big)^{(2,2)}}^{-1}
\end{align*}
zu ${\cal S}_{0,q,[-\alpha,\infty)}[\tjm$,""$\leq]$ gehört. Hieraus folgt wegen \fref{drlth3bw1} und \fref{drlth3bw2} nun, dass wegen
\begin{align*}
	&\ \big(\Salm^{\sklam{s}}\big)^{(2,1)}(z)\tF(z)+\big(\Salm^{\sklam{s}}\big)^{(2,2)}(z) \\
	&= -\big(\Sigma^{\sklam{t}}_{-\arm}\big)^{(2,1)}(-z)\big(\negthinspace-\widetilde{G}\big)(-z)+\big(\Sigma^{\sklam{t}}_{-\arm}\big)^{(2,2)}(-z) \\
	&= \big(\Sigma^{\sklam{t}}_{-\arm}\big)^{(2,1)}(-z)\widetilde{G}(-z)+\big(\Sigma^{\sklam{t}}_{-\arm}\big)^{(2,2)}(-z)
\end{align*}
für alle $z\in\Pm\cup(\alpha,\infty)\cup(\Pp\setminus\MF)$ dann $\det\beklam{(\Salm^{\sklam{s}})^{(2,1)}\tF+(\Salm^{\sklam{s}})^{(2,2)}}$ nicht die Nullfunktion ist und eine diskrete Teilmenge $\D$ von $\Pm\cup(\alpha,\infty)\cup(\Pp\setminus\MF)$ existiert, so dass
\begin{align}	\label{drlth3bw3}
	S(z) &= \eklam{\big(\Salm^{\sklam{s}}\big)^{(1,1)}(z)\tF(z)+\big(\Salm^{\sklam{s}}\big)^{(1,2)}(z)} \notag \\
	&\quad \cdot\eklam{\big(\Salm^{\sklam{s}}\big)^{(2,1)}(z)\tF(z)+\big(\Salm^{\sklam{s}}\big)^{(2,2)}(z)}^{-1} \notag \\
	&= \eklam{\big(\Sigma^{\sklam{t}}_{-\arm}\big)^{(1,1)}(-z)\big(\negthinspace-\widetilde{G}\big)(-z)-\big(\Sigma^{\sklam{t}}_{-\arm}\big)^{(1,2)}(-z)} \notag \\
	&\quad \cdot\eklam{-\big(\Sigma^{\sklam{t}}_{-\arm}\big)^{(2,1)}(-z)\big(\negthinspace-\widetilde{G}\big)(-z)+\big(\Sigma^{\sklam{t}}_{-\arm}\big)^{(2,2)}(-z)}^{-1} \notag \\
	&= -\eklam{\big(\Sigma^{\sklam{t}}_{-\arm}\big)^{(1,1)}(-z)\widetilde{G}(-z)+\big(\Sigma^{\sklam{t}}_{-\arm}\big)^{(1,2)}(-z)} \notag \\
	&\quad \cdot\eklam{\big(\Sigma^{\sklam{t}}_{\arm}\big)^{(2,1)}(-z)\widetilde{G}(-z)+\big(\Sigma^{\sklam{t}}_{-\arm}\big)^{(2,2)}(-z)}^{-1} = -T(-z)
\end{align}
für alle $z\in(\Pm\cup(\alpha,\infty)\cup(\Pp\setminus\MF))\setminus\D$ erfüllt ist. Sei $\widetilde{\D} = \D\cup\MF$. Unter Beachtung von $\C\setminus(-\infty,\alpha] = \Pm\cup(\alpha,\infty)\cup\Pp$ und der Tatsache, dass $\MF$ eine diskrete Teilmenge von $\Pp$ ist, ist dann $\widetilde{\D}$ eine diskrete Teilmenge von $\C\setminus(-\infty,\alpha]$ und wegen \fref{drlth3bw3} gilt
\begin{align*}
	S(z) = -T(-z)
\end{align*}
für alle $z\in\C\setminus\big((-\infty,\alpha]\cup\widetilde{D}\big)$. Hieraus folgt wegen des Identitätssatzes für meromorphe Funktionen (vergleiche z.\,B. im skalaren Fall \cite[Satz 10.3.2]{Funk}; im matriziellen Fall betrachtet man die einzelnen Einträge der Matrixfunktion) und Teil (d) von \thref{asmbm5} dann $S\in\Sqmasmu$.

Zu (b): Sei $\widecheck{S}:[-\alpha,\infty)\rightarrow\Cqq$ definiert gemäß $\widecheck{S}(z):=-S(-z)$. Wegen Teil (d) von \thref{asmbm5} gilt dann $\widecheck{S}\in{\cal S}_{0,q,[-\alpha,\infty)}[\tjm$,""$\leq]$. Hieraus folgt wegen Teil (b) von \thref{drth3} dann, dass genau eine in $\C\setminus[-\alpha,\infty)$ meromorphe \textit{q}$\times$\textit{q}-Matrixfunktion $\widetilde{G}$ mit $\Rstr_{\Pp}\widetilde{G}\in\Dqma$ derart existiert, dass $\det\beklam{(\Sigma^{\sklam{t}}_{-\arm})^{(2,1)}\widetilde{G}+(\Sigma^{\sklam{t}}_{-\arm})^{(2,2)}}$ nicht die Nullfunktion ist und
\begin{align*}
    \widecheck{S} := \eklam{\big(\Sigma^{\sklam{t}}_{-\arm}\big)^{(1,1)}\widetilde{G}+\big(\Sigma^{\sklam{t}}_{-\arm}\big)^{(1,2)}}\eklam{\big(\Sigma^{\sklam{t}}_{-\arm}\big)^{(2,1)}\widetilde{G}+\big(\Sigma^{\sklam{t}}_{-\arm}\big)^{(2,2)}}^{-1}
\end{align*}
erfüllt ist. Sei nun $\widetilde{F}:\C\setminus(-\infty,\alpha]$ definiert gemäß $\widetilde{F}(z) = -\widetilde{G}(-z)$. Unter Beachtung von $\Rstr_{\Pp}\widetilde{G}\in\Dqma$ gilt wegen Teil (b) von \thref{spbm4} dann $\Rstr_{\Pm}\widetilde{F}\in\Eqa$. Wegen \fref{drlth3bw1} gilt
\begin{align*}
	&\ \big(\Salm^{\sklam{s}}\big)^{(2,1)}(z)\tF(z)+\big(\Salm^{\sklam{s}}\big)^{(2,2)}(z) \\
	&= -\big(\Sigma^{\sklam{t}}_{-\arm}\big)^{(2,1)}(-z)\big({-\widetilde{G}}\big)(-z)+\big(\Sigma^{\sklam{t}}_{-\arm}\big)^{(2,2)}(-z) \\
	&= \big(\Sigma^{\sklam{t}}_{-\arm}\big)^{(2,1)}(-z)\widetilde{G}(-z)+\big(\Sigma^{\sklam{t}}_{-\arm}\big)^{(2,2)}(-z)
\end{align*}
für alle $z\in\Pm$. Unter Beachtung, dass $\det\beklam{(\Sigma^{\sklam{t}}_{-\arm})^{(2,1)}\widetilde{G}+(\Sigma^{\sklam{t}}_{-\arm})^{(2,2)}}$ nicht die Nullfunktion ist, ist dann auch $\det\beklam{(\Salm^{\sklam{s}})^{(2,1)}\tF+(\Salm^{\sklam{s}})^{(2,2)}}$ nicht die Nullfunktion und es existiert eine diskrete Teilmenge $\D$ von $\Pm$ mit
\begin{align*}
	S(z) &= -\widecheck{S}(-z) \\
	&= -\eklam{\big(\Sigma^{\sklam{t}}_{-\arm}\big)^{(1,1)}(-z)\widetilde{G}(-z)+\big(\Sigma^{\sklam{t}}_{-\arm}\big)^{(1,2)}(-z)} \\
	&\quad \cdot\eklam{\big(\Sigma^{\sklam{t}}_{\arm}\big)^{(2,1)}(-z)\widetilde{G}(-z)+\big(\Sigma^{\sklam{t}}_{-\arm}\big)^{(2,2)}(-z)}^{-1} \\
	&= \eklam{\big(\Sigma^{\sklam{t}}_{-\arm}\big)^{(1,1)}(-z)\big({-\widetilde{G}}\big)(-z)-\big(\Sigma^{\sklam{t}}_{-\arm}\big)^{(1,2)}(-z)} \\
	&\quad \cdot\eklam{-\big(\Sigma^{\sklam{t}}_{-\arm}\big)^{(2,1)}(-z)\big({-\widetilde{G}}\big)(-z)+\big(\Sigma^{\sklam{t}}_{-\arm}\big)^{(2,2)}(-z)}^{-1} \\
	&= \eklam{\big(\Salm^{\sklam{s}}\big)^{(1,1)}(z)\tF(z)+\big(\Salm^{\sklam{s}}\big)^{(1,2)}(z)} \\
	&\quad \cdot\eklam{\big(\Salm^{\sklam{s}}\big)^{(2,1)}(z)\tF(z)+\big(\Salm^{\sklam{s}}\big)^{(2,2)}(z)}^{-1}
\end{align*}
für alle $z\in\Pm\setminus\D$. Unter Beachtung, dass $\Rstr_{\Pm}\widetilde{F}$ eine \textit{q}$\times$\textit{q}-meromorphe Matrixfunktion auf $\Pm$ ist, folgt hieraus wegen des Identitätssatzes für meromorphe Funktionen (vergleiche z.\,B. im skalaren Fall \cite[Satz 10.3.2]{Funk}; im matriziellen Fall betrachtet man die einzelnen Einträge der Matrixfunktion) dann
\begin{align*}
	S = \eklam{\big(\Salm^{\sklam{s}}\big)^{(1,1)}\tF+\big(\Salm^{\sklam{s}}\big)^{(1,2)}} \eklam{\big(\Salm^{\sklam{s}}\big)^{(2,1)}\tF+\big(\Salm^{\sklam{s}}\big)^{(2,2)}}^{-1}.
\end{align*}
Die Einzigartigkeit von $\widetilde{F}$ folgt unmittelbar aus der Einzigartigkeit von $\widetilde{G}$, Teil (b) von \thref{spbm4} und dem eben erwähnten Identitätssatz für meromorphe Funktionen. \bwend

Seien $m\in\N$, $\alpha\in\R$ und $\sjm\in\Lpqma$. Weiterhin sei $\Salmmin$ bzw. $\Salmmax$ das untere bzw. obere Extremalelement von $\Sqmasmu$. Aus dem Beweis von \thref{drlbsp1} geht hervor, dass $\Salmmin$ bzw. $\Salmmax$ die in Teil (a) von \thref{drlth2} zum \linebreak \textit{q}$\times$\textit{q}-Stieltjes-Paar $\binom{\cal I}{\cal O}$ bzw. $\binom{\cal O}{\cal I}$ in $\C\setminus(-\infty,\alpha]$ zugehörige Funktion aus \linebreak $\Sqmasmu$ darstellt, wobei ${\cal I}$ bzw. ${\cal O}$ die in $\C\setminus(-\infty,\alpha]$ konstante Matrixfunktion mit dem Wert $\Iq$ bzw. $\Oq$ bezeichnet. Es folgt nun, dass $\Salmmin$ bzw. $\Salmmax$ auch in Teil (a) von \thref{drlth3} mithilfe der Funktion ${\cal I}$ bzw. $-{\cal I}$ gewonnen werden kann. Mit den Bezeichnungen aus \thref{drlth3} und \thref{drldef4} gelten wegen
\begin{align*}
	\Salm = \frac{1}{\sqrt{2}} \begin{pmatrix} -i\Ualm^{(1,1)}-\Ualm^{(1,2)} & -i\Ualm^{(1,1)}+\Ualm^{(1,2)} \\ -i\Ualm^{(2,1)}-\Ualm^{(2,2)} & -i\Ualm^{(2,1)}+\Ualm^{(2,2)} \end{pmatrix}
\end{align*}
und \thref{drldef3} nämlich
\begin{align*}
	&\ \brklam{\Salm^{(1,1)}{\cal I}+\Salm^{(1,2)}}\brklam{\Salm^{(2,1)}{\cal I}+\Salm^{(2,2)}}^{-1} \\
	&= \brklam{{-2i}\Ualm^{(1,1)}}\brklam{{-2i}\Ualm^{(2,1)}}^{-1}
	= \Aalfm\Calfm^{-1} = \Salmmin
\end{align*}
und
\begin{align*}
	&\ \brklam{\Salm^{(1,1)}[-{\cal I}]+\Salm^{(1,2)}}\brklam{\Salm^{(2,1)}[-{\cal I}]+\Salm^{(2,2)}}^{-1} \\
	&= \brklam{2\Ualm^{(1,2)}}\brklam{2\Ualm^{(2,2)}}^{-1}
	= \Balfm\Dalfm^{-1} = \Salmmax.
\end{align*}

\newpage
\section[Das \texorpdfstring{$\alpha$}{a}-Stieltjes-Quadrupel bezüglich \texorpdfstring{$\alpha$}{a}-Stieltjes-positiv definiter Folgen]{Das \texorpdfstring{$\alpha$}{a}-Stieltjes-Quadrupel bezüglich \texorpdfstring{$\alpha$}{a}-Stieltjes- \\ positiv definiter Folgen} \label{chapsq}

Im Mittelpunkt des vorliegenden Abschnitts steht eine weitere Untersuchung der in Kapitel \ref{chapdr} eingeführten Folge von $2$\textit{q}$\times2$\textit{q}-$\alpha$-Dyukarev-Matrixpolynomen bezüglich einer $\alpha$-Stieltjes-positiv definiten Folge. Wir zeigen nun, dass diese in enger Beziehung zu dem in Abschnitt \ref{chapmp} eingeführten monischen links-orthogonalen System von Matrixpolynomen bezüglich einer Hankel-positiv definiten Folge (vergleiche \thref{mpdef2} und \thref{mpsa2}) sowie dem zugehörigen linken System von Matrixpolynomen zweiter Art (vergleiche \thref{mpdef3}) stehen. Unter Verwendung dieser beiden Systeme führen wir das $\alpha$-Stieltjes-Quadrupel bezüglich einer $\alpha$-Stieltjes-positiv definiten Folge ein und stellen neben den bereits angesprochenen Beziehungen zur Folge von $2$\textit{q}$\times2$\textit{q}-$\alpha$-Dyukarev-Matrixpolynomen bezüglich jener Folge einige weitere Beziehungen zu dem $\alpha$-Dyukarev-Quadrupel, der $\alpha$-Stieltjes-Parametrisierung und dem Favard-Paar jener Folge zusammen.

Im rechtsseitigen Fall orientieren wir uns an der Vorgehensweise von \cite[Chapter 4]{CR1}, wo eine gegebene Folge aus ${\cal K}^{>}_{q,\infty,0}$ zugrunde gelegt wurde. Im linksseitigen Fall greifen wir wieder hauptsächlich auf die Resultate für den rechtsseitigen Fall zurück.

\subsection{Der rechtsseitige Fall}

Unsere folgende Überlegung ist auf die Einführung eines Quadrupels von Folgen von \textit{q}$\times$\textit{q}-Matrixfunktionen ausgerichtet, welches der in Abschnitt \ref{chapmp} eingeführten in unserem Fall mit einer rechtsseitig $\alpha$-Stieltjes-positiv definiten Folge und derer durch rechtsseitige $\alpha$-Verschiebung generierten Folge assoziierten monischen links-orthogonalen Systeme von Matrixpolynomen sowie der hiermit assoziierten linken Systeme von Matrixpolynomen zweiter Art zugrunde liegt. Zuvor wollen wir einige Resultate für besagte Systeme rekapitulieren und erweitern. Dafür benötigen wir zunächst einige weitere Bezeichnungen.

\begin{bez}	\thlabel{sqbz1}
	Seien $\kappa \in \Na$, $\alpha \in \R$ und $\sjk$ eine Folge aus $\Cqq$. Dann bezeichnet für $n \in \Zokm$
	\begin{align*}
  		\Sarn := \begin{cases} s_{\aro} & \text{falls } n=0 \\
  		\begin{pmatrix} s_{\aro} & \Oq & \ldots & \Oq \\ s_{\ar 1} & s_{\aro} & \ddots & \vdots \\ \vdots & \ddots & \ddots & \Oq \\ s_{\arn} & \ldots & s_{\ar 1} & s_{\aro} \end{pmatrix} & \text{falls } n>0 \end{cases}
	\end{align*}
	die untere n-te Blockdreiecksmatrix von $\sarjk$. Weiterhin seien $\dSaro := s_0$ und für $n \in \Zek$
	\begin{align*}
  		\dSarn := \begin{pmatrix} s_0 & \Oq & \ldots & \Oq \\ s_{\ar o} & s_0 & \ddots & \vdots \\ \vdots & \ddots & \ddots & \Oq \\ s_{\arn-1} & \ldots & s_{\aro} & s_0 \end{pmatrix} 
  = S_n-\alpha \begin{pmatrix} \Oqn & \Oq \\ S_{n-1} & \Onq \end{pmatrix}.
	\end{align*}
\end{bez}

\begin{satz}	\thlabel{sqsa1}
  Seien $\kappa \in \Na$, $\alpha \in \R$ und $\sjk \in \Kpqka$. \dgfa
  \begin{itemize}
    \item [\rm{(a)}] Sei $\Pnfk$ das monische links-orthogonale System von Matrixpolynomen bezüglich $\sjk$. Dann gelten $P_0 \equiv \Iq$ und
      \begin{align*}
	P_n(z) = \begin{pmatrix} -z_{n,2n-1}\Hnm^{-1} & \Iq \end{pmatrix}E_n(z)
      \end{align*}
      für alle $z \in \C$ und $n \in \Zefkp$.
    \item [\rm{(b)}] Sei $\Psnfk$ das linke System von Matrixpolynomen zweiter Art bezüglich \linebreak $\sjk$. Dann gelten $\Ps_0 \equiv \Oq$ und
      \begin{align*}
	\Ps_n(z) = \begin{pmatrix} -z_{n,2n-1}\Hnm^{-1} & \Iq \end{pmatrix} \begin{pmatrix} \Oqn \\ S_{n-1} \end{pmatrix} E_{n-1}(z)
      \end{align*}
      für alle $z \in \C$ und $n \in \Zefkp$.
    \item [\rm{(c)}] Sei $\Parnfk$ das monische links-orthogonale System von Matrixpolynomen bezüglich $\sarjk$. Dann gelten $P_{\aro} \equiv \Iq$ und im Fall $\kappa\geq2$
      \begin{align*}
	P_{\arn}(z) = \begin{pmatrix} -z_{\arn,2n-1}\Harnm^{-1} & \Iq \end{pmatrix}E_n(z)
      \end{align*}
      für alle $z \in \C$ und $n \in \Zefk$.
    \item [\rm{(d)}] Sei $\Psarnfk$ das linke System von Matrixpolynomen zweiter Art bezüglich \linebreak $\sarjk$. Dann gelten $\Pars_{\aro} \equiv \Oq$ und im Fall $\kappa\geq2$
      \begin{align*}
	\Pars_{\arn}(z) = \begin{pmatrix} -z_{\arn,2n-1}\Harnm^{-1} & \Iq \end{pmatrix} \begin{pmatrix} \Oqn \\ S_{\arn-1} \end{pmatrix} E_{n-1}(z)
      \end{align*}
      für alle $z \in \C$ und $n \in \Zefk$.
    \item [\rm{(e)}] Seien $\widehat{P}_{\arn} := (z-\alpha)P_{\arn}(z)$ für alle $z \in \C$ und $n \in \Zofk$.
      Weiterhin sei $\dPs_{\arn}$ das zu $\sjk$ gehörige Matrixpolynom bezüglich $\widehat{P}_{\arn}$ für alle $n \in \Zofk$. 
      Dann gelten $\dPs_{\aro} \equiv s_0$ und im Fall $\kappa\geq2$
      \begin{align*}
	\dPs_{\arn}(z) = \begin{pmatrix} -z_{\arn,2n-1}\Harnm^{-1} & \Iq \end{pmatrix} \dSarn E_n(z)
      \end{align*}
      für alle $z \in \C$ und $n \in \Zefk$ sowie
      \begin{align*}
	\dPs_{\arn}(z) = \Pars_{\arn}(z) + P_{\arn}(z) s_0
      \end{align*}
      für alle $z \in \C$ und $n \in \Zofk$.
  \end{itemize}
\end{satz}

\bwanf Zu (a) - (d): Unter Beachtung von \thref{asmbm1} gilt wegen $\sjk \in \Kpqka$ und Teil (b) von \thref{asmdef2} $\sjfkm \in \Hpqfkm$ bzw. $(s_{\arj})^{2\fklam{\kappa-2}} \in \Hpqfkzm$, falls $\kappa \geq 2$.
Hieraus folgt wegen \thref{mpsa2} dann die Behauptung von (a) bzw. (c) und wegen \thref{mpbem1} die Behauptung von (b) bzw. (d).

Zu (e): Sei $n \in \Zofk$. Unter Beachtung von (c) gelten dann $\deg \dPs_{\arn} = n+1$ und
\begin{align}	\label{sqsa1bw1}
  \dP^{[j]}_{\arn} = \begin{cases} -\alpha P^{[0]}_{\arn} & \text{falls } j = 0 \\
  P^{[j-1]}_{\arn}-\alpha P^{[j]}_{\arn} & \text{falls } 1 \leq j \leq n \\
  P^{[n]}_{\arn} & \text{falls } j = n+1 \end{cases}
\end{align}
für alle $j \in \Zonp$. Sei nun $z \in \C$. Wegen \fref{sqsa1bw1} und (c) gelten dann
\begin{align*}
  \dPs_{\aro}(z) & = \begin{pmatrix} \dP^{[0]}_{\aro} & \dP^{[1]}_{\aro} \end{pmatrix} \begin{pmatrix} \Oq \\ S_0 \end{pmatrix} E_0(z) 
  = \begin{pmatrix} -\alpha P^{[0]}_{\aro} & P^{[0]}_{\aro} \end{pmatrix} \begin{pmatrix} \Oq \\ s_0 \end{pmatrix} \\
  & = \begin{pmatrix} -\alpha \Iq & \Iq \end{pmatrix} \begin{pmatrix} \Oq \\ s_0 \end{pmatrix} 
  = s_0
\end{align*}
und im Fall $\kappa\geq2$
\begin{align*}
  \dPs_{\arn}(z) & = \begin{pmatrix} \dP^{[0]}_{\arn} & \ldots & \dP^{[n+1]}_{\arn} \end{pmatrix} \begin{pmatrix} 0_{q\times(n+1)q} \\ S_n \end{pmatrix} \En(z) \\ 
  & = \eklam{\begin{pmatrix} \Oq & P^{[0]}_{\arn} & \ldots & P^{[n]}_{\arn} \end{pmatrix} - \alpha \begin{pmatrix} P^{[0]}_{\arn} & \ldots & P^{[n]}_{\arn} & \Oq \end{pmatrix}} \begin{pmatrix} 0_{q\times(n+1)q} \\ S_n \end{pmatrix} \En(z) \\ 
  & = \eklam{\begin{pmatrix} P^{[0]}_{\arn} & \ldots & P^{[n]}_{\arn} \end{pmatrix} S_n - \alpha \begin{pmatrix} P^{[0]}_{\arn} & \ldots & P^{[n]}_{\arn} \end{pmatrix} \begin{pmatrix} \Oqn & \Oq \\ S_{n-1} & \Onq \end{pmatrix}} \En(z) \\
  & = \begin{pmatrix} P^{[0]}_{\arn} & \ldots & P^{[n]}_{\arn} \end{pmatrix} \dSarn \En(z)
\end{align*}
für alle $n \in \Zefk$. Wegen (c) und (d) gelten weiterhin
\begin{align*}
  \Pars_{\aro}(z) + P_{\aro}(z)s_0 = s_0 = \dPs_{\aro}(z)
\end{align*}
sowie im Fall $\kappa\geq2$
\begin{align*}
  &\ \Pars_{\arn}(z) + P_{\arn}(z)s_0 \\
  & = \begin{pmatrix} -z_{\arn,2n-1}\Harnm^{-1} & \Iq \end{pmatrix} \begin{pmatrix} \Oqn \\ S_{\arn-1} \end{pmatrix} \Enm(z) + \begin{pmatrix} -z_{\arn,2n-1}\Harnm^{-1} & \Iq \end{pmatrix}\En(z)s_0 \\
  & = \begin{pmatrix} -z_{\arn,2n-1}\Harnm^{-1} & \Iq \end{pmatrix} \eklam{\begin{pmatrix} \Oqn & \Oq \\ S_{\arn-1} & \Onq \end{pmatrix} \En(z) + \diag(s_0,...,s_0)\En(z)} \\
  & = \begin{pmatrix} -z_{\arn,2n-1}\Harnm^{-1} & \Iq \end{pmatrix} \dSarn \En(z) = \dPs_{\arn}(z)
\end{align*}
für alle $n \in \Zefk$. \bwend

\begin{defi}	\thlabel{sqdef1}
  Seien $\kappa \in \Na$, $\alpha \in \R$ und $\sjk \in \Kpqka$. 
  Weiterhin seien $(P_{n,s})^{\fklam{\kappa+1}}_{n=0}$ das monische links-orthogonale System von Matrixpolynomen bezüglich $\sjk$, \linebreak
  $(\Ps_{n,s})^{\fklam{\kappa+1}}_{n=0}$ das linke System von Matrixpolynomen zweiter Art bezüglich $\sjk$, \linebreak
  $(P_{\arn,s})^{\fklam{\kappa}}_{n=0}$ das monische links-orthogonale System von Matrixpolynomen bezüglich \linebreak $\sarjk$,
  $\widehat{P}_{\arn,s} := (z-\alpha)P_{\arn,s}(z)$ für alle $z \in \C$ und $n \in \Zofk$
  sowie $\dPs_{\arn,s}$ das zu $\sjk$-gehörige Matrixpolynom bezüglich $\widehat{P}_{\arn,s}$ für alle $n \in \Zofk$.
  Dann heißt $\SQsark$ das \textbf{rechtsseitige $\alpha$-Stieltjes-Qua\-dru\-pel} bezüglich $\sjk$.
  Falls klar ist, von welchem $\sjk$ die Rede ist, lassen wir das \anf{$s$} im unteren Index weg.
\end{defi}

Unsere nachfolgenden Überlegungen sind nun auf die Herleitung von Zusammenhängen zwischen dem gerade eingeführten rechtsseitigen $\alpha$-Stieltjes-Quadrupel bezüglich einer rechtsseitig $\alpha$-Stieltjes-positiv definiten Folge und dem rechtsseitigen $\alpha$-Dyukarev-Quadrupel bezüglich jener Folge orientiert. Wir werden zunächst das rechtsseitige $\alpha$-Stieltjes-Quadrupel auf eine Form bringen, die mit den Bezeichnungen aus Kapitel \ref{chapdr} konform ist. Hierfür wird folgendes Lemma uns eine wichtige Hilfestellung leisten.

\begin{lemma}	\thlabel{sqlm1}
  Seien $\kappa \in \Na$, $\alpha \in \R$ und $\sjk$ eine Folge aus $\Cqq$. \dgfa
  \begin{itemize}
    \item [\rm{(a)}] Für alle $z \in \C$ und $n \in \No$ gilt
      \begin{align*}
	R_n(z)v_n = \En(z).
      \end{align*}    
    \item [\rm{(b)}] Für alle $z \in \C$ und $n \in \Zok$ gilt
      \begin{align*}
	R_n(z)y_{0,n} = S_n\En(z).
      \end{align*}
    \item [\rm{(c)}] Für alle $z \in \C$ und $n \in \Zekp$ gilt
      \begin{align*}
	R_n(z)u_n = \begin{pmatrix} \Oqn \\ S_{n-1} \end{pmatrix} \Enm(z).
      \end{align*}
    \item [\rm{(d)}] Für alle $z \in \C$ und $n \in \Zokm$ gilt
      \begin{align*}
	R_n(z)y_{\aro,n} = \Sarn\En(z).
      \end{align*}
    \item [\rm{(e)}] Für alle $z \in \C$ und $n \in \Zok$ gilt
      \begin{align*}
	R_n(z)u_{\arn} = \dSarn\En(z).
      \end{align*}
  \end{itemize}
\end{lemma}

\bwanf Zu (a): Dies folgt sogleich aus der Definition der beteiligten Größen (vergleiche \thref{drbz1} und \thref{adpbz1}).

Zu (b): Der Fall $n = 0$ folgt sogleich aus der Definition der beteiligten Größen (vergleiche \thref{drbz1}, \thref{asmbz1}, \thref{mpbz3} und \thref{adpbz1}). Seien nun $z \in \C$ und $n \in \Zek$. Dann gilt
\begin{align*}
  \Rn(z)\yn = \begin{pmatrix} s_0 \\ zs_0+s_1 \\ \vdots \\ \sum^{n}_{j=0} z^{n-j}s_j \end{pmatrix} = S_n\En(z).
\end{align*}

Zu (c): Seien $z \in \C$ und $n \in \Zekp$. Wegen \thref{drbz3} und (b) gilt dann
\begin{align*}
  \Rn(z)\un = \begin{pmatrix} \Oq \\ \Rnm(z)\ynm \end{pmatrix} = \begin{pmatrix} \Oq \\ S_{n-1}\Enm(z) \end{pmatrix} = \begin{pmatrix} \Oqn \\ S_{n-1} \end{pmatrix} \Enm(z).
\end{align*}

Zu (d): Der Fall $n = 0$ folgt sogleich aus der Definition der beteiligten Größen (vergleiche \thref{drbz1}, Teil (a) von \thref{aspdef1}, \thref{sqbz1} und \thref{adpbz1}). Seien nun $z \in \C$ und $n \in \Zek$. Dann gilt
\begin{align*}
  \Rn(z)\yarn = \begin{pmatrix} s_{\aro} \\ zs_{\aro}+s_{\ar1} \\ \vdots \\ \sum^{n}_{j=0} z^{n-j}s_{\arj} \end{pmatrix} = \Sarn\En(z).
\end{align*}

Zu (e): Der Fall $n = 0$ folgt sogleich aus der Definition der beteiligten Größen (vergleiche \thref{drbz1}, \thref{drbz3}, \thref{sqbz1} und \thref{adpbz1}). Seien nun $z \in \C$ und $n \in \Zekp$. Wegen (b) und (c) gilt dann
\begin{align*}
  \Rn(z)\uarn & = \Rn(z)\eklam{\yn-\alpha\un} = S_n\En(z)-\alpha \begin{pmatrix} \Oqn \\ S_{n-1} \end{pmatrix} \Enm(z) \\
  & = \eklam{S_n-\alpha \begin{pmatrix} \Oqn & \Oq \\ S_{n-1} & \Onq \end{pmatrix}}\En(z) = \dSarn\En(z). \tag*{$\Box$}
\end{align*}

\begin{bem}	\thlabel{sqbm1}
  Seien $\kappa \in \Na$, $\alpha \in \R$, $\sjk \in \Kpqka$ und $\SQark$ das rechtsseitige $\alpha$-Stieltjes-Quadrupel bezüglich $\sjk$.
  Dann gelten $P_0 \equiv \Iq$, $\Ps_0 \equiv \Oq$, $P_{\aro} \equiv \Iq$ und $\dPs_{\aro} \equiv s_0$ sowie
  \begin{align*}
    P_n(z) = \begin{pmatrix} -z_{n,2n-1}\Hnm^{-1} & \Iq \end{pmatrix}R_n(z)v_n, \\
    \Ps_n(z) = \begin{pmatrix} -z_{n,2n-1}\Hnm^{-1} & \Iq \end{pmatrix}R_n(z)u_n
  \end{align*}
  für alle $z \in \C$ und $n \in \Zefkp$ sowie im Fall $\kappa\geq2$
  \begin{align*}
    P_{\arn}(z) = \begin{pmatrix} -z_{\arn,2n-1}\Harnm^{-1} & \Iq \end{pmatrix}R_n(z)v_n, \\
    \dPs_{\arn}(z) = \begin{pmatrix} -z_{\arn,2n-1}\Harnm^{-1} & \Iq \end{pmatrix}R_n(z)u_{\arn}
  \end{align*}
  für alle $z \in \C$ und $n \in \Zefk$.
\end{bem}

\bwanf Dies folgt aus \thref{sqsa1} und \thref{sqlm1}. \bwend

Mithilfe von \thref{sqbm1} erkennen wir, dass die in \thref{sqdef1} eingeführten Größen unter Beachtung von \thref{sqsa1} mit denen aus \cite[Definition 4.1]{CR1} in Verbindung mit \cite[Proposition 4.2]{CR1} im Fall $\alpha=0$ und $\kappa=\infty$ übereinstimmen.

Wir kommen nun zu einem ersten Zusammenhang zwischen rechtsseitigem $\alpha$-Stieltjes-Quadrupel und rechtsseitigem $\alpha$-Dyukarev-Quadrupel bezüglich einer rechtsseitigen $\alpha$-Stieltjes-positiv definiten Folge (vergleiche \cite[Proposition 4.9(a)]{CR1} für den Fall $\alpha=0$ und $\kappa=\infty$, wobei dort in den Formeln (4.45) und (4.48) das falsche Vorzeichen gewählt wurde). 

\begin{satz}	\thlabel{sqsa2}
  Seien $\kappa \in \Na$, $\alpha \in \R$ und $\sjk \in \Kpqka$. Weiterhin seien $\SQark$ das rechtsseitige $\alpha$-Stieltjes-Quadrupel bezüglich $\sjk$
  und $\ABCDark$ das rechtsseitige $\alpha$-Dyukarev-Qua\-dru\-pel bezüglich $\sjk$. Dann gelten
  \begin{align*}
    \Aarn(z) & = \Iq + (z-\alpha) \sum^{n}_{j=0} \beklam{\Ps_j(\za)}^{\ast}\dH^{-1}_jP_j(\alpha), \\
    \Carn(z) & = - (z-\alpha) \sum^{n}_{j=0} \beklam{P_j(\za)}^{\ast}\dH^{-1}_jP_j(\alpha)
  \end{align*}
  für alle $z \in \C$ und $n \in \Zofk$ und
  \begin{align*}
    \Barn(z) & = \sum^{n-1}_{j=0} \beklam{\dPs_{\arj}(\za)}^{\ast}\dH^{-1}_{\arj}\dPs_{\arj}(\alpha), \\
    \Darn(z) & = \Iq - (z-\alpha) \sum^{n-1}_{j=0} \beklam{P_{\arj}(\za)}^{\ast}\dH^{-1}_{\arj}\dPs_{\arj}(\alpha)
  \end{align*}
  für alle $z \in \C$ und $n \in \Zefkp$.
\end{satz}

\bwanf Wegen \thref{drdef1} und \thref{sqbm1} gelten
\begin{align*}
  \Aaro(z) = \Iq = \Iq + (z-\alpha) \beklam{\Ps_0(\za)}^{\ast}\dH^{-1}_0P_0(\alpha)
\end{align*}
für alle $z \in \C$ und unter Beachtung von \thref{adplm2}, den Teilen (b) und (d) von \thref{drlm2} sowie Teil (c) von \thref{drlm3} weiterhin
\begin{align*}
  \Aarn(z) & = \Iq + (z-\alpha)u^{\ast}_nR^{\ast}_n(\za)H^{-1}_nR_n(\alpha)v_n \\
  & = \Iq + (z-\alpha)u^{\ast}_nR^{\ast}_n(\za) \begin{pmatrix} \Hnm^{-1} & \Onq \\ \Oqn & \Oq \end{pmatrix} R_n(\alpha)v_n \\
  &\quad + (z-\alpha)u^{\ast}_nR^{\ast}_n(\za) \begin{pmatrix} -\Hnm^{-1}y_{n,2n-1} \\ \Iq \end{pmatrix} \dHn^{-1} \begin{pmatrix} -z_{n,2n-1}\Hnm^{-1} & \Iq \end{pmatrix} R_n(\alpha)v_n \\
  & = \Iq + (z-\alpha)u^{\ast}_nR^{\ast}_n(\za)\dL_n\Hnm^{-1}\dL^{\ast}_nR_n(\alpha)v_n + (z-\alpha)\beklam{\Ps_n(\za)}^{\ast}\dHn^{-1}P_n(\alpha) \\
  & = \Iq + (z-\alpha)u^{\ast}_{n-1}R^{\ast}_{n-1}(\za)\Hnm^{-1}R_{n-1}(\alpha)v_{n-1} + (z-\alpha)\beklam{\Ps_n(\za)}^{\ast}\dHn^{-1}P_n(\alpha) \\
  & = \Aarnm(z) + (z-\alpha)\beklam{\Ps_n(\za)}^{\ast}\dHn^{-1}P_n(\alpha)
\end{align*}
für alle $z \in \C$ und $n \in \Zefk$. Hieraus folgt iterativ dann
\begin{align*}
  \Aarn(z) = \Iq + (z-\alpha) \sum^{n}_{j=0} \beklam{\Ps_j(\za)}^{\ast}\dH^{-1}_jP_j(\alpha)
\end{align*}
für alle $z \in \C$ und $n \in \Zofk$. Wegen \thref{drdef1} und \thref{sqbm1} gelten
\begin{align*}
  \Bare(z) = s_0s^{-1}_{\aro}s_0 = \beklam{\dPs_{\aro}(\za)}^{\ast}\dHaro^{-1}\dPs_{\aro}(\alpha)
\end{align*}
für alle $z \in \C$ und unter Beachtung von \thref{adplm2}, Teil (d) von \thref{drlm2} sowie den Teilen (b), (d) und (e) von \thref{drlm3} weiterhin
\begin{align*}
  \Barnp(z) & = \uarn^{\ast}R^{\ast}_n(\za)\Harn^{-1}y_{0,n} \\
  & = \uarn^{\ast}R^{\ast}_n(\za) \begin{pmatrix} \Harnm^{-1} & \Onq \\ \Oqn & \Oq \end{pmatrix} y_{0,n} \\
  &\quad + \uarn^{\ast}R^{\ast}_n(\za) \begin{pmatrix} -\Harnm^{-1}y_{\arn,2n-1} \\ \Iq \end{pmatrix} \dHarn^{-1} \begin{pmatrix} -z_{\arn,2n-1}\Harnm^{-1} & \Iq \end{pmatrix} R^{\ast}_n(\alpha)\uarn \\
  & = \uarn^{\ast}R^{\ast}_n(\za)\dL_n\Harnm^{-1}\dL^{\ast}_ny_{0,n} + \beklam{\dPs_{\arn}(\za)}^{\ast}\dHarn^{-1}\dPs_{\arn}(\alpha) \\
  & = \uarnm^{\ast}R^{\ast}_{n-1}(\za)\Harnm^{-1}y_{0,n-1} + \beklam{\dPs_{\arn}(\za)}^{\ast}\dHarn^{-1}\dPs_{\arn}(\alpha) \\
  & = \Barn(z) + \beklam{\dPs_{\arn}(\za)}^{\ast}\dHarn^{-1}\dPs_{\arn}(\alpha)  
\end{align*}
für alle $z \in \C$ und $n \in \Z{1}{\fklam{\kappa+3}}$. Hieraus folgt iterativ dann
\begin{align*}
  \Barn(z) = \sum^{n-1}_{j=0} \beklam{\dPs_{\arj}(\za)}^{\ast}\dH^{-1}_{\arj}\dPs_{\arj}(\alpha)
\end{align*}
für alle $z \in \C$ und $n \in \Zefkp$. Wegen \thref{drdef1} und \thref{sqbm1} gelten
\begin{align*}
  \Caro(z) = -(z-\alpha)s^{-1}_0 = -(z-\alpha)\beklam{P_0(\za)}^{\ast}\dH^{-1}_0P_0(\alpha)
\end{align*}
für alle $z \in \C$ und unter Beachtung von \thref{adplm2} sowie den Teilen (b) und (d) von \thref{drlm2} weiterhin
\begin{align*}
  \Carn(z) & = -(z-\alpha)v^{\ast}_nR^{\ast}_n(\za)H^{-1}_nR_n(\alpha)v_n \\
  & = -(z-\alpha)v^{\ast}_nR^{\ast}_n(\za) \begin{pmatrix} \Hnm^{-1} & \Onq \\ \Oqn & \Oq \end{pmatrix} R_n(\alpha)v_n \\
  &\quad -(z-\alpha)v^{\ast}_nR^{\ast}_n(\za) \begin{pmatrix} -\Hnm^{-1}y_{n,2n-1} \\ \Iq \end{pmatrix} \dHn^{-1} \begin{pmatrix} -z_{n,2n-1}\Hnm^{-1} & \Iq \end{pmatrix} R_n(\alpha)v_n \\
  & = -(z-\alpha)v^{\ast}_nR^{\ast}_n(\za)\dL_n\Hnm^{-1}\dL^{\ast}_nR_n(\alpha)v_n - (z-\alpha)\beklam{P_n(\za)}^{\ast}\dHn^{-1}P_n(\alpha) \\
  & = -(z-\alpha)v^{\ast}_{n-1}R^{\ast}_{n-1}(\za)\Hnm^{-1}R_{n-1}(\alpha)v_{n-1} - (z-\alpha)\beklam{P_n(\za)}^{\ast}\dHn^{-1}P_n(\alpha) \\
  & = \Carnm(z) - (z-\alpha)\beklam{P_n(\za)}^{\ast}\dHn^{-1}P_n(\alpha)
\end{align*}
für alle $z \in \C$ und $n \in \Zefk$. Hieraus folgt iterativ dann
\begin{align*}
  \Carn(z) = - (z-\alpha) \sum^{n}_{j=0} \beklam{P_j(\za)}^{\ast}\dH^{-1}_jP_j(\alpha)
\end{align*}
für alle $z \in \C$ und $n \in \Zofk$. Wegen \thref{drdef1} und \thref{sqbm1} gelten
\begin{align*}
  \Dare(z) = \Iq - (z-\alpha)s^{-1}_{\aro}s_0 = \Iq - (z-\alpha)\beklam{P_{\aro}(\za)}^{\ast}\dHaro^{-1}\dPs_{\aro}(\alpha)
\end{align*}
für alle $z \in \C$ und unter Beachtung von \thref{adplm2}, den Teilen (b) und (d) von \thref{drlm2} sowie den Teilen (b) und (e) von \thref{drlm3} weiterhin
\begin{align*}
  \Darnp(z) & = \Iq - (z-\alpha)v^{\ast}_nR^{\ast}_n(\za)\Harn^{-1}y_{0,n} \\
  & = \Iq - (z-\alpha)v^{\ast}_nR^{\ast}_n(\za) \begin{pmatrix} \Harnm^{-1} & \Onq \\ \Oqn & \Oq \end{pmatrix} y_{0,n} - (z-\alpha)v^{\ast}_nR^{\ast}_n(\za) \\
  &\quad \cdot \begin{pmatrix} -\Harnm^{-1}y_{\arn,2n-1} \\ \Iq \end{pmatrix} \dHarn^{-1} \begin{pmatrix} -z_{\arn,2n-1}\Harnm^{-1} & \Iq \end{pmatrix} R^{\ast}_n(\alpha)\uarn \\
  & = \Iq - (z-\alpha)v^{\ast}_nR^{\ast}_n(\za)\dL_n\Harnm^{-1}\dL^{\ast}_ny_{0,n} - (z-\alpha)\beklam{P_{\arn}(\za)}^{\ast}\dHarn^{-1}\dPs_{\arn}(\alpha) \\
  & = \Iq - (z-\alpha)v^{\ast}_{n-1}R^{\ast}_{n-1}(\za)\Harnm^{-1}y_{0,n-1} - (z-\alpha)\beklam{P_{\arn}(\za)}^{\ast}\dHarn^{-1}\dPs_{\arn}(\alpha) \\
  & = \Darn(z) - (z-\alpha)\beklam{P_{\arn}(\za)}^{\ast}\dHarn^{-1}\dPs_{\arn}(\alpha)
\end{align*}
für alle $z \in \C$ und $n \in \Z{1}{\fklam{\kappa+3}}$. Hieraus folgt iterativ dann
\begin{align*}
  \Darn(z) = \Iq - (z-\alpha) \sum^{n-1}_{j=0} \beklam{P_{\arj}(\za)}^{\ast}\dH^{-1}_{\arj}\dPs_{\arj}(\alpha)
\end{align*}
für alle $z \in \C$ und $n \in \Zefkp$. \bwend

Wir zeigen nun ein weiteres Zusammenspiel zwischen rechtsseitigem $\alpha$"=Stieltjes"=Quadrupel und rechtsseitigem $\alpha$"=Dyukarev"=Quadrupel bezüglich einer rechtsseitig $\alpha$"=Stieltjes"=positiv definiten Folge, das für unsere nächsten Überlegungen von Wichtigkeit ist, wir aber später noch einmal konkretisieren wollen (vergleiche \cite[Lemma 4.4]{CR1} im Fall $\alpha=0$ und $\kappa=\infty$).

\begin{bem}	\thlabel{sqbm2}
  Seien $\kappa \in \Na$, $\alpha \in \R$ und $\sjk \in \Kpqka$. Weiterhin seien $\SQark$ das rechtsseitige $\alpha$-Stieltjes-Quadrupel bezüglich $\sjk$
  und $\ABCDark$ das rechtsseitige $\alpha$-Dyukarev-Qua\-dru\-pel bezüglich $\sjk$. Dann gelten
  \begin{align*}
    \Aarn(z)\beklam{\dPs_{\arn}(\alpha)}^{\ast} & = \beklam{\dPs_{\arn}(\za)}^{\ast}, \\
    \Carn(z)\beklam{\dPs_{\arn}(\alpha)}^{\ast} & =  -(z-\alpha)\beklam{P_{\arn}(\za)}^{\ast}
  \end{align*}
  für alle $z \in \C$ und $n \in \Zofk$ sowie
  \begin{align*}
    \Barn(z)\beklam{P_n(\alpha)}^{\ast} & = -\beklam{\Ps_n(\za)}^{\ast}, \\
    \Darn(z)\beklam{P_n(\alpha)}^{\ast} & = \beklam{P_n(\za)}^{\ast}
  \end{align*}
  für alle $z \in \C$ und $n \in \Zofkp$.
\end{bem}

\bwanf Sei $z \in \C$. Wegen \thref{drdef1} und \thref{sqbm1} gelten dann
\begin{align*}
  \Aaro(z)\beklam{\dPs_{\aro}(\alpha)}^{\ast} & = s_0 = \beklam{\dPs_{\aro}(\za)}^{\ast}, \\
  \Baro(z)\beklam{P_0(\alpha)}^{\ast} & = \Oq = -\beklam{\Ps_0(\za)}^{\ast}, \\
  \Caro(z)\beklam{\dPs_{\aro}(\alpha)}^{\ast} & = -(z-\alpha)\Iq = -(z-\alpha)\beklam{P_{\aro}(\za)}^{\ast}, \\
  \Daro(z)\beklam{P_0(\alpha)}^{\ast} & = \Iq = \beklam{P_0(\za)}^{\ast}.
\end{align*}
Sei nun $n \in \Zefkp$. Sei $\tHn$ definiert wie in Teil (f) von \thref{drbm1}. Wegen Teil (d) von \thref{drlm2} und Teil (f) von \thref{drbm1} gilt dann
\begin{align}	\label{sqbm2bw4}
  y_{0,n-1}v^{\ast}_n\Rna(\alpha) + \Harnm R^{\ast}_{n-1}(\alpha)L^{\ast}_n = \dL^{\ast}_n\tHn.
\end{align}
Weiterhin gilt
\begin{align}	\label{sqbm2bw5}
  \dL^{\ast}_n\tHn \begin{pmatrix} -\Hnm^{-1}y_{n,2n-1} \\ \Iq \end{pmatrix} 
  = \begin{pmatrix} \Hnm & y_{n,2n-1} \end{pmatrix} \begin{pmatrix} -\Hnm^{-1}y_{n,2n-1} \\ \Iq \end{pmatrix} = \Onq.
\end{align}
Wegen \thref{drdef1}, \thref{sqbm1}, der Teile (e) und (b) von \thref{drlm3}, Teil (d) von \thref{drlm1}, Teil (c) von \thref{drlm1}, \fref{sqbm2bw4} und \fref{sqbm2bw5} gilt nun
\begin{align*}
  &\ \Barn(z)\beklam{P_n(\alpha)}^{\ast} + \beklam{\Ps_n(\za)}^{\ast} \\
  & = \uarnm^{\ast}\Rnma(\za)\Harnm^{-1}y_{0,n-1} \eklam{\begin{pmatrix} -z_{n,2n-1}\Hnm^{-1} & \Iq \end{pmatrix} \Rn(\alpha)v_n}^{\ast} \\
  &\quad +  \eklam{\begin{pmatrix} -z_{n,2n-1}\Hnm^{-1} & \Iq \end{pmatrix} \Rn(\za)u_n}^{\ast} \\
  & = \eklam{\uarnm^{\ast}\Rnma(\za)\Harnm^{-1}y_{0,n-1}v^{\ast}_n\Rna(\alpha) + u^{\ast}_n\Rna(\za)} \begin{pmatrix} -\Hnm^{-1}y_{n,2n-1} \\ \Iq \end{pmatrix} \\
  & = \eklam{y^{\ast}_{0,n-1}\Rnma(\za)\Rnm^{-\ast}(\alpha)\Harnm^{-1}y_{0,n-1}v^{\ast}_n\Rna(\alpha) + y^{\ast}_{0,n-1}\Rnma(\za)L^{\ast}_n} \begin{pmatrix} -\Hnm^{-1}y_{n,2n-1} \\ \Iq \end{pmatrix} \\
  & = y^{\ast}_{0,n-1}\Rnma(\za)\Rnm^{-\ast}(\alpha)\Harnm^{-1} \eklam{y_{0,n-1}v^{\ast}_n\Rna(\alpha) + \Harnm\Rnma(\alpha)L^{\ast}_n} \begin{pmatrix} -\Hnm^{-1}y_{n,2n-1} \\ \Iq \end{pmatrix} \\
  & = y^{\ast}_{0,n-1}\Rnma(\za)\Rnm^{-\ast}(\alpha)\Harnm^{-1} \dL^{\ast}_n\tHn \begin{pmatrix} -\Hnm^{-1}y_{n,2n-1} \\ \Iq \end{pmatrix} = \Oq.
\end{align*}
Wegen der Teile (d) und (a) von \thref{drlm1} sowie der Teile (b), (a) und (d) von \thref{drlm2} gilt weiterhin
\begin{align}	\label{sqbm2bw6}
  -\vna\eklam{\Rna(\alpha)-\Rna(\za)} & = -(\alpha-z)\vna\Rna(\za)T^{\ast}_n\Rna(\alpha) \notag \\
  & = (z-\alpha)\vnma\dL^{\ast}_nT^{\ast}_n\Rna(\za)\Rna(\alpha) \notag \\
  & = (z-\alpha)\vnma L^{\ast}_n\Rna(\za)\Rna(\alpha) \notag \\
  & = (z-\alpha)\vnma\Rnma(\za)\Rnma(\alpha)L^{\ast}_n.
\end{align}
Wegen \thref{drdef1}, \thref{sqbm1}, \fref{sqbm2bw6}, \fref{sqbm2bw4} und \fref{sqbm2bw5} gilt nun
\begin{align*}
  &\ \Darn(z)\beklam{P_n(\alpha)}^{\ast} - \beklam{P_n(\za)}^{\ast} \\
  & = \eklam{\Iq-(z-\alpha)\vnma\Rnma(\za)\Harnm^{-1}y_{0,n-1}} \eklam{\begin{pmatrix} -z_{n,2n-1}\Hnm^{-1} & \Iq \end{pmatrix} \Rn(\alpha)v_n}^{\ast} \\
  &\quad - \eklam{\begin{pmatrix} -z_{n,2n-1}\Hnm^{-1} & \Iq \end{pmatrix} \Rn(\za)v_n}^{\ast} \\
  & = \eklam{\vna\eklam{\Rna(\alpha)-\Rna(\za)} - (z-\alpha)\vnma\Rnma(\za)\Harnm^{-1}y_{0,n-1}\vna\Rna(\alpha)} \begin{pmatrix} -\Hnm^{-1}y_{n,2n-1} \\ \Iq \end{pmatrix} \\
  & = \eklam{-(z-\alpha)\vnma\Rnma(\za)\Rnma(\alpha)L^{\ast}_n - (z-\alpha)\vnma\Rnma(\za)\Harnm^{-1}y_{0,n-1}\vna\Rna(\alpha)} \\
  &\quad \cdot \begin{pmatrix} -\Hnm^{-1}y_{n,2n-1} \\ \Iq \end{pmatrix} \\
  & = -(z-\alpha)\vnma\Rnma(\za)\Harnm^{-1} \eklam{\Harnm\Rnma(\alpha)L^{\ast}_n + y_{0,n-1}\vna\Rna(\alpha)} \begin{pmatrix} -\Hnm^{-1}y_{n,2n-1} \\ \Iq \end{pmatrix} \\
  & = -(z-\alpha)\vnma\Rnma(\za)\Harnm^{-1} \dL^{\ast}_n\tHn \begin{pmatrix} -\Hnm^{-1}y_{n,2n-1} \\ \Iq \end{pmatrix} = \Oq.
\end{align*}
Sei nun $\kappa\geq2$ und $n \in \Zefk$. Wegen Teil (d) von \thref{drlm1} sowie der Teile (e) und (a) von \thref{drlm3} gilt
\begin{align}	\label{sqbm2bw1}
  \uarn^{\ast}\eklam{\Rna(\za)-\Rna(\alpha)} & = (z-\alpha)\uarn^{\ast}\Rna(\alpha)T^{\ast}_n\Rna(\za) \notag \\
  & = (z-\alpha)y^{\ast}_{0,n}T^{\ast}_n\Rna(\za) \notag \\
  & = (z-\alpha)u^{\ast}_n\Rna(\za).
\end{align}
Sei $\tHarn$ definiert wie in Teil (e) von \thref{drbm1}. Wegen Teil (e) von \thref{drbm1} gilt dann
\begin{align}	\label{sqbm2bw2}
  \Rn(\alpha)v_ny^{\ast}_{0,n} - H_n = \Rn(\alpha)\beklam{\Rn^{-1}(\alpha)H_n-T_n\tHarn}-H_n = -\Rn(\alpha)T_n\tHarn.
\end{align}
Weiterhin gilt
\begin{align}	\label{sqbm2bw3}
  T_n\tHarn \begin{pmatrix} -\Harnm^{-1}y_{\arn,2n-1} \\ \Iq \end{pmatrix} 
  = \begin{pmatrix} \Oqn & \Oq \\ \Harnm & y_{\arn,2n-1} \end{pmatrix} \begin{pmatrix} -\Harnm^{-1}y_{\arn,2n-1} \\ \Iq \end{pmatrix} = 0_{(n+1)q\times q}.
\end{align}
Wegen \thref{drdef1}, \thref{sqbm1}, \fref{sqbm2bw1}, Teil (e) von \thref{drlm3}, \fref{sqbm2bw2} und \fref{sqbm2bw3} gilt nun
\begin{align*}
  &\ \Aarn(z)\beklam{\dPs_{\arn}(\alpha)}^{\ast} - \beklam{\dPs_{\arn}(\za)}^{\ast} \\
  & = \eklam{\Iq + (z-\alpha)u^{\ast}_n\Rna(\za)\Hn^{-1}\Rn(\alpha)v_n} \eklam{\begin{pmatrix} -z_{\arn,2n-1}\Harnm^{-1} & \Iq \end{pmatrix} \Rn(\alpha)\uarn}^{\ast} \\
  &\quad - \eklam{ \begin{pmatrix} -z_{\arn,2n-1}\Harnm^{-1} & \Iq \end{pmatrix} \Rn(\za)\uarn}^{\ast} \\
  & = \eklam{-\uarn^{\ast}\eklam{\Rna(\za)-\Rna(\alpha)} + (z-\alpha)u^{\ast}_n\Rna(\za)\Hn^{-1}\Rn(\alpha)v_n \uarn^{\ast}\Rna(\alpha)} \\
  &\quad \cdot \begin{pmatrix} -\Harnm^{-1}y_{\arn,2n-1} \\ \Iq \end{pmatrix} \\
  & = \eklam{-(z-\alpha)u^{\ast}_n\Rna(\za) + (z-\alpha)u^{\ast}_n\Rna(\za)\Hn^{-1}\Rn(\alpha)v_n y^{\ast}_{0,n}} \begin{pmatrix} -\Harnm^{-1}y_{\arn,2n-1} \\ \Iq \end{pmatrix} \\
  & = (z-\alpha)u^{\ast}_n\Rna(\za)\Hn^{-1}\eklam{-\Hn + \Rn(\alpha)v_ny^{\ast}_{0,n}} \begin{pmatrix} -\Harnm^{-1}y_{\arn,2n-1} \\ \Iq \end{pmatrix} \\
  & = -(z-\alpha)u^{\ast}_n\Rna(\za)\Hn^{-1}\Rn(\alpha)T_n\tHarn \begin{pmatrix} -\Harnm^{-1}y_{\arn,2n-1} \\ \Iq \end{pmatrix} = \Oq.
\end{align*}
Wegen \thref{drdef1}, \thref{sqbm1}, Teil (e) von \thref{drlm3}, \fref{sqbm2bw2} und \fref{sqbm2bw3} gilt weiterhin
\begin{align*}
  &\ \Carn(z)\beklam{\dPs_{\arn}(\alpha)}^{\ast} + (z-\alpha)\beklam{P_{\arn}(\za)}^{\ast} \\
  & = -(z-\alpha)v^{\ast}_n\Rna(\za)\Hn^{-1}\Rn(\alpha)v_n \eklam{\begin{pmatrix} -z_{\arn,2n-1}\Harnm^{-1} & \Iq \end{pmatrix} \Rn(\alpha)\uarn}^{\ast} \\
  &\quad + (z-\alpha) \eklam{ \begin{pmatrix} -z_{\arn,2n-1}\Harnm^{-1} & \Iq \end{pmatrix} \Rn(\za)v_n}^{\ast} \\
  & = \eklam{ -(z-\alpha)v^{\ast}_n\Rna(\za)\Hn^{-1}\Rn(\alpha)v_n \uarn^{\ast}\Rna(\alpha) + (z-\alpha)v^{\ast}_n\Rna(\za)} \\
  &\quad \cdot \begin{pmatrix} -\Harnm^{-1}y_{\arn,2n-1} \\ \Iq \end{pmatrix} \\
  & = -(z-\alpha)v^{\ast}_n\Rna(\za)\Hn^{-1} \rklam{\Rn(\alpha)v_ny^{\ast}_{0,n}-\Hn} \begin{pmatrix} -\Harnm^{-1}y_{\arn,2n-1} \\ \Iq \end{pmatrix} \\
  & = (z-\alpha)v^{\ast}_n\Rna(\za)\Hn^{-1} \Rn(\alpha)T_n\tHarn \begin{pmatrix} -\Harnm^{-1}y_{\arn,2n-1} \\ \Iq \end{pmatrix} = \Oq. \tag*{$\Box$}
\end{align*}

Folgendes Resultat liefert uns eine Darstellung der Matrixpolynome des rechtsseitigen $\alpha$-Stieltjes-Quadrupels bezüglich einer rechtsseitig $\alpha$-Stieltjes-positiv definiten Folge an der Stelle $\alpha$ mithilfe der rechtsseitigen $\alpha$-Dyukarev-Stieltjes-Parametrisierung jener Folge (vergleiche \cite[Proposition 4.5]{CR1} im Fall $\alpha=0$ und $\kappa=\infty$).

\begin{satz}	\thlabel{sqsa3}
	Seien $\kappa \in \Na$, $\alpha \in \R$ und $\sjk \in \Kpqka$. Weiterhin seien $\SQark$ das rechtsseitige $\alpha$-Stieltjes-Quadrupel bezüglich $\sjk$ und $\LMkarn$  die rechtsseitige $\alpha$-Dyukarev-Stieltjes-Parametrisierung \linebreak von $\sjk$. \dgfa
	\begin{itemize}
		\item [\rm{(a)}] Es sind $\brklam{P_n(\alpha)}^{\fklam{\kappa+1}}_{n=0}$, $\brklam{\Ps_n(\alpha)}^{\fklam{\kappa+1}}_{n=1}$, $\brklam{P_{\arn}(\alpha)}^{\fklam{\kappa}}_{n=0}$ und $\brklam{\dPs_{\arn}(\alpha)}^{\fklam{\kappa}}_{n=0}$ Folgen von regulären Matrizen aus $\Cqq$.
		\item [\rm{(b)}] Es gelten
		\begin{align*}
			P_n(\alpha) &= (-1)^n\prodr^{n-1}_{j=0}\rklam{\Marj^{-1}\Larj^{-1}}, \\
			\Ps_n(\alpha) &= (-1)^{n+1}\prodr^{n-1}_{j=0}\rklam{\Marj^{-1}\Larj^{-1}}\sum^{n-1}_{j=0}\Larj
		\end{align*}
		für alle $n\in\Zefkp$ und im Fall $\kappa\geq2$
		\begin{align*}
			P_{\arn}(\alpha) &= (-1)^n\prodr^{n-1}_{j=0}\rklam{\Marj^{-1}\Larj^{-1}}\Marn^{-1}\sum^n_{j=0}\Marj, \\
			\dPs_{\arn}(\alpha) &= (-1)^n\prodr^{n-1}_{j=0}\rklam{\Marj^{-1}\Larj^{-1}}\Marn^{-1}
		\end{align*}
		für alle $n\in\Zefk$.
	\end{itemize}
\end{satz}

\bwanf Zu (a): Unter Beachtung von \thref{adpbm1} sind wegen $P_0 \equiv \Iq$, $P_{\aro} \equiv \Iq$ und $\dPs_{\aro} \equiv s_0$ (vergleiche \thref{sqbm1}) die Matrizen $P_0(\alpha)$, $P_{\aro}(\alpha)$ und $\dPs_{\aro}(\alpha)$ regulär. Sei $\ABCDark$ das rechtsseitige $\alpha$-Dyukarev-Qua\-dru\-pel bezüglich $\sjk$. Wegen \thref{sqbm2} gelten dann
\begin{align*}
    \Carn(z)\beklam{\dPs_{\arn}(\alpha)}^{\ast} & =  -(z-\alpha)\beklam{P_{\arn}(\za)}^{\ast}
\end{align*}
für alle $z \in \C$ und $n \in \Zofk$ sowie
\begin{align*}
    \Barn(z)\beklam{P_n(\alpha)}^{\ast} & = -\beklam{\Ps_n(\za)}^{\ast}, \\
    \Darn(z)\beklam{P_n(\alpha)}^{\ast} & = \beklam{P_n(\za)}^{\ast}
\end{align*}
für alle $z \in \C$ und $n \in \Zofkp$. 
Hieraus folgen wegen \thref{drsa7} nun 
\begin{align}
	\beklam{P_n(\za)}^{\ast} &= \eklam{(-1)^n(z-\alpha)^n\prodr^{n-1}_{j=0}\rklam{\Marj\Larj} + \ldots + \Iq} \beklam{P_n(\alpha)}^{\ast}, \label{sqsa3bw1}\\
	\beklam{\Ps_n(\za)}^{\ast} &= \eklam{(-1)^n(z-\alpha)^{n-1}\Laro\prodr^{n-1}_{j=1}\rklam{\Marj\Larj} + \ldots - \sum^{n-1}_{j=0}\Larj} \beklam{P_n(\alpha)}^{\ast} \label{sqsa3bw2}
\end{align}
für alle $z \in \C$ und $n \in \Zefkp$ sowie im Fall $\kappa\geq2$
\begin{align} \label{sqsa3bw3}
	&\ (z-\alpha)\beklam{P_{\arn}(\za)}^{\ast} \notag \\
	&= \eklam{(-1)^n(z-\alpha)^{n+1}\Maro\prodr^{n-1}_{j=0}\rklam{\Larj\Marjp} + \ldots + (z-\alpha)\sum^n_{j=0}\Marj} \beklam{\dPs_{\arn}(\alpha)}^{\ast}
\end{align}
für alle $z \in \C$ und $n \in \Zefk$.
Wegen \thref{sqdef1} und \thref{mpdef2} ist der Leitkoeffizient von $P_n$ für alle $n \in \Zofkp$ und $P_{\arn}$ für alle $n \in \Zofk$ gleich $\Iq$. Hieraus folgt durch Koeffizientenvergleich in \fref{sqsa3bw1} bzw. \fref{sqsa3bw3} dann
\begin{align} \label{sqsa3bw4}
	\Iq = (-1)^n\prodr^{n-1}_{j=0}\rklam{\Marj\Larj} \beklam{P_n(\alpha)}^{\ast}
\end{align} 
für alle $n \in \Zefkp$ bzw. im Fall $\kappa\geq2$
\begin{align} \label{sqsa3bw5}
	\Iq = (-1)^n\Maro\prodr^{n-1}_{j=0}\rklam{\Larj\Marjp} \beklam{\dPs_{\arn}(\alpha)}^{\ast}
\end{align}
für alle $n \in \Zefk$.
Wegen \thref{adpbm3} sind $(\Larn)^{\fklam{\kappa-1}}_{n=0}$ und $(\Marn)^{\fklam{\kappa}}_{n=0}$ Folgen von regulären Matrizen. 
Hieraus folgt wegen \fref{sqsa3bw4} bzw. \fref{sqsa3bw5} dann, dass auch $P_n(\alpha)$ für alle $n \in \Zefkp$ bzw. im Fall $\kappa\geq2$ auch $\dPs_{\arn}(\alpha)$ für alle $n \in \Zefk$ reguläre Matrizen sind.
Wegen \fref{sqsa3bw2} gilt
\begin{align} \label{sqsa3bw6}
	\beklam{\Ps_n(\alpha)}^{\ast} =  - \rklam{\sum^{n-1}_{j=0}\Larj} \beklam{P_n(\alpha)}^{\ast}
\end{align}
für alle $n \in \Zefkp$. Unter Beachtung von \thref{adpbm3} ist $\sum^{n-1}_{j=0}\Larj$ für alle \linebreak $n \in \Zefkp$ eine positiv hermitesche und insbesondere reguläre Matrix, also ist wegen \fref{sqsa3bw6} auch $\Ps_n(\alpha)$ für alle $n \in \Zefkp$ eine reguläre Matrix.

Seien nun $\kappa\geq2$ und $\tHarn$ für alle $n \in \Zefk$ definiert wie in Teil (e) von \thref{drbm1}. Wegen Teil (e) von \thref{drbm1} gilt dann
\begin{align}	\label{sqsa3bw8}
  \Rn(\alpha)v_ny^{\ast}_{0,n} - H_n = \Rn(\alpha)\beklam{\Rn^{-1}(\alpha)H_n-T_n\tHarn}-H_n = -\Rn(\alpha)T_n\tHarn
\end{align}
für alle $n \in \Zefk$. Weiterhin gilt
\begin{align}	\label{sqsa3bw9}
  T_n\tHarn \begin{pmatrix} -\Harnm^{-1}y_{\arn,2n-1} \\ \Iq \end{pmatrix} 
  = \begin{pmatrix} \Oqn & \Oq \\ \Harnm & y_{\arn,2n-1} \end{pmatrix} \begin{pmatrix} -\Harnm^{-1}y_{\arn,2n-1} \\ \Iq \end{pmatrix} = 0_{(n+1)q\times q}
\end{align}
für alle $n \in \Zefk$.
Wegen \thref{adpdef1} und Teil (a) von \thref{sqlm1} gilt
\begin{align*}
	\sum^n_{j=0}\Marj &= \sum^n_{j=1}\eklam{\Ej^{\ast}(\alpha)\Hj^{-1}\Ej(\alpha)-\Ejm^{\ast}(\alpha)\Hjm^{-1}\Ejm(\alpha)}+s^{-1}_0 \\
	&= \En^{\ast}(\alpha)\Hn^{-1}\En(\alpha) \\
	&= \vna\Rna(\alpha)\Hn^{-1}\Rn(\alpha)\vn
\end{align*}
für alle $n \in \Zefk$.
Hieraus folgt wegen \thref{sqbm1}, Teil (e) von \thref{drlm3}, \fref{sqsa3bw8} und \fref{sqsa3bw9} dann
\begin{align*}
	& \rklam{\sum^n_{j=0}\Marj} \beklam{\dPs_{\arn}(\alpha)}^{\ast}-\beklam{P_{\arn}(\alpha)}^{\ast} \\
	&= \vna\Rna(\alpha)\Hn^{-1}\Rn(\alpha)\vn \uarn^{\ast}\Rna(\alpha)\begin{pmatrix} -\Harnm^{-1}y_{\arn,2n-1} \\ \Iq \end{pmatrix} \\
	&\quad -\vna\Rna(\alpha)\begin{pmatrix} -\Harnm^{-1}y_{\arn,2n-1} \\ \Iq \end{pmatrix} \\
	&= \vna\Rna(\alpha)\Hn^{-1}\eklam{\Rn(\alpha)\vn\yna - \Hn}\begin{pmatrix} -\Harnm^{-1}y_{\arn,2n-1} \\ \Iq \end{pmatrix} \\
	&= -\vna\Rna(\alpha)\Hn^{-1}\Rn(\alpha)T_n\tHarn\begin{pmatrix} -\Harnm^{-1}y_{\arn,2n-1} \\ \Iq \end{pmatrix} = \Oq
\end{align*}
für alle $n \in \Zefk$.
Hieraus folgt nun
\begin{align} \label{sqsa3bw7}
	\beklam{P_{\arn}(\alpha)}^{\ast} = \rklam{\sum^n_{j=0}\Marj} \beklam{\dPs_{\arn}(\alpha)}^{\ast}
\end{align}
für alle $n \in \Zefk$. Unter Beachtung von \thref{adpbm3} ist $\sum^n_{j=0}\Marj$ für alle $n \in \Zefk$ eine positiv hermitesche und insbesondere reguläre Matrix, also ist wegen \fref{sqsa3bw7} auch $P_{\arn}(\alpha)$ für alle $n \in \Zefk$ eine reguläre Matrix.

Zu (b): Dies folgt wegen \thref{adpbm3} aus \fref{sqsa3bw4}, \fref{sqsa3bw5}, \fref{sqsa3bw6} und \fref{sqsa3bw7}. \bwend

Wir können nun die Aussage von \thref{sqbm2} konkretisieren und die rechtsseitigen $2$\textit{q}$\times2$\textit{q}-$\alpha$-Dyukarev-Matrixpolynome bezüglich einer rechtsseitig $\alpha$-Stieltjes-positiv definiten Folge mithilfe des rechtsseitigen $\alpha$-Stieltjes-Quadrupels bezüglich jener Folge darstellen (vergleiche \cite[Theorem 4.6]{CR1}, \cite[Theorem 4.7]{CR1} und \cite[Theorem 4.8]{CR1} im Fall $\alpha=0$ und $\kappa=\infty$).

\begin{satz}	\thlabel{sqsa4}
	Seien $\kappa \in \Na$, $\alpha \in \R$, $\sjk \in \Kpqka$ und $\SQark$ das rechtsseitige $\alpha$-Stieltjes-Quadrupel bezüglich $\sjk$. \dgfa
	\begin{itemize}
		\item [\rm{(a)}] Sei $\ABCDark$ das rechtsseitige $\alpha$-Dyukarev-Qua\-dru\-pel bezüglich $\sjk$. Dann sind $P_n(\alpha)$ für alle $n \in \Zofkp$ sowie $\dPs_{\arn}(\alpha)$ für alle $n \in \Zofk$ reguläre Matrizen und es gelten
		\begin{align*}
			\Aarn(z) & = \beklam{\dPs_{\arn}(\za)}^{\ast}\beklam{\dPs_{\arn}(\alpha)}^{-\ast}, \\
    		\Carn(z) & =  -(z-\alpha)\beklam{P_{\arn}(\za)}^{\ast}\beklam{\dPs_{\arn}(\alpha)}^{-\ast}
  		\end{align*}
  		für alle $z \in \C$ und $n \in \Zofk$ sowie
  		\begin{align*}
   			\Barn(z) & = -\beklam{\Ps_n(\za)}^{\ast}\beklam{P_n(\alpha)}^{-\ast}, \\
    		\Darn(z) & = \beklam{P_n(\za)}^{\ast}\beklam{P_n(\alpha)}^{-\ast}
  		\end{align*}
  		für alle $z \in \C$ und $n \in \Zofkp$.
  		\item [\rm{(b)}] Sei $\Uarmk$ die Folge von rechtsseitigen $2$\textit{q}$\times2$\textit{q}-$\alpha$-Dyukarev-Matrixpolynomen bezüglich $\sjk$. Weiterhin seien $z \in \C$ und $m \in \Zok$. Dann gilt
  		\begin{align}	\label{sqsa4bw1}
  			\Uarm(z) = \begin{pmatrix} \beklam{\dPs_{\ar\fklam{m}}(\za)}^{\ast}\beklam{\dPs_{\ar\fklam{m}}(\alpha)}^{-\ast} & -\beklam{\Ps_{\fklam{m+1}}(\za)}^{\ast}\beklam{P_{\fklam{m+1}}(\alpha)}^{-\ast} \\ -(z-\alpha)\beklam{P_{\ar\fklam{m}}(\za)}^{\ast}\beklam{\dPs_{\ar\fklam{m}}(\alpha)}^{-\ast} & \beklam{P_{\fklam{m+1}}(\za)}^{\ast}\beklam{P_{\fklam{m+1}}(\alpha)}^{-\ast} \end{pmatrix}.
  		\end{align}
  		Außerdem ist $\Uarm(z)$ eine reguläre Matrix und es gilt
  		\begin{align*}
  			\Uarm^{-1}(z) = \begin{pmatrix} \beklam{P_{\fklam{m+1}}(\alpha)}^{-1}P_{\fklam{m+1}}(z) & \beklam{P_{\fklam{m+1}}(\alpha)}^{-1}\Ps_{\fklam{m+1}}(z) \\ (z-\alpha)\beklam{\dPs_{\ar\fklam{m}}(\alpha)}^{-1}P_{\ar\fklam{m}}(z) & \beklam{\dPs_{\ar\fklam{m}}(\alpha)}^{-1}\dPs_{\ar\fklam{m}}(z) \end{pmatrix}.
  		\end{align*}
  		\item [\rm{(c)}] Seien $m\in\Zek$ und $\Sarmmin$ bzw. $\Sarmmax$ das untere bzw. obere Extremalelement von $\Sqasmu$. Dann gelten
  		\begin{align*}  			
  			\Sarmmin(z) &= -\beklam{\Ps_{\fklam{m+1}}(\za)}^{\ast}\beklam{P_{\fklam{m+1}}(\za)}^{-\ast}, \\
  			\Sarmmax(z) &= -\frac{1}{z-\alpha}\beklam{\dPs_{\ar\fklam{m}}(\za)}^{\ast}\beklam{P_{\ar\fklam{m}}(\za)}^{-\ast}
  		\end{align*}
  		für alle $z \in \C\setminus[\alpha,\infty)$.  		
  	\end{itemize}		
\end{satz}

\bwanf Zu (a): Dies folgt wegen Teil (a) von \thref{sqsa3} aus \thref{sqbm2}.

Zu (b): \fref{sqsa4bw1} folgt wegen (a) aus \thref{drdef2}. Wegen \thref{spbsp1}, Teil (c) von \thref{drbm5}, Teil (a) von \thref{splm1} und \fref{sqsa4bw1} ist $\Uarm(z)$ eine reguläre Matrix und es gilt
\begin{align*}
	&\ \Uarm^{-1}(z) = \tJq\beklam{\Uarm(\za)}^{\ast}\tJq \\
	&= \begin{pmatrix} \Oq & -i\Iq \\ i\Iq & \Oq \end{pmatrix} 
	\begin{pmatrix} \beklam{\dPs_{\ar\fklam{m}}(\alpha)}^{-1}\dPs_{\ar\fklam{m}}(z) & -(z-\alpha)\beklam{\dPs_{\ar\fklam{m}}(\alpha)}^{-1}P_{\ar\fklam{m}}(z) \\ -\beklam{P_{\fklam{m+1}}(\alpha)}^{-1}\Ps_{\fklam{m+1}}(z) & \beklam{P_{\fklam{m+1}}(\alpha)}^{-1}P_{\fklam{m+1}}(z) \end{pmatrix} \\
	&\quad \cdot\begin{pmatrix} \Oq & -i\Iq \\ i\Iq & \Oq \end{pmatrix} \\
	&= \begin{pmatrix} i\beklam{P_{\fklam{m+1}}(\alpha)}^{-1}\Ps_{\fklam{m+1}}(z) & -i\beklam{P_{\fklam{m+1}}(\alpha)}^{-1}P_{\fklam{m+1}}(z) \\ i\beklam{\dPs_{\ar\fklam{m}}(\alpha)}^{-1}\dPs_{\ar\fklam{m}}(z) & -i(z-\alpha)\beklam{\dPs_{\ar\fklam{m}}(\alpha)}^{-1}P_{\ar\fklam{m}}(z) \end{pmatrix} \begin{pmatrix} \Oq & -i\Iq \\ i\Iq & \Oq \end{pmatrix} \\
	&= \begin{pmatrix} \beklam{P_{\fklam{m+1}}(\alpha)}^{-1}P_{\fklam{m+1}}(z) & \beklam{P_{\fklam{m+1}}(\alpha)}^{-1}\Ps_{\fklam{m+1}}(z) \\ (z-\alpha)\beklam{\dPs_{\ar\fklam{m}}(\alpha)}^{-1}P_{\ar\fklam{m}}(z) & \beklam{\dPs_{\ar\fklam{m}}(\alpha)}^{-1}\dPs_{\ar\fklam{m}}(z) \end{pmatrix}.
\end{align*}

Zu (c): Wegen \thref{drdef4} und (a) gelten
\begin{align*}
	\Sarmmin(z) &= \Barfm(z)\Darfm^{-1}(z) \\
	&= -\beklam{\Ps_{\fklam{m+1}}(\za)}^{\ast}\beklam{P_{\fklam{m+1}}(\alpha)}^{-\ast}\beklam{P_{\fklam{m+1}}(\alpha)}^{\ast}\beklam{P_{\fklam{m+1}}(\za)}^{-\ast} \\
	&= -\beklam{\Ps_{\fklam{m+1}}(\za)}^{\ast}\beklam{P_{\fklam{m+1}}(\za)}^{-\ast}	
\end{align*}
und
\begin{align*}
	\Sarmmax(z) &= \Aarfm(z)\Carfm^{-1}(z) \\
	&= \beklam{\dPs_{\ar\fklam{m}}(\za)}^{\ast}\beklam{\dPs_{\ar\fklam{m}}(\alpha)}^{-\ast}\frac{-1}{z-\alpha}\beklam{\dPs_{\ar\fklam{m}}(\alpha)}^{\ast}\beklam{P_{\ar\fklam{m}}(\za)}^{-\ast} \\
	&= -\frac{1}{z-\alpha}\beklam{\dPs_{\ar\fklam{m}}(\za)}^{\ast}\beklam{P_{\ar\fklam{m}}(\za)}^{-\ast}
\end{align*}
für alle $z \in \C\setminus[\alpha,\infty)$. \bwend

Im folgenden Satz betrachten wir die Lokalisierung der  Nullstellen der Determinanten der einzelnen \textit{q}$\times$\textit{q}"=Matrixpolynome des rechtsseitigen $\alpha$-Stieltjes-Quadrupels bezüglich einer rechtsseitig $\alpha$-Stieltjes-positiv definiten Folge.

\begin{satz}	\thlabel{sqsa7}
	Seien $\kappa \in \Na$, $\alpha \in \R$ und $\sjk \in \Kpqka$. Weiterhin seien $\SQark$ das rechtsseitige $\alpha$-Stieltjes-Quadrupel bezüglich $\sjk$ und $z\in\C\setminus(\alpha,\infty)$. Dann gelten $\det P_n(z) \neq 0$ für alle $n\in\Zofkp$, $\det \Ps_n(z) \neq 0$ für alle $n\in\Zefkp$, $\det P_{\arn}(z) \neq 0$ für alle $n\in\Zofk$ und $\det\dPs_{\arn}(z) \neq 0$ für alle $n\in\Zofk$.	
\end{satz}

\bwanf Dies folgt aus Teil (a) von \thref{sqsa3} und Teil (b) von \thref{drsa12} in Verbindung mit Teil (a) von \thref{sqsa4}. \bwend

Im Fall $q=1$ ist jener Teil von \thref{sqsa7}, welcher die Polynomenfolgen $(P_n)^{\fklam{\kappa+1}}$ und $(P_{\arn})^{\fklam{\kappa}}$ betrifft, bereits ein klassisches Resultat der Theorie orthogonaler Polynome mit reellen Koeffizienten auf der reellen Achse (vergleiche z.\,B. \cite[Theorem 5.2]{Ch} oder \cite[Satz 2.24]{Jung}). Die dortigen Beweise werden aber mit völlig verschiedenen Methoden geführt und sind nicht auf den Matrixfall übertragbar.

Wir wollen nun die rechtsseitige $\alpha$-Stieltjes-Parametrisierung einer rechtsseitig $\alpha$"=Stieltjes"=positiv definiten Folge mithilfe des rechtsseitigen $\alpha$-Stieltjes-Quadrupels bezüglich jener Folge darstellen (vergleiche \cite[Proposition 4.9(c), (e)]{CR1} für den Fall $\alpha=0$ und $\kappa=\infty$). 

\begin{satz}	\thlabel{sqsa5}
	Seien $\kappa \in \Na$, $\alpha \in \R$ und $\sjk \in \Kpqka$. Weiterhin seien $\SQark$ das rechtsseitige $\alpha$-Stieltjes-Quadrupel bezüglich $\sjk$ und $\Qarjk$ die rechtsseitige $\alpha$-Stieltjes-Parametrisierung von $\sjk$. \dgfa
	\begin{itemize}
		\item [\rm{(a)}] Es gelten
		\begin{align*}
			Q_{\ar 2n} = P_n(\alpha)\beklam{\dPs_{\arn}(\alpha)}^{\ast}
		\end{align*}
		für alle $n\in\Zofk$ und
		\begin{align*}
			Q_{\ar 2n+1} = -\dPs_{\arn}(\alpha)\beklam{P_{n+1}(\alpha)}^{\ast}
		\end{align*}
		für alle $n\in\Zofkm$.
		\item [\rm{(b)}] Es gelten
		\begin{align*}
			P_n(\alpha) = (-1)^n\prodl^{n-1}_{j=0}\rklam{Q_{\ar2j+1}Q^{-1}_{\ar2j}}
		\end{align*}
		für alle $n\in\Zefkp$ und im Fall $\kappa\geq2$
		\begin{align*}
			\dPs_{\arn}(\alpha) = (-1)^nQ_{\ar2n}\prodl^{n-1}_{j=0}\rklam{Q^{-1}_{\ar2j+1}Q_{\ar2j}}
		\end{align*}
		für alle $n\in\Zefk$.
	\end{itemize}
\end{satz}

\bwanf Zu (a): Wegen \thref{sqbm1} und Teil (a) von \thref{aspdef2} gilt
\begin{align*}
	P_0(\alpha)\beklam{\dPs_{\aro}(\alpha)}^{\ast} = \Iq \cdot s^{\ast}_0 = \dH_0 = Q_{\aro}.
\end{align*}
Sei nun $\kappa\geq2$ und $n\in\Zefk$. Wegen \thref{sqsa2} und Teil (a) von \thref{sqsa4} gilt dann
\begin{align}	\label{sqsa5bw1}
	\beklam{P_{\arn}(\za)}^{\ast}\beklam{\dPs_{\arn}(\alpha)}^{-\ast} = \sum^{n}_{j=0} \beklam{P_j(\za)}^{\ast}\dH^{-1}_jP_j(\alpha)
\end{align}
für alle $z \in \C\setminus\gklam{\alpha}$. Wegen \thref{sqdef1} und \thref{mpdef2} ist der Leitkoeffizient von $P_n$ sowie $P_{\arn}$ gleich $\Iq$ und es gilt $\deg P_{j} = j$ für alle $j \in \Zofkp$. Hieraus folgt durch Koeffizientenvergleich in \fref{sqsa5bw1} dann
\begin{align*}
	\beklam{\dPs_{\arn}(\alpha)}^{-\ast} = \dH^{-1}_nP_n(\alpha).
\end{align*}
Wegen Teil (a) von \thref{aspdef2} folgt hieraus
\begin{align*}
	Q_{\ar 2n} = \dH_n = P_n(\alpha)\beklam{\dPs_{\arn}(\alpha)}^{\ast}.
\end{align*}
Wegen \thref{sqsa2} und Teil (a) von \thref{sqsa4} gilt
\begin{align}	\label{sqsa5bw2}
	\beklam{P_n(\za)}^{\ast}\beklam{P_n(\alpha)}^{-\ast} = \Iq - (z-\alpha) \sum^{n-1}_{j=0} \beklam{P_{\arj}(\za)}^{\ast}\dH^{-1}_{\arj}\dPs_{\arj}(\alpha)
\end{align}
für alle $z \in \C$ und $n \in \Zefkp$. Wegen \thref{sqdef1} und \thref{mpdef2} ist der Leitkoeffizient von $P_n$ sowie $P_{\arn-1}$ gleich $\Iq$ und es gilt $\deg P_{\arj} = j$ für alle $j \in \Zofk$. Hieraus folgt durch Koeffizientenvergleich in \fref{sqsa5bw2} dann
\begin{align*}
	\beklam{P_n(\alpha)}^{-\ast} = -\dH^{-1}_{\arn-1}\dPs_{\arn-1}(\alpha)
\end{align*}
für alle $n \in \Zefkp$.
Wegen Teil (a) von \thref{aspdef2} folgt hieraus
\begin{align*}
	Q_{\ar 2n+1} = \dH_{\arn} = -\dPs_{\arn}(\alpha)\beklam{P_{n+1}(\alpha)}^{\ast}
\end{align*}
für alle $n \in \Zofkm$.

Zu (b): Wegen \thref{sqbm1} gilt $P_0(\alpha) = \Iq$. Unter Beachtung von Teil (c) von \thref{aspsa1} und Teil (a) von \thref{sqsa3} gilt wegen (a) nun
\begin{align*}
	\beklam{\dPs_{\arn}(\alpha)}^{-\ast} = Q^{-1}_{\ar2n}P_n(\alpha)
\end{align*}
für alle $n\in\Zofk$ und
\begin{align*}
	P_{n+1}(\alpha) = -Q_{\ar2n+1}\beklam{\dPs_{\arn}(\alpha)}^{-\ast}
\end{align*}
für alle $n\in\Zofkm$. Hieraus folgt
\begin{align*}
	P_{n+1}(\alpha) = -Q_{\ar2n+1}Q^{-1}_{\ar2n}P_n(\alpha)
\end{align*}
für alle $n\in\Zofkm$. Iterativ folgt dann
\begin{align*}
	P_n(\alpha) = (-1)^n\prodl^{n-1}_{j=0}\rklam{Q_{\ar2j+1}Q^{-1}_{\ar2j}}
\end{align*}
für alle $n\in\Zefkp$. Wegen \thref{sqbm1} und Teil (a) von \thref{aspdef2} gilt
\begin{align*}
	\dPs_{\aro}(\alpha) = s_0 = \dH_0 = Q_{\aro}.
\end{align*}
Unter Beachtung von Teil (c) von \thref{aspsa1} und Teil (a) von \thref{sqsa3} gilt wegen (a) nun
\begin{align*}
	\dPs_{\arn}(\alpha) = Q_{\ar2n}\beklam{P_n(\alpha)}^{-\ast}
\end{align*}
für alle $n\in\Zofk$ und
\begin{align*}
	\beklam{P_{n+1}(\alpha)}^{-\ast} = -Q^{-1}_{\ar2n+1}\dPs_{\arn}(\alpha)
\end{align*}
für alle $n\in\Zofkm$. Hieraus folgt im Fall $\kappa\geq2$ nun
\begin{align*}
	\dPs_{\arn}(\alpha) = -Q_{\ar2n}Q^{-1}_{\ar2n-1}\dPs_{\arn-1}(\alpha)
\end{align*}
für alle $n\in\Zefk$. Iterativ folgt dann
\begin{align*}
	\dPs_{\arn}(\alpha) = (-1)^nQ_{\ar2n}\prodl^{n-1}_{j=0}\rklam{Q^{-1}_{\ar2j+1}Q_{\ar2j}}
\end{align*}
für alle $n\in\Zefk$. \bwend

Teil (b) von \thref{sqsa5} wurde auf andere Weise auch schon in \thref{adplm1} bewiesen. Umgekehrt finden wir einen alternativen Beweis von \thref{adpsa2} durch Verwendung von Teil (a) von \thref{sqsa5} in Verbindung mit Teil (b) von \thref{sqsa3}.

Abschließend zeigen wir noch einen Zusammenhang zwischen rechtsseitigem $\alpha$-Stieltjes-Quadrupel bezüglich einer rechtsseitig $\alpha$-Stieltjes-positiv definiten Folge und den Favard"=Paaren bezüglich jener Folge und derer durch rechtsseitige $\alpha$-Verschiebung generierten Folge (vergleiche \cite[Proposition 9.1]{CR1} für den Fall $\alpha=0$ und $\kappa=\infty$).

\begin{satz}	\thlabel{sqsa6}
	Seien $\kappa \in \Na$, $\alpha \in \R$ und $\sjk \in \Kpqka$. Weiterhin seien $\SQark$ das rechtsseitige $\alpha$-Stieltjes-Quadrupel bezüglich $\sjk$, $\ABfk$ das Favard-Paar bezüglich $\sjk$ und $(B_{\arj})^0_{j=0}$ bzw. im Fall \linebreak $\kappa\geq2$ $\ABarfk$ das Favard-Paar bezüglich $\sarjk$. \dgfa
	\begin{itemize}
		\item [\rm{(a)}] Es gelten $P_{0} \equiv \Iq$,
    	\begin{align*}
			P_{1}(z) = (z\Iq-A_0)P_{0}(z)
   		\end{align*}
    	für alle $z \in \C$ und im Fall $\kappa \geq 3$
    	\begin{align*}
			P_{n}(z) = (z\Iq-A_{n-1})P_{n-1}(z)-B^{\ast}_{n-1}P_{n-2}(z)
    	\end{align*}
   		für alle $z \in \C$ und $n \in \Zzfkp$.
   		\item [\rm{(b)}] Es gelten $\Ps_{0} \equiv 0_{\qq}$, $\Ps_{1} \equiv B_{0}$ und im Fall $\kappa \geq 3$
      	\begin{align*}
			\Ps_{n}(z) = (z\Iq-A_{n-1})\Ps_{n-1}(z) - B^{\ast}_{n-1}\Ps_{n-2}(z)
     	\end{align*}
      	für alle $z \in \C$ und $n \in \Zzfkp$.
      	\item [\rm{(c)}] Es gelten $P_{\aro} \equiv \Iq$, im Fall $\kappa \geq 2$
    	\begin{align*}
			P_{\ar1}(z) = (z\Iq-A_{\aro})P_{\aro}(z)
   		\end{align*}
    	für alle $z \in \C$ und im Fall $\kappa \geq 4$
    	\begin{align*}
			P_{\arn}(z) = (z\Iq-A_{\arn-1})P_{\arn-1}(z)-B^{\ast}_{\arn-1}P_{\arn-2}(z)
    	\end{align*}
   		für alle $z \in \C$ und $n \in \Zzfk$.
   		\item [\rm{(d)}] Es gelten $\dPs_{\aro} \equiv B_0$, im Fall $\kappa \geq 2$
    	\begin{align*}
			\dPs_{\ar1}(z) = B_{\aro}+(z\Iq-A_{\aro})\dPs_{\aro}(z)
   		\end{align*}
    	für alle $z \in \C$ und im Fall $\kappa \geq 4$
		\begin{align*}
			\dPs_{\arn}(z) = (z\Iq-A_{\arn-1})\dPs_{\arn-1}(z)-B^{\ast}_{\arn-1}\dPs_{\arn-2}(z)
    	\end{align*}
   		für alle $z \in \C$ und $n \in \Zzfk$.	
	\end{itemize}	
\end{satz}

\bwanf Wegen Teil (b) von \thref{asmdef2} und \thref{asmbm1} gelten \linebreak $\sjfkm\in\Hpqfkm$ und $(s_{\arj})^{2\fklam{\kappa-2}}_{j=0} \in \Hpqfkzm$, falls $\kappa\geq2$.

Zu (a) und (c): Dies folgt aus \thref{sqdef1} und \thref{mpfo2}.

Zu (b): Dies folgt aus \thref{sqdef1} und \thref{mpfo3}.

Zu (d): Sei $\Psarnfk$ das linke System von Matrixpolynomen zweiter Art bezüglich $\sarjk$. Wegen Teil (e) von \thref{sqsa1} und \thref{fpdef1} gilt dann
\begin{align}	\label{sqsa6bw1}
	\dPs_{\arn}(z) = \Pars_{\arn}(z) + P_{\arn}(z) B_0
\end{align}
für alle $z \in \C$ und $n \in \Zofk$. Wegen \thref{mpfo3} gelten weiterhin $\Pars_{\aro} \equiv 0_{\qq}$, im Fall $\kappa\geq2$ $\Pars_{\ar1} \equiv B_{\aro}$ und im Fall $\kappa \geq 4$
\begin{align*}
	\Pars_{\arn}(z) = (z\Iq-A_{\arn-1})\Pars_{\arn-1}(z) - B^{\ast}_{\arn-1}\Pars_{\arn-2}(z)
\end{align*}
für alle $z \in \C$ und $n \in \Zzfk$. Hieraus folgt wegen \fref{sqsa6bw1} und (c) dann die Behauptung. \bwend

\subsection{Der linksseitige Fall}

Unsere folgende Überlegung ist auf die Einführung eines Quadrupels von Folgen von \textit{q}$\times$\textit{q}-Matrixfunktionen ausgerichtet, welches der in Kapitel \ref{chapmp} eingeführten in unserem Fall mit einer linksseitig $\alpha$-Stieltjes-positiv definiten Folge und derer durch linksseitige $\alpha$-Verschiebung generierten Folge assoziierten monischen links-orthogonalen Systeme von Matrixpolynomen sowie der hiermit assoziierten linken Systeme von Matrixpolynomen zweiter Art zugrunde liegt. Zuvor wollen wir einige Resultate für besagte Systeme rekapitulieren und erweitern. Dafür benötigen wir zunächst einige weitere Bezeichnungen.

\begin{bez}	\thlabel{sqlbz1}
	Seien $\kappa \in \Na$, $\alpha \in \R$ und $\sjk$ eine Folge aus $\Cqq$. Dann bezeichnet für $n \in \Zokm$
	\begin{align*}
  		\Saln := \begin{cases} s_{\alo} & \text{falls } n=0 \\
  		\begin{pmatrix} s_{\alo} & \Oq & \ldots & \Oq \\ s_{\al 1} & s_{\alo} & \ddots & \vdots \\ \vdots & \ddots & \ddots & \Oq \\ s_{\aln} & \ldots & s_{\al 1} & s_{\alo} \end{pmatrix} & \text{falls } n>0 \end{cases}
	\end{align*}
	die untere n-te Blockdreiecksmatrix von $\saljk$. Weiterhin seien $\dSalo := s_0$ und für $n \in \Zek$
	\begin{align*}
  		\dSaln := \begin{pmatrix} -s_0 & \Oq & \ldots & \Oq \\ s_{\alo} & -s_0 & \ddots & \vdots \\ \vdots & \ddots & \ddots & \Oq \\ s_{\aln-1} & \ldots & s_{\alo} & -s_0 \end{pmatrix} 
  = \alpha \begin{pmatrix} \Oqn & \Oq \\ S_{n-1} & \Onq \end{pmatrix} - S_n.
	\end{align*}
\end{bez}

\begin{satz}	\thlabel{sqlsa1}
  Seien $\kappa \in \Na$, $\alpha \in \R$ und $\sjk \in \Lpqka$. \dgfa
  \begin{itemize}
    \item [\rm{(a)}] Sei $\Pnfk$ das monische links-orthogonale System von Matrixpolynomen bezüglich $\sjk$. Dann gelten $P_0 \equiv \Iq$ und
      \begin{align*}
	P_n(z) = \begin{pmatrix} -z_{n,2n-1}\Hnm^{-1} & \Iq \end{pmatrix}E_n(z)
      \end{align*}
      für alle $z \in \C$ und $n \in \Zefkp$.
    \item [\rm{(b)}] Sei $\Psnfk$ das linke System von Matrixpolynomen zweiter Art bezüglich \linebreak $\sjk$. Dann gelten $\Ps_0 \equiv \Oq$ und
      \begin{align*}
	\Ps_n(z) = \begin{pmatrix} -z_{n,2n-1}\Hnm^{-1} & \Iq \end{pmatrix} \begin{pmatrix} \Oqn \\ S_{n-1} \end{pmatrix} E_{n-1}(z)
      \end{align*}
      für alle $z \in \C$ und $n \in \Zefkp$.
    \item [\rm{(c)}] Sei $\Palnfk$ das monische links-orthogonale System von Matrixpolynomen bezüglich $\saljk$. Dann gelten $P_{\alo} \equiv \Iq$ und im Fall $\kappa\geq2$
      \begin{align*}
	P_{\aln}(z) = \begin{pmatrix} -z_{\aln,2n-1}\Halnm^{-1} & \Iq \end{pmatrix}E_n(z)
      \end{align*}
      für alle $z \in \C$ und $n \in \Zefk$.
    \item [\rm{(d)}] Sei $\Psalnfk$ das linke System von Matrixpolynomen zweiter Art bezüglich \linebreak $\saljk$. Dann gelten $\Pals_{\alo} \equiv \Oq$ und im Fall $\kappa\geq2$
      \begin{align*}
	\Pals_{\aln}(z) = \begin{pmatrix} -z_{\aln,2n-1}\Halnm^{-1} & \Iq \end{pmatrix} \begin{pmatrix} \Oqn \\ S_{\aln-1} \end{pmatrix} E_{n-1}(z)
      \end{align*}
      für alle $z \in \C$ und $n \in \Zefk$.
    \item [\rm{(e)}] Seien $\widehat{P}_{\aln} := (\alpha-z)P_{\aln}(z)$ für alle $z \in \C$ und $n \in \Zofk$.
      Weiterhin sei $\dPs_{\aln}$ das zu $\sjk$ gehörige Matrixpolynom bezüglich $\widehat{P}_{\aln}$ für alle $n \in \Zofk$. 
      Dann gelten $\dPs_{\alo} \equiv -s_0$ und im Fall $\kappa\geq2$
      \begin{align*}
	\dPs_{\aln}(z) = \begin{pmatrix} -z_{\aln,2n-1}\Halnm^{-1} & \Iq \end{pmatrix} \dSaln E_n(z)
      \end{align*}
      für alle $z \in \C$ und $n \in \Zefk$ sowie
      \begin{align*}
	\dPs_{\aln}(z) = \Pals_{\aln}(z) - P_{\aln}(z) s_0
      \end{align*}
      für alle $z \in \C$ und $n \in \Zofk$.
  \end{itemize}
\end{satz}

\bwanf Zu (a) - (d): Unter Beachtung von \thref{asmbm1} gilt wegen $\sjk \in \Lpqka$ und Teil (b) von \thref{asmdef3} $\sjfkm \in \Hpqfkm$ bzw. $(s_{\alj})^{2\fklam{\kappa-2}} \in \Hpqfkzm$, falls $\kappa \geq 2$.
Hieraus folgt wegen \thref{mpsa2} dann die Behauptung von (a) bzw. (c) sowie wegen \thref{mpbem1} die Behauptung von (b) bzw. (d).

Zu (e): Sei $n \in \Zofk$. Unter Beachtung von (c) gelten dann $\deg \dPs_{\aln} = n+1$ und
\begin{align}	\label{sqlsa1bw1}
  \dP^{[j]}_{\aln} = \begin{cases} \alpha P^{[0]}_{\aln} & \text{falls } j = 0 \\
  \alpha P^{[j]}_{\aln} - P^{[j-1]}_{\aln} & \text{falls } 1 \leq j \leq n \\
  -P^{[n]}_{\aln} & \text{falls } j = n+1 \end{cases}
\end{align}
für alle $j \in \Zonp$. Sei nun $z \in \C$. Wegen \fref{sqlsa1bw1} und (c) gelten dann
\begin{align*}
  \dPs_{\alo}(z) & = \begin{pmatrix} \dP^{[0]}_{\alo} & \dP^{[1]}_{\alo} \end{pmatrix} \begin{pmatrix} \Oq \\ S_0 \end{pmatrix} E_0(z) 
  = \begin{pmatrix} \alpha P^{[0]}_{\alo} & -P^{[0]}_{\alo} \end{pmatrix} \begin{pmatrix} \Oq \\ s_0 \end{pmatrix} \\
  & = \begin{pmatrix} \alpha \Iq & -\Iq \end{pmatrix} \begin{pmatrix} \Oq \\ s_0 \end{pmatrix} 
  = -s_0
\end{align*}
sowie im Fall $\kappa\geq2$
\begin{align*}
  \dPs_{\aln}(z) & = \begin{pmatrix} \dP^{[0]}_{\aln} & \ldots & \dP^{[n+1]}_{\aln} \end{pmatrix} \begin{pmatrix} 0_{q\times(n+1)q} \\ S_n \end{pmatrix} \En(z) \\ 
  & = \eklam{\alpha \begin{pmatrix} P^{[0]}_{\aln} & \ldots & P^{[n]}_{\aln} & \Oq \end{pmatrix} - \begin{pmatrix} \Oq & P^{[0]}_{\aln} & \ldots & P^{[n]}_{\aln} \end{pmatrix}} \begin{pmatrix} 0_{q\times(n+1)q} \\ S_n \end{pmatrix} \En(z) \\ 
  & = \eklam{\alpha \begin{pmatrix} P^{[0]}_{\aln} & \ldots & P^{[n]}_{\aln} \end{pmatrix} \begin{pmatrix} \Oqn & \Oq \\ S_{n-1} & \Onq \end{pmatrix} - \begin{pmatrix} P^{[0]}_{\aln} & \ldots & P^{[n]}_{\aln} \end{pmatrix} S_n} \En(z) \\
  & = \begin{pmatrix} P^{[0]}_{\aln} & \ldots & P^{[n]}_{\aln} \end{pmatrix} \dSaln \En(z)
\end{align*}
für alle $n \in \Zefk$. Wegen (c) und (d) gelten weiterhin
\begin{align*}
  \Pals_{\alo}(z) - P_{\alo}(z)s_0 = -s_0 = \dPs_{\alo}(z)
\end{align*}
sowie im Fall $\kappa\geq2$
\begin{align*}
  &\ \Pals_{\aln}(z) - P_{\aln}(z)s_0 \\
  & = \begin{pmatrix} -z_{\aln,2n-1}\Halnm^{-1} & \Iq \end{pmatrix} \begin{pmatrix} \Oqn \\ S_{\aln-1} \end{pmatrix} \Enm(z) - \begin{pmatrix} -z_{\aln,2n-1}\Halnm^{-1} & \Iq \end{pmatrix}\En(z)s_0 \\
  & = \begin{pmatrix} -z_{\aln,2n-1}\Halnm^{-1} & \Iq \end{pmatrix} \eklam{\begin{pmatrix} \Oqn & \Oq \\ S_{\aln-1} & \Onq \end{pmatrix} \En(z) - \diag(s_0,...,s_0)\En(z)} \\
  & = \begin{pmatrix} -z_{\aln,2n-1}\Halnm^{-1} & \Iq \end{pmatrix} \dSaln \En(z) = \dPs_{\aln}(z)
\end{align*}
für alle $n \in \Zefk$. \bwend

\begin{defi}	\thlabel{sqldef1}
  Seien $\kappa \in \Na$, $\alpha \in \R$ und $\sjk \in \Lpqka$. 
  Weiterhin seien $(P_{n,s})^{\fklam{\kappa+1}}_{n=0}$ das monische links-orthogonale System von Matrixpolynomen bezüglich $\sjk$, \linebreak
  $(\Ps_{n,s})^{\fklam{\kappa+1}}_{n=0}$ das linke System von Matrixpolynomen zweiter Art bezüglich $\sjk$, \linebreak
  $(P_{\aln,s})^{\fklam{\kappa}}_{n=0}$ das monische links-orthogonale System von Matrixpolynomen bezüglich \linebreak $\saljk$,
  $\widehat{P}_{\aln,s} := (\alpha-z)P_{\aln,s}(z)$ für alle $z \in \C$ und $n \in \Zofk$
  sowie $\dPs_{\aln,s}$ das zu $\sjk$-gehörige Matrixpolynom bezüglich $\widehat{P}_{\aln,s}$ für alle $n \in \Zofk$.
  Dann heißt $\SQsalk$ das \textbf{linksseitige $\alpha$-Stieltjes"=Quadrupel} bezüglich $\sjk$.
  Falls klar ist, von welchem $\sjk$ die Rede ist, lassen wir das \anf{$s$} im Index weg.
\end{defi}

Unsere nachfolgenden Überlegungen sind nun auf die Herleitung von Zusammenhängen zwischen dem gerade eingeführten linksseitigen $\alpha$-Stieltjes-Quadrupel bezüglich einer linksseitig $\alpha$-Stieltjes-positiv definiten Folge und dem linksseitigen $\alpha$-Dyukarev-Quadrupel bezüglich jener Folge orientiert. Wir werden zunächst das linksseitige $\alpha$-Stieltjes-Quadrupel auf eine Form bringen, die mit den Bezeichnungen aus Kapitel \ref{chapdr} konform ist. Hierfür werden \thref{sqlm1} und folgendes Lemma uns wichtige Hilfestellungen leisten.

\begin{lemma}	\thlabel{sqllm1}
  Seien $\kappa \in \Na$, $\alpha \in \R$ und $\sjk$ eine Folge aus $\Cqq$. \dgfa
  \begin{itemize}
    \item [\rm{(a)}] Für alle $z \in \C$ und $n \in \Zokm$ gilt
      \begin{align*}
	R_n(z)y_{\alo,n} = \Saln\En(z).
      \end{align*}
    \item [\rm{(b)}] Für alle $z \in \C$ und $n \in \Zok$ gilt
      \begin{align*}
	R_n(z)u_{\aln} = \dSaln\En(z).
      \end{align*}
  \end{itemize}
\end{lemma}

\bwanf Zu (a): Der Fall $n = 0$ folgt sogleich aus der Definition der beteiligten Größen (vergleiche \thref{drbz1}, Teil (b) von \thref{aspdef1}, \thref{sqlbz1} und \thref{adpbz1}). Seien nun $z \in \C$ und $n \in \Zek$. Dann gilt
\begin{align*}
  \Rn(z)\yaln = \begin{pmatrix} s_{\alo} \\ zs_{\alo}+s_{\al1} \\ \vdots \\ \sum^{n}_{j=0} z^{n-j}s_{\alj} \end{pmatrix} = \Saln\En(z).
\end{align*}

Zu (b): Der Fall $n = 0$ folgt sogleich aus der Definition der beteiligten Größen (vergleiche \thref{drbz1}, \thref{drlbz1}, \thref{sqlbz1} und \thref{adpbz1}). Seien nun $z \in \C$ und $n \in \Zekp$. Wegen Teil (b) und (c) von \thref{sqlm1} gilt dann
\begin{align*}
  \Rn(z)\ualn & = \Rn(z)\eklam{\alpha\un-\yn} = \alpha \begin{pmatrix} \Oqn \\ S_{n-1} \end{pmatrix} \Enm(z) - S_n\En(z)\\
  & = \eklam{\alpha \begin{pmatrix} \Oqn & \Oq \\ S_{n-1} & \Onq \end{pmatrix} - S_n}\En(z) = \dSaln\En(z).
\end{align*}
\bwend

\begin{bem}	\thlabel{sqlbm1}
  Seien $\kappa \in \Na$, $\alpha \in \R$, $\sjk \in \Lpqka$ und $\SQalk$ das linksseitige $\alpha$-Stieltjes-Quadrupel bezüglich $\sjk$.
  Dann gelten $P_0 \equiv \Iq$, $\Ps_0 \equiv \Oq$, $P_{\alo} \equiv \Iq$ und $\dPs_{\alo} \equiv -s_0$ sowie
  \begin{align*}
    P_n(z) = \begin{pmatrix} -z_{n,2n-1}\Hnm^{-1} & \Iq \end{pmatrix}\Rn(z)\vn, \\
    \Ps_n(z) = \begin{pmatrix} -z_{n,2n-1}\Hnm^{-1} & \Iq \end{pmatrix}\Rn(z)\un
  \end{align*}
  für alle $z \in \C$ und $n \in \Zefkp$ sowie im Fall $\kappa\geq2$
  \begin{align*}
    P_{\aln}(z) = \begin{pmatrix} -z_{\aln,2n-1}\Halnm^{-1} & \Iq \end{pmatrix}\Rn(z)\vn, \\
    \dPs_{\aln}(z) = \begin{pmatrix} -z_{\aln,2n-1}\Halnm^{-1} & \Iq \end{pmatrix}\Rn(z)\ualn
  \end{align*}
  für alle $z \in \C$ und $n \in \Zefk$.
\end{bem}

\bwanf Dies folgt aus \thref{sqlsa1}, \thref{sqlm1} und \thref{sqllm1}. \bwend

Für unsere weiteren Betrachtungen wird uns folgendes Lemma dazu dienen, die Resultate für den rechtsseitigen Fall auf den linksseitigen Fall zu übertragen.

\begin{lemma}	\thlabel{sqllm2}
	Seien $\kappa \in \Na$, $\alpha \in \R$, $\sjk$ eine Folge aus $\Cqq$ und $t_j:=(-1)^js_j$ für alle $j\in\Zok$. Weiterhin sei $\tjk\in\Lpqkma$ oder $\sjk\in\Kpqka$ erfüllt. Dann gelten $\sjk\in\Kpqka$ bzw. $\tjk\in\Lpqkma$. Seien nun $\SQsark$ das rechtsseitige $\alpha$-Stieltjes-Quadrupel bezüglich $\sjk$ und $[(P_{n,t})^{\fklam{\kappa+1}}_{n=0}$,""$(P^{\sklam{t}}_{n,t})^{\fklam{\kappa+1}}_{n=0}$,""$(P_{-\aln,t})^{\fklam{\kappa}}_{n=0}$,""$(\dP^{\sklam{t}}_{-\aln,t})^{\fklam{\kappa}}_{n=0}]$ das linksseitige $-\alpha$-Stieltjes-Quadrupel bezüglich $\tjk$. Dann gelten
	\begin{align*}
		P_{n,t}(-z) = (-1)^nP_{n,s}(z), \\
		P^{\sklam{t}}_{n,t}(-z) = (-1)^{n+1}\Ps_{n,s}(z)
	\end{align*}
	für alle $z \in \C$ und $n \in \Zofkp$ sowie
	\begin{align*}
		P_{-\aln,t}(-z) = (-1)^nP_{\arn,s}(z), \\
		\dP^{\sklam{t}}_{-\aln,t}(-z) = (-1)^{n+1}\dPs_{\arn,s}(z)
	\end{align*}
	für alle $z \in \C$ und $n \in \Zofk$.
\end{lemma}

\bwanf Wegen Teil (a) von \thref{asmbm3} gilt $\sjk\in\Kpqka$ bzw. $\tjk\in\Lpqkma$. Wegen \thref{sqlbm1} und \thref{sqbm1} gelten
\begin{align*}
	P_{0,t}(-z) &= \Iq = P_{0,s}(z), \\
	P^{\sklam{t}}_{0,t}(-z) &= \Oq = -\Ps_{0,s}(z), \\
	P_{-\alo,t}(-z) &= \Iq = P_{\aro,s}(z), \\
	\dP^{\sklam{t}}_{-\alo,t}(-z) &= -t_0 = -s_0 = -\dPs_{\aro,s}(z)
\end{align*}
für alle $z\in\C$.

Seien nun $z \in \C$ und $n \in \Zefkp$. Wegen \thref{sqlbm1}, der Teile (a), (b) und (d) von \thref{asmlm1}, der Teile (a), (b) und (c) von \thref{drllm1} sowie \thref{sqbm1} gelten dann
\begin{align*}
	P_{n,t}(-z) &= \begin{pmatrix} -z^{\sklam{t}}_{n,2n-1}\big(\Hnm^{\sklam{t}}\big)^{-1} & \Iq \end{pmatrix}\Rn(-z)\vn \\
	&= \begin{pmatrix} -(-1)^nz^{\sklam{s}}_{n,2n-1}\Vnma\Vnm\big(\Hsnm\big)^{-1}\Vnma & \Iq \end{pmatrix}\Vn\Rn(z)\Vna\Vn\vn \\
	&= (-1)^n\begin{pmatrix} -z^{\sklam{s}}_{n,2n-1}\Vnma\Vnm\big(\Hsnm\big)^{-1} & \Iq \end{pmatrix}\Vna\Vn\Rn(z)\Vna\Vn\vn \\
	&= (-1)^n\begin{pmatrix} -z^{\sklam{s}}_{n,2n-1}\big(\Hsnm\big)^{-1} & \Iq \end{pmatrix}\Rn(z)\vn = (-1)^nP_{n,s}(z)
\end{align*}
und
\begin{align*}
	P^{\sklam{t}}_{n,t}(-z) &= \begin{pmatrix} -z^{\sklam{t}}_{n,2n-1}\big(\Hnm^{\sklam{t}}\big)^{-1} & \Iq \end{pmatrix}\Rn(-z)\un^{\sklam{t}} \\
	&= -\begin{pmatrix} -(-1)^nz^{\sklam{s}}_{n,2n-1}\Vnma\Vnm\big(\Hsnm\big)^{-1}\Vnma & \Iq \end{pmatrix}\Vn\Rn(z)\Vna\Vn\usn \\
	&= (-1)^{n+1}\begin{pmatrix} -z^{\sklam{s}}_{n,2n-1}\Vnma\Vnm\big(\Hsnm\big)^{-1} & \Iq \end{pmatrix}\Vna\Vn\Rn(z)\Vna\Vn\usn \\
	&= (-1)^{n+1}\begin{pmatrix} -z^{\sklam{s}}_{n,2n-1}\big(\Hsnm\big)^{-1} & \Iq \end{pmatrix}\Rn(z)\usn = (-1)^{n+1}\Ps_{n,s}(z).
\end{align*}

Seien nun $z \in \C$, $\kappa\geq2$ und $n \in \Zefk$. Wegen \thref{sqlbm1}, Teil (a) von \thref{asmlm1}, der Teile (b) und (c) von \thref{asplm1}, der Teile (a), (b) und (c) von \thref{drllm1} sowie \thref{sqbm1} gelten dann
\begin{align*}
	P_{-\aln,t}(-z) &= \begin{pmatrix} -z^{\sklam{t}}_{-\aln,2n-1}\big(H^{\sklam{t}}_{-\aln-1}\big)^{-1} & \Iq \end{pmatrix}\Rn(-z)\vn \\
	&= \begin{pmatrix} -(-1)^nz^{\sklam{s}}_{\arn,2n-1}\Vnma\Vnm\big(\Hsarnm\big)^{-1}\Vnma & \Iq \end{pmatrix}\Vn\Rn(z)\Vna\Vn\vn \\
	&= (-1)^n\begin{pmatrix} -z^{\sklam{s}}_{\arn,2n-1}\Vnma\Vnm\big(\Hsarnm\big)^{-1} & \Iq \end{pmatrix}\Vna\Vn\Rn(z)\Vna\Vn\vn \\
	&= (-1)^n\begin{pmatrix} -z^{\sklam{s}}_{\arn,2n-1}\big(\Hsarnm\big)^{-1} & \Iq \end{pmatrix}\Rn(z)\vn = (-1)^nP_{\arn,s}(z)
\end{align*}
und
\begin{align*}
	\dP^{\sklam{t}}_{-\aln,t}(-z) &= \begin{pmatrix} -z^{\sklam{t}}_{-\aln,2n-1}\big(H^{\sklam{t}}_{-\aln-1}\big)^{-1} & \Iq \end{pmatrix}\Rn(-z)u^{\sklam{t}}_{-\aln} \\
	&= -\begin{pmatrix} -(-1)^nz^{\sklam{s}}_{\arn,2n-1}\Vnma\Vnm\big(\Hsarnm\big)^{-1}\Vnma & \Iq \end{pmatrix}\Vn\Rn(z)\Vna\Vn\usarn \\
	&= (-1)^{n+1}\begin{pmatrix} -z^{\sklam{s}}_{\arn,2n-1}\Vnma\Vnm\big(\Hsarnm\big)^{-1} & \Iq \end{pmatrix}\Vna\Vn\Rn(z)\Vna\Vn\usarn \\
	&= (-1)^{n+1}\begin{pmatrix} -z^{\sklam{s}}_{\arn,2n-1}\big(\Hsarnm\big)^{-1} & \Iq \end{pmatrix}\Rn(z)\usarn = (-1)^{n+1}\dPs_{\arn,s}(z). \tag*{$\Box$}
\end{align*}

Wir kommen nun zu einem ersten Zusammenhang zwischen linksseitigem $\alpha$-Stieltjes-Quadrupel und linksseitigem $\alpha$-Dyukarev-Quadrupel bezüglich einer linksseitigen $\alpha$-Stieltjes-positiv definiten Folge. 

\begin{satz}	\thlabel{sqlsa2}
  	Seien $\kappa \in \Na$, $\alpha \in \R$ und $\sjk \in \Lpqka$. Weiterhin seien $\SQalk$ das linksseitige $\alpha$-Stieltjes-Quadrupel bezüglich $\sjk$ und $\ABCDalk$ das linksseitige $\alpha$-Dyukarev-Qua\-dru\-pel bezüglich $\sjk$. Dann gelten
  \begin{align*}
    \Aaln(z) & = \Iq - (\alpha-z) \sum^{n}_{j=0} \beklam{\Ps_j(\za)}^{\ast}\dH^{-1}_jP_j(\alpha), \\
    \Caln(z) & = (\alpha-z) \sum^{n}_{j=0} \beklam{P_j(\za)}^{\ast}\dH^{-1}_jP_j(\alpha)
  \end{align*}
  für alle $z \in \C$ und $n \in \Zofk$ und
  \begin{align*}
    \Baln(z) & = -\sum^{n-1}_{j=0} \beklam{\dPs_{\alj}(\za)}^{\ast}\dH^{-1}_{\alj}\dPs_{\alj}(\alpha), \\
    \Daln(z) & = \Iq + (\alpha-z) \sum^{n-1}_{j=0} \beklam{P_{\alj}(\za)}^{\ast}\dH^{-1}_{\alj}\dPs_{\alj}(\alpha)
  \end{align*}
  für alle $z \in \C$ und $n \in \Zefkp$.
\end{satz}

\bwanf Sei $t_j := (-1)^js_j$ für alle $j\in\Zok$. Wegen Teil (a) von \thref{asmbm3} gilt dann $\tjk\in\Kpqkma$. Seien $[(P_{n,t})^{\fklam{\kappa+1}}_{n=0}$,""$(P^{\sklam{t}}_{n,t})^{\fklam{\kappa+1}}_{n=0}$,""$(P_{-\arn,t})^{\fklam{\kappa}}_{n=0}$,""$(\dP^{\sklam{t}}_{-\arn,t})^{\fklam{\kappa}}_{n=0}]$ das rechtsseitige $-\alpha$-Stieltjes-Quadrupel bezüglich $\tjk$ und $[(\mathbf{A}^{\sklam{t}}_{-\arn})^{\fklam{\kappa}}_{n=0}$,""$(\mathbf{B}^{\sklam{t}}_{-\arn})^{\fklam{\kappa+1}}_{n=0}$,""$(\mathbf{C}^{\sklam{t}}_{-\arn})^{\fklam{\kappa}}_{n=0}$,""$(\mathbf{D}^{\sklam{t}}_{-\arn})^{\fklam{\kappa+1}}_{n=0}]$ das rechtsseitige $-\alpha$-Dyukarev-Quadrupel bezüglich $\tjk$. Wegen \thref{drllm5}, \thref{sqsa2}, \thref{sqllm2} und Teil (e) von \thref{asmlm1} gelten dann
\begin{align*}
	\Aaln^{\sklam{s}}(z) &= \mathbf{A}^{\sklam{t}}_{-\arn}(-z) \\
	&= \Iq + (-z-(-\alpha)) \sum^{n}_{j=0} \beklam{P^{\sklam{t}}_{j,t}(-\za)}^{\ast}\big(\dH^{\sklam{t}}_j\big)^{-1}P_{j,t}(-\alpha) \\
	&= \Iq - (\alpha-z) \sum^{n}_{j=0} \beklam{\Ps_{j,s}(\za)}^{\ast}\big(\dH^{\sklam{s}}_j\big)^{-1}P_{j,s}(\alpha)
\end{align*}
und
\begin{align*}
	\Caln^{\sklam{s}}(z) &= -\mathbf{C}^{\sklam{t}}_{-\arn}(-z) \\
	&= (-z-(-\alpha)) \sum^{n}_{j=0} \beklam{P_{j,t}(-\za)}^{\ast}\big(\dH^{\sklam{t}}_j\big)^{-1}P_{j,t}(-\alpha) \\
	&= (\alpha-z) \sum^{n}_{j=0} \beklam{P_{j,s}(\za)}^{\ast}\big(\dH^{\sklam{s}}_j\big)^{-1}P_{j,s}(\alpha)
\end{align*}
für alle $z \in \C$ und $n \in \Zofk$. Wegen \thref{drllm5}, \thref{sqsa2}, \thref{sqllm2} und Teil (d) von \thref{asplm1} gelten weiterhin
\begin{align*}
	\Baln^{\sklam{s}}(z) &= -\mathbf{B}^{\sklam{t}}_{-\arn}(-z) \\
	&= -\sum^{n-1}_{j=0} \beklam{\dP^{\sklam{t}}_{-\arj,t}(-\za)}^{\ast}\big(\dH^{\sklam{t}}_{-\arj}\big)^{-1}\dP^{\sklam{t}}_{-\arj,t}(-\alpha) \\
	&= -\sum^{n-1}_{j=0} \beklam{\dPs_{\alj,s}(\za)}^{\ast}\big(\dH^{\sklam{s}}_{\alj}\big)^{-1}\dPs_{\alj,s}(\alpha)
\end{align*}
und
\begin{align*}
	\Daln^{\sklam{s}}(z) &= \mathbf{D}^{\sklam{t}}_{-\arn}(-z) \\
	&= \Iq - (-z-(-\alpha)) \sum^{n-1}_{j=0} \beklam{P_{-\arj,t}(-\za)}^{\ast}\big(\dH^{\sklam{t}}_{-\arj}\big)^{-1}\dP^{\sklam{t}}_{-\arj,t}(-\alpha) \\
	&= \Iq + (\alpha-z) \sum^{n-1}_{j=0} \beklam{P_{\alj,s}(\za)}^{\ast}\big(\dH^{\sklam{s}}_{\alj}\big)^{-1}\dPs_{\alj,s}(\alpha)
\end{align*}
für alle $z \in \C$ und $n \in \Zefkp$. \bwend

Folgendes Resultat liefert uns eine Darstellung der Matrixpolynome des linksseitigen $\alpha$-Stieltjes-Quadrupels bezüglich einer linksseitig $\alpha$-Stieltjes-positiv definiten Folge an der Stelle $\alpha$ mithilfe der linksseitigen $\alpha$-Dyukarev-Stieltjes-Parametrisierung jener Folge.

\begin{satz}	\thlabel{sqlsa3}
	Seien $\kappa \in \Na$, $\alpha \in \R$ und $\sjk \in \Lpqka$. Weiterhin seien $\SQalk$ das linksseitige $\alpha$-Stieltjes-Quadrupel bezüglich $\sjk$ und $\LMkaln$ die linksseitige $\alpha$-Dyukarev-Stieltjes-Parametrisierung von $\sjk$. \dgfa
	\begin{itemize}
		\item [\rm{(a)}] Es sind $\brklam{P_n(\alpha)}^{\fklam{\kappa+1}}_{n=0}$, $\brklam{\Ps_n(\alpha)}^{\fklam{\kappa+1}}_{n=1}$, $\brklam{P_{\aln}(\alpha)}^{\fklam{\kappa}}_{n=0}$ und $\brklam{\dPs_{\aln}(\alpha)}^{\fklam{\kappa}}_{n=0}$ Folgen von regulären Matrizen aus $\Cqq$.
		\item [\rm{(b)}] Es gelten
		\begin{align*}
			P_n(\alpha) &= \prodr^{n-1}_{j=0}\rklam{\Malj^{-1}\Lalj^{-1}}, \\
			\Ps_n(\alpha) &= \prodr^{n-1}_{j=0}\rklam{\Malj^{-1}\Lalj^{-1}}\sum^{n-1}_{j=0}\Lalj
		\end{align*}
		für alle $n\in\Zefkp$ und im Fall $\kappa\geq2$
		\begin{align*}
			P_{\aln}(\alpha) &= \prodr^{n-1}_{j=0}\rklam{\Malj^{-1}\Lalj^{-1}}\Maln^{-1}\sum^n_{j=0}\Malj, \\
			\dPs_{\aln}(\alpha) &= -\prodr^{n-1}_{j=0}\rklam{\Malj^{-1}\Lalj^{-1}}\Maln^{-1}
		\end{align*}
		für alle $n\in\Zefk$.
	\end{itemize}
\end{satz}

\bwanf Sei $t_j := (-1)^js_j$ für alle $j\in\Zok$. Wegen Teil (a) von \thref{asmbm3} gilt dann $\tjk\in\Kpqkma$. Seien $[(P_{n,t})^{\fklam{\kappa+1}}_{n=0}$,""$(P^{\sklam{t}}_{n,t})^{\fklam{\kappa+1}}_{n=0}$,""$(P_{-\arn,t})^{\fklam{\kappa}}_{n=0}$,""$(\dP^{\sklam{t}}_{-\arn,t})^{\fklam{\kappa}}_{n=0}]$ das rechtsseitige $-\alpha$-Stieltjes-Quadrupel bezüglich $\tjk$ und $[(\mathbf{L}^{\sklam{t}}_{-\arn})^{\fklam{\kappa-1}}_{n=0}$,""$(\mathbf{M}^{\sklam{t}}_{-\arn})^{\fklam{\kappa}}_{n=0}]$ die rechtsseitige $-\alpha$-Dyukarev-Stieltjes-Parametrisierung von $\tjk$.

Zu (a): Wegen Teil (a) von \thref{sqsa3} sind $\big(P_{n,t}(-\alpha)\big)^{\fklam{\kappa+1}}_{n=0}$, $\big(P^{\sklam{t}}_{n,t}(-\alpha)\big)^{\fklam{\kappa+1}}_{n=1}$, \linebreak $\big(P_{-\arn,t}(-\alpha)\big)^{\fklam{\kappa}}_{n=0}$ und $\big(\dP^{\sklam{t}}_{-\arn,t}(-\alpha)\big)^{\fklam{\kappa}}_{n=0}$ Folgen von regulären Matrizen aus $\Cqq$. Hier\-aus folgt wegen \thref{sqllm2} dann die Behauptung.

Zu (b): Wegen \thref{sqllm2}, Teil (b) von \thref{sqsa3} und Teil (c) von \thref{adpllm1} gelten
\begin{align*}
	P_{n,s}(\alpha) &= (-1)^n P_{n,t}(-\alpha) = \prodr^{n-1}_{j=0}\eklam{\big(\mathbf{M}^{\sklam{t}}_{-\arj}\big)^{-1}\big(\mathbf{L}^{\sklam{t}}_{-\arj}\big)^{-1}} = \prodr^{n-1}_{j=0}\eklam{\big(\Malj^{\sklam{s}}\big)^{-1}\big(\Lalj^{\sklam{s}}\big)^{-1}}
\end{align*}
und
\begin{align*}
	\Ps_{n,s}(\alpha) &= (-1)^{n+1} P^{\sklam{t}}_{n,t}(-\alpha) \\
	&= \prodr^{n-1}_{j=0}\eklam{\big(\mathbf{M}^{\sklam{t}}_{-\arj}\big)^{-1}\big(\mathbf{L}^{\sklam{t}}_{-\arj}\big)^{-1}}\sum^{n-1}_{j=0}\mathbf{L}^{\sklam{t}}_{-\arj} \\
	&= \prodr^{n-1}_{j=0}\eklam{\big(\Malj^{\sklam{s}}\big)^{-1}\big(\Lalj^{\sklam{s}}\big)^{-1}}\sum^{n-1}_{j=0}\Lalj^{\sklam{s}}
\end{align*}
für alle $n\in\Zefkp$ sowie im Fall $\kappa\geq2$
\begin{align*}
	P_{\aln,s}(\alpha) &= (-1)^n P_{-\arn,t}(-\alpha) \\
	&= \prodr^{n-1}_{j=0}\eklam{\big(\mathbf{M}^{\sklam{t}}_{-\arj}\big)^{-1}\big(\mathbf{L}^{\sklam{t}}_{-\arj}\big)^{-1}}\big(\mathbf{M}^{\sklam{t}}_{-\arn}\big)^{-1}\sum^n_{j=0}\mathbf{M}^{\sklam{t}}_{-\arj} \\
	&= \prodr^{n-1}_{j=0}\eklam{\big(\Malj^{\sklam{s}}\big)^{-1}\big(\Lalj^{\sklam{s}}\big)^{-1}}\big(\Maln^{\sklam{s}}\big)^{-1}\sum^n_{j=0}\Malj^{\sklam{s}}
\end{align*}
und
\begin{align*}
	\dPs_{\aln,s}(\alpha) &= (-1)^{n+1} \dP^{\sklam{t}}_{-\arn,t}(-\alpha) \\
	&= -\prodr^{n-1}_{j=0}\eklam{\big(\mathbf{M}^{\sklam{t}}_{-\arj}\big)^{-1}\big(\mathbf{L}^{\sklam{t}}_{-\arj}\big)^{-1}}\big(\mathbf{M}^{\sklam{t}}_{-\arn}\big)^{-1} \\
	&= -\prodr^{n-1}_{j=0}\eklam{\big(\Malj^{\sklam{s}}\big)^{-1}\big(\Lalj^{\sklam{s}}\big)^{-1}}\big(\Maln^{\sklam{s}}\big)^{-1}
\end{align*}
für alle $n\in\Zefk$. \bwend

Wir können nun die linksseitigen $2$\textit{q}$\times2$\textit{q}-$\alpha$-Dyukarev-Matrixpolynome bezüglich einer linksseitig $\alpha$-Stieltjes-positiv definiten Folge mithilfe des linksseitigen $\alpha$-Stieltjes-Quadrupels bezüglich jener Folge darstellen.

\begin{satz}	\thlabel{sqlsa4}
	Seien $\kappa \in \Na$, $\alpha \in \R$, $\sjk \in \Lpqka$ und $\SQalk$ das linksseitige $\alpha$-Stieltjes-Quadrupel bezüglich $\sjk$. \dgfa
	\begin{itemize}
		\item [\rm{(a)}] Sei $\ABCDalk$ das linksseitige $\alpha$"=Dyukarev"=\linebreak Quadrupel bezüglich $\sjk$. Dann sind $P_n(\alpha)$ für alle $n \in \Zofkp$ sowie $\dPs_{\aln}(\alpha)$ für alle $n \in \Zofk$ reguläre Matrizen und es gelten
		\begin{align*}
			\Aaln(z) & = \beklam{\dPs_{\aln}(\za)}^{\ast}\beklam{\dPs_{\aln}(\alpha)}^{-\ast}, \\
    		\Caln(z) & =  -(\alpha-z)\beklam{P_{\aln}(\za)}^{\ast}\beklam{\dPs_{\aln}(\alpha)}^{-\ast}
  		\end{align*}
  		für alle $z \in \C$ und $n \in \Zofk$ sowie
  		\begin{align*}
   			\Baln(z) & = -\beklam{\Ps_n(\za)}^{\ast}\beklam{P_n(\alpha)}^{-\ast}, \\
    		\Daln(z) & = \beklam{P_n(\za)}^{\ast}\beklam{P_n(\alpha)}^{-\ast}
  		\end{align*}
  		für alle $z \in \C$ und $n \in \Zofkp$.
  		\item [\rm{(b)}] Sei $\Ualmk$ die Folge von linksseitigen $2$\textit{q}$\times2$\textit{q}-$\alpha$-Dyukarev-Matrixpolynomen bezüglich $\sjk$. Weiterhin seien $z \in \C$ und $m \in \Zok$. Dann gilt
  		\begin{align}	\label{sqlsa4bw1}
  			\Ualm(z) = \begin{pmatrix} \beklam{\dPs_{\al\fklam{m}}(\za)}^{\ast}\beklam{\dPs_{\al\fklam{m}}(\alpha)}^{-\ast} & -\beklam{\Ps_{\fklam{m+1}}(\za)}^{\ast}\beklam{P_{\fklam{m+1}}(\alpha)}^{-\ast} \\ -(\alpha-z)\beklam{P_{\al\fklam{m}}(\za)}^{\ast}\beklam{\dPs_{\al\fklam{m}}(\alpha)}^{-\ast} & \beklam{P_{\fklam{m+1}}(\za)}^{\ast}\beklam{P_{\fklam{m+1}}(\alpha)}^{-\ast} \end{pmatrix}.
  		\end{align}
  		Außerdem ist $\Ualm(z)$ eine reguläre Matrix und es gilt
  		\begin{align*}
  			\Ualm^{-1}(z) = \begin{pmatrix} \beklam{P_{\fklam{m+1}}(\alpha)}^{-1}P_{\fklam{m+1}}(z) & \beklam{P_{\fklam{m+1}}(\alpha)}^{-1}\Ps_{\fklam{m+1}}(z) \\ (\alpha-z)\beklam{\dPs_{\al\fklam{m}}(\alpha)}^{-1}P_{\al\fklam{m}}(z) & \beklam{\dPs_{\al\fklam{m}}(\alpha)}^{-1}\dPs_{\al\fklam{m}}(z) \end{pmatrix}.
  		\end{align*}
  		\item [\rm{(c)}] Seien $m\in\Zek$ und $\Salmmin$ bzw. $\Salmmax$ das untere bzw. obere Extremalelement von $\Sqmasmu$. Dann gelten
  		\begin{align*}
  			\Salmmin(z) &= -\frac{1}{\alpha-z}\beklam{\dPs_{\al\fklam{m}}(\za)}^{\ast}\beklam{P_{\al\fklam{m}}(\za)}^{-\ast}, \\
  			\Salmmax(z) &= -\beklam{\Ps_{\fklam{m+1}}(\za)}^{\ast}\beklam{P_{\fklam{m+1}}(\za)}^{-\ast}
  		\end{align*}
  		für alle $z \in \C\setminus(-\infty,\alpha]$.  		
  	\end{itemize}		
\end{satz}

\bwanf Seien $t_j = (-1)^js_j$ für alle $j\in\Zok$. Wegen Teil (a) von \thref{asmbm3} gilt dann $\tjk\in\Kpqkma$. Sei $[(P_{n,t})^{\fklam{\kappa+1}}_{n=0}$,""$(P^{\sklam{t}}_{n,t})^{\fklam{\kappa+1}}_{n=0}$,""$(P_{-\arn,t})^{\fklam{\kappa}}_{n=0}$,""$(\dP^{\sklam{t}}_{-\arn,t})^{\fklam{\kappa}}_{n=0}]$ das rechtsseitige $-\alpha$-Stieltjes-Quadrupel bezüglich $\tjk$.

Zu (a): Wegen Teil (a) von \thref{sqlsa3} sind $P_n(\alpha)$ für alle $n \in \Zofkp$ sowie $\dPs_{\aln}(\alpha)$ für alle $n \in \Zofk$ reguläre Matrizen. Sei $[(\mathbf{A}^{\sklam{t}}_{-\arn})^{\fklam{\kappa}}_{n=0}$,""$(\mathbf{B}^{\sklam{t}}_{-\arn})^{\fklam{\kappa+1}}_{n=0}$,""$(\mathbf{C}^{\sklam{t}}_{-\arn})^{\fklam{\kappa}}_{n=0}$,""$(\mathbf{D}^{\sklam{t}}_{-\arn})^{\fklam{\kappa+1}}_{n=0}]$ das rechtsseitige $-\alpha$-Dyukarev-Quadrupel bezüglich $\tjk$. Wegen \thref{drllm5}, Teil (a) von \thref{sqsa4} und \thref{sqllm2} gelten
\begin{align*}
	\Aaln^{\sklam{s}}(z) = \mathbf{A}^{\sklam{t}}_{-\arn}(-z) 
	= \beklam{\dP^{\sklam{t}}_{-\arn,t}(-\za)}^{\ast}\beklam{\dP^{\sklam{t}}_{-\arn,t}(-\alpha)}^{-\ast} 
	= \beklam{\dPs_{\aln,s}(\za)}^{\ast}\beklam{\dPs_{\aln,s}(\alpha)}^{-\ast}
\end{align*}
und
\begin{align*}
	\Caln^{\sklam{s}}(z) &= -\mathbf{C}^{\sklam{t}}_{-\arn}(-z) \\
	&= (-z-(-\alpha))\beklam{P_{-\arn,t}(-\za)}^{\ast}\beklam{\dP^{\sklam{t}}_{-\arn,t}(-\alpha)}^{-\ast} \\
	&= -(\alpha-z)\beklam{P_{\aln,s}(\za)}^{\ast}\beklam{\dPs_{\aln,s}(\alpha)}^{-\ast}
\end{align*}
für alle $z \in \C$ und $n \in \Zofk$ sowie
\begin{align*}
	\Baln^{\sklam{s}}(z) = -\mathbf{B}^{\sklam{t}}_{-\arn}(-z)
	= \beklam{P^{\sklam{t}}_{n,t}(-\za)}^{\ast}\beklam{P_{n,t}(-\alpha)}^{-\ast}
	= -\beklam{\Ps_{n,s}(\za)}^{\ast}\beklam{P_{n,s}(\alpha)}^{-\ast}
\end{align*}
und
\begin{align*}
	\Daln^{\sklam{s}}(z) = \mathbf{D}^{\sklam{t}}_{-\arn}(-z)
	= \beklam{P_{n,t}(-\za)}^{\ast}\beklam{P_{n,t}(-\alpha)}^{-\ast}
	= \beklam{P_{n,s}(\za)}^{\ast}\beklam{P_{n,s}(\alpha)}^{-\ast}
\end{align*}
für alle $z \in \C$ und $n \in \Zofkp$.

Zu (b): \fref{sqlsa4bw1} folgt wegen (a) aus \thref{drldef3}. Wegen \thref{spbsp1}, Teil (c) von \thref{drlbm2}, Teil (a) von \thref{splm1} und \fref{sqlsa4bw1} ist $\Ualm(z)$ eine reguläre Matrix und es gilt
\begin{align*}
	&\ \Ualm^{-1}(z) = \tJq\beklam{\Ualm(\za)}^{\ast}\tJq \\
	&= \begin{pmatrix} \Oq & -i\Iq \\ i\Iq & \Oq \end{pmatrix} 
	\begin{pmatrix} \beklam{\dPs_{\al\fklam{m}}(\alpha)}^{-1}\dPs_{\al\fklam{m}}(z) & -(\alpha-z)\beklam{\dPs_{\al\fklam{m}}(\alpha)}^{-1}P_{\al\fklam{m}}(z) \\ -\beklam{P_{\fklam{m+1}}(\alpha)}^{-1}\Ps_{\fklam{m+1}}(z) & \beklam{P_{\fklam{m+1}}(\alpha)}^{-1}P_{\fklam{m+1}}(z) \end{pmatrix} \\
	&\quad \cdot\begin{pmatrix} \Oq & -i\Iq \\ i\Iq & \Oq \end{pmatrix} \\
	&= \begin{pmatrix} i\beklam{P_{\fklam{m+1}}(\alpha)}^{-1}\Ps_{\fklam{m+1}}(z) & -i\beklam{P_{\fklam{m+1}}(\alpha)}^{-1}P_{\fklam{m+1}}(z) \\ i\beklam{\dPs_{\al\fklam{m}}(\alpha)}^{-1}\dPs_{\al\fklam{m}}(z) & -i(\alpha-z)\beklam{\dPs_{\al\fklam{m}}(\alpha)}^{-1}P_{\al\fklam{m}}(z) \end{pmatrix} \begin{pmatrix} \Oq & -i\Iq \\ i\Iq & \Oq \end{pmatrix} \\
	&= \begin{pmatrix} \beklam{P_{\fklam{m+1}}(\alpha)}^{-1}P_{\fklam{m+1}}(z) & \beklam{P_{\fklam{m+1}}(\alpha)}^{-1}\Ps_{\fklam{m+1}}(z) \\ (\alpha-z)\beklam{\dPs_{\al\fklam{m}}(\alpha)}^{-1}P_{\al\fklam{m}}(z) & \beklam{\dPs_{\al\fklam{m}}(\alpha)}^{-1}\dPs_{\al\fklam{m}}(z) \end{pmatrix}.
\end{align*}

Zu (c): Wegen \thref{drldef4} und (a) gelten
\begin{align*}
	\Salmmin(z) &= \Aalfm(z)\Calfm^{-1}(z) \\
	&= \beklam{\dPs_{\al\fklam{m}}(\za)}^{\ast}\beklam{\dPs_{\al\fklam{m}}(\alpha)}^{-\ast}\frac{-1}{\alpha-z}\beklam{\dPs_{\al\fklam{m}}(\alpha)}^{\ast}\beklam{P_{\al\fklam{m}}(\za)}^{-\ast} \\
	&= -\frac{1}{\alpha-z}\beklam{\dPs_{\al\fklam{m}}(\za)}^{\ast}\beklam{P_{\al\fklam{m}}(\za)}^{-\ast}
\end{align*}
und
\begin{align*}
	\Salmmax(z) &= \Balfm(z)\Dalfm^{-1}(z) \\
	&= -\beklam{\Ps_{\fklam{m+1}}(\za)}^{\ast}\beklam{P_{\fklam{m+1}}(\alpha)}^{-\ast}\beklam{P_{\fklam{m+1}}(\alpha)}^{\ast}\beklam{P_{\fklam{m+1}}(\za)}^{-\ast} \\
	&= -\beklam{\Ps_{\fklam{m+1}}(\za)}^{\ast}\beklam{P_{\fklam{m+1}}(\za)}^{-\ast}	
\end{align*}
für alle $z \in \C\setminus(-\infty,\alpha]$. \bwend

Im folgenden Satz betrachten wir die Lokalisierung der Nullstellen der Determinanten der einzelnen \textit{q}$\times$\textit{q}"=Matrixpolynome des linksseitigen $\alpha$-Stieltjes-Quadrupels bezüglich einer linksseitig $\alpha$-Stieltjes-positiv definiten Folge.

\begin{satz}	\thlabel{sqlsa7}
	Seien $\kappa \in \Na$, $\alpha \in \R$ und $\sjk \in \Lpqka$. Weiterhin seien $\SQalk$ das linksseitige $\alpha$-Stieltjes-Quadrupel bezüglich $\sjk$ und $z\in\C\setminus(-\infty,\alpha)$. Dann gelten $\det P_n(z) \neq 0$ für alle $n\in\Zofkp$, $\det \Ps_n(z) \neq 0$ für alle $n\in\Zefkp$, $\det P_{\aln}(z) \neq 0$ für alle $n\in\Zofk$ und $\det\dPs_{\aln}(z) \neq 0$ für alle $n\in\Zofk$.	
\end{satz}

\bwanf Dies folgt aus \thref{sqllm2} und \thref{sqsa7}. \bwend

Wir wollen nun die linksseitige $\alpha$-Stieltjes-Parametrisierung einer linksseitig $\alpha$-Stieltjes-positiv definiten Folge mithilfe des linksseitigen $\alpha$-Stieltjes-Quadrupels bezüglich jener Folge darstellen. 

\begin{satz}	\thlabel{sqlsa5}
	Seien $\kappa \in \Na$, $\alpha \in \R$ und $\sjk \in \Lpqka$. Weiterhin seien $\SQalk$ das linksseitige $\alpha$-Stieltjes-Quadrupel bezüglich $\sjk$ und $\Qaljk$ die linksseitige $\alpha$-Stieltjes-Parametrisierung von $\sjk$. \dgfa
	\begin{itemize}
		\item [\rm{(a)}] Es gelten
		\begin{align*}
			Q_{\al 2n} = -P_n(\alpha)\beklam{\dPs_{\aln}(\alpha)}^{\ast}
		\end{align*}
		für alle $n\in\Zofk$ und
		\begin{align*}
			Q_{\al 2n+1} = \dPs_{\aln}(\alpha)\beklam{P_{n+1}(\alpha)}^{\ast}
		\end{align*}
		für alle $n\in\Zofkm$.
		\item [\rm{(b)}] Es gelten
		\begin{align*}
			P_n(\alpha) = \prodl^{n-1}_{j=0}\rklam{Q_{\al2j+1}Q^{-1}_{\al2j}}
		\end{align*}
		für alle $n\in\Zefkp$ und im Fall $\kappa\geq2$
		\begin{align*}
			\dPs_{\aln}(\alpha) = -Q_{\al2n}\prodl^{n-1}_{j=0}\rklam{Q^{-1}_{\al2j+1}Q_{\al2j}}
		\end{align*}
		für alle $n\in\Zefk$.
	\end{itemize}
\end{satz}

\bwanf Seien $t_j = (-1)^js_j$ für alle $j\in\Zok$. Wegen Teil (a) von \thref{asmbm3} gilt dann $\tjk\in\Kpqkma$. Seien $[(P_{n,t})^{\fklam{\kappa+1}}_{n=0}$,""$(P^{\sklam{t}}_{n,t})^{\fklam{\kappa+1}}_{n=0}$,""$(P_{-\arn,t})^{\fklam{\kappa}}_{n=0}$,""$(\dP^{\sklam{t}}_{-\arn,t})^{\fklam{\kappa}}_{n=0}]$ das rechtsseitige $-\alpha$-Stieltjes-Quadrupel bezüglich $\tjk$ und $(Q^{\sklam{t}}_{-\arj})^{\kappa}_{j=0}$ die rechtsseitige $-\alpha$-Stieltjes-Parametrisierung von $\tjk$.

Zu (a): Wegen \thref{aspbm4}, Teil (a) von \thref{sqsa5} und \thref{sqllm2} gelten
\begin{align*}
	Q^{\sklam{s}}_{\al 2n} = Q^{\sklam{t}}_{-\ar 2n}
	= P_{n,t}(-\alpha)\beklam{\dP^{\sklam{t}}_{-\arn,t}(-\alpha)}^{\ast}
	= -P_{n,s}(\alpha)\beklam{\dPs_{\aln,s}(\alpha)}^{\ast}
\end{align*}
für alle $n\in\Zofk$ und
\begin{align*}
	Q^{\sklam{s}}_{\al 2n+1} = Q^{\sklam{t}}_{-\ar 2n+1}
	= -\dP^{\sklam{t}}_{-\arn,t}(-\alpha)\beklam{P_{n+1,t}(-\alpha)}^{\ast}
	= \dPs_{\aln,s}(\alpha)\beklam{P_{n+1,s}(\alpha)}^{\ast}
\end{align*}
für alle $n\in\Zofkm$.

Zu (b): Wegen \thref{sqllm2}, Teil (b) von \thref{sqsa5} und \thref{aspbm4} gelten
\begin{align*}
	P_{n,s}(\alpha) = (-1)^n P_{n,t}(-\alpha)
	= \prodl^{n-1}_{j=0}\eklam{Q^{\sklam{t}}_{-\ar2j+1}\big(Q^{\sklam{t}}_{-\ar2j}\big)^{-1}}
	= \prodl^{n-1}_{j=0}\eklam{Q^{\sklam{s}}_{\al2j+1}\big(Q^{\sklam{s}}_{\al2j}\big)^{-1}}
\end{align*}
für alle $n\in\Zefkp$ und
\begin{align*}
	\dPs_{\aln,s}(\alpha) &= (-1)^{n+1} \dP^{\sklam{t}}_{-\arn,t}(-\alpha) \\
	&= -Q^{\sklam{t}}_{-\ar2n}\prodl^{n-1}_{j=0}\eklam{\big(Q^{\sklam{t}}_{-\ar2j+1}\big)^{-1}Q^{\sklam{t}}_{-\ar2j}} \\
	&= -Q^{\sklam{s}}_{\al2n}\prodl^{n-1}_{j=0}\eklam{\big(Q^{\sklam{s}}_{\al2j+1}\big)^{-1}Q^{\sklam{s}}_{\al2j}}
\end{align*}
für alle $n\in\Zefk$. \bwend

Abschließend zeigen wir noch einen Zusammenhang zwischen linksseitigem $\alpha$-Stieltjes-Quadrupel bezüglich einer linksseitig $\alpha$-Stieltjes-positiv definiten Folge und den Favard-Paaren bezüglich jener Folge und derer durch linksseitige $\alpha$-Verschiebung generierten Folge.

\begin{satz}	\thlabel{sqlsa6}
	Seien $\kappa \in \Na$, $\alpha \in \R$ und $\sjk \in \Lpqka$. Weiterhin seien $\SQalk$ das linksseitige $\alpha$-Stieltjes-Quadrupel bezüglich $\sjk$, $\ABfk$ das Favard-Paar bezüglich $\sjk$ und $(B_{\alj})^0_{j=0}$ bzw. im Fall \linebreak $\kappa\geq2$ $\ABalfk$ das Favard-Paar bezüglich $\saljk$. \dgfa
	\begin{itemize}
		\item [\rm{(a)}] Es gelten $P_{0} \equiv \Iq$,
    	\begin{align*}
			P_{1}(z) = (z\Iq-A_0)P_{0}(z)
   		\end{align*}
    	für alle $z \in \C$ und im Fall $\kappa \geq 3$
    	\begin{align*}
			P_{n}(z) = (z\Iq-A_{n-1})P_{n-1}(z)-B^{\ast}_{n-1}P_{n-2}(z)
    	\end{align*}
   		für alle $z \in \C$ und $n \in \Zzfkp$.
   		\item [\rm{(b)}] Es gelten $\Ps_{0} \equiv 0_{\qq}$, $\Ps_{1} \equiv B_{0}$ und im Fall $\kappa \geq 3$
      	\begin{align*}
			\Ps_{n}(z) = (z\Iq-A_{n-1})\Ps_{n-1}(z) - B^{\ast}_{n-1}\Ps_{n-2}(z)
     	\end{align*}
      	für alle $z \in \C$ und $n \in \Zzfkp$.
      	\item [\rm{(c)}] Es gelten $P_{\alo} \equiv \Iq$, im Fall $\kappa \geq 2$
    	\begin{align*}
			P_{\al1}(z) = (z\Iq-A_{\alo})P_{\alo}(z)
   		\end{align*}
    	für alle $z \in \C$ und im Fall $\kappa \geq 4$
    	\begin{align*}
			P_{\aln}(z) = (z\Iq-A_{\aln-1})P_{\aln-1}(z)-B^{\ast}_{\aln-1}P_{\aln-2}(z)
    	\end{align*}
   		für alle $z \in \C$ und $n \in \Zzfk$.
   		\item [\rm{(d)}] Es gelten $\dPs_{\alo} \equiv -B_0$, im Fall $\kappa \geq 2$
    	\begin{align*}
			\dPs_{\al1}(z) = B_{\alo}+(z\Iq-A_{\alo})\dPs_{\alo}(z)
   		\end{align*}
    	für alle $z \in \C$ und im Fall $\kappa \geq 4$
		\begin{align*}
			\dPs_{\aln}(z) = (z\Iq-A_{\aln-1})\dPs_{\aln-1}(z)-B^{\ast}_{\aln-1}\dPs_{\aln-2}(z)
    	\end{align*}
   		für alle $z \in \C$ und $n \in \Zzfk$.	
	\end{itemize}	
\end{satz}

\bwanf Wegen Teil (b) von \thref{asmdef3} und \thref{asmbm1} gelten \linebreak $\sjfkm\in\Hpqfkm$ und $(s_{\alj})^{2\fklam{\kappa-2}}_{j=0} \in \Hpqfkzm$, falls $\kappa\geq2$.

Zu (a) und (c): Dies folgt aus \thref{sqldef1} und \thref{mpfo2}.

Zu (b): Dies folgt aus \thref{sqldef1} und \thref{mpfo3}.

Zu (d): Sei $\Psalnfk$ das linke System von Matrixpolynomen zweiter Art bezüglich $\saljk$. Wegen Teil (e) von \thref{sqlsa1} und \thref{fpdef1} gilt dann
\begin{align}	\label{sqlsa6bw1}
	\dPs_{\aln}(z) = \Pals_{\aln}(z) - P_{\aln}(z) B_0
\end{align}
für alle $z \in \C$ und $n \in \Zofk$. Wegen \thref{mpfo3} gelten weiterhin $\Pals_{\alo} \equiv 0_{\qq}$, im Fall $\kappa\geq2$ $\Pals_{\al1} \equiv B_{\alo}$ und im Fall $\kappa \geq 4$
\begin{align*}
	\Pals_{\aln}(z) = (z\Iq-A_{\aln-1})\Pals_{\aln-1}(z) - B^{\ast}_{\aln-1}\Pals_{\aln-2}(z)
\end{align*}
für alle $z \in \C$ und $n \in \Zzfk$. Hieraus folgt wegen \fref{sqlsa6bw1} und (c) dann die Behauptung. \bwend

%% file: sm5.tex
\newpage
\section{Weitere Zusammenhänge zwischen einigen Parametrisierungen \texorpdfstring{$\alpha$}{a}"=Stieltjes"=positiv definiter Folgen} \label{chapwz}

Nach der Einführung aller Parametrisierungen und ersten Zusammenhängen zwischen ihnen, können wir nun in diesem Kapitel durch Verknüpfen jener Ergebnisse neue Zusammenhänge schaffen. Insbesondere werden wir die Favard-Paare der beiden einer $\alpha$-Stieltjes-positiv definiten Folge zugrundeliegenden Hankel-positiv definiten Folgen sowohl mit der $\alpha$-Dyukarev-Stieltjes-Parametrisierung als auch mit der $\alpha$-Stieltjes-Parametrisierung jener Folge in Verbindung bringen.

Hierbei verfolgen wir im rechtsseitigen Fall die gleiche Vorgehensweise wie in \cite[Chapter 9]{CR1}, wo von einer gegebenen Folge aus ${\cal K}^{>}_{q,\infty,0}$ ausgegangen wurde. Im linksseitigen Fall werden wir wieder auf die Resultate für den rechtsseitigen Fall zurückgreifen.

\subsection{Der rechtsseitige Fall}

Bevor wir uns einem Zusammehang zwischen der rechtsseitgen $\alpha$-Dyukarev-Stieltjes-Paremetrisierung einer rechtsseitig $\alpha$-Stieltjes-positiv definiten Folge und den Favard-Paaren bezüglich jener Folge und deren durch rechtsseitige $\alpha$-Verschiebung generierten Folge widmen (vergleiche \cite[Theorem 9.3]{CR1} für den Fall $\alpha=0$ und $\kappa=\infty$), benötigen wir ein Hilfsresultat, das zusätzlich noch das rechtsseitige $\alpha$-Stieltjes-Quadrupel bezüglich jener Folge umfasst (vergleiche \cite[Lemma 9.2]{CR1} für den Fall $\alpha=0$ und $\kappa=\infty$).

\begin{lemma}	\thlabel{wzlm1}
	Seien $\kappa \in \Na$, $\alpha \in \R$ und $\sjk \in \Kpqka$. Weiterhin seien $\SQark$ das rechtsseitige $\alpha$-Stieltjes-Quadrupel bezüglich $\sjk$ und $\LMkarn$ die rechtsseitige $\alpha$-Dyukarev-Stieltjes-Parametrisierung \linebreak von $\sjk$. \dgfa
	\begin{itemize}
	\item [\rm{(a)}] Sei $\ABfk$ das Favard-Paar bezüglich $\sjk$. Dann gelten $B_0 = s_0$, $A_0-\alpha\Iq = -P_1(\alpha)$ und im Fall $\kappa\geq3$
	\begin{align*}
		B_n &= \beklam{P_{n-1}(\alpha)}^{-\ast}\Larnm^{-1}\Larn \beklam{P_{n+1}(\alpha)}^{\ast}, \\
		A_n-\alpha\Iq &= -P_{n+1}(\alpha)\rklam{\Larnm+\Larn}\Larnm^{-1}\beklam{P_n(\alpha)}^{-1}
	\end{align*}
	für alle $n\in\Zefkm$.
	\item [\rm{(b)}] Sei $(B_{\arn})^0_{n=0}$ bzw. im Fall $\kappa\geq2$ $\ABarfk$ das Favard-Paar bezüglich $\sarjk$. Dann gelten $B_{\aro} = s_{\aro}$, im Fall $\kappa\geq2$
	\begin{align*}
		A_{\arn}-\alpha\Iq &= -\dPs_{\arn+1}(\alpha)\rklam{\Marn+\Marnp}\Marn^{-1}\beklam{\dPs_{\arn}(\alpha)}^{-1}
	\end{align*} 
	für alle $n\in\Z{0}{\fklam{\kappa-2}}$ und im Fall $\kappa\geq4$
	\begin{align*}
		B_{\arn} &= \beklam{\dPs_{\arn-1}(\alpha)}^{-\ast}\Marn^{-1}\Marnp\beklam{\dPs_{\arn+1}(\alpha)}^{\ast}
	\end{align*}
	für alle $n\in\Z{1}{\fklam{\kappa-2}}$.
	\end{itemize}
\end{lemma}

\bwanf Zu (a): Aus \thref{fpdef1} folgt $B_0 = s_0$. Wegen Teil (a) von \thref{sqsa6} gilt
\begin{align*}
	A_0-\alpha\Iq = -(\alpha\Iq-A_0)P_0(\alpha) = -P_1(\alpha).
\end{align*}
Wegen Teil (b) von \thref{sqsa3} gilt
\begin{align}	\label{wzlm1bw1}
	\Ps_n(\alpha) = -P_n(\alpha)\sum^{n-1}_{j=0}\Larj
\end{align}
für alle $n\in\Zefkp$. 

Sei nun $\kappa\geq3$. Wegen Teil (b) von \thref{sqsa6} und \fref{wzlm1bw1} gilt dann
\begin{align*}
	\Oq &= \Ps_2(\alpha)-(\alpha\Iq-A_1)\Ps_1(\alpha)+B^{\ast}_1\Ps_0(\alpha) \\
	&= -P_2(\alpha)(\Laro+\Lare)+(\alpha\Iq-A_1)P_1(\alpha)\Laro
\end{align*}
und somit unter Beachtung von Teil (a) von \thref{sqsa3} und \thref{adpbm3} dann
\begin{align}	\label{wzlm1bw4}
	A_1-\alpha\Iq &= -P_2(\alpha)(\Laro+\Lare)\Laro^{-1}\beklam{P_1(\alpha)}^{-1}.
\end{align}
Wegen Teil (a) von \thref{sqsa6} gilt
\begin{align}	\label{wzlm1bw5}
	\Oq = P_{n+1}(\alpha)-(\alpha\Iq-A_n)P_n(\alpha)+B^{\ast}_nP_{n-1}(\alpha)
\end{align}
für alle $n\in\Zefkm$. Hieraus folgt
\begin{align}	\label{wzlm1bw3}
	\Oq = P_{n+1}(\alpha)\sum^{n-2}_{j=0}\Larj-(\alpha\Iq-A_n)P_n(\alpha)\sum^{n-2}_{j=0}\Larj+B^{\ast}_nP_{n-1}(\alpha)\sum^{n-2}_{j=0}\Larj
\end{align}
für alle $n\in\Zzfkm$. Wegen Teil (b) von \thref{sqsa6} und \fref{wzlm1bw1} gilt weiterhin
\begin{align*}
	\Oq &= \Ps_{n+1}(\alpha)-(\alpha\Iq-A_n)\Ps_n(\alpha)+B^{\ast}_n\Ps_{n-1}(\alpha) \\
	&= -P_{n+1}(\alpha)\sum^n_{j=0}\Larj+(\alpha\Iq-A_n)P_n(\alpha)\sum^{n-1}_{j=0}\Larj-B^{\ast}_nP_{n-1}(\alpha)\sum^{n-2}_{j=0}\Larj
\end{align*}
für alle $n\in\Zzfkm$. Hieraus folgt wegen \fref{wzlm1bw3} nun
\begin{align*}
	\Oq = -P_{n+1}(\alpha)(\Larnm+\Larn)+(\alpha\Iq-A_n)P_n(\alpha)\Larnm
\end{align*}
für alle $n\in\Zzfkm$ und somit unter Beachtung von Teil (a) von \thref{sqsa3} und \thref{adpbm3} dann
\begin{align*}
	A_n-\alpha\Iq &= -P_{n+1}(\alpha)(\Larnm+\Larn)\Larnm^{-1}\beklam{P_n(\alpha)}^{-1}
\end{align*}
für alle $n\in\Zzfkm$. Hieraus folgt wegen \fref{wzlm1bw5} und \fref{wzlm1bw4} weiterhin
\begin{align*}
	\Oq &= P_{n+1}(\alpha)-P_{n+1}(\alpha)(\Larnm+\Larn)\Larnm^{-1}+B^{\ast}_nP_{n-1}(\alpha) \\
	&= -P_{n+1}(\alpha)\Larn\Larnm^{-1}+B^{\ast}_nP_{n-1}(\alpha)
\end{align*}
für alle $n\in\Zefkm$ und somit unter Beachtung von Teil (a) von \thref{sqsa3} und \thref{adpbm3} dann
\begin{align*}
	B_n &= \beklam{P_{n-1}(\alpha)}^{-\ast}\Larnm^{-1}\Larn \beklam{P_{n+1}(\alpha)}^{\ast}
\end{align*}
für alle $n\in\Zefkm$. 

Zu (b): Aus \thref{fpdef1} folgt $B_{\aro} = s_{\aro}$. 

Sei nun $\kappa\geq2$. 
Wegen Teil (b) von \thref{sqsa3} gilt dann
\begin{align}	\label{wzlm1bw2}
	P_{\arn}(\alpha) = \dPs_{\arn}\sum^n_{j=0}\Marj
\end{align}
für alle $n\in\Zefk$. Hieraus folgt wegen \thref{adpdef1}, \thref{sqbm1} und Teil (c) von \thref{sqsa6} dann
\begin{align*}
	& -\dPs_{\ar1}(\alpha)(\Maro+\Mare)\Maro^{-1}\beklam{\dPs_{\aro}(\alpha)}^{-1} \\
	&= -P_{\ar1}(\alpha)s_0s^{-1}_0 
	= -(\alpha\Iq-A_{\aro})P_{\aro}(\alpha) \\
	&= A_{\aro}-\alpha\Iq.
\end{align*}

Sei nun $\kappa\geq4$. Wegen Teil (d) von \thref{sqsa6} gilt dann
\begin{align}	\label{wzlm1bw6}
	\Oq = \dPs_{\arn+1}(\alpha)-(\alpha\Iq-A_{\arn})\dPs_{\arn}(\alpha)+B^{\ast}_{\arn}\dPs_{\arn-1}(\alpha)
\end{align}
für alle $n\in\Z{1}{\fklam{\kappa-2}}$. Hieraus folgt
\begin{align}	\label{wzlm1bw7}
	\Oq = \dPs_{\arn+1}(\alpha)\sum^{n-1}_{j=0}\Marj-(\alpha\Iq-A_{\arn})\dPs_{\arn}(\alpha)\sum^{n-1}_{j=0}\Marj+B^{\ast}_{\arn}\dPs_{\arn-1}(\alpha)\sum^{n-1}_{j=0}\Marj
\end{align}
für alle $n\in\Z{1}{\fklam{\kappa-2}}$. Wegen Teil (c) von \thref{sqsa6} und \fref{wzlm1bw2} gilt weiterhin
\begin{align*}
	\Oq &= P_{\arn+1}(\alpha)-(\alpha\Iq-A_{\arn})P_{\arn}(\alpha)+B^{\ast}_{\arn}P_{\arn-1}(\alpha) \\
	&= \dPs_{\arn+1}(\alpha)\sum^{n+1}_{j=0}\Marj-(\alpha\Iq-A_{\arn})\dPs_{\arn}(\alpha)\sum^n_{j=0}\Marj+B^{\ast}_{\arn}\dPs_{\arn-1}(\alpha)\sum^{n-1}_{j=0}\Marj
\end{align*}
für alle $n\in\Z{1}{\fklam{\kappa-2}}$. Hieraus folgt wegen \fref{wzlm1bw7} nun
\begin{align*}
	\Oq = \dPs_{\arn+1}(\alpha)(\Marn+\Marnp)-(\alpha\Iq-A_{\arn})\dPs_{\arn}(\alpha)\Marn
\end{align*}
für alle $n\in\Z{1}{\fklam{\kappa-2}}$ und somit unter Beachtung von Teil (a) von \thref{sqsa3} und \thref{adpbm3} dann
\begin{align*}
	A_{\arn}-\alpha\Iq &= -\dPs_{\arn+1}(\alpha)(\Marn+\Marnp)\Marn^{-1}\beklam{\dPs_{\arn}(\alpha)}^{-1}
\end{align*}
für alle $n\in\Z{1}{\fklam{\kappa-2}}$. Hieraus folgt wegen \fref{wzlm1bw6} weiterhin
\begin{align*}
	\Oq &= \dPs_{\arn+1}(\alpha)-\dPs_{\arn+1}(\alpha)(\Marn+\Marnp)\Marn^{-1}+B^{\ast}_{\arn}\dPs_{\arn-1}(\alpha) \\
	&= -\dPs_{\arn+1}(\alpha)\Marnp\Marn^{-1}+B^{\ast}_{\arn}\dPs_{\arn-1}(\alpha)
\end{align*}
für alle $n\in\Z{1}{\fklam{\kappa-2}}$ und somit unter Beachtung von Teil (a) von \thref{sqsa3} und \thref{adpbm3} dann
\begin{align*}
	B_{\arn} &= \beklam{\dPs_{\arn-1}(\alpha)}^{-\ast}\Marn^{-1}\Marnp\beklam{\dPs_{\arn+1}(\alpha)}^{\ast}
\end{align*}
für alle $n\in\Z{1}{\fklam{\kappa-2}}$. \bwend

\begin{satz}	\thlabel{wzsa1}
	Seien $\kappa \in \Na$, $\alpha \in \R$, $\sjk \in \Kpqka$ und $\LMkarn$  die rechtsseitige $\alpha$-Dyukarev-Stieltjes-Parametrisierung von $\sjk$. \dgfa
	\begin{itemize}
	\item [\rm{(a)}] Sei $\ABfk$ das Favard-Paar bezüglich $\sjk$. Dann gelten \linebreak  $B_0 = \Maro^{-1}$, $A_0-\alpha\Iq = \Maro^{-1}\Laro^{-1}$, im Fall $\kappa\geq2$
	\begin{align*}
		B_1 &= \Laro^{-1}\Mare^{-1}\Laro^{-1}\Maro^{-1},
	\end{align*}
	im Fall $\kappa\geq3$
	\begin{align*}
		A_1-\alpha\Iq = \Maro^{-1}\Laro^{-1}\Mare^{-1}\Lare^{-1}(\Laro+\Lare)\Maro,
	\end{align*}
	im Fall $\kappa\geq4$
	\begin{align*}
		B_n = \prodr^{n-2}_{j=0}(\Marj\Larj)\Larnm^{-1}\Marn^{-1}\prodl^{n-1}_{j=0}\rklam{\Larj^{-1}\Marj^{-1}}
	\end{align*}
	für alle $n\in\Zzfk$ und im Fall $\kappa\geq5$
	\begin{align*}
		A_n-\alpha\Iq = \prodr^n_{j=0}\rklam{\Marj^{-1}\Larj^{-1}}(\Larnm+\Larn)\Marnm\prodl^{n-2}_{j=0}(\Larj\Marj)
	\end{align*}
	für alle $n\in\Zzfkm$. 
	\item [\rm{(b)}] Sei $(B_{\arn})^0_{n=0}$ bzw. im Fall $\kappa\geq2$ $\ABarfk$ das Favard-Paar bezüglich $\sarjk$. Dann gelten $B_{\aro} = \Maro^{-1}\Laro^{-1}\Maro^{-1}$, im Fall $\kappa\geq2$
	\begin{align*}
		A_{\aro}-\alpha\Iq = \Maro^{-1}\Laro^{-1}\Mare^{-1}(\Maro+\Mare),
	\end{align*}
	im Fall $\kappa\geq3$
	\begin{align*}
		B_{\ar1} = \Maro\Mare^{-1}\Lare^{-1}\Mare^{-1}\Laro^{-1}\Maro^{-1},
	\end{align*}
	im Fall $\kappa\geq4$
	\begin{align*}
		A_{\arn}-\alpha\Iq = \prodr^n_{j=0}\rklam{\Marj^{-1}\Larj^{-1}}\Marnp^{-1}(\Marn+\Marnp)\prodl^{n-1}_{j=0}(\Larj\Marj)
	\end{align*}
	für alle $n\in\Z{1}{\fklam{\kappa-2}}$ und im Fall $\kappa\geq5$
	\begin{align*}
		B_{\arn} = \prodr^{n-2}_{j=0}(\Marj\Larj)\Marnm\Marn^{-1}\prodl^n_{j=0}\rklam{\Larj^{-1}\Marj^{-1}}
	\end{align*}
	für alle $n\in\Zzfkm$.
	\end{itemize}
\end{satz}

\bwanf Sei $\SQark$ das rechtsseitige $\alpha$-Stieltjes-Qua\-drupel bezüglich $\sjk$.

Zu (a): Wegen \thref{adpdef1} und \thref{fpdef1} gelten
\begin{align*}
	\Maro^{-1} = s_0 = B_0
\end{align*}
und
\begin{align*}
	\Maro^{-1}\Laro^{-1} = s_0(s_0s^{-1}_{\aro}s_0)^{-1} = s_{\aro}s^{-1}_0 = (-\alpha s_0+s_1)s^{-1}_0 = -\alpha\Iq+A_0.
\end{align*}

Sei nun $\kappa\geq2$. Unter Beachtung von \thref{adpbm3} gilt wegen Teil (a) von \thref{wzlm1}, \thref{sqbm1} und Teil (b) von \thref{sqsa3} dann
\begin{align*}
	B_1 &= \beklam{P_{0}(\alpha)}^{-\ast}\Laro^{-1}\Lare \beklam{P_{2}(\alpha)}^{\ast}
	= \Laro^{-1}\Lare\rklam{\Maro^{-1}\Laro^{-1}\Mare^{-1}\Lare^{-1}}^{\ast} \\
	&= \Laro^{-1}\Mare^{-1}\Laro^{-1}\Maro^{-1}
\end{align*}
(im Fall $\kappa=2$ existiert wegen Teil (c) von \thref{asmsa1} ein $s_3\in\Cqq$ mit $(s_j)^3_{j=0}\in{\cal K}^{>}_{q,3,\alpha}$, so dass die Formel für $B_1$ in Teil (a) von \thref{wzlm1} bzw. für $P_2(\alpha)$ in Teil (b) von \thref{sqsa3} benutzt werden kann; am Ende der Rechnung spielt durch das Wegfallen von $\Lare$ dann die Matrix $s_3$ keine Rolle mehr). 

Sei nun $\kappa\geq3$. Wegen Teil (a) von \thref{wzlm1} und Teil (b) von \thref{sqsa3} gilt dann
\begin{align*}
	A_1-\alpha\Iq &= -P_{2}(\alpha)(\Laro+\Lare)\Laro^{-1}\beklam{P_1(\alpha)}^{-1} \\
	&= -\Maro^{-1}\Laro^{-1}\Mare^{-1}\Lare^{-1}(\Laro+\Lare)\Laro^{-1}\rklam{-\Maro^{-1}\Laro^{-1}}^{-1} \\
	&= \Maro^{-1}\Laro^{-1}\Mare^{-1}\Lare^{-1}(\Laro+\Lare)\Maro.
\end{align*}

Sei nun $\kappa\geq4$. Unter Beachtung von \thref{adpbm3} gilt wegen Teil (a) von \thref{wzlm1} und Teil (b) von \thref{sqsa3} dann
\begin{align*}
	B_n &= \beklam{P_{n-1}(\alpha)}^{-\ast}\Larnm^{-1}\Larn \beklam{P_{n+1}(\alpha)}^{\ast} \\
	&= \eklam{(-1)^{n-1}\prodr^{n-2}_{j=0}\rklam{\Marj^{-1}\Larj^{-1}}}^{-\ast}\Larnm^{-1}\Larn\eklam{(-1)^{n+1}\prodr^n_{j=0}\rklam{\Marj^{-1}\Larj^{-1}}}^{\ast} \\
	&= \prodr^{n-2}_{j=0}(\Marj\Larj)\Larnm^{-1}\Marn^{-1}\prodl^{n-1}_{j=0}\rklam{\Larj^{-1}\Marj^{-1}}
\end{align*}
für alle $n\in\Zzfk$ (im Fall $2n=\kappa$ existiert wegen Teil (c) von \thref{asmsa1} ein $s_{2n+1}\in\Cqq$ mit $(s_j)^{2n+1}_{j=0}\in{\cal K}^{>}_{q,2n+1,\alpha}$, so dass die Formel für $B_n$ in Teil (a) von \thref{wzlm1} bzw. für $P_{n+1}(\alpha)$ in Teil (b) von \thref{sqsa3} benutzt werden kann; am Ende der Rechnung spielt durch das Wegfallen von $\Larn$ dann die Matrix $s_{2n+1}$ keine Rolle mehr). 

Sei nun $\kappa\geq5$. Wegen Teil (a) von \thref{wzlm1} und Teil (b) von \thref{sqsa3} gilt dann
\begin{align*}
	A_n-\alpha\Iq &= -P_{n+1}(\alpha)(\Larnm+\Larn)\Larnm^{-1}\beklam{P_n(\alpha)}^{-1} \\
	&= (-1)^{n+2}\prodr^n_{j=0}\rklam{\Marj^{-1}\Larj^{-1}}(\Larnm+\Larn)\Larnm^{-1}\eklam{(-1)^{n}\prodr^{n-1}_{j=0}\rklam{\Marj^{-1}\Larj^{-1}}}^{-1} \\
	&= \prodr^n_{j=0}\rklam{\Marj^{-1}\Larj^{-1}}(\Larnm+\Larn)\Marnm\prodl^{n-2}_{j=0}(\Larj\Marj)	
\end{align*}
für alle $n\in\Zzfkm$.

Zu (b): Wegen \thref{adpdef1} und \thref{fpdef1} gilt
\begin{align*}
	\Maro^{-1}\Laro^{-1}\Maro^{-1} = s_0(s_0s^{-1}_{\aro}s_0)^{-1}s_0 = s_{\aro} = B_{\aro}.
\end{align*}

Sei nun $\kappa\geq2$. Wegen Teil (b) von \thref{wzlm1}, Teil (b) von \thref{sqsa3}, \thref{adpdef1} und \thref{sqbm1} gilt dann
\begin{align*}
	A_{\aro}-\alpha\Iq &= -\dPs_{\ar1}(\alpha)(\Maro+\Mare)\Maro^{-1}\beklam{\dPs_{\aro}(\alpha)}^{-1} \\
	&= \Maro^{-1}\Laro^{-1}\Mare^{-1}(\Maro+\Mare)s_0s^{-1}_0 \\
	&= \Maro^{-1}\Laro^{-1}\Mare^{-1}(\Maro+\Mare).
\end{align*}

Sei nun $\kappa\geq3$. Unter Beachtung von \thref{adpbm3} gilt wegen Teil (b) von \thref{wzlm1}, \thref{sqbm1}, \thref{adpdef1} und Teil (b) von \thref{sqsa3} dann
\begin{align*}
	B_{\ar1} &= \beklam{\dPs_{\aro}(\alpha)}^{-\ast}\Mare^{-1}\Marz\beklam{\dPs_{\ar2}(\alpha)}^{\ast} 
	= s^{-\ast}_0\Mare^{-1}\Marz\rklam{\Maro^{-1}\Laro^{-1}\Mare^{-1}\Lare^{-1}\Marz^{-1}}^{\ast} \\
	&= \Maro\Mare^{-1}\Lare^{-1}\Mare^{-1}\Laro^{-1}\Maro^{-1}
\end{align*}
(im Fall $\kappa=3$ existiert wegen Teil (c) von \thref{asmsa1} ein $s_4\in\Cqq$ mit $(s_j)^4_{j=0}\in{\cal K}^{>}_{q,4,\alpha}$, so dass die Formel für $B_{\ar1}$ in Teil (b) von \thref{wzlm1} bzw. für $\dPs_{\ar2}(\alpha)$ in Teil (b) von \thref{sqsa3} benutzt werden kann; am Ende der Rechnung spielt durch das Wegfallen von $\Marz$ dann die Matrix $s_4$ keine Rolle mehr). 

Sei nun $\kappa\geq4$. Wegen Teil (b) von \thref{wzlm1} und Teil (b) von \thref{sqsa3} gilt dann
\begin{align*}
	A_{\arn}-\alpha\Iq &= -\dPs_{\arn+1}(\alpha)(\Marn+\Marnp)\Marn^{-1}\beklam{\dPs_{\arn}(\alpha)}^{-1} \\
	&= (-1)^{n+2}\prodr^n_{j=0}\rklam{\Marj^{-1}\Larj^{-1}}\Marnp^{-1}(\Marn+\Marnp)\Marn^{-1} \\ &\quad \cdot\eklam{(-1)^n\prodr^{n-1}_{j=0}\rklam{\Marj^{-1}\Larj^{-1}}\Marn^{-1}}^{-1} \\
	&= \prodr^n_{j=0}\rklam{\Marj^{-1}\Larj^{-1}}\Marnp^{-1}(\Marn+\Marnp)\prodl^{n-1}_{j=0}(\Larj\Marj)
\end{align*}
für alle $n\in\Z{1}{\fklam{\kappa-2}}$. 

Sei nun $\kappa\geq5$. Unter Beachtung von \thref{adpbm3} gilt wegen Teil (b) von \thref{wzlm1} und Teil (b) von \thref{sqsa3} dann
\begin{align*}
	B_{\arn} &= \beklam{\dPs_{\arn-1}(\alpha)}^{-\ast}\Marn^{-1}\Marnp\beklam{\dPs_{\arn+1}(\alpha)}^{\ast} \\
	&= \eklam{(-1)^{n-1}\prodr^{n-2}_{j=0}\rklam{\Marj^{-1}\Larj^{-1}}\Marnm^{-1}}^{-\ast}\Marn^{-1}\Marnp \\ &\quad \cdot\eklam{(-1)^{n+1}\prodr^n_{j=0}\rklam{\Marj^{-1}\Larj^{-1}}\Marnp^{-1}}^{\ast} \\
	&= \prodr^{n-2}_{j=0}(\Marj\Larj)\Marnm\Marn^{-1}\prodl^n_{j=0}\rklam{\Larj^{-1}\Marj^{-1}}
\end{align*}
für alle $n\in\Zzfkm$ (im Fall $2n+1=\kappa$ existiert wegen Teil (c) von \thref{asmsa1} ein $s_{2n+2}\in\Cqq$ mit $(s_j)^{2n+2}_{j=0}\in{\cal K}^{>}_{q,2n+2,\alpha}$, so dass die Formel für $B_{\arn}$ in Teil (b) von \thref{wzlm1} bzw. für $\dPs_{\arn+1}(\alpha)$ in Teil (b) von \thref{sqsa3} benutzt werden kann; am Ende der Rechnung spielt durch das Wegfallen von $\Marnp$ dann die Matrix $s_{2n+2}$ keine Rolle mehr). \bwend

Nun zeigen wir einen Zusammenhang zwischen der rechtsseitigen $\alpha$-Stieltjes"=Parametrisierung einer rechtsseitig $\alpha$-Stieltjes-positiv definiten Folge und den Favard-Paaren bezüglich jener Folge und deren durch rechtsseitige $\alpha$-Verschiebung generierten Folge (vergleiche Teil (a) mit \cite[Proposition 9.4]{CR1} für den Fall $\alpha=0$ und $\kappa=\infty$).

\begin{satz}	\thlabel{wzsa2}
	Seien $\kappa \in \Na$, $\alpha \in \R$, $\sjk \in \Kpqka$ und $\Qarjk$ die rechtsseitige $\alpha$-Stieltjes-Parametrisierung von $\sjk$. \dgfa
	\begin{itemize}
		\item [\rm{(a)}] Sei $\ABfk$ das Favard-Paar bezüglich $\sjk$. Dann gelten \linebreak $B_0 = Q_{\aro}$, $A_0-\alpha\Iq = Q_{\ar1}Q^{-1}_{\aro}$, im Fall $\kappa\geq2$
		\begin{align*}
			B_n = Q^{-1}_{\ar2n-2}Q_{\ar2n}
		\end{align*}
		für alle $n\in\Zefk$ und im Fall $\kappa\geq3$
		\begin{align*}
			A_n-\alpha\Iq = Q_{\ar2n+1}Q^{-1}_{\ar2n}+Q_{\ar2n}Q^{-1}_{\ar2n-1}
		\end{align*}
		für alle $n\in\Zefkm$.
		\item [\rm{(b)}] Sei $(B_{\arn})^0_{n=0}$ bzw. im Fall $\kappa\geq2$ $\ABarfk$ das Favard-Paar bezüglich $\sarjk$. Dann gelten $B_{\aro} = Q_{\ar1}$, im Fall $\kappa\geq2$
		\begin{align*}
			A_{\arn}-\alpha\Iq = Q_{\ar2n+2}Q^{-1}_{\ar2n+1}+Q_{\ar2n+1}Q^{-1}_{\ar2n}
		\end{align*}
		für alle $n\in\Z{0}{\fklam{\kappa-2}}$ und im Fall $\kappa\geq3$
		\begin{align*}
			B_{\arn} = Q^{-1}_{\ar2n-1}Q_{\ar2n+1}
		\end{align*}
		für alle $n\in\Zefkm$.
	\end{itemize}
\end{satz}

\bwanf Seien $\SQark$ das rechtsseitige $\alpha$-Stieltjes-Quadrupel bezüglich $\sjk$ und $\LMkarn$ die rechtsseitige $\alpha$-Dyukarev-Stieltjes-Parametrisierung von $\sjk$.

Zu (a): Wegen \thref{fpdef1} und Teil (a) von \thref{aspdef2} gelten
\begin{align*}
	B_0 = s_0 = \dH_0 = Q_{\aro}
\end{align*}
und
\begin{align*}
	A_0-\alpha\Iq = s_1s^{-1}_0-\alpha s_0s^{-1}_0 
	= s_{\aro}s^{-1}_0 = \dH_{\aro}\dH^{-1}_0 = Q_{\ar1}Q^{-1}_{\aro}.
\end{align*}

Sei nun $\kappa\geq2$. Wegen \thref{fpdef1} und Teil (a) von \thref{aspdef2} gilt dann
\begin{align*}
	B_n = \dH^{-1}_{n-1}\dH_{n} = Q^{-1}_{\ar2n-2}Q_{\ar2n}
\end{align*}
für alle $n\in\Zefk$. 

Sei nun $\kappa\geq3$. Wegen \thref{sqbm1}, \thref{adpdef1} und Teil (b) von \thref{sqsa3} gelten
\begin{align*}
	\dPs_{\aro}(\alpha) = s_0 = \Maro^{-1} = -P_1(\alpha)\Laro
\end{align*}
und
\begin{align*}
	\dPs_{\arn}(\alpha) = -P_{n+1}(\alpha)\Larn
\end{align*}
für alle $n\in\Zefkm$. Hieraus folgt wegen Teil (a) von \thref{wzlm1} und Teil (b) von \thref{sqsa5} nun
\begin{align*}
	A_n-\alpha\Iq &= -P_{n+1}(\alpha)(\Larnm+\Larn)\Larnm^{-1}\beklam{P_n(\alpha)}^{-1} \\
	&= -P_{n+1}(\alpha)\beklam{P_n(\alpha)}^{-1}-P_{n+1}(\alpha)\Larn\Larnm^{-1}\beklam{P_n(\alpha)}^{-1} \\
	&= -P_{n+1}(\alpha)\beklam{P_n(\alpha)}^{-1}-\dPs_{\arn}(\alpha)\beklam{\dPs_{\arn-1}(\alpha)}^{-1} \\
	&= \prodl^n_{j=0}\rklam{Q_{\ar2j+1}Q^{-1}_{\ar2j}}\prodr^{n-1}_{j=0}\rklam{Q_{\ar2j}Q^{-1}_{\ar2j+1}} \\ &\quad +Q_{\ar2n}\prodl^{n-1}_{j=0}\rklam{Q^{-1}_{\ar2j+1}Q_{\ar2j}}\prodr^{n-2}_{j=0}\rklam{Q^{-1}_{\ar2j}Q_{\ar2j+1}}Q^{-1}_{\ar2n-2} \\
	&= Q_{\ar2n+1}Q^{-1}_{\ar2n}+Q_{\ar2n}Q^{-1}_{\ar2n-1}
\end{align*}
für alle $n\in\Zefkm$.

Zu (b): Wegen \thref{fpdef1} und Teil (a) von \thref{aspdef2} gilt
\begin{align*}
	B_{\aro} = s_{\aro} = \dH_{\aro} = Q_{\ar1}.
\end{align*}

Sei nun $\kappa\geq2$. Wegen \thref{sqbm1} und \thref{adpdef1} gilt dann
\begin{align*}
	P_0(\alpha) = \Iq = s_0s^{-1}_0 = \dPs_{\aro}(\alpha)\Maro
\end{align*}
und wegen Teil (b) von \thref{sqsa3} gilt weiterhin
\begin{align*}
	P_n(\alpha) = \dPs_{\arn}\Marn
\end{align*}
für alle $n\in\Zefk$. Hieraus folgt wegen Teil (b) von \thref{wzlm1} und Teil (b) von \thref{sqsa5} nun
\begin{align*}
	A_{\arn}-\alpha\Iq &= -\dPs_{\arn+1}(\alpha)(\Marn+\Marnp)\Marn^{-1}\beklam{\dPs_{\arn}(\alpha)}^{-1} \\\
	&= -\dPs_{\arn+1}(\alpha)\beklam{\dPs_{\arn}(\alpha)}^{-1}-\dPs_{\arn+1}(\alpha)\Marnp\Marn^{-1}\beklam{\dPs_{\arn}(\alpha)}^{-1} \\
	&= -\dPs_{\arn+1}(\alpha)\beklam{\dPs_{\arn}(\alpha)}^{-1}-P_{n+1}(\alpha)\beklam{P_n(\alpha)}^{-1} \\
	&= Q_{\ar2n+2}\prodl^n_{j=0}\rklam{Q^{-1}_{\ar2j+1}Q_{\ar2j}}\prodr^{n-1}_{j=0}\rklam{Q^{-1}_{\ar2j}Q_{\ar2j+1}}Q^{-1}_{\ar2n} \\ &\quad +\prodl^n_{j=0}\rklam{Q_{\ar2j+1}Q^{-1}_{\ar2j}}\prodr^{n-1}_{j=0}\rklam{Q_{\ar2j}Q^{-1}_{\ar2j+1}} \\
	&= Q_{\ar2n+2}Q^{-1}_{\ar2n+1}+Q_{\ar2n+1}Q^{-1}_{\ar2n}
\end{align*}
für alle $n\in\Z{0}{\fklam{\kappa-2}}$. 

Sei nun $\kappa\geq3$. Wegen \thref{fpdef1} und Teil (a) von \thref{aspdef2} gilt dann
\begin{align*}
	B_{\arn} = \dH^{-1}_{\arn-1}\dH_{\arn} = Q^{-1}_{\ar2n-1}Q_{\ar2n+1}
\end{align*}
für alle $n\in\Zefkm$. \bwend

\subsection{Der linksseitige Fall}

Für unser weiteres Vorgehen benötigen wir folgendes Lemma, das die Favard-Paare für den linksseitigen und rechtsseitigen Fall verbindet. Somit können wir im linksseitigen Fall auf entsprechende Resultate für den rechtsseitgen Fall zurückgreifen.

\begin{lemma}	\thlabel{wzllm1}
	Seien $\kappa \in \Noa$, $\sjk$ eine Folge aus $\Cqq$ und $t_j:=(-1)^js_j$ für alle $j\in\Zok$. \dgfa
	\begin{itemize}
		\item [\rm{(a)}] Sei im Fall $\kappa\geq1$ $\sjfkm\in\Hpqfkm$ oder $(t_j)^{2\fklam{\kappa-1}}_{j=0}\in\Hpqfkm$ erfüllt. Dann gilt $(t_j)^{2\fklam{\kappa-1}}_{j=0}\in\Hpqfkm$ bzw. $\sjfkm\in\Hpqfkm$, falls $\kappa\geq1$. Weiterhin seien $(B^{\sklam{s}}_{n})^0_{n=0}$ bzw. im Fall $\kappa\geq1$ $\ABsfk$ das Favard-Paar bezüglich $\sjk$ und $(B^{\sklam{t}}_{n})^0_{n=0}$ bzw. im Fall $\kappa\geq1$ $[(A^{\sklam{t}}_{n})^{\fklam{\kappa-1}}_{n=0}$,""$(B^{\sklam{t}}_{n})^{\fklam{\kappa}}_{n=0}]$ das Favard-Paar bezüglich $\tjk$. Dann gelten
		\begin{align*}
			B^{\sklam{t}}_{n} = B^{\sklam{s}}_{n}
		\end{align*}
		für alle $n\in\Zofk$ und im Fall $\kappa\geq1$
		\begin{align*}
			A^{\sklam{t}}_{n} = -A^{\sklam{s}}_{n}
		\end{align*}
		für alle $n\in\Zofkm$.
		\item [\rm{(b)}] Seien nun $\kappa\geq1$, $\alpha\in\C$ und im Fall $\kappa\geq2$ $(s_{\arj})^{2\fklam{\kappa-2}}_{j=0}\in\Hpqfkzm$ oder $(t_{-\alj})^{2\fklam{\kappa-2}}_{j=0}\in\Hpqfkzm$ erfüllt. Dann gilt $(t_{-\alj})^{2\fklam{\kappa-2}}_{j=0}\in\Hpqfkzm$ bzw. \linebreak $(s_{\arj})^{2\fklam{\kappa-2}}_{j=0}\in\Hpqfkzm$, falls $\kappa\geq2$. Weiterhin seien $(B^{\sklam{s}}_{\arn})^0_{n=0}$ bzw. im Fall $\kappa\geq2$ $\ABsarfk$ das Favard-Paar bezüglich $\sarjk$ und $(B^{\sklam{t}}_{-\aln})^0_{n=0}$ bzw. im Fall $\kappa\geq2$ $[(A^{\sklam{t}}_{-\aln})^{\fklam{\kappa-2}}_{n=0}$,""$(B^{\sklam{t}}_{-\aln})^{\fklam{\kappa-1}}_{n=0}]$ das Favard-Paar bezüglich $(t_{-\alj})^{\kappa-1}_{j=0}$. Dann gelten
		\begin{align*}
			B^{\sklam{t}}_{-\aln} = B^{\sklam{s}}_{\arn}
		\end{align*}
		für alle $n\in\Zofkm$ und im Fall $\kappa\geq2$
		\begin{align*}
			A^{\sklam{t}}_{-\aln} = -A^{\sklam{s}}_{\arn}
		\end{align*}
		für alle $n\in\Z{0}{\fklam{\kappa-2}}$.
	\end{itemize}
\end{lemma}

\bwanf Zu (a): Wegen Teil (a) von \thref{asmdef1} und der Teile (a) und (b) von \thref{asmlm1} gelten im Fall $\kappa\geq1$ $(t_j)^{2\fklam{\kappa-1}}_{j=0}\in\Hpqfkm$ und $\sjfkm\in\Hpqfkm$. Wegen \thref{fpdef1} gelten weiterhin
\begin{align*}
	B^{\sklam{t}}_{0} &= t_0 = s_0 = B^{\sklam{s}}_{0}
\end{align*}
und im Fall $\kappa\geq1$
\begin{align*}
	A^{\sklam{t}}_{0} &= t_1t^{-1}_0 = -s_1s^{-1}_0 = -A^{\sklam{s}}_{0}.
\end{align*}

Sei nun $\kappa\geq2$. Wegen \thref{fpdef1} und Teil (e) von \thref{asmlm1} gilt dann
\begin{align*}
	B^{\sklam{t}}_{n} = \brklam{\dH^{\sklam{t}}_{n-1}}^{-1}\dH^{\sklam{t}}_{n} 
	= \brklam{\dH^{\sklam{s}}_{n-1}}^{-1}\dH^{\sklam{s}}_{n} = B^{\sklam{s}}_{n}
\end{align*}
für alle $n\in\Zefk$. 

Sei nun $\kappa\geq3$. Wegen \thref{fpdef1} und \thref{asmlm1} gilt dann
\begin{align*}
	A^{\sklam{t}}_{n} &= \begin{pmatrix} -z^{\sklam{t}}_{n,2n-1}\brklam{\Hnm^{\sklam{t}}}^{-1} & \Iq \end{pmatrix} \Kn^{\sklam{t}} \begin{pmatrix} -\brklam{\Hnm^{\sklam{t}}}^{-1}y^{\sklam{t}}_{n,2n-1} \\ \Iq \end{pmatrix} \big(\dHn^{\sklam{t}}\big)^{-1} \\
	&= \begin{pmatrix} -(-1)^nz^{\sklam{s}}_{n,2n-1}\Vnma\Vnm\brklam{\Hsnm}^{-1}\Vnma & \Iq \end{pmatrix} \brklam{-\Vn\Ksn\Vna} \\
	& \quad \cdot\begin{pmatrix} -\Vnm\brklam{\Hsnm}^{-1}\Vnma(-1)^n\Vnm y^{\sklam{s}}_{n,2n-1} \\ \Iq \end{pmatrix} \big(\dHn^{\sklam{s}}\big)^{-1} \\
	&= (-1)^{2n+1}\begin{pmatrix} -z^{\sklam{s}}_{n,2n-1}\Vnma\Vnm\brklam{\Hsnm}^{-1} & \Iq \end{pmatrix} \Vna\Vn\Ksn\Vna\Vn \\
	& \quad \cdot\begin{pmatrix} -\brklam{\Hsnm}^{-1}\Vnma\Vnm y^{\sklam{s}}_{n,2n-1} \\ \Iq \end{pmatrix} \big(\dHn^{\sklam{s}}\big)^{-1} \\
	&= -\begin{pmatrix} -z^{\sklam{s}}_{n,2n-1}\brklam{\Hsnm}^{-1} & \Iq \end{pmatrix} \Ksn \begin{pmatrix} -\brklam{\Hsnm}^{-1}y^{\sklam{s}}_{n,2n-1} \\ \Iq \end{pmatrix} \big(\dHn^{\sklam{s}}\big)^{-1}
	= -A^{\sklam{s}}_{n}
\end{align*}
für alle $n\in\Zefkm$.

Zu (b): Dies folgt wegen Teil (a) von \thref{asplm1} aus (a). \bwend

Wir behandeln zuerst einen Zusammenhang zwischen der linksseitigen $\alpha$"=Dyukarev"=Stieltjes"=Paremetrisierung einer linksseitig $\alpha$-Stieltjes-positiv definiten Folge und den Favard-Paaren bezüglich jener Folge und deren durch linksseitige $\alpha$-Verschiebung generierten Folge.

\begin{satz}	\thlabel{wzlsa1}
	Seien $\kappa \in \Na$, $\alpha \in \R$, $\sjk \in \Lpqka$ und $\LMkaln$  die linksseitige $\alpha$-Dyukarev-Stieltjes-Parametrisierung von $\sjk$. \dgfa
	\begin{itemize}
	\item [\rm{(a)}] Sei $\ABfk$ das Favard-Paar bezüglich $\sjk$. Dann gelten \linebreak $B_0 = \Malo^{-1}$, $A_0-\alpha\Iq = -\Malo^{-1}\Lalo^{-1}$, im Fall $\kappa\geq2$
	\begin{align*}
		B_1 &= \Lalo^{-1}\Male^{-1}\Lalo^{-1}\Malo^{-1},
	\end{align*}
	im Fall $\kappa\geq3$
	\begin{align*}
		A_1-\alpha\Iq = -\Malo^{-1}\Lalo^{-1}\Male^{-1}\Lale^{-1}(\Lalo+\Lale)\Malo,
	\end{align*}
	im Fall $\kappa\geq4$
	\begin{align*}
		B_n = \prodr^{n-2}_{j=0}(\Malj\Lalj)\Lalnm^{-1}\Maln^{-1}\prodl^{n-1}_{j=0}\rklam{\Lalj^{-1}\Malj^{-1}}
	\end{align*}
	für alle $n\in\Zzfk$ und im Fall $\kappa\geq5$
	\begin{align*}
		A_n-\alpha\Iq = -\prodr^n_{j=0}\rklam{\Malj^{-1}\Lalj^{-1}}(\Lalnm+\Laln)\Malnm\prodl^{n-2}_{j=0}(\Lalj\Malj)
	\end{align*}
	für alle $n\in\Zzfkm$. 
	\item [\rm{(b)}] Sei $(B_{\aln})^0_{n=0}$ bzw. im Fall $\kappa\geq2$ $\ABalfk$ das Favard-Paar bezüglich $\saljk$. Dann gelten $B_{\alo} = \Malo^{-1}\Lalo^{-1}\Malo^{-1}$, im Fall $\kappa\geq2$
	\begin{align*}
		A_{\alo}-\alpha\Iq = -\Malo^{-1}\Lalo^{-1}\Male^{-1}(\Malo+\Male),
	\end{align*}
	im Fall $\kappa\geq3$
	\begin{align*}
		B_{\al1} = \Malo\Male^{-1}\Lale^{-1}\Male^{-1}\Lalo^{-1}\Malo^{-1},
	\end{align*}
	im Fall $\kappa\geq4$
	\begin{align*}
		A_{\aln}-\alpha\Iq = -\prodr^n_{j=0}\rklam{\Malj^{-1}\Lalj^{-1}}\Malnp^{-1}(\Maln+\Malnp)\prodl^{n-1}_{j=0}(\Lalj\Malj)
	\end{align*}
	für alle $n\in\Z{1}{\fklam{\kappa-2}}$ und im Fall $\kappa\geq5$
	\begin{align*}
		B_{\aln} = \prodr^{n-2}_{j=0}(\Malj\Lalj)\Malnm\Maln^{-1}\prodl^n_{j=0}\rklam{\Lalj^{-1}\Malj^{-1}}
	\end{align*}
	für alle $n\in\Zzfkm$.
	\end{itemize}
\end{satz}

\bwanf Sei $t_j := (-1)^js_j$ für alle $j\in\Zok$. Wegen Teil (a) von \thref{asmbm3} gilt dann $\tjk\in\Kpqkma$. Weiterhin sei $[(\mathbf{L}^{\sklam{t}}_{-\arn})^{\fklam{\kappa-1}}_{n=0}$,""$(\mathbf{M}^{\sklam{t}}_{-\arn})^{\fklam{\kappa}}_{n=0}]$ die rechtsseitige $-\alpha$-Dyukarev-Stieltjes-Parametrisierung von $\tjk$.

Zu (a): Sei $[(A^{\sklam{t}}_{n})^{\fklam{\kappa-1}}_{n=0}$,""$(B^{\sklam{t}}_{n})^{\fklam{\kappa}}_{n=0}]$ das Favard-Paar von $\tjk$. Wegen Teil (a) von \thref{wzllm1}, Teil (a) von \thref{wzsa1} und Teil (c) von \thref{adpllm1} gelten dann
\begin{align*}
	B^{\sklam{s}}_0 &= B^{\sklam{t}}_0 = \brklam{\Mmaro^{\sklam{t}}}^{-1} = \brklam{\Malo^{\sklam{s}}}^{-1}, \\
	A^{\sklam{s}}_0-\alpha\Iq &= -\eklam{A^{\sklam{t}}_0-(-\alpha)\Iq} = -\brklam{\Mmaro^{\sklam{t}}}^{-1}\brklam{\Lmaro^{\sklam{t}}}^{-1} = -\brklam{\Malo^{\sklam{s}}}^{-1}\brklam{\Lalo^{\sklam{s}}}^{-1},
\end{align*}
im Fall $\kappa\geq2$
\begin{align*}
	B^{\sklam{s}}_1 &= B^{\sklam{t}}_1 = \brklam{\Lmaro^{\sklam{t}}}^{-1}\brklam{\Mmare^{\sklam{t}}}^{-1}\brklam{\Lmaro^{\sklam{t}}}^{-1}\brklam{\Mmaro^{\sklam{t}}}^{-1} \\
	&= \brklam{\Lalo^{\sklam{s}}}^{-1}\brklam{\Male^{\sklam{s}}}^{-1}\brklam{\Lalo^{\sklam{s}}}^{-1}\brklam{\Malo^{\sklam{s}}}^{-1},
\end{align*}
im Fall $\kappa\geq3$
\begin{align*}
	A^{\sklam{s}}_1-\alpha\Iq &= -\eklam{A^{\sklam{t}}_1-(-\alpha)\Iq} \\
	&= -\brklam{\Mmaro^{\sklam{t}}}^{-1}\brklam{\Lmaro^{\sklam{t}}}^{-1}\brklam{\Mmare^{\sklam{t}}}^{-1}\brklam{\Lmare^{\sklam{t}}}^{-1}\brklam{\Lmaro^{\sklam{t}}+\Lmare^{\sklam{t}}}\Mmaro^{\sklam{t}} \\
	&= -\brklam{\Malo^{\sklam{s}}}^{-1}\brklam{\Lalo^{\sklam{s}}}^{-1}\brklam{\Male^{\sklam{s}}}^{-1}\brklam{\Lale^{\sklam{s}}}^{-1}\brklam{\Lalo^{\sklam{s}}+\Lale^{\sklam{s}}}\Malo^{\sklam{s}},
\end{align*}
im Fall $\kappa\geq4$
\begin{align*}
	B^{\sklam{s}}_n &= B^{\sklam{t}}_n  
	= \prodr^{n-2}_{j=0}\beklam{\Mmarj^{\sklam{t}}\Lmarj^{\sklam{t}}}\brklam{\Lmarnm^{\sklam{t}}}^{-1}\brklam{\Mmarn^{\sklam{t}}}^{-1}\prodl^{n-1}_{j=0}\eklam{\brklam{\Lmarj^{\sklam{t}}}^{-1}\brklam{\Mmarj^{\sklam{t}}}^{-1}}	\\
	&= \prodr^{n-2}_{j=0}\beklam{\Malj^{\sklam{s}}\Lalj^{\sklam{s}}}\brklam{\Lalnm^{\sklam{s}}}^{-1}\brklam{\Maln^{\sklam{s}}}^{-1}\prodl^{n-1}_{j=0}\eklam{(\Lalj^{\sklam{s}})^{-1}\brklam{\Malj^{\sklam{s}}}^{-1}}
\end{align*}
für alle $n\in\Zzfk$ und im Fall $\kappa\geq5$
\begin{align*}
	&\ A^{\sklam{s}}_n-\alpha\Iq = -\beklam{A^{\sklam{t}}_n-(-\alpha)\Iq} \\
	&= -\prodr^n_{j=0}\eklam{\brklam{\Mmarj^{\sklam{t}}}^{-1}\brklam{\Lmarj^{\sklam{t}}}^{-1}}\rklam{\Lmarnm^{\sklam{t}}+\Lmarn^{\sklam{t}}}\Mmarnm^{\sklam{t}}\prodl^{n-2}_{j=0}\beklam{\Lmarj^{\sklam{t}}\Mmarj^{\sklam{t}}} \\
	&= -\prodr^n_{j=0}\eklam{\brklam{\Malj^{\sklam{s}}}^{-1}\brklam{\Lalj^{\sklam{s}}}^{-1}}\brklam{\Lalnm^{\sklam{s}}+\Laln^{\sklam{s}}}\Malnm^{\sklam{s}}\prodl^{n-2}_{j=0}\beklam{\Lalj^{\sklam{s}}\Malj^{\sklam{s}}}
\end{align*}
für alle $n\in\Zzfkm$. 

Zu (b): Sei $(B^{\sklam{t}}_{-\arj})^0_{j=0}$ bzw. im Fall $\kappa\geq2$ $[(A^{\sklam{t}}_{-\arn})^{\fklam{\kappa-2}}_{n=0}$,""$(B^{\sklam{t}}_{-\arn})^{\fklam{\kappa-1}}_{n=0}]$ das Favard-Paar bezüglich $(t_{-\arj})^{\kappa-1}_{j=0}$. Wegen Teil (b) von \thref{wzllm1}, Teil (b) von \thref{wzsa1} und Teil (c) von \thref{adpllm1} gelten dann
\begin{align*}
	B^{\sklam{s}}_{\alo} = B^{\sklam{t}}_{-\aro}
	= \brklam{\Mmaro^{\sklam{t}}}^{-1}\brklam{\Lmaro^{\sklam{t}}}^{-1}\brklam{\Mmaro^{\sklam{t}}}^{-1}
	= \brklam{\Malo^{\sklam{s}}}^{-1}\brklam{\Lalo^{\sklam{s}}}^{-1}\brklam{\Malo^{\sklam{s}}}^{-1}, 
\end{align*}
im Fall $\kappa\geq2$
\begin{align*}
	A^{\sklam{s}}_{\alo}-\alpha\Iq &= -\beklam{A^{\sklam{t}}_{\alo}-(-\alpha)\Iq} \\
	&= -\brklam{\Mmaro^{\sklam{t}}}^{-1}\brklam{\Lmaro^{\sklam{t}}}^{-1}\brklam{\Mmare^{\sklam{t}}}^{-1}\brklam{\Mmaro^{\sklam{t}}+\Mmare^{\sklam{t}}} \\
	&= -\brklam{\Malo^{\sklam{s}}}^{-1}\brklam{\Lalo^{\sklam{s}}}^{-1}\brklam{\Male^{\sklam{s}}}^{-1}\brklam{\Malo^{\sklam{s}}+\Male^{\sklam{s}}},
\end{align*}
im Fall $\kappa\geq3$
\begin{align*}
	B^{\sklam{s}}_{\al1} &= B^{\sklam{t}}_{-\ar1} 
	= \Mmaro^{\sklam{t}}\brklam{\Mmare^{\sklam{t}}}^{-1}\brklam{\Lmare^{\sklam{t}}}^{-1}\brklam{\Mmare^{\sklam{t}}}^{-1}\brklam{\Lmaro^{\sklam{t}}}^{-1}\brklam{\Mmaro^{\sklam{t}}}^{-1} \\
	&= \Malo^{\sklam{s}}\brklam{\Male^{\sklam{s}}}^{-1}\brklam{\Lale^{\sklam{s}}}^{-1}\brklam{\Male^{\sklam{s}}}^{-1}\brklam{\Lalo^{\sklam{s}}}^{-1}\brklam{\Malo^{\sklam{s}}}^{-1},
\end{align*}
im Fall $\kappa\geq4$
\begin{align*}
	&\ A^{\sklam{s}}_{\aln}-\alpha\Iq = -\beklam{A^{\sklam{t}}_{\aln}-(-\alpha)\Iq} \\
	&= -\prodr^n_{j=0}\eklam{\brklam{\Mmarj^{\sklam{t}}}^{-1}\brklam{\Lmarj^{\sklam{t}}}^{-1}}\brklam{\Mmarnp^{\sklam{t}}}^{-1}\brklam{\Mmarn^{\sklam{t}}+\Mmarnp^{\sklam{t}}}\prodl^{n-1}_{j=0}\beklam{\Lmarj^{\sklam{t}}\Mmarj^{\sklam{t}}} \\
	&= -\prodr^n_{j=0}\eklam{\brklam{\Malj^{\sklam{s}}}^{-1}\brklam{\Lalj^{\sklam{s}}}^{-1}}\brklam{\Malnp^{\sklam{s}}}^{-1}\brklam{\Maln^{\sklam{s}}+\Malnp^{\sklam{s}}}\prodl^{n-1}_{j=0}\beklam{\Lalj^{\sklam{s}}\Malj^{\sklam{s}}}
\end{align*}
für alle $n\in\Z{1}{\fklam{\kappa-2}}$ und im Fall $\kappa\geq5$
\begin{align*}
	B^{\sklam{s}}_{\aln} &= B^{\sklam{t}}_{-\arn} 
	= \prodr^{n-2}_{j=0}\beklam{\Mmarj^{\sklam{t}}\Lmarj^{\sklam{t}}}\Mmarnm^{\sklam{t}}\brklam{\Mmarn^{\sklam{t}}}^{-1}\prodl^n_{j=0}\eklam{\brklam{\Lmarj^{\sklam{t}}}^{-1}\brklam{\Mmarj^{\sklam{t}}}^{-1}} \\
	&= \prodr^{n-2}_{j=0}\beklam{\Malj^{\sklam{s}}\Lalj^{\sklam{s}}}\Malnm^{\sklam{s}}\brklam{\Maln^{\sklam{s}}}^{-1}\prodl^n_{j=0}\eklam{\brklam{\Lalj^{\sklam{s}}}^{-1}\brklam{\Malj^{\sklam{s}}}^{-1}}
\end{align*}
für alle $n\in\Zzfkm$. \bwend

Nun zeigen wir einen Zusammenhang zwischen der linksseitigen $\alpha$-Stieltjes"=Parametrisierung einer linksseitig $\alpha$-Stieltjes-positiv definiten Folge und den Favard-Paaren bezüglich jener Folge und deren durch linksseitige $\alpha$-Verschiebung generierten Folge.

\begin{satz}	\thlabel{wzlsa2}
	Seien $\kappa \in \Na$, $\alpha \in \R$, $\sjk \in \Lpqka$ und $\Qaljk$ die linksseitige $\alpha$-Stieltjes-Parametrisierung von $\sjk$. \dgfa
	\begin{itemize}
		\item [\rm{(a)}] Sei $\ABfk$ das Favard-Paar bezüglich $\sjk$. Dann gelten \linebreak $B_0 = Q_{\alo}$, $A_0-\alpha\Iq = -Q_{\al1}Q^{-1}_{\alo}$, im Fall $\kappa\geq2$
		\begin{align*}
			B_n = Q^{-1}_{\al2n-2}Q_{\al2n}
		\end{align*}
		für alle $n\in\Zefk$ und im Fall $\kappa\geq3$
		\begin{align*}
			A_n-\alpha\Iq = -Q_{\al2n+1}Q^{-1}_{\al2n}-Q_{\al2n}Q^{-1}_{\al2n-1}
		\end{align*}
		für alle $n\in\Zefkm$.
		\item [\rm{(b)}] Sei $(B_{\aln})^0_{n=0}$ bzw. im Fall $\kappa\geq2$ $\ABalfk$ das Favard-Paar bezüglich $\saljk$. Dann gelten $B_{\alo} = Q_{\al1}$, im Fall $\kappa\geq2$
		\begin{align*}
			A_{\aln}-\alpha\Iq = -Q_{\al2n+2}Q^{-1}_{\al2n+1}-Q_{\al2n+1}Q^{-1}_{\al2n}
		\end{align*}
		für alle $n\in\Z{0}{\fklam{\kappa-2}}$ und im Fall $\kappa\geq3$
		\begin{align*}
			B_{\aln} = Q^{-1}_{\al2n-1}Q_{\al2n+1}
		\end{align*}
		für alle $n\in\Zefkm$.
	\end{itemize}
\end{satz}

\bwanf Sei $t_j := (-1)^js_j$ für alle $j\in\Zok$. Wegen Teil (a) von \thref{asmbm3} gilt dann $\tjk\in\Kpqkma$. Weiterhin sei $(Q^{\sklam{t}}_{-\arj})^{\kappa}_{j=0}$ die rechtsseitige $-\alpha$-Stieltjes-Parametrisierung von $\tjk$.

Zu (a): Sei $[(A^{\sklam{t}}_{n})^{\fklam{\kappa-1}}_{n=0}$,""$(B^{\sklam{t}}_{n})^{\fklam{\kappa}}_{n=0}]$ das Favard-Paar von $\tjk$. Wegen Teil (a) von \thref{wzllm1}, Teil (a) von \thref{wzsa2} und \thref{aspbm4} gelten dann
\begin{align*}
	B^{\sklam{s}}_0 &= B^{\sklam{t}}_0 
	= Q^{\sklam{t}}_{-\aro} = Q^{\sklam{s}}_{\alo}, \\
	A^{\sklam{s}}_0-\alpha\Iq &= -\eklam{A^{\sklam{t}}_0-(-\alpha)\Iq}
	= -Q^{\sklam{t}}_{-\ar1}\brklam{Q^{\sklam{t}}_{-\aro}}^{-1} = -Q^{\sklam{s}}_{\al1}\brklam{Q^{\sklam{s}}_{\alo}}^{-1}, 
\end{align*}
im Fall $\kappa\geq2$
\begin{align*}
	B^{\sklam{s}}_n &= B^{\sklam{t}}_n 
	= \brklam{Q^{\sklam{t}}_{-\ar2n-2}}^{-1}Q^{\sklam{t}}_{-\ar2n} = \brklam{Q^{\sklam{s}}_{\al2n-2}}^{-1}Q^{\sklam{s}}_{\al2n}
\end{align*}
für alle $n\in\Zefk$ und im Fall $\kappa\geq3$
\begin{align*}
	A^{\sklam{s}}_n-\alpha\Iq &= -\eklam{A^{\sklam{t}}_n-(-\alpha)\Iq}
	= -Q^{\sklam{t}}_{-\ar2n+1}\brklam{Q^{\sklam{t}}_{-\ar2n}}^{-1}-Q^{\sklam{t}}_{-\ar2n}\brklam{Q^{\sklam{t}}_{-\ar2n-1}}^{-1} \\
	&= -Q^{\sklam{s}}_{\al2n+1}\brklam{Q^{\sklam{s}}_{\al2n}}^{-1}-Q^{\sklam{s}}_{\al2n}\brklam{Q^{\sklam{s}}_{\al2n-1}}^{-1}
\end{align*}
für alle $n\in\Zefkm$.

Zu (b): Sei $(B^{\sklam{t}}_{-\arn})^0_{n=0}$ bzw. im Fall $\kappa\geq2$ $[(A^{\sklam{t}}_{-\arn})^{\fklam{\kappa-2}}_{n=0}$,""$(B^{\sklam{t}}_{-\arn})^{\fklam{\kappa-1}}_{n=0}]$ das Favard-Paar bezüglich $(t_{-\arj})^{\kappa-1}_{j=0}$. Wegen Teil (b) von \thref{wzllm1}, Teil (b) von \thref{wzsa2} und \thref{aspbm4} gelten dann
\begin{align*}
	B^{\sklam{s}}_{\alo} = B^{\sklam{t}}_{-\aro} 
	= Q^{\sklam{t}}_{-\ar1} = Q^{\sklam{s}}_{\al1}, 
\end{align*}
im Fall $\kappa\geq2$
\begin{align*}
	A^{\sklam{s}}_{\aln}-\alpha\Iq &= -\eklam{A^{\sklam{t}}_{-\arn}-(-\alpha)\Iq }
	= -Q^{\sklam{t}}_{-\ar2n+2}\brklam{Q^{\sklam{t}}_{-\ar2n+1}}^{-1}-Q^{\sklam{t}}_{-\ar2n+1}\brklam{Q^{\sklam{t}}_{-\ar2n}}^{-1} \\
	&= -Q^{\sklam{s}}_{\al2n+2}\brklam{Q^{\sklam{s}}_{\al2n+1}}^{-1}-Q^{\sklam{s}}_{\al2n+1}\brklam{Q^{\sklam{s}}_{\al2n}}^{-1}
\end{align*}
für alle $n\in\Z{0}{\fklam{\kappa-2}}$ und im Fall $\kappa\geq3$
\begin{align*}
	B^{\sklam{s}}_{\aln} = B^{\sklam{t}}_{-\arn} 
	= \brklam{Q^{\sklam{t}}_{-\ar2n-1}}^{-1}Q^{\sklam{t}}_{-\ar2n+1} = \brklam{Q^{\sklam{s}}_{\al2n-1}}^{-1}Q^{\sklam{s}}_{\al2n+1}
\end{align*}
für alle $n\in\Zefkm$. \bwend

\newpage
\section[Einige Zusammenhänge zum matriziellen Hausdorffschen Momentenproblem]{Einige Zusammenhänge zum matriziellen \\ Hausdorffschen Momentenproblem} \label{chapma}

Seien $m\in\Noa$, $\sjk$ eine Folge aus $\Cqq$ und $\alpha,\beta\in\R$ mit $\alpha<\beta$. Dann betrachten wir in diesem Abschnitt das Momentenproblem $\Mabskg$, auch Hausdorffsches Momentenproblem genannt. Der Fall $\kappa=2n+1$ mit $n\in\No$ wurde in \cite{C06} und der Fall $\kappa=2n$ mit $n\in\No$ wurde in \cite{C07} behandelt.
Wir wollen hier einige Zusammenhänge zu unseren bisherigen Ergebnissen für das rechtsseitige bzw. linksseitige $\alpha$-Stieltjes Momentenproblem finden, wobei wir im linksseitgen Fall hierfür $\beta$ statt $\alpha$ verwenden. 

Mithilfe der $[\alpha,\beta]$-Stieltjes-Transformation (vergleiche \thref{masa1} und \thref{madef1}) können wir nun das matrizielle Hausdorffsche Momentenproblem wie folgt umformulieren:
\label{Sss}
\begin{itemize}
  	\item $\Sabskg$: Seien $\alpha,\beta \in \R$, $\kappa \in \Noa$ und $\sjk$ eine Folge aus $\Cqq$. 
  	Beschreibe die Menge $\Rqabsjk$ aller $S \in \Rqab$, deren zugehöriges Stieltjes-Maß zu $\Mqabskg$ gehört.
\end{itemize}

Die Ansätze von \thref{mabm2} bringen uns auf folgendes Resultat, wo wir speziell das Hausdorffsche Momentenproblem auf dem Intervall $[\alpha,\beta]$ sowie das rechtsseitige $\alpha$-Stieltjes Momentenproblem und das linksseitige $\beta$-Stieltjes Momentenproblem betrachten. Eine allgemeinere Form dieser Aussage findet man z.\,B. in \cite[Remark 2.25]{12}.

\begin{satz}	\thlabel{masa2}
	Seien $\kappa\in\Noa$, $\alpha,\beta\in\R$ mit $\alpha<\beta$, $\sjk$ eine Folge aus $\Cqq$ und $\Mqabskg$ nichtleer. Weiterhin seien $\mu\in\Mqabskg$ sowie $\mu_{[\alpha,\infty)}:\B_{[\alpha,\infty)}\rightarrow\Cqq$ und $\mu_{(-\infty,\beta]}:\B_{(-\infty,\beta]}\rightarrow\Cqq$ definiert gemäß
	\begin{align*}
		\mu_{[\alpha,\infty)}(A) := \mu(A\cap[\alpha,\beta])
	\end{align*}
	und
	\begin{align*}
		\mu_{(-\infty,\beta]}(A) := \mu(A\cap[\alpha,\beta]).
	\end{align*}
	Dann gelten $\mu_{[\alpha,\infty)}\in\Mqaskg$ und $\mu_{(-\infty,\beta]}\in\Mqmbskg$. Insbesondere sind jeweils $\Mqaskg$ und $\Mqmbskg$ nichtleer.
\end{satz}

\bwanf Unter Beachtung von $\mu\in\Mqab$ gelten wegen Teil (a) von \thref{hmlbm2} dann $\mu_{[\alpha,\infty)}\in\Mqa$ und $\mu_{(-\infty,\beta]}\in\Mqmb$. Wegen Teil (b) von \thref{hmlbm2} gelten weiterhin
\begin{align*}
	s^{\eklam{\mu_{[\alpha,\infty)}}}_j 
	= \int_{[\alpha,\infty)} t^j \;\mu_{[\alpha,\infty)}(dt)
	= \int_{[\alpha,\beta]} t^j \;\mu(dt)
	= s^{\eklam{\mu}}_j
\end{align*}
und
\begin{align*}
	s^{\eklam{\mu_{(-\infty,\beta]}}}_j 
	= \int_{(-\infty,\beta]} t^j \;\mu_{(-\infty,\beta]}(dt)
	= \int_{[\alpha,\beta]} t^j \;\mu(dt)
	= s^{\eklam{\mu}}_j
\end{align*}
für alle $j\in\Zok$. Hieraus folgen unter Beachtung von $s^{\eklam{\mu}}_j = s_j$ für alle $j\in\Zok$ dann $\mu_{[\alpha,\infty)}\in\Mqaskg$ und $\mu_{(-\infty,\beta]}\in\Mqmbskg$. \bwend

Wir kommen nun zu den Lösbarkeitsbedingungen des Hausdorffschen Momentenproblems auf dem Intervall $[\alpha,\beta]$ und betrachten anschließend einen Zusammenhang zum rechtsseitigen $\alpha$-Stieltjes Momentenproblem und linksseitigen $\beta$-Stieltjes Momentenproblem. Hierzu benötigen wir noch folgende Bezeichnung.

\begin{bez}	\thlabel{mabz2}
	Seien $\kappa\in\Na\setminus{\gklam{1}}$. Dann sei für alle $n\in\Z{0}{\fklam{\kappa-2}}$
	\begin{align*}
		\wKsn := (s_{j+k+2})^n_{j,k=0}.
	\end{align*}
	Falls klar ist, von welchem $\sjk$ die Rede ist, lassen wir das \anf{$\sklam{s}$} als oberen Index weg.
\end{bez}

\begin{satz}	\thlabel{masa3}
	Seien $\alpha,\beta\in\R$ mit $\alpha<\beta$. \dgfa
	\begin{itemize}
		\item [\rm{(a)}] Seien $n\in\No$ und $\sjne$ eine Folge aus $\Cqq$. Dann ist $\Mqabsneg$ genau dann nichtleer, wenn $\Harn$ und $\Hbln$ nichtnegativ hermitesch sind.
		\item [\rm{(b)}] Seien $n\in\No$ und $\sjn$ eine Folge aus $\Cqq$. Dann ist $\Mqabsng$ genau dann nichtleer, wenn $\Hn$ und im Fall $n\in\N$ auch 
		\begin{align*}
			-\alpha\beta\Hnm+(\alpha+\beta)\Knm-\wKnm
		\end{align*}
		nichtnegativ hermitesch sind.
	\end{itemize}
\end{satz}

\bwanf Zu (a): Siehe \cite[Theorem 1.3]{C06}.

Zu (b): Der Fall $n\in\N$ wurde in \cite[Theorem 1.3]{C07} behandelt. Im Fall $n=0$ ist $\Mqabsng$ die Menge aller $\mu \in \Mqab$ mit $\mu([\alpha,\beta])=s_0$. Somit ist $\Mqabsog$ genau dann nichtleer, wenn $H_0 = s_0$ nichtnegativ hermitesch ist. \bwend

\begin{satz}	\thlabel{masa4}
	Seien $n\in\No$, $\alpha,\beta\in\R$ mit $\alpha<\beta$ und $\sjne$ eine Folge aus $\Cqq$. Dann ist $\Mqabsneg$ genau dann nichtleer, wenn $\Mqasneu$ und $\Mqmbsneuu$ nichtleer sind.
\end{satz}

\bwanf Sei $\Mqabsneg$ nichtleer. Wegen Teil (a) von \thref{masa3} sind dann $\Harn$ und $\Hbln$ nichtnegativ hermitesch. Es gilt
\begin{align*}
	\frac{1}{\beta-\alpha}(\Harn+\Hbln) = \frac{1}{\beta-\alpha}(-\alpha\Hn+\Kn+\beta\Hn-\Kn) = \Hn.
\end{align*}
Unter Beachtung von $\alpha<\beta$ folgt hieraus dann, dass $\Hn$ nichtnegativ hermitesch ist. Wegen Teil (c) von \thref{asmth2} und Teil (b) von \thref{aspbm1} ist dann \linebreak $\Mqasneu$ nichtleer. Weiterhin ist wegen Teil (d) von \thref{asmth3} und Teil (b) von \thref{asplbm1} auch $\Mqmbsneuu$ nichtleer.

Seien nun $\Mqasneu$ und $\Mqmbsneuu$ nichtleer. Wegen Teil (c) von \thref{asmth2} und Teil (b) von \thref{aspbm1} ist dann $\Harn$ nichtnegativ hermitesch. Weiterhin ist wegen Teil (d) von \thref{asmth3} und Teil (b) von \thref{asplbm1} auch $\Hbln$ nichtnegativ hermitesch. Hieraus folgt wegen Teil (a) von \thref{masa3} dann, dass $\Mqabsneg$ nichtleer ist. \bwend

Ein zu \thref{masa3} ähnliches Resultat findet man in der zu dieser Dissertation parallel entstandenen Arbeit \cite{HMS} unter \cite[Proposition 9.3]{HMS}.

%% file: anhang.tex
\newpage
\section[Einige Aussagen zur Integrationstheorie nichtnegativ hermitescher Maße]{Einige Aussagen zur Integrationstheorie \\ nichtnegativ hermitescher Maße} \label{chapAA}

In diesem Abschnitt behandeln wir einige ausgewählte Resultate der Integrationstheorie nichtnegativ hermitescher Maße. Es sei hier nochmal auf die Einleitung verwiesen, wo wir den Integrationsbegriff nichtnegativ hermitescher Maße einführen (vergleiche \sref{Mq}). Da wir hier dieses Thema nicht detailliert ausarbeiten werden, kann der Leser in \cite[Anhang M]{MP} und \cite{Pe} eine ausführlichere Beschreibung finden.

Unser erstes Resultat dieses Abschnitts ist von zentraler Bedeutung und ermöglicht uns das darauffolgende Lemma, das uns dazu dient, die Lösbarkeitsbedingungen für das linksseitige $\alpha$-Stieltjes Momentenproblem aus dem rechtsseitigen Fall zu folgern. In diesem Satz betrachten wir die Integration einer messbaren Abbildung unter beliebiger Transformation.

\begin{satz}	\thlabel{hmsa1}
  Seien $(\Omega_1,\A_1)$ und $(\Omega_2,\A_2)$ messbare Räume sowie $\mu \in \Mqoae$.
  Weiterhin sei $T: \Omega_1 \rightarrow \Omega_2$ eine $\A_1$-$\A_2$-messbare Abbildung. \dgfa
  \begin{itemize}
    \item [\rm{(a)}] Es gilt $T(\mu) \in \Mqoaz$.
    \item [\rm{(b)}] Sei $f: \Omega_2 \rightarrow \C$ eine $\A_2$-$\B_\C$-messbare Abbildung. \dsfaa
    \begin{itemize}
      \item [\rm{(i)}] Es gilt $f \in {\cal L}^1\rklam{\Omega_2,\A_2,T(\mu);\C}$.
      \item [\rm{(ii)}] Es gilt $f\circ T \in {\cal L}^1\rklam{\Omega_1,\A_1,\mu;\C}$.
    \end{itemize}
    \item [\rm{(c)}] Sei $f \in {\cal L}^1\rklam{\Omega_2,\A_2,T(\mu);\C}$. Dann gilt
    \begin{align*}
      \int_A f \;d\eklam{T(\mu)} = \int_{T^{-1}(A)} (f\circ T) \;d\mu
    \end{align*}
    für alle $A \in \A_2$.
  \end{itemize}
\end{satz}

\bwanf Dies ist ein wohlbekanntes Resultat der Maßtheorie. Einen Beweis findet man z.\,B. unter \cite[Proposition B.1]{ot226} oder \cite[Satz M.47]{MP}. \bwend

Wir werden nun die Aussage von \thref{hmsa1} für einen für uns relevanten Spezialfall betrachten (vergleiche \cite[Lemma B.2]{ot226}). Dafür verwenden wir als zugrundeliegende Transformation eine Spiegelung und führen zunächst einige Bezeichnungen ein.

\begin{bez}	\thlabel{hmbz1}
	Sei $\Omega$ eine Teilmenge von $\C$. Dann bezeichne $\widecheck{\Omega} := \gklam{-\omega \;|\; \omega \in \Omega}$.
	Seien nun $\Omega$ eine nichtleere Borel-Teilmenge von $\R$ und $\mu \in \Mqo$. Für $a,b \in \R$ sei weiterhin $T_{a,b}: \R \rightarrow \R$ definiert gemäß $T_{a,b}(x) := ax+b$. Dann bezeichne $\widecheck{\mu} :=  T_{-1,0}(\mu)$ das nach Teil (a) von \thref{hmsa1} zu $\Mqco$ zugehörige Maß.
\end{bez}

\begin{lemma}	\thlabel{hmlm1}
  Seien $\Omega$ eine nichtleere Borel-Teilmenge von $\R$, $\mu \in \Mqo$ und $\kappa \in \Noa$. \dgfa
  \begin{itemize}
    \item [\rm{(a)}] \esfaa
    \begin{itemize}
      \item [\rm{(i)}] Es gilt $\mu \in \Mqko$.
      \item [\rm{(ii)}] Es gilt $\widecheck{\mu} \in \Mqkco$.
    \end{itemize}
    \item [\rm{(b)}] Sei $\mu \in \Mqko$. Dann gilt $s^{(\widecheck{\mu})}_j = (-1)^j s^{(\mu)}_j$ für alle $j \in \Zok$.
  \end{itemize}
\end{lemma}

\bwanf Offensichtlich ist $T_{-1,0}$ eine $\B_\R$-$\B_\R$-messbare Abbildung und es gilt \linebreak $T^{-1}_{-1,0}(\widecheck{\Omega}) = \Omega$.
Sei $j \in \Zok$. Weiterhin sei $f_j: \R \rightarrow \C$ definiert gemäß $f_j(x) = x^j$. 
Dann gelten $f_j \in {\cal L}^1\rklam{\R,\B_\R,\widecheck{\mu};\C}$ und $f_j\circ T_{-1,0} = (-1)^j f_j$. 
Hieraus folgt mithilfe von Teil (c) von \thref{hmsa1} dann
\begin{align*}
  \int_{\widecheck{\Omega}} f_j \;d\widecheck{\mu} = (-1)^j \int_\Omega f_j \;d\mu.
\end{align*}
Hieraus folgen wiederum dann alle Behauptungen. \bwend

Nun betrachten wir ein nichtnegativ hermitesches Maß, das auf einer nichtleeren Teilmenge einer Grundmenge definiert ist, und wollen es auf ein nichtnegativ hermitesches Maß, das auf der ganzen Grundmenge definiert ist, erweitern.

\begin{bem}	\thlabel{hmlbm2}
	Seien $(\Omega,\A)$ ein messbarer Raum, $B\in\A\setminus\gklam{\emptyset}$ und $\mu\in\M^q_{\geq}(B,\linebreak\A\cap B)$. Weiterhin sei $\mu_B:\A\rightarrow\Cqq$ definiert gemäß $\mu_B := \mu(A\cap B)$. \dgfa
	\begin{itemize}
		\item [\rm{(a)}] Es gelten $\mu_B\in\Mqoa$ und $\Rstr_{\A\cap B}\mu_B = \mu$.
		\item [\rm{(b)}] Sei $f:\Omega\rightarrow\C$ eine $\A$-$\B_\C$-messbare Abbildung. Dann gilt $f\in{\cal L}^1(\Omega,\A,\mu_B;\C)$ genau dann, wenn $\Rstr_B f\in{\cal L}^1(B,\A\cap B,\mu;\C)$ erfüllt ist. In diesem Fall gilt weiterhin
			\begin{align*}
				\int_\Omega f \;d\mu_B = \int_B f \;d\mu.
			\end{align*}
	\end{itemize}
\end{bem}

\bwanf Dies ist ein wohlbekanntes Resultat und kann mit Standardmethoden der Maß- und Integrationstheorie gezeigt werden. \bwend

Wir wollen nun eine spezielle Abschätzung der Spektralnorm eines Integrals einer messbaren Abbildung mithilfe des Spurmaßes des zugrundeliegenden nichtnegativ hermiteschen Maßes behandeln.

\begin{bem}	\thlabel{hmlbm1}
	Seien $(\Omega,\A)$ ein messbarer Raum, $\mu\in\Mqoa$ und $\tau$ das Spurmaß von $\mu$. Weiterhin sei $f:\Omega\rightarrow\C$ eine $\A$-$\B_\C$-messbare Abbildung. \dgfa
	\begin{itemize}
		\item [\rm{(a)}] \esfaa
		\begin{itemize}
      		\item [\rm{(i)}] Es gilt $f\in{\cal L}^1(\Omega,\A,\mu;\C)$.
      		\item [\rm{(ii)}] Es gilt $\abs{f}\in{\cal L}^1(\Omega,\A,\tau;\C)$.
      	\end{itemize}
      	\item [\rm{(b)}] Sei $f\in{\cal L}^1(\Omega,\A,\mu;\C)$. Dann gilt
      	\begin{align*}
      		\snorm{\int_\Omega f \;d\mu} \leq 4q \int_\Omega \abs{f} d\tau.
      	\end{align*}
    \end{itemize}
\end{bem}

\bwanf Einen Beweis findet man unter \cite[Lemma M.39]{MP}, wobei dort die euklidische Norm $\enorm{\cdot}$ verwendet wird, d.\,h. es gilt
\begin{align*}
	\enorm{\int_\Omega f \;d\mu} \leq 4q \int_\Omega \abs{f} d\tau.
\end{align*}
Unter Beachtung von
\begin{align*}
	\snorm{A} \leq \enorm{A} \leq \sqrt{pq} \snorm{A}
\end{align*}
für alle $A\in\Cpq$ folgt dann die Behauptung. \bwend

Unser letztes Resultat für diesen Abschnitt widmet sich der Adjungierten von einem Integral einer messbaren Abbildung.

\begin{lemma}	\thlabel{hmllm2}
	Seien $(\Omega,\A)$ ein messbarer Raum, $\mu\in\Mqoa$ und $f\in{\cal L}^1(\Omega,\A,\mu;\C)$. Dann gelten $\overline{f}\in{\cal L}^1(\Omega,\A,\mu;\C)$ und
	\begin{align*}
		\int_\Omega \overline{f} \;d\mu = \rklam{\int_\Omega f \;d\mu}^{\ast}.
	\end{align*}
\end{lemma}

\bwanf Siehe z.\,B. \cite[Lemma 2.23]{Pe}. \bwend

\section{Über die Stieltjes-Transformation von nichtnegativ hermiteschen Maßen} \label{chapAS}

In diesem Abschnitt betrachten wir nichtnegativ hermitesche \textit{q}$\times$\textit{q}-Maße auf der Borelschen $\sigma$-Algebra eines reellen Intervalls des Typs $[\alpha,\infty)$, $(-\infty,\alpha]$ oder $[\alpha,\beta]$ für beliebige $\alpha,\beta\in\R$ mit $\alpha<\beta$. Über eine spezielle Integraltransformation ordnen wir einem solchen nichtnegativ hermiteschen \textit{q}$\times$\textit{q}-Maß eine in $\C\setminus[\alpha,\infty)$, $\C\setminus(-\infty,\alpha]$ bzw. $\C\setminus[\alpha,\beta]$ holomorphe \textit{q}$\times$\textit{q}-Matrixfunktion zu. Hierdurch wird ein bijektiver Zusammenhang zwischen der Menge von nichtnegativ hermiteschen \textit{q}$\times$\textit{q}-Maßen auf der Borelschen $\sigma$-Algebra von $[\alpha,\infty)$, $(-\infty,\alpha]$ bzw. $[\alpha,\beta]$ und einer Menge von in $\C\setminus[\alpha,\infty)$, $\C\setminus(-\infty,\alpha]$ bzw. $\C\setminus[\alpha,\beta]$ holomorphen \textit{q}$\times$\textit{q}-Matrixfunktionen hergestellt.

Diese Vorgehensweise geht im klassischen Fall $q=1$ für das Intervall $[0,\infty)$ auf T.-J. Stieltjes \cite{Sur} zurück. Man spricht deshalb auch von der Stieltjes-Transformation. Via Stieltjes-Transformation wird es dann möglich, die ursprünglich für Maße formulierte Probleme (wie z.\,B. Momentenprobleme) in äquivalente Probleme zur Bestimmung von holomorphen Matrixfunktionen mit gewissen Eigenschaften zu überführen. Hierbei sei bemerkt, dass Momentenprobleme vielfach in Probleme der Bestimmung von in einer Halbebene holomorphen Matrixfunktionen mit durch die vorgegebenen Momente festgelegten asymptotischen Entwicklungen übergehen. Wir beginnen nun mit den konkreten Detailbetrachtungen.


Zuerst widmen wir uns Maßen auf dem Intervall $[\alpha,\infty)$, die in Verbindung mit dem rechtsseitigen $\alpha$-Stieltjes Momentenproblem stehen. Wir stellen zunächst eine für unsere weiteren Betrachtungen wichtige Klasse holomorpher Matrixfunktionen bereit.

\begin{bez}	\thlabel{asmbz3}
	Sei $\alpha \in \R$. Dann bezeichne $\Sqa$ die Menge aller in $\C\setminus[\alpha,\infty)$ \linebreak holomorphen \textit{q}$\times$\textit{q}-Matrixfunktionen $S$ mit 
$\im S(z) \geq \Oq$ für alle $z \in \Pp$ und \linebreak $S(x) \geq \Oq$ für alle $x \in (-\infty,\alpha)$. Weiterhin sei
	\begin{align*}
  		\Soqa := \{ S \in \Sqa \;|\; \sup{}_{y \in [1,\infty)} y \snorm{S(iy)} < \infty \}.
	\end{align*}
\end{bez}

Die Klasse $\Sqa$ bzw. $\Soqa$ wurde in \cite[Chapter 3 and 4]{SF} bzw. \cite[Chapter 5]{SF} näher beleuchtet.

\begin{satz}[$\lbrack\alpha,\infty)$-Stieltjes-Transformation]	\thlabel{asmth4}   Sei $\alpha \in \R$. \dgfa
\begin{itemize}
  \item [\rm{(a)}] Sei $ S \in \Soqa $. Dann existiert genau ein $ \mu \in \Mqa $ mit
  \begin{align}		\label{asmth4bw1}
    S(z) = \int_{[\alpha,\infty)} \frac{1}{t-z} \;\mu(dt)
  \end{align}
  für alle $z \in \C\setminus[\alpha,\infty)$.
  \item [\rm{(b)}] Sei $ \mu \in \Mqa $. Dann definiert $ S: \C\setminus[\alpha,\infty) \rightarrow \Cqq $ gemäß \fref{asmth4bw1} eine
  wohldefinierte Matrixfunktion aus $ \Soqa $.
\end{itemize}
\end{satz}

\bwanf Siehe \cite[Theorem 5.1]{SF} oder \cite[Theorem 1.73]{MP}. \bwend

\begin{defi}	\thlabel{asmdef4}
  Sei $\alpha \in \R$. Weiterhin sei $ \mu \in \Mqa $. Dann heißt die in Teil (b) von \thref{asmth4} definierte Abbildung $ S $ \textbf{Stieltjes-Transformierte} von $ \mu $.
  Sei nun umgekehrt $ S \in \Soqa $. Dann heißt das in Teil (a) von \thref{asmth4} definierte Maß $ \mu $ das zu $ S $ gehörige \textbf{Stieltjes-Maß}.
\end{defi}

Das nachfolgende Resultat behandelt die Spalten- und Nullräume für die Funktionswerte einer Matrixfunktion aus $\Soqa$. Es stellt sich heraus, dass diese Räume sogar konstant sind und mithilfe des zugehörigen Stieltjes-Maßes angegeben werden können.

\begin{bem}	\thlabel{mabm3}
	Seien $\alpha\in\R$ und $S\in\Soqa$. Weiterhin bezeichne $\mu$ das zu $S$ gehörige Stieltjes-Maß aus $\Mqa$. Dann gelten
	\begin{align*}
		{\cal R}\brklam{S(z)} = {\cal R}\brklam{\mu([\alpha,\infty))}
		\quad \text{und} \quad
		{\cal N}\brklam{S(z)} = {\cal N}\brklam{\mu([\alpha,\infty))}
	\end{align*}
	für alle $z\in\C\setminus[\alpha,\infty)$.
\end{bem}

\bwanf Siehe \cite[Proposition 5.3]{SF}. \bwend

Wir wollen nun die Stieltjes-Transformation für Maße auf dem Intervall $(-\infty,\alpha]$, die in Verbindung mit dem linksseitigen $\alpha$-Stieltjes Momentenproblem stehen, einführen. Hierbei werden wir die Stieltjes-Transformation mithilfe von \thref{hmsa1} aus \thref{asmth4} gewinnen.

\begin{bez}	\thlabel{asmbz4}
	Sei $\alpha \in \R$. Dann bezeichne $\Sqma$ die Menge aller in $\C\setminus(-\infty,\alpha]$ holomorphen \textit{q}$\times$\textit{q}-Matrixfunktionen $S$ mit 
$-\im S(z) \geq \Oq$ für alle $z \in \Pm$ und $-S(x) \geq \Oq$ für alle $x \in (\alpha,\infty)$. Weiterhin sei
	\begin{align*}
  		\Soqma := \{ S \in \Sqma \;|\; \sup{}_{y \in (-\infty,-1]} -y \snorm{S(iy)} < \infty \}.
	\end{align*}
\end{bez}

Eine äquivalente Formulierung der Klassen $\Sqma$ bzw. $\Soqma$ wurde in \cite[Chapter 7 and 8]{SF} bzw. \cite[Chapter 9]{SF} näher beleuchtet. Man erhält die dort vorgenommene Definition durch die für $S\in\Sqma$ gültige Beziehung $S^{\ast}(\za)=S(z)$ für alle $z\in\C\setminus(-\infty,\alpha]$. Diese Eigenschaft geht aus der Nevanlinna-Parametrisierung oder speziell für $S\in\Soqma$ auch aus der Stieltjes-Transformation hervor (vergleiche Beweis von \thref{mabm1}).

\begin{bem}	\thlabel{asmbm4}
	Sei $S: \C\setminus[\alpha,\infty) \rightarrow \Cqq$ eine Matrixfunktion. 
	Weiterhin sei $\widecheck{S}: \C\setminus(-\infty,-\alpha] \rightarrow \Cqq$ definiert gemäß $\widecheck{S}(z) := -S(-z)$. Dann gilt $S\in\Sqa$ bzw. $S\in\Soqa$ genau dann, wenn $\widecheck{S}\in\Sqmma$ bzw. $\widecheck{S}\in\Soqmma$ erfüllt ist.
\end{bem}

\bwanf Aus der Definition von $\widecheck{S}$ folgen sogleich
\begin{itemize}
	\item [\rm{(I)}] Es ist $S$ holomorph in $\C\setminus[\alpha,\infty)$ genau dann, wenn $\widecheck{S}$ in $\C\setminus(-\infty,-\alpha]$ holomorph ist.
	\item [\rm{(II)}] Es gilt $S(x) \geq \Oq$ für alle $x \in (-\infty,\alpha)$ genau dann, wenn $-\widecheck{S}(x) \geq \Oq$ für alle $x \in (-\alpha,\infty)$ erfüllt ist.
	\item [\rm{(III)}] Es gilt $\im S(z) \geq \Oq$ für alle $z \in \Pp$ genau dann, wenn $-\im \widecheck{S}(z) \geq \Oq$ für alle $z \in \Pm$ erfüllt ist.
\end{itemize}
Es gilt
\begin{align*}
	\sup_{y \in (-\infty,-1]} -y \bsnorm{\widecheck{S}(iy)} 
	= \sup_{y \in (-\infty,-1]} (-y) \abs{-1}\snorm{S(i(-y))}
	= \sup_{y \in [1,\infty)} y \snorm{S(iy)}.
\end{align*}
Hieraus folgt
\begin{itemize}
	\item [\rm{(IV)}] Es gilt 
	$\sup_{y \in [1,\infty)} \negthinspace y ||S(iy)||\negthinspace < \negthinspace\infty$
	genau dann, wenn 
	$\sup_{y \in (-\infty,-1]} \negthinspace-y ||\widecheck{S}(iy)||\negthinspace < \negthinspace\infty$ 
	erfüllt ist.
\end{itemize}
Wegen (I)-(IV) folgt dann die Behauptung (vergleiche \thref{asmbz3} und \thref{asmbz4}). \bwend

\begin{satz}[$(-\infty,\alpha\rbrack$-Stieltjes-Transformation]	\thlabel{asmth5} Sei $\alpha \in \R$. \dgfa
\begin{itemize}
  \item [\rm{(a)}] Sei $ S \in \Soqma $. Dann existiert genau ein $ \mu \in \Mqma $ mit
  \begin{align} \label{asmth4lbw1}
    S(z) = \int_{(-\infty,\alpha]} \frac{1}{t-z} \;\mu(dt)
  \end{align}
  für alle $z \in \C\setminus(-\infty,\alpha]$.
  \item [\rm{(b)}] Sei $ \mu \in \Mqma $. Dann definiert $ S: \C\setminus(-\infty,\alpha] \rightarrow \Cqq $ gemäß \fref{asmth4lbw1} eine
  wohldefinierte Matrixfunktion aus $ \Soqma $.
\end{itemize}
\end{satz}

\bwanf Zu (a): Sei $\widecheck{S}: \C\setminus[-\alpha,\infty) \rightarrow \Cqq$ definiert gemäß $\widecheck{S}(z) := -S(-z)$. Wegen \thref{asmbm4} ist dann $\widecheck{S}\in{\cal S}_{0,q,[-\alpha,\infty)}$. Wegen Teil (a) von \thref{asmth4} existiert genau ein $\tau \in \M^q_{\geq}([-\alpha,\infty))$ mit
\begin{align*}
	\widecheck{S}(z) = \int_{[-\alpha,\infty)} \frac{1}{t-z} \;\tau(dt)
\end{align*}
für alle $z \in \C\setminus[-\alpha,\infty)$. Sei $\mu := \widecheck{\tau}$. Wegen Teil (a) von \thref{hmsa1} ist dann \linebreak $\mu\in\Mqma$ und wegen Teil (c) von \thref{hmsa1} gilt
\begin{align*}
	\int_{(-\infty,\alpha]} \frac{1}{t-z} \;\mu(dt) 
	&= \int_{[-\alpha,\infty)} \frac{1}{-t-z} \;\tau(dt) \\
	&= -\int_{[-\alpha,\infty)} \frac{1}{t-(-z)} \;\tau(dt) \\
	&= -\widecheck{S}(-z) = S(z)
\end{align*}
für alle $z \in \C\setminus(-\infty,\alpha]$. 

Bleibt nur noch die Eindeutigkeit von $\mu$ zu zeigen. Sei hierfür $\sigma\in\Mqma$ mit
\begin{align*}
    S(z) = \int_{(-\infty,\alpha]} \frac{1}{t-z} \;\sigma(dt)
\end{align*}
für alle $z \in \C\setminus(-\infty,\alpha]$. Hieraus folgt wegen Teil (c) von \thref{hmsa1} dann
\begin{align*}
	\widecheck{S}(z) &= -S(-z) = -\int_{(-\infty,\alpha]} \frac{1}{t-(-z)} \;\sigma(dt) \\
	&= \int_{(-\infty,\alpha]} \frac{1}{-t-z} \;\sigma(dt)
	= \int_{[-\alpha,\infty)} \frac{1}{t-z} \;\widecheck{\sigma}(dt)
\end{align*}
für alle $z \in \C\setminus[-\alpha,\infty)$. Hieraus folgt wegen der Eindeutigkeit von $\tau$ dann $\tau = \widecheck{\sigma}$. Hieraus folgt nun $\mu = \sigma$.

Zu (b): Wegen Teil (a) von \thref{hmsa1} ist $\widecheck{\mu}\in\Mqma$. Hieraus folgt wegen Teil (b) von \thref{asmth4} dann, dass $ T: \C\setminus[-\alpha,\infty) \rightarrow \Cqq $ definiert gemäß 
\begin{align*}
	T(z) := \int_{[-\alpha,\infty)} \frac{1}{t-z} \;\widecheck{\mu}(dt)
\end{align*}
eine wohldefinierte Matrixfunktion aus ${\cal S}_{0,q,[-\alpha,\infty)}$ ist. Sei $S: \C\setminus(-\infty,\alpha] \rightarrow \Cqq$ definiert gemäß $S(z) := -T(-z)$. Wegen \thref{asmbm4} ist dann $S\in\Soqma$ und wegen Teil (c) von \thref{hmsa1} gilt
\begin{align*}
	S(z) &= -T(-z)
	= -\int_{[-\alpha,\infty)} \frac{1}{t-(-z)} \;\widecheck{\mu}(dt) \\
	&= \int_{[-\alpha,\infty)} \frac{1}{-t-z} \;\widecheck{\mu}(dt)
	= \int_{(-\infty,\alpha]} \frac{1}{t-z} \;\mu(dt)
\end{align*}
für alle $z \in \C\setminus(-\infty,\alpha]$. \bwend

\begin{defi}	\thlabel{asmdef5}
  Sei $\alpha \in \R$. Weiterhin sei $ \mu \in \Mqma $. Dann heißt die in Teil (b) von \thref{asmth5} definierte Abbildung $ S $ \textbf{Stieltjes-Transformierte} von $ \mu $.
  Sei nun umgekehrt $ S \in \Soqma $. Dann heißt das in Teil (a) von \thref{asmth5} definierte Maß $ \mu $ das zu $ S $ gehörige \textbf{Stieltjes-Maß}.
\end{defi}

Das nachfolgende Resultat behandelt die Spalten- und Nullräume für die Funktionswerte einer Matrixfunktion aus $\Soqma$. Es stellt sich heraus, dass diese Räume sogar konstant sind und mithilfe des zugehörigen Stieltjes-Maßes angegeben werden können.

\begin{bem}	\thlabel{mabm5}
	Seien $\alpha\in\R$ und $S\in\Soqma$. Weiterhin bezeichne $\mu$ das zu $S$ gehörige Stieltjes-Maß aus $\Mqma$. Dann gelten
	\begin{align*}
		{\cal R}\brklam{S(z)} = {\cal R}\brklam{\mu((-\infty,\alpha])}
		\quad \text{und} \quad
		{\cal N}\brklam{S(z)} = {\cal N}\brklam{\mu((-\infty,\alpha])}
	\end{align*}
	für alle $z\in\C\setminus(-\infty,\alpha]$.
\end{bem}

\bwanf Sei $\widecheck{S}: \C\setminus[-\alpha,\infty) \rightarrow \Cqq$ definiert gemäß $\widecheck{S}(z) := -S(-z)$. Wegen \thref{asmbm4} ist dann $\widecheck{S}\in{\cal S}_{0,q,[-\alpha,\infty)}$. Sei $\tau$ das zu $\widecheck{S}$ gehörige Stieltjes-Maß aus $\M^q_{\geq}([-\alpha,\infty))$. Wegen Teil (a) von \thref{asmth4} und Teil (a) von \thref{asmth5} gilt dann
\begin{align*}
	\int_{[-\alpha,\infty)} \frac{1}{t-z} \;\tau(dt) &= \widecheck{S}(z) = -S(-z) \\
	&= -\int_{(-\infty,\alpha]} \frac{1}{t-(-z)} \;\mu(dt) \\
	&= \int_{(-\infty,\alpha]} \frac{1}{-t-z} \;\mu(dt)
\end{align*}
für alle $z\in\C\setminus[-\alpha,\infty)$. Hieraus folgt wegen Teil (c) von \thref{hmsa1} sowie der Eindeutigkeit von $\mu$ und $\tau$ dann $\tau=\widecheck{\mu}$. Hieraus folgen wegen der Linearität der Spalten- und Nullräume sowie \thref{mabm3} nun
\begin{align*}
	{\cal R}\brklam{S(z)} = {\cal R}\brklam{{-\widecheck{S}(-z)}} = {\cal R}\brklam{\widecheck{S}(-z)} = {\cal R}\brklam{\tau([-\alpha,\infty))} = {\cal R}\brklam{\mu((-\infty,\alpha])}
\end{align*}
und
\begin{align*}
	{\cal N}\brklam{S(z)} = {\cal N}\brklam{{-\widecheck{S}(-z)}} = {\cal N}\brklam{\widecheck{S}(-z)} = {\cal N}\brklam{\tau([-\alpha,\infty))} = {\cal N}\brklam{\mu((-\infty,\alpha])}
\end{align*}
für alle $z\in\C\setminus(-\infty,\alpha]$. \bwend

Wir kommen nun zur Stieltjes-Transformation für Maße auf dem Intervall $[\alpha,\beta]$, die in Verbindung zum matriziellen Hausdorffschen Momentenproblem stehen.

\begin{bez}	\thlabel{mabz1}
	Seien $\alpha,\beta \in \R$ mit $\alpha<\beta$. Dann bezeichne $\Rqab$ die Menge aller in $\C\setminus[\alpha,\beta]$ holomorphen \textit{q}$\times$\textit{q}"=Matrixfunktionen $S$ mit
	\begin{itemize}
		\item [\rm{(i)}] Es gilt $\im S(z)\geq\Oq$ für alle $z\in\Pp$.
		\item [\rm{(ii)}] Es gilt $S(x)\geq\Oq$ für alle $x\in(-\infty,\alpha)$.
		\item [\rm{(iii)}] Es gilt $-S(x)\geq\Oq$ für alle $x\in(\beta,\infty)$.
	\end{itemize}
\end{bez}

\begin{satz}[$\lbrack\alpha,\beta\rbrack$-Stieltjes-Transformation]	\thlabel{masa1}   
Seien $\alpha,\beta \in \R$ mit $\alpha<\beta$. \dgfa
\begin{itemize}
  	\item [\rm{(a)}] Sei $S \in \Rqab$. Dann existiert genau ein $\mu \in \Mqab$ mit
  	\begin{align}	\label{masabw1}
    	S(z) = \int_{[\alpha,\beta]} \frac{1}{t-z} \;\mu(dt)
  	\end{align}
  	für alle $z \in \C\setminus[\alpha,\beta]$.
  	\item [\rm{(b)}] Sei $\mu \in \Mqab$. Dann definiert $S: \C\setminus[\alpha,\beta] \rightarrow \Cqq$ gemäß \fref{masabw1} eine wohldefinierte Matrixfunktion aus $\Rqab$.
\end{itemize}
\end{satz}

\bwanf Siehe \cite[Theorem 1.1]{C06}.

\begin{defi}	\thlabel{madef1}
  	Seien $\alpha,\beta \in \R$ mit $\alpha<\beta$. Weiterhin sei $\mu \in \Mqab$. Dann heißt die in Teil (b) von \thref{masa1} definierte Abbildung $S$ \textbf{Stieltjes-Transformierte} von $\mu$.
  	Sei nun umgekehrt $S \in \Rqab$. Dann heißt das in Teil (a) von \thref{masa1} definierte Maß $\mu$ das zu $S$ gehörige \textbf{Stieltjes-Maß}.
\end{defi}

Durch die Stieltjes-Transformation erhalten wir einen Symmetrie-Effekt zwischen der oberen und unteren offenen Halbebene von $\C$ für die Imaginärteile der für uns relevanten Stieltjes-Transformierten.

\begin{bem}	\thlabel{mabm1}	Seien $\alpha,\beta\in\R$ mit $\alpha<\beta$. \dgfa
\begin{itemize} 
	\item [\rm{(a)}] Sei $S\in\Soqa$. 
	Dann gilt $-\im S(z)\geq\Oq$ für alle $z\in\Pm$.
	\item [\rm{(b)}] Sei $S\in\Soqmb$. 
	Dann gilt $\im S(z)\geq\Oq$ für alle $z\in\Pp$.
	\item [\rm{(c)}] Sei $S\in\Rqab$. 
	Dann gilt $-\im S(z)\geq\Oq$ für alle $z\in\Pm$.	
\end{itemize}
\end{bem}

\bwanf Zu (a): Wegen Teil (a) von \thref{asmth4} existiert genau ein $\mu \in \Mqa$ mit
\begin{align*}
	S(z) = \int_{[\alpha,\infty)} \frac{1}{t-z} \;\mu(dt)
\end{align*}
für alle $z \in \C\setminus[\alpha,\infty)$. Wegen \thref{hmllm2} folgt hieraus dann
\begin{align}	\label{mabm1bw1}
	S^{\ast}(z) = \int_{[\alpha,\infty)} \overline{\rklam{\frac{1}{t-z}}} \;\mu(dt)
	= \int_{[\alpha,\infty)} \frac{1}{t-\za} \;\mu(dt)
	= S(\za)
\end{align}
für alle $z \in \C\setminus[\alpha,\infty)$. Sei nun $z\in\Pm$. Dann gelten $\za\in\Pp$ sowie wegen \fref{mabm1bw1} und $S\in\Soqa$ (vergleiche \thref{asmbz3}) weiterhin
\begin{align*}
	-\im S(z) = \frac{-1}{2i}\eklam{S(z)-S^{\ast}(z)} = \im S^{\ast}(z) = \im S(\za)\geq\Oq.
\end{align*}

Zu (b): Wegen Teil (a) von \thref{asmth5} existiert genau ein $\mu \in \Mqmb$ mit
\begin{align*}
	S(z) = \int_{(-\infty,\beta]} \frac{1}{t-z} \;\mu(dt)
\end{align*}
für alle $z \in \C\setminus(-\infty,\beta]$. Wegen \thref{hmllm2} folgt hieraus dann
\begin{align}	\label{mabm1bw2}
	S^{\ast}(z) = \int_{(-\infty,\beta]} \overline{\rklam{\frac{1}{t-z}}} \;\mu(dt)
	= \int_{(-\infty,\beta]} \frac{1}{t-\za} \;\mu(dt)
	= S(\za)
\end{align}
für alle $z \in \C\setminus(-\infty,\beta]$. Sei nun $z\in\Pp$. Dann gelten $\za\in\Pm$ und wegen \fref{mabm1bw2} und $S\in\Soqmb$ (vergleiche \thref{asmbz4}) weiterhin
\begin{align*}
	\im S(z) = \frac{1}{2i}\eklam{S(z)-S^{\ast}(z)} = -\im S^{\ast}(z) = -\im S(\za)\geq\Oq.
\end{align*}

Zu (c): Wegen Teil (a) von \thref{masa1} existiert genau ein $\mu \in \Mqab$ mit
\begin{align*}
	S(z) = \int_{[\alpha,\beta]} \frac{1}{t-z} \;\mu(dt)
\end{align*}
für alle $z \in \C\setminus[\alpha,\beta]$. Wegen \thref{hmllm2} folgt hieraus dann
\begin{align}	\label{mabm1bw3}
	S^{\ast}(z) = \int_{[\alpha,\beta]} \overline{\rklam{\frac{1}{t-z}}} \;\mu(dt)
	= \int_{[\alpha,\beta]} \frac{1}{t-\za} \;\mu(dt)
	= S(\za)
\end{align}
für alle $z \in \C\setminus[\alpha,\beta]$. Sei nun $z\in\Pm$. Dann gelten $\za\in\Pp$ und wegen \fref{mabm1bw3} und $S\in\Rqab$ (vergleiche \thref{mabz1}) weiterhin
\begin{align*}
	-\im S(z) = \frac{-1}{2i}\eklam{S(z)-S^{\ast}(z)} = \im S^{\ast}(z) = \im S(\za)\geq\Oq. \tag*{$\Box$}
\end{align*}

Abschließend zeigen wir, wie wir für eine Funktion aus $\Rqab$ Funktionen aus $\Soqa$ und $\Soqmb$ gewinnen können und wie die zugehörigen Stieltjes-Maße aussehen.

\begin{bem}	\thlabel{mabm2}	
	Seien $\alpha,\beta\in\R$ mit $\alpha<\beta$ und $S:\C\setminus[\alpha,\beta]$. \dgfa
	\begin{itemize}
	\item [\rm{(a)}] \esfaa
		\begin{itemize}
			\item [\rm{(i)}] Es gilt $S\in\Rqab$.
			\item [\rm{(ii)}] Es gelten $\Rstr_{\C\setminus[\alpha,\infty)}S\in\Soqa$ und $\Rstr_{\C\setminus(-\infty,\beta]}S\in\Soqmb$.
		\end{itemize}
	\item [\rm{(a)}] Sei (i) erfüllt. Weiterhin seien $\mu\in\Mqab$ das zu $S$ gehörige Stieltjes-Maß und $\mu_{[\alpha,\infty)}:\B_{[\alpha,\infty)}\rightarrow\Cqq$ bzw. $\mu_{(-\infty,\beta]}:\B_{(-\infty,\beta]}\rightarrow\Cqq$ definiert gemäß
		\begin{align*}
			\mu_{[\alpha,\infty)}(A) := \mu(A\cap[\alpha,\beta])
		\end{align*}
		bzw.
		\begin{align*}
			\mu_{(-\infty,\beta]}(A) := \mu(A\cap[\alpha,\beta]).
		\end{align*}
		Dann gelten $\mu_{[\alpha,\infty)}\in\Mqa$ sowie $\mu_{(-\infty,\beta]}\in\Mqmb$ und es ist $\mu_{[\alpha,\infty)}$ bzw. $\mu_{(-\infty,\beta]}$ das zu $\Rstr_{\C\setminus[\alpha,\infty)}S$ bzw. $\Rstr_{\C\setminus(-\infty,\beta]}S$ gehörige Stieltjes-Maß.
	\end{itemize}
\end{bem}

\bwanf Zu (a): Sei zunächst (i) erfüllt. Hieraus (vergleiche \thref{mabz1}) folgt dann, dass $S$ eine in $\C\setminus[\alpha,\beta]$ holomorphe \textit{q}$\times$\textit{q}-Matrixfunktion ist. Hieraus folgt nun, dass $\Rstr_{\C\setminus[\alpha,\infty)}S$ eine in $\C\setminus[\alpha,\infty)$ bzw. $\Rstr_{\C\setminus(-\infty,\beta]}S$ eine in $\C\setminus(-\infty,\beta]$ holomorphe \textit{q}$\times$\textit{q}-Matrixfunktion ist. Aus $S\in\Rqab$ (vergleiche \thref{mabz1}) folgt weiterhin
\begin{align*}
	\im \eklam{\Rstr_{\C\setminus[\alpha,\infty)}S(z)} = \im S(z) \geq \Oq
\end{align*}
für alle $z\in\Pp$ und
\begin{align*}
	\Rstr_{\C\setminus[\alpha,\infty)}S(x) = S(x) \geq \Oq
\end{align*}
für alle $x\in(-\infty,\alpha)$. Hieraus folgt nun $\Rstr_{\C\setminus[\alpha,\infty)}S\in\Sqa$ (vergleiche \thref{asmbz3}). Wegen Teil (a) von \thref{masa1} und Teil (b) von \thref{hmlbm1} existiert ein $\mu \in \Mqab$ mit Spurmaß $\tau$ und
\begin{align*}
	y\snorm{\Rstr_{\C\setminus[\alpha,\infty)}S(iy)} &= y\snorm{S(iy)} 
	= y\snorm{\int_{[\alpha,\beta]} \frac{1}{t-iy} \;\mu(dt)} \\
	&\leq 4qy\int_{[\alpha,\beta]} \abs{\frac{1}{t-iy}} \tau(dt) 
	= 4qy\int_{[\alpha,\beta]} \frac{1}{\sqrt{t^2+(-y)^2}} \;\tau(dt) \\
	&\leq 4qy\int_{[\alpha,\beta]} \frac{1}{\abs{y}} \;\tau(dt) 
	= 4q\tau([\alpha,\beta])
\end{align*}
für alle $y\in[1,\infty)$, also unter Beachtung, dass $\tau$ ein endliches Maß ist, folgt nun
\begin{align*}
	\sup_{y \in [1,\infty)} y \snorm{\Rstr_{\C\setminus[\alpha,\infty)}S(iy)} \leq 4q\tau([\alpha,\beta]) < \infty.
\end{align*} 
Hieraus folgt dann $\Rstr_{\C\setminus[\alpha,\infty)}S\in\Soqa$ (vergleiche \thref{asmbz3}).
Wegen Teil (c) von \thref{mabm1} bzw. $S\in\Rqab$ (vergleiche \thref{mabz1}) gilt
\begin{align*}
	-\im \eklam{\Rstr_{\C\setminus(-\infty,\beta]}S(z)} = -\im S(z) \geq \Oq
\end{align*}
für alle $z\in\Pm$ bzw.
\begin{align*}
	-\Rstr_{\C\setminus(-\infty,\beta]}S(x) = -S(x) \geq \Oq
\end{align*}
für alle $x\in(\beta,\infty)$. Hieraus folgt dann $\Rstr_{\C\setminus(-\infty,\beta]}S\in\Sqmb$ (vergleiche \thref{asmbz4}). Wegen Teil (a) von \thref{masa1} und Teil (b) von \thref{hmlbm1} existiert ein $\mu \in \Mqab$ mit Spurmaß $\tau$ und
\begin{align*}
	& -y\snorm{\Rstr_{\C\setminus(-\infty,\beta]}S(iy)} = -y\snorm{S(iy)} 
	= -y\snorm{\int_{[\alpha,\beta]} \frac{1}{t-iy} \;\mu(dt)} \\
	&\leq -4qy\int_{[\alpha,\beta]} \abs{\frac{1}{t-iy}} \tau(dt) 
	= -4qy\int_{[\alpha,\beta]} \frac{1}{\sqrt{t^2+(-y)^2}} \;\tau(dt) \\
	&\leq -4qy\int_{[\alpha,\beta]} \frac{1}{\abs{y}} \;\tau(dt) 
	= 4q\tau([\alpha,\beta])
\end{align*}
für alle $y\in(-\infty,-1]$, also unter Beachtung, dass $\tau$ ein endliches Maß ist, folgt nun
\begin{align*}
	\sup_{y \in (-\infty,-1]} -y \snorm{\Rstr_{\C\setminus(-\infty,\beta])}S(iy)} \leq 4q\tau([\alpha,\beta]) < \infty.
\end{align*} 
Hieraus folgt dann $\Rstr_{\C\setminus(-\infty,\beta]}S\in\Soqmb$ (vergleiche \thref{asmbz4}). Somit gilt also (ii).

Sei nun umgekehrt (ii) erfüllt. Hieraus (vergleiche Bezeichnungen \ref{asmbz3} und \ref{asmbz4}) folgt dann, dass $\Rstr_{\C\setminus[\alpha,\infty)}S$ eine in $\C\setminus[\alpha,\infty)$ bzw. $\Rstr_{\C\setminus(-\infty,\beta]}S$ eine in $\C\setminus(-\infty,\beta]$ holomorphe \textit{q}$\times$\textit{q}-Matrixfunktion ist. Unter Beachtung von 
\begin{align*}
	\C\setminus[\alpha,\beta] = \brklam{\C\setminus[\alpha,\infty)} \cup \brklam{\C\setminus(-\infty,\beta]}
\end{align*}
folgt, dass $S$ eine in $\C\setminus[\alpha,\beta]$ holomorphe \textit{q}$\times$\textit{q}-Matrixfunktion ist. Aus \linebreak $\Rstr_{\C\setminus[\alpha,\infty)}S\in\Soqa$ (vergleiche \thref{asmbz3}) folgt weiterhin
\begin{align}	\label{mabm2bw1}
	\im S(z) = \im \eklam{\Rstr_{\C\setminus[\alpha,\infty)}S(z)} \geq \Oq
\end{align}
für alle $z\in\Pp$ und
\begin{align}	\label{mabm2bw2}
	S(x) = \Rstr_{\C\setminus[\alpha,\infty)}S(x) \geq \Oq
\end{align}
für alle $x\in(-\infty,\alpha)$.
Aus $\Rstr_{\C\setminus(-\infty,\beta]}S\in\Soqmb$ (vergleiche \thref{asmbz4}) folgt weiterhin
\begin{align}	\label{mabm2bw3}
	-S(x) = -\Rstr_{\C\setminus(-\infty,\beta]}S(x) \geq \Oq
\end{align}
für alle $x\in(\beta,\infty)$. Wegen \fref{mabm2bw1}, \fref{mabm2bw2} und \fref{mabm2bw3} folgt dann $S\in\Rqab$ (vergleiche \thref{mabz1}).

Zu (b): Wegen Teil (a) von \thref{hmlbm2} gelten $\mu_{[\alpha,\infty)}\in\Mqa$ und \linebreak $\mu_{(-\infty,\beta]}\in\Mqmb$. Wegen Teil (b) von \thref{hmlbm2} und \thref{madef1} gelten weiterhin
\begin{align*}
	\int_{[\alpha,\infty)} \frac{1}{t-z} \;\mu_{[\alpha,\infty)}(dt) 
	= \int_{[\alpha,\beta]} \frac{1}{t-z} \;\mu(dt)
	= S(z)
	= \Rstr_{\C\setminus[\alpha,\infty)}S(z)
\end{align*}
für alle $z\in\C\setminus[\alpha,\infty)$ und
\begin{align*}
	\int_{(-\infty,\beta]} \frac{1}{t-z} \;\mu_{(-\infty,\beta]}(dt) 
	= \int_{[\alpha,\beta]} \frac{1}{t-z} \;\mu(dt)
	= S(z)
	= \Rstr_{\C\setminus(-\infty,\beta]}S(z)
\end{align*}
für alle $z\in\C\setminus(-\infty,\beta]$. Hieraus folgt wegen \thref{madef1} dann, dass $\mu_{[\alpha,\infty)}$ bzw. $\mu_{(-\infty,\beta]}$ das zu $\Rstr_{\C\setminus[\alpha,\infty)}S$ bzw. $\Rstr_{\C\setminus(-\infty,\beta]}S$ gehörige Stieltjes-Maß ist. \bwend

\section{Einige Aussagen der Matrizentheorie}

In diesem Abschnitt geben wir einige ausgewählte Resultate der Matrizentheorie an. 

Zunächst betrachten wir zwei bestimmte Produkte von positiv hermiteschen Matrizenfolgen, die für die in Kapitel \ref{chapadp} behandelten Zusammenhänge zwischen der $\alpha$-Dyukarev-Stieltjes-Parametrisierung und der $\alpha$-Stieltjes-Parametrisierung von Wichtigkeit sind.

\begin{lemma}	\thlabel{amlm1}
  Seien $\kappa \in \Noa$ sowie $(X_{n})^{\kappa}_{n=0}$ bzw. $(X_{n})^{\kappa+1}_{n=0}$ und $(Y_{n})^{\kappa}_{n=0}$ Folgen aus $\Cqq_{>}$. 
  Weiterhin seien
  \begin{align*}
    A_{n} := \rklam{\prodr^{n}_{j=0}X_{j}Y^{-1}_{j}}Y_{n}\rklam{\prodr^{n}_{j=0}X_{j}Y^{-1}_{j}}^{\ast}
  \end{align*}
  für alle $n \in \Zok$ und
  \begin{align*}
    B_{n} := \begin{cases} X^{-1}_{0} & \text{falls } n=0 \\
    \rklam{\prodr^{n-1}_{j=0}X^{-1}_{j}Y_{j}}X^{-1}_{n}\rklam{\prodr^{n-1}_{j=0}X^{-1}_{j}Y_{j}}^{\ast} & \text{falls } n \geq 1 \end{cases}
  \end{align*}
  für alle $n \in \Zok$ bzw. $n \in \Zokp$. Dann sind $(A_{n})^{\kappa}_{n=0}$ sowie $(B_{n})^{\kappa}_{n=0}$ bzw. $(B_{n})^{\kappa+1}_{n=0}$ Folgen aus $\Cqq_{>}$ und es gelten
  \begin{align*}
    X_{n} = \begin{cases} B^{-1}_{0} & \text{falls } n=0 \\
    \rklam{\prodr^{n-1}_{j=0}B_{j}A_{j}}^{-\ast}B^{-1}_{n}\rklam{\prodr^{n-1}_{j=0}B_{j}A_{j}}^{-1} & \text{falls } n \geq 1 \end{cases}
  \end{align*}
  für alle $n \in \Zok$ bzw. $n \in \Zokp$ sowie
  \begin{align*}
    Y_{n} = \rklam{\prodr^{n}_{j=0}B_{j}A_{j}}^{-\ast}A_{n}\rklam{\prodr^{n}_{j=0}B_{j}A_{j}}^{-1}
  \end{align*}
  für alle $n \in \Zok$.
\end{lemma}

\bwanf Der Fall $\kappa = \infty$ wurde in \cite[Lemma 8.7]{Trans} behandelt. Der Beweis im endlichen Fall verläuft analog. \bwend

\begin{lemma}	\thlabel{amlm3}
  Seien $\kappa \in \Noa$ sowie $(A_{n})^{\kappa}_{n=0}$ und $(B_{n})^{\kappa}_{n=0}$ bzw. $(B_{n})^{\kappa+1}_{n=0}$ Folgen aus $\Cqq_{>}$. Weiterhin seien
  \begin{align*}
    X_{n} := \begin{cases} B^{-1}_{0} & \text{falls } n=0 \\
    \rklam{\prodr^{n-1}_{j=0}B_{j}A_{j}}^{-\ast}B^{-1}_{n}\rklam{\prodr^{n-1}_{j=0}B_{j}A_{j}}^{-1} & \text{falls } n \geq 1 \end{cases}
  \end{align*}
  für alle $n \in \Zok$ bzw. $n \in \Zokp$ und
  \begin{align*}
    Y_{n} := \rklam{\prodr^{n}_{j=0}B_{j}A_{j}}^{-\ast}A_{n}\rklam{\prodr^{n}_{j=0}B_{j}A_{j}}^{-1}
  \end{align*}
  für alle $n \in \Zok$. Dann sind $(X_{n})^{\kappa}_{n=0}$ bzw. $(X_{n})^{\kappa+1}_{n=0}$ sowie $(Y_{n})^{\kappa}_{n=0}$ Folgen aus $\Cqq_{>}$ und es gelten
  \begin{align*}
    A_{n} = \rklam{\prodr^{n}_{j=0}X_{j}Y^{-1}_{j}}Y_{n}\rklam{\prodr^{n}_{j=0}X_{j}Y^{-1}_{j}}^{\ast}
  \end{align*}
  für alle $n \in \Zok$ sowie
  \begin{align*}
    B_{n} = \begin{cases} X^{-1}_{0} & \text{falls } n=0 \\
    \rklam{\prodr^{n-1}_{j=0}X^{-1}_{j}Y_{j}}X^{-1}_{n}\rklam{\prodr^{n-1}_{j=0}X^{-1}_{j}Y_{j}}^{\ast} & \text{falls } n \geq 1 \end{cases}
  \end{align*}
  für alle $n \in \Zok$ bzw. $n \in \Zokp$.  
\end{lemma}

\bwanf Der Fall $\kappa = \infty$ wurde in \cite[Lemma 8.5]{Trans} behandelt. Der Beweis im endlichen Fall verläuft analog. \bwend

Unser nächstes Resultat behandelt Blockmatrizen in Verbindung mit dem Begriff der nichtnegativen Hermitizität. Hierbei spielt das linke bzw. rechte Schurkomplement eine tragende Rolle. Diese Aussage geht auf Albert \cite{A} und Efimov/Potapov \cite{EP} zurück.

\begin{lemma}	\thlabel{amlm2}
  Seien $A \in \Cpp$, $B \in \Cpq$, $C \in \Cqp$, $D \in \Cqq$ und
  \begin{align*}
    E := \begin{pmatrix} A & B \\ C & D \end{pmatrix}.
  \end{align*}
  \dsfaa
  \begin{itemize}
    \item [\rm{(i)}] Es gilt $E \in \C^{(p+q)\times(p+q)}_{\geq}$.
    \item [\rm{(ii)}] Es gelten $A \in \Cpp_{\geq}$, $AA^{+}B = B$, $C = B^{\ast}$ und $D-CA^{+}B \in \Cqq_{\geq}$.
    \item [\rm{(iii)}] Es gelten $D \in \Cqq_{\geq}$, $DD^{+}C = C$, $B = C^{\ast}$ und $A-BD^{+}C \in \Cpp_{\geq}$.
  \end{itemize}
\end{lemma}

\bwanf Siehe \cite[Lemma 1.1.9]{DFK}. \bwend

Aus \thref{amlm2} lässt sich ein analoges Resultat für den positiv hermiteschen Fall herleiten.

\begin{lemma}	\thlabel{amlm4}
  Seien $A \in \Cpp$, $B \in \Cpq$, $C \in \Cqp$, $D \in \Cqq$ und
  \begin{align*}
    E := \begin{pmatrix} A & B \\ C & D \end{pmatrix}.
  \end{align*}
  \dsfaa
  \begin{itemize}
    \item [\rm{(i)}] Es gilt $E \in \C^{(p+q)\times(p+q)}_{>}$.
    \item [\rm{(ii)}] Es gelten $A \in \Cpp_{>}$, $C = B^{\ast}$ und $D-CA^{-1}B \in \Cqq_{>}$.
    \item [\rm{(iii)}] Es gelten $D \in \Cqq_{>}$, $B = C^{\ast}$ und $A-BD^{-1}C \in \Cpp_{>}$.
  \end{itemize}
\end{lemma}

\bwanf Siehe \cite[Lemma 1.1.9]{DFK}. \bwend

Abschließend widmen wir uns den Matrixintervallen bezüglich der Löwner-Halbordnung und zeigen einen Zusammenhang zu Matrizenkreisen. Ein ähnliches Resultat findet man in \cite[Lemma 10.1]{HMS}.

\begin{bem}	\thlabel{ambm1}
	Seien $A,B\in\Cqq_H$ mit $A\leq B$. \dgfa
	\begin{itemize}
		\item [\rm{(a)}] Es gilt
		\begin{align*}
			[A,B] = \bgklam{A+\sqrt{B-A}K\sqrt{B-A} \;\big|\; K\in[\Oq,\Iq]}.
		\end{align*}
		\item [\rm{(b)}] Es gilt
		\begin{align*}
			[A,B] = \bgklam{B-\sqrt{B-A}K\sqrt{B-A} \;\big|\; K\in[\Oq,\Iq]}.
		\end{align*}
	\end{itemize}
\end{bem}

\bwanf Zu (a): Seien zunächst $K\in[\Oq,\Iq]$ und $X:=A+\sqrt{B-A}K\sqrt{B-A}$. Dann gelten $X\in\Cqq_H$ sowie
\begin{align*}
	X-A = \sqrt{B-A}K\sqrt{B-A} \geq \Oq
\end{align*}
und
\begin{align*}
	B-X = B-A-\sqrt{B-A}K\sqrt{B-A} = \sqrt{B-A}(\Iq-K)\sqrt{B-A} \geq \Oq,
\end{align*}
also $X\in[A,B]$.

Seien nun $X\in[A,B]$ und $K:= \sqrt{B-A}^{+}(X-A)\sqrt{B-A}^{+}$. Dann gelten $K\in\Cqq_H$ und
\begin{align}	\label{ambm1bw1}
	B-A\geq X-A\geq \Oq,
\end{align}
also $K\geq\Oq$. Aus \fref{ambm1bw1} folgen weiterhin ${\cal R}(X-A)\subseteq{\cal R}(B-A)$ und \linebreak
${\cal N}(B-A)\subseteq{\cal N}(X-A)$. Hieraus folgt dann
\begin{align*}
	\sqrt{B-A}K\sqrt{B-A} = \sqrt{B-A}\sqrt{B-A}^{+}(X-A)\sqrt{B-A}^{+}\sqrt{B-A} = X-A,
\end{align*}
also $X = A+\sqrt{B-A}K\sqrt{B-A}$. Sei $P_{{\cal R}(\sqrt{B-A})}$ die Orthoprojektionsmatrix von $\Cq$ auf ${\cal R}(\sqrt{B-A})$. Dann gelten $P_{{\cal R}(\sqrt{B-A})}=\sqrt{B-A}\sqrt{B-A}^{+}$ und $P_{{\cal R}(\sqrt{B-A})}\leq\Iq$. Hieraus folgt wegen \fref{ambm1bw1} und der Definition der Moore-Penrose-Inverse nun
\begin{align*}
	K \leq \sqrt{B-A}^{+}(B-A)\sqrt{B-A}^{+} = \sqrt{B-A}\sqrt{B-A}^{+} \leq \Iq,
\end{align*}
also $K\in[\Oq,\Iq]$.

Zu (b): Seien zunächst $K\in[\Oq,\Iq]$ und $X:=B-\sqrt{B-A}K\sqrt{B-A}$. Dann gelten $X\in\Cqq_H$ sowie
\begin{align*}
	B-X = \sqrt{B-A}K\sqrt{B-A} \geq \Oq
\end{align*}
und
\begin{align*}
	X-A = B-A-\sqrt{B-A}K\sqrt{B-A} = \sqrt{B-A}(\Iq-K)\sqrt{B-A} \geq \Oq,
\end{align*}
also $X\in[A,B]$.

Seien nun $X\in[A,B]$ und $K:= \sqrt{B-A}^{+}(B-X)\sqrt{B-A}^{+}$. Dann gelten $K\in\Cqq_H$ und
\begin{align}	\label{ambm1bw2}
	B-A\geq B-X\geq \Oq,
\end{align}
also $K\geq\Oq$. Aus \fref{ambm1bw2} folgen weiterhin ${\cal R}(B-X)\subseteq{\cal R}(B-A)$ und \linebreak
${\cal N}(B-A)\subseteq{\cal N}(B-X)$. Hieraus folgt dann
\begin{align*}
	\sqrt{B-A}K\sqrt{B-A} = \sqrt{B-A}\sqrt{B-A}^{+}(B-X)\sqrt{B-A}^{+}\sqrt{B-A} = B-X,
\end{align*}
also $X = B-\sqrt{B-A}K\sqrt{B-A}$. Sei $P_{{\cal R}(\sqrt{B-A})}$ die Orthoprojektionsmatrix von $\Cq$ auf ${\cal R}(\sqrt{B-A})$. Dann gelten $P_{{\cal R}(\sqrt{B-A})}=\sqrt{B-A}\sqrt{B-A}^{+}$ und $P_{{\cal R}(\sqrt{B-A})}\leq\Iq$. Hieraus folgt wegen \fref{ambm1bw2} und der Definition der Moore-Penrose-Inverse nun
\begin{align*}
	K \leq \sqrt{B-A}^{+}(B-A)\sqrt{B-A}^{+} = \sqrt{B-A}\sqrt{B-A}^{+} \leq \Iq,
\end{align*}
also $K\in[\Oq,\Iq]$. \bwend

\section{Einige Aussagen der \textit{J}-Theorie}  \label{chapsm}

In diesem Abschnitt stellen wir einige für unsere Betrachtungen bedeutsame Resultate der $J$-Theorie zusammen.

Zuerst führen wir den Begriff der Signaturmatrix ein und wollen dann die Definitionen von kontraktiven, expansiven und unitären Matrizen mithilfe der Signaturmatrizen verallgemeinern. Detailliertere Ausführungen hierzu findet man z.\,B. in \cite[Chapter 1.3.]{DFK}.

\begin{defi}	\thlabel{spdef1}
  Sei $J \in \Cpp$. Dann heißt $J$ \textbf{p$\times$p-Signaturmatrix}, falls $J = J^{\ast}$ und $J^2 = \Ip$ erfüllt sind.
\end{defi}

Aus \thref{spdef1} erkennt man sogleich, dass die Einheitsmatrix $\Ip$ eine \textit{p}$\times$\textit{p}"=Signaturmatrix ist.

\begin{defi}	\thlabel{spdef2} Seien $J$ eine \textit{p}$\times$\textit{p}-Signaturmatrix und $A \in \Cpp$.
\begin{itemize}
  \item [\rm{(a)}] $A$ heißt \textbf{J-kontraktiv} bzw. \textbf{J-expansiv}, falls $J-A^{\ast}JA \geq \Op$ bzw. \linebreak $A^{\ast}JA-J \geq \Op$ erfüllt ist.
  \item [\rm{(b)}] $A$ heißt \textbf{J-unitär}, falls $J-A^{\ast}JA = \Op$ erfüllt ist.
\end{itemize}
\end{defi}

Es sei bemerkt, dass man in \thref{spdef2} mit der Setzung $J=\Iq$ die Definitionen für kontraktive, expansive und unitäre Matrizen erhält. Wir betrachten nun einige grundlegende Eigenschaften für die so eben eingeführten Begriffe.

\begin{bem}	\thlabel{spbm1} Seien $J$ eine \textit{p}$\times$\textit{p}-Signaturmatrix und $A,B \in \Cpp$. \dgfa
\begin{itemize}
  \item [\rm{(a)}] \esfaa
  \begin{itemize}
    \item [\rm{(i)}]  $A$ ist $J$-kontraktiv bzw. $J$-expansiv.
    \item [\rm{(ii)}]  $A^{\ast}$ ist $J$-kontraktiv bzw. $J$-expansiv.
  \end{itemize}
  \item [\rm{(b)}] Sei $A$ regulär. \dsfaa
  \begin{itemize}
    \item [\rm{(iii)}]  $A$ ist $J$-kontraktiv.
    \item [\rm{(iv)}]  $A^{-1}$ ist $J$-expansiv.
  \end{itemize}
  \item [\rm{(c)}] Sei $A$ $J$-unitär. Dann ist $A$ regulär sowie $A^{\ast}$ und $A^{-1}$ $J$-unitär. 
  \item [\rm{(d)}] Seien $A$ und $B$ $J$-kontraktiv bzw. $J$-expansiv bzw. $J$-unitär. Dann ist $AB$ $J$-kontraktiv bzw. $J$-expansiv bzw. $J$-unitär.
\end{itemize}
\end{bem}

\bwanf Zu (a): Siehe Teile (a) und (c) von \cite[Theorem 1.3.3]{DFK}.

Zu (b): Siehe Teil (a) und (b) von \cite[Lemma 1.3.15]{DFK}.

Zu (c): Siehe Teil (e) von \cite[Theorem 1.3.3]{DFK} und Teil (c) von \cite[Lemma 1.3.15]{DFK}.

Zu (d): Siehe \cite[Lemma 1.3.13]{DFK}. \bwend

In dieser Arbeit verwenden wir oftmals folgende speziell für unser Vorgehen relevante Signaturmatrix.

\begin{beispiel}	\thlabel{spbsp1} Sei
\begin{align*}
  \tJq := \begin{pmatrix} \Oq & -i\Iq \\ i\Iq & \Oq \end{pmatrix}.
\end{align*}
Dann ist $\tJq$ eine $2$\textit{q}$\times2$\textit{q}-Signaturmatrix und es gilt
\begin{align*}
  \begin{pmatrix} A \\ B \end{pmatrix}^{\ast} \brklam{{-\tJq}} \begin{pmatrix} A \\ B \end{pmatrix} = 2 \im \eklam{B^{\ast}A}
\end{align*}
für alle $A,B \in \Cqq$.
\end{beispiel}

\bwanf Aus der Definition von $\tJq$ folgen sogleich $\tJq = \tJq^{\ast}$ und $\tJq^2 = I_{2q}$. Somit ist wegen \thref{spdef1} dann $\tJq$ eine $2$\textit{q}$\times2$\textit{q}-Signaturmatrix. Weiterhin gilt
\begin{align*}
  	\begin{pmatrix} A \\ B \end{pmatrix}^{\ast} \brklam{{-\tJq}} \begin{pmatrix} A \\ B \end{pmatrix} 
  	= \begin{pmatrix} A \\ B \end{pmatrix}^{\ast} \begin{pmatrix} iB \\ -iA \end{pmatrix}  
  	= iA^{\ast}B-iB^{\ast}A
  	= \frac{1}{i}\eklam{B^{\ast}A-A^{\ast}B}
  	= 2 \im \eklam{B^{\ast}A}
\end{align*}
für alle $A,B \in \Cqq$. \bwend

Wir betrachten nun gewisse Blockmatrizen von Dreiecksgestalt und deren Verhalten im Bezug zur so eben eingeführten Signaturmatrix $\tJq$ (vergleiche \cite[Remark 6.6]{C06}).

\begin{bem}	\thlabel{spbm2}
  Seien $A \in \Cqq$ sowie
  \begin{align*}
    M_1 := \begin{pmatrix} \Iq & \Oq \\ A & \Iq \end{pmatrix} \quad \text{und} \quad M_2 := \begin{pmatrix} \Iq & A \\ \Oq & \Iq \end{pmatrix}.
  \end{align*}
  \dgfa
  \begin{itemize}
    \item [\rm{(a)}] Es gilt $M^{\ast}_1\tJq M_1 = \tJq + \diag(i(A^{\ast}-A), \Oq)$.
    \item [\rm{(b)}] Es gilt $M^{\ast}_2\tJq M_2 = \tJq + \diag(\Oq, i(A-A^{\ast}))$.
    \item [\rm{(c)}] \esfaa
    \begin{itemize}
      \item [\rm{(i)}] Es gilt $A^{\ast} = A$.
      \item [\rm{(ii)}] $M_1$ ist $\tJq$-unitär.
      \item [\rm{(iii)}] $M_2$ ist $\tJq$-unitär.
    \end{itemize}
  \end{itemize}
\end{bem}

\bwanf 
Zu (a): Es gilt
\begin{align*}
	M^{\ast}_1\tJq M_1 
	&= \begin{pmatrix} \Iq & A^{\ast} \\ \Oq & \Iq \end{pmatrix} \begin{pmatrix} -iA & -i\Iq \\ i\Iq & \Oq \end{pmatrix}
	= \begin{pmatrix} i(A^{\ast}-A) & -i\Iq \\ i\Iq & \Oq \end{pmatrix} \\
	&= \tJq + \diag(i(A^{\ast}-A), \Oq).
\end{align*}

Zu (b): Es gilt
\begin{align*}
	M^{\ast}_2\tJq M_2 
	&= \begin{pmatrix} \Iq & \Oq \\ A^{\ast} & \Iq \end{pmatrix} \begin{pmatrix} \Oq & -i\Iq \\ i\Iq & iA \end{pmatrix}
	= \begin{pmatrix} \Oq & -i\Iq \\ i\Iq & i(A-A^{\ast}) \end{pmatrix} \\
	&= \tJq + \diag(\Oq, i(A-A^{\ast})).
\end{align*}

Zu (c): Dies folgt wegen \thref{spbsp1} und Teil (b) von \thref{spdef2} sogleich aus (a) und (b). \bwend

Wir kommen nun auf eine weitere Signaturmatrix zu sprechen, welche in dieser Arbeit von Bedeutung ist.

\begin{beispiel}	\thlabel{spbsp4} Sei
\begin{align*}
  \Jq := \begin{pmatrix} \Oq & -\Iq \\ -\Iq & \Oq \end{pmatrix}.
\end{align*}
Dann ist $\Jq$ eine $2$\textit{q}$\times2$\textit{q}-Signaturmatrix und es gilt
\begin{align*}
  \begin{pmatrix} A \\ B \end{pmatrix}^{\ast} \rklam{-\Jq} \begin{pmatrix} A \\ B \end{pmatrix} = 2 \re \eklam{B^{\ast}A}
\end{align*}
für alle $A,B \in \Cqq$.
\end{beispiel}

\bwanf Aus der Definition von $\Jq$ folgen sogleich $\Jq = \Jq^{\ast}$ und $\Jq^2 = I_{2q}$. Somit ist wegen \thref{spdef1} dann $\Jq$ eine $2$\textit{q}$\times2$\textit{q}-Signaturmatrix. Weiterhin gilt
\begin{align*}
  	\begin{pmatrix} A \\ B \end{pmatrix}^{\ast} \rklam{-\Jq} \begin{pmatrix} A \\ B \end{pmatrix} 
  	= \begin{pmatrix} A \\ B \end{pmatrix}^{\ast} \begin{pmatrix} B \\ A \end{pmatrix}  
  	= A^{\ast}B+B^{\ast}A
  	= 2 \re \eklam{B^{\ast}A}
\end{align*}
für alle $A,B \in \Cqq$. \bwend

Folgend betrachten wir ein Analogon von \thref{spbm2} für die so eben eingeführte Signaturmatrix $\Jq$.

\begin{bem}	\thlabel{spbm5}
  Seien $A \in \Cqq$ sowie
  \begin{align*}
    M_1 := \begin{pmatrix} \Iq & \Oq \\ A & \Iq \end{pmatrix} \quad \text{und} \quad M_2 := \begin{pmatrix} \Iq & A \\ \Oq & \Iq \end{pmatrix}.
  \end{align*}
  \dgfa
  \begin{itemize}
    \item [\rm{(a)}] Es gilt $M^{\ast}_1\Jq M_1 = \Jq - \diag(A+A^{\ast}, \Oq)$.
    \item [\rm{(b)}] Es gilt $M^{\ast}_2\Jq M_2 = \Jq - \diag(\Oq, A+A^{\ast})$.
    \item [\rm{(c)}] \esfaa
    \begin{itemize}
      \item [\rm{(i)}] Es gilt $A\in\Cqq_\geq$.
      \item [\rm{(ii)}] $M_1$ ist $\Jq$-kontraktiv.
      \item [\rm{(iii)}] $M_2$ ist $\Jq$-kontraktiv.
    \end{itemize}
  \end{itemize}
\end{bem}

\bwanf 
Zu (a): Es gilt
\begin{align*}
	M^{\ast}_1\Jq M_1 
	&= \begin{pmatrix} \Iq & A^{\ast} \\ \Oq & \Iq \end{pmatrix} \begin{pmatrix} -A & -\Iq \\ -\Iq & \Oq \end{pmatrix}
	= \begin{pmatrix} -A-A^{\ast} & -\Iq \\ -\Iq & \Oq \end{pmatrix} \\
	&= \Jq - \diag(A+A^{\ast}, \Oq).
\end{align*}

Zu (b): Es gilt
\begin{align*}
	M^{\ast}_2\Jq M_2 
	&= \begin{pmatrix} \Iq & \Oq \\ A^{\ast} & \Iq \end{pmatrix} \begin{pmatrix} \Oq & -\Iq \\ -\Iq & -A \end{pmatrix}
	= \begin{pmatrix} \Oq & -\Iq \\ -\Iq & -A^{\ast}-A \end{pmatrix} \\
	&= \Jq - \diag(\Oq, A+A^{\ast}).
\end{align*}

Zu (c): Wegen (a) bzw. (b) gilt
\begin{align*}
	\Jq-M^{\ast}_1\Jq M_1=\diag(A+A^{\ast}, \Oq)
\end{align*}
bzw.
\begin{align*}
	\Jq-M^{\ast}_2\Jq M_2=\diag(\Oq, A+A^{\ast}).
\end{align*}
Hieraus folgt wegen \thref{spbsp4} und Teil (a) von \thref{spdef2} dann die Behauptung. \bwend

Wir kommen nun auf eine letzte Signaturmatrix zu sprechen, die in dieser Arbeit Verwendung findet.

\begin{beispiel}	\thlabel{spbsp5} Sei $\jqq := \diag(\Iq,-\Iq)$.
Dann ist $\jqq$ eine $2$\textit{q}$\times2$\textit{q}-Signaturmatrix.
\end{beispiel}

\bwanf Aus der Definition von $\jqq$ folgen sogleich $\jqq = \jqq^{\ast}$ und $\jqq^2 = I_{2q}$. Somit ist wegen \thref{spdef1} dann $\jqq$ eine $2$\textit{q}$\times2$\textit{q}-Signaturmatrix. \bwend

Folgende Definition verallgemeinert den Begriff der $J$-kontraktiven, $J$-expansiven und $J$-unitären Matrizen unter Verwendung einer zweiten Signaturmatrix.

\begin{defi}	\thlabel{spdef7} Seien $J^{(1)}$ und $J^{(2)}$ jeweils eine \textit{p}$\times$\textit{p}-Signaturmatrix sowie $A \in \Cpp$.
\begin{itemize}
  \item [\rm{(a)}] $A$ heißt \textbf{J$\mathbf{^{(2)}}$-J$\mathbf{^{(1)}}$-kontraktiv} bzw. \textbf{J$\mathbf{^{(2)}}$-J$\mathbf{^{(1)}}$-expansiv}, falls $J^{(2)}-A^{\ast}J^{(1)}A \geq \Op$ bzw. $A^{\ast}J^{(1)}A-J^{(2)} \geq \Op$ erfüllt ist.
  \item [\rm{(b)}] $A$ heißt \textbf{J$\mathbf{^{(2)}}$-J$\mathbf{^{(1)}}$-unitär}, falls $J^{(2)}-A^{\ast}J^{(1)}A = \Op$ erfüllt ist.
\end{itemize}
\end{defi}

Nimmt man in \thref{spdef7} die Setzung $J^{(1)} = J^{(2)}$ vor, so erhält man die Begriffsbildungen aus \thref{spdef2}. Wir betrachten nun zwei spezielle Resultate für die so eben eingeführten Begriffsbildungen, die die in Kapitel \ref{chapar} vorliegende Situation umfassen.

\begin{bem}	\thlabel{spbm7}
	Sei
	\begin{align*}
		E:=\frac{1}{\sqrt{2}}\begin{pmatrix} -i\Iq & i\Iq \\ \Iq & \Iq \end{pmatrix}.
	\end{align*}
	\dgfa
	\begin{itemize}
		\item [\rm{(a)}] Es ist $E$ unitär.
		\item [\rm{(b)}] Es gilt $\jqq = E^{\ast}\tJq E$.
		\item [\rm{(c)}] Seien $A\in\C^{2q\times2q}$ und $B:=AE$. \dgfa
		\begin{itemize}
			\item [\rm{(c1)}] Es gilt
			\begin{align*}
				\jqq-B^{\ast}\tJq B = E^{\ast}\brklam{\tJq-A^{\ast}\tJq A}E.
			\end{align*}
			\item [\rm{(c2)}] Es ist $B$ genau dann $\jqq$-$\tJq$-kontraktiv bzw. $\jqq$-$\tJq$-expansiv, wenn $A$ eine $\tJq$-kontraktive bzw. $\tJq$-expansive Matrix ist.
			\item [\rm{(c3)}] Es ist $B$ genau dann $\jqq$-$\tJq$-unitär, wenn $A$ eine $\tJq$-unitäre Matrix ist.
		\end{itemize}
	\end{itemize}
\end{bem}

\bwanf Zu (a): Es gilt
\begin{align*}
	E^{\ast}E = \frac{1}{\sqrt{2}}\begin{pmatrix} i\Iq & \Iq \\ -i\Iq & \Iq \end{pmatrix}\frac{1}{\sqrt{2}}\begin{pmatrix} -i\Iq & i\Iq \\ \Iq & \Iq \end{pmatrix} 
	= \frac{1}{2}\begin{pmatrix} 2\Iq & \Oq \\ \Oq & 2\Iq \end{pmatrix} = I_{2q\times2q},
\end{align*}
also ist $E$ unitär.

Zu (b): Es gilt
\begin{align*}
	E^{\ast}\tJq E &= \frac{1}{\sqrt{2}}\begin{pmatrix} i\Iq & \Iq \\ -i\Iq & \Iq \end{pmatrix}\begin{pmatrix} \Oq & -i\Iq \\ i\Iq & \Oq \end{pmatrix}\frac{1}{\sqrt{2}}\begin{pmatrix} -i\Iq & i\Iq \\ \Iq & \Iq \end{pmatrix} \\
	&= \frac{1}{2}\begin{pmatrix} i\Iq & \Iq \\ -i\Iq & \Iq \end{pmatrix}\begin{pmatrix} -i\Iq & -i\Iq \\ \Iq & -\Iq \end{pmatrix}
	= \frac{1}{2}\begin{pmatrix} 2\Iq & \Oq \\ \Oq & -2\Iq \end{pmatrix} = \jqq.
\end{align*}

Zu (c1): Wegen (b) gilt
\begin{align*}
	E^{\ast}\brklam{\tJq-A^{\ast}\tJq A}E
	= E^{\ast}\tJq E-(AE)^{\ast}\tJq AE
	= \jqq-B^{\ast}\tJq B.
\end{align*}

Zu (c2): Dies folgt wegen (a) und (c1) aus \thref{spbsp1}, \thref{spbsp5}, Teil (a) von \thref{spdef2} und Teil (a) von \thref{spdef7}.

Zu (c3): Dies folgt wegen (a) und (c1) aus \thref{spbsp1}, \thref{spbsp5}, Teil (b) von \thref{spdef2} und Teil (b) von \thref{spdef7}. \bwend

\begin{bem}	\thlabel{spbm10}
	Sei
	\begin{align*}
		\widetilde{E}:=\frac{1}{\sqrt{2}}\begin{pmatrix} -i\Iq & -i\Iq \\ -\Iq & \Iq \end{pmatrix}.
	\end{align*}
	\dgfa
	\begin{itemize}
		\item [\rm{(a)}] Es ist $\widetilde{E}$ unitär.
		\item [\rm{(b)}] Es gilt $\jqq = \widetilde{E}^{\ast}\tJq \widetilde{E}$.
		\item [\rm{(c)}] Seien $A\in\C^{2q\times2q}$ und $B:=A\widetilde{E}$. \dgfa
		\begin{itemize}
			\item [\rm{(c1)}] Es gilt
			\begin{align*}
				\jqq-B^{\ast}\tJq B = \widetilde{E}^{\ast}\brklam{\tJq-A^{\ast}\tJq A}\widetilde{E}.
			\end{align*}
			\item [\rm{(c2)}] Es ist $B$ genau dann $\jqq$-$\tJq$-kontraktiv bzw. $\jqq$-$\tJq$-expansiv, wenn $A$ eine $\tJq$-kontraktive bzw. $\tJq$-expansive Matrix ist.
			\item [\rm{(c3)}] Es ist $B$ genau dann $\jqq$-$\tJq$-unitär, wenn $A$ eine $\tJq$-unitäre Matrix ist.
		\end{itemize}
	\end{itemize}
\end{bem}

\bwanf Zu (a): Es gilt
\begin{align*}
	\widetilde{E}^{\ast}\widetilde{E} = \frac{1}{\sqrt{2}}\begin{pmatrix} i\Iq & -\Iq \\ i\Iq & \Iq \end{pmatrix}\frac{1}{\sqrt{2}}\begin{pmatrix} -i\Iq & -i\Iq \\ -\Iq & \Iq \end{pmatrix} 
	= \frac{1}{2}\begin{pmatrix} 2\Iq & \Oq \\ \Oq & 2\Iq \end{pmatrix} = I_{2q\times2q},
\end{align*}
also ist $E$ unitär.

Zu (b): Es gilt
\begin{align*}
	\widetilde{E}^{\ast}\tJq \widetilde{E} &= \frac{1}{\sqrt{2}}\begin{pmatrix} i\Iq & -\Iq \\ i\Iq & \Iq \end{pmatrix}\begin{pmatrix} \Oq & -i\Iq \\ i\Iq & \Oq \end{pmatrix}\frac{1}{\sqrt{2}}\begin{pmatrix} -i\Iq & -i\Iq \\ -\Iq & \Iq \end{pmatrix} \\
	&= \frac{1}{2}\begin{pmatrix} i\Iq & -\Iq \\ i\Iq & \Iq \end{pmatrix}\begin{pmatrix} i\Iq & -i\Iq \\ \Iq & \Iq \end{pmatrix}
	= \frac{1}{2}\begin{pmatrix} 2\Iq & \Oq \\ \Oq & -2\Iq \end{pmatrix} = \jqq.
\end{align*}

Zu (c1): Wegen (b) gilt
\begin{align*}
	\widetilde{E}^{\ast}\brklam{\tJq-A^{\ast}\tJq A}\widetilde{E}
	= \widetilde{E}^{\ast}\tJq \widetilde{E}-(A\widetilde{E})^{\ast}\tJq A\widetilde{E}
	= \jqq-B^{\ast}\tJq B.
\end{align*}

Zu (c2): Dies folgt wegen (a) und (c1) aus \thref{spbsp1}, \thref{spbsp5}, Teil (a) von \thref{spdef2} und Teil (a) von \thref{spdef7}.

Zu (c3): Dies folgt wegen (a) und (c1) aus \thref{spbsp1}, \thref{spbsp5}, Teil (b) von \thref{spdef2} und Teil (b) von \thref{spdef7}. \bwend

\section[Einige Aussagen über ganze Funktionen aus \textit{J}"=Potapov"=Klassen bezüglich Halbebenen]{Einige Aussagen über ganze Funktionen aus \\ \textit{J}-Potapov-Klassen bezüglich Halbebenen}  \label{chappf}

In diesem Anhang betrachten wir ganze Matrixfunktionen, welche in spezieller Weise mit Signaturmatrizen verknüpft sind. Derartige Matrixfunktionen wurden systematisch von V.\,P. Potapov untersucht, welcher in \cite{Po} ein tiefliegendes Faktorisierungstheorem für Funktionen dieser Klasse formulierte.

Wir kommen nun zur Definition der Potapov-Funktion bezüglich der Halbebene $\Pp$ und betrachten anschließend einige grundlegende Eigenschaften.

\begin{defi}	\thlabel{spdef3} Seien $J$ eine \textit{p}$\times$\textit{p}-Signaturmatrix und $W$ eine in $\C$ holomorphe \textit{p}$\times$\textit{p}-Matrixfunktion.
\begin{itemize}
  \item [\rm{(a)}] Es heißt $W$ \textbf{J-Potapov-Funktion} bezüglich $\Pp$, falls $W(z)$ für alle $z \in \Pp$ eine $J$-kontraktive Matrix ist. Mit $\bJPp$ bezeichnen wir die Menge aller $J$-Potapov-Funktionen bezüglich $\Pp$.
  \item [\rm{(b)}] Sei $W \in \bJPp$. Dann heißt $W$ \textbf{J-innere Funktion} aus $\bJPp$, falls $W(x)$ für alle $x \in \R$ eine $J$-unitäre Matrix ist. Mit $\tbJPp$ bezeichnen wir die Menge aller $J$-inneren Funktionen aus $\bJPp$.
\end{itemize}
\end{defi}

\begin{lemma}	\thlabel{splm1} Seien $J$ eine \textit{p}$\times$\textit{p}-Signaturmatrix und $W \in \tbJPp$. \dgfa
\begin{itemize}
  \item [\rm{(a)}] Sei $z \in \C$. Dann ist $W(z)$ regulär und es gilt $W^{-1}(z) = JW^{\ast}(\za)J$.
  \item [\rm{(b)}] Seien $z, \omega \in \C$. Dann gilt 
  \begin{align*}
    J-W^{-\ast}(z)JW^{-1}(\omega) = J\eklam{J-W(\za)JW^{\ast}(\omegaa)}J.
  \end{align*}
  \item [\rm{(c)}] Sei $z \in \Pm$. Dann ist $W(z)$ eine $J$-expansive Matrix.
  \item [\rm{(d)}] Sei $z \in \C\setminus\R$. Dann gilt
  \begin{align*}
    \frac{W^{\ast}(z)JW(z)-J}{i(z-\za)} \geq \Op.
  \end{align*}
\end{itemize}
\end{lemma}

\bwanf Zu (a): Siehe Teil (a) von \cite[Lemma 5.1]{C06}.

Zu (b): Wegen \thref{spdef1} und (a) gilt
\begin{align*}
	J\eklam{J-W(\za)JW^{\ast}(\omegaa)}J
	&= J-JW(\za)JW^{\ast}(\omegaa)J \\
	&= J-\eklam{JW(\za)J}J\eklam{JW^{\ast}(\omegaa)J} \\
	&= J-W^{-\ast}(z)JW^{-1}(\omega).
\end{align*}

Zu (c): Dies folgt unter Beachtung von Teil (a) von \thref{spdef2} aus Teil (b) von \cite[Lemma 5.1]{C06}.

Zu (d): Siehe Teil (c) von \cite[Lemma 5.1]{C06}. \bwend

\thref{splm1} erlaubt uns folgende Schlussfolgerung speziell für Matrixpolynome.

\begin{folg}	\thlabel{spfo1}
	Seien $J$ eine \textit{p}$\times$\textit{p}-Signaturmatrix und $W \in \tbJPp$ ein \textit{p}$\times$\textit{p}"=Matrixpolynom. Dann ist $\det W$ eine konstante Funktion auf $\C$ und deren Wert von Null verschieden.
\end{folg}

\bwanf Da $W$ ein \textit{p}$\times$\textit{p}-Matrixpolynom ist, ist $\det W$ ein skalares Matrixpolynom. Wegen Teil (a) von \thref{splm1} verschwindet $\det W$ nirgends in $\C$. Der Fundamentalsatz der Algebra liefert dann, dass $\det W$ eine konstante Funktion auf $\C$ ist. \bwend

Folgendes Resultat veranschaulicht die Multiplikativität der in \thref{spdef3} eingeführten Klassen von ganzen Matrixfunktionen.

\begin{bem}	\thlabel{spbm6}
	Seien $J$ eine \textit{p}$\times$\textit{p}-Signaturmatrix sowie $W_1$ und $W_2$ Funktionen aus $\bJPp$ bzw. $\tbJPp$. Dann ist auch $W_1W_2$ eine Funktion aus $\bJPp$ bzw. $\tbJPp$.
\end{bem}

\bwanf Dies folgt sogleich aus Teil (d) von \thref{spbm1} und \thref{spdef3}. \bwend

Wir kommen nun zur Definition der Potapov-Funktion bezüglich der Halbebene $\Lam$ für beliebige reelle $\alpha$.

\begin{defi}	\thlabel{spdef8} Seien $\alpha\in\R$, $J$ eine \textit{p}$\times$\textit{p}-Signaturmatrix und $W$ eine in $\C$ holomorphe \textit{p}$\times$\textit{p}-Matrixfunktion.
\begin{itemize}
  \item [\rm{(a)}] Es heißt $W$ \textbf{J-Potapov-Funktion} bezüglich $\Lam$, falls $W(z)$ für alle $z \in \Lam$ eine $J$-kontraktive Matrix ist. Mit $\bJLam$ bezeichnen wir die Menge aller $J$-Potapov-Funktionen bezüglich $\Lam$.
  \item [\rm{(b)}] Sei $W \in \bJLam$. Dann heißt $W$ \textbf{J-innere Funktion} aus $\bJLam$, falls $W(z)$ für alle $z \in \C$ mit $\re z=\alpha$ eine $J$-unitäre Matrix ist. Mit $\tbJLam$ bezeichnen wir die Menge aller $J$-inneren Funktionen aus $\bJLam$.
\end{itemize}
\end{defi}

\begin{bem}	\thlabel{spbm8}
	Seien $\alpha\in\R$, $J$ eine \textit{p}$\times$\textit{p}-Signaturmatrix sowie $W_1$ und $W_2$ Funktionen aus $\bJLam$ bzw. $\tbJLam$. Dann ist auch $W_1W_2$ eine Funktion aus $\bJLam$ bzw. $\tbJLam$.
\end{bem}

\bwanf Dies folgt sogleich aus Teil (d) von \thref{spbm1} und \thref{spdef8}. \bwend

Wir wollen nun den Begriff der Potapov-Funktionen bezüglich $\Pp$ aus \thref{spdef3} erweitern, indem wir statt nur einer Signaturmatrix nun zwei verwenden werden. Diese Überlegung basiert auf die in \thref{spdef7} eingeführten Begriffsbildungen.

\begin{defi}	\thlabel{spdef9}
	Seien $J^{(1)}$ und $J^{(2)}$ jeweils eine \textit{p}$\times$\textit{p}-Signaturmatrix sowie $W$ eine in $\C$ holomorphe \textit{p}$\times$\textit{p}-Matrixfunktion.
	\begin{itemize}
 		\item [\rm{(a)}] Es heißt $W$ \textbf{J$\mathbf{^{(2)}}$-J$\mathbf{^{(1)}}$-Potapov-Funktion} bezüglich $\Pp$, falls $W(z)$ für alle \linebreak $z \in \Pp$ eine $J^{(2)}$-$J^{(1)}$-kontraktive Matrix ist. Mit $\bJJPp$ bezeichnen wir die Menge aller $J^{(2)}$-$J^{(1)}$-Potapov-Funktionen bezüglich $\Pp$.
  		\item [\rm{(b)}] Sei $W \in \bJJPp$. Dann heißt $W$ \textbf{J$\mathbf{^{(2)}}$-J$\mathbf{^{(1)}}$-innere Funktion} aus \linebreak $\bJJPp$, falls $W(x)$ für alle $x \in \R$ eine $J^{(2)}$-$J^{(1)}$-unitäre Matrix ist. Mit $\tbJJPp$ bezeichnen wir die Menge aller $J^{(2)}$-$J^{(1)}$-inneren Funktionen aus $\bJJPp$.
	\end{itemize}
\end{defi} 

Nimmt man in \thref{spdef9} die Setzung $J^{(1)} = J^{(2)}$ vor, so erhält man die Begriffsbildungen aus \thref{spdef3}. Wir betrachten nun ein spezielles Resultat für die so eben eingeführten Begriffsbildungen, das die in Kapitel \ref{chapar} vorliegende Situation umfasst.

\begin{bem}	\thlabel{spbm9}
	Sei $W_1$ eine in $\C$ holomorphe $2$\textit{q}$\times2$\textit{q}-Matrixfunktion.
	\dgfa
	\begin{itemize}
		\item [\rm{(a)}] Seien
		\begin{align*}
			E:=\frac{1}{\sqrt{2}}\begin{pmatrix} -i\Iq & i\Iq \\ \Iq & \Iq \end{pmatrix}
		\end{align*}
		und $W_2:=W_1E$. \dgfa
		\begin{itemize}
			\item [\rm{(a1)}] Es ist $W_1\in\btJqPp$ genau dann, wenn $W_2\in\btJJPp$ erfüllt ist.
			\item [\rm{(a2)}] Es ist $W_1\in\tbtJqPp$ genau dann, wenn $W_2\in\tbtJJPp$ erfüllt ist.
		\end{itemize}
		\item [\rm{(b)}] Seien
		\begin{align*}
			\widetilde{E}:=\frac{1}{\sqrt{2}}\begin{pmatrix} -i\Iq & -i\Iq \\ -\Iq & \Iq \end{pmatrix}
		\end{align*}
		und $W_2:=W_1\widetilde{E}$. \dgfa
		\begin{itemize}
			\item [\rm{(b1)}] Es ist $W_1\in\btJqPp$ genau dann, wenn $W_2\in\btJJPp$ erfüllt ist.
			\item [\rm{(b2)}] Es ist $W_1\in\tbtJqPp$ genau dann, wenn $W_2\in\tbtJJPp$ erfüllt ist.
		\end{itemize}
	\end{itemize}
\end{bem}

\bwanf Zu (a): Dies folgt unter Beachtung von \thref{spdef3} und \thref{spdef9} aus \thref{spbsp1}, \thref{spbsp5} sowie den Teilen (c2) und (c3) von \thref{spbm7}. 

Zu (b): Dies folgt unter Beachtung von \thref{spdef3} und \thref{spdef9} aus \thref{spbsp1}, \thref{spbsp5} sowie den Teilen (c2) und (c3) von \thref{spbm10}. \bwend

\section[Einige Aussagen über Stieltjes"=Paare von meromorphen Matrixfunktionen]{Einige Aussagen über Stieltjes"=Paare von \\ meromorphen Matrixfunktionen}  \label{chapsp}

Im Mittelpunkt dieses Anhangs steht die Diskussion spezieller Klassen geordneter Paare von meromorphen Matrixfunktionen, welche in dieser Arbeit die Rolle der Parametermenge bei der Parametrisierung der Lösungsmenge eines nichtdegenerierten matriziellen Momentenproblems vom $\alpha$-Stieltjes-Typ spielen (vergleiche Kapitel \ref{chapdr}).

Zunächst führen wir den Begriff des \textit{q}$\times$\textit{q}-Stieltjes-Paares in $\C\setminus[\alpha,\infty)$ ein. Eine grundlegende Beschreibung dieses Themas findet man z.\,B. in \cite[Kapitel 10]{Maka}.

\begin{defi}	\thlabel{spdef4} 
	Sei $\alpha \in \R$. 
	\begin{itemize}
  		\item [\rm{(a)}] Seien $\phi$ und $\psi$ in $\C\setminus[\alpha,\infty)$ meromorphe \textit{q}$\times$\textit{q}-Matrixfunktionen.
  		Dann heißt $\phipsi$ \textbf{q$\times$q-Stieltjes-Paar in $\C\setminus[\alpha,\infty)$}, falls eine diskrete Teilmenge $\D$ von $\C\setminus[\alpha,\infty)$ existiert mit
 		\begin{itemize}
    		\item [\rm{(i)}] $\phi$ und $\psi$ sind in $\C\setminus\rklam{[\alpha,\infty)\cup\D}$ holomorph.
    		\item [\rm{(ii)}] Es gilt $\rank \phipsiz = q$ für alle $z \in \C\setminus\rklam{[\alpha,\infty)\cup\D}$.
    		\item [\rm{(iii)}] Es gelten
    		\begin{align*}
     			\binom{(z-\alpha)\phi(z)}{\psi(z)}^{\ast}\bbrklam{\frac{-\tJq}{2\im z}}\binom{(z-\alpha)\phi(z)}{\psi(z)} \geq \Oq
    		\end{align*}
    		und
    		\begin{align*}
      			\phipsiz^{\ast}\bbrklam{\frac{-\tJq}{2\im z}}\phipsiz \geq \Oq
    		\end{align*}
    		für alle $z \in \C\setminus\rklam{\R\cup\D}$.
  		\end{itemize}
  		Mit $\PtJqCa$ bezeichnen wir die Menge aller q$\times$q-Stieltjes-Paare in \linebreak $\C\setminus[\alpha,\infty)$. Gelten {\rm (i)}-{\rm (iii)} sogar für $\D=\emptyset$, so schreiben wir $\phipsi\in\dPtJqCa$.
  		\item [\rm{(b)}] Sei $\phipsi\in\PtJqCa$. Dann heißt das Paar $\phipsi$ \textbf{eigentlich}, falls $\det\psi$ nicht die Nullfunktion ist.  		
	\end{itemize}
\end{defi}  

Im Folgenden wollen wir nun eine Äquivalenzrelation auf $\PtJqCa$ einführen. Hierfür benötigen wir zunächst noch folgendes Lemma (vergleiche \cite[Lemma 1.11]{Maka2}).

\begin{lemma}	\thlabel{splm6}
	Seien $\alpha\in\R$ und $\phipsi\in\PtJqCa$. Weiterhin sei $g$ eine in \linebreak $\C\setminus[\alpha,\infty)$ meromorphe \textit{q}$\times$\textit{q}-Matrixfunktion derart, dass $\det g$ nicht die Nullfunktion ist. Dann gilt $\binom{\phi g}{\psi g}\in\PtJqCa$.
\end{lemma}

\bwanf Sei $\D$ die gemäß Teil (a) von \thref{spdef4} existierende diskrete Teilmenge von $\C\setminus[\alpha,\infty)$. Da $\det g$ nicht die Nullfunktion ist, existiert unter Beachtung des Identitätssatzes für meromorphe Funktionen (vergleiche z.\,B. im skalaren Fall \cite[Satz 10.3.2]{Funk}; im matriziellen Fall betrachtet man die einzelnen Einträge der Matrixfunktion) eine diskrete Teilmenge $\D_1$ von $\C\setminus[\alpha,\infty)$ mit $\det g(z)\neq0$ für alle $z \in \C\setminus\rklam{[\alpha,\infty)\cup\D_1}$. Sei $\D_2 := \D \cup \D_1$. Dann kann man sich leicht davon überzeugen, dass die Bedingungen (i)-(iv) von Teil (a) von \thref{spdef4} für das Paar $\binom{\phi g}{\psi g}$ und der diskreten Teilmenge $\D_2$ erfüllt sind. \bwend

\thref{splm6} bringt uns nun auf folgende Definition.

\begin{defi}	\thlabel{spdef4b} 
	Seien $\alpha \in \R$ und $\binom{\phi_1}{\psi_1}, \binom{\phi_2}{\psi_2} \in \PtJqCa$.
  	Dann heißen $\binom{\phi_1}{\psi_1}$ und $\binom{\phi_2}{\psi_2}$ äquivalent, falls eine in $\C\setminus[\alpha,\infty)$ meromorphe \textit{q}$\times$\textit{q}-Matrixfunktion $g$ und eine diskrete Teilmenge $\D$ von $\C\setminus[\alpha,\infty)$ existieren mit
  	\begin{itemize}
    	\item [\rm{(i)}] $\phi_1, \psi_1, \phi_2, \psi_2$ und $g$ sind in $\C\setminus\rklam{[\alpha,\infty)\cup\D}$ holomorph.
    	\item [\rm{(ii)}] Es gelten $\det g(z) \neq 0$ und
    	\begin{align*}
      		\binom{\phi_2(z)}{\psi_2(z)} = \binom{\phi_1(z)}{\psi_1(z)}g(z)
    	\end{align*}
    	für alle $z \in \C\setminus\rklam{[\alpha,\infty)\cup\D}$.
  	\end{itemize}
\end{defi}

Man kann sich leicht davon überzeugen, dass folgende Bemerkung richtig ist (vergleiche auch \cite[Bemerkung 10.12]{Maka}).

\begin{bem}	\thlabel{spbm12}
	Der gemäß \thref{spdef4b} definierte Äquivalenzbegriff ist eine Äquivalenzrelation auf $\PtJqCa$. Mit $\aphipsi$ bezeichnen wir die durch $\phipsi\in\PtJqCa$ erzeugte Äquivalenzklasse.
\end{bem}

Wir zeigen nun, dass die in Teil (a) von \thref{spdef4} eingeführte Klasse von \textit{q}$\times$\textit{q}-Stieltjes-Paaren in $\C\setminus[\alpha,\infty)$ als eine projektive Erweiterung der Funktionenklasse $\Sqa$ (vergleiche \thref{asmbz3}) aufgefasst werden kann (vergleiche Teil (a) mit \cite[Lemma 10.16]{Maka}).

\begin{bem}	\thlabel{spbm11}
	Sei $\alpha\in\R$. \dgfa
	\begin{itemize}
		\item [\rm{(a)}] Sei $S\in\Sqa$. Weiterhin bezeichne ${\cal I}$ die in $\C\setminus[\alpha,\infty)$ konstante Matrixfunktion mit dem Wert $\Iq$. Dann ist $\binom{S}{\cal I}$ ein eigentliches Paar aus $\dPtJqCa$.
		\item [\rm{(b)}] Seien $S_1,S_2\in\Sqa$ so beschaffen, dass die Paare $\binom{S_1}{\cal I}$ und $\binom{S_2}{\cal I}$ äquivalent sind. Dann gilt $S_1 = S_2$.
	\end{itemize}
\end{bem}

\bwanf Zu (a): Wegen \cite[Proposition 4.4]{SF} gelten
\begin{itemize}
	\item [\rm{(i)}] Es ist $S$ in $\C\setminus[\alpha,\infty)$ holomorph.
	\item [\rm{(ii)}] Es gilt $\im S(z) \geq \Oq$ für alle $z\in\Pp$.
	\item [\rm{(iii)}] Es gilt $-\im S(z) \geq \Oq$ für alle $z\in\Pm$.
\end{itemize}
Wegen (ii), (iii) und \thref{spbsp1} gilt dann
\begin{align}	\label{spbm11bw1}
	\binom{S(z)}{\Iq}^{\ast}\bbrklam{\frac{-\tJq}{2\im z}}\binom{S(z)}{\Iq} = \frac{1}{\im z}\im S(z) \geq \Oq
\end{align}
für alle $z\in\C\setminus\R$. Wegen \thref{spbsp1} und \cite[Lemma 4.2]{SF} gilt weiterhin
\begin{align}	\label{spbm11bw2}
	\binom{(z-\alpha)S(z)}{\Iq}^{\ast}\bbrklam{\frac{-\tJq}{2\im z}}\binom{(z-\alpha)S(z)}{\Iq} = \frac{1}{\im z}\im\eklam{(z-\alpha)S(z)} \geq \Oq
\end{align}
für alle $z\in\C\setminus\R$. Offensichtlich ist ${\cal I}$ in $\C\setminus[\alpha,\infty)$ holomorph und $\rank{\cal I}(z)=q$ für alle $z\in\C\setminus[\alpha,\infty)$. Hieraus folgt wegen (i), \fref{spbm11bw1}, \fref{spbm11bw2} und \thref{spdef4} dann, dass $\binom{S}{\cal I}$ ein eigentliches Paar aus $\dPtJqCa$ ist. 

Zu (b): Wegen $S_1,S_2\in\Sqa$ (vergleiche \thref{asmbz3}) sind $S_1$ und $S_2$ in \linebreak $\C\setminus[\alpha,\infty)$ holomorph. Wegen (a) und \thref{spdef4b} existiert eine in $\C\setminus[\alpha,\infty)$ meromorphe \textit{q}$\times$\textit{q}-Matrixfunktion $g$ und eine diskrete Teilmenge $\D$ von $\C\setminus[\alpha,\infty)$ mit
\begin{itemize}
   	\item [\rm{(iv)}] $g$ ist in $\C\setminus\rklam{[\alpha,\infty)\cup\D}$ holomorph.
   	\item [\rm{(v)}] Es gelten $\det g(z) \neq 0$ und
    \begin{align*}
    	\binom{S_2(z)}{\Iq} = \binom{S_1(z)}{\Iq}g(z)
    \end{align*}
    für alle $z \in \C\setminus\rklam{[\alpha,\infty)\cup\D}$.
\end{itemize}
Aus (v) folgt sogleich $g(z) = \Iq$ für alle $z \in \C\setminus\rklam{[\alpha,\infty)\cup\D}$ und somit $S_2(z)=S_1(z)$ für alle $z \in \C\setminus\rklam{[\alpha,\infty)\cup\D}$. Wegen des Identitätssatzes für meromorphe Funktionen (vergleiche z.\,B. im skalaren Fall \cite[Satz 10.3.2]{Funk}; im matriziellen Fall betrachtet man die einzelnen Einträge der Matrixfunktion) folgt nun $S_2=S_1$. \bwend

Wir betrachten nun den umgekehrten Fall von \thref{spbm11}, indem wir aus einem gegebenen \textit{q}$\times$\textit{q}-Stieltjes-Paar in $\C\setminus[\alpha,\infty)$ eine zu $\Sqa$ gehörige Matrixfunktion konstruieren (vergleiche \cite[Lemma 10.18]{Maka}).

\begin{bem}	\thlabel{spbm13}
	Seien $\alpha\in\R$ und $\phipsi$ ein eigentliches Paar aus $\PtJqCa$. \dgfa
	\begin{itemize}
		\item [\rm{(a)}] Bezeichne ${\cal I}$ die in $\C\setminus[\alpha,\infty)$ konstante Matrixfunktion mit dem Wert $\Iq$. Dann ist $\binom{\phi\psi^{-1}}{\cal I}$ ein eigentliches Paar aus $\PtJqCa$, welches zu $\phipsi$ äquivalent ist.
		\item [\rm{(b)}] Es gilt $\phi\psi^{-1}\in\Sqa$.
	\end{itemize}
\end{bem}

\bwanf Zu (a): Unter Beachtung von $\binom{\phi\psi^{-1}}{\cal I}=\phipsi \psi^{-1}$ gilt wegen \thref{splm6}, Teil (b) von \thref{spdef4} sowie $\rank{\cal I}(z)=q$ für alle $z\in\C\setminus[\alpha,\infty)$, dass $\binom{\phi\psi^{-1}}{\cal I}$ ein eigentliches Paar aus $\PtJqCa$ ist, und wegen \thref{spdef4b} weiterhin, dass $\binom{\phi\psi^{-1}}{\cal I}$ äquivalent ist zu $\phipsi$.

Zu (b): Wegen (a) und Teil (a) von \thref{spdef4} existiert eine diskrete Teilmenge $\D$ von $\C\setminus[\alpha,\infty)$ mit
\begin{itemize}
	\item [\rm{(i)}] $\phi\psi^{-1}$ ist in $\C\setminus\rklam{[\alpha,\infty)\cup\D}$ holomorph.
	\item [\rm{(ii)}] Es gelten
	\begin{align*}
    	\binom{(z-\alpha)\phi(z)\psi^{-1}(z)}{\Iq}^{\ast}\bbrklam{\frac{-\tJq}{2\im z}}\binom{(z-\alpha)\phi(z)\psi^{-1}(z)}{\Iq} \geq \Oq
    \end{align*}
    und
    \begin{align*}
    	\binom{\phi(z)\psi^{-1}(z)}{\Iq}^{\ast}\bbrklam{\frac{-\tJq}{2\im z}}\binom{\phi(z)\psi^{-1}(z)}{\Iq} \geq \Oq
    \end{align*}
    für alle $z \in \C\setminus\rklam{\R\cup\D}$.
\end{itemize}
Aus \thref{spbsp1} und (ii) folgt dann
\begin{align*}
	\frac{1}{\im z}\im\eklam{\phi(z)\psi^{-1}(z)} \geq \Oq
\end{align*}
für alle $z \in \C\setminus\rklam{\R\cup\D}$. Hieraus folgen dann
\begin{align}	\label{spbm13bw2}
	\im\eklam{\phi(z)\psi^{-1}(z)} \geq \Oq
\end{align}
für alle $z\in\Pp\setminus\D$ und
\begin{align}	\label{spbm13bw3}
	-\im\eklam{\phi(z)\psi^{-1}(z)} \geq \Oq
\end{align}
für alle $z\in\Pm\setminus\D$. Wegen \thref{spbsp1} und (ii) gilt weiterhin
\begin{align}	\label{spbm13bw1}
	&\ \frac{\alpha-\re z}{\im z}\im\eklam{\phi(z)\psi^{-1}(z)} \notag \\
	&= (\alpha-\re z)\binom{\phi(z)\psi^{-1}(z)}{\Iq}^{\ast}\bbrklam{\frac{-\tJq}{2\im z}}\binom{\phi(z)\psi^{-1}(z)}{\Iq}
	\geq \Oq
\end{align}
für alle $z\in \Lam\setminus(\R\cup\D)$. Unter Beachtung von
\begin{align*}
	\im\eklam{zA} = \re z \im A + \im z \re A
\end{align*}
für beliebige $z\in\C$ und $A\in\Cqq$ (vergleiche z.\,B. \cite[Bemerkung A.23]{Sch}) gilt wegen \fref{spbm13bw1}, \thref{spbsp1} und (ii) dann
\begin{align*}
	&\ \re\eklam{\phi(z)\psi^{-1}(z)}
	= \frac{1}{\im z}\im\eklam{z\phi(z)\psi^{-1}(z)} - \frac{\re z}{\im z}\im\eklam{\phi(z)\psi^{-1}(z)} \notag \\
	&= \frac{1}{\im z}\im\eklam{z\phi(z)\psi^{-1}(z)} + \frac{\alpha-\re z}{\im z}\im\eklam{\phi(z)\psi^{-1}(z)} - \frac{\alpha}{\im z}\im\eklam{\phi(z)\psi^{-1}(z)} \notag \\
	&\geq \frac{1}{\im z}\im\eklam{z\phi(z)\psi^{-1}(z)} - \frac{1}{\im z}\im\eklam{\alpha\phi(z)\psi^{-1}(z)} \notag \\
	&= \frac{1}{\im z}\im\eklam{(z-\alpha)\phi(z)\psi^{-1}(z)} \notag \\
	&= \binom{(z-\alpha)\phi(z)\psi^{-1}(z)}{\Iq}^{\ast}\bbrklam{\frac{-\tJq}{2\im z}}\binom{(z-\alpha)\phi(z)\psi^{-1}(z)}{\Iq} 
	\geq \Oq
\end{align*}
für alle $z\in \Lam\setminus(\R\cup\D)$. Hieraus folgt wegen der Stetigkeit der Funktion $\phi\psi^{-1}$ in $\C\setminus([\alpha,\infty)\cup\D)$ nun
\begin{align}	\label{spbm13bw4}
	\re\eklam{\phi(z)\psi^{-1}(z)} \geq \Oq
\end{align}
für alle $z\in \Lam\setminus\D$. Wegen \fref{spbm13bw2}, \fref{spbm13bw3}, \fref{spbm13bw4} sowie unter Beachtung von
\begin{align*}
	\re\eklam{iA} = \frac{1}{2}\brklam{iA+\eklam{iA}^{\ast}}=\frac{1}{2i}\brklam{{-A-\eklam{-A}^{\ast}}} = \im\eklam{-A}
\end{align*}
für alle $A\in\Cqq$ und \cite[Lemma 2.19]{DFK} oder \cite[Lemma 4.23]{Sch} ist jedes $z\in\D$ eine hebbare Singularität von $\phi\psi^{-1}$ und $\phi\psi^{-1}$ lässt sich zu einer auf $\C\setminus[\alpha,\infty)$ holomorphen Matrixfunktion $F$ fortsetzen mit
\begin{itemize}
   	\item [\rm{(iii)}] Es gilt $\im F(z) \geq \Oq$ für alle $z \in \Pp$.
   	\item [\rm{(iv)}] Es gilt $-\im F(z) \geq \Oq$ für alle $z \in \Pm$.
    \item [\rm{(v)}] Es gilt $\re F(z) \geq \Oq$ für alle $z\in\Lam$.
\end{itemize}
Hieraus folgt wegen \cite[Proposition 4.4]{SF} dann $F\in\Sqa$. \bwend

Mithilfe von \thref{spbm11} und \thref{spbm13} erkennen wir, dass eine Bijektion zwischen $\Sqa$ und der Menge der Äquivalenzklassen von eigentlichen \textit{q}$\times$\textit{q}-Stieltjes-Paaren in $\C\setminus[\alpha,\infty)$ besteht.

Wir geben nun Beispiele zweier bemerkenswerter \textit{q}$\times$\textit{q}-Stieltjes-Paare in $\C\setminus[\alpha,\infty)$ an, welche in unseren Überlegungen für Abschnitt \ref{EE} eine wesentliche Rolle spielen (vergleiche \cite[Example 1.17]{CR1} für den Fall $\alpha=0$).

\begin{beispiel}	\thlabel{spbsp2}
	Sei $\alpha\in\R$. Bezeichne ${\cal I}$ bzw. ${\cal O}$ die in $\C\setminus[\alpha,\infty)$ konstante Matrixfunktion mit dem Wert $\Iq$ bzw. $\Oq$. Dann gilt $\binom{\cal I}{\cal O}, \binom{\cal O}{\cal I} \in \dPtJqCa$. Insbesondere ist das Paar $\binom{\cal O}{\cal I}$ eigentlich.
\end{beispiel}

\bwanf Offensichtlich sind ${\cal I}$ und ${\cal O}$ in $\C\setminus[\alpha,\infty)$ holomorph und es gilt $\rank{\cal I}(z) = q$ für alle $z\in\C\setminus[\alpha,\infty)$. Hieraus folgt unter Beachtung von \thref{spbsp1} und \thref{spdef4} dann die Behauptung. \bwend	

Wir führen nun den Begriff des \textit{q}$\times$\textit{q}-Stieltjes-Paares in $\C\setminus(-\infty,\alpha]$ ein.

\begin{defi}	\thlabel{spdef5} 
Sei $\alpha \in \R$. 
	\begin{itemize}
  	\item [\rm{(a)}] Seien $\phi$ und $\psi$ in $\C\setminus(-\infty,\alpha]$ meromorphe \textit{q}$\times$\textit{q}-Matrixfunktionen.
  	Dann heißt $\phipsi$ \textbf{q$\times$q-Stieltjes-Paar in $\C\setminus(-\infty,\alpha]$}, falls eine diskrete Teilmenge $\D$ von $\C\setminus(-\infty,\alpha]$ existiert mit
  \begin{itemize}
    \item [\rm{(i)}] $\phi$ und $\psi$ sind in $\C\setminus\rklam{(-\infty,\alpha]\cup\D}$ holomorph.
    \item [\rm{(ii)}] Es gilt $\rank \phipsiz = q$ für alle $z \in \C\setminus\rklam{(-\infty,\alpha]\cup\D}$.
    \item [\rm{(iii)}] Es gelten
    \begin{align*}
      \binom{(\alpha-z)\phi(z)}{\psi(z)}^{\ast}\bbrklam{\frac{-\tJq}{2\im z}}\binom{(\alpha-z)\phi(z)}{\psi(z)} \geq \Oq
    \end{align*}
    und
    \begin{align*}
      \phipsiz^{\ast}\bbrklam{\frac{-\tJq}{2\im z}}\phipsiz \geq \Oq
    \end{align*}
    für alle $z \in \C\setminus\rklam{\R\cup\D}$.
  \end{itemize}
  Mit $\PtJqCma$ bezeichnen wir die Menge aller q$\times$q-Stieltjes-Paare in $\C\setminus(-\infty,\alpha]$. Gelten {\rm (i)}-{\rm (iii)} sogar für $\D=\emptyset$, so schreiben wir $\phipsi\in\dPtJqCma$.  
  \item [\rm{(b)}] Sei $\phipsi\in\PtJqCma$. Dann heißt das Paar $\phipsi$ \textbf{eigentlich}, falls $\det\psi$ nicht die Nullfunktion ist.  
\end{itemize}
\end{defi}

Folgende Bemerkung liefert uns einen Zusammenhang zwischen den in \thref{spdef4} und \thref{spdef5} eingeführten Klassen.

\begin{bem}	\thlabel{spbm3}
	Seien $\alpha\in\R$, $j\in\gklam{0,1}$ sowie $\phi, \psi:\C\setminus[\alpha,\infty)\rightarrow\Cqq$ und $\widecheck{\phi}, \widecheck{\psi}:\C\setminus(-\infty,-\alpha]\rightarrow\Cqq$ definiert gemäß $\widecheck{\phi}(z) := (-1)^j\phi(-z)$ bzw. $\widecheck{\psi}(z) := (-1)^{j+1}\psi(-z)$. Weiterhin seien $\phi$ und $\psi$ in $\C\setminus[\alpha,\infty)$ oder $\widecheck{\phi}$ und $\widecheck{\psi}$ in $\C\setminus(-\infty,-\alpha]$ meromorphe \textit{q}$\times$\textit{q}-Matrixfunktionen. Dann sind $\widecheck{\phi}$ und $ \widecheck{\psi}$ in $\C\setminus(-\infty,-\alpha]$ bzw. $\phi$ und $\psi$ in $\C\setminus[\alpha,\infty)$ meromorphe \textit{q}$\times$\textit{q}-Matrixfunktionen und es gelten folgende Aussagen:
	\begin{itemize}
		\item [\rm{(a)}] \esfaa
		\begin{itemize}
			\item [\rm{(i)}] Es gilt $\phipsi\in\PtJqCa$.
			\item [\rm{(ii)}] Es gilt $\bbinom{\widecheck{\phi}}{\widecheck{\psi}}\in{\cal P}^{(q,q)}_{-\tJq,\geq}(\C$,""$(-\infty,-\alpha])$.
		\end{itemize}
		\item [\rm{(b)}] \esfaa
		\begin{itemize}
			\item [\rm{(iii)}] Es gilt $\phipsi\in\dPtJqCa$.
			\item [\rm{(iv)}] Es gilt $\bbinom{\widecheck{\phi}}{\widecheck{\psi}}\in\widehat{\cal P}^{(q,q)}_{-\tJq,\geq}(\C$,""$(-\infty,-\alpha])$.
		\end{itemize}
		\item [\rm{(c)}] Sei {\rm (i)} erfüllt. Dann ist $\phipsi$ ein eigentliches Paar genau dann, wenn $\bbinom{\widecheck{\phi}}{\widecheck{\psi}}$ ein eigentliches Paar ist.
	\end{itemize}
\end{bem}

\bwanf Zu (a): Seien $\D$ eine diskrete Teilmenge von $\C\setminus[\alpha,\infty)$ und $\widecheck{\D} := \gklam{-\omega \;|\; \omega \in \D}$. Offensichtlich gilt dann
\begin{itemize}
	\item [\rm{(I)}] Es sind $\phi$ und $\psi$ genau dann in $\C\setminus\rklam{[\alpha,\infty)\cup\D}$ holomorph, wenn $\widecheck{\phi}$ und $\widecheck{\psi}$ in $\C\setminus\brklam{(-\infty,-\alpha]\cup\widecheck{\D}}$ holomorph sind.
\end{itemize}
Weiterhin gilt
\begin{align}	\label{spbm3bw1}
	\rank \phipsiz = \rank \begin{pmatrix}\widecheck{\phi}(-z)\\\widecheck{\psi}(-z)\end{pmatrix}
\end{align}
für alle $z\in\C\setminus\rklam{[\alpha,\infty)\cup\D}$. Wegen \thref{spbsp1} gelten
\begin{align*}
	&\ \binom{(z-\alpha)\phi(z)}{\psi(z)}^{\ast}\bbrklam{\frac{-\tJq}{2\im z}}\binom{(z-\alpha)\phi(z)}{\psi(z)} \\
	&= \frac{1}{\im z}\im\eklam{(z-\alpha)\psi^{\ast}(z)\phi(z)} \\
	&= \frac{1}{\im\eklam{-z}}\im\beklam{(-\alpha-(-z))\widecheck{\psi}^{\ast}(-z)\widecheck{\phi}(-z)} \\
	&= \begin{pmatrix}(-\alpha-(-z))\widecheck{\phi}(-z)\\\widecheck{\psi}(-z)\end{pmatrix}^{\ast}\bbrklam{\frac{-\tJq}{2\im\eklam{-z}}}\begin{pmatrix}(-\alpha-(-z))\widecheck{\phi}(-z)\\\widecheck{\psi}(-z)\end{pmatrix}
\end{align*}
und
\begin{align*}
	&\ \binom{\phi(z)}{\psi(z)}^{\ast}\bbrklam{\frac{-\tJq}{2\im z}}\binom{\phi(z)}{\psi(z)}
	= \frac{1}{\im z}\im\eklam{\psi^{\ast}(z)\phi(z)} \\
	&= \frac{1}{\im\eklam{-z}}\im\beklam{\widecheck{\psi}^{\ast}(-z)\widecheck{\phi}(-z)}
	= \begin{pmatrix}\widecheck{\phi}(-z)\\\widecheck{\psi}(-z)\end{pmatrix}^{\ast}\bbrklam{\frac{-\tJq}{2\im\eklam{-z}}}\begin{pmatrix}\widecheck{\phi}(-z)\\\widecheck{\psi}(-z)\end{pmatrix}
\end{align*}
für alle $z\in\C\setminus\rklam{\R\cup\D}$. Unter Beachtung von Teil (a) von \thref{spdef4} und Teil (a) von \thref{spdef5} folgt hieraus wegen (I) und \fref{spbm3bw1} dann die Äquivalenz von (i) und (ii).

Zu (b): Dies folgt unter Beachtung, dass im Fall $\D=\emptyset$ dann $\widecheck{\D}=\emptyset$ gilt, aus dem Beweis von (a) sowie Teil (a) von \thref{spdef4} und Teil (a) von \thref{spdef5}.

Zu (c): Offensichtlich ist $\det\psi$ nicht die Nullfunktion genau dann, wenn $\det\widecheck{\psi}$ nicht die Nullfunktion ist. Hieraus folgt wegen Teil (b) von \thref{spdef4} und Teil (b) von \thref{spdef5} dann die Behauptung. \bwend

Im Folgenden wollen wir nun eine Äquivalenzrelation auf $\PtJqCma$ einführen. Hierzu benötigen wir zunächst noch folgendes Lemma.

\begin{lemma}	\thlabel{splm7}
	Seien $\alpha\in\R$ und $\phipsi\in\PtJqCma$. Weiterhin sei $g$ eine in $\C\setminus(-\infty,\alpha]$ meromorphe \textit{q}$\times$\textit{q}-Matrixfunktion derart, dass $\det g$ nicht die Nullfunktion ist. Dann gilt $\binom{\phi g}{\psi g}\in\PtJqCma$.
\end{lemma}

\bwanf Seien $\widecheck{\phi}, \widecheck{\psi}:\C\setminus[-\alpha,\infty)\rightarrow\Cqq$ definiert gemäß $\widecheck{\phi}(z) = -\phi(-z)$ bzw. $\widecheck{\psi}(z) = \psi(-z)$. Wegen Teil (a) von \thref{spbm3} gilt dann $\bbinom{\widecheck{\phi}}{\widecheck{\psi}}\in{\cal P}^{(q,q)}_{-\tJq,\geq}(\C$,""$[-\alpha,\infty))$. Weiterhin sei $\widecheck{g}(z)=g(-z)$ für alle $z\in\C\setminus[-\alpha,\infty)$. Dann ist $\widecheck{g}$ eine in $\C\setminus[-\alpha,\infty)$ meromorphe \textit{q}$\times$\textit{q}-Matrixfunktion derart, dass $\det \widecheck{g}$ nicht die Nullfunktion ist. Hieraus folgt wegen \thref{splm6} dann $\bbinom{\widecheck{\phi}\widecheck{g}}{\widecheck{\psi}\widecheck{g}}\in{\cal P}^{(q,q)}_{-\tJq,\geq}(\C$,""$[-\alpha,\infty))$. Hieraus folgt wegen Teil (a) von \thref{spbm3} wiederum $\binom{\phi g}{\psi g}\in\PtJqCma$. \bwend

\thref{splm7} bringt uns nun auf folgende Definition.
  
\begin{defi}	\thlabel{spdef5b}
  	Seien $\alpha\in\R$ und $\binom{\phi_1}{\psi_1}, \binom{\phi_2}{\psi_2} \in \PtJqCma$.
  	Dann heißen $\binom{\phi_1}{\psi_1}$ und $\binom{\phi_2}{\psi_2}$ äquivalent, falls eine in $\C\setminus(-\infty,\alpha]$ meromorphe \textit{q}$\times$\textit{q}-Matrixfunktion $g$ und eine diskrete Teilmenge $\D$ von $\C\setminus(-\infty,\alpha]$ existieren mit
  	\begin{itemize}
    	\item [\rm{(i)}] $\phi_1, \psi_1, \phi_2, \psi_2$ und $g$ sind in $\C\setminus\rklam{(-\infty,\alpha]\cup\D}$ holomorph.
    	\item [\rm{(ii)}] Es gelten $\det g(z) \neq 0$ und
    	\begin{align*}
      		\binom{\phi_2(z)}{\psi_2(z)} = \binom{\phi_1(z)}{\psi_1(z)}g(z)
    	\end{align*}
    	für alle $z \in \C\setminus\rklam{(-\infty,\alpha]\cup\D}$.
  	\end{itemize}
\end{defi}

\begin{bem}	\thlabel{spbm14}
	Man kann sich leicht davon überzeugen, dass der gemäß \thref{spdef5b} definierte Äquivalenzbegriff eine Äquivalenzrelation auf $\PtJqCma$ ist. Mit $\aphipsi$ bezeichnen wir die durch $\phipsi\in\PtJqCma$ erzeugte Äquivalenzklasse.
\end{bem}

Wir zeigen nun, dass die in Teil (a) von \thref{spdef5} eingeführte Klasse von \textit{q}$\times$\textit{q}-Stieltjes-Paaren in $\C\setminus(-\infty,\alpha]$ als eine projektive Erweiterung der Funktionenklasse $\Sqma$ (vergleiche \thref{asmbz4}) aufgefasst werden kann.

\begin{bem}	\thlabel{spbm15}
	Sei $\alpha\in\R$. \dgfa
	\begin{itemize}
		\item [\rm{(a)}] Sei $S\in\Sqma$. Weiterhin bezeichne ${\cal I}$ die in $\C\setminus(-\infty,\alpha]$ konstante Matrixfunktion mit dem Wert $\Iq$. Dann ist $\binom{S}{\cal I}$ ein eigentliches Paar aus $\dPtJqCma$.
		\item [\rm{(b)}] Seien $S_1,S_2\in\Sqma$ so beschaffen, dass die Paare $\binom{S_1}{\cal I}$ und $\binom{S_2}{\cal I}$ äquivalent sind. Dann gilt $S_1 = S_2$.
	\end{itemize}
\end{bem}

\bwanf Bezeichne $\widecheck{\cal I}$ die in $\C\setminus[-\alpha,\infty)$ konstante Matrixfunktion mit dem Wert $\Iq$.

Zu (a): Sei $\widecheck{S}:\C\setminus[-\alpha,\infty)\rightarrow\Cqq$ definiert gemäß $\widecheck{S}(z) := -S(-z)$. Wegen \thref{asmbm4} gilt dann $\widecheck{S}\in{\cal S}_{q,[-\alpha,\infty)}$. Wegen Teil (a) von \thref{spbm11} ist dann $\bbinom{\widecheck{S}}{\widecheck{\cal I}}$ ein eigentliches Paar aus $\widehat{\cal P}^{(q,q)}_{-\tJq,\geq}(\C$,""$[-\alpha,\infty))$. Wegen der Teile (b) und (c) von \thref{spbm3} gilt nun, dass $\binom{S}{\cal I}$ ein eigentliches Paar aus $\dPtJqCma$ ist.

Zu (b): Seien $\widecheck{S}_1,\widecheck{S}_2:\C\setminus[-\alpha,\infty)\rightarrow\Cqq$ definiert gemäß $\widecheck{S}_1(z) := -S_1(-z)$ bzw. $\widecheck{S}_2(z) := -S_2(-z)$. Wegen \thref{asmbm4} gelten dann $\widecheck{S}_1,\widecheck{S}_2\in{\cal S}_{q,[-\alpha,\infty)}$. Wegen \thref{spdef5b} existieren eine in $\C\setminus(-\infty,\alpha]$ meromorphe \textit{q}$\times$\textit{q}-Matrixfunktion $g$ und eine diskrete Teilmenge $\D$ von $\C\setminus(-\infty,\alpha]$ mit
\begin{itemize}
	\item [\rm{(i)}] $g$ ist in $z\in\C\setminus((-\infty,\alpha]\cup\D)$ holomorph.
	\item [\rm{(ii)}] Es gelten $\det g(z)\neq0$ und
	\begin{align*}
		\begin{pmatrix}-\widecheck{S}_2(-z)\\\Iq\end{pmatrix} = \binom{S_2(z)}{\Iq} = \binom{S_1(z)}{\Iq}g(z) = \begin{pmatrix}-\widecheck{S}_1(-z)\\\Iq\end{pmatrix}g(z)
	\end{align*}
	für alle $z\in\C\setminus((-\infty,\alpha]\cup\D)$.
\end{itemize}
Hieraus folgt wegen \thref{spdef4b} dann, dass $\bbinom{-\widecheck{S}_1}{\widecheck{\cal I}}$ und $\bbinom{-\widecheck{S}_2}{\widecheck{\cal I}}$ äquivalent sind. Hieraus folgt wegen Teil (b) von \thref{spbm11} nun $-\widecheck{S}_1 = -\widecheck{S}_2$, also sogar $S_1 = S_2$. \bwend

Wir betrachten nun den umgekehrten Fall von \thref{spbm15}, indem wir aus einem gegebenen \textit{q}$\times$\textit{q}-Stieltjes-Paar in $\C\setminus(-\infty,\alpha]$ eine zu $\Sqma$ gehörige Matrixfunktion konstruieren.

\begin{bem}	\thlabel{spbm16}
	Seien $\alpha\in\R$ und $\phipsi$ ein eigentliches Paar aus $\PtJqCma$. \dgfa
	\begin{itemize}
		\item [\rm{(a)}] Bezeichne ${\cal I}$ die in $\C\setminus(-\infty,\alpha]$ konstante Matrixfunktion mit dem Wert $\Iq$. Dann ist $\binom{\phi\psi^{-1}}{\cal I}$ ein eigentliches Paar aus $\PtJqCma$, welches zu $\phipsi$ äquivalent ist.
		\item [\rm{(b)}] Es gilt $\phi\psi^{-1}\in\Sqma$.
	\end{itemize}
\end{bem}

\bwanf Zu (a): Unter Beachtung von $\binom{\phi\psi^{-1}}{\cal I}=\phipsi \psi^{-1}$ gilt wegen \thref{splm7}, Teil (b) von \thref{spdef5} sowie $\rank{\cal I}(z)=q$ für alle $z\in\C\setminus(-\infty,\alpha]$, dass $\binom{\phi\psi^{-1}}{\cal I}$ ein eigentliches Paar aus $\PtJqCma$ ist, und wegen \thref{spdef5b} weiterhin, dass $\binom{\phi\psi^{-1}}{\cal I}$ äquivalent ist zu $\phipsi$.

Zu (b): Seien $\widecheck{\phi}, \widecheck{\psi}:\C\setminus[-\alpha,\infty)\rightarrow\Cqq$ definiert gemäß $\widecheck{\phi}(z) = -\phi(-z)$ bzw. \linebreak $\widecheck{\psi}(z) = \psi(-z)$. Wegen der Teile (a) und (c) von \thref{spbm3} ist dann $\bbinom{\widecheck{\phi}}{\widecheck{\psi}}$ ein eigentliches Paar aus ${\cal P}^{(q,q)}_{-\tJq,\geq}(\C$,""$[-\alpha,\infty))$. Hieraus folgt wegen Teil (b) von \thref{spbm13} dann $\widecheck{\phi}\widecheck{\psi}^{-1}\in{\cal S}_{q,[-\alpha,\infty)}$. Hieraus folgt wegen \thref{asmbm4} nun $\phi\psi^{-1}\in\Sqma$. \bwend

Mithilfe von \thref{spbm15} und \thref{spbm16} erkennen wir, dass eine Bijektion zwischen $\Sqma$ und der Menge der Äquivalenzklassen von eigentlichen \textit{q}$\times$\textit{q}-Stieltjes-Paaren in $\C\setminus(-\infty,\alpha]$ besteht.

Wir geben nun Beispiele zweier bemerkenswerter \textit{q}$\times$\textit{q}-Stieltjes-Paare in $\C\setminus(-\infty,\alpha]$ an, welche in unseren Überlegungen für Abschnitt \ref{EEl} eine wesentliche Rolle spielen.

\begin{beispiel}	\thlabel{spbsp3}
	Sei $\alpha\in\R$. Bezeichne ${\cal I}$ bzw. ${\cal O}$ die in $\C\setminus(-\infty,\alpha]$ konstante Matrixfunktion mit dem Wert $\Iq$ bzw. $\Oq$. Dann gilt $\binom{\cal I}{\cal O}, \binom{\cal O}{\cal I} \in \dPtJqCma$. Insbesondere ist das Paar $\binom{\cal O}{\cal I}$ eigentlich.
\end{beispiel}

\bwanf Dies folgt sogleich aus \thref{spbsp2} und den Teilen (b) und (c) von \thref{spbm3}. \bwend

Unsere folgenden Untersuchungen sind einer alternativen Beschreibung der Menge $\PtJqCa$ gewidmet. Statt der ersten Bedingung von Teil (iii) von \thref{spdef4} wollen wir eine äquivalente Ungleichung mithilfe der Löwner-Halbordnung bezüglich der Signaturmatrix $\Jq$ (vergleiche \thref{spbsp4}) auf der Halbebene $\Lam$ finden. Wir verwenden die gleiche Vorgehensweise wie in \cite[Abschnitt 8.4]{Wb}. Hierfür benötigen wir zunächst noch folgende Lemmas (vergleiche \thref{splm5} mit \cite[Lemma 10.15]{Maka} oder \cite[Lemma 8.16]{Wb}, \thref{splm4} mit \cite[Lemma 8.38]{Wb} und \thref{splm3} mit \cite[Lemma 8.35]{Wb}).

\begin{lemma}	\thlabel{splm5}
	Seien $\alpha\in\R$ und $\phipsi\in\PtJqCa$. Weiterhin sei $\D$ die gemäß Teil (a) von \thref{spdef4} existierende diskrete Teilmenge von $\C\setminus[\alpha,\infty)$. Dann gilt
	\begin{align*}
		\phipsiz^{\ast}\rklam{-\Jq}\phipsiz \geq \Oq
	\end{align*}
	für alle $z\in \Lam\setminus\D$.
\end{lemma}

\bwanf Wegen Teil (a) von \thref{spdef4} gelten
\begin{align}	\label{splm5bw1a}
	\binom{(z-\alpha)\phi(z)}{\psi(z)}^{\ast}\bbrklam{\frac{-\tJq}{2\im z}}\binom{(z-\alpha)\phi(z)}{\psi(z)} \geq \Oq
\end{align}
und
\begin{align}	\label{splm5bw1b}
	\binom{\phi(z)}{\psi(z)}^{\ast}\bbrklam{\frac{-\tJq}{2\im z}}\binom{\phi(z)}{\psi(z)} \geq \Oq
\end{align}
für alle $z \in \C\setminus\rklam{\R\cup\D}$. Wegen \thref{spbsp1} und \fref{splm5bw1b} gilt dann
\begin{align}	\label{splm5bw2}
	\frac{\alpha-\re z}{\im z}\im\eklam{\psi^{\ast}(z)\phi(z)}
	= (\alpha-\re z)\phipsiz^{\ast}\bbrklam{\frac{-\tJq}{2\im z}}\phipsiz
	\geq \Oq
\end{align}
für alle $z\in \Lam\setminus(\R\cup\D)$. Unter Beachtung von
\begin{align*}
	\im\eklam{zA} = \re z \im A + \im z \re A
\end{align*}
für beliebige $z\in\C$ und $A\in\Cqq$ (vergleiche z.\,B. \cite[Bemerkung A.23]{Sch}) gilt wegen \fref{splm5bw2}, \thref{spbsp1} und \fref{splm5bw1a} dann
\begin{align*}
	\re\eklam{\psi^{\ast}(z)\phi(z)}
	&= \frac{1}{\im z}\im\eklam{z\psi^{\ast}(z)\phi(z)} - \frac{\re z}{\im z}\im\eklam{\psi^{\ast}(z)\phi(z)} \notag \\
	&= \frac{1}{\im z}\im\eklam{z\psi^{\ast}(z)\phi(z)} + \frac{\alpha-\re z}{\im z}\im\eklam{\psi^{\ast}(z)\phi(z)} - \frac{\alpha}{\im z}\im\eklam{\psi^{\ast}(z)\phi(z)} \notag \\
	&\geq \frac{1}{\im z}\im\eklam{z\psi^{\ast}(z)\phi(z)} - \frac{1}{\im z}\im\eklam{\alpha\psi^{\ast}(z)\phi(z)} \notag \\
	&= \frac{1}{\im z}\im\eklam{(z-\alpha)\psi^{\ast}(z)\phi(z)} \notag \\
	&= \binom{(z-\alpha)\phi(z)}{\psi(z)}^{\ast}\bbrklam{\frac{-\tJq}{2\im z}}\binom{(z-\alpha)\phi(z)}{\psi(z)} 
	\geq \Oq
\end{align*}
für alle $z\in \Lam\setminus(\R\cup\D)$. Hieraus folgt wegen der Stetigkeit der Funktionen $\phi$ und $\psi$ in $\C\setminus([\alpha,\infty)\cup\D)$ nun
\begin{align*}
	\re\eklam{\psi^{\ast}(z)\phi(z)} \geq \Oq
\end{align*}
für alle $z\in \Lam\setminus\D$. Hieraus folgt wegen \thref{spbsp4} dann die Behauptung. \bwend

\begin{lemma}	\thlabel{splm4}
	Seien $\alpha\in\R$ sowie $\phi$ und $\psi$ in $\C\setminus[\alpha,\infty)$ meromorphe Matrixfunktionen derart, dass $\det\psi$ nicht die Nullfunktion ist und eine diskrete Teilmenge $\D$ von \linebreak $\C\setminus[\alpha,\infty)$ mit folgenden Eigenschaften existiert:
	\begin{itemize}
		\item [\rm{(i)}] $\phi$ und $\psi$ sind in $\C\setminus\rklam{[\alpha,\infty)\cup\D}$ holomorph.
    	\item [\rm{(ii)}] Es gilt $\rank \phipsiz = q$ für alle $z \in \C\setminus\rklam{[\alpha,\infty)\cup\D}$.
    	\item [\rm{(iii)}] Es gilt
    	\begin{align*}
      		\phipsiz^{\ast}\bbrklam{\frac{-\tJq}{2\im z}}\phipsiz \geq \Oq
    	\end{align*}
    	für alle $z \in \C\setminus\rklam{\R\cup\D}$.
    	\item [\rm{(iv)}] Es gilt 
    	\begin{align*}
    		\phipsiz^{\ast}\rklam{-\Jq}\phipsiz \geq \Oq
    	\end{align*}
    	für alle $z\in\Lam\setminus\D$.
	\end{itemize}
	Dann gilt
	\begin{align*}
		\binom{(z-\alpha)\phi(z)}{\psi(z)}^{\ast}\bbrklam{\frac{-\tJq}{2\im z}}\binom{(z-\alpha)\phi(z)}{\psi(z)} \geq \Oq
	\end{align*}
	für alle $z \in \C\setminus\rklam{\R\cup\D}$.
\end{lemma}

\bwanf Da $\det\psi$ nicht die Nullfunktion ist, existiert unter Beachtung des Identitätssatzes für meromorphe Funktionen (vergleiche z.\,B. im skalaren Fall \cite[Satz 10.3.2]{Funk}; im matriziellen Fall betrachtet man die einzelnen Einträge der Matrixfunktion) eine diskrete Teilmenge $\D_1$ von $\C\setminus[\alpha,\infty)$ mit 
\begin{align}	\label{splm4bw1}
	\det\psi(z)\neq0
\end{align}
für alle $z \in \C\setminus\rklam{[\alpha,\infty)\cup\D_1}$. Sei nun $\D_2=\D\cup\D_1$. Wegen \rm{(i)} und \fref{splm4bw1} ist dann $\phi\psi^{-1}$ in $\C\setminus\rklam{[\alpha,\infty)\cup\D_2}$ holomorph. Unter Beachtung von
\begin{align*}
	\im\eklam{B^{\ast}AB} = \frac{1}{2i}\rklam{B^{\ast}AB-B^{\ast}A^{\ast}B }
	= B^{\ast}\frac{1}{2i}\rklam{A-A^{\ast}}B = B^{\ast}(\im A)B
\end{align*}
und
\begin{align*}
	\re\eklam{B^{\ast}AB} = \frac{1}{2}\rklam{B^{\ast}AB+B^{\ast}A^{\ast}B }
	= B^{\ast}\frac{1}{2}\rklam{A+A^{\ast}}B = B^{\ast}(\re A)B
\end{align*}
für beliebige $A,B\in\Cqq$ gelten
\begin{align}	\label{splm4bw2}
	\im\eklam{\psi^{\ast}(z)\phi(z)} = \im\eklam{\psi^{\ast}(z)\phi(z)\psi^{-1}(z)\psi(z)} = \psi^{\ast}(z)\im\eklam{\phi(z)\psi^{-1}(z)}\psi(z)
\end{align}
und
\begin{align}	\label{splm4bw3}
	\re\eklam{\psi^{\ast}(z)\phi(z)} = \re\eklam{\psi^{\ast}(z)\phi(z)\psi^{-1}(z)\psi(z)} = \psi^{\ast}(z)\re\eklam{\phi(z)\psi^{-1}(z)}\psi(z)
\end{align}
für alle $z \in \C\setminus\rklam{[\alpha,\infty)\cup\D_1}$. Wegen \rm{(iii)} und \thref{spbsp1} gilt
\begin{align*}
	\frac{1}{\im z}\im\eklam{\psi^{\ast}(z)\phi(z)} = \phipsiz^{\ast}\bbrklam{\frac{-\tJq}{2\im z}}\phipsiz \geq \Oq
\end{align*}
für alle $z \in \C\setminus\rklam{\R\cup\D}$. Hieraus folgt wegen \fref{splm4bw2} und \fref{splm4bw1} dann
\begin{align*}
	\frac{1}{\im z}\im\eklam{\phi(z)\psi^{-1}(z)} \geq \Oq
\end{align*}
für alle $z \in \C\setminus\rklam{\R\cup\D_2}$. Hieraus folgen nun
\begin{align}	\label{splm4bw4}
	\im\eklam{\phi(z)\psi^{-1}(z)} \geq \Oq
\end{align}
für alle $z \in \Pp\setminus\D_2$ und
\begin{align}	\label{splm4bw5}
	-\im\eklam{\phi(z)\psi^{-1}(z)} \geq \Oq
\end{align}
für alle $z \in \Pm\setminus\D_2$. Wegen \rm{(iv)} und \thref{spbsp4} gilt
\begin{align*}
	\re\eklam{\psi^{\ast}(z)\phi(z)} \geq \Oq 
\end{align*}
für alle $z\in\Lam\setminus\D$. Hieraus folgt wegen \fref{splm4bw3} und \fref{splm4bw1} dann
\begin{align}	\label{splm4bw6}
	\re\eklam{\phi(z)\psi^{-1}(z)} \geq \Oq
\end{align}
für alle $z\in\Lam\setminus\D_2$. Wegen \fref{splm4bw4}, \fref{splm4bw5}, \fref{splm4bw6} sowie unter Beachtung von
\begin{align*}
	\re\eklam{iA} = \frac{1}{2}\brklam{iA+\eklam{iA}^{\ast}}=\frac{1}{2i}\brklam{{-A-\eklam{-A}^{\ast}}} = \im\eklam{-A}
\end{align*}
für alle $A\in\Cqq$ und \cite[Lemma 2.19]{DFK} oder \cite[Lemma 4.23]{Sch} ist jedes $z\in\D_2$ eine hebbare Singularität von $\phi\psi^{-1}$ und $\phi\psi^{-1}$ lässt sich zu einer auf $\C\setminus[\alpha,\infty)$ holomorphen Matrixfunktion $F$ fortsetzen mit
\begin{itemize}
   	\item [\rm{(v)}] Es gilt $\im F(z) \geq \Oq$ für alle $z \in \Pp$.
   	\item [\rm{(vi)}] Es gilt $-\im F(z) \geq \Oq$ für alle $z \in \Pm$.
    \item [\rm{(iv)}] Es gilt $\re F(z) \geq \Oq$ für alle $z\in\Lam$.
\end{itemize}
Hieraus folgt wegen \cite[Proposition 4.4]{SF} dann $F\in\Sqa$. Hieraus folgt wegen \cite[Lemma 4.2]{SF} nun
\begin{align*}
	\frac{1}{\im z}\im\eklam{(z-\alpha)F(z)} \geq \Oq
\end{align*}
für alle $z \in \C\setminus\R$. Hieraus folgt wegen \fref{splm4bw1} dann
\begin{align*}
	\frac{1}{\im z}\im\eklam{(z-\alpha)\psi^{\ast}(z)\phi(z)}
	&= \frac{1}{\im z}\im\eklam{(z-\alpha)\psi^{\ast}(z)\phi(z)\psi^{-1}(z)\psi(z)} \\
	&= \frac{1}{\im z}\psi^{\ast}(z)\im\eklam{(z-\alpha)F(z)}\psi(z) \geq \Oq
\end{align*}
für alle $z \in \C\setminus\rklam{\R\cup\D_2}$. Hieraus folgt wegen der Stetigkeit von $\phi$ und $\psi$ in \linebreak $\C\setminus\rklam{[\alpha,\infty)\cup\D}$ nun
\begin{align*}
	\binom{(z-\alpha)\phi(z)}{\psi(z)}^{\ast}\bbrklam{\frac{-\tJq}{2\im z}}\binom{(z-\alpha)\phi(z)}{\psi(z)} = \frac{1}{\im z}\im\eklam{(z-\alpha)\psi^{\ast}(z)\phi(z)} \geq \Oq
\end{align*}
für alle $z \in \C\setminus\rklam{\R\cup\D}$. \bwend

\begin{lemma}	\thlabel{splm3}
	Seien $\alpha\in\R$ sowie $\phi$ und $\psi$ in $\C\setminus[\alpha,\infty)$ meromorphe Matrixfunktionen derart, dass eine diskrete Teilmenge $\D$ von $\C\setminus[\alpha,\infty)$ mit folgenden Eigenschaften existiert:
	\begin{itemize}
		\item [\rm{(i)}] $\phi$ und $\psi$ sind in $\C\setminus\rklam{[\alpha,\infty)\cup\D}$ holomorph.
    	\item [\rm{(ii)}] Es gilt $\rank \phipsiz = q$ für alle $z \in \C\setminus\rklam{[\alpha,\infty)\cup\D}$.
    	\item [\rm{(iii)}] Es gilt
    	\begin{align*}
      		\phipsiz^{\ast}\bbrklam{\frac{-\tJq}{2\im z}}\phipsiz \geq \Oq
    	\end{align*}
    	für alle $z \in \C\setminus\rklam{\R\cup\D}$.
    	\item [\rm{(iv)}] Es gilt 
    	\begin{align*}
    		\phipsiz^{\ast}\rklam{-\Jq}\phipsiz \geq \Oq
    	\end{align*}
    	für alle $z\in\Lam\setminus\D$.
	\end{itemize}
	\dgfa
	\begin{itemize}
		\item [\rm{(a)}] Sei $\eps\in(0,\infty)$. \dgfa
	\begin{itemize}
		\item [\rm{(v)}] $\phi$ und $\eps\phi+\psi$ sind in $\C\setminus\rklam{[\alpha,\infty)\cup\D}$ holomorph.
    	\item [\rm{(vi)}] Es gilt $\rank \binom{\phi(z)}{\eps\phi(z)+\psi(z)} = q$ für alle $z \in \C\setminus\rklam{[\alpha,\infty)\cup\D}$.
    	\item [\rm{(vii)}] Es gilt
    	\begin{align*}
      		\binom{\phi(z)}{\eps\phi(z)+\psi(z)}^{\ast}\bbrklam{\frac{-\tJq}{2\im z}}\binom{\phi(z)}{\eps\phi(z)+\psi(z)} \geq \Oq
    	\end{align*}
    	für alle $z \in \C\setminus\rklam{\R\cup\D}$.
    	\item [\rm{(viii)}] Es gilt
    	\begin{align*} 
    		\binom{\phi(z)}{\eps\phi(z)+\psi(z)}^{\ast}\rklam{-\Jq}\binom{\phi(z)}{\eps\phi(z)+\psi(z)} \geq \Oq
    	\end{align*}
    	für alle $z\in\Lam\setminus\D$.
    	\item [\rm{(ix)}] Es gilt $\det[\eps\phi(z)+\psi(z)]\neq0$ für alle $z\in\Lam\setminus([\alpha,\infty)\cup\D)$.
	\end{itemize}
		\item [\rm{(b)}] Es gilt
		\begin{align*}
			\binom{(z-\alpha)\phi(z)}{\psi(z)}^{\ast}\bbrklam{\frac{-\tJq}{2\im z}}\binom{(z-\alpha)\phi(z)}{\psi(z)} \geq \Oq
		\end{align*}
		für alle $z \in \C\setminus\rklam{\R\cup\D}$.
	\end{itemize}
\end{lemma}

\bwanf Zu \rm{(a)}: Aus \rm{(i)} folgt unmittelbar \rm{(v)}. 
Unter Beachtung von
\begin{align}	\label{splm3bw1}
	\binom{\phi(z)}{\eps\phi(z)+\psi(z)} = \begin{pmatrix} \Iq & \Oq \\ \eps\Iq & \Iq \end{pmatrix} \phipsiz
\end{align}
für alle $z \in \C\setminus\rklam{[\alpha,\infty)\cup\D}$ und
\begin{align*}
	\det\begin{pmatrix} \Iq & \Oq \\ \eps\Iq & \Iq \end{pmatrix} = 1 \neq 0
\end{align*}
folgt aus \rm{(ii)} dann \rm{(vi)}.
Unter Beachtung von \fref{splm3bw1} und 
\begin{align*}
	\begin{pmatrix} \Iq & \Oq \\ \eps\Iq & \Iq \end{pmatrix}^{\ast} \brklam{{-\tJq}} \begin{pmatrix} \Iq & \Oq \\ \eps\Iq & \Iq \end{pmatrix} 
	= \begin{pmatrix} \Iq & \eps\Iq \\ \Oq & \Iq \end{pmatrix} \begin{pmatrix} i\eps\Iq & i\Iq \\ -i\Iq & \Oq \end{pmatrix} = -\tJq
\end{align*}
folgt aus \rm{(iii)} dann \rm{(vii)}. 
Wegen \rm{(iv)} und \thref{spbsp4} gilt
\begin{align}	\label{splm3bw2}
	\psi^{\ast}(z)\phi(z) + \phi^{\ast}(z)\psi(z) = 2\re\eklam{\psi^{\ast}(z)\phi(z)} = \phipsiz^{\ast}\rklam{-\Jq}\phipsiz \geq \Oq
\end{align}
für alle $z\in\Lam\setminus\D$.
Hieraus folgt unter Beachtung von $\eps\phi^{\ast}(z)\phi(z)\geq\Oq$ für alle $z\in\C\setminus([\alpha,\infty)\cup\D)$ und \thref{spbsp4} dann
\begin{align*}
	\binom{\phi(z)}{\eps\phi(z)+\psi(z)}^{\ast}\rklam{-\Jq}\binom{\phi(z)}{\eps\phi(z)+\psi(z)} 
	&= 2\re\brklam{[\eps\phi(z)+\psi(z)]^{\ast}\phi(z)} \\
	&= \brklam{[\eps\phi(z)+\psi(z)]^{\ast}\phi(z) + \phi^{\ast}(z)[\eps\phi(z)+\psi(z)]} \\
	&= 2\eps\phi^{\ast}(z)\phi(z) + \brklam{\psi^{\ast}(z)\phi(z) + \phi^{\ast}(z)\psi(z)} \\
	&\geq \Oq
\end{align*}
für alle $z\in\Lam\setminus\D$. Wegen \rm{(ii)} gilt
\begin{align*}
	\phi^{\ast}(z)\phi(z)+\psi^{\ast}(z)\psi(z) = \phipsiz^{\ast}\phipsiz > \Oq
\end{align*}
für alle $z \in \C\setminus\rklam{[\alpha,\infty)\cup\D}$. Hieraus folgt wegen \fref{splm3bw2} nun
\begin{align}	\label{splm3bw3}
	&\ \eklam{\eps\phi(z)+\psi(z)}^{\ast}\eklam{\phi(z)+\eps\psi(z)} + \eklam{\phi(z)+\eps\psi(z)}^{\ast}\eklam{\eps\phi(z)+\psi(z)} \notag \\
	&= 2\eps\eklam{\phi^{\ast}(z)\phi(z)+\psi^{\ast}(z)\psi(z)} + (1+\eps)\eklam{\phi^{\ast}(z)\psi(z)+\psi^{\ast}(z)\phi(z)} > \Oq
\end{align}
für alle $z\in\Lam\setminus([\alpha,\infty)\cup\D)$.
Seien nun $z\in\Lam\setminus([\alpha,\infty)\cup\D)$ und \linebreak $u\in{\cal N}(\eps\phi(z)+\psi(z))$. Dann gilt
\begin{align*}
	&\ u^{\ast}\brklam{\eklam{\eps\phi(z)+\psi(z)}^{\ast}\eklam{\phi(z)+\eps\psi(z)} + \eklam{\phi(z)+\eps\psi(z)}^{\ast}\eklam{\eps\phi(z)+\psi(z)}}u \\
	&= 0^{\ast}_{q\times1}\eklam{\phi(z)+\eps\psi(z)}u + u^{\ast}\eklam{\phi(z)+\eps\psi(z)}^{\ast}0_{q\times1} = 0.
\end{align*}
Hieraus folgt wegen \fref{splm3bw3} dann $u=0_{q\times1}$ und somit $\det[\eps\phi(z)+\psi(z)]\neq0$. 

Zu \rm{(b)}: Wegen \rm{(a)} und \thref{splm4} gilt
\begin{align*}
	\binom{(z-\alpha)\phi(z)}{\eps\phi(z)+\psi(z)}^{\ast}\bbrklam{\frac{-\tJq}{2\im z}}\binom{(z-\alpha)\phi(z)}{\eps\phi(z)+\psi(z)} \geq \Oq
\end{align*}
für alle $z \in \C\setminus\rklam{\R\cup\D}$ und $\eps\in(0,\infty)$. Mit der Grenzwertbetrachtung $\eps\rightarrow0$ folgt dann die Behauptung. \bwend

Es folgt nun das Hauptresultat dieses Abschnitts (vergleiche \cite[Satz 8.39]{Wb}).

\begin{satz}	\thlabel{spsa2}
	Seien $\alpha\in\R$ sowie $\phi$ und $\psi$ in $\C\setminus[\alpha,\infty)$ meromorphe Matrixfunktionen. Dann gilt $\phipsi\in\PtJqCa$ genau dann, wenn eine diskrete Teilmenge $\D$ von $\C\setminus[\alpha,\infty)$ mit folgenden Eigenschaften existiert:
	\begin{itemize}
		\item [\rm{(i)}] $\phi$ und $\psi$ sind in $\C\setminus\rklam{[\alpha,\infty)\cup\D}$ holomorph.
    	\item [\rm{(ii)}] Es gilt $\rank \phipsiz = q$ für alle $z \in \C\setminus\rklam{[\alpha,\infty)\cup\D}$.
    	\item [\rm{(iii)}] Es gilt
    	\begin{align*}
      		\phipsiz^{\ast}\bbrklam{\frac{-\tJq}{2\im z}}\phipsiz \geq \Oq
    	\end{align*}
    	für alle $z \in \C\setminus\rklam{\R\cup\D}$.
    	\item [\rm{(iv)}] Es gilt 
    	\begin{align*}
    		\phipsiz^{\ast}\rklam{-\Jq}\phipsiz \geq \Oq
    	\end{align*}
    	für alle $z\in\Lam\setminus\D$.
	\end{itemize}
	Insbesondere gilt $\phipsi\in\dPtJqCa$ genau dann, wenn {\rm (i)}-{\rm (iv)} mit $\D=\emptyset$ erfüllt sind.
\end{satz}

\bwanf Sei zunächst $\phipsi\in\PtJqCa$. Dann sind wegen Teil (a) von \thref{spdef4} die Aussagen (i)-(iii) erfüllt. Wegen \thref{splm5} gilt dann auch (iv).

Sei nun umgekehrt (i)-(iv) erfüllt. Wegen Teil (b) von \thref{splm3} gilt dann
\begin{align*}
	\binom{(z-\alpha)\phi(z)}{\psi(z)}^{\ast}\bbrklam{\frac{-\tJq}{2\im z}}\binom{(z-\alpha)\phi(z)}{\psi(z)} \geq \Oq
\end{align*}
für alle $z \in \C\setminus\rklam{\R\cup\D}$. Unter Beachtung von Teil (a) von \thref{spdef4} folgt hieraus wegen (i)-(iii) dann $\phipsi\in\PtJqCa$.

Da die gewählte diskrete Teilmenge $\D$ mit derjenigen aus Teil (a) von \thref{spdef4} übereinstimmt, gilt $\phipsi\in\dPtJqCa$ genau dann, wenn {\rm (i)}-{\rm (iv)} mit $\D=\emptyset$ erfüllt sind. \bwend

\section[Einige Aussagen über Teilklassen von Schur"=Funktionen auf Halbebenen]{Einige Aussagen über Teilklassen von Schur- \\ Funktionen auf Halbebenen}  \label{chapsf}

In diesem Abschnitt behandeln wir zwei spezielle Teilklassen von Schur-Funktion auf der oberen bzw. unteren offenen Halbebene von $\C$, die in Kapitel \ref{chapar} für eine weitere Beschreibung der Lösungsmenge des vollständig nichtdegenerierten matriziellen $\alpha$-Stieltjes-Momentenproblems als Parametermenge Verwendung finden. Hierbei legen wir besonderen Augenmerk auf eine Verbindung zu den \textit{q}$\times$\textit{q}-Stieltjes-Paaren aus Kapitel \ref{chapsp}. Wir beginnen nun mit den zentralen Begriffsbildungen für diesen Abschnitt.

\begin{defi}	\thlabel{spdef6} 
	Sei $\G$ ein Gebiet von $\C$.
	\begin{itemize}
		\item [\rm{(a)}] Es heißt $S:\G \rightarrow \Cpq$ \textbf{p$\times$q-Schur-Funktion} auf $\G$, falls $S$ in $\G$ holomorph und $\Iq - S^{\ast}(z)S(z) \geq \Oq$ für alle $z \in \G$ erfüllt ist.
  		Mit $\SpqG$ bezeichnen wir die Menge aller \textit{p}$\times$\textit{q}-Schur-Funktionen auf $\G$.
  		\item [\rm{(b)}] Sei $\alpha\in\R$. Dann bezeichne $\Dqa$ die Menge aller $S\in\SqPp$ mit
  		\begin{itemize}
  			\item [\rm{(i)}] Es gilt $\im S(z) \geq \Oq$ für alle $z\in\Pp$ mit $\re z\in(-\infty,\alpha)$.
  			\item [\rm{(ii)}] Für alle $x\in(-\infty,\alpha)$ existiert $U_x := \lim_{z\rightarrow x} S(z)$ und es ist $U_x$ unitär.
  		\end{itemize}
  		\item [\rm{(c)}] Sei $\alpha\in\R$. Dann bezeichne $\Eqa$ die Menge aller $S\in\SqPm$ mit
  		\begin{itemize}
  			\item [\rm{(iii)}] Es gilt $-\im S(z) \geq \Oq$ für alle $z\in\Pm$ mit $\re z\in(\alpha,\infty)$.
  			\item [\rm{(iv)}] Für alle $x\in(\alpha,\infty)$ existiert $U_x := \lim_{z\rightarrow x} S(z)$ und es ist $U_x$ unitär.
  		\end{itemize}
  	\end{itemize}
\end{defi}

Folgendes Resultat findet man in ähnlicher Form auch unter \cite[Lemma 4.21]{Sch} oder \cite[Lemma 1.4]{Maka}.

\begin{lemma}	\thlabel{splm8}
	Seien $\G$ ein Gebiet von $\C$ und $\D$ eine nichtleere diskrete Teilmenge von $\G$. Weiterhin sei $S:\G\setminus\D\rightarrow\Cpq$ eine holomorphe Matrixfunktion derart, dass
	\begin{align*}
		\Iq-S^{\ast}(z)S(z) \geq \Oq
	\end{align*}
	für alle $z\in\G\setminus\D$ erfüllt ist. Dann lässt sich $S$ zu einer Funktion aus $\SpqG$ fortsetzen.
\end{lemma}

\bwanf Es gilt $S^{\ast}(z)S(z) \leq \Iq$ für alle $z\in\G\setminus\D$. Unter Beachtung von \linebreak $\enorm{A}^2 = \tr (A^{\ast}A)$ für alle $A\in\Cpq$ folgt dann
\begin{align*}
	\enorm{S(z)} \leq \enorm{\Iq} = \sqrt{q}
\end{align*}
für alle $z\in\G\setminus\D$. Sei nun $S:=(S_{ij})_{i\in\Zep,j\in\Zeq}$. Dann gilt $\abs{S_{ij}(z)}\leq\sqrt{q}$ für alle $i\in\Zep$, $j\in\Zeq$ und $z\in\G\setminus\D$. Wegen des Riemannschen Fortsetzungssatzes (vergleiche z.\,B. \cite[Satz 7.3.3]{Funk}) lässt sich dann $S_{ij}$ für alle $i\in\Zep$ und $j\in\Zeq$ auf eine in $\G$ holomorphe Funktion fortsetzen. Somit lässt sich auch $S$ auf eine in $\G$ holomorphe Funktion fortsetzen. Aus der Stetigkeit jener Fortsetzung und Teil (a) von \thref{spdef6} folgt dann die Behauptung. \bwend

Das folgende Resultat schafft eine Verbindung zwischen den in den Teilen (b) und (c) von \thref{spdef6} definierten speziellen Teilklassen von Schur-Funktionen auf der oberen und unteren offenen Halbebene von $\C$.

\begin{bem}	\thlabel{spbm4}
	\egfa
	\begin{itemize}
		\item [\rm{(a)}] Seien $\G$ ein Gebiet von $\C$ und $\widecheck{\G} := \gklam{-\omega \;|\; \omega \in \G}$. Weiterhin seien $S:\linebreak\G \rightarrow \Cpq$ und $\widecheck{S}:\widecheck{\G} \rightarrow \Cpq$ definiert gemäß $\widecheck{S}(z) := -S(-z)$. \dsfaa
		\begin{itemize}
			\item [\rm{(i)}] Es gilt $S\in\SpqG$.
			\item [\rm{(ii)}] Es gilt $\widecheck{S}\in\SpqcG$.
		\end{itemize}
		\item [\rm{(b)}] Seien $\alpha\in\R$, $S:\Pp \rightarrow \Cqq$ und $\widecheck{S}:\Pm \rightarrow \Cqq$ definiert gemäß \linebreak $\widecheck{S}(z) := -S(-z)$. \dsfaa
		\begin{itemize}
			\item [\rm{(iii)}] Es gilt $S\in\Dqa$.
			\item [\rm{(iv)}] Es gilt $\widecheck{S}\in\Eqma$.
		\end{itemize}
	\end{itemize}
\end{bem}

\bwanf Zu (a): Offensichtlich ist $S$ genau dann in $\G$ holomorph, wenn $\widecheck{S}$ in $\widecheck{G}$ holomorph ist. Weiterhin gilt
\begin{align*}
	\Iq - S^{\ast}(z)S(z) = \Iq - \widecheck{S}^{\ast}(-z)\widecheck{S}(-z)
\end{align*}
für alle $z\in\G$. Hieraus folgt wegen Teil (a) von \thref{spdef6} dann die Behauptung.

Zu (b): Es gilt
\begin{align}	\label{spbmbw1}
	\widecheck{\Pi}_{+} = \gklam{-\omega \;|\; \omega \in \Pp}
	= \gklam{\omega\in\C \;| -\im\omega > 0}
	= \gklam{\omega\in\C \;|\; \im\omega < 0} = \Pm.
\end{align}
Weiterhin gelten
\begin{align*}
	\im S(z) = -\im \widecheck{S}(-z)
\end{align*}
für alle $z\in\Pp$ und
\begin{align*}
	\lim_{z\rightarrow x}S(z) 
	= -\lim_{z\rightarrow x}\widecheck{S}(-z)
	= -\lim_{\omega\rightarrow -x}\widecheck{S}(\omega)
\end{align*}
für alle $x\in(\infty,\alpha)$. Hieraus folgt unter Beachtung von \fref{spbmbw1} und (a) wegen \thref{spdef6} dann die Behauptung. \bwend

Eine ausführlichere Version des nächsten Satzes, der eine Verbindung zwischen der Menge von \textit{q}$\times$\textit{q}-Stieltjes-Paaren in $\C\setminus[\alpha,\infty)$ und der in Teil (b) von \thref{spdef6} definierten speziellen Teilmenge $\Dqa$ der Menge von Schur-Funktionen auf der oberen offenen Halbebene von $\C$ schafft, findet man unter \cite[Satz 8.41]{Wb} in Verbindung mit \cite[Satz 8.39]{Wb}. Für unsere Zwecke reicht die einfachere Version aus. Zur besseren Anschauung geben wir einen detaillierten Beweis an. Hierfür wird uns die alternative Definition der \textit{q}$\times$\textit{q}-Stieltjes-Paare in $\C\setminus[\alpha,\infty)$ von \thref{spsa2} und folgendes Lemma besondere Dienste erweisen (vergleiche \cite[Lemma 8.40]{Wb}).

\begin{lemma}	\thlabel{splm2}
	Seien $\G$ ein Gebiet von $\C$, $\D$ eine diskrete Teilmenge von $\G$ mit \linebreak $(\G\cap\R)\setminus\D\neq\emptyset$ und $f:\G\setminus\D\rightarrow\Cqq$ stetig mit
	\begin{align*}
		\frac{1}{\im z}f(z) \geq \Oq
	\end{align*}
	für alle $z\in\G\setminus(\R\cup\D)$. Dann gilt $f(x) = \Oq$ für alle $x\in(\G\cap\R)\setminus\D$.
\end{lemma}

\bwanf Seien $x\in(\G\cap\R)\setminus\D$ und $(y_n)^{\infty}_{n=0}$ bzw. $(z_n)^{\infty}_{n=0}$ Folgen aus $(\Pp\cap\G)\setminus\D$ bzw. $(\Pm\cap\G)\setminus\D$ mit
\begin{align}	\label{splm2bw1}
	\lim_{n\rightarrow\infty} y_n = x = \lim_{n\rightarrow\infty} z_n.
\end{align}
Wegen
\begin{align*}
	\frac{1}{\im y_n}f(y_n) \geq \Oq \quad \text{und} \quad \im y_n > 0
\end{align*}
für alle $n\in\No$ gilt $f(y_n)\geq\Oq$ für alle $n\in\No$, also wegen \fref{splm2bw1} und der Stetigkeit von $f$ in $\G\setminus\D$ dann $f(x)\geq\Oq$. Wegen
\begin{align*}
	\frac{1}{\im z_n}f(z_n) \geq \Oq \quad \text{und} \quad -\im z_n > 0
\end{align*}
für alle $n\in\No$ gilt $-f(z_n)\geq\Oq$ für alle $n\in\No$, also wegen \fref{splm2bw1} und der Stetigkeit von $f$ in $\G\setminus\D$ dann $-f(x)\geq\Oq$. Somit folgt $f(x)=\Oq$. \bwend

\begin{satz}	\thlabel{spsa1}
Sei $\alpha\in\R$. \dgfa	
	\begin{itemize}
		\item [\rm{(a)}] Sei $\phipsi\in\PtJqCa$. Dann ist $\det(\psi-i\phi)$ nicht die Nullfunktion. Weiterhin sei $\tF := (\psi+i\phi)(\psi-i\phi)^{-1}$. Dann gelten $\Rstr_{\Pp}\tF\in\Dqa$ sowie
    	\begin{align*}
    		\phi = \frac{i}{2}\brklam{\Iq-\tF}\rklam{\psi-i\phi} \quad \text{und} \quad
    		\psi = \frac{1}{2}\brklam{\Iq+\tF}\rklam{\psi-i\phi}.
    	\end{align*}
    	
		\item [\rm{(b)}] Sei $F\in\Dqa$. Weiterhin seien
		\begin{align*}
			\MF := \gklam{z\in\Pm \;|\; \det F(\za) = 0},
		\end{align*}
		$\tF: \Pp\cup(-\infty,\alpha)\cup(\Pm\setminus{\cal M}_F) \rightarrow \Cqq$ definiert gemäß
		\begin{align*}
			\tF(z) := \begin{cases} F(z) & \text{falls } z\in\Pp \\
			\lim_{\omega\rightarrow z}F(\omega) & \text{falls } z\in(-\infty,\alpha) \\
			F^{-\ast}(\za) & \text{falls } z\in\Pm\setminus\MF \end{cases}
		\end{align*}
		sowie $\phi := i\big(\Iq-\tF\big)$ und $\psi := \Iq+\tF$.
		Dann gelten
		\begin{itemize}
			\item [\rm{(i)}] Es ist $\tF$ eine in $\C\setminus[\alpha,\infty)$ meromorphe \textit{q}$\times$\textit{q}-Matrixfunktion und es gilt $\phipsi\in\PtJqCa$. 
			\item [\rm{(ii)}] Es ist $\det(\psi-i\phi)$ nicht die Nullfunktion und es gilt
		\begin{align*}
			\tF = (\psi+i\phi)(\psi-i\phi)^{-1}.
		\end{align*}
		\end{itemize}
		
		\item [\rm{(c)}] Seien $\binom{\phi_1}{\psi_1},\binom{\phi_2}{\psi_2}\in\PtJqCa$. Dann ist $\det(\psi_j-i\phi_j)$ für alle $j\in\gklam{1,2}$ nicht die Nullfunktion. Weiterhin sei $\tF_j:=(\psi_j+i\phi_j)(\psi_j-i\phi_j)^{-1}$ für alle \linebreak $j\in\gklam{1,2}$. \dsfaa
		\begin{itemize} 
			\item [\rm{(iii)}] Es gilt $\Rstr_{\Pp}\tF_1=\Rstr_{\Pp}\tF_2$.
			\item [\rm{(iv)}] Die Paare $\binom{\phi_1}{\psi_1}$ und $\binom{\phi_2}{\psi_2}$ sind äquivalent.
		\end{itemize}
	\end{itemize}
\end{satz}

\bwanf Zu (a): Wegen \thref{spsa2} existiert eine diskrete Teilmenge $\D$ von $\C\setminus[\alpha,\infty)$ mit
\begin{itemize}
    \item [\rm{(v)}] $\phi$ und $\psi$ sind in $\C\setminus([\alpha,\infty)\cup\D)$ holomorph.
    \item [\rm{(vi)}] Es gilt $\rank \phipsiz = q$ für alle $z \in \C\setminus([\alpha,\infty)\cup\D)$.
    \item [\rm{(vii)}] Es gilt
    \begin{align*}
      	\phipsiz^{\ast}\bbrklam{\frac{-\tJq}{2\im z}}\phipsiz \geq \Oq
    \end{align*}
    für alle $z \in \C\setminus(\R\cup\D)$.
    \item [\rm{(viii)}] Es gilt 
    \begin{align*}
    	\phipsiz^{\ast}\rklam{-\Jq}\phipsiz \geq \Oq
    \end{align*}
   	für alle $z\in\Lam\setminus\D$.
\end{itemize}
Wegen {\rm (vi)} gilt
\begin{align}	\label{spsa1bw1}
	\phipsiz^{\ast}\phipsiz > \Oq
\end{align}
für alle $z\in\C\setminus([\alpha,\infty)\cup\D)$. Wegen {\rm (vii)} gilt weiterhin
\begin{align}	\label{spsa1bw6}
    \phipsiz^{\ast}\brklam{{-\tJq}}\phipsiz \geq \Oq
\end{align}
für alle $z\in\Pp\setminus\D$. Wegen \thref{spbsp1} gilt nun
\begin{align}	\label{spsa1bw2}
	[\psi(z)-i\phi(z)]^{\ast}[\psi(z)-i\phi(z)]
	&= \phi^{\ast}(z)\phi(z)+\psi^{\ast}(z)\psi(z)+i\eklam{\phi^{\ast}(z)\psi(z)-\psi^{\ast}(z)\phi(z)} \notag \\
	&= \phipsiz^{\ast}\phipsiz + 2\im \eklam{\psi^{\ast}(z)\phi(z)} \notag \\
	&= \phipsiz^{\ast}\phipsiz + \phipsiz^{\ast}\brklam{{-\tJq}}\phipsiz
\end{align}
für alle $z\in\C\setminus([\alpha,\infty)\cup\D)$. Hieraus folgt wegen \fref{spsa1bw6} und \fref{spsa1bw1} nun
\begin{align*}
	[\psi(z)-i\phi(z)]^{\ast}[\psi(z)-i\phi(z)]
	&\geq \phipsiz^{\ast}\phipsiz > \Oq
\end{align*}
für alle $z\in\Pp\setminus\D$. Somit gilt $\det(\psi(z)-i\phi(z))\neq0$ für alle $z\in\Pp\setminus\D$. Hieraus folgt dann, dass $\det(\psi-i\phi)$ nicht die Nullfunktion ist, und wegen {\rm (v)} gilt weiterhin, dass $F:=\Rstr_{\Pp}\tF$ in $\Pp\setminus\D$ holomorph ist. Insbesondere ist wegen der Wahl von $\phi$ und $\psi$ die Matrixfunktion $\tF$ meromorph in $\C\setminus[\alpha,\infty)$. Wegen \thref{spbsp1} gilt
\begin{align*}
	[\psi(z)+i\phi(z)]^{\ast}[\psi(z)+i\phi(z)] 
	&= \phi^{\ast}(z)\phi(z)+\psi^{\ast}(z)\psi(z)-i\eklam{\phi^{\ast}(z)\psi(z)-\psi^{\ast}(z)\phi(z)} \\
	&= \phipsiz^{\ast}\phipsiz - 2\im \eklam{\psi^{\ast}(z)\phi(z)} \\
	&= \phipsiz^{\ast}\phipsiz - \phipsiz^{\ast}\brklam{{-\tJq}}\phipsiz
\end{align*}
für alle $z\in\C\setminus([\alpha,\infty)\cup\D)$. Hieraus folgt wegen $\det(\psi(z)-i\phi(z))\neq0$ für alle $z\in\Pp\setminus\D$, \fref{spsa1bw2} und \fref{spsa1bw6} dann
\begin{align}	\label{spsa1bw10}
	\Iq-F^{\ast}(z)F(z)
	&= \Iq-\eklam{\psi(z)-i\phi(z)}^{-\ast}\eklam{\psi(z)+i\phi(z)}^{\ast}\eklam{\psi(z)+i\phi(z)}\eklam{\psi(z)-i\phi(z)}^{-1} \notag \\
	&= \eklam{\psi(z)-i\phi(z)}^{-\ast}\big(\eklam{\psi(z)-i\phi(z)}^{\ast}\eklam{\psi(z)-i\phi(z)} \notag \\
	&\quad -\eklam{\psi(z)+i\phi(z)}^{\ast}\eklam{\psi(z)+i\phi(z)}\big)\eklam{\psi(z)-i\phi(z)}^{-1} \notag \\
	&= 2 \eklam{\psi(z)-i\phi(z)}^{-\ast} \phipsiz^{\ast}\brklam{{-\tJq}}\phipsiz \eklam{\psi(z)-i\phi(z)}^{-1} \geq \Oq
\end{align}
für alle $z\in\Pp\setminus\D$. Hieraus folgt wegen Teil (a) von \thref{spdef6} im Fall $\D=\emptyset$ sogleich $F\in\SqPp$. Andererseits lässt sich $F$ wegen \thref{splm8} im Fall $\D\neq\emptyset$ auf eine Funktion aus $\SqPp$ fortsetzen. Wegen (viii) und \thref{spbsp4} gilt
\begin{align}	\label{spsa1bw4}
	\re\eklam{\psi^{\ast}(z)\phi(z)} \geq \Oq
\end{align}
für alle $z\in \Lam\setminus\D$. Weiterhin gilt
\begin{align}	\label{spsa1bw5}
	& \im\brklam{\eklam{\psi(z)-i\phi(z)}^{\ast}\eklam{\psi(z)+i\phi(z)}} \notag \\
	&= \im\brklam{\psi^{\ast}(z)\psi(z)-\phi^{\ast}(z)\phi(z)+i\eklam{\psi^{\ast}(z)\phi(z)+\phi^{\ast}(z)\psi(z)}} \notag \\
	&= \frac{1}{2i}\brklam{2i\eklam{\psi^{\ast}(z)\phi(z)+\phi^{\ast}(z)\psi(z)}}
	= 2\re\eklam{\psi^{\ast}(z)\phi(z)}
\end{align}
für alle $z\in\C\setminus([\alpha,\infty)\cup\D)$. Unter Beachtung von
\begin{align*}
	\im B^{\ast}AB = \frac{1}{2i}\rklam{B^{\ast}AB-B^{\ast}A^{\ast}B }
	= B^{\ast}\frac{1}{2i}\rklam{A-A^{\ast}}B = B^{\ast}(\im A)B
\end{align*}
für beliebige $A,B\in\Cqq$ folgt wegen $\det(\psi(z)-i\phi(z))\neq0$ für alle $z\in\Pp\setminus\D$, \fref{spsa1bw5} und \fref{spsa1bw4} nun
\begin{align}	\label{spsa1bw25}
	\im F(z)
	&= \im\brklam{\eklam{\psi(z)+i\phi(z)}\eklam{\psi(z)-i\phi(z)}^{-1}} \notag \\
	&= \im\brklam{\eklam{\psi(z)-i\phi(z)}^{-\ast}\eklam{\psi(z)-i\phi(z)}^{\ast}\eklam{\psi(z)+i\phi(z)}\eklam{\psi(z)-i\phi(z)}^{-1}} \notag \\
	&= \eklam{\psi(z)-i\phi(z)}^{-\ast} \eklam{\im\brklam{\eklam{\psi(z)-i\phi(z)}^{\ast}\eklam{\psi(z)+i\phi(z)}}} \eklam{\psi(z)-i\phi(z)}^{-1} \notag \\
	&= 2\eklam{\psi(z)-i\phi(z)}^{-\ast} \brklam{\re\eklam{\psi^{\ast}(z)\phi(z)}} \eklam{\psi(z)-i\phi(z)}^{-1}
	\geq \Oq
\end{align}
für alle $z\in\Pp\setminus\D$ mit $\re z \in (-\infty,\alpha)$. Unter Beachtung von
\begin{align*}
	\re\eklam{iA} = \frac{1}{2}\rklam{iA+\eklam{iA}^{\ast}}=\frac{1}{2i}\rklam{-A-\eklam{-A}^{\ast}} = \im\eklam{-A}
\end{align*}
für alle $A\in\Cqq$ ist wegen \fref{spsa1bw25} und \cite[Lemma 2.19]{DFK} oder \cite[Lemma 4.23]{Sch} jedes $z\in\Pp\cap\D$ mit $\re z \in (-\infty,\alpha)$ eine hebbare Singularität von $F$ und für die Fortsetzung von $F$ aus $\SqPp$ gilt 
\begin{align}	\label{spsa1bw7}
	\im F(z)\geq\Oq
\end{align}
für alle $z\in\Pp\setminus\D$ mit $\re z \in (-\infty,\alpha)$. Wegen {\rm (v)}, {\rm (vii)} und \thref{splm2} gilt
\begin{align}	\label{spsa1bw8}
    \binom{\phi(x)}{\psi(x)}^{\ast}\brklam{{-\tJq}}\binom{\phi(x)}{\psi(x)} = \Oq
\end{align}
für alle $x\in(-\infty,\alpha)\setminus\D$. Hieraus folgt wegen \fref{spsa1bw1} und \fref{spsa1bw2} dann
\begin{align*}
	[\psi(x)-i\phi(x)]^{\ast}[\psi(x)-i\phi(x)]
	&\geq \binom{\phi(x)}{\psi(x)}^{\ast}\binom{\phi(x)}{\psi(x)} > \Oq
\end{align*}
für alle $x\in(-\infty,\alpha)\setminus\D$. Somit gilt $\det(\psi(x)-i\phi(x))\neq0$ für alle $x\in(-\infty,\alpha)\setminus\D$. Hieraus folgt wegen {\rm (v)} dann, dass $\tF$ in $(-\infty,\alpha)\setminus\D$ holomorph ist. Wir zeigen nun im Fall $(-\infty,\alpha)\cap\D\neq\emptyset$ indirekt, dass $\tF$ in jedem Punkt $x\in(-\infty,\alpha)\cap\D$ hebbar ist. Angenommen, $x$ ist eine Polstelle der in $\C\setminus[\alpha,\infty)$ meromorphen Matrixfunktion $\tF$. Sei $\tF := \brklam{\tF_{ij}}_{i,j\in\Zeq}$. Dann ist $x$ eine Polstelle von $\tF_{ij}$ für alle $i,j\in\Zeq$.  Wegen \cite[Korollar 10.1.1]{Funk} ist dann $\lim_{z\rightarrow x}\babs{\tF_{ij}(z)}=\infty$ für alle $i,j\in\Zeq$, d.\,h. für eine Folge $(y_n)^{\infty}_{n=0}$ aus $\Pp\setminus\D$ existiert ein $N\in\No$ mit \begin{align}	\label{spsa1bw24}
	\babs{\tF_{ij}(y_N)} > \sqrt{q}
\end{align}
für alle $i,j\in\Zeq$. Wegen \fref{spsa1bw10} gilt $\tF^{\ast}(z)\tF(z) \leq \Iq$ für alle $z\in\Pp\setminus\D$. Unter Beachtung von $\enorm{A}^2 = \tr (A^{\ast}A)$ für alle $A\in\Cpq$ folgt dann
\begin{align*}
	\benorm{\tF(z)} \leq \enorm{\Iq} = \sqrt{q}
\end{align*}
für alle $z\in\Pp\setminus\D$. Somit gilt $\babs{\tF_{ij}(z)}\leq\sqrt{q}$ für alle $i,j\in\Zeq$ und $z\in\Pp\setminus\D$, im Widerspruch zu \fref{spsa1bw24}. Somit ist $x$ keine Polstelle, sondern eine hebbare Singularität, d.\,h. $\tF$ ist auf $(-\infty,\alpha)$ fortsetzbar. Es existiert nun
\begin{align}	\label{spsa1bw9}
	U_x := \lim_{z\rightarrow x}F(x)
\end{align}
für alle $x\in(-\infty,\alpha)$. Wegen \fref{spsa1bw10}, \fref{spsa1bw9} und \fref{spsa1bw8} gilt
\begin{align*}
	\Iq-U^{\ast}_xU_x = 2 \eklam{\psi(x)-i\phi(x)}^{-\ast} \binom{\phi(x)}{\psi(x)}^{\ast}\brklam{{-\tJq}}\binom{\phi(x)}{\psi(x)} \eklam{\psi(x)-i\phi(x)}^{-1} = \Oq
\end{align*}
für alle $x\in(-\infty,\alpha)\setminus\D$ und. Hieraus folgt wegen der Stetigkeit der Fortsetzung von $\tF$ in $(-\infty,\alpha)$ sogar $U^{\ast}_xU_x = \Iq$ für alle $x\in(-\infty,\alpha)$, also ist $U_x$ für alle \linebreak $x\in(-\infty,\alpha)$ eine unitäre Matrix. Hieraus folgt wegen $F\in\SqPp$, \fref{spsa1bw7}, \fref{spsa1bw9} und Teil (b) von \thref{spdef6} dann $F\in\Dqa$. Weiterhin gelten
\begin{align*}
	\frac{i}{2}\brklam{\Iq-\tF}(\psi-i\phi)
	&= \frac{i}{2}\eklam{\Iq-(\psi+i\phi)(\psi-i\phi)^{-1}}(\psi-i\phi) \\
	&= \frac{i}{2}\eklam{(\psi-i\phi)-(\psi+i\phi)}
	= \phi
\end{align*}
und
\begin{align*}
	\frac{1}{2}\brklam{\Iq+\tF}(\psi-i\phi)
	&= \frac{1}{2}\eklam{\Iq+(\psi+i\phi)(\psi-i\phi)^{-1}}(\psi-i\phi) \\
	&= \frac{1}{2}\eklam{(\psi-i\phi)+(\psi+i\phi)}
	= \psi.
\end{align*}

Zu (b): Wir zeigen zunächst {\rm (i)}. Wegen der Teile (a) und (b) von \thref{spdef6} ist $F$ in $\Pp$ holomorph und es gilt
\begin{itemize}
	\item [\rm{(ix)}] Für alle $x\in(-\infty,\alpha)$ existiert $U_x := \lim_{z\rightarrow x} F(z)$ und es ist $U_x$ unitär.
\end{itemize}
Angenommen, $\det F$ ist die Nullfunktion auf $\Pp$. Wegen der Stetigkeit von $F$ in $\Pp$ und {\rm (ix)} gilt dann
\begin{align*}
	\det U_x = \det \eklam{\lim_{z\rightarrow x} F(z)} = \lim_{z\rightarrow x} \eklam{\det F(z)} = 0
\end{align*}
für alle $x\in(-\infty,\alpha)$, im Widerspruch zu {\rm (ix)}, dass $U_x$ für alle $x\in(-\infty,\alpha)$ unitär ist. Somit ist $\det F$ nicht die Nullfunktion sowie unter Beachtung des Identitätssatzes für meromorphe Funktionen (vergleiche z.\,B. im skalaren Fall \cite[Satz 10.3.2]{Funk}; im matriziellen Fall betrachtet man die einzelnen Einträge der Matrixfunktion) ist $\MF$ eine diskrete Teilmenge von $\Pm$. Unter Beachtung von $\Pm=\gklam{\za \;|\; z\in\Pp}$ ist $\Rstr_{\Pm\setminus\MF}\tF$ wohldefiniert und in $\Pm\setminus\MF$ holomorph. Sei
\begin{align*}
	\widecheck{\M}_F := \gklam{\za \;|\; z\in\MF} = \gklam{z\in\Pp \;|\; \det F(z)=0}
\end{align*}
Wegen {\rm (ix)} gilt dann
\begin{align*}
	U_x = U^{-\ast}_x 
	= \Beklam{\lim_{\substack{z\rightarrow x \\ z\in\Pp\setminus\widecheck{\M}_F}}F(z)}^{-\ast} 
	= \lim_{\substack{z\rightarrow x \\ z\in\Pp\setminus\widecheck{\M}_F}}F^{-\ast}(z)
	= \lim_{\substack{z\rightarrow x \\ z\in\Pm\setminus\MF}}F^{-\ast}(\za)
	= \lim_{\substack{z\rightarrow x \\ z\in\Pm\setminus\MF}}\tF(z)
\end{align*}
für alle $x\in(-\infty,\alpha)$. Hieraus folgt wegen {\rm (ix)} und der Tatsache, dass $F=\Rstr_{\Pp}\tF$ bzw. $\Rstr_{\Pm\setminus\MF}\tF$ in $\Pp$ bzw. $\Pm\setminus\MF$ holomorph ist, dann, dass $\tF$ in
\begin{align*}
	\Pp\cup(-\infty,\alpha)\cup(\Pm\setminus\MF) = \C\setminus([\alpha,\infty)\cup\MF)
\end{align*}
holomorph und somit in $\C\setminus[\alpha,\infty)$ meromorph ist. Hieraus folgt dann, dass $\phi$ und $\psi$ in $\C\setminus([\alpha,\infty)\cup\MF)$ holomorph sind. Unter Beachtung von $\tF^{\ast}(z)\tF(z)\geq\Oq$ für alle $z\in\C\setminus([\alpha,\infty)\cup\MF)$ gilt
\begin{align*}
	&\ \phipsiz^{\ast}\phipsiz = \phi^{\ast}(z)\phi(z)+\psi^{\ast}(z)\psi(z) \\
	&= \brklam{i\beklam{\Iq-\tF(z)}}^{\ast}\brklam{i\beklam{\Iq-\tF(z)}}+\beklam{\Iq+\tF(z)}^{\ast}\beklam{\Iq+\tF(z)} \\
	&= 2\beklam{\Iq+\tF^{\ast}(z)\tF(z)}
	\geq 2\Iq > \Oq
\end{align*}
für alle $z\in\C\setminus([\alpha,\infty)\cup\MF)$. Hieraus folgt dann
\begin{align}	\label{spsa1bw11}
	\rank\phipsiz = q
\end{align}
für alle $z\in\C\setminus([\alpha,\infty)\cup\MF)$. Wegen \thref{spbsp1} gilt
\begin{align}	\label{spsa1bw12}
	&\ \phipsiz^{\ast}\brklam{{-\tJq}}\phipsiz
	= 2\im\eklam{\psi^{\ast}(z)\phi(z)}
	= \frac{1}{i}\eklam{\psi^{\ast}(z)\phi(z)-\phi^{\ast}(z)\psi(z)} \notag \\
	&= \frac{1}{i}\Beklam{\brklam{\Iq+\tF(z)}^{\ast}\brklam{i\beklam{\Iq-\tF(z)}}-\brklam{i\beklam{\Iq-\tF(z)}}^{\ast}\brklam{\Iq+\tF(z)}} \notag \\
	&= 2\beklam{\Iq-\tF^{\ast}(z)\tF(z)}
\end{align}
für alle $z\in\C\setminus([\alpha,\infty)\cup\MF)$. Hieraus folgt wegen $F=\Rstr_{\Pp}\tF\in\SqPp$ (vergleiche Teil (a) von \thref{spdef6}) nun
\begin{align*}
	\phipsiz^{\ast}\brklam{{-\tJq}}\phipsiz \geq \Oq
\end{align*}
für alle $z\in\Pp$. Hieraus folgt dann
\begin{align}	\label{spsa1bw13}
	\phipsiz^{\ast}\bbrklam{\frac{-\tJq}{2\im z}}\phipsiz \geq \Oq
\end{align}
für alle $z\in\Pp$. Wegen $F\in\SqPp$ (vergleiche Teil (a) von \thref{spdef6}) gilt weiterhin
\begin{align*}
	F^{\ast}(\za)F(\za) \leq \Iq
\end{align*}
für alle $z\in\Pm$. Hieraus folgt dann
\begin{align*}
	\tF^{\ast}(z)\tF(z) = F^{-1}(\za)F^{-\ast}(\za) = \beklam{F^{\ast}(\za)F(\za)}^{-1}
	\geq \Iq^{-1}=\Iq
\end{align*}
für alle $z\in\Pm\setminus\MF$. Hieraus folgt wegen \fref{spsa1bw12} nun
\begin{align*}
	-\phipsiz^{\ast}\brklam{{-\tJq}}\phipsiz \geq \Oq
\end{align*}
für alle $z\in\Pm\setminus\MF$. Hieraus folgt dann
\begin{align}	\label{spsa1bw14}
	\phipsiz^{\ast}\bbrklam{\frac{-\tJq}{2\im z}}\phipsiz \geq \Oq
\end{align}
für alle $z\in\Pm\setminus\MF$, also wegen \fref{spsa1bw13} sogar für alle $z\in\C\setminus(\R\cup\MF)$. Es gilt
\begin{align}	\label{spsa1bw23}
	\re\eklam{\psi^{\ast}(z)\phi(z)} 
	&= \frac{1}{2}\eklam{\psi^{\ast}(z)\phi(z)+\phi^{\ast}(z)\psi(z)} \notag \\
	&= \frac{1}{2}\Beklam{\brklam{\Iq+\tF(z)}^{\ast}\brklam{i\beklam{\Iq-\tF(z)}}+\brklam{i\beklam{\Iq-\tF(z)}}^{\ast}\brklam{\Iq+\tF(z)}} \notag \\
	&= i\beklam{\tF^{\ast}(z)-\tF(z)} = 2\im\tF(z)
\end{align}
für alle $z\in\C\setminus([\alpha,\infty)\cup\MF)$. Hieraus folgen wegen $F=\Rstr_{\Pp}\tF\in\Dqa$ (vergleiche Teil (b) von \thref{spdef6}) dann
\begin{align}	\label{spsa1bw17}
	\re\eklam{\psi^{\ast}(z)\phi(z)} = 2\im F(z) \geq \Oq
\end{align}
für alle $z\in\Pp$ mit $\re z\in(-\infty,\alpha)$ und unter Beachtung von
\begin{align*}
	\im A^{-\ast} &= \frac{1}{2i}\rklam{A^{-\ast}-A^{-1}}
	= \frac{1}{2i}\rklam{A^{-\ast}AA^{-1}-A^{-\ast}A^{\ast}A^{-1}} \\
	&= A^{-\ast}\frac{1}{2i}\rklam{A-A^{\ast}} A^{-1} = A^{-\ast}(\im A)A^{-1}
\end{align*} 
für alle reguläre $A\in\Cqq$ weiterhin
\begin{align}	\label{spsa1bw18}
	\re\eklam{\psi^{\ast}(z)\phi(z)} = 2F^{-\ast}(\za)\eklam{\im F(\za)}F^{-1}(\za) \geq \Oq
\end{align}
für alle $z\in\Pm\setminus\MF$ mit $\re z\in(-\infty,\alpha)$. Wegen der Stetigkeit von $F$ in $\Pp$, \fref{spsa1bw23} und \fref{spsa1bw17} gilt
\begin{align*}
	\re\eklam{\psi^{\ast}(x)\phi(x)} = 2\im\Beklam{\lim_{z\rightarrow x}F(z)} \geq \Oq
\end{align*}
für alle $x\in(-\infty,\alpha)$. Hieraus folgt wegen \thref{spbsp4}, \fref{spsa1bw17} und \fref{spsa1bw18} dann
\begin{align*}
	\phipsiz^{\ast}\rklam{-\Jq}\phipsiz = 2\re\eklam{\psi^{\ast}(z)\phi(z)} \geq \Oq
\end{align*}
für alle $z\in\Lam\setminus\MF$. Hieraus folgt unter Beachtung, dass $\phi$ und $\psi$ in $\C\setminus([\alpha,\infty)\cup\MF)$ holomorph sind, \fref{spsa1bw11} und \fref{spsa1bw14} wegen \thref{spsa2} dann $\phipsi\in\PtJqCa$. 

Es bleibt noch {\rm (ii)} zu zeigen. 
Es gelten
\begin{align}	\label{spsa1bw21}
	\psi-i\phi = \Iq+\tF-i\beklam{i\brklam{\Iq-\tF}} = 2\Iq
\end{align}
und
\begin{align}	\label{spsa1bw22}
	\psi+i\phi = \Iq+\tF+i\beklam{i\brklam{\Iq-\tF}} = 2\tF.
\end{align}
Wegen \fref{spsa1bw21} ist $\det(\psi-i\phi)$ nicht die Nullfunktion. Hieraus folgt wegen \fref{spsa1bw21} und \fref{spsa1bw22} weiterhin
\begin{align*}
	(\psi+i\phi)(\psi-i\phi)^{-1} = \brklam{2\tF}\brklam{2\Iq}^{-1} = \tF.
\end{align*}

Zu (c): Wegen (a) ist $\det(\psi_j-i\phi_j)$ für alle $j\in\gklam{1,2}$ nicht die Nullfunktion.

Sei zunächst (iii) erfüllt. Wegen des Identitätssatzes für meromorphe Funktionen (vergleiche z.\,B. im skalaren Fall \cite[Satz 10.3.2]{Funk}; im matriziellen Fall betrachtet man die einzelnen Einträge der Matrixfunktion) gilt dann $\tF_1 = \tF_2$. Hieraus folgt wegen (a) dann
\begin{align*}
	\phi_j = \frac{i}{2}\brklam{\Iq-\tF_1}\rklam{\psi_j-i\phi_j} \quad \text{und} \quad
    \psi_j = \frac{1}{2}\brklam{\Iq+\tF_1}\rklam{\psi_j-i\phi_j}
\end{align*}
für alle $j\in\gklam{1,2}$. Hieraus folgt nun
\begin{align*}
	\binom{\phi_2}{\psi_2} &= \begin{pmatrix} \frac{i}{2}\brklam{\Iq-\tF_1} \\ \frac{1}{2}\brklam{\Iq+\tF_1} \end{pmatrix} \rklam{\psi_2-i\phi_2} \\
	&= \begin{pmatrix} \frac{i}{2}\brklam{\Iq-\tF_1} \\ \frac{1}{2}\brklam{\Iq+\tF_1} \end{pmatrix} \rklam{\psi_1-i\phi_1}\rklam{\psi_1-i\phi_1}^{-1}\rklam{\psi_2-i\phi_2} \\
	&= \binom{\phi_1}{\psi_1}\rklam{\psi_1-i\phi_1}^{-1}\rklam{\psi_2-i\phi_2}.
\end{align*}
Sei $g:=\rklam{\psi_1-i\phi_1}^{-1}\rklam{\psi_2-i\phi_2}$. Dann existiert eine Teilmenge $\D$ von $\C\setminus[\alpha,\infty)$ (als Vereinigung endlich vieler diskreter Teilmengen von $\C\setminus[\alpha,\infty)$), so dass $\phi_1, \psi_1, \phi_2, \psi_2$ und $g$ in $\C\setminus([\alpha,\infty)\cup\D$ holomorph sind und $\det g(z) \neq 0$ sowie 
\begin{align*}
	\binom{\phi_2(z)}{\psi_2(z)} = \binom{\phi_1(z)}{\psi_1(z)}g(z)
\end{align*}
für alle $z \in \C\setminus([\alpha,\infty)\cup\D$ erfüllt sind. Wegen \thref{spdef4b} gilt dann $\sklam{\binom{\phi_1}{\psi_1}} = \sklam{\binom{\phi_2}{\psi_2}}$.

Sei nun (iv) erfüllt. Wegen \thref{spdef4b} existieren eine in $\C\setminus[\alpha,\infty)$ meromorphe \textit{q}$\times$\textit{q}"=Matrixfunktion $g$ und eine diskrete Teilmenge $\D$ von $\C\setminus[\alpha,\infty)$, so dass $\phi_1, \psi_1, \phi_2, \psi_2$ und $g$ in $\C\setminus([\alpha,\infty)\cup\D)$ holomorph sind sowie $\det g(z) \neq 0$ und
\begin{align*}
  \binom{\phi_2(z)}{\psi_2(z)} = \binom{\phi_1(z)}{\psi_1(z)}g(z)
\end{align*}
für alle $z \in \C\setminus([\alpha,\infty)\cup\D)$ erfüllt sind. Hieraus folgt dann
\begin{align*}
	\tF_2 &= (\psi_2+i\phi_2)(\psi_2-i\phi_2)^{-1} \\
	&= (\psi_1g+i\phi_1g)(\psi_1g-i\phi_1g)^{-1} \\
	&= (\psi_1+i\phi_1)gg^{-1}(\psi_1-i\phi_1)^{-1} \\
	&= (\psi_1+i\phi_1)(\psi_1-i\phi_1)^{-1} = \tF_1,
\end{align*}
also insbesondere $\Rstr_{\Pp}\tF_2 = \Rstr_{\Pp}\tF_1$. \bwend

Mithilfe von \thref{spsa1} erkennen wir, dass eine Bijektion zwischen der Teilklasse $\Dqa$ von \textit{q}$\times$\textit{q}-Schur-Funktionen auf $\Pp$ und der Menge der Äquivalenzklassen von \textit{q}$\times$\textit{q}-Stieltjes-Paaren in $\C\setminus[\alpha,\infty)$ besteht. 
Mithilfe von Teil (a) von \thref{spbm3} und Teil (b) von \thref{spbm4} lässt sich weiterhin zeigen, dass eine Bijektion zwischen der Teilklasse $\Eqa$ von \textit{q}$\times$\textit{q}-Schur-Funktionen auf $\Pm$ und der Menge der Äquivalenzklassen von \textit{q}$\times$\textit{q}-Stieltjes-Paaren in $\C\setminus(-\infty,\alpha]$ besteht. 

Die Ausführungen von \thref{spsa1} bringen uns auf folgende Bemerkung.

\begin{bem}	\thlabel{spbm17}
	Sei $\alpha\in\R$. Bezeichne ${\cal I}$ die auf $\C\setminus[\alpha,\infty)$ konstante Matrixfunktion mit dem Wert $\Iq$. \dgfa
	\begin{itemize}
		\item [\rm{(a)}] Sei $\phipsi\in\PtJqCa$. Dann existiert ein zu $\phipsi$ äquivalentes Paar \linebreak $\tphipsi\in\PtJqCa$ mit $\tpsi-i\tphi=2{\cal I}$.
		\item [\rm{(b)}] Sei $\phipsi\in\PtJqCa$ derart, dass $\psi-i\phi=2{\cal I}$ erfüllt ist. Dann existiert eine diskrete Teilmenge $\D$ von $\C\setminus[\alpha,\infty)$ mit 
		\begin{align*}
			\enorm{\binom{\phi(x)}{\psi(x)}} = 2\sqrt{q}
		\end{align*}
		für alle $x\in(-\infty,\alpha)\setminus\D$.
	\end{itemize}
\end{bem}

\bwanf Zu (a): Dies folgt sogleich aus \thref{spsa1} (vergleiche hierzu insbesondere Formel \fref{spsa1bw21} im zugehörigen Beweis).

Zu (b): Sei $\tF := (\psi+i\phi)(\psi-i\phi)^{-1}$. Wegen Teil (a) von \thref{spsa1} gilt dann \linebreak $\Rstr_{\Pp}\tF\in\Dqa$. Wegen Teil (a) von \thref{spdef4} existiert eine diskrete Teilmenge $\D$ von $\C\setminus[\alpha,\infty)$, so dass $\phi$ und $\psi$ in $\C\setminus([\alpha,\infty)\cup\D)$ holomorph sind. Somit ist $\tF$ in $\C\setminus([\alpha,\infty)\cup\D)$ stetig und wegen Teil (b) von \thref{spdef6} ist dann $\tF(x)$ für alle $x\in(-\infty,\alpha)\setminus\D$ unitär. Hieraus folgt wegen Teil (a) von \thref{spsa1} dann
\begin{align*}
	\binom{\phi(x)}{\psi(x)}^{\ast}\binom{\phi(x)}{\psi(x)} = \begin{pmatrix} i\beklam{\Iq-\tF(x)} \\ \Iq+\tF(x) \end{pmatrix}^{\ast} \begin{pmatrix} i\beklam{\Iq-\tF(x)} \\ \Iq+\tF(x) \end{pmatrix} = 2\beklam{\Iq+\tF^{\ast}(x)\tF(x)} = 4\Iq
\end{align*}
für alle $x\in(-\infty,\alpha)\setminus\D$. Hieraus folgt unter Beachtung von $\enorm{A}^2 = \tr (A^{\ast}A)$ für alle $A\in\Cpq$ nun
\begin{align*}
	\enorm{\binom{\phi(x)}{\psi(x)}} = \enorm{2\Iq} = 2\sqrt{q}
\end{align*}
für alle $x\in(-\infty,\alpha)\setminus\D$. \bwend

%% file: begriffe.tex
\renewcommand{\refname}{Begriffsverzeichnis}
\fancyhead[RO,LE]{Begriffsverzeichnis}
{\Large\sffamily\bfseries  Begriffsverzeichnis}
\addcontentsline{toc}{section}{Begriffsverzeichnis}

\begin{tabular} [t]{p{11.5cm}l}	\\

durch rechts- bzw. linksseitige $\alpha$-Verschiebung generierte Folge	&	\thref{aspdef1}	\\
Favard-Paar	&	\thref{fpdef1}	\\
Folge von linksseitigen $2$\textit{q}$\times2$\textit{q}-$\alpha$-Dyukarev-Matrixpolynomen	&	\thref{drldef3}	\\
Folge von rechtsseitigen $2$\textit{q}$\times2$\textit{q}-$\alpha$-Dyukarev-Matrixpolynomen	&	\thref{drdef2}	\\
Hankel-nichtnegativ bzw. -positiv Definitheit	&	\thref{asmdef1}	\\
$J$-innere Funktion aus $\bJLam$	&	\thref{spdef8}	\\
$J$-innere Funktion aus $\bJPp$	&	\thref{spdef3}	\\
$J^{(2)}$-$J^{(1)}$-innere Funktion aus $\bJJPp$	&	\thref{spdef9}	\\
$J$-kontraktiv bzw. $J$-expansiv bzw. $J$-unitär	&	\thref{spdef2}	\\
$J^{(2)}$-$J^{(1)}$-kontraktiv bzw. $J^{(2)}$-$J^{(1)}$-expansiv bzw. $J^{(2)}$-$J^{(1)}$-unitär	&	\thref{spdef7}	\\
$J$-Potapov-Funktion bezüglich $\Lam$	&	\thref{spdef8}	\\
$J$-Potapov-Funktion bezüglich $\Pp$	&	\thref{spdef3}	\\
$J^{(2)}$-$J^{(1)}$-Potapov-Funktion bezüglich $\Pp$	&	\thref{spdef9}	\\
kanonische Hankel-Parametrisierung	&	\thref{mpdef1}	\\
linke Schur-Komplement	&	\thref{asmbz1}	\\
linke System von Matrixpolynomen zweiter Art	&	\thref{mpdef3}	\\
linksseitig $\alpha$-Stieltjes-nichtnegativ bzw. -positiv Definitheit	&	\thref{asmdef3}	\\
linksseitige $2$\textit{q}$\times2$\textit{q}-$\alpha$-Dyukarev-Matrixpolynom	&	\thref{drldef3}	\\
linksseitige $\alpha$-Dyukarev-Quadrupel	&	\thref{drldef2}	\\
linksseitige $\alpha$"=Dyukarev"=Stieltjes"=Parametrisierung	&	\thref{adpldef1}	\\
linksseitige $\alpha$-Stieltjes-Quadrupel	&	\thref{sqldef1}	\\
monische links-orthogonale System von Matrixpolynomen	&	\thref{mpdef2}	\\
n-te Block-Hankel-Matrix	&	\thref{asmbz1}	\\
Potapovsche Fundamentalmatrizen	im linksseitigen Fall & \thref{drldef1}	\\
Potapovsche Fundamentalmatrizen	im rechtsseitigen Fall & \thref{drdef3} \\
\textit{p}$\times$\textit{q}-Schur-Funktion	&	\thref{spdef6}	\\
\textit{p}$\times$\textit{p}-Signaturmatrix	&	\thref{spdef1}	\\
\textit{q}$\times$\textit{q}-Stieltjes-Paar in $\C\setminus[\alpha,\infty)$	&	\thref{spdef4}	\\
\textit{q}$\times$\textit{q}-Stieltjes-Paar in $\C\setminus(-\infty,\alpha]$	&	\thref{spdef5}	\\
rechtsseitig $\alpha$-Stieltjes-nichtnegativ bzw. -positiv Definitheit	&	\thref{asmdef2}	\\
rechtsseitige $2$\textit{q}$\times2$\textit{q}-$\alpha$-Dyukarev-Matrixpolynom	&	\thref{drdef2}	\\
rechtsseitige $\alpha$-Dyukarev-Quadrupel	&	\thref{drdef1}	\\
rechtsseitige $\alpha$"=Dyukarev"=Stieltjes"=Parametrisierung	&	\thref{adpdef1}	\\
rechtsseitige bzw. linksseitige $\alpha$-Stieltjes-Parametrisierung	&	\thref{aspdef2}	\\
rechtsseitige $\alpha$-Stieltjes-Quadrupel	&	\thref{sqdef1}	\\
Resolventenmatrix des Momentenproblems	&	\thref{drdef0}	\\
Stieltjes-Maß	&	\thref{asmdef5}	\\
Stieltjes-Transformierte	&	\thref{asmdef5}	\\
untere n-te Blockdreiecksmatrix	&	\thref{mpbz3}	\\
unteres bzw. oberes Extremalelement von $\Sqasmu$	&	\thref{drdef4} \\
unteres bzw. oberes Extremalelement von $\Sqmasmu$	&	\thref{drldef4} \\

\end{tabular}

\newpage
\begin{tabular} [t]{p{11.5cm}l}	\\

Weylsche Intervall für eine linksseitig $\alpha$-Stieltjes-positiv definite Folge	& \thref{drldef5}	\\
Weylsche Intervall für eine rechtsseitig $\alpha$-Stieltjes-positiv definite Folge	& \thref{drdef5}	\\
zu einer Matrizenfolge gehörige Matrixpolynom	&	\thref{mpbz3}	\\

\end{tabular}

%% file: symbole.tex
\newpage
\renewcommand{\refname}{Symbolverzeichnis}
\fancyhead[RO,LE]{Symbolverzeichnis}
{\Large\bf  Symbolverzeichnis} \\
\addcontentsline{toc}{section}{Symbolverzeichnis}

Symbole ohne \anf{$\sklam{s}$} als oberen Index oder ohne \anf{s} im unteren Index wurden an der selben Stelle eingeführt, wie Symbole mit jenem Index. Sie werden verwendet, wenn klar ist, von welcher Folge von Matrizen die Rede ist.

\begin{minipage}[t]{0.49\textwidth}
\begin{tabular}[t]{l l}

$\enorm{A}$						&	\sref{norm}			\\
$\snorm{A}$						&	\sref{norm}			\\
$[A,B]$							&	\sref{MI}			\\
$\fklamo{\kappa}$				&	\sref{fklam}		\\

$A^{\ast}$						&	\sref{A}			\\
$A^{+}$							&	\sref{A}			\\
$A^{-1}$						&	\sref{A2}			\\
$A^{-\ast}$						&	\sref{A2}			\\
$\sqrt{A}$						&	\sref{MI}			\\

$A^{\sklam{s}}_n$				&	\thref{fpdef1}		\\
$\Asarn$						&	\thref{drdef1}		\\
$\Asaln$						&	\thref{drldef2}		\\

$B^{\sklam{s}}_n$				&	\thref{fpdef1}		\\
$\Bsarn$						&	\thref{drdef1}		\\
$\Bsaln$						&	\thref{drldef2}		\\

$\bJLam$						&	\thref{spdef8}		\\
$\tbJLam$						&	\thref{spdef8}		\\
$\bJPp$							&	\thref{spdef3}		\\
$\tbJPp$						&	\thref{spdef3}		\\
$\bJJPp$						&	\thref{spdef9}		\\
$\tbJJPp$						&	\thref{spdef9}		\\

$C_{n}$							&	\thref{mpdef1}		\\
$\Csarn$						&	\thref{drdef1}		\\
$\Csaln$						&	\thref{drldef2}		\\

$\C$							&	\sref{Mengen}		\\
$\Cpq$							&	\sref{Cpq}			\\
$\Cq$							&	\sref{Cpq}			\\
$\Cqq_H$						&	\sref{CH}			\\
$\Cqq_>$						&	\sref{CH}			\\
$\Cqq_{\geq}$					&	\sref{CH}			\\

$\Lap$							&	\sref{La}			\\
$\Lam$							&	\sref{La}			\\

$\deg P$						&	\thref{mpbz2}		\\
$\det A$						&	\sref{A}			\\
$\det f$						&	\sref{f}			\\
$\diag\rklam{A_0, \ldots, A_n}$	&	\sref{diag}			\\

$D_{n}$							&	\thref{mpdef1}		\\
$\Dsarn$						&	\thref{drdef1}		\\
$\Dsaln$						&	\thref{drldef2}		\\

\end{tabular}
\end{minipage}
\hfill
\begin{minipage}[t]{0.49\textwidth}
\begin{tabular}[t]{l l}

$\Dqa$							&	\thref{spdef6}		\\

$\delta_{j,k}$					&	\sref{djk}			\\

$\En$							&	\thref{adpbz1}		\\

$\Eqa$							&	\thref{spdef6}		\\

$f^{\ast}$						&	\sref{f}			\\
$f^{-1}$						&	\sref{f}			\\
$f^{-\ast}$						&	\sref{f}			\\

$\Ffns$							&	\thref{drdef3}		\\
$\Ffarns$						&	\thref{drdef3}		\\
$\Ffalns$						&	\thref{drldef1}		\\
$\dFfns$						&	\thref{drdef3}		\\
$\dFfarns$						&	\thref{drdef3}		\\
$\dFfalns$						&	\thref{drldef1}		\\

$\Hsn$							&	\thref{asmbz1}		\\	
$\Hsarn$						&	\thref{aspdef1}		\\	
$\Hsaln$						&	\thref{aspdef1}		\\	

$\dHsn$							&	\thref{asmbz1}		\\	
$\dHsarn$						&	\thref{aspdef1}		\\	
$\dHsaln$						&	\thref{aspdef1}		\\	

$\Hqn$							&	\thref{asmdef1}		\\	
$\Heqm$							&	\thref{asmdef1}		\\
$\Hpqk$							&	\thref{asmdef1}		\\

$\im A$							&	\sref{reim}			\\

$\Iq$							&	\sref{Cpq}			\\

$\jqq$							&	\thref{spbsp5}		\\
$\Jq$							&	\thref{spbsp4}		\\
$\tJq$							&	\thref{spbsp1}		\\
 
$\Ksn$							&	\thref{asmbz1}		\\
$\wKsn$							&	\thref{mabz2}		\\

$\Kqma$							&	\thref{asmdef2}		\\
$\Keqma$						&	\thref{asmdef2}		\\
$\Kpqma$						&	\thref{asmdef2}		\\ 

$L_n$							&	\thref{drbz1}		\\
$\dL_n$							&	\thref{drbz1}		\\

$\Larn^{\sklam{s}}$				&	\thref{adpdef1}		\\

${\cal L}^1(\Omega,\A,\mu;\C)$	&	\sref{L1}			\\

$\Lqma$							&	\thref{asmdef3}		\\
$\Leqma$						&	\thref{asmdef3}		\\

\end{tabular}
\end{minipage}

\newpage

\begin{minipage}[t]{0.49\textwidth}
\begin{tabular}[t]{l l}

$\Lpqma$						&	\thref{asmdef3}		\\

$\Lambda^{\sklam{s}}_{n}$		&	\thref{mpbz1}		\\

$\Mskg$							&	\sref{Ms}			\\
$\Msmu$							&	\sref{Ms}			\\

$M^{\sklam{s}}_{\arn}$			&	\thref{drbz5}		\\
$M^{\sklam{s}}_{\aln}$			&	\thref{drlbz2}		\\
$\widetilde{M}^{\sklam{s}}_{\arn}$	&	\thref{drbz5}		\\
$\widetilde{M}^{\sklam{s}}_{\aln}$	&	\thref{drlbz2}		\\

$\Marn^{\sklam{s}}$				&	\thref{adpdef1}		\\
$\Maln^{\sklam{s}}$				&	\thref{adpldef1}	\\

$\Mqo$							&	\sref{Mq}			\\
$\Mqko$							&	\sref{Mqko}			\\
$\Mqoa$							&	\sref{Mq}			\\

$\widecheck{\mu}$				&	\thref{hmbz1}		\\

$\Opq$							&	\sref{Cpq}			\\

$\widecheck{\Omega}$			&	\thref{hmbz1}		\\

$\N$							&	\sref{Mengen}		\\
$\No$							&	\sref{Mengen}		\\
$\Na$							&	\sref{Mengen}		\\
$\Noa$							&	\sref{Mengen}		\\

$P^{[j]}$						&	\thref{mpbz2}		\\
$\Ps$							&	\thref{mpbz3}		\\
$P_n$							&	\thref{mpdef2}		\\
$\Ps_n$							&	\thref{mpdef3}		\\

$P_{n,s}$						&	\thref{sqdef1}		\\
$\Ps_{n,s}$						&	\thref{sqdef1}		\\
$P_{\arn,s}$					&	\thref{sqdef1}		\\
$P_{\aln,s}$					&	\thref{sqldef1}		\\
$\widehat{P}_{\arn,s}$			&	\thref{sqdef1}		\\
$\widehat{P}_{\aln,s}$			&	\thref{sqldef1}		\\
$\dPs_{\arn,s}$					&	\thref{sqdef1}		\\
$\dPs_{\aln,s}$					&	\thref{sqldef1}		\\

$\PtJqCa$						&	\thref{spdef4}		\\
$\dPtJqCa$						&	\thref{spdef4}		\\
$\PtJqCma$						&	\thref{spdef5}		\\
$\dPtJqCma$						&	\thref{spdef5}		\\

$\Pp$							& 	\sref{P}			\\
$\Pm$							& 	\sref{P}			\\ 
 
$\Yarn^{\sklam{s}}$				&	\thref{drbz7}		\\
$\Yaln^{\sklam{s}}$				&	\thref{drlbz2}		\\
$\tYarn^{\sklam{s}}$			&	\thref{drbz7}		\\
$\tYaln^{\sklam{s}}$			&	\thref{drlbz2}		\\

$Q^{\sklam{s}}_{\ar m}$			&	\thref{aspdef2}		\\

\end{tabular}
\end{minipage}
\hfill
\begin{minipage}[t]{0.55\textwidth}
\begin{tabular}[t]{l l}

$Q^{\sklam{s}}_{\al m}$			&	\thref{aspdef2}		\\

$\rank A$						&	\sref{A}			\\

$\re A$							&	\sref{reim}			\\

$R_n$							&	\thref{drbz1}		\\

$\Rstr_Z f$						&	\sref{Rstr}			\\

$\Rqab$							&	\thref{mabz1}		\\
$\Rqabsjk$						&	\sref{Sss}			\\

$\R$							&	\sref{Mengen}		\\

$s^{(\mu)}_{j}$					&	\sref{sjs}			\\
$s_{\arj}$						&	\thref{aspdef1}		\\
$s_{\alj}$						&	\thref{aspdef1}		\\

$\Sarmmaxs$						&	\thref{drdef4}		\\
$\Salmmaxs$						&	\thref{drldef4}		\\
$\Sarmmins$						&	\thref{drdef4}		\\
$\Salmmins$						&	\thref{drldef4}		\\

$\Saskg$						&	\sref{Ss}			\\
$\Sasmu$  						&	\sref{Ss}			\\
$\Sabskg$						&	\sref{Sss}			\\
$\Smaskg$						&	\sref{Ssl}			\\
$\Smasmu$  						&	\sref{Ssl}			\\

$S_n$							&	\thref{mpbz3}		\\

$\Sarn$							&	\thref{sqbz1}		\\
$\dSarn$						&	\thref{sqbz1}		\\
$\Saln$							&	\thref{sqlbz1}		\\
$\dSaln$						&	\thref{sqlbz1}		\\

$\SpqG$							&	\thref{spdef6}		\\

$\Sqa$							&	\thref{asmbz3}		\\
$\Sqma$							&	\thref{asmbz4}		\\
$\Soqa$							&	\thref{asmbz3}		\\
$\Sqaskg$						&	\sref{Ss}			\\
$\Sqasmu$						&	\sref{Ss}			\\ 
$\Soqma$						&	\thref{asmbz4}		\\
$\Sqmaskg$						&	\sref{Ssl}			\\
$\Sqmasmu$						&	\sref{Ssl}			\\ 

$\Sarm$							&	\thref{drth3}		\\
$\Salm$							&	\thref{drlth3}		\\

$T_{a,b}$						&	\thref{hmbz1}		\\
 
$T_n$							&	\thref{drbz1}		\\

$\Tarn^{\sklam{s}}$				&	\thref{drbz6}		\\
$\Taln^{\sklam{s}}$				&	\thref{drlbz2}		\\
$\tTarn^{\sklam{s}}$			&	\thref{drbz6}		\\
$\tTaln^{\sklam{s}}$			&	\thref{drlbz2}		\\

$\usn$							&	\thref{drbz3}		\\
$\usarn$						&	\thref{drbz3}		\\
 
\end{tabular}
\end{minipage}

\newpage

\begin{minipage}[t]{0.49\textwidth}
\begin{tabular}[t]{l l}

$\usaln$						&	\thref{drlbz1}		\\

$\Usarm$						&	\thref{drdef2}		\\  
$\Usalm$						&	\thref{drldef3}		\\  
$\tUarm$						&	\thref{drbm5}		\\ 

$\tUalm$					&	\thref{drlbm2}		\\ 

$\vn$						& 	\thref{drbz1}			\\

$\Vn$						&	\thref{asmbz2}			\\

$\Varm^{\sklam{s}}$			&	\thref{drbz4}			\\
$\Valm^{\sklam{s}}$			&	\thref{drlbz2}			\\


$\ysjk$						&	\thref{asmbz1}			\\
$\ysarjk$					&	\thref{aspdef1}		\\
$\ysaljk$					&	\thref{aspdef1}		\\

${\mathbb Z}$				&	\sref{Mengen}		\\
$\Z{a}{b}$					&	\sref{Mengen}		\\

$\zsjk$						&	\thref{asmbz1}			\\
$\zsarjk$					&	\thref{aspdef1}		\\
$\zsaljk$					&	\thref{aspdef1}		\\

\end{tabular}
\end{minipage}
\hfill
\begin{minipage}[t]{0.49\textwidth}
\begin{tabular}[t]{l l}

\end{tabular}
\end{minipage}

%% file: literatur.tex
\newpage
\renewcommand{\refname}{Literaturverzeichnis}
\fancyhead[RO,LE]{Literaturverzeichnis}

%% file: sm.bbl
\begin{thebibliography}{100} 

\addcontentsline{toc}{section}{Literaturverzeichnis}

\newcounter{bib}

\stepcounter{bib}
\bibitem[\thebib]{ATU1}
Adamyan, V.\,M.; Tkachenko, I.\,M.:
\textit{Truncated Hamburger moment problems with constraints}.
In: Recent progress in functional analysis (Valencia, 2000), North-Holland Math. Stud., 189, North-Holland, Amsterdam (2001), S. 321-333

\stepcounter{bib}
\bibitem[\thebib]{ATU2}
Adamyan, V.\,M.; Tkachenko, I.\,M.:
\textit{Solution of the Stieltjes truncated matrix moment problem}.
In: Opuscula Math., 25 (2005), Nr. 1, S. 5-24

\stepcounter{bib}
\bibitem[\thebib]{ATU3}
Adamyan, V.\,M.; Tkachenko, I.\,M.:
\textit{General solution of the Stieltjes truncated matrix moment problem}.
In: Operator theory and indefinite inner product spaces, 
Oper. Theory: Adv. and Appl., vol. 163, 
Birkhäuser, Basel (2006), S. 1-22

\stepcounter{bib}
\bibitem[\thebib]{ATU4}
Adamyan, V.\,M.; Tkachenko, I.\,M.; Urrea, M.:
\textit{Solution of the Stieltjes truncated moment problem}.
In: J. Appl. Anal., 9 (2003), Nr. 1, S. 57-74

\stepcounter{bib}
\bibitem[\thebib]{AK}
Akhiezer, N.\,I.; Krein, M.\,G.:
\textit{Some Questions in Theory of Moments},
(Russ.) Gos. Nauchn.-Tehn. Izd-vo. Ukr., Kharkov (1983);
Englische Übersetzung in: Translations of Mathematical Monographs, Vol. 2, Amer. Math. Soc., Providence, R.\,I. (1962)

\stepcounter{bib}
\bibitem[\thebib]{A}
Albert, A.:
\textit{Conditions for positive and nonnegative definiteness in terms of pseudoinverses},
SIAM J. Appl. Math. 17 (1969), S. 434-440


\stepcounter{bib}
\bibitem[\thebib]{Bana}
Banachiewicz, T.: 
\textit{Zur Berechnung der Determinanten, wie auch der Inversen und zur darauf basierten Auflösung der Systeme linearer Gleichungen}
In: Acta Astronom., Ser. C, 3 (1937), S. 41-67

\stepcounter{bib}
\bibitem[\thebib]{Bo1}
Bolotnikov, V.\,A.:
\textit{Descriptions of solutions of a degenerate moment problem on the axis and the halfaxis}. 
In: (Russ.) Teor. Funktsii, Funktsional. Anal. i Prilozhen., 50 (1988), S. 25-31; 
Englische Übersetzung in: J. Soviet Math., 49 (1990), Nr. 6, S. 1253-1258

\stepcounter{bib}
\bibitem[\thebib]{Bo2}
Bolotnikov, V.\,A.:
\textit{Degenerate Stieltjes moment problem and associated J-inner polynomials}.
In: Z. Anal. Anwendungen, 14 (1995), Nr. 3, S. 441-468




\stepcounter{bib}
\bibitem[\thebib]{CH3}
Chen, G.-N.; Hu, Y.-J.:
\textit{The truncated Hamburger matrix moment problems in the nondegenerate and degenerate cases, and matrix continued fractions}.
In: Linear Algebra Appl., 277 (1998), S. 199-236

\stepcounter{bib}
\bibitem[\thebib]{CH1}
Chen, G.-N.; Hu, Y.-J.:
\textit{A unified treatment for the matrix Stieltjes moment problem in both nondegenerate and degenerate cases}.
In: J. Math. Anal. Appl., 254 (2001), S. 23-34

\stepcounter{bib}
\bibitem[\thebib]{CH2}
Chen, G.-N.; Hu, Y.-J.:
\textit{Unified treatment for the matrix Stieltjes moment problem}.
In: Linear Algebra Appl., 380 (2004), S. 227-239


\stepcounter{bib}
\bibitem[\thebib]{Ch}
Chihara, T.\,S.:
\textit{An Introduction to Orthogonal Polynomials},
Gordon and Breach, New York (1978) 

\stepcounter{bib}
\bibitem[\thebib]{CR}
Choque Rivero, A.\,E.:
\textit{Ein finites Matrixmomentenproblem auf einem endlichen Intervall},
Dissertation, Universität Leipzig, Leipzig (2001)

\stepcounter{bib}
\bibitem[\thebib]{CR1}
Choque Rivero, A.\,E.:
\textit{On Dyukarev's resolvent matrix for a truncated Stieltjes matrix moment problem under the view of orthogonal matrix polynomials}.
In: Linear Algebra Appl., 474 (2015), S. 44-109

\stepcounter{bib}
\bibitem[\thebib]{C06}
Choque Rivero, A.\,E.; Dyukarev, Y.; Fritzsche, B.; Kirstein, B.: 
\textit{A Truncated Matricial Moment Problem on a Finite Interval}. 
In: Interpolation, Schur Functions and Moment Problems, 
Oper. Theory: Adv. and Appl., vol. 165, 
Birkhäuser, Basel (2006), S. 121-173

\stepcounter{bib}
\bibitem[\thebib]{C07}
Choque Rivero, A.\,E.; Dyukarev, Y.; Fritzsche, B.; Kirstein, B.: 
\textit{A Truncated Matricial Moment Problem on a Finite Interval. The Case of an Odd Number of Prescribed Moments}. 
In: System Theory, the Schur Algorithm and Multidimensional Analysis, 
Oper. Theory: Adv. and Appl., vol. 176, 
Birkhäuser, Basel (2007), S. 99-164

\stepcounter{bib}
\bibitem[\thebib]{Pert}
Choque Rivero, A.\,E.; Mädler, C.:
\textit{On Hankel Positive Definite Pertubations of Hankel Positive Definite Sequences and Interrelations of Orthogonal Matrix Polynomials}.
In: Complex Anal. Oper. Theory, 8 (2014), S. 1645-1698

\stepcounter{bib}
\bibitem[\thebib]{Du1}
Dubovoj, V.\,K:
\textit{Indefinite metric in the interpolation problem of Schur for analytic matrix functions}.
In: (Russ.) Theor. Funkcii, Funkcional. Anal. i Prilozen, 
Part I: 37 (1982), S. 14-26;
Part II: 38 (1982), S. 32-39;
Part III: 41 (1984), S. 55-64;
Part IV: 42 (1984), S. 46-57;
Part V: 45 (1986), S. 16-21;
Part VI: 47 (1987), S. 112-119

\stepcounter{bib}
\bibitem[\thebib]{Du2}
Dubovoj, V.\,K:
\textit{Parametrization of multiple elementary factor of nonfull rank}.
In: (Russ.) Analysis in Infinite Dimensional Spaces and Operator Theory (Ed.: V.\,A. Marcenko),
Naukova Dumka, Kiev (1983), S. 54-68

\stepcounter{bib}
\bibitem[\thebib]{DFK}
Dubovoj, V.\,K.; Fritzsche, B.; Kirstein, B.: 
\textit{Matricial version of the classical Schur Problem}. 
In: Teubner-Texte zur Mathematik, Band 129, B. G. Teubner, Stuttgart-Leipzig (1992)

\stepcounter{bib}
\bibitem[\thebib]{dyu2}
Dyukarev, Yu.\,M.:
\textit{The Stieltjes matrix moment problem},
Manuskript, deponiert in VINITI (Moskau) am 22.3.1981, Nr. 2628-81 (1981)

\stepcounter{bib}
\bibitem[\thebib]{dyu1}
Dyukarev, Yu.\,M.:
\textit{Interpolation problems in the Stieltjes class},
Dissertation, Universität im. A.\,M. Gorkogo, Kharkov (1982)

\stepcounter{bib}
\bibitem[\thebib]{dyu3}
Dyukarev, Yu.\,M.:
\textit{A general scheme for solving interpolation problems in the Stieltjes class that is based on consistant integral representations of pairs of nonnegativ operators},
Mat. Fiz. Anal. Geom. 6, Nr. 1-2 (1999), S. 30-54

\stepcounter{bib}
\bibitem[\thebib]{dyu}
Dyukarev, Yu.\,M.:
\textit{Indeterminacy criteria for the Stieltjes matrix moment problem}.
In: (Russ.) Mat. Zametki, 75 (2004), S. 71-88;
Englische Übersetzung in: Math. Notes, 75 (2004), S. 66-82

\stepcounter{bib}
\bibitem[\thebib]{dyu4}
Dyukarev, Yu.\,M.:
\textit{Theory of interpolation problems in the Stieltjes class and some related questions in analysis},
Habilitation, Universität Charkiw, Charkiw (2006)

\stepcounter{bib}
\bibitem[\thebib]{13}
Dyukarev, Yu.\,M.; Fritzsche, B.; Kirstein, B.; Mädler, C.; Thiele, H.\,C.:
\textit{On distinguished solutions of truncated matricial Hamburger moment problems}.
In: Complex Anal. Oper. Theory, 3 (2009), Nr. 4, S. 759-834

\stepcounter{bib}
\bibitem[\thebib]{12}
Dyukarev, Yu.\,M.; Fritzsche, B.; Kirstein, B.; Mädler, C.:
\textit{On truncated matricial Stieltjes type moment problems}. 
In: Complex Anal. Oper. Theory, 4 (2010), Nr. 4, S. 905-951

\stepcounter{bib}
\bibitem[\thebib]{DK}
Dyukarev, Yu.\,M.; Katsnelson, V.\,E.:
\textit{Multiplicative and additive classes
of analytic matrix functions of Stieljes type and associated interpolation
problems connected with them}.
In: (Russ.) Teor. Funktsii, Funktsional. Anal. i Prilozhen., Teil I: 36, (1981), S. 13-27; Teil II: 38 (1982), S. 40-48; Teil III: 41 (1984), S. 64-70

\stepcounter{bib}
\bibitem[\thebib]{EP}
Efimov, A.\,V.; Potapov, V.\,P.: 
\textit{J-expansive matrix-valued functions and their role in the analytic theory of electrical circuits}. 
In: (Russ.) Uspekhi Mat. Nauk, 28 (1973), Nr. 1 (169), S. 65-130; 
Englische Übersetzung in: Russ. Math. Surv., 28 (1973), Nr. 1, S. 69-140

\stepcounter{bib}
\bibitem[\thebib]{FR}
Friese, N.; Reichel, M.:
\textit{Matrixgleichungen im Zusammenhang mit Potenzmomentenproblemen vom Stieltjes und Hausdorff-Typ},
Diplomarbeit, Universität Leipzig, Leipzig (2015)

\stepcounter{bib}
\bibitem[\thebib]{14}
Fritzsche, B.; Kirstein, B.:
\textit{Schwache Konvergenz nichtnegativ hermitescher Borelmaße}, 
Wiss. Z. Karl-Marx-Univ. Leipzig Math.-Natur., 37 (1988), Nr. 4, S. 375-398

\stepcounter{bib}
\bibitem[\thebib]{15}
Fritzsche, B.; Kirstein, B.; Mädler, C.:
\textit{On Hankel nonnegative definite sequences, the canonical Hankel parametrization, and orthogonal matrix polynomials}. 
In: Complex Anal. Oper. Theory, 5 (2011), Nr. 2, S. 447-511

\stepcounter{bib}
\bibitem[\thebib]{ot226}
Fritzsche, B.; Kirstein, B.; Mädler, C.:
\textit{On a Special Parametrization of Matricial $\alpha$-Stieltjes One-sided Non-negative Definite Sequences}
In: Interpolation, Schur Functions and Moment Problems II, 
Oper. Theory: Adv. and Appl., vol. 226, 
Birkhäuser/Springer Basel AG, Basel (2012), S. 211-250

\stepcounter{bib}
\bibitem[\thebib]{Trans}
Fritzsche, B.; Kirstein, B.; Mädler, C.:
\textit{Transformations of matricial $\alpha$-Stieltjes non-negative definite sequences}.
In: Linear Algebra and its Applications, 439 (2013), S. 3893-3933

\stepcounter{bib}
\bibitem[\thebib]{SimH}
Fritzsche, B.; Kirstein, B.; Mädler, C.:
\textit{On a simultaneous approach to the even and odd truncated matricial Hamburger moment problems}.
In: Recent Advances in Inverse Scattering, Schur Analysis and Stochastic Processes, 
Oper. Theory: Adv. and Appl., vol. 244, 
Birkhäuser/Springer Basel AG, Basel (2015), S. 181-285

\stepcounter{bib}
\bibitem[\thebib]{SF}
Fritzsche, B.; Kirstein, B.; Mädler, C.:
\textit{On matrix-valued Stieltjes functions with an emphasis on particular subclasses}.
In: arXiv:1506.01600v1 [math.CV] (2015); Erscheint in: Large Truncated Toeplitz Matrices, Toeplitz Operators, and Related Topics, Oper. Theory: Adv. and Appl., vol. 259, Birkhäuser/Springer Basel AG, Basel (2017)

\stepcounter{bib}
\bibitem[\thebib]{Sim}
Fritzsche, B.; Kirstein, B.; Mädler, C.:
\textit{On a Simultaneous Approach to the Even and Odd Truncated Matricial Stieltjes Moment Problem}.
In: Part I, arXiv:1604.07240v1 [math.CV] (2016); Part II, arXiv:1604.07629v1 [math.CV] (2016);
Eingereicht in: Linear Algebra and its Applications

\stepcounter{bib}
\bibitem[\thebib]{HMS}
Fritzsche, B.; Kirstein, B.; Mädler, C.:
\textit{On the structure of Hausdorff moment sequences of complex matrices}.
In: arXiv:1701.04246v1 [math.CV] (2017)

\stepcounter{bib}
\bibitem[\thebib]{GL}
Gerecke, U.; Lorenz, J.:
\textit{Grundlegende Aussagen über nichtnegativ hermitesche Maße, Maße mit orthogonalen Werten sowie projektorwertige Maße},
Diplomarbeit, Universität Leipzig, Leipzig (1996)

\stepcounter{bib}
\bibitem[\thebib]{GK}
Gantmacher, F.\,R.; Krein, M.\,G.:
\textit{Oszillationsmatrizen, Oszillationskerne und kleine Schwingungen mechanischer Systeme},
Wisssenschaftliche Bearbeitung der deutschen Ausgabe: Alfred Stöhr, Math. Lehrbüch. Monogr. Abt. 1, Band V, Akademie-Verlag, Berlin (1960);
Englische Version: \textit{Oscillation Matrices and Kernels and Small Vibrations of Mechanical Systems}, revised ed., AMS Chelsea Publishing, Providence, RI (2002), Übersetzung basiert auf dem 1941 erschienenen Original in Russisch, bearbeitet von Alex Eremenko

\stepcounter{bib}
\bibitem[\thebib]{Ham}
Hamburger, H.\,L.:
\textit{Über eine Erweiterung des Stieltjesschen Momentenproblems}.
In: Math. Ann., Teil I: 81 (1920), S. 235-319; Teil II: 82 (1921), S. 120-164; Teil III: 82 (1921), S. 168-187

\stepcounter{bib}
\bibitem[\thebib]{HJ}
Horn, R.\,A.; Johnson, C.\,R.:
\textit{Matrix analysis},
second edition, Cambridge University Press, Cambridge (2013)


\stepcounter{bib}
\bibitem[\thebib]{Jung}
Junghanns, P.:
\textit{EAGLE-GUIDE Orthogonale Polynome},
Edition am Gutenbergplatz, Leipzig (2009)

\stepcounter{bib}
\bibitem[\thebib]{KS}
Karlin, S.; Studden, W.\,J.:
\textit{Tchebycheff systems: with applications in analysis and statistics},
Interscience Publishers, New-London-Sydney (1966)

\stepcounter{bib}
\bibitem[\thebib]{Ka}
Kats, I.\,S.:
\textit{On Hilbert spaces generated by Hermitian monotone matrix functions}.
In: (Russ.) Zap. Mat. Otd. Fiz.-Mat. Fak. i Kharkov. Mat. Obsh., 22 (1950), S. 95-113

\stepcounter{bib}
\bibitem[\thebib]{Kat1}
Katsnelson, V.\,E.:
\textit{Continual analogues of the Hamburger-Nevanlinna theorem and fundamental matrix inequalities of classical problems}.
In: (Russ.) Teor. Funkcii, Funkcional. Anal. i Prilozen,
Part I: 36 (1981), S. 31-48;
Part II: 37 (1982), S. 31-48;
Part III: 39 (1983), S. 61-73;
Part IV: 40 (1983), S. 79-90

\stepcounter{bib}
\bibitem[\thebib]{Ka1}
Katsnelson, V.\,E.:
\textit{Methods of J-theory in continuous interpolation problems of analysis},
(Russ.) hinterlegt in VINITI (1983).
Englische Übersetzung des Part I von T. Ando, Hokkaido University, Sapporo (1985)

\stepcounter{bib}
\bibitem[\thebib]{Kat2}
Katsnelson, V.\,E.:
\textit{Integral representation of Hemitian positive kernels of mixed type and the generalized Nehari problem}.
In: (Russ.) Teor. Funkcii, Funkcional. Anal. i Prilozen, 43 (1985), S. 54-70

\stepcounter{bib}
\bibitem[\thebib]{Kat3}
Katsnelson, V.\,E.:
\textit{Extremal and factorization properties of the radii in the problem of representations of positive Hermitian matrix functions}.
In: (Russ.) Mathematical Physics and Functional Analysis (Ed.: V. A. Marcenko),
Naukova Dumka, Kiev (1986), S. 80-94

\stepcounter{bib}
\bibitem[\thebib]{Ka2}
Katsnelson, V.\,E.:
\textit{On transformations of Potapov's fundamental matrix inequality}.
In: Topics in interpolation theory, 
Oper. Theory: Adv. and Appl., vol. 95, 
Birkhäuser, Basel (1997), S. 253-281

\stepcounter{bib}
\bibitem[\thebib]{Ko1}
Kovalishina, I.\,V.:
\textit{J-expansive matrix-valued functions and the classical moment problem}.
In: (Russ.) Akad. Nauk Armjan. SSR Dokl., 60 (1975), Nr. 1, S. 3-10

\stepcounter{bib}
\bibitem[\thebib]{Ko2}
Kovalishina, I.\,V.:
\textit{Analytic theory of a class of interpolation problems}.
In: (Russ.) Izv. Akad. Nauk SSSR Ser. Mat., 47 (1983), Nr. 3, S. 455-497

\stepcounter{bib}
\bibitem[\thebib]{Kr1}
Krein, M.\,G.:
\textit{On a generalized problem of moments}.
In: (Russ.) Dokl. Akad. Nauk SSSR (N.S.), 44 (1944), S. 219-222

\stepcounter{bib}
\bibitem[\thebib]{Kr2}
Krein, M.\,G.:
\textit{Infinite J-matrices and a matrix-moment problem}.
In: (Russ.) Dokl. Akad. Nauk SSSR (N.S.), 69 (1949), S. 125-128

\stepcounter{bib}
\bibitem[\thebib]{Kr3}
Krein, M.\,G.:
\textit{The fundamental propositions of the theory of representations
of Hermitian operators with deficiency index (m,m)}.
In: (Russ.) Ukrain. Mat. Zurnal, 1 (1949), no. 2, S. 3-66

\stepcounter{bib}
\bibitem[\thebib]{Kr5}
Krein, M.\,G.:
\textit{The ideas of P. L. \v{C}eby\v{s}ev and A. A. Markov in the theory of limited values of integrals and their further development}.
In: (Russ.) Usp. Mat. Nauk 6 (1951), no. 4, S. 3-66 (with the redactional participation of P. G. Rekhtman);
Englische Übersetzung: Amer. Math. Soc. Transl., Series 2, 12 (1959), S. 1-121

\stepcounter{bib}
\bibitem[\thebib]{Kr4}
Krein, M.\,G.; Krasnoselskii, M.\,A.:
\textit{Fundamental theorems on the extension of Hermitian operators and certain of their applications to the theory of orthogonal polynomials and the problem of moments}.
In: (Russ.) Uspekhi Mat. Nauk (N.S.), 2 (1947), Nr. 3 (19), S. 60-106

\stepcounter{bib}
\bibitem[\thebib]{KN}
Krein, M.\,G.; Nudelman, A.\,A.:
\textit{The Markov Moment Problem and Extremal Problems},
(Russ.) Nauka, Moscow (1973);
Englische Übersetzung in: Translations of Mathematical Monographs, Vol. 50, Amer. Math. Soc., Providence, R. I. (1977)

\stepcounter{bib}
\bibitem[\thebib]{Maka2}
Makarevich, T.:
\textit{Darstellungen der Lösungsmengen finiter matrizieller
Hamburgerscher Potenzmomentenprobleme}, 
Diplomarbeit, Universität Leipzig, Leipzig (2009)

\stepcounter{bib}
\bibitem[\thebib]{Maka}
Makarevich, T.:
\textit{Ein matrizielles finites Momentenproblem vom Stieltjes-Typ}, 
Dissertation, Universität Leipzig, Leipzig (2014)

\stepcounter{bib}
\bibitem[\thebib]{MP}
Mühling, Ch.; Pfeffing, C.:
\textit{Über ein matrizielles Potenzmomentenproblem vom Stieltjes-Typ},
Diplomarbeit, Universität Leipzig, Leipzig (2011)

\stepcounter{bib}
\bibitem[\thebib]{Nev}
Nevanlinna, R.:
\textit{Asymptotische Entwicklungen beschränkter Funktionen und das Stieltjessche Momentenproblem}.
In: Ann. Acad. Sci. Fenn., A 18 (1922), 5, S. 1-53

\stepcounter{bib}
\bibitem[\thebib]{Nm1}
Nudelman, A.\,A.:
\textit{The work of M. G. Krein on the moment problem}.
In: (Russ.) Ukrain. Mat. Zh., 46 (1994), Nr. 1-2, S. 62-74

\stepcounter{bib}
\bibitem[\thebib]{Nm2}
Nudelman, A.\,A.:
\textit{On M. G. Krein's contribution to the moment problem}.
In: Operator theory and related topics II (Odessa, 1997), Oper. Theory Adv. Appl., 118, Birkhäuser, Basel (2000), S. 17-32

\stepcounter{bib}
\bibitem[\thebib]{Pe}
Petzel, O.:
\textit{Beiträge zur Integrationstheorie nichtnegativ hermitescher Maße mit Anwendungen zu matriziellen Momentenproblemen}, 
Diplomarbeit, Universität Leipzig, Leipzig (2008)

\stepcounter{bib}
\bibitem[\thebib]{Po}
Potapov, V.\,P.: 
\textit{The multiplicative structure of J-contractive matrix functions}.
In: (Russ.) Tr. Mosk. Mat. Obs., 1955, Volume 4, S. 125-236

\stepcounter{bib}
\bibitem[\thebib]{Funk}
Remmert, R.; Schumacher, G.:
\textit{Funktionentheorie 1},
fünfte Auflage, Springer Verlag, Berlin (2002)

\stepcounter{bib}
\bibitem[\thebib]{Ro}
Rosenberg, M.:
\textit{The square integrability of matrix-valued functions with respect to a non-negative Hermitian measure}.
In: Duke Math. J., 31 (1964), Nr. 2, S. 291-298


\stepcounter{bib}
\bibitem[\thebib]{Sh}
Scheithauer, P.:
\textit{Über ein matrizielles Potenzmomentenproblem vom Stieltjes-Typ},
Diplomarbeit, Universität Leipzig, Leipzig (2011)

\stepcounter{bib}
\bibitem[\thebib]{Sch}
Schröder, T.:
\textit{Beweise für notwendige und hinreichende Lösbarkeitsbedingungen finiter matrizieller Hamburgerscher Potenzmomentenprobleme}, 
Diplomarbeit, Universität Leipzig, Leipzig (2013)

\stepcounter{bib}
\bibitem[\thebib]{Sur}
Stieltjes, T.-J.:
\textit{Recherches sur les fractions continues},
(Franz.) Annales de la Faculté des sciences de Toulouse: Mathématiques, vol. 8 iss. 4 (1894), vol. 9 iss. 1 (1895)

\stepcounter{bib}
\bibitem[\thebib]{Wb}
Wöhlbier, P.-E.:
\textit{Einige Beiträge zu einem matriziellen Momentenproblem vom Stieltjes-Typ},
Diplomarbeit, Universität Leipzig, Leipzig (2015)

\end{thebibliography}
